\documentclass[11pt]{amsart}
\usepackage[T1]{fontenc}
\usepackage[utf8]{inputenc}
\usepackage{MnSymbol} %For llangle, rrangle
\DeclareMathSymbol{\subsetneq} {\mathrel}{AMSb}{"28}
\DeclareMathSymbol{\supsetneq} {\mathrel}{AMSb}{"29}

\usepackage{dsfont}
\usepackage[noadjust]{cite}
\usepackage{booktabs}
\usepackage{hyphenat}
\usepackage{mathtools}
\usepackage{mdwtab}
\usepackage{enumitem}
\usepackage{color}
\usepackage{soul}
\usepackage{mathrsfs}
\usepackage[pdftex, lmargin=1in, rmargin=1in, tmargin=1in, bmargin=1in]{geometry}
\usepackage[bookmarks=true, bookmarksopen=true,%
bookmarksdepth=3,bookmarksopenlevel=2,%
colorlinks=true,%
linkcolor=blue,%
citecolor=blue,%
filecolor=blue,%
menucolor=blue,%
urlcolor=blue]{hyperref}
\hypersetup{pdftitle={Diagonals and A-infinity tensor products}}
\hypersetup{pdfauthor={Robert Lipshitz, Peter S. Ozsváth and Dylan
    P. Thurston}}
\usepackage{letltxmacro}
\usepackage{color}
\usepackage{tikz}
\usetikzlibrary{matrix,arrows,calc}
\newcommand\tikzsetnextfilename[1]{}
\tikzset{cdlabel/.style={above,sloped,
   execute at begin node=$\scriptstyle,execute at end node=$}}
\tikzset{alga/.style={->, thick}}%Alg arrow
\tikzset{stumpa/.style={|->, thick}}%Alg arrow
\tikzset{blga/.style={->, thick}}%Blg arrow
\tikzset{taa/.style={->, double}}
\tikzset{moda/.style={->, dashed}}
\tikzset{dmoda/.style={->, dashed}}
\tikzset{damoda/.style={->, dashed}}
\tikzset{tbb/.style={->,double}} %Tensor on Blg.
\tikzset{smallpic/.style={x=0.75cm,y=0.75cm,node font=\small}}
\tikzset{smallpicwide/.style={x=1.25cm,y=0.75cm,node font=\small}}

\allowdisplaybreaks

\newcommand{\RR}{\mathbb R}
\newcommand{\CC}{\mathbb C}
\newcommand{\ZZ}{\mathbb Z}

\newcommand{\FF}{\mathbb F}
\newcommand{\NN}{\mathbb N}

\newcommand{\llbracket}{\lsem}
\newcommand{\rrbracket}{\rsem}

\newcommand\HHH{\mathbb{H}}

\newcommand{\co}{\nobreak\mskip2mu\mathpunct{}\nonscript
  \mkern-\thinmuskip{:}\penalty300\mskip6muplus1mu\relax}

\newcommand{\into}{\hookrightarrow}

\newcommand{\bdy}{\partial}

\newcommand{\lbracket}{[}
\newcommand{\rbracket}{]}

\DeclareMathOperator{\Sym}{Sym}

\DeclareMathOperator{\Hom}{Hom}
\DeclareMathOperator{\Aut}{Aut}

\DeclareMathOperator{\gr}{gr}

\DeclareMathOperator{\wgr}{wt}

\newcommand{\Barop}{{\mathrm{Bar}}}

\theoremstyle{plain}

\numberwithin{equation}{section}
\newtheorem{theorem}[equation]{Theorem}

\newtheorem{proposition}[equation]{Proposition}
\newtheorem{lemma}[equation]{Lemma}
\newtheorem{corollary}[equation]{Corollary}

\newtheorem{convention}[equation]{Convention}
\newtheorem{definition}[equation]{Definition}
\newtheorem{construction}[equation]{Construction}

\theoremstyle{definition}

\theoremstyle{remark}
\newtheorem{example}[equation]{Example}
\newtheorem{remark}[equation]{Remark}
\newtheorem{warning}[equation]{Warning}

\newcommand{\HF}{\mathit{HF}}
\newcommand{\HFa}{\widehat {\HF}}

\newcommand{\HFm}{{\HF}^-}

\newcommand{\z}{\mathbf z}
\newcommand{\w}{\mathbf w}

\newcommand\HH{\mathit{HH}}

\newcommand\Hochschild\HH

\newcommand{\grs}[1]{\llangle{#1}\rrangle} %Notation for a grading shift

\newcommand{\Ainf}{A_\infty}
\newcommand{\Alg}{\mathcal{A}}

\newcommand\Blg{\mathcal{B}}
\newcommand\Clg{\mathcal{C}}

\newcommand{\cA}{{\mathcal{A}}}
\newcommand{\cB}{{\mathcal{B}}}

\newcommand{\alphas}{{\boldsymbol{\alpha}}}
\newcommand{\betas}{{\boldsymbol{\beta}}}
\newcommand{\gammas}{{\boldsymbol{\gamma}}}

\newcommand{\cM}{\mathcal{M}}

\newcommand{\DD}{\textit{DD}}
\newcommand{\DA}{\textit{DA}}

\newcommand{\AAm}{\textit{AA}} %\AA conflicts with something
\newcommand{\CFD}{\mathit{CFD}}

\newcommand{\CFA}{\mathit{CFA}}

\newcommand{\CFDa}{\widehat{\CFD}}

\newcommand{\CFAa}{\widehat{\CFA}}

\newcommand{\dg}{\textit{dg} }

\newcommand{\cModule}{\mathcal{M}}
\newcommand{\cNodule}{\mathcal{N}}

\newcommand\Id{\mathbb{I}}
\newcommand\Ground{\mathds{k}}
\newcommand\DT{\boxtimes}

\newcommand\Tensor{\mathcal T}

\newcommand{\Field}{{\FF_2}}

\newcommand{\Ring}{R}

\renewcommand{\th}{^\text{th}}
\renewcommand{\st}{^\text{st}}

\newcommand{\ModCat}{\mathsf{Mod}}

\newcommand{\Cat}{\mathscr{C}}
\newcommand{\Dat}{\mathscr{D}}

\DeclareMathOperator{\Mor}{Mor}

\DeclareMathOperator{\ob}{Ob}
\DeclareMathOperator\suMor{\Mor^{u}}
\DeclareMathOperator\wMor{wMor}
\DeclareMathOperator\uwMor{uMor}

\newcommand{\op}{\mathrm{op}}

\makeatletter
\newcommand\honestalg[3]{\bigl\lbracket
\begin{smallmatrix} #1\@ifempty{#3}{}{&#3} \\ #2 \end{smallmatrix}
\bigr\rbracket}
\makeatother
\newcommand{\lab}[1]{$\scriptstyle #1$}

\newcommand{\lsub}[2]{{}_{#1}#2}
\newcommand{\lsup}[2]{{}^{#1}\mskip-.6\thinmuskip#2}

\newcommand{\cell}{\mathrm{cell}}

\newcommand\unit{\mathbf 1}

\newcommand{\corolla}[1]{\Psi_{#1}}%Fairly pretty corolla
\newcommand{\wcorolla}[2]{\corolla{#1}^{#2}}

\newcommand{\DegenTree}{\mathord{\downarrow}}
\newcommand{\IdTree}{\DegenTree} %was $T_\Id$
\newcommand{\trunkvertex}{\tau}

\newcommand{\kotimes}[1]{\otimes}%To allow us to suppress ground field from certain tensor products.
\newcommand{\rotimes}[1]{\otimes_{#1}}%To allow us to suppress ground field from certain tensor products.

\newcommand{\AsDiag}{{\boldsymbol{\Gamma}}} %A diagonal for the associahedron.
\newcommand{\TrDiag}{\boldsymbol{\gamma}}
\newcommand{\wTrDiag}[2]{\boldsymbol{\gamma}^{#1,#2}}
\newcommand{\MDiag}{\mathbf{M}}
\newcommand{\wADiag}{\boldsymbol{\Gamma}}
\newcommand{\wADiagCell}{\boldsymbol{\gamma}}
\newcommand{\wMDiag}{\mathbf{M}}
\newcommand{\TrMPrim}{\mathbf{p}}
\newcommand{\TrMDiag}{\mathbf{m}}

\newcommand{\TrPMDiag}{\mathbf{p}}
\newcommand{\TrPMorDiag}{\mathbf{q}}
\newcommand{\TrDADD}{\mathbf{r}}
\newcommand{\MDtp}[1][\MDiag]{\otimes_{#1}}
\newcommand{\ADtp}[1][\AsDiag]{\otimes_{#1}}
\newcommand{\wADtp}[1][\wADiag]{\otimes_{#1}}
\newcommand{\wMDtp}[1][\MDiag]{\otimes_{#1}}

\newcommand{\sing}{\mathrm{sing}}
\newcommand{\Trees}{\mathcal{T}}
\newcommand{\ModTransTrees}[1]{{\mathcal{T}^{MT}_{#1}}}
\newcommand{\gModTransTrees}[1]{{g\mathcal{T}^{MT}_{#1}}}
\newcommand{\MulDiag}{{\boldsymbol{\Theta}}}
\newcommand{\TrMulDiag}{\boldsymbol{\theta}}
\newcommand{\wTrMulDiag}[2]{\boldsymbol{\theta}^{#1,#2}}
\newcommand{\wTrDADD}[2]{\mathbf{r}^{#1,#2}}

\newcommand{\ModMulDiag}{\mathbf{L}}
\newcommand{\PartTrModMulDiag}{\mathbf{k}}
\newcommand{\TrModMulDiag}{\boldsymbol{\ell}}
\newcommand{\wDiag}[2]{\wADiag^{#1,#2}}
\newcommand{\wDiagCell}[2]{\wADiagCell^{#1,#2}}
\newcommand{\wDiagNS}{\wADiag}
\newcommand{\wModDiag}[2]{\mathbf{M}^{#1,#2}}
\newcommand{\wModDiagNS}{\mathbf{M}}
\newcommand{\wModDiagCell}[2]{\mathbf{m}^{#1,#2}}
\newcommand{\wTrPMDiag}[2]{\mathbf{p}^{#1,#2}}
\newcommand{\wTrPMDiagNS}{\mathbf{p}}
\newcommand{\stump}{\top}
\newcommand{\wMulDiag}[2]{\boldsymbol{\Theta}^{#1,#2}}
\newcommand{\wMulDiagNS}{\boldsymbol{\Theta}}
\newcommand{\wModMulDiag}[2]{\mathbf{L}^{#1,#2}}
\newcommand{\wModMulDiagNS}{\mathbf{L}}
\newcommand{\wTrModMulDiag}[2]{\boldsymbol{\ell}^{#1,#2}}

\newcommand{\wTrPMorDiag}[2]{\mathbf{q}^{#1,#2}}
\newcommand{\wTrPMorDiagNS}{\mathbf{q}}
\newcommand{\wPartTrModMulDiag}[2]{\mathbf{k}^{#1,#2}}
\newcommand{\wPartTrModMulDiagNS}{\mathbf{k}}
\newcommand{\uwDiag}[2]{u\wADiag^{#1,#2}}
\newcommand{\uwMDiag}[2]{u\wMDiag^{#1,#2}}
\newcommand{\uwMulDiag}[2]{u\boldsymbol{\Theta}^{#1,#2}}
\newcommand{\uwModMulDiag}[2]{u\mathbf{L}^{#1,#2}}
\newcommand{\uwDiagNS}{u\wADiag}
\newcommand{\uwMDiagNS}{u\wMDiag}
\newcommand{\uwMulDiagNS}{u\boldsymbol{\Theta}}
\newcommand{\uwModMulDiagNS}{u\mathbf{L}}

\DeclareMathOperator{\RootJoin}{RoJ} %Root joining
\DeclareMathOperator{\LeftJoin}{LeJ} %Left joining
\DeclareMathOperator{\LRjoin}{LR} %Left join on left, root join on right.
\newcommand\LRjoinW{\LRjoin} %weighted version of LRjoin

\DeclareMathOperator{\Spl}{Spl} %Splitting a module transformation tree

\DeclareMathOperator{\wRootJoin}{\RootJoin} %Root joining

\newcommand\Vertices{\mathrm{Vert}}
\newcommand\Edges{\mathrm{Edge}}

\newcommand{\Filt}{\mathcal{F}}

\newcommand{\pt}{\mathit{pt}}
\newcommand{\cellC}[1]{C^{\mathit{cell}}_{#1}}
\newcommand{\singC}[1]{C^{\mathit{sing}}_{#1}}

\newcommand{\wAlg}{\mathscr A}
\newcommand{\wBlg}{\mathscr B}
\newcommand{\wClg}{\mathscr C}
\newcommand{\wBlgop}{{\mathscr B}^{\op}}

\newcommand{\wMod}{\mathscr M}
\newcommand{\wNod}{\mathscr N}
\newcommand{\wP}{P}
\newcommand{\wQ}{Q}
\newcommand{\One}{\boldsymbol{1}}

\newcommand\dpm{\mathrm{dpm}}

\newcommand{\uTreesCx}[2][*]{u\!X_{#1}^{#2}}
\newcommand{\uTransCx}[2][*]{u\!J_{#1}^{#2}}
\newcommand{\wTreesCx}[3][*]{X_{#1}^{#2,#3}}
\newcommand{\wTransCx}[3][*]{J_{#1}^{#2,#3}}
\newcommand{\xwTransCx}[3][*]{\widetilde{J}_{#1}^{#2,#3}}
\newcommand{\wMTransCx}[3][*]{J\!M_{#1}^{#2,#3}}
\newcommand{\gwMTransCx}[3][*]{gJ\!M_{#1}^{#2,#3}}
\newcommand{\wMTreesCx}[3][*]{X\!M_{#1}^{#2,#3}}
\newcommand{\uwTreesCx}[3][*]{u\!X_{#1}^{#2,#3}}
\newcommand{\uwMTreesCx}[3][*]{u\!X\!M_{#1}^{#2,#3}}
\newcommand{\uwTransCx}[3][*]{u\!J_{#1}^{#2,#3}}
\newcommand{\uwMTransCx}[3][*]{u\!J\!M_{#1}^{#2,#3}}
\newcommand{\wTrees}[2]{\mathcal{T}_{#1,#2}}
\newcommand{\wMTrees}[2]{\mathcal{T\!M}_{#1,#2}}
\newcommand{\xwTreesCx}[3][*]{\widetilde{X}_{#1}^{#2,#3}}
\newcommand{\xwTrees}[2]{\widetilde{\mathcal{T}}_{#1,#2}}
\newcommand{\xwMTreesCx}[2]{\widetilde{XM}_*^{#1,#2}}
\newcommand{\wSeed}{\mathscr{S}}
\newcommand{\rcorolla}[1]{\textcolor{red}{\mathrm{r}\hspace{-1pt}\Psi_{#1}}}
\newcommand{\bcorolla}[1]{\textcolor{blue}{\mathrm{b}\hspace{-1pt}\Psi_{#1}}} 
\newcommand{\pcorolla}[1]{\textcolor{purple}{\mathrm{p}\hspace{-1pt}\Psi_{#1}}}
\newcommand{\bRootJoin}{\textcolor{blue}{\mathrm{bRoJ}}}
\newcommand{\wrcorolla}[2]{\textcolor{red}{\mathrm{r}\hspace{-1pt}\Psi_{#1}^{#2}}}
\newcommand{\wbcorolla}[2]{\textcolor{blue}{\mathrm{b}\hspace{-1pt}\Psi_{#1}^{#2}}} 
\newcommand{\wpcorolla}[2]{\textcolor{purple}{\mathrm{p}\hspace{-1pt}\Psi_{#1}^{#2}}}
\newcommand{\Forget}{\mathcal{F}}

\newcommand{\dfs}{\mathrm{dfs}}
\newcommand{\depth}{\mathrm{dfsind}}

\newcommand{\signissue}{ }

\DeclareMathOperator{\valence}{val}

\newcommand{\seq}[2][j]{\{#2\}_{#1=1}^\infty}

\newread\testin

\graphicspath{{draws/}{mpdraws/}{}}
\makeatletter
\def\input@path{{}{draws/}}
\makeatother

\def\mathcenter#1{%
  \vcenter{\hbox{$#1$}}%
}

\makeatletter
\newcommand\mi@kern[1]{%
  \settowidth\@tempdima{$\mi@obj^{#1}$}
  \kern-\@tempdima
  #1
  \settowidth\@tempdima{$\mi@obj$}
  \kern\@tempdima
}

\newtoks\mi@toksp
\newtoks\mi@toksb
\DeclareRobustCommand{\manyindices}[5]{
  \def\mi@obj{#5}
  \mi@toksp\expandafter{\mi@kern{#2}}
  \mi@toksb\expandafter{\mi@kern{#1}}
  \@mathmeasure4\textstyle{#5_{#1}^{#2}}
  \@mathmeasure6\textstyle{#5_{#3}^{#4}}
  \dimen0-\wd6 \advance\dimen0\wd4
  \@mathmeasure8\textstyle{\hphantom{{}_{#1}^{#2}}#5^{\the\mi@toksp#4}_{\the\mi@toksb#3}}
  \hbox to \dimen0{}{\kern-\dimen0\box8}
}
\makeatother 

\usepackage{ifpdf}
\ifpdf
  \let\textalt\texorpdfstring
\else
  \newcommand{\textalt}[2]{#1}
\fi

\begin{document}
\title{Diagonals and A-infinity Tensor Products}

\author[Lipshitz]{Robert Lipshitz} \thanks{\texttt{RL\ \ was supported by
    NSF grant\ \ DMS-1810893.}}
\address{Department of Mathematics, University of Oregon\\
  Eugene, OR 97403} \email{lipshitz@uoregon.edu}

\author[Ozsv\'ath]{Peter Ozsv\'ath}
\thanks{\texttt{PSO was supported by NSF grants DMS-1708284 and DMS-2104536.}}
\address {Department of Mathematics, Princeton University\\ New
  Jersey, 08544}
\email {petero@math.princeton.edu}

\author[Thurston]{Dylan~P.~Thurston}
\thanks{\texttt{DPT was supported by NSF grants DMS-1507244 and DMS-2110143.}}
\address{Department of Mathematics\\
         Indiana University,
         Bloomington, Indiana 47405\\
         USA}
\email{dpthurst@indiana.edu}

\begin{abstract}
  Extending work of Saneblidze-Umble and others, we use diagonals for
  the associahedron and multiplihedron to define tensor products of
  $\Ainf$-algebras, modules, algebra homomorphisms, and module
  morphisms, as well as to define a bimodule analogue of twisted
  complexes (type \DD\
  structures, in the language of bordered Heegaard Floer homology) and
  their one- and two-sided tensor products. We then give analogous
  definitions for 1-parameter deformations of $\Ainf$-algebras; this
  involves another collection of complexes. These constructions are
  relevant to bordered Heegaard Floer homology.
\end{abstract} 

\date{\today}

% Floer homology

\maketitle 
%Post to Math.RA (rings and algebras), cross-list to math.GT
%(Geometric Topology), math.SG (Symplectic Geometry)

\tableofcontents

% \newpage
% \listoftodos

\newpage
\listoffigures

\section{Introduction}

\subsection{Context and motivation}

This paper develops the homological algebra underpinning construction
of bordered Floer homology for the full Heegaard Floer homology
package~\cite{OS04:HolomorphicDisks}. While the paper is purely about
homological algebra, we start with a brief discussion of the
motivation from topology. This discussion is not needed for the rest
of the paper.

Bordered Floer homology~\cite{LOT1}
associates a \dg algebra $\Alg(F)$ to a surface $F$, an $\Ainf$-module
$\CFAa(Y_1)$ over $\Alg(F)$ to a three-manifold $Y_1$ equipped with an
identification $\partial Y_1=F$, and a twisted complex, or type $D$
structure, $\CFDa(Y_2)$ over $\Alg(F)$ to a three-manifold $Y_2$ with
$\partial Y_2=-F$.  A pairing theorem then describes the
$U=0$-specialized Heegaard Floer complex of the closed three-manifold
$Y=Y_1\cup_F Y_2$, as a certain type of tensor product of $\CFAa(Y_1)$
with $\CFDa(Y_2)$.  Moreover, the modules $\CFAa(Y_2)$ and
$\CFDa(Y_2)$ are related by tensoring with certain bimodules, called
the type \AAm\ and \DD\ bimodules of the identity map;
see~\cite{LOT2}.

The $U=0$ specialization of Heegaard Floer homology, denoted $\HFa$,
is useful for certain three-dimensional applications, such as
detecting the Thurston norm of the
three-manifold~\cite{OS04:ThurstonNorm,Ni13:spheres}; it is, however,
insufficient to access many the four-dimensional aspects, which were
the original motivation for the construction of the theory;
see~\cite{OS06:HolDiskFour}. To gain access to the full power of
Heegaard Floer homology, one must use  the $U$-unspecialized
version of Heegaard Floer homology, $\HFm$.

Extending the bordered theory to $\HFm$ presents new challenges,
arising from the existence of holomorphic disks which are now allowed to cover
the entire Heegaard surface. This phenomenon has both
geometric and algebraic manifestations. 

The aim of this work is to treat the algebraic aspects. In particular,
the basic algebraic structures appearing in bordered Floer homology
need to be generalized for the extension.  Rather than associating a
\dg algebra to a surface, we associate a certain kind of curved
$\Ainf$-algebra, which we call a \emph{weighted $\Ainf$-algebra}.
Similarly, the modules associated to three-manifolds with boundary
have the form of weighted $\Ainf$-modules, and weighted type $D$
structures. The definitions of these algebraic objects are fairly
straightforward. By contrast, the notions of bimodules and the
generalizations of the tensor product to them are more complicated to
define.  Specifically, they depend on a generalization of Stasheff's
associahedron which we call the \emph{weighted associaplex} and
various notions of ``diagonals'', generalizing the ones used by
Saneblidze-Umble to construct tensor products of
$\Ainf$-algebras~\cite{SU04:Diagonals} (see
also~\cite{MS06:AssociahedraProdAinf,Loday11:DiagonalStasheff}).

Although bordered Heegaard Floer homology was our motivation for this
work, we hope that some of the formal aspects explored here will be
of independent interest. For example, it is quite likely that other
gauge-theoretic invariants, such as instanton
homology~\cite{Floer88:instanton} and Seiberg-Witten
theory~\cite{KronheimerMrowka} have bordered analogues, with
similar algebraic structure; see
especially~\cite{DWang}. It is also conceivable that the
algebraic structures may be of  interest elsewhere in symplectic
geometry (compare, for instance,~\cite{Seidel02:FukayaDef,Auroux10:Bordered,Amorim16:tensor}).

\subsection{Statement of results in the unweighted case}\label{sec:results}

As a warm-up, we start with the now well understood case of tensor
products of $\Ainf$-algebras.  This operation was described
by Saneblidze and Umble; we review their
construction, with a view towards subsequent generalizations.

Fix a commutative, unital $\Field$-algebra $\Ring$ (possibly graded).
Throughout, by a module (respectively chain complex) we mean a free
or, more generally, projective $\Ring$-module (respectively chain complex of
projective $\Ring$-modules). Unless otherwise specified, in this
section undecorated tensor products are over~$\Ring$. (For example, in
bordered Floer theory one might take $\Ring=\FF_2$.)

Let $K_n$ denote the $(n-2)$-dimensional {\em associahedron},
introduced by Stasheff in~\cite{Stasheff63:associahedron1}. Cells of
$K_n$ are in a natural one-to-one correspondence with planar, rooted
trees with $n$ inputs and no $2$-valent vertices.
An {\em{associahedron diagonal}} is a
collection of chain maps
$$
\AsDiag=\{\AsDiag_n\co \cellC{*}(K_n)\to \cellC{*}(K_{n}\times K_{n})\}_{n=2}^{\infty}
$$ 
satisfying compatibility and non-degeneracy conditions
(Definition~\ref{def:AssociahedronDiagonal}).

Given $\Ainf$-algebras $\Alg_1$ and $\Alg_2$ over $\Ring$, we use the
associahedron diagonal to define a tensor product over $\Ring$, denoted
$\Alg_1\ADtp\Alg_2$ (Definition~\ref{def:Alg-tp}), giving the 
following reformulation of results from Saneblidze-Umble~\cite{SU04:Diagonals}:

\begin{theorem}
  \label{thm:Algebras}
  There exist associahedron diagonals $\AsDiag$,
  which can be used to define a tensor
  product $\Alg_1\ADtp\Alg_2$ of $\Ainf$-algebras
  $\Alg_1$ and $\Alg_2$, with the following properties:
  \begin{enumerate}[label=(A\arabic*),ref=(A\arabic*)]
  \item\label{item:Alg-thm-vec-space} The underlying $\Ring$-module of
    $\Alg_1\ADtp\Alg_2$ is the tensor product (over $\Ring$)
    $A_1\otimes A_2$ of the underlying $\Ring$-modules of $\Alg_1$ and
    $\Alg_2$.
  \item\label{item:Alg-thm-dg} If $\Alg_1$ and $\Alg_2$ are
    differential graded algebras, then $\Alg_1 \ADtp\Alg_2$ is a
    \dg~algebra and there is an isomorphism of \dg~algebras
    \[\Alg_1\ADtp\Alg_2\cong\Alg_1\otimes \Alg_2,\]
    where the right-hand-side is the usual tensor product of
    \dg~algebras.
  \item\label{item:Alg-thm-qi} If $\Alg_1$ is quasi-isomorphic to
    $\Alg_1'$ and $\Alg_2$ is quasi-isomorphic to $\Alg_2'$, then
    $\Alg_1\ADtp\Alg_2$ is quasi-isomorphic to
    $\Alg_1'\ADtp\Alg_2'$.
  \item\label{item:Alg-thm-assoc} The tensor product is associative up
    to isomorphism: if $\Alg_1$, $\Alg_2$, and $\Alg_3$ are
    $\Ainf$-algebras, there is an isomorphism of $\Ainf$-algebras
    $$\Alg_1\ADtp(\Alg_2\ADtp \Alg_3)
    \cong(\Alg_1\ADtp \Alg_2)\ADtp \Alg_3.$$
  \item\label{item:Alg-thm-change-diag} If $\AsDiag$ and $\AsDiag'$ are
    both associahedron diagonals, then there is an isomorphism of $\Ainf$-algebras
    \[
      \Alg_1\ADtp\Alg_2 \cong \Alg_1\ADtp[\AsDiag']\Alg_2.
    \]
  \end{enumerate}
\end{theorem}

Since any $\Ainf$-algebra is quasi-isomorphic to a \dg algebra, the
following corollary characterizes $\Alg_1\ADtp\Alg_2$ up to
quasi-isomorphism:
\begin{corollary}
  \label{cor:DGresolution}
  Let $\Alg_1$ and $\Alg_2$ be $\Ainf$-algebras and $\Blg_1$ and
  $\Blg_2$ \dg algebras so that $\Blg_i$ is quasi-isomorphic to
  $\Alg_i$. Then for any associahedron diagonal $\AsDiag$,
  $\Alg_1\ADtp\Alg_2$ is quasi-isomorphic to
  $\Blg_1\otimes\Blg_2$.
\end{corollary}

\begin{remark}
  Corollary~\ref{cor:DGresolution} follows from
  Properties~\ref{item:Alg-thm-dg} and~\ref{item:Alg-thm-qi} of
  Theorem~\ref{thm:Algebras}.  Since any $\Ainf$-algebra is
  quasi-isomorphic to a \dg algebra, we can view
  Properties~\ref{item:Alg-thm-assoc}
  and~\ref{item:Alg-thm-change-diag} as easy consequences of
  Corollary~\ref{cor:DGresolution}.  We have chosen instead the
  formulation from Theorem~\ref{thm:Algebras}, as it generalizes more
  readily to the weighted context; see
  Remark~\ref{rmk:WeightedResolution}.
\end{remark}

Going beyond tensor products of $\Ainf$-algebras, we consider next the
case of $\Ainf$-modules.  Fix an associahedron diagonal $\AsDiag$.  A
(right) \emph{module diagonal} $\MDiag$ compatible with $\AsDiag$ is a
collection of chain maps
$$\MDiag=\{\MDiag_n\co \cellC{*}(K_n)\to \cellC{*}(K_{n}\times K_{n})\}_{n=2}^{\infty}$$
that satisfy a compatibility condition with respect to $\AsDiag$ and
a non-degeneracy condition (Definition~\ref{def:ModuleDiagonal}).

Given $\Ainf$-algebras $\Alg_1$ and $\Alg_2$ and (right)
$\Ainf$-modules $\cModule_1$ and $\cModule_2$ over $\Alg_1$ and
$\Alg_2$, we use $\MDiag$ to define a tensor product
$\cModule_1\MDtp \cModule_2$ of $\cModule_1$ and $\cModule_2$ over
$\Ring$ (Definition~\ref{def:Mod-tp}).

\begin{theorem}
  \label{thm:ModuleDiagonalExists}
  For any associahedron diagonal $\AsDiag$, there is as a
  compatible module diagonal $\MDiag$ which, given $\Ainf$-modules
  $\cModule_1$ and $\cModule_2$ over $\Alg_1$ and $\Alg_2$ respectively, can
  be used to define their tensor product $\cModule_1\MDtp\cModule_2$.
  The tensor product $\cModule_1\MDtp\cModule_2$
  is an $\Ainf$-module over $\Alg_1\ADtp \Alg_2$, and
  this construction has the following properties:
  \begin{enumerate}[label=(M\arabic*),ref=(M\arabic*)]
  \item\label{item:Mod-thm-dg} If $\cModule_1$ and $\cModule_2$ are
    differential graded modules over \dg algebras $\Alg_1$ and $\Alg_2$, then
    \[\cModule_1\MDtp\cModule_2\cong \cModule_1\otimes \cModule_2,\] where the
    right-hand-side is the usual (external) tensor product of \dg modules.
  \item\label{item:Mod-thm-qi} If $\cModule_1$ is $\Ainf$-homotopy equivalent to $\cNodule_1$
    and $\cModule_2$ is $\Ainf$-homotopy equivalent to $\cNodule_2$, then
    $\cModule_1\MDtp\cModule_2$ is $\Ainf$-homotopy equivalent to
    $\cNodule_1\MDtp\cNodule_2$.
  \item\label{item:Mod-thm-change-diag} If $\MDiag$ and $\MDiag'$ are
    both module diagonals compatible with $\AsDiag$, then
    $\cModule_1\MDtp\cModule_2$ is isomorphic to
    $\cModule_1\MDtp[\MDiag']\cModule_2$ as modules over $\Alg_1\ADtp
    \Alg_2$.
  \end{enumerate}
\end{theorem}

The proofs of parts~\ref{item:Alg-thm-qi}
and~\ref{item:Alg-thm-change-diag} of Theorem~\ref{thm:Algebras}
involve another sort of diagonal,
a \emph{multiplihedron diagonal}
(Definition~\ref{def:multiplihedron-diag}), which more generally allows
one to tensor $\Ainf$-algebra homomorphisms together; see
Lemma~\ref{lem:Alg-map-tp}.  Similarly, the proofs of
parts~\ref{item:Mod-thm-qi} and~\ref{item:Mod-thm-change-diag}
of Theorem~\ref{thm:ModuleDiagonalExists} involve
module-map diagonals (Definition~\ref{def:mod-map-diag}), which more
generally allow one to tensor morphisms of modules together; see
Lemma~\ref{lem:tens-mod-maps}, Proposition~\ref{prop:dg-bifunctor},
and Corollary~\ref{cor:mod-id-tens-id}.

Our interest in diagonals stems from the fact that they can be used to
define a bimodule analogue of special types of modules, which we called~type $D$
structures~\cite{LOT1}. Type $D$ structures are essentially the same
as twisted complexes over $\Alg$
(see~\cite{BondalKapranov,Kontsevich}) or comodules over
the bar complex $\Barop(\Alg)$ (compare~\cite{LefevreAInfinity,KellerLefevre}); see
also~\cite[Remarks 2.2.36
and 2.2.37]{LOT2}.

In more detail, let $\Ground$ be a commutative $\Ring$-algebra,
possibly graded, but free (or projective) as an $\Ring$-module. (For example,
$\Ground$ might be a finite direct sum of copies of $\Ring$ or $\Ring[Y]$. See also Convention~\ref{conv:Ring}.) Fix an
$\Ainf$-algebra $\Alg=(A,\{\mu_n\})$ over $\Ground$, and let
\[
  {\Tensor}^*A=\bigoplus_{n=0}^\infty A^{\otimes_\Ground n}
  \qquad{\text{and}}\qquad 
  {\overline\Tensor}^*A=\prod_{n=0}^\infty A^{\otimes_\Ground n}.
\]
Let $P$ be a projective,
graded $\Ground$-module, equipped with a map
$\delta^1_P\co P \to A\otimes_\Ground P\grs{1}$,
where $\grs{1}$ denotes a grading shift (Convention~\ref{conv:grading-shift}).
The map $\delta^1_P$ can be iterated to construct a map
$\delta\co P \to {\overline\Tensor}^*A\otimes P,$
and 
the $\Ainf$ operations on $\Alg$ fit together to give a map
$\mu \co {\Tensor}^*A \to A.$
We say that the pair $(P,\delta^1_P\co P \to A\otimes P\grs{1})$ is a
{\em type $D$ structure} if the following compatibility condition
holds:
\begin{equation}\label{eq:typeD}
(\mu\otimes \Id_P)\circ \delta = 0.
\end{equation}
This formula makes sense if
either large enough iterates of $\delta^1_P$ vanish or
the $\mu_i$ on $\Alg$ vanish for $i$ sufficiently large.  There is a
tensor product pairing right $\Ainf$-modules $\cModule_{\Alg}$ with left
type $D$ structures, written
$\cModule_{\Alg}\DT\lsup{\Alg}P$. (See~\cite{LOT1} or
Section~\ref{sec:box} below.)

Given a diagonal, there is a bimodule analogue of type $D$ structures:

\begin{definition}\label{def:DD-intro}
  Fix $\Ainf$-algebras $\Alg$ and $\Blg$. A (left-left)
  \emph{type~\DD~structure} over $\Alg$ and $\Blg$ is a type $D$
  structure over $\Alg\ADtp\Blg$.
\end{definition}

(We adopt the convention from our earlier papers~\cite{LOT2} that an
algebra as a subscript, as in $\cModule_\Alg$, denotes an $\Ainf$-module,
while an algebra as a superscript, as in $\lsup{\Alg}P$, denotes a
type $D$ structure. An $\Ainf$-algebra $\Alg$ has an \emph{opposite
  algebra} $\Alg^\op$, and a right $\Ainf$-module $\cModule_\Alg$ over $\Alg$
is the same as a left $\Ainf$-module $\lsub{\Alg^\op}\cModule$ over
$\Alg^\op$; and similarly for type $D$ structures. In particular, a
left-right type \DD\ structure $\lsup{\Alg}P^\Blg$ over $\Alg$ and
$\Blg$ is the same as a left-left type \DD\ structure
$\lsup{\Alg,\Blg^\op}P$ over $\Alg$ and $\Blg^\op$.)

Fix an associahedron diagonal $\AsDiag$ and corresponding module
diagonal $\MDiag$.  If $\lsup{\Alg,\Blg}P$ is a type~\DD\ structure
with respect to $\AsDiag$ and $\cModule_\Alg$, $\cNodule_\Blg$ are
$\Ainf$-modules over $\Alg$ and $\Blg$, we can define the triple tensor product
$[\cModule_{\Alg}\DT \lsup{\Alg}P^{\Blg^\op}\DT \lsub{\Blg^\op}\cNodule]_{\MDiag}$ to be
the tensor product of the $(\Alg\ADtp\Blg)$-module
$\cModule_{\Alg}\MDtp \cNodule_{\Blg}$ with the
$\Alg\ADtp\Blg$ type $D$ structure
$\lsup{\Alg,\Blg}P$.
\begin{theorem}     
  \label{thm:TripleTensorProduct}
  Assuming that either $P$ is bounded or $\cModule$, $\cNodule$, $\Alg$
  and $\Blg$ are all bonsai
  (as defined in Section~\ref{sec:algebra}), 
  the triple tensor product is a well-defined chain complex. Moreover,
  different choices of $\MDiag$ give rise to homotopy equivalent triple
  tensor products.
\end{theorem}

Fix an associahedron diagonal. We will define the notion of a compatible 
\emph{module diagonal primitive}, $\TrMPrim$ (Definition~\ref{def:M-prim}).  The key
results are the following:

\begin{proposition}\label{prop:prim-exist}
  Given an associahedron diagonal $\AsDiag$, there is a module
  diagonal primitive $\TrMPrim$ compatible with $\AsDiag$.
\end{proposition}

\begin{lemma}\label{lem:prim-gives-diag}
  Given an associahedron diagonal $\AsDiag$, a module diagonal
  primitive $\TrMPrim$ compatible with $\AsDiag$ gives rise to an
  associated module diagonal $\MDiag=\MDiag_{\TrMPrim}$.
\end{lemma}

The utility of this construction is the following:

\begin{theorem}\label{thm:prim-DT}
  Under appropriate boundedness hypotheses, a module diagonal
  primitive gives a way to form a type $D$ structure
  $\lsup{\Alg}P^{\Blg^\op}\DT^\TrMPrim\lsub{\Blg^\op}\cNodule$ over $\Alg$.
  If $\cModule_{\Alg}$ is an $\Ainf$-module over $\Alg$, then there is
  an isomorphism
  \[
  \cModule_{\Alg}\DT (\lsup{\Alg}P^{\Blg^\op}\DT^\TrMPrim \lsub{\Blg^\op}\cNodule) \cong
  [\cModule_{\Alg}\DT \lsup{\Alg}P^{\Blg^\op} \DT
  \lsub{\Blg^\op}\cNodule]_{\MDiag_{\TrMPrim}},
  \] 
  where here the right-hand-side denotes the triple tensor product of
  Theorem~\ref{thm:TripleTensorProduct}.
\end{theorem}

The operation $\DT^{\TrMPrim}$ is functorial in an appropriate sense; see Proposition~\ref{prop:DT-funct} and Lemma~\ref{lem:Id-DT-Id}.

\subsection{Statement of results in the weighted case}\label{sec:wresults}
We start with some definitions. Let $\Ring$ and $\Ground$ be as above.
\begin{definition}
  A \emph{weighted $\Ainf$-algebra} or \emph{$w$-algebra} over $\Ground$ is a curved
  $\Ainf$-algebra $\wAlg=(A\llbracket t\rrbracket,\mu_n)$ over
  $\Ground\llbracket t\rrbracket$, free over
  $\Ring\llbracket t\rrbracket$, with $t$ central, and so that the curvature $\mu_0$ lies in
  $tA\llbracket t\rrbracket$. In particular,
  $A=A\llbracket t\rrbracket/tA\llbracket t\rrbracket$ inherits the
  structure of an uncurved $\Ainf$-algebra, called the
  \emph{undeformed $\Ainf$-algebra of $\wAlg$}.

  A \emph{homomorphism} $f\co \wAlg\to\wBlg$ of weighted
  $\Ainf$-algebras is a homomorphism of curved $\Ainf$-algebras so
  that $f_0\in t B\llbracket t\rrbracket$. A homomorphism $f$ is a
  \emph{quasi-isomorphism} if the induced map of undeformed
  $\Ainf$-algebras is a quasi-isomorphism.

  A \emph{weighted $\Ainf$-module} or \emph{$w$-module} $\wMod$ over
  $\wAlg$ is a curved $\Ainf$-module
  $\wMod=(M\llbracket t\rrbracket,m_n)$ over
  $\Ground\llbracket t\rrbracket$, free over
  $\Ring\llbracket t\rrbracket$, so that the curvature
  $m_0$ lies in $tM\llbracket t\rrbracket$.
\end{definition}
See Definitions~\ref{def:wAinfty},~\ref{def:wAlg-homo},
and~\ref{def:wmod} for equivalent, more explicit definitions. The
variable $t$ is allowed to have an arbitrary grading, and different
gradings of $t$ lead to different notions. (Again, see
Section~\ref{sec:wAinfty} for further discussion.)  For the complex of
morphisms between weighted $\Ainf$-modules, see
Section~\ref{sec:wmod-morph}.

An ordinary $\Ainf$-algebra $\Alg$ can be viewed as a weighted
$\Ainf$-algebra by extension of scalars. We call such algebras
\emph{weighted trivially}. (Equivalently, in the notation of
Definition~\ref{def:wAinfty}, an algebra is weighted trivially if the
maps $\mu_n^w$ vanish when $w>0$.)

The analogue of the cellular chain complex of the associahedron in the
weighted case is the \emph{weighted trees complex} $\wTreesCx{*}{*}$,
defined in Section~\ref{sec:wAlgs}, and a slight extension of
it, the \emph{extended weighted trees complex} $\xwTreesCx{*}{*}$
(Section~\ref{sec:walg-diag}). The weighted trees complex is the
cellular chain complex of the \emph{weighted associaplex}, a cell
complex introduced in Section~\ref{sec:Associaplex}. Unlike the
associahedron, the associaplex is not a polyhedron; but like the associahedron,
it is
contractible (Theorems~\ref{thm:AssociaplexAcyclic} and~\ref{thm:assoc-ball}).

A \emph{weighted algebra diagonal} is a map
\[
\wDiag{n}{w}\co \wTreesCx{n}{w}\to \bigoplus_{w_1,w_2\leq
      w}Y_1^{w-w_1}Y_2^{w-w_2}\xwTreesCx{n}{w_1}\otimes\xwTreesCx{n}{w_2}\otimes\Ring[Y_1,Y_2]
\]
satisfying certain compatibility and non-degeneracy conditions
(Definition~\ref{def:w-alg-diag}). Unlike the unweighted case, not all
weighted algebra diagonals are equivalent: rather, the homotopy class
is determined by the \emph{seed} $\wDiag{0}{1}$.

A weighted algebra $\wAlg$ is \emph{strictly unital} if there is an
element $\One\in A$ so that for each $a\in A$,
$\mu_2(\One,a)=\mu_2(a,\One)=a$ and for any $n\neq 2$ and
$a_1,\dots,a_{n-1}\in A$,
$\mu_n(a_1,\dots,a_{i-1},\One,a_i,\dots,a_{n-1})=0$. Strict unitality
for weighted modules over strictly unital algebras and homomorphisms
of strictly unital algebras is defined similarly; see
Section~\ref{sec:w-units}. Unlike in the unweighted case, we 
define tensor products only of strictly unital weighted algebras and
modules (or, more generally, homotopy unital ones, in Section~\ref{sec:w-htpy-unital}). 

We can now state the weighted analogue of Theorem~\ref{thm:Algebras}:
\begin{theorem}
  \label{thm:wAlgebras}
  There is a weighted algebra diagonal $\wDiag{*}{*}$ with any given
  seed. Further, given a weighted algebra diagonal $\wDiag{*}{*}$, there is an associated
  tensor product of strictly unital weighted algebras $\wAlg_1$ and $\wAlg_1$ over
  $\Ring$, $\wAlg_1\wADtp\wAlg_2$, satisfying the following
  properties:
  \begin{enumerate}[label=(wA\arabic*),ref=(wA\arabic*)]
  \item\label{item:wAlg-thm-Ainf} If $\wAlg_1$ and $\wAlg_2$ are
    weighted trivially then $\wAlg_1\wADtp\wAlg_2$ coincides with the tensor product 
    of undeformed $\Ainf$-algebras (with respect to a diagonal obtained by
    restricting $\wDiag{*}{*}$).
  \item\label{item:wAlg-thm-qi} If $\wAlg_1$ is quasi-isomorphic to
    $\wBlg_1$ and $\wAlg_2$ is quasi-isomorphic to $\wBlg_2$, then
    $\wAlg_1\wADtp\wAlg_2$ is quasi-isomorphic to
    $\wBlg_1\wADtp\wBlg_2$.
  \item\label{item:wAlg-thm-change-diag} If $\wDiag{*}{*}_1$ and $\wDiag{*}{*}_2$ are
    weighted algebra diagonals with the same seed, then $\wAlg_1\wADtp[\wADiag_1]\wAlg_2$
    is isomorphic to $\wAlg_1\wADtp[\wADiag_2]\wAlg_2$.
  \end{enumerate}
\end{theorem}
The tensor product of strictly unital weighted algebras is not itself
strictly unital, but only weakly unital
(Definition~\ref{def:w-weakly-unital}). Any weakly unital weighted algebra is
isomorphic to a strictly unital one
(Theorem~\ref{thm:UnitalIsUnitalW}), but this extra isomorphism makes
associativity of
the tensor product in the weighted case more subtle than the
unweighted case; see Section~\ref{sec:w-htpy-unital}.

\begin{remark}
  \label{rmk:WeightedResolution}
  While every $\Ainf$-algebra over a field is quasi-isomorphic to a \dg algebra,
  we are not aware of a completely analogous result for weighted
  $\Ainf$-algebras (but see~\cite[Corollary 3.4]{Amorim16:tensor}). So, while one can define
  the tensor product of two $\Ainf$-algebras by resolving and then taking the
  ordinary tensor product, we do not know a weighted analogue of this operation.

  For a particular seed, Amorim has constructed a tensor product as in
  Theorem~\ref{thm:wAlgebras} by other methods~\cite{Amorim16:tensor}; see also
  Remark~\ref{rem:Amorim}.
\end{remark}

Turning to modules, the weighted analogue of Theorem~\ref{thm:ModuleDiagonalExists} is:
\begin{theorem}
  \label{thm:wMods}
  For any weighted algebra diagonal $\wDiag{*}{*}$, there is as a
  corresponding weighted module diagonal $\wModDiag{*}{*}$ which, given
  strictly unital weighted $\Ainf$-modules $\wMod_1$ and $\wMod_2$ over $\wAlg_1$ and
  $\wAlg_2$ respectively, can be used to define their tensor product
  $\wMod_1\wMDtp\wMod_2$, a weighted $\Ainf$-module over
  $\wAlg_1\wADtp\wAlg_2$. This tensor product has the
  following properties:
  \begin{enumerate}[label=(M\arabic*),ref=(M\arabic*)]
  \item\label{item:w-Mod-thm-unweighted} If $\wMod_1$ and $\wMod_2$ are
    trivially weighted $\Ainf$-modules over trivially weighted
    $\Ainf$-algebras then their tensor product is the same as in
    Theorem~\ref{thm:ModuleDiagonalExists}.
  \item\label{item:w-Mod-thm-qi} If $\wMod_1$ is homotopy equivalent to $\wNod_1$
    and $\wMod_2$ is homotopy equivalent to $\wNod_2$, then
    $\wMod_1\wMDtp\wMod_2$ is homotopy equivalent to
    $\wNod_1\wMDtp\wNod_2$.
  \item\label{item:w-Mod-thm-change-diag} If $\wModDiag{*}{*}_1$ and $\wModDiag{*}{*}_2$ are
    weighted module diagonals compatible with $\wDiag{*}{*}$, then
    $\wMod_1\wMDtp[\wModDiagNS_1]\wMod_2$ is homotopy equivalent to
    $\wMod_1\wMDtp[\wModDiagNS_2]\wMod_2$ as weighted $\Ainf$-modules over $\wAlg_1\wADtp
    \wAlg_2$.
  \end{enumerate}
\end{theorem}

Again, the tensor product does not preserve strict unitality, though
any weighted $\Ainf$-module is isomorphic to a strictly unital one
(Theorem~\ref{thm:UnitalIsUnitalMw}).

Similarly to the unweighted case, functoriality in
Theorem~\ref{thm:wAlgebras} uses weighted map diagonals
(Definition~\ref{def:w-map-diag}), the existence of which is guaranteed by
Lemma~\ref{lem:wmul-diag-exists}. See
Definition~\ref{def:w-alg-map-tensor} for the definition of the tensor
product of two weighted algebra homomorphisms.  Similarly,
functoriality in Theorem~\ref{thm:wMods} uses weighted module-map
diagonals (Definition~\ref{def:w-mod-map-diag}), existence and
uniqueness (up to homotopy) of which is
Lemma~\ref{lem:w-mod-map-diag-exists-unique}. Functoriality properties
of the tensor product of weighted modules are spelled out in
Proposition~\ref{prop:w-dg-bifunctor} and
Lemma~\ref{lem:w-mod-id-tens-id}.

We turn next to weighted type $D$ and \DD\ structures and their box
tensor products. With $\Ring$ and~$\Ground$ as above, let
$\wAlg=(A\llbracket t\rrbracket, \mu)$ be a
weighted $\Ainf$-algebra over $\Ground$ and $X$ an element of
$\Ground$ which acts centrally on $A$, in the sense that $aX=Xa$ for
all $a\in A$. Write
$\mu_n=\sum \mu_n^w t^w$. A \emph{weighted type $D$ structure with
  charge $X$} over $\wAlg$ is a projective $\Ground$-module $P$ and a map
$\delta^1\co P\to A\otimes P$ satisfying
\begin{equation}\label{eq:intro-D-str-rel}
  \sum_{w=0}^\infty\sum_{n=0}^\infty X^w(\mu_n^w\otimes\Id)\circ \delta^n=0\in A\otimes P.
\end{equation}
Here, $\delta^n\co P\to A^{\otimes n}\otimes P$ is the result of
iterating $\delta^1$ $n$ times. This structure equation makes sense
under suitable boundedness hypotheses on $\wAlg$ and/or $P$; see
Section~\ref{sec:wD}.  There is a box tensor product of a weighted
type $D$ structure and a weighted $\Ainf$-module, $\wMod\DT P$,
defined similarly to the unweighted case; see
Equation~\eqref{eq:w-DT-def} (Section~\ref{sec:wD}). (One could also
consider charges $X\in\Ground\llbracket t\rrbracket$.
Then, for charge
$X=t$, Equation~\eqref{eq:intro-D-str-rel} is the usual type $D$
structure relation over $\wAlg$ considered as a curved $\Ainf$-algebra.)

Parallel to the unweighted case, given weighted algebras $\wAlg$ and $\wBlg$
and a weighted algebra diagonal $\wDiag{*}{*}$, a weighted type \DD\
structure over $\wAlg$ and $\wBlg$ with respect to $\wDiag{*}{*}$ is a
type $D$ structure over $\wAlg\wADtp\wBlg$.
Given weighted modules $\wMod$ and $\wNod$ over $\wAlg$ and $\wBlg$,
under suitable boundedness hypotheses there is a corresponding triple
tensor product (Definition~\ref{def:wtriple-prod})
\[
  [\wMod_{\wAlg}\DT \lsup{\wAlg}P^{\wBlgop}\DT \lsub{\wBlgop}\wNod]_{\wModDiagNS}
  \coloneqq
  (\wMod\otimes_\wModDiagNS \wNod)_{\wAlg\otimes_{\wDiagNS}\wBlg} \,\DT\, \lsup{\wAlg\otimes_{\wDiagNS} \wBlg}P.
\]

\begin{theorem}
  \label{thm:wTripleTensorProduct}
  Assume that either $\lsup{\wAlg}P^{\wBlgop}$ is operationally bounded
  or $\wAlg$, $\wBlg$, $\wMod$, and $\wNod$ are bonsai. Then the
  triple tensor product is a well-defined chain complex. Moreover,
  different choices of weighted module diagonal~$\MDiag$ give rise to
  homotopy equivalent
  triple tensor products.
\end{theorem}

Again, for one-sided box tensor products of type \DD\ structures we use
a \emph{weighted module diagonal primitive} $\wTrPMDiag{*}{*}$
(Definition~\ref{def:wM-prim}). Some key results are:
\begin{itemize}
\item Given any weighted algebra diagonal, there is a compatible
  weighted module diagonal primitive
  (Proposition~\ref{prop:wprim-exist}).
\item A weighted module diagonal primitive gives rise to a weighted
  module diagonal (Lemma~\ref{lem:wprim-gives-diag}).
\item Weighted module diagonal primitives give a way of defining a
  one-sided box tensor product $\wMod\DT^{\wTrPMDiag{*}{*}}P$
  (Definition~\ref{def:w-one-sided-DT}).
\item The triple box tensor product agrees, up to homotopy, with
  taking the one-sided box tensor product twice
  (Proposition~\ref{prop:one-sided-DT-works}).
\end{itemize}
Functoriality of the one-sided box tensor product uses module-map
diagonal primitives (Definition~\ref{def:w-mod-map-prim}), and is
spelled out in Section~\ref{sec:w-one-sided-DT-func}, particularly
Proposition~\ref{prop:w-one-side-DT-bifunc} and
Lemma~\ref{lem:w-Id-DT-Id}.

\subsection{Additional topics}
To prove associativity of the tensor product of weighted algebras and
to prove that the category of type \DD\ structures is independent of
the choice of diagonal, we need to relax the unitality assumption, and
work with homotopy unital weighted algebras. We generally give these
generalizations at the end of each section.

As mentioned above, many of the algebraic constructions in this paper
require some boundedness properties to make sense. In most of the
paper, we work with boundedness properties which are easy to state and
make sense quite generally, but are not, in fact, sufficiently general
for the application to bordered Heegaard Floer homology. We relax
these conditions, over particular ground rings, in
Section~\ref{sec:boundedness}, to conditions that are satisfied in
bordered Floer theory.

For the reader's convenience, Appendix~\ref{sec:pictures} collects the
first few terms in various kinds of diagonals, and
Appendix~\ref{sec:tables} has a table of the various objects (Table 1)
introduced in this paper and a table of conventions (Table 2).

\subsection*{Acknowledgments.} We thank the referee for helpful comments.

%%% Local Variables: 
%%% mode: latex
%%% TeX-master: "AbstractDiagonal"
%%% TeX-command-extra-options: "-shell-escape"
%%% End: 

\section{Diagonals: definitions, existence}\label{sec:diags}
\begin{convention}\label{conv:diag-tens-F2}
  Let $\Ring$ be a commutative ring (possibly graded).
  In this section, undecorated tensor products are over $\Ring$.
\end{convention}

\subsection{Trees and the associahedron}\label{sec:associahedron}
\begin{convention}\label{conv:trees}
  In this paper, all trees will be planar trees, with a distinguished
  root leaf, which we call the \emph{output} (so our planar trees are
  \emph{rooted}); we call the other leaves the \emph{inputs}. (This will
  change in Section~\ref{sec:wAinfty}, where some leaves are neither
  inputs nor the output.) Non-leaf vertices are \emph{internal vertices}.

  We admit as a degenerate tree $\IdTree$ a single edge with one
  input, one output, and no internal vertices, which we sometimes call the
  \emph{identity tree}.

  We draw our trees so that the inputs come in the top and the root
  exits through the bottom, and each edge points downward.
  Consequently, each internal vertex has a leftmost and a rightmost parent.
\end{convention}

Each vertex $v$ of a tree $T$ has a \emph{valence} $d(v)$, the number of edges
incident to the vertex. A tree is called {\em binary} if each internal
vertex
has two parents (a left parent and a right parent), or equivalently if
each internal vertex has valence $3$.

We denote the free $\Ring$-module spanned by trees with $n$ inputs and
no $2$-valent vertices by $\Trees_{n}$. Let $\Trees$ be the graded
$\Ring$-module $\bigoplus_{n=1}^\infty\Trees_n$.

Note that the special case $\Trees_1$ is $1$-dimensional, spanned by
the tree $\downarrow$ with one input, one output, and no internal
vertices. (The module $\Trees_2$ is also $1$-dimensional.)

\begin{figure}
  \centering
  \input{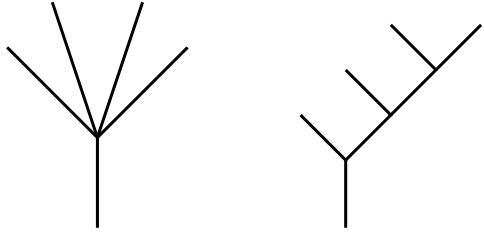}
  \caption[Some trees]{\textbf{Some trees}. 
    \label{fig:SomeTrees}
    At the left is the corolla $\corolla{4}$ with four inputs, which
    corresponds to the parenthesization $(1,2,3,4)$. To its right is
    the right-associated binary tree with four inputs, which corresponds to
    the parenthesization $(1,(2,(3,4)))$.}
\end{figure}

\begin{definition}
  A rooted, planar tree with a single internal vertex is called a
  \emph{corolla}. Let $\corolla{n}$ denote the corolla with $n$
  inputs (i.e., $n+1$ leaves); see
  Figure~\ref{fig:SomeTrees}.
\end{definition}

Given trees $S$ and $T$, let $T\circ_i S$ be the result of gluing the
output leaf of $S$ to the $i\th$ input leaf of $T$. (Note that
$\circ_i$ is not associative, because of how the indexing behaves.)
Extend $\circ_i$ bilinearly to a map
$\Trees_n\otimes \Trees_m\to \Trees_{m+n-1}$. Given
$S\in\Trees_m$, $T\in \Trees_n$ let $T\circ S=\sum_{i=1}^n T\circ_iS$.
It will also sometimes be convenient to use a multi-linear composition
map $\Trees_k\otimes \Trees_{m_1}\otimes\cdots\otimes\Trees_{m_k}\to
\Trees_{m_1+\cdots+m_k}$,
\begin{equation}\label{eq:multi-comp}
  (T,S_1,\dots,S_k)\mapsto T\circ(S_1,\dots,S_k),
\end{equation}
obtained by gluing the output of $S_i$ to the $i\th$ input of
$T$. That is,
\[
  T\circ (S_1,\dots,S_k)=(((T\circ_k
  S_k)\circ_{k-1}S_{k-1})\circ_{k-2}\cdots) \circ_1 S_1.
\]

We will often be interested in pairs of trees with $n$ inputs (or
elements of $\Trees_n\otimes \Trees_n$). Define
$(T_1,T_2)\circ_i (S_1,S_2)=(T_1\circ_i S_1,T_2\circ_{i} S_2)$ and
\begin{equation}\label{eq:circ}
  (T_1,T_2)\circ(S_1,S_2)=\sum_i(T_1,T_2)\circ_i(S_1,S_2)
  =\sum_i(T_1\circ_i S_1,T_2\circ_{i} S_2).
\end{equation}

We also use a multi-composition for pairs of trees, analogous to
Formula~\eqref{eq:multi-comp}. Specifically, given a pair of trees
$(T,T')$ with $k$ inputs and pairs of trees
$(S_1,S'_1),\dots,(S_k,S'_k)$, let
\begin{equation}\label{eq:multi-comp-pair}
  (T,T')\circ \bigl((S_1,S'_1),\dots,(S_k,S'_k)\bigr)=\bigl(T\circ(S_1,\dots,S_k),T'\circ(S'_1,\dots,S'_k)\bigr).
\end{equation}
The similarity in notation between Formulas~\eqref{eq:circ}
and~\eqref{eq:multi-comp-pair} should not cause confusion: the only
overlap is when $(T,T')$ are trees with a single input, in which case
the two notations agree.

The {\em associahedron} $K_n$, introduced
by Stasheff~\cite{Stasheff63:associahedron1}, is a CW complex which is a
natural compactification of the configuration space of $n\geq 2$ distinct
points in
$\RR$ modulo translation and scaling~\cite{FukayaOh97:open-string,Devadoss98:tessellations}.  The
cells in the associahedron are in a natural one-to-one correspondence
with trees with $n$ inputs and no $2$-valent vertices. (We will
sometimes call such trees \emph{associahedron trees}.) 
In turn, cells correspond to parenthesizations of
the inputs; see Figure~\ref{fig:SomeTrees}.
In
particular, for each $n\geq 2$, the cellular chain complex
$\cellC{*}(K_n;\Ring)$ is identified with $\Trees_n$.

\begin{convention}\label{conv:cellC}
  By $\cellC{*}(\thinspace\cdot\thinspace)=\cellC{*}(\thinspace\thinspace\cdot\thinspace\thinspace;\Ring)$ we mean the cellular chain complex with
  $\Ring$-coefficients.
\end{convention}

\begin{definition}\label{def:dim}
  The \emph{dimension} of a tree $T \in \Trees_n$ with $k$ internal
  vertices is the dimension of the corresponding cell in $K_n$; namely,
  $\dim(T) = n - k - 1$.
\end{definition}

\begin{remark}
  We use the word {\em dimension} to refer to the grading on $\Trees_n$,
  $\cellC{*}(K_n)$, and so on, to distinguish it from another grading which appears
  in Section~\ref{sec:wDiagApps}.
\end{remark}

The differential in the cellular chain complex $\cellC{*}(X_n)$ of (the cell
corresponding to) such a tree $T$ is the sum over all the ways of
inserting an edge into $T$, so as to preserve the property that no
vertex has valence two.\signissue
For example, if we write a tree as a bracketing of its inputs then
\[
  \partial(1,2,3,4)=((1,2,3),4)+(1,(2,3,4))+((1,2),3,4)+(1,(2,3),4)+(1,2,(3,4)).
\]

The following fundamental result guides many of the results in this paper:
\begin{proposition}\cite{Stasheff63:associahedron1}
  \label{prop:AssociahedronIsContractible}
  The complex $\cellC{*}(K_n)$ is contractible; i.e.,
  \[
  H_i(\cellC{*}(K_n))=\left\{\begin{array}{ll}
  \Ring &\text{if $i=0$} \\
  0 &{\text{otherwise.}}\end{array}\right.\]
\end{proposition}
(The proposition also follows from the fact that $K_n$ 
can be realized as a convex
polytope; see~\cite{CeballosSantosZiegler} for a discussion of various
realizations of the associahedron and further citations.)

\begin{figure}
  \centering
  \input{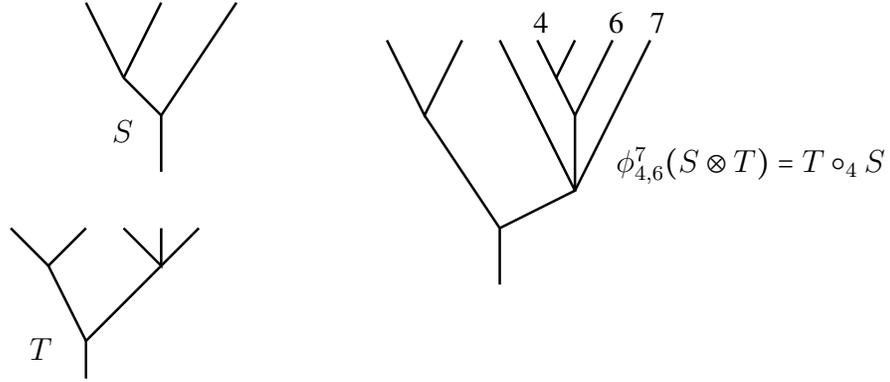}
  \caption[Examples of composition of trees]{\textbf{Examples of compositions of trees}. 
    An illustration of the stacking of
    trees $S$ and $T$, $\phi^7_{4,6}(S\otimes T)$.}
  \label{fig:SomeStackingTrees}
\end{figure}

Given any $1\leq i < j \leq n$, there is a codimension-1 face
$F^n_{i,j}$ such that
\begin{equation}\label{eq:assoc-face-include}
K_{j-i+1}\times K_{n+i-j} \cong F^n_{i,j}\subset K_{n},
\end{equation}
and these are all the codimension-1 faces.
These face maps induce chain maps
\begin{equation}\label{eq:face-maps}
\phi^n_{i,j}\co \cellC{*}(K_{j-i+1})\otimes \cellC{*}(K_{n+i-j})
\to \cellC{*}(K_n),
\end{equation}
which can be interpreted as compositions of trees: 
\[
\phi^n_{i,j}(S,T)=T\circ_iS.
\]
See Figure~\ref{fig:SomeStackingTrees}.

\subsection{Algebra diagonals}
For any space~$X$, 
the diagonal map $X\to X\times X$ induces a chain map 
on the singular chain complexes
% there is a natural diagonal map of the singular chain complex
\[
\AsDiag^X_{\sing} \co \singC{*}(X) \to \singC{*}(X\times X).
\]
Our first goal is to define an {\em associahedron diagonal}, which is a 
cellular approximation
$$\AsDiag^n\co \cellC{*}(K_n) \to \cellC{*}(K_n\times K_n)$$
to the diagonal and which is natural under the face maps, in a sense to be
made precise presently. 

\begin{definition}
  \label{def:AssociahedronDiagonal}
  An {\em associahedron diagonal} is a sequence of (degree-preserving) chain maps
  $$
  \bigl\{\,\AsDiag^n\co \cellC{*}(K_n)\to \cellC{*}(K_n)\otimes \cellC{*}(K_n)\,\bigr\}_{n=2}^\infty
  $$
  satisfying the following conditions:
  \begin{itemize}
  \item \textbf{Compatibility under stacking}:\signissue
    \begin{equation}
      \label{eq:CompatibilityCondition}
      \AsDiag^n\circ \phi^n_{i,j} = (\phi^n_{i,j}\otimes \phi^n_{i,j})\circ (\AsDiag^{j-i+1}\otimes \AsDiag^{n+i-j}),
    \end{equation}
    with the understanding that the composition on the right-hand side
    of Formula~\eqref{eq:CompatibilityCondition} involves a shuffling
    of tensor factors, i.e., it implicitly uses the identification
    \begin{multline}\label{eq:as-diag-shuf}
      \cellC{*}(K_{j-i+1})\otimes \cellC{*}(K_{n+i-j})\otimes \cellC{*}(K_{j-i+1})\otimes \cellC{*}(K_{n+i-j}) \\
      \cong \cellC{*}(K_{j-i+1})\otimes \cellC{*}(K_{j-i+1})\otimes \cellC{*}(K_{n+i-j}) \otimes \cellC{*}(K_{n+i-j}).
    \end{multline}
  \item \textbf{Non-degeneracy}: $\AsDiag^2\co \cellC{0}(K_2)\to\cellC{0}(K_2)\otimes \cellC{0}(K_2)$ is the standard isomorphism
    ${\Ring}\cong {\Ring}\rotimes{\Ring}{\Ring}$. (Note that $K_2$ is a single point.)
  \end{itemize}
\end{definition}

The term {\em associahedron diagonal} is used to distinguish these
maps from other diagonals (such as
module diagonals and multiplihedron diagonals) defined below.
An associahedron diagonal was constructed by Saneblidze and
Umble~\cite{SU04:Diagonals} to define the tensor products of
$\Ainf$-algebras; see
also~\cite{MS06:AssociahedraProdAinf,Loday11:DiagonalStasheff}.

An associahedron diagonal is encoded conveniently in the following:

\begin{definition}
  \label{def:DiagonalCells}
  An {\em associahedron tree diagonal} consists of $(n-2)$-chains
  $\TrDiag_n\in \cellC{*}(K_n)\otimes \cellC{*}(K_n)$ for each $n\geq 2$ satisfying:
  \begin{itemize}
  \item \textbf{Compatibility}:\signissue
    \begin{equation}
      \label{eq:DiagonalCellCompatibility}
      \partial \TrDiag_n = \sum_{1\leq i<j\leq n}(\phi^n_{i,j}\otimes \phi^n_{i,j})(\TrDiag_{j-i+1}\otimes \TrDiag_{n+i-j}).
    \end{equation}
    (Here, we have the same shuffling of factors as in
    Formula~\eqref{eq:CompatibilityCondition}.)
  \item \textbf{Non-degeneracy}: $\TrDiag_2\in C_0(K_2)\otimes
    C_0(K_2)$ is the generator. 
  \end{itemize}

  Equivalently, we can view an associahedron tree diagonal as elements
  $\TrDiag_n\in \Trees_n\otimes\Trees_n$ of dimension $n-2$,
  satisfying:
  \begin{itemize}
  \item \textbf{Compatibility}:\signissue
    \begin{equation}
      \label{eq:DiagonalCellCompatibilityTree}
      \partial \TrDiag_n = \sum_{i+j=n+1}\TrDiag_i\circ \TrDiag_j,
    \end{equation}
    where $\circ$ is the sum of all ways of composing, as in Equation~\eqref{eq:circ}.
  \item \textbf{Non-degeneracy}: $\TrDiag_2=\corolla{2}\otimes\corolla{2}$.
  \end{itemize}
\end{definition}

Schematically, the compatibility condition is:
\[
  \tikzsetnextfilename{def-DiagonalCells-1}
\bdy\Biggl(\mathcenter{
  \begin{tikzpicture}[smallpic]
    \node at (0,0) (tl) {};
    \node at (1,0) (tr) {};
    \node at (0,-1) (delta1) {$\TrDiag$};
    \node at (1,-1) (delta2) {$\TrDiag$};
    \node at (0,-2) (bl) {};
    \node at (1,-2) (br) {};
    \draw[taa] (tl) to (delta1);
    \draw[taa] (tr) to (delta2);
    \draw[alga] (delta1) to (bl);
    \draw[alga] (delta2) to (br);
  \end{tikzpicture}
}\Biggr)
=
\mathcenter{
  \tikzsetnextfilename{def-DiagonalCells-2}
  \begin{tikzpicture}[smallpic]
    \node at (-2.5,0) (tll) {};
    \node at (-1.5,0) (tl) {};
    \node at (-.5,0) (tlr) {};
    \node at (.5,0) (trl) {};
    \node at (1.5,0) (tr) {};
    \node at (2.5,0) (trr) {};
    \node at (-1.5,-1) (deltaL) {$\TrDiag$};
    \node at (1.5,-1) (deltaR) {$\TrDiag$};
    \node at (-1.5,-2) (deltaL2) {$\TrDiag$};
    \node at (1.5,-2) (deltaR2) {$\TrDiag$};
    \node at (-1.5,-3) (bl) {};
    \node at (1.5,-3) (br) {};
    \draw[taa] (tl) to (deltaL);
    \draw[taa] (tr) to (deltaR);
    \draw[taa] (tll) to (deltaL2);
    \draw[taa] (tlr) to (deltaL2);
    \draw[taa] (trl) to (deltaR2);
    \draw[taa] (trr) to (deltaR2);
    \draw[alga] (deltaL) to (deltaL2);
    \draw[alga] (deltaR) to (deltaR2);
    \draw[alga] (deltaL2) to (bl);
    \draw[alga] (deltaR2) to (br);
  \end{tikzpicture}
}.
\]

The correspondence between associahedron diagonals and associahedron tree diagonals
comes from the fact that associahedron diagonals are
determined by their values on the corollas, which correspond to the
top-dimensional cells of the associahedra.

\begin{lemma}
  \label{lem:DiagonalCells}
  There is a one-to-one correspondence between associahedron diagonals
  and associahedron tree diagonals, which associates to an
  associahedron diagonal $\AsDiag$ the associahedron tree diagonal
  $\AsDiag^n(\corolla{n})$.
\end{lemma}
\begin{proof}
  First, we show that if $\AsDiag^n$ is an associahedron diagonal, then the cells
  $\TrDiag_n=\AsDiag^n(\corolla{n})$ form an associahedron tree diagonal, in the sense of
  Definition~\ref{def:DiagonalCells}.
  Start from the identity from Equation~\eqref{eq:assoc-face-include}
  $$\partial \corolla{n}=\sum_{1\leq i<j\leq n} \phi^n_{i,j}(\corolla{j-i+1}\otimes \corolla{n+i-j}).$$
  Apply $\AsDiag^n$ to both sides and use
  Equation~\eqref{eq:CompatibilityCondition} to conclude that
  Equation~\eqref{eq:DiagonalCellCompatibility} holds.

  Conversely, given an associahedron tree diagonal
  $\TrDiag\in \Trees_n$ define
  $\AsDiag^n\co \cellC{*}(K_n)\to \cellC{*}(K_n)\otimes
  \cellC{*}(K_n)$
  as follows. Identifying $\cellC{*}(K_n)$ and $\Trees_n$, it suffices
  to define the value of $\AsDiag^n$ on any tree $T\in\Trees_n$. We
  can write $T$ as a composition of corollas,
  \[
    T=\corolla{n_k}\circ_{i_{k-1}}\corolla{n_{k-1}}\circ_{i_{k-2}}\cdots\circ_{i_1}\corolla{n_1}.
  \]
  Define 
  \[
    \AsDiag(T)=\TrDiag_{n_k}\circ_{i_{k-1}}\TrDiag_{n_{k-1}}\circ_{i_{k-2}}\cdots\circ_{i_1}\TrDiag_{n_1}.
  \]
  (See Figure~\ref{fig:diag-cells}.)
  It is clear that this definition is independent of the decomposition
  of $T$ and, consequently, that $\AsDiag$ satisfies
  Equation~\eqref{eq:CompatibilityCondition}. The non-degeneracy
  conditions on $\TrDiag$ and on~$\AsDiag$ also obviously correspond.

  \begin{figure}
    \centering
    %Font is 12 point.
    \includegraphics{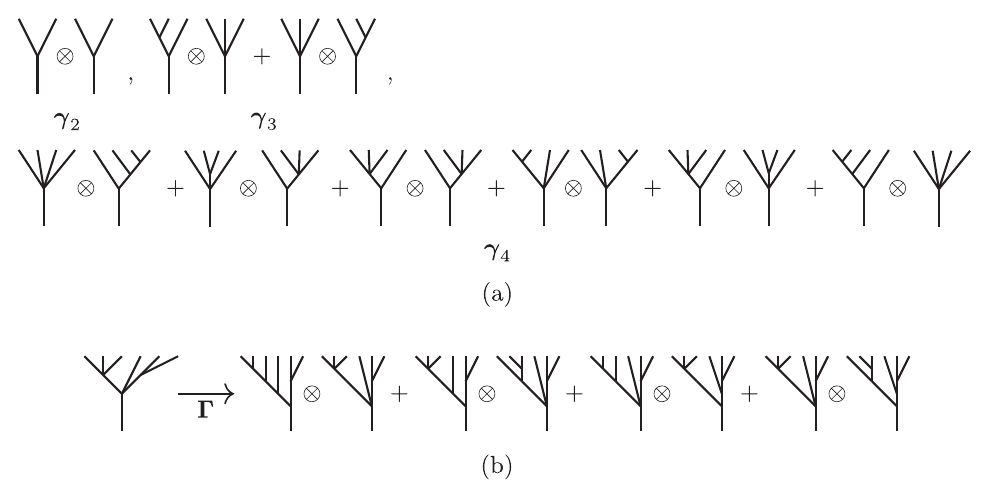}
    \caption[From an associahedron tree diagonal to an associahedron diagonal]{\textbf{From an associahedron tree diagonal to an associahedron diagonal.} (a) The first few terms in a particular associahedron tree diagonal. (b) Substituting these terms into a more complicated tree.}
    \label{fig:diag-cells}
  \end{figure}
  
  It remains to verify that the extension is a chain map. Given a tree
  $T\in \Trees_n=\cellC{*}(K_n)$, the cellular boundary
  $\bdy(T)$ is the sum of all ways of inserting an edge at some
  internal vertex~$v$ in~$T$.
  The terms in
  $\AsDiag(\bdy(T))$ coming from $v$ correspond to breaking $v$ into two
  vertices, with valences $i,j>2$, and inserting $\TrDiag_i$
  and $\TrDiag_j$ at these two vertices. This corresponds to the
  right-hand side of Formula~(\ref{eq:DiagonalCellCompatibility}).
  The terms in $\bdy (\AsDiag(T))$ coming from $v$ correspond to
  inserting $\bdy\TrDiag_n$ (rather than $\TrDiag_n$) at $v$. This
  corresponds to the left-hand side of
  Formula~(\ref{eq:DiagonalCellCompatibility}).
  The
  result follows.\signissue
\end{proof}

\begin{lemma}
  \label{lem:ExistAssociahedronDiagonal}
  There is an associahedron diagonal.
\end{lemma}

\begin{proof}
  We describe the inductive procedure for constructing an
  associahedron tree diagonal. The base case $n=2$ is specified by the
  non-degeneracy condition, and a solution for $n=3$ is shown in
  Figure~\ref{fig:diag-cells}. For the inductive step, suppose we
  already have $\TrDiag_n$ for $n=2,\dots,m$ satisfying
  Equation~\eqref{eq:DiagonalCellCompatibility}. To find
  $\TrDiag_{m+1}$, $m\geq 3$, we observe that the right-hand-side of
  Equation~\eqref{eq:DiagonalCellCompatibility} determines a cycle:
  \signissue
  \[ 
  \partial \sum_{1\leq i<j\leq {m+1}}(\phi^{m+1}_{i,j}\otimes \phi^{m+1}_{i,j})(\TrDiag_{j-i+1}\otimes \TrDiag_{{m+1}+i-j})=0.
  \]
  This is a straightforward consequence of the facts that
  $\phi^{m+1}_{i,j}$ are chain maps and the $\TrDiag_{n}$ satisfy
  Equation~\eqref{eq:DiagonalCellCompatibility} for $n \le m$: the
  terms cancel in pairs.  According to Stasheff, $X$
  is contractible~\cite[Proposition 3]{Stasheff63:associahedron1}; in
  particular, $H_*(X)=0$ for
  $*>0$. It follows that $H_*(X\times X)=0$ for $*>0$ as well. Since the right hand side of Equation~\eqref{eq:DiagonalCellCompatibility} lies in dimension $n-3$ we can find some $\TrDiag_{m+1}$ satisfying
  Equation~\eqref{eq:DiagonalCellCompatibility} for $n=m+1\geq 4$, as needed.
\end{proof}

\begin{figure}
  \centering
   \input{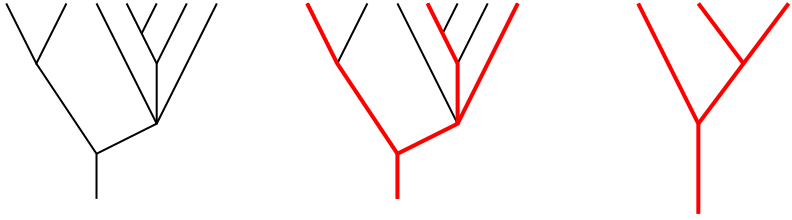} 
   \caption[A profile tree]{\textbf{A profile tree.} The tree on the left has $n=7$. At the
     right we have drawn the profile tree corresponding to the triple of inputs
     $1$, $4$, and $7$.}
  \label{fig:profile-tree}
\end{figure}

\begin{example}\label{eg:explicit-diag}
  Given a tree $T$ with $n$ inputs and a subset $V\subset \{1,\dots,n\}$ of the inputs to $T$, the
  \emph{profile tree} of $T$ with profile $V$, denoted $T(V)$, is the result of deleting the inputs to $T$ which are not in $V$ and collapsing any resulting $2$-valent vertices. See Figure~\ref{fig:profile-tree}.
  
  Call a pair of three-input trees \emph{right-moving} if it is one of the following five pairs:
  \begin{center}
    \includegraphics{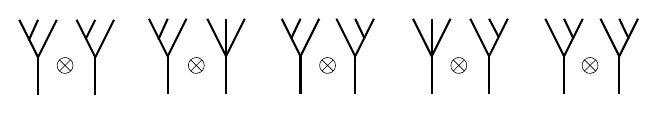}
  \end{center}
  Call a pair $(S,T)$ of $n$-input trees \emph{right-moving} if for
  any triple $v_1,v_2,v_3\in \{1,\dots,n\}$, the pair of three-input trees
  $(S(v_1,v_2,v_3),T(v_1,v_2,v_3))$ is right-moving.

  Masuda-Thomas-Tonks-Vallette~\cite{MTTV:explicit-diag} have shown
  that the sum of all trees which are right-moving and in dimension
  $n-2$ (where $n$ is the number of inputs to each tree) forms an
  associahedron tree diagonal. This is the associahedron tree diagonal
  for which the first three terms are shown in
  Figure~\ref{fig:diag-cells}. (Their technique uses a geometric
  procedure which is applicable to more general convex polytopes.)
\end{example}

\subsection{Module diagonals}

\begin{definition}
  \label{def:ModuleDiagonal}
  Fix an associahedron diagonal $\AsDiag$.  A \emph{right module diagonal
  compatible with $\AsDiag$} is a collection of chain
  maps $\MDiag^n \co \cellC{*}(K_n)
  \to \cellC{*}(K_n)\otimes \cellC{*}(K_n)$, $n\geq 2$, satisfying the
  following conditions:
  \begin{itemize}
  \item \textbf{Compatibility}:\signissue
    \begin{equation}\label{eq:MDiagCompat}
    \MDiag^n\circ \phi^n_{i,j}=\begin{cases}
      (\phi^n_{1,j}\otimes \phi^n_{1,j})\circ
      (\MDiag^{j}\otimes\MDiag^{n+1-j}) & {\text{if $i=1$}}\\
      (\phi^n_{i,j}\otimes \phi^n_{i,j})\circ
      (\AsDiag^{j-i+1}\otimes\MDiag^{n+i-j}) & {\text{if $i>1$.}} 
    \end{cases}
    \end{equation}
  \item \textbf{Non-degeneracy}: $\MDiag^2\co \cellC{0}(K_2)\to
    \cellC{0}(K_2)\otimes \cellC{0}(K_2)$ is the standard isomorphism
    ${\Ring}\cong {\Ring}\rotimes{\Ring}{\Ring}$.
  \end{itemize}

  Equivalently, we can view $\MDiag$ as a formal linear combination of
  pairs of trees $\TrMDiag_n\in \Trees_n\otimes\Trees_n$ in dimension $n-2$,
  corresponding to the cells $\MDiag_n(\corolla{n})$. For module
  diagonals, the leftmost strands are distinguished.
  These trees must satisfy:
  \begin{itemize}
  \item \textbf{Compatibility}: 
    \begin{equation}
      \label{eq:MDiagCompatTree}
      \partial \TrMDiag_n = \sum_{k+\ell=n+1}\Bigl(\TrMDiag_k\circ_1\TrMDiag_\ell+\sum_{i=2}^k\TrMDiag_k\circ_i\TrDiag_\ell\Bigr).
    \end{equation}
  \item \textbf{Non-degeneracy}: $\TrMDiag_2=\corolla{2}\otimes\corolla{2}$.
  \end{itemize}
  We will sometimes call the data $\{\TrMDiag_n\}$ a \emph{module tree diagonal}. 

  \emph{Left module (tree) diagonals} are defined similarly, except that the
  rightmost strands are distinguished. That is, the compatibility
  condition~(\ref{eq:MDiagCompatTree}) is replaced by
  \[
  \partial \TrMDiag_n = \sum_{k+\ell=n+1}\Bigl(\TrMDiag_k\circ_k\TrMDiag_\ell+\sum_{i=1}^{k-1}\TrMDiag_k\circ_i\TrDiag_\ell\Bigr).
  \]
  
  By a \emph{module (tree) diagonal} we will always mean a right module (tree) diagonal.
\end{definition}
The proof of the equivalence of the above two definitions of 
module diagonals is similar to the
proof of Lemma~\ref{lem:DiagonalCells}. Schematically, the
compatibility condition for a right module diagonal is:
\[
\bdy\Biggl(
\mathcenter{
  \tikzsetnextfilename{def-ModuleDiagonal-1}
\begin{tikzpicture}[smallpic]
\node at (0,0) (tl) {};
\node at (1,0) (tc) {};
\node at (1.5,0) (tr) {};
\node at (2.5,0) (trr) {};
\node at (0,-1) (DeltaL) {$\TrMDiag$};
\node at (1.5,-1) (DeltaR) {$\TrMDiag$};
\node at (0,-2) (bl) {};
\node at (1.5,-2) (br) {};
\draw[moda] (tl) to (DeltaL);
\draw[moda] (tr) to (DeltaR);
\draw [moda](DeltaL) to (bl);
\draw [moda](DeltaR) to (br);
\draw[taa] (tc) to (DeltaL);
\draw[taa] (trr) to (DeltaR);
\end{tikzpicture}}\Biggr)
=
\left(\rule{0cm}{1.15cm}\right.
\mathcenter{
  \tikzsetnextfilename{def-ModuleDiagonal-2}
\begin{tikzpicture}[smallpic]
\node at (0,0) (tl) {};
\node at (.5,0) (tcl) {};
\node at (1,0) (tc) {};
\node at (1.5,0) (tr) {};
\node at (2,0) (trrl) {};
\node at (2.5,0) (trr) {};
\node at (0,-1) (DeltaL) {$\TrMDiag$};
\node at (1.5,-1) (DeltaR) {$\TrMDiag$};
\node at (0,-2) (DeltaL2) {$\TrMDiag$};
\node at (1.5,-2) (DeltaR2) {$\TrMDiag$};
\node at (0,-3) (bl) {};
\node at (1.5,-3) (br) {};
\draw[moda] (tl) to (DeltaL);
\draw[moda] (tr) to (DeltaR);
\draw[moda] (DeltaL) to (DeltaL2);
\draw[moda] (DeltaR) to (DeltaR2);
\draw [moda](DeltaL2) to (bl);
\draw [moda](DeltaR2) to (br);
\draw[taa] (tcl) to (DeltaL);
\draw[taa] (trrl) to (DeltaR);
\draw[taa] (tc) to (DeltaL2);
\draw[taa] (trr) to (DeltaR2);
\end{tikzpicture}}
\left.\rule{0cm}{1.15cm}\right)
+
\left(\rule{0cm}{1.15cm}\right.
\mathcenter{
  \tikzsetnextfilename{def-ModuleDiagonal-3}
\begin{tikzpicture}[smallpic]
\node at (0,0) (tl) {};
\node at (.5,0) (tcl) {};
\node at (1.5,0) (tcr) {};
\node at (1,0) (tc) {};
\node at (2,0) (tr) {};
\node at (3.5,0) (trrr) {};
\node at (2.5,0) (trrl) {};
\node at (3,0) (trr) {};
\node at (.5,-1) (ADL) {$\TrDiag$};
\node at (2.5,-1) (ADR) {$\TrDiag$};
\node at (0,-2) (DeltaL) {$\TrMDiag$};
\node at (2,-2) (DeltaR) {$\TrMDiag$};
\node at (0,-3) (bl) {};
\node at (2,-3) (br) {};
\draw[moda] (tl) to (DeltaL);
\draw[moda] (tr) to (DeltaR);
\draw [moda](DeltaL) to (bl);
\draw [moda](DeltaR) to (br);
\draw[taa] (tcl) to (DeltaL);
\draw[taa] (tcr) to (DeltaL);
\draw[taa] (trrl) to (DeltaR);
\draw[taa] (trrr) to (DeltaR);
\draw[taa] (tc) to (ADL);
\draw[taa] (trr) to (ADR);
\draw[alga] (ADL) to (DeltaL);
\draw[alga] (ADR) to (DeltaR);
\end{tikzpicture}}
\left.\rule{0cm}{1.15cm}\right)
\]
(For a left module diagonal one reflects the pictures horizontally.)

\begin{figure}
  \centering
  \includegraphics{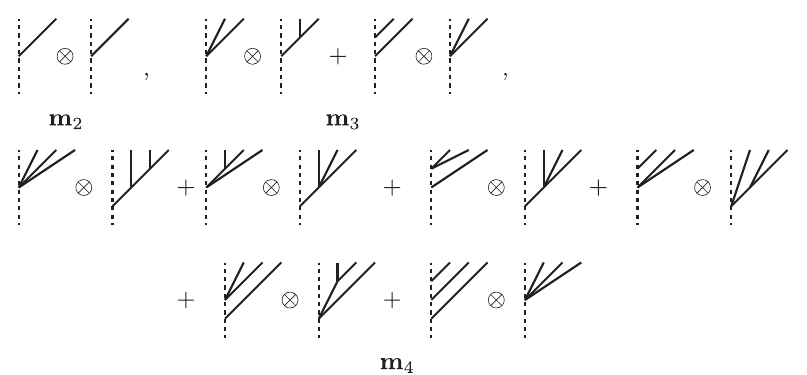}
  \caption[The first three terms in a module diagonal associated
      to a module diagonal primitive]{\textbf{The first three terms in a module diagonal associated
      to a module diagonal primitive.} This is the module diagonal
    associated to the primitive in Figure~\ref{fig:mprim-terms}.}
  \label{fig:prim-to-diag}
\end{figure}
The first few terms in a particular module tree diagonal are shown in
Figure~\ref{fig:prim-to-diag}.

\begin{example}\label{eg:M-diag}
  Any associahedron diagonal can be viewed as a module diagonal; i.e.,
  given an associahedron diagonal $\AsDiag$, the collection of chain
  maps $\MDiag^n=\AsDiag^n$ satisfies the hypotheses of
  Definition~\ref{def:ModuleDiagonal}.
\end{example}

\subsection{Module diagonal primitives}
The goals of this section are to define module diagonal primitives and
relate them to module diagonals.

\begin{definition}
  The {\em leftmost  strand} of a rooted, planar tree is the path
  which is composed of a sequence of edges
  $e_1,\dots,e_n$, where $e_n$ is the terminal edge (incident to the
  output leaf), $e_i$ is the leftmost edge into the initial vertex of
  $e_{i+1}$, and $e_1$ is edge incident to the leftmost input leaf.
\end{definition}

\begin{figure}
  \centering
  \includegraphics{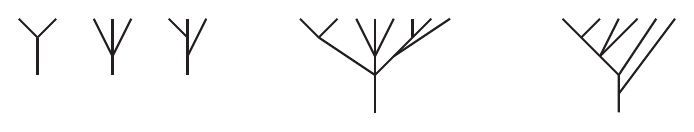}
  \caption[Joining trees]{\textbf{Joining trees.} Left: trees $T_1$, $T_2$, and
    $T_3$. Center: the root joining $\RootJoin(T_1,T_2,T_3)$ of
    $T_1,T_2,T_3$. Right: the left joining $\LeftJoin(T_1,T_2,T_3)$ of
    $T_1,T_2,T_3$.}
  \label{fig:joiningTrees}
\end{figure}

\begin{definition}\label{def:join-trees}
  Given a sequence $S_1,\dots, S_n$ of trees, the \emph{root joining}
  $\RootJoin(S_1,\dots,S_n)$ of $S_1,\dots,S_n$ is obtained by joining
  the outputs of $S_1,\dots,S_n$ into a single new node. That is, the
  root joining is 
  \[
  \RootJoin(S_1,\dots,S_n)= ((((\corolla{n}\circ_n
  S_n)\circ_{n-1}S_{n-1})\circ_{n-2})\cdots) \circ_1 S_1=\corolla{n}\circ(S_1,\dots,S_n).
  \]
  See Figure~\ref{fig:joiningTrees}.

  Given a sequence $T_1,\dots, T_n$ of trees, the \emph{left joining}
  $\LeftJoin(T_1,\dots,T_n)$ of $T_1,\dots, T_n$ is given by
  attaching the output of $T_i$ (for $i=1,\dots,n-1$) to the leftmost
  input of $T_{i+1}$. That is, 
  \[
  \LeftJoin(T_1,\dots,T_n)=T_n\circ_{1}(T_{n-1}\circ_{1}(\cdots\circ_{1}T_1)).
  \]
  Again, see Figure~\ref{fig:joiningTrees}.
\end{definition}

When working with primitives, we will be interested in pairs of trees
$(T,S)\in\Trees_{n}\otimes \Trees_{n-1}$.

\begin{definition}\label{def:shift-join}
  Given $(T_1,S_1)\in\Trees_n\otimes\Trees_{n-1}$ and
  $(T_2,S_2)\in\Trees_m\otimes\Trees_m$, define a shifted version of
  composition by
  $(T_1,S_1)\circ'_i(T_2,S_2)=(T_1\circ_{i+1} T_2,S_1\circ_i S_2)$, and
  let
  $(T_1,S_1)\circ'
  (T_2,S_2)=\sum_{i=1}^{n-1}(T_1,S_1)\circ'_i(T_2,S_2)$.

  For $n \ge 2$, given a sequence
  $(T_1,S_1),(T_2,S_2),\dots,(T_n,S_n)$ of pairs of trees, with each
  $T_i$ having one more input than $S_i$, define the
  \emph{left-root joining} of the sequence to be
  \[
  \LRjoin((T_1,S_1),\dots,(T_n,S_n))=\LeftJoin(T_1,\dots,T_n)\otimes\RootJoin(S_1,\dots,S_n),
  \]
  where the left join $\LeftJoin$ and root join $\RootJoin$ are
  defined in Definition~\ref{def:join-trees}.
  Extend $\LRjoin$ multi-linearly to a function
  \[
  \LRjoin\co (\Trees\otimes\Trees)^{\otimes n}\to (\Trees\otimes\Trees)
  \]
  for $n \ge 2$.
  We define $\LRjoin$ to vanish on
  $(\Trees\otimes\Trees)^{\otimes 1}$ and on
  $(\Trees\otimes\Trees)^{\otimes 0}$.
\end{definition}

\begin{definition}\label{def:M-prim}
  Fix an associahedron tree diagonal~$\TrDiag$.
  A \emph{(right) module diagonal primitive compatible with $\TrDiag$}
  consists of a linear combination of trees for $n \ge 2$
  \begin{equation}\label{eq:M-prim-1}
  \TrPMDiag_n=\sum_{(S,T)}n_{S,T}(S,T)\in \Trees_n\otimes \Trees_{n-1},
  \end{equation}
  in dimension $n-2$, satisfying the following conditions:
  \begin{itemize}
  \item \textbf{Compatibility}:
    \begin{equation}
      \label{eq:M-prim-compat}
      \partial \TrPMDiag = 
      \LRjoin(\TrPMDiag^{\otimes\bullet})+\TrPMDiag\circ'\TrDiag.
    \end{equation}
    Here, $\TrPMDiag=\sum_n\TrPMDiag_n\in\prod_n\Trees_n\otimes\Trees_{n-1}$ and $\TrPMDiag^{\otimes\bullet} = \sum_{n \ge 0}\TrPMDiag^{\otimes n}$.\signissue
  \item \textbf{Non-degeneracy}: $\TrPMDiag_2=\corolla{2}\otimes \IdTree$.
  \end{itemize}
\end{definition}

(The reader familiar with \DD\ bimodules may want to skip ahead to
Definition~\ref{def:one-side-DT} for the motivation for module
diagonal primitives.)

\begin{figure}
  \centering
  %Font is 12 point.
  \includegraphics{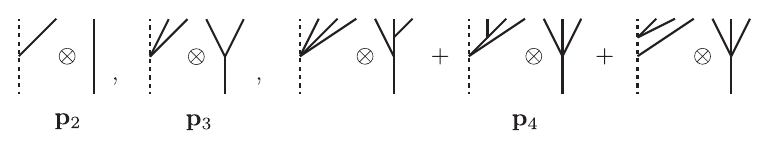}
  \caption[The first three terms in a module diagonal primitive]{\textbf{The first three terms in a module diagonal
      primitive.} The terms $\TrPMDiag_2$, $\TrPMDiag_3$, and
    $\TrPMDiag_4$ are shown. This primitive is compatible with the
    associahedron diagonal from Figure~\ref{fig:diag-cells}.}
  \label{fig:mprim-terms}
\end{figure}

Schematically, the module diagonal primitive compatibility condition is:
\[
\bdy\Biggl(
\mathcenter{
  \tikzsetnextfilename{def-M-prim-1}
\begin{tikzpicture}[smallpic]
  \node at (2,1) (tr) {};
  \node at (2.5,1) (tc) {};
  \node at (2,0) (PR) {$\TrMPrim$};
  \node at (2,-1) (br) {};
  \draw[taa] (tc) to (PR);
  \draw[moda] (tr) to (PR);
  \draw[moda] (PR) to (br);
\end{tikzpicture}}
\mathcenter{
  \tikzsetnextfilename{def-M-prim-2}
\begin{tikzpicture}[smallpic]
  \node at (1.5,1) (tl) {};
  \node at (1.5,0) (PL) {$\TrMPrim$};
  \node at (1.5,-1) (bl) {};
  \draw[taa] (tl) to (PL);
  \draw[alga] (PL) to (bl);
\end{tikzpicture}}
\Biggr)
=
\mathcenter{
  \tikzsetnextfilename{def-M-prim-3}
  \begin{tikzpicture}[smallpic]
  \node at (1,1) (tr) {};
  \node at (2,1) (tc) {};
  \node at (2.5,1) (tcr) {};
  \node at (1,0) (PR) {$\TrMPrim$};
  \node at (1,-1) (dotsr) {$\vdots$};
  \node at (1,-2) (PR2) {$\TrMPrim$};
  \node at (1,-4) (br) {};
  \draw[taa] (tc) to (PR);
  \draw[taa] (tcr) to (PR2);
  \draw[moda] (tr) to (PR);
  \draw[moda] (PR) to (dotsr);
  \draw[moda] (dotsr) to (PR2);
  \draw[moda] (PR2) to (br);
  \end{tikzpicture}
}
\mathcenter{
  \tikzsetnextfilename{def-M-prim-4}
  \begin{tikzpicture}[smallpic]
  \node at (-2,1) (tl) {};
  \node at (0,1) (tlr) {};
  \node at (-2,-1) (PL) {$\TrMPrim$};
  \node at (-1,-1) (dotsl) {$\cdots$};
  \node at (0,-1) (PL2) {$\TrMPrim$};
  \node at (-1,-2) (mul) {$\corolla{}$};
  \node at (-1,-4) (bl) {};
  \draw[taa] (tl) to (PL);
  \draw[taa] (tlr) to (PL2);
  \draw[alga] (PL) to (mul);    
  \draw[alga] (PL2) to (mul);    
  \draw[alga] (mul) to (bl);    
  \end{tikzpicture}
}
+
\mathcenter{
  \tikzsetnextfilename{def-M-prim-5}
  \begin{tikzpicture}[smallpic]
  \node at (0,1) (tr) {};
  \node at (1,1) (tc) {};
  \node at (.5,1) (tcl) {};
  \node at (1.5,1) (tcr) {};
  \node at (.5,0) (DeltaR) {$\TrDiag$};
  \node at (0,-1) (PR) {$\TrMPrim$};
  \node at (0,-2) (br) {};
  \draw[taa] (tcl) to (PR);
  \draw[taa] (tcr) to (PR);
  \draw[taa] (tc) to (DeltaR);
  \draw[alga] (DeltaR) to (PR);
  \draw[moda] (tr) to (PR);
  \draw[moda] (PR) to (br);
  \end{tikzpicture}
}
\mathcenter{
  \tikzsetnextfilename{def-M-prim-6}
  \begin{tikzpicture}[smallpic]
  \node at (-1.5,1) (tl) {};
  \node at (-2.5,1) (tll) {};
  \node at (-.5,1) (tlr) {};
  \node at (-1.5,0) (DeltaL) {$\TrDiag$};
  \node at (-1.5,-1) (PL) {$\TrMPrim$};
  \node at (-1.5,-2) (bl) {};
  \draw[taa] (tll) to (PL);
  \draw[taa] (tlr) to (PL);
  \draw[taa] (tl) to (DeltaL);
  \draw[alga] (DeltaL) to (PL);
  \draw[alga] (PL) to (bl);    
  \end{tikzpicture}
}
\]
The first few terms in one module diagonal primitive, compatible with
the associahedron diagonal from Figure~\ref{fig:diag-cells}, are shown
in Figure~\ref{fig:mprim-terms}.

\begin{example}\label{eg:explicit-prim}
  Given a tree $S$, a \emph{left inflation} of $S$ is a tree obtained
  from $S$ by adding one more input, to the left of the other inputs
  of $S$, and connecting that input to some (existing) internal vertex on the
  left-most strand of $S$.
  Let $\Delta=\{(S,T)\}$ be the set of pairs
  of trees appearing in the diagonal in
  Example~\ref{eg:explicit-diag}. One can show that
  \[
    \{(S',T)\mid \exists (S,T)\in\Delta\text{ such that }S'\text{ is a left inflation of }S\}
  \]
  is a module diagonal primitive compatible with $\Delta$.
  This is the module diagonal primitive shown in
  Figure~\ref{fig:mprim-terms}. (See also the proof of
  Lemma~\ref{lem:prim-unique}.)
\end{example}

\begin{lemma}\label{lem:join-diff}
  Given $x_1\otimes\cdots\otimes x_n\in (\Trees\otimes\Trees)^{\otimes n}$,
  the operation $\LRjoin$ satisfies
  \begin{multline*}
  \bdy(\LRjoin(x_1\otimes\cdots\otimes x_n))=\sum_{i=1}^n \LRjoin(x_1\otimes\cdots\otimes \bdy(x_i)\otimes\cdots\otimes x_n) \\
  + \sum_{j=2}^{n-1}\sum_{i=0}^{n-j}\LRjoin(x_1\otimes\cdots\otimes x_i\otimes \LRjoin(x_{i+1}\otimes\cdots\otimes x_{i+j})\otimes x_{i+j+1}\otimes\cdots\otimes x_n).
  \end{multline*}
  The operation $\circ'$ satisfies $\bdy(x\circ'y)=\bdy(x)\circ'
  y+x\circ'\bdy(y)$.
  Finally,
  \[ 
    \LRjoin(x_1\otimes\dots\otimes x_n)\circ'y=\sum_{i=1}^n \LRjoin(x_1,\dots,x_{i-1}, x_i\circ' y,x_{i+1},\dots,x_n).
  \]
\end{lemma}
\begin{proof}
  The claim about $\bdy\circ\LRjoin$ is immediate from the definition: the
  second sum comes from differentiating the corolla in the root
  joining. The Leibniz rule for $\circ'$ and distributivity of
  $\circ'$ and $\LRjoin$ are also immediate from the definitions.
\end{proof}

Recall that Proposition~\ref{prop:prim-exist} asserts the existence of module diagonal primitives.
\begin{proof}[Proof of Proposition~\ref{prop:prim-exist}]
  As in Lemma~\ref{lem:ExistAssociahedronDiagonal}, we proceed by
  induction on the number of inputs. The case $n=2$ is determined by
  the non-degeneracy condition. For the case $n=3$, the compatibility
  condition is
  \[
    \bdy\TrPMDiag_3=
    \mathcenter{
  \tikzsetnextfilename{lem-join-diff-1}
      \begin{tikzpicture}[smallpic]
        % First left tree
        \draw[->, dashed] (1,0) to (1,-1);
        \draw (1.25,0) to (1,-.25);
        \draw (1.5, 0) to (1,-.5);
      \end{tikzpicture}
    }
    \otimes\ 
    \mathcenter{
  \tikzsetnextfilename{lem-join-diff-2}
      \begin{tikzpicture}[smallpic]
        %First right tree
        \draw (-.25,0) to (0,-.5);
        \draw (.25,0) to (0,-.5);
        \draw[->] (0,-.5) to (0,-1);
      \end{tikzpicture}
    }\ +\ 
    \mathcenter{
  \tikzsetnextfilename{lem-join-diff-3}
      \begin{tikzpicture}[smallpic]
        % Second left tree
        \node at (.9,0) {};
        \draw[->, dashed] (1,0) to (1,-1);
        \draw (1.25,0) to (1.25,-.25);
        \draw (1.5, 0) to (1,-.5);
      \end{tikzpicture}
    }
    \otimes
    \mathcenter{
  \tikzsetnextfilename{lem-join-diff-4}
      \begin{tikzpicture}[smallpic]
        %Second right tree
        \draw (-.25,0) to (0,-.5);
        \draw (.25,0) to (0,-.5);
        \draw[->] (0,-.5) to (0,-1);
      \end{tikzpicture}.
    }
  \]
  (The two terms come from the two terms on the right of
  Formula~(\ref{eq:M-prim-compat}), respectively.) This is solved by
  \[
    \TrPMDiag_3=
    \mathcenter{
  \tikzsetnextfilename{lem-join-diff-5}
      \begin{tikzpicture}[smallpic]
        % Left tree
        \draw[->, dashed] (1,0) to (1,-1);
        \draw (1.25,0) to (1,-.5);
        \draw (1.5, 0) to (1,-.5);
      \end{tikzpicture}
    }
    \otimes\ 
    \mathcenter{
  \tikzsetnextfilename{lem-join-diff-6}
      \begin{tikzpicture}[smallpic]
        % Right tree
        \draw (-.25,0) to (0,-.5);
        \draw (.25,0) to (0,-.5);
        \draw[->] (0,-.5) to (0,-1);
      \end{tikzpicture}.
    }
  \]

  If we have constructed
  primitives with fewer than $n$ leaves, $n>3$,
  it suffices to verify that the expression on the right of
  Formula~(\ref{eq:M-prim-compat}) is a cycle: the fact that it is a
  boundary then follows from contractibility of $K_n\times K_{n-1}$
  and the fact that a primitive lies in dimension
  $n-2$. Using Lemma~\ref{lem:join-diff}, we have
  \begin{align*}
  \bdy(\LRjoin(\TrPMDiag^{\otimes\bullet})+\TrPMDiag\circ'\TrDiag)
  &=
  \LRjoin(\TrPMDiag^{\otimes\bullet} \otimes \LRjoin(\TrPMDiag^{\otimes\bullet})\otimes\TrPMDiag^{\otimes\bullet})
  + 
  \LRjoin(\TrPMDiag^{\otimes\bullet} \otimes \bdy(\TrPMDiag)\otimes\TrPMDiag^{\otimes\bullet})\\
  &\qquad\qquad+ 
  \bdy(\TrPMDiag)\circ'\TrDiag
  +\TrPMDiag\circ'\bdy(\TrDiag).
  \end{align*}
  Applying the primitive compatibility condition
  (Formula~(\ref{eq:M-prim-compat})) to the second and third terms and
  the associahedron diagonal compatibility condition
  (Formula~(\ref{eq:DiagonalCellCompatibilityTree})) to the fourth
  term gives
  \begin{multline*}
  \LRjoin(\TrPMDiag^{\otimes\bullet}\otimes \LRjoin(\TrPMDiag^{\otimes\bullet})\otimes\TrPMDiag^{\otimes\bullet})
  + 
  \LRjoin(\TrPMDiag^{\otimes\bullet} \otimes\LRjoin(\TrPMDiag^{\otimes\bullet})\otimes\TrPMDiag^{\otimes\bullet})
  +
  \LRjoin(\TrPMDiag^{\otimes\bullet}\otimes (\TrPMDiag\circ' \TrDiag)\otimes\TrPMDiag^{\otimes\bullet})\\
  +
  \LRjoin(\TrPMDiag^{\otimes\bullet})\circ'\TrDiag
  +
  (\TrPMDiag\circ'\TrDiag)\circ'\TrDiag
  +
  \TrPMDiag\circ'(\TrDiag\circ\TrDiag).    
  \end{multline*}
  The first two terms cancel, the third and fourth terms cancel, and
  the fifth and sixth terms cancel. The result follows.\signissue
\end{proof}

We can construct a module diagonal from a primitive.
\begin{definition}\label{def:assoc-mod-diag}
  Given a sequence $(S_1,T_1),\dots,(S_k,T_k)$ of pairs of trees,
  define a modified left-root joining by
  \[
  \LRjoin'((S_1,T_1),\dots,(S_k,T_k)) =
    \LeftJoin(S_1,\dots,S_k) \otimes
      \RootJoin(\IdTree,T_1,\dots,T_k).
  \]
  If each $S_i$ has one more input than the corresponding $T_i$, then
  in the result of $\LRjoin'$ the two trees have the same number of
  inputs. Extend $\LRjoin'$ multi-linearly to a function
  \[
  \LRjoin'\co (\Trees\otimes\Trees)^{\otimes n}\to (\Trees\otimes\Trees)
  \]
  for $n \ge 1$.
  Define $\LRjoin'$ to vanish on $(\Trees\otimes\Trees)^{\otimes 0}$.
  Then, if $\TrPMDiag = \sum_n \TrPMDiag_n$ is a linear combination of
  trees as in Equation~\eqref{eq:M-prim-1}, the \emph{module tree
    diagonal associated to $\TrPMDiag$} is
  \begin{equation}\label{eq:prim-to-diag}
    \TrMDiag_\TrPMDiag \coloneqq \LRjoin'(\TrPMDiag^{\otimes\bullet}).
  \end{equation}
\end{definition}

Schematically, the associated module tree diagonal is:
\[
\mathcenter{
  \tikzsetnextfilename{def-assoc-mod-diag-1}
  \begin{tikzpicture}[smallpic]
    \node at (0,0) (tl) {};
    \node at (.5,0) (tc) {};
    \node at (1,0) (tr) {};
    \node at (1.5,0) (trr) {};
    \node at (0,-1) (DeltaL) {$\TrMDiag$};
    \node at (1,-1) (DeltaR) {$\TrMDiag$};
    \node at (0,-2) (bl) {};
    \node at (1,-2) (br) {};
    \draw[moda] (tl) to (DeltaL);
    \draw[moda] (tr) to (DeltaR);
    \draw [moda](DeltaL) to (bl);
    \draw [moda](DeltaR) to (br);
    \draw[taa] (tc) to (DeltaL);
    \draw[taa] (trr) to (DeltaR);
  \end{tikzpicture}}
\coloneqq
\mathcenter{
  \tikzsetnextfilename{def-assoc-mod-diag-2}
  \begin{tikzpicture}[smallpic]
    \node at (2,0) (tr) {};
    \node at (2.5,0) (trr) {};
    \node at (3,0) (trrr) {};
    \node at (2,-1) (PR) {$\TrMPrim$};
    \node at (2,-2) (rdots) {$\vdots$};
    \node at (2,-3) (PR2) {$\TrMPrim$};
    \node at (2,-5) (br) {};
    \draw[moda] (tr) to (PR);
    \draw[moda] (PR) to (rdots);
    \draw[moda] (rdots) to (PR2);
    \draw[moda] (PR2) to (br);
    \draw[taa] (trr) to (PR);
    \draw[taa] (trrr) to (PR2);
  \end{tikzpicture}
}
\mathcenter{
  \tikzsetnextfilename{def-assoc-mod-diag-3}
  \begin{tikzpicture}[smallpic]
    \node at (.5,0) (tll) {};
    \node at (1,0) (tl) {};
    \node at (1.5,0) (tc) {};
    \node at (2,0) (tcr) {};
    \node at (1,-1) (PL) {$\TrMPrim$};
    \node at (1,-2) (ldots) {$\vdots$};
    \node at (1,-3) (PL2) {$\TrMPrim$};
    \node at (.5,-4) (mul) {$\corolla{}$};
    \node at (.5,-5) (bll) {};
    \draw[moda] (tll) to (mul);
    \draw[moda] (mul) to (bll);
    \draw[taa] (tc) to (PL);
    \draw[taa] (tcr) to (PL2);
    \draw[alga] (PL) to (mul);
    \draw[alga] (PL2) to (mul);
  \end{tikzpicture}
}.
\]
An example is shown in Figure~\ref{fig:prim-to-diag}. In that figure,
the trees in each $\TrMDiag_n$ are listed in order of the number $k$
of terms in $\TrPMDiag^{\otimes k}$ which are joined together.

\begin{lemma}\label{lem:mod-join-diff}
  Given $x_1\otimes\cdots\otimes x_n\in (\Trees\otimes\Trees)^{\otimes
    n}$,
  the operation $\LRjoin'$ satisfies
  \begin{align*}
  \bdy(\LRjoin'(x_1\otimes\cdots\otimes x_n))
   &=\sum_{i=1}^n \LRjoin'(x_1\otimes\cdots\otimes \bdy(x_i)\otimes\cdots\otimes x_n) \\
  &+ \sum_{i=1}^{n-1} \LRjoin'(x_1\otimes\dots\otimes x_i) \circ_1 \LRjoin'(x_{i+1}\otimes\dots\otimes x_n)\\
  &+ \sum_{j=2}^n\sum_{i=0}^{n-j}\LRjoin'(x_1\otimes\cdots\otimes \LRjoin(x_{i+1}\otimes\cdots\otimes x_{i+j})\otimes\cdots\otimes x_n).
  \end{align*}
\end{lemma}

\begin{proof}
  This is immediate from the definitions. The second and third terms come from
  differentiating the corolla in the root joining, depending on whether
  the new edge is on the left or not.
\end{proof}

Lemma~\ref{lem:prim-gives-diag} asserts that the associated module
tree diagonal is, in fact, a module tree diagonal.

\begin{proof}[Proof of Lemma~\ref{lem:prim-gives-diag}]
  Let $\TrMDiag$ be the module tree diagonal associated
  to~$\TrPMDiag$. Non-degeneracy of $\TrMDiag$ is immediate. From
  Lemma~\ref{lem:mod-join-diff} and the module diagonal primitive
  compatibility Equation~\eqref{eq:M-prim-compat}, we have
  \begin{align*}
    \bdy \TrMDiag  &= \LRjoin'(\TrPMDiag^{\otimes\bullet} \otimes
      \bdy\TrPMDiag \otimes \TrPMDiag^{\otimes\bullet})
      + \LRjoin'(\TrPMDiag^{\otimes\bullet}) \circ_1
                     \LRjoin'(\TrPMDiag^{\otimes\bullet})
      + \LRjoin'(\TrPMDiag^{\otimes\bullet}\otimes\LRjoin(\TrPMDiag^{\otimes\bullet})
         \otimes\TrPMDiag^{\otimes\bullet})\\
    &= \LRjoin'(\TrPMDiag^{\otimes\bullet} \otimes (\TrPMDiag \circ' \TrDiag)
         \otimes \TrPMDiag^{\otimes\bullet})
      + \LRjoin'(\TrPMDiag^{\otimes\bullet}) \circ_1
                     \LRjoin'(\TrPMDiag^{\otimes\bullet})\\
    &= \sum_{i \ge 2}\TrMDiag \circ_i \TrDiag + \TrMDiag \circ_1 \TrMDiag,
  \end{align*}
  agreeing with the module diagonal compatibility
  Equation~\eqref{eq:MDiagCompatTree}.
\end{proof}

Module diagonal primitives, if they exist, are unique, and satisfy the
module primitive compatibility condition:
\begin{lemma}\label{lem:prim-unique}
  If $\TrPMDiag$ and $\TrPMDiag'$ are module diagonal primitives with
  $\LRjoin'(\TrPMDiag)=\LRjoin'(\TrPMDiag')$ then
  $\TrPMDiag=\TrPMDiag'$. Further, if $\TrMDiag$ is a module tree
  diagonal so that $\TrMDiag=\LRjoin'(\TrPMDiag)$ for some collection
  of pairs of trees $\TrPMDiag$ then $\TrPMDiag$
  is a module diagonal primitive.
\end{lemma}
\begin{proof}
  For the first statement, note that there is a bijection between
  pairs of trees $(S,T)$ in $\TrPMDiag_n$ and pairs of trees $(S',T')$
  in $\TrMDiag_n$ for which
  the trunk vertex $T'$ (the internal vertex closest to the root)
  has valence $3$. Given $(S',T')$, the corresponding tree $(S,T)$ has
  $S=S'$ and $T$ obtained from $T'$ by deleting the left-most
  edge. (Compare Figures~\ref{fig:mprim-terms}
  and~\ref{fig:prim-to-diag} and Example~\ref{eg:explicit-prim}.)

  For the second statement, if $\TrMDiag$ satisfies the module
  diagonal compatibility condition then, as in the proof of
  Lemma~\ref{lem:prim-gives-diag}, we must have
  \begin{multline*}
    \LRjoin'(\TrPMDiag^{\otimes\bullet} \otimes
    \bdy\TrPMDiag \otimes \TrPMDiag^{\otimes\bullet})
    + \LRjoin'(\TrPMDiag^{\otimes\bullet}) \circ_1
    \LRjoin'(\TrPMDiag^{\otimes\bullet})
    + \LRjoin'(\TrPMDiag^{\otimes\bullet}\otimes\LRjoin(\TrPMDiag^{\otimes\bullet})
    \otimes\TrPMDiag^{\otimes\bullet})\\
    = \LRjoin'(\TrPMDiag^{\otimes\bullet} \otimes (\TrPMDiag \circ' \TrDiag)
    \otimes \TrPMDiag^{\otimes\bullet})
    + \LRjoin'(\TrPMDiag^{\otimes\bullet}) \circ_1
    \LRjoin'(\TrPMDiag^{\otimes\bullet}).
  \end{multline*}
  Thus,
  \[
    \LRjoin'(\TrPMDiag^{\otimes\bullet} \otimes
    \bdy\TrPMDiag \otimes \TrPMDiag^{\otimes\bullet})
    + \LRjoin'(\TrPMDiag^{\otimes\bullet}\otimes\LRjoin(\TrPMDiag^{\otimes\bullet})
    \otimes\TrPMDiag^{\otimes\bullet})
    +\LRjoin'(\TrPMDiag^{\otimes\bullet} \otimes (\TrPMDiag \circ' \TrDiag)
    \otimes \TrPMDiag^{\otimes\bullet})
    =0.    
  \]
  Restricting to the terms where the trunk vertex
  of the left tree has valence $3$ then gives
  \[
    \LRjoin'(
    \bdy\TrPMDiag)
    + \LRjoin'(\LRjoin(\TrPMDiag^{\otimes\bullet}))
    +\LRjoin'(\TrPMDiag \circ' \TrDiag)
    =0.
  \]
  Considering the terms where the trunk vertex of the right tree also has
  valence~$3$ shows that
  \[
    \bdy\TrPMDiag
    + \LRjoin(\TrPMDiag^{\otimes\bullet})
    +\TrPMDiag \circ' \TrDiag
    =0,
  \]
  as desired. Finally, non-degeneracy of $\TrMDiag$ implies
  non-degeneracy of $\TrPMDiag$.
\end{proof}

\begin{remark}
  It would be nice to have a cellular definition of a module diagonal
  primitive, along the lines of the (non-tree versions of)
  Definitions~\ref{def:AssociahedronDiagonal}
  and~\ref{def:ModuleDiagonal}.
\end{remark}

\begin{remark}
  Fix an associahedron tree diagonal $\TrDiag$. We can define a
  differential algebra $\Alg$ as follows. Let $\Alg_n$ by the free
  $\Ring$-module spanned by the pairs of trees with $n$ inputs (and no
  $2$-valent vertices), and let $\Alg=\prod_{n=1}^\infty \Alg_n$. The
  multiplication on $\Alg$ is given by $(T_1,T_2)\cdot
  (S_1,S_2)=(T_1\circ_1 S_1, T_2\circ_1 S_2)$. The differential on
  $\Alg_n$ is given by
  \[
  \bdy(S,T)=(\bdy(S),T)+(S,\bdy(T))+\sum_{i=2}^n(S,T)\circ_i \TrDiag
  \]
  (where $\bdy(S)$ denotes the usual differential on trees, i.e., the
  sum of all ways of adding an edge to $S$). It is straightforward to
  verify that $\Alg$ is a differential algebra.

  A module diagonal is a rank $1$ type $D$ structure (see
  Section~\ref{sec:typeD}) over $\Alg$
  satisfying the non-degeneracy condition that
  $\delta^1$ of the generator contains $\corolla{2}\otimes\corolla{2}$.

  Similarly, define an $\Ainf$-algebra $\Blg$ as follows. Let $\Blg_n$
  be the free $\Ring$-module spanned by all pairs of trees $(S,T)$ where
  $S$ has $n$ inputs, $T$ has $n-1$ inputs, and neither $S$ nor $T$
  has any $2$-valent vertices. Let $\Blg=\prod_{n=2}^\infty
  \Blg_n$. Define products $\mu_k$ on $\Blg$, $k\geq 2$, by
  \[
  \mu_k((S_1,T_1),\dots,(S_k,T_k))=\LRjoin((S_1,T_1),\dots,(S_k,T_k)).
  \]
  Define a differential on $\Blg_n$ by
  \[
    \bdy(S,T)=(\bdy(S),T)+(S,\bdy(T))+(S,T)\circ'\TrDiag.
  \]
  It is straightforward to verify that $\Blg$ is an $\Ainf$-algebra.

  A module diagonal primitive is a rank $1$ type $D$ structure over $\Blg$.

  There is an $\Ainf$-homomorphism $F\co \Blg\to\Alg$ defined by 
  \[
    F_k((S_1,T_1),\dots,(S_n,T_n))=\LRjoin'((S_1,T_1),\dots,(S_k,T_k)).
  \]
  Viewing $F$ as a rank-one type \DA\ structure $\lsup{\Alg}[F]_{\Blg}$ over
  $\Alg$ and $\Blg$ (see~\cite[Definition 2.2.48]{LOT2}), the
  operation of turning a module diagonal
  primitive into a module tree diagonal is taking the $\DT$ product,
  over $\Blg$, with $\lsup{\Alg}[F]_{\Blg}$ (see~\cite[Section
  2.3.2]{LOT2}), i.e., applying the induction functor associated to
  $F$ to the module diagonal primitive.
\end{remark}

\subsection{Multiplihedra and multiplihedron diagonals}\label{sec:multiplihedron}
The associahedron is relevant to the definition of an
$\Ainf$-algebra.
There is a different contractible CW complex, the
\emph{multiplihedron} $J_n$, which is relevant to morphisms of
$\Ainf$-algebras~\cite{IwaseMimura89:multiplihedron}. (Like the
associahedron, the multiplihedron can be understood in terms of disks
with boundary marked points~\cite{MauWoodward10:multiplihedra}.)

\begin{definition}
\label{def:TransformationTree}
By a \emph{transformation tree} we mean a tree whose edges are colored by one
of two colors, {\em red} and {\em blue}, such that:
\begin{enumerate}
\item The edges adjacent to input leaves are red.
\item The edge adjacent to the output leaf is blue.
\item For each vertex $v$, all of the inputs of $v$ have the same color (red or blue).
\item If a vertex $v$ has a red output, then all of the inputs of $v$ are red;
  if a vertex has blue inputs, then its output is also blue.
\end{enumerate}
Call an internal vertex \emph{red} (respectively \emph{blue}) if all of its
inputs and its output are red (respectively blue), and \emph{purple}
if its inputs are red
but its output is blue. We also require that
\begin{enumerate}[resume]
\item every $2$-valent vertex is purple.
\end{enumerate}
\end{definition}

\begin{figure}
  \centering
  \includegraphics{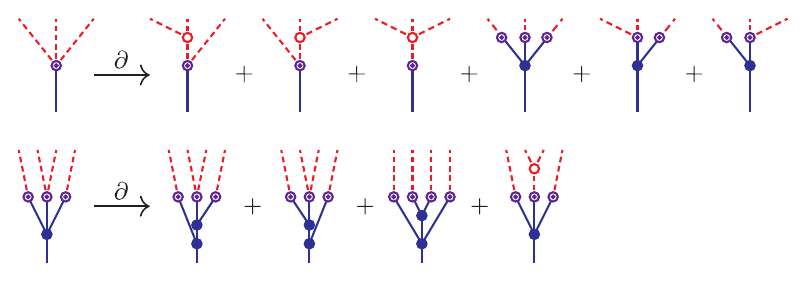}
  \caption[The differential on the multiplihedron]{\textbf{The differential on $\cellC{*}(J_n)$.} Top: the
    differential of a purple corolla. Bottom: the differential of a
    tree with four internal vertices. Blue segments are solid, red segments are
    dashed, blue vertices are solid, red vertices are hollow, and
    purple vertices are solid and hollow.}
  \label{fig:multiplihedron-diff}
\end{figure}

The cells in the multiplihedron are in natural correspondence with
transformation trees.
The differential on $\cellC{*}(J_n)$ is defined as follows.
Given a blue (respectively red) corolla $\bcorolla{n}$ (respectively
$\rcorolla{n}$), define the differential of $\bcorolla{n}$
(respectively $\rcorolla{n}$) as in Section~\ref{sec:associahedron},
by summing over all ways of inserting an edge.
To define the differential of a purple corolla $\pcorolla{n}$, we use
the \emph{blue root joining} $\bRootJoin(S_1,\dots,S_n)$ of a collection of
transformation trees, which is just their root joining at a vertex which is colored blue.
Given a purple corolla $\pcorolla{n}$, define $\bdy(\pcorolla{n})$ to
be
\begin{equation}\label{eq:purple-diff}
\bdy(\pcorolla{n})=\sum_{k+\ell=n+1}\pcorolla{k}\circ\rcorolla{\ell}
+
\sum_{k_1+\cdots+k_\ell=n}
\bRootJoin(\pcorolla{k_1},\dots,\pcorolla{k_{\ell}})
\end{equation}
The differential of a transformation tree $T$ is the sum over the
internal vertices $v$ of $T$ of the result of replacing $v$ by $\bdy(v)$.  See
Figure~\ref{fig:multiplihedron-diff}.

The dimension of a transformation tree~$T$ with $n$~inputs, $r$~red
internal vertices, and $b$~blue internal vertices is
\[
  \dim(T) = n - r - b - 1.
\]

The terms in Formula~\eqref{eq:purple-diff} correspond
to the codimension-1 faces of $J_n$. In particular, for any $1\leq
i<j\leq n$ and any $0=i_0< i_1<\dots<i_\ell=n$ there are
codimension-1 faces $G_{i,j}^n$ and $H_{i_0,\dots,i_\ell}^n$ of $J_n$,
identifications
% XXX: Maybe try to line up $\cong$ marks as well
\begin{align*}
  K_{j-i+1}\times J_{n+i-j} \cong G_{i,j}^n&\subset J_n\\
  J_{i_1-i_0}\times\cdots\times J_{i_{\ell}-i_{\ell-1}}\times K_{\ell}\cong H_{i_1,\dots,i_\ell}^n&\subset J_n
\end{align*}
and corresponding maps
\begin{align*}
  \psi_{i,j}^n\co \cellC{*}(K_{j-i+1})\otimes \cellC{*}(J_{n+i-j})&\to \cellC{*}(J_n)\\
  \xi^n_{i_0,\dots,i_\ell}\co \cellC{*}(J_{i_1-i_0})\otimes\cdots\otimes \cellC{*}(J_{i_{\ell}-i_{\ell-1}})\otimes \cellC{*}(K_\ell)&\to \cellC{*}(J_n).
\end{align*}

In terms of trees, the map $\psi^n_{i,j}$ corresponds to $(S,T)\mapsto T\circ_iS$ (where $T$ is a transformation tree and $S$ is a (red) associahedron tree). The map $\xi^n_{i_0,\dots,i_\ell}$ corresponds to 
\[
(T_1,\dots,T_\ell,S)\mapsto ((S\circ_{\ell}T_\ell)\circ_{\ell-1}\cdots)\circ_1 T_1=S\circ(T_1,\dots,T_\ell)
\]
(where the $T_i$ are transformation trees and $S$ is a (blue)
associahedron tree).

\begin{lemma}
  The multiplihedron $J_n$ is contractible. In particular,
  $H_0(J_n)=\Ring$ and $H_i(J_n)=0$ for $i\neq 0$.
\end{lemma}
\begin{proof}
  This is clear from the construction of $J_n$ as a subspace of
  $\RR^{n-1}$~\cite[Section 2]{IwaseMimura89:multiplihedron}.  Since
  it will be useful to us later, we also give a direct proof.

  Given a (weighted) rooted, planar tree $T$, we can order the
  internal vertices of $T$ by \emph{(pre-order) depth-first search}. If $v$ is
  the trunk vertex of $T$ and the sub-trees feeding into $v$ are
  $T_1,\dots,T_k$ (so $v$ has valence $k+1$) then the depth-first
  search is the sequence of internal vertices of $T$
  \begin{equation}
    \label{eq:DFS}
    \dfs(T)=(v,\dfs(T_1),\dfs(T_2),\dots,\dfs(T_k)).
  \end{equation}
  Given an internal vertex $v'$ of $T$,
  $\depth(v')\in\NN$ denote the index
  (location) of $v'$ in $\dfs(T)$.  See Figure~\ref{fig:dfs}.
  \begin{figure}
    \centering
    % Font is 12 point
    \includegraphics{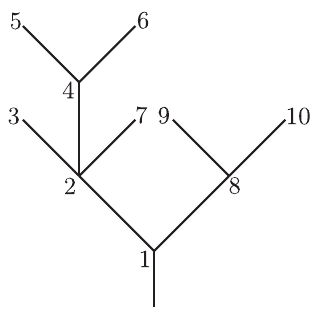}
    \caption[Depth-first search]{\textbf{Depth-first search.} The ordering of the vertices, not the weights of the vertices, is indicated by the numbering.}
    \label{fig:dfs}
  \end{figure}

  Consider the operator $H\co J_n\to J_n$ defined as follows.  Fix
  $T\in J_n$. Using depth-first search, look for the first purple vertex~$v$
  in~$T$, skipping 2-valent vertices whose input is an input to~$T$.
  If $v$ is $2$-valent, contract the input edge into $v$, and call the result
  $H(T)$. If $v$ is not $2$-valent or if there is no such vertex, then $H(T)=0$. 
  
  Consider the filtration $F(T)=\sum_{\{u\in\Vertices(T)\mid u~\text{is not
      blue}\}} (\valence(u)-2)$, where $\valence(u)$ denotes the
  valence of $u$. Clearly, $F(T)\geq 0$ for all $T$, and 
  the differential respects the filtration determined by
  $F$, i.e., if $S$ appears in $\bdy(T)$ then $F(S)\leq F(T)$.
  We claim that if $F(T)>0$, then the terms in 
  \[
    (\partial\circ H + H \circ\partial + \Id)(T)
  \]
  are in filtration level strictly less than $F(T)$. Indeed, the terms
  in $\bdy(H(T))$ and $H(\bdy(T))$ clearly cancel in pairs except:
  \begin{itemize}
  \item Terms in $\bdy(H(T))$ where $H$ contracts the input edge to a
    $2$-valent vertex creating a new purple vertex $w$, and $\bdy$ splits the
    vertex $w$. These occur only if the vertex $v$ in $T$ found by
    depth-first search is $2$-valent. All of the resulting terms in
    $\bdy(H(T))$ have filtration strictly less than $F(T)$ (because they
    have more blue vertices) except the term $T$ itself, which cancels
    against $\Id(T)$.
  \item Terms in $H(\bdy(T))$ where $\bdy(T)$ creates a new $2$-valent
    purple vertex $v'$ which is the first vertex found by depth-first
    search. This occurs only if the first purple vertex $v$ in $T$
    found by depth-first search has valence $>2$, so $H(T)=0$.
    All but
    one of the resulting terms in $H(\bdy(T))$ have filtration less
    than $F(T)$ (because they have more blue vertices); the remaining
    term in $H(\bdy(T))$ is $T$ itself, which cancels with $\Id(T)$.
  \end{itemize}

  If $F(T)=0$, all the non-blue nodes in the tree $T$ are $2$-valent
  (and hence they are purple, and all their inputs are inputs to
  $T$). In that case,
  \[
    (\partial\circ H + H \circ\partial)(T) =0.
  \]
  It follows that $J_n$ is chain homotopy equivalent to the subcomplex
  in filtration level~$0$, which in turn is identified with
  the associahedron with $n$ inputs.  The result now follows from the
  corresponding fact for the associahedron.
\end{proof}

\begin{figure}
  \centering
  \includegraphics{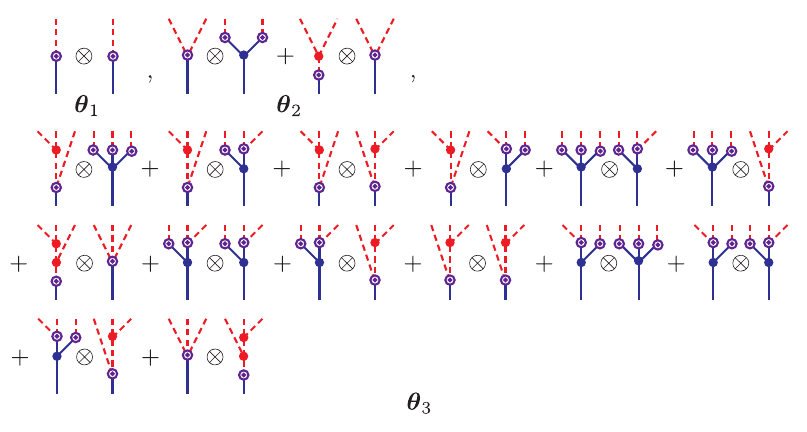}
  \caption[The first three terms in a multiplihedron
      diagonal]{\textbf{The first three terms in a multiplihedron
      diagonal.} This multiplihedron diagonal is compatible with the
    associahedron diagonal from Figure~\ref{fig:diag-cells} (on both sides).}
  \label{fig:mult-diag}
\end{figure}

Tensoring maps of $\Ainf$-algebras will use multiplihedron diagonals:
\begin{definition}\label{def:multiplihedron-diag}
  Fix associahedron diagonals $\AsDiag_1$ and $\AsDiag_2$.  A
  \emph{multiplihedron diagonal} compatible with $\AsDiag_1$ and
  $\AsDiag_2$ consists of a sequence of (degree-preserving) chain maps
  \[
  \{\MulDiag^n\co \cellC{*}(J_n)\to \cellC{*}(J_n)\otimes\cellC{*}(J_n)\}_{n=1}^\infty
  \]
  satisfying the following conditions:
  \begin{itemize}
  \item \textbf{Compatibility under stacking}: 
    \begin{align*}
      \MulDiag^n\circ \psi_{i,j}^n&=(\psi_{i,j}^{n}\otimes\psi_{i,j}^{n})\circ(\AsDiag_1^{j-i+1}\otimes\MulDiag^{n+i-j}) \\
      \MulDiag^n\circ \xi^n_{i_0,\dots,i_\ell}&=
      (\xi^n_{i_0,\dots,i_\ell}\otimes \xi^n_{i_0,\dots,i_\ell})\circ(\MulDiag^{i_1-i_0}\otimes \cdots\otimes \MulDiag^{i_\ell-i_{\ell-1}}\otimes \AsDiag_2^\ell),
    \end{align*}
    with the understanding that the compositions on the right-hand
    side involve shuffling of factors (compare
    Formula~\eqref{eq:as-diag-shuf}).
  \item \textbf{Non-degeneracy}: $\MulDiag^1\co \cellC{0}(J_1)\to
    \cellC{0}(J_1)\otimes\cellC{0}(J_1)$ is the standard isomorphism
    $\Ring\cong \Ring\otimes_\Ring \Ring$. (Note that $J_1$ is a single point.)
  \end{itemize}

  In terms of trees, if $\AsDiag_1$ and $\AsDiag_2$ correspond to the
  associahedron tree diagonals $\TrDiag^1$ and~$\TrDiag^2$, a
  multiplihedron diagonal consists of elements
  $\TrMulDiag_n=\MulDiag^n(\pcorolla{n})\in
  \cellC{*}(J_n)\otimes\cellC{*}(J_n)$ in dimension $n-1$ satisfying:
  \begin{itemize}
  \item \textbf{Compatibility under stacking}: 
  \begin{equation}\label{eq:MulDiagCompatTree}
    \bdy(\TrMulDiag_n)=\sum_{i+j=n+1}\TrMulDiag_i\circ \TrDiag^1_j +
    \sum_k\sum_{m_1+\cdots+m_k=n}\!\!\!\!\TrDiag^2_k\circ (\TrMulDiag_{m_1},\cdots,\TrMulDiag_{m_k}).
  \end{equation}
\item \textbf{Non-degeneracy}:
  $\TrMulDiag_1\in\cellC{*}(J_1)\otimes\cellC{*}(J_1)$ is the (unique)
  pair of $1$-input transformation trees (with $1$ internal vertex
  each).
  \end{itemize}
\end{definition}
(The reader might find it helpful to compare
Formulas~(\ref{eq:purple-diff}) and~(\ref{eq:MulDiagCompatTree}). The
first three terms in a multiplihedron diagonal are shown in Figure~\ref{fig:mult-diag}.)

\begin{lemma}\label{lem:mul-diag-exists}
  Given any associahedron diagonals $\AsDiag_1$ and $\AsDiag_2$ there
  is a multiplihedron diagonal $\MulDiag$ compatible with $\AsDiag_1$
  and $\AsDiag_2$.
\end{lemma}
\begin{proof}
  This follows from an inductive argument similar to the proof of
  Lemma~\ref{lem:ExistAssociahedronDiagonal}.
\end{proof}

A particular multiplihedron diagonal was constructed by
Saneblidze-Umble~\cite{SU04:Diagonals}.

\subsection{Module-map diagonals}\label{sec:mod-map-diag}

\begin{figure}
  \centering
  \includegraphics{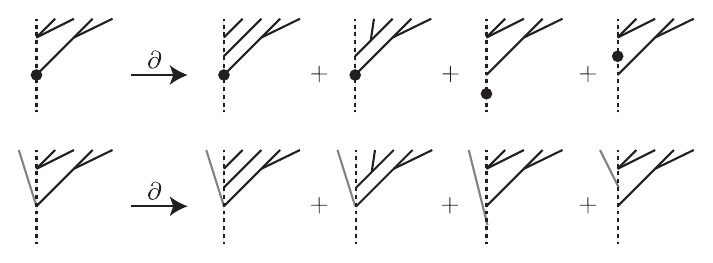}
  \caption[Module transformation trees]{\textbf{Module transformation trees.} Top: a module
    transformation tree with five inputs and its differential, where
    the distinguished vertex in each tree is indicated by the
    dot. Bottom: the same trees, but viewed as associahedron trees
    with $6$ inputs.}
  \label{fig:mod-trans-tree}
\end{figure}

To tensor maps of $\Ainf$-modules we will use
module-map diagonals. By a \emph{module transformation
  tree} we mean a tree $T$ together with a distinguished internal vertex $v$ on
the leftmost strand of $T$. We require that no vertex of $T$ except
perhaps $v$ have valence $2$. (Alternately, we can think of a module
transformation tree as a tree $T$ together with a coloring of each
internal vertex on the left-most strand of $T$ by red, blue or purple, so that
all of the red vertices come above the purple vertex and the purple
vertex comes above the blue vertices.) Define the differential of a
module transformation tree $(T,v)$ to be the sum of all ways of
inserting an edge in $T$ so as to get a new module transformation
tree. If the edge is inserted at a vertex other than $v$ then $v$
specifies the distinguished vertex in the new tree; if one inserts the
new edge at $v$ then there are sometimes two choices of distinguished vertex in
the new tree, and if so we take the sum of both choices. See
Figure~\ref{fig:mod-trans-tree}. The \emph{dimension} of a module
transformation tree is the number of inputs minus the number of
internal vertices (including $v$).  Let $\ModTransTrees{n}$ be the free
$\Ring$-module generated by the module transformation trees with $n$
inputs. The differential just defined makes
$\ModTransTrees{n}=\ModTransTrees{n,*}$ into a chain complex.

There is a correspondence between module transformation trees with $n$
inputs and associahedron trees with $n+1$ inputs, as follows. Given a
module transformation tree $(T,v)$, define a new tree $T'$ by adding a
new input leaf to the left of all inputs in $T$ and connecting the new
input to $v$. It is easy to see that this gives a bijection, and
moreover that it respects the differentials. In particular, the set of
module transformation trees with $n$ inputs corresponds to the cells
in $K_{n+1}$, and the differential on module transformation trees
corresponds to the cellular differential on $K_{n+1}$.
Thus:
\begin{corollary}\label{cor:ModTransTrees-contrac}
  For each $n$, the chain complex $\ModTransTrees{n,*}$ has homology
  $R$ in dimension $0$ and trivial homology in all other dimensions.
\end{corollary}

There are face inclusions
\begin{align*}
  \phi^n_{i,j}&\co \cellC{*}(K_{j-i+1})\otimes \ModTransTrees{n+i-j}
                \to \ModTransTrees{n}\\
  \phi^n_{1,j}&\co \ModTransTrees{j}\otimes \cellC{*}(K_{n+1-j})
                \to \ModTransTrees{n}\\
  \phi^n_{1,j}&\co \cellC{*}(K_j)\otimes \ModTransTrees{n+1-j}
                \to \ModTransTrees{n}
\end{align*}
which, at the level of trees, are given by
\[
\phi^n_{i,j}(S,T)=T\circ_iS.
\]
(If we identify $\ModTransTrees{n}\cong \cellC{*}(K_{1+n})$ then these
maps $\phi$ undergo an indexing shift: $\phi^n_{i,j}$ as defined here
corresponds to $\phi^{n+1}_{i+1,j+1}$ for the associahedron.)
\begin{definition}\label{def:mod-map-diag}
  Fix an associahedron diagonal $\AsDiag$ and module diagonals
  $\MDiag^1$ and $\MDiag^2$ compatible with $\AsDiag$. A
  \emph{module-map diagonal} compatible with $\MDiag^1$ and $\MDiag^2$
  is a collection of (degree-preserving) chain maps $\ModMulDiag_n\co
  \ModTransTrees{n,*}\to\ModTransTrees{n,*}\otimes\ModTransTrees{n,*}$
  satisfying the following conditions:
  \begin{itemize}
  \item \textbf{Compatibility}:
    \[
    \ModMulDiag_n\circ \phi_{i,j}^n=
    \begin{cases}
      (\phi_{1,j}^n\otimes \phi_{1,j}^n)\circ (\ModMulDiag_j\otimes\MDiag^2_{n-j+1} + \MDiag^1_j\otimes\ModMulDiag_{n-j+1}) & i=1\\
      (\phi_{i,j}^n\otimes\phi_{i,j}^n)\circ (\AsDiag_{j-i+1}\otimes\ModMulDiag_{n+i-j}) & i\ge 2.
    \end{cases}
  \]
  \item \textbf{Non-degeneracy}: $\ModMulDiag_1\co
    \ModTransTrees{1,*}\to\ModTransTrees{1,*}\otimes\ModTransTrees{1,*}$ is the
    standard isomorphism $\Ring\cong \Ring\otimes\Ring$.
  \end{itemize}

  Equivalently, we can view $\ModMulDiag_n$ as a formal linear
  combination of pairs of module transformation trees
  \[
    \TrModMulDiag_n\in\bigoplus_{i+j=n-1}\ModTransTrees{n,i}\otimes\ModTransTrees{n,j}
  \]
  in dimension $n-1$ (i.e., a \emph{module-map tree diagonal}). These
  trees must satisfy:
  \begin{itemize}
  \item \textbf{Compatibility}:
    \[
    \bdy(\TrModMulDiag_n)=\sum_{k+\ell=n+1}\Bigl(\TrModMulDiag_k\circ_1\TrMDiag_\ell^1+\TrMDiag_k^2\circ_1\TrModMulDiag_\ell + \sum_{i=2}^{k}\TrModMulDiag_k\circ_i\TrDiag_\ell\Bigr).
    \]
  \item \textbf{Non-degeneracy}: $\TrModMulDiag_1$ is the (unique)
    pair of module transformation trees with one input.
  \end{itemize}
\end{definition}

\begin{figure}
  \centering
  %Font is 12 point
  \includegraphics{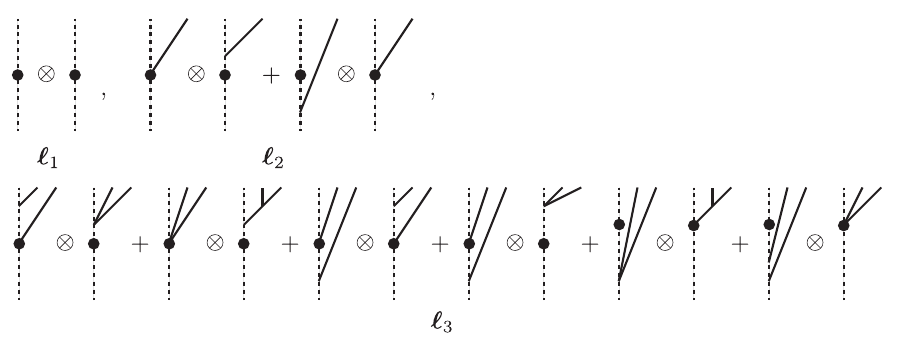}
  \caption[The first few terms in a module-map tree diagonal]{\textbf{The first few terms in a module-map tree diagonal.}
    This module-map tree diagonal is compatible with $\TrMDiag$ and
    $\TrMDiag$, where $\TrMDiag$ is the module tree diagonal of
    Figure~\ref{fig:prim-to-diag}.}
  \label{fig:m-tree-diag}
\end{figure}

The first few terms in a module-map tree diagonal are shown in Figure~\ref{fig:m-tree-diag}.

\begin{lemma}\label{lem:mod-map-diag-exists}
  Given any associahedron diagonal $\AsDiag$ and module diagonals
  $\MDiag^1$ and $\MDiag^2$ there is a module-map
  diagonal $\ModMulDiag$ compatible with $\MDiag^1$ and
  $\MDiag^2$.
\end{lemma}
\begin{proof}
  Again, this follows from an inductive argument similar to the proof
  of Lemma~\ref{lem:ExistAssociahedronDiagonal}.
\end{proof}

We turn next to uniqueness of module-map diagonals. Fix an
associahedron tree diagonal $\TrDiag$ and compatible module tree
diagonals $\TrMDiag^1$ and $\TrMDiag^2$. Given a formal linear
combination of pairs of module transformation trees
\[
  \eta_n\in \bigoplus_{i+j=k+n-1}\ModTransTrees{n,i}\otimes\ModTransTrees{n,j},
\]
for $n\geq 1$ and some given $k$,
we define $d(\eta)$, where $\eta=\{\eta_n\}$, to be the sequence whose $n\th$ term is
\begin{equation}\label{eq:deformed-mm-d}
  d(\eta)_n=\bdy(\eta_n)+\sum_{p+m=n+1}\Bigl(\eta_p\circ_1\TrMDiag_m^1+\TrMDiag_p^2\circ_1\eta_m + \sum_{i=2}^{p} \eta_p\circ_i\TrDiag_m\Bigr).
\end{equation}
\begin{lemma}\label{lem:deform-mod-map-cx}
  The operation $d$ makes the space of sequences $\eta_n$ into a chain
  complex, i.e., satisfies $d^2=0$. (Here, the grading on the chain
  complex is given by the parameter~$k$.)
\end{lemma}
\begin{proof}
  The verification is a straightforward computation, using 
  the compatibility conditions for the algebra diagonal $\gammas$ 
  (Equation~\eqref{eq:DiagonalCellCompatibilityTree}),
  the module diagonals $\MDiag^1$ and $\MDiag^2$
  (Equation~\eqref{eq:MDiagCompatTree}), and the fact that 
  $\sum_{i>1}(\cdot \circ_i \gamma_q)$
  is a derivation over $\circ_1$ in the sense that, for any~$\eta$,
  \[
    \sum_{i=2}^{p+m-1} (\TrMDiag^2_p\circ_1 \eta_m)\circ_i \TrDiag_q
    = \sum_{i=2}^{p} (\TrMDiag^2_p\circ_i \TrDiag_q)\circ_1 \eta_m
    + \sum_{i=2}^{m} \TrMDiag^2_p\circ_1 (\eta_m\circ_i \TrDiag_q).\qedhere
  \]
\end{proof}
The Compatibility condition in Definition~\ref{def:mod-map-diag} is
the statement that $\TrModMulDiag$ is a cycle in this chain
complex (with grading $0$). 

\begin{definition}\label{def:homotopy-module-map-diag}
  Given module-map diagonals $\TrModMulDiag$ and
  $\TrModMulDiag'$, a \emph{homotopy} from $\TrModMulDiag$ to
  $\TrModMulDiag'$ is a chain $\eta$ so that
  $d(\eta)=\TrModMulDiag-\TrModMulDiag'$.
\end{definition}

\begin{lemma}
  The chain complex from Lemma~\ref{lem:deform-mod-map-cx} has
  homology $\Ring$ in dimension $0$ and trivial homology in all other
  dimensions.
\end{lemma}
\begin{proof}
  The number of inputs~$n$ gives a descending filtration on this chain
  complex. The associated graded object is
  \[
    \prod_{n\geq 1}\ModTransTrees{n,*}\otimes\ModTransTrees{n,*},
  \]
  with the usual differential $\bdy$. By Corollary~\ref{cor:ModTransTrees-contrac}, the homology of the associated graded---i.e., the $E^1$-page of the associated spectral sequence---is 
  \[
    \prod_{n\geq 1}\Ring\otimes\Ring,
  \]
  The $i\th$ summand is generated by any pair of trees with $i$ inputs so that the
  purple vertices are $2$-valent and all other internal vertices are
  $3$-valent; this lies in grading $k=-n+1$. The $d_1$-differential in
  the spectral sequence comes
  from splicing in $\corolla{2}\otimes\corolla{2}$ at any of the $i$
  inputs or feeding the output into $\corolla{2}\otimes\corolla{2}$.
  This has $i+1$ terms, so the $E^1$-page is
  \[
    0\longrightarrow \Ring\stackrel{2}{\longrightarrow}\Ring\stackrel{3}{\longrightarrow}\Ring\stackrel{4}{\longrightarrow}\cdots.
  \]
  The resulting homology is $\Ring$, in dimension $0$ (corresponding to the pair of $1$-input trees).
\end{proof}

\begin{proposition}\label{prop:mod-map-diag-homotopic}
  All module-map diagonals are homotopic.
\end{proposition}
\begin{proof}
  Fix module-map diagonals $\TrModMulDiag$ and $\TrModMulDiag'$; we
  will construct a homotopy $\eta=\{\eta_n\}$ inductively in $n$.
  From the non-degeneracy condition,
  $\TrModMulDiag_1=\TrModMulDiag'_1$. So, we can define
  $\eta_1=0$. 
  Suppose inductively that we have defined all $\eta_i$ 
  satisfying $d(\eta)_i = \TrModMulDiag_i-\TrModMulDiag'_i$ for all $i<n$, for some $n>1$.
  Then, using Lemma~\ref{lem:deform-mod-map-cx},
  \[
    \bdy\left(\TrModMulDiag_n'-\TrModMulDiag_n+
    \sum_{p+m=n+1}\Bigl(\eta_p\circ_1\TrMDiag_m^1+\TrMDiag_p^2\circ_1\eta_m + \sum_{i=2}^{p} \eta_p\circ_i\TrDiag_m\Bigr)\right)=0.
  \]
  Further, the terms in
  \[
    \TrModMulDiag_n'-\TrModMulDiag_n+
  \sum_{p+m=n+1}\Bigl(\eta_p\circ_1\TrMDiag_m^1+\TrMDiag_p^2\circ_1\eta_m + \sum_{i=2}^{p} \eta_p\circ_i\TrDiag_m\Bigr)
  \]
  are in dimension $n-1>0$ in
  $\ModTransTrees{n,*}\otimes \ModTransTrees{n,*}$. So the result
  follows from Corollary~\ref{cor:ModTransTrees-contrac}.
\end{proof}

\subsection{Module-map primitives and partial module-map diagonals}\label{sec:mod-map-prim}
To define primitives of module-map diagonals, we use two degenerate
trees. One is the identity tree $\IdTree$ with $1$ input, $1$ output,
and $0$ internal vertices; this eventually corresponds to the identity
map from the algebra to itself.
The other is the \emph{stump} $\stump$, which is a $0$-input,
$1$-output tree with $0$ internal vertices; this eventually
corresponds to the unit element of the algebra. Composing $\IdTree$
(in any way) is the identity map. Composing a corolla with the stump
is $0$ except for the following cases:
\begin{align*}
  \IdTree\circ_1\stump &= \stump & \corolla{2}\circ_1\stump &= \IdTree & \corolla{2}\circ_2\stump &=\IdTree.
\end{align*}
These composition maps are chain maps; see Lemma~\ref{lem:extend-circ}.

\begin{definition}\label{def:mod-map-prim}
  Fix an associahedron diagonal $\AsDiag$ and two module diagonal
  primitives $\TrPMDiag^1$ and $\TrPMDiag^2$ compatible with
  $\AsDiag$. A \emph{module-map primitive} compatible with
  $\TrPMDiag^1$ and $\TrPMDiag^2$ is a linear combination
  \begin{equation}\label{eq:ModMap-prim-1}
    \TrPMorDiag_n=\sum_{(S,T)}n_{S,T}(S,T)\in
    \bigoplus_{i+j=n-1}\ModTransTrees{n,i}\otimes\cellC{j}(K_{n-1})
  \end{equation}
  of pairs of trees for each $n \ge 1$.  Here, the tree $S$ is viewed
  as a module transformation tree with $n$ inputs, and the tree $T$
  is an ordinary associahedron tree with $n-1$ inputs.
  In the case
  $n=1$, $\cellC{j}(K_{0})$ is interpreted as one-dimensional,
  generated by $\stump$, while for $n=2$, $\cellC{j}(K_{1})$ is
  interpreted as one-dimensional, generated by $\IdTree$.
  These elements $\TrPMorDiag_n$ are required to satisfy the following
  conditions:
  \begin{itemize}
  \item \textbf{Compatibility}:
    \begin{equation}
      \label{eq:Mor-prim-compat}
      \partial \TrPMorDiag = 
      \LRjoin((\TrPMDiag^1)^{\otimes\bullet}
      \otimes {\TrPMorDiag}\otimes
      (\TrPMDiag^2)^{\otimes\bullet})+\TrPMorDiag\circ'\TrDiag.
    \end{equation}
    Here, $\LRjoin$ from Definition~\ref{def:shift-join} is extended
    in the obvious way to module
    transformation trees, i.e.,
    \[
      \LRjoin((S_1,T_1),\dots,(S_k,T_k))=\LeftJoin(S_1,\dots,S_k)\otimes\RootJoin(T_1,\dots,T_k)
    \]
    (which is the tensor product of a module transformation
    tree and an associahedron tree).
    We have also invoked the usual convention that $\TrPMorDiag=\sum_n\TrPMorDiag_n\in\prod_n\ModTransTrees{n,*}\otimes\cellC{*}(K_{n-1})$.
    \signissue
  \item \textbf{Non-degeneracy}: $\TrPMorDiag_1=\pcorolla{1}\otimes\stump$.
  \end{itemize}
\end{definition}

\begin{figure}
  \centering
  %Font is 12 point.
  \includegraphics{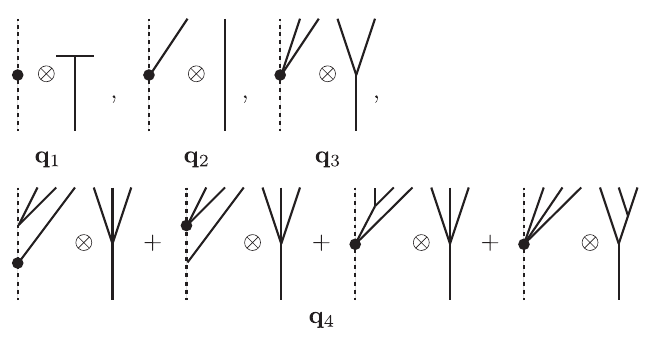}
  \caption[The first few terms in a module-map diagonal primitive]{\textbf{The first few terms in a module-map primitive.}
    This module-map primitive is compatible with $\TrPMDiag$ and
    $\TrPMDiag$, where $\TrPMDiag$ is the module diagonal primitive of
    Figure~\ref{fig:mprim-terms}.}
  \label{fig:m-map-prim}
\end{figure}

Note that the operation $\LRjoin$ used in
Definition~\ref{def:mod-map-prim} satisfies the formulas for $\LRjoin$
in Lemma~\ref{lem:join-diff}.

\begin{lemma}\label{lem:m-map-prim-exist}
  Given any associahedron diagonal $\AsDiag$ and module diagonal
  primitives $\TrPMDiag^1$ and $\TrPMDiag^2$, there is a module-map primitive
  $\TrPMorDiag$ compatible with $\TrPMDiag^1$ and $\TrPMDiag^2$.
\end{lemma}

\begin{proof}
  This follows as in the proof of Proposition~\ref{prop:prim-exist}.
\end{proof}

\begin{definition}\label{def:homotopy-module-map-primitive}
Given an associahedron tree diagonal, compatible module diagonal primitives
$\TrPMDiag^1$ and $\TrPMDiag^2$, and compatible module-map primitives $\TrPMorDiag^1$ and $\TrPMorDiag^2$,
a {\em homotopy of module-map primitives} is a linear combination of trees
\[ \xi_n \in \bigoplus_{i+j=n-1} \ModTransTrees{n,i}\otimes \cellC{j}(K_{n}) \]
so that 
\[
  \partial \xi_n + \LRjoin((\TrPMDiag^1)^{\otimes\bullet}\otimes \xi\otimes (\TrPMDiag^2)^{\otimes\bullet})+\xi\circ'\TrDiag
  = \TrPMorDiag^1-\TrPMorDiag^2.
\]
Two module-map primitives are called {\em homotopic} if there exists a homotopy between them.
\end{definition}

\begin{proposition}\label{prop:mod-map-prim-homotopic}
  All module-map primitives compatible with $\TrPMDiag^1$ and $\TrPMDiag^2$ are homotopic.
\end{proposition}

\begin{proof}
  Construct the homotopy via the same inductive procedure used to
  prove Proposition~\ref{prop:mod-map-diag-homotopic}, using
  Corollary~\ref{cor:ModTransTrees-contrac} and contractibility of the
  associahedron (Proposition~\ref{prop:AssociahedronIsContractible})
  to guarantee that the relevant homology groups
  vanish.
\end{proof}

In a special case, we can explicitly describe a module-map primitive in terms of a module diagonal primitive:
\begin{lemma}\label{lem:mor-prim-example}
  Let $\TrPMDiag$ be a module diagonal primitive. Let
  $\TrPMorDiag_1=\corolla{1}\otimes\stump$ and for $n\geq 2$, let
  $\TrPMorDiag_n$ be the sum over $(S,T)$ in $\TrPMDiag_n$ of all
  pairs $(S',T)$ where $S'$ is obtained from $S$ by making one
  internal vertex
  on the leftmost strand of $S$ distinguished (purple). Then $\TrPMorDiag$ is a
  module-map primitive compatible with $\TrPMDiag$ and $\TrPMDiag$.
\end{lemma}
\begin{proof}
  The differential $\bdy\TrPMorDiag$ has four kinds of terms, based on
  the distinguished vertex in the left tree:
  \begin{enumerate}[label=(\arabic*)]
  \item\label{item:mor-prim-eg-term-1} Terms appearing in $\bdy
    \TrPMDiag$ but with one internal vertex on the leftmost strand (of the left tree) distinguished.
  \item\label{item:mor-prim-eg-term-2} Terms where there is a 2-valent distinguished vertex at the top or the bottom of the leftmost strand (of the left tree).
  \item\label{item:mor-prim-eg-term-3} Terms where there is a 2-valent distinguished vertex on the leftmost strand (of the left tree), but not at the top or bottom.
  \end{enumerate}
  The third kind of terms cancel in pairs. By the compatibility
  condition for $\TrPMDiag$ (i.e., the right hand side of 
  Equation~\eqref{eq:M-prim-compat}), the
  first and second types of terms are of the form on the right side of
  compatibility
  condition for $\TrPMorDiag$,
  Equation~\eqref{eq:Mor-prim-compat}. The other possible terms on the right of
  Equation~\eqref{eq:Mor-prim-compat} are of the form
  \[
    \LRjoin((\TrPMDiag)^{\otimes k} \otimes {\TrPMorDiag_1}\otimes
    (\TrPMDiag)^{\otimes \ell})
  \]
  where $k$ and $\ell$ are each $\geq 1$. Since the right tree of $\TrPMorDiag_1$ is a
  stump, however, these terms vanish.
\end{proof}

The relationship between module-map primitives and module-map diagonals is more complicated than the relationship between module diagonal primitives and module diagonals. To state it, we need an auxiliary notion:
\begin{definition}\label{def:part-mod-map-diag}
  Let $\Filt_{\geq 3}\ModTransTrees{n,*}\subset \ModTransTrees{n,*}$ be the subspace spanned by
  module transformation trees in which the purple vertex has valence
  $\geq 3$ (i.e., is not 2-valent). This is not a
  subcomplex of $\ModTransTrees{n,*}$. Let
  \[
    \gModTransTrees{n,*}=\ModTransTrees{n,*}/\bigl(\Filt_{\geq 3}\ModTransTrees{n,*}+\bdy (\Filt_{\geq 3}\ModTransTrees{n,*+1})\bigr).
  \]

  A \emph{partial module-map diagonal} consists of elements
  \[
    \PartTrModMulDiag_n\in\bigoplus_{i+j=n-1}\ModTransTrees{n,i}\otimes\gModTransTrees{n,j}
  \]
  satisfying:
  \begin{itemize}
  \item \textbf{Compatibility}:
    \begin{equation}
      \bdy(\PartTrModMulDiag_n)-\sum_{k+\ell=n+1}\Bigl(\PartTrModMulDiag_k\circ_1\TrMDiag_\ell^1+\TrMDiag_k^2\circ_1\PartTrModMulDiag_\ell
      + \sum_{i=2}^{k}\PartTrModMulDiag_k\circ_i\TrDiag_\ell\Bigr)
      =0\in\ModTransTrees{n,*}\otimes\gModTransTrees{n,*}.
    \end{equation}
  \item \textbf{Non-degeneracy}: $\PartTrModMulDiag_1$ is the (unique)
    pair of module transformation trees with one input.
  \end{itemize}
  
  Partial module-map diagonals $\PartTrModMulDiag^1$ and $\PartTrModMulDiag^2$ are \emph{homotopic} if there is a collection of elements $\zeta_n\in \bigoplus_{i+j=n}\ModTransTrees{n,i}\otimes\gModTransTrees{n,j}$ satisfying
  \begin{equation}\label{eq:part-mod-map-htpy}
    \bdy(\zeta_n)-\sum_{k+\ell=n+1}\Bigl(\zeta_k\circ_1\TrMDiag_\ell^1+\TrMDiag_k^2\circ_1\zeta_\ell
    + \sum_{i=2}^{k}\zeta_k\circ_i\TrDiag_\ell\Bigr)+\PartTrModMulDiag^1_n-\PartTrModMulDiag^2_n=0\in\ModTransTrees{n,*}\otimes\gModTransTrees{n,*}.
  \end{equation}
  
\end{definition}
Explicitly, elements of $\gModTransTrees{n,*}$ are module
transformation trees where the purple vertex is $2$-valent, and if
two trees differ only in the location of the purple vertex they are
equivalent in $\gModTransTrees{n,*}$.

Given a module-map diagonal~$\TrModMulDiag$, the
image~$\PartTrModMulDiag$ of~$\TrModMulDiag$ under the quotient map
\[
  \ModTransTrees{}\otimes\ModTransTrees{}\to \ModTransTrees{}\otimes
  \gModTransTrees{}
\]
is a partial module-map diagonal.
In this case, we say that $\TrModMulDiag$ is a module-map
diagonal \emph{extending}~$\PartTrModMulDiag$.

\begin{remark}
  As we will see in Section~\ref{sec:box}, a partial module-map
  diagonal is the data needed to define $\Id\otimes g$, the tensor
  product of the identity map of $M$ and an arbitrary $\Ainf$-module
  map $g\co N\to N'$.
\end{remark}

\begin{lemma}\label{lem:mod-trans-subquot}
  The complex $\gModTransTrees{n,*}$ is contractible, i.e., has
  homology $\Ring$ in dimension $0$ and homology~$0$ in all other dimensions.
\end{lemma}
\begin{proof}
  Forgetting the purple vertex gives an isomorphism between
  $\gModTransTrees{n,*}$ and $\cellC{*}(K_n)$.
\end{proof}

\begin{lemma}\label{lem:part-mod-maps-htpic}
  All partial module-map diagonals compatible with $\TrMDiag^1$ and $\TrMDiag^2$
  are homotopic.
\end{lemma}
\begin{proof}
  We build the maps $\zeta_n$ satisfying
  Equation~(\ref{eq:part-mod-map-htpy}) inductively. Existence for
  $n=1$ is immediate from the non-degeneracy condition for partial
  module-map diagonals. For the inductive step, we claim that the element
  \begin{equation}\label{eq:zeta-elt}
    \sum_{k+\ell=n+1}\Bigl(\zeta_k\circ_1\TrMDiag_\ell^1+\TrMDiag_k^2\circ_1\zeta_\ell
    + \sum_{i=2}^{k}\zeta_k\circ_i\TrDiag_\ell\Bigr)+\PartTrModMulDiag^1_n-\PartTrModMulDiag^2_n
  \end{equation}
  in
  $\bigoplus_{i+j=n-1}\ModTransTrees{n,i}\otimes\gModTransTrees{n,j}$
  is a cycle. Indeed,
  \begin{align*}
    \bdy\biggl[\sum_{k+\ell=n+1}\Bigl(\zeta_k
    &\circ_1\TrMDiag_\ell^1+\TrMDiag_k^2\circ_1\zeta_\ell
    +
      \sum_{i=2}^{k}\zeta_k\circ_i\TrDiag_\ell\Bigr)+\PartTrModMulDiag^1_n-\PartTrModMulDiag^2_n\biggr]\\
    &=
      \sum_{k+\ell=n+1}\Bigl((\alpha^1_k+\bdy\beta^1_k)\circ_1\TrMDiag_\ell^1+\TrMDiag_k^2\circ_1(\alpha^1_\ell+\bdy\beta^1_\ell)
      +
      \sum_{i=2}^{k}(\alpha^1_k+\bdy\beta^1_k)\circ_i\TrDiag_\ell\Bigr)+\alpha^2_n+\bdy\beta^2_n
  \end{align*}
  for some
  $\alpha^1,\alpha^2,\beta^1,\beta^2\in\ModTransTrees{}\otimes\Filt_{\geq
    3}\ModTransTrees{}$. The terms involving $\alpha^1$ and
  $\alpha^2$ are in $\Filt_{\geq 3}\ModTransTrees{n,*}$, and
  $\bdy\beta^2_n\in\bdy\Filt_{\geq 3}\ModTransTrees{n,*+1}$. Further,
  \[
    \left[\sum_{k+\ell=n+1}\Bigl((\bdy\beta^1_k)\circ_1\TrMDiag_\ell^1+\TrMDiag_k^2\circ_1(\bdy\beta^1_\ell)
      +
      \sum_{i=2}^{k}(\bdy\beta^1_k)\circ_i\TrDiag_\ell\Bigr)\right]
    -
    \bdy\left[\sum_{k+\ell=n+1}\Bigl(\beta^1_k\circ_1\TrMDiag_\ell^1+\TrMDiag_k^2\circ_1\beta^1_\ell
      +
      \sum_{i=2}^{k}\beta^1_k\circ_i\TrDiag_\ell\Bigr)\right]
  \]
  lies in $\Filt_{\geq 3}\ModTransTrees{n,*}$, so
  \[
    \left[\sum_{k+\ell=n+1}\Bigl((\bdy\beta^1_k)\circ_1\TrMDiag_\ell^1+\TrMDiag_k^2\circ_1(\bdy\beta^1_\ell)
      +
      \sum_{i=2}^{k}(\bdy\beta^1_k)\circ_i\TrDiag_\ell\Bigr)\right]
    \in \Filt_{\geq 3}\ModTransTrees{n,*}+\bdy\Filt_{\geq 3}\ModTransTrees{n,*+1}.
  \]
  Thus, the
  element~\eqref{eq:zeta-elt} is a cycle in
  $\ModTransTrees{n,*}\otimes\gModTransTrees{n,*}$.

  Since the element~\eqref{eq:zeta-elt} is a cycle, from
  Corollary~\ref{cor:ModTransTrees-contrac} and
  Lemma~\ref{lem:mod-trans-subquot}, it is therefore a boundary unless
  $i=0$ and $j=0$, which only occurs for the base case $n=1$. Define
  $\zeta_n$ to be an element with boundary~\eqref{eq:zeta-elt} and
  continue the induction.
\end{proof}

Next we use an analogue of the operation $\LRjoin'$ from
Definition~\ref{def:assoc-mod-diag} to build a partial
module-map diagonal from a module-map primitive. Specifically, given a
sequence $(S_1,T_1),\dots,(S_k,T_k)$ of pairs of trees, with
$(S_i,T_i)\in \Trees_{n_i+1}\otimes\Trees_{n_i}$ for $i\neq j$ and
$(S_j,T_j)\in \ModTransTrees{n_j+1}\otimes \Trees_{n_j}$ define 
\begin{equation}\label{eq:RLjoinP-mor}
    \LRjoin'((S_1,T_1),\dots,(S_k,T_k))=
    \LeftJoin(S_1,\dots,S_k)\otimes
  \RootJoin(\IdTree,T_1,\dots,T_k),
\end{equation}
i.e., the obvious analogue of the definition of $\LRjoin'$ used for
module trees. Then $\LRjoin'((S_1,T_1),\dots,(S_k,T_k))$ is the tensor
product of a module transformation tree and an ordinary tree.

Note that this extension of $\LRjoin'$ still satisfies
Lemma~\ref{lem:mod-join-diff}.

We can turn $\LRjoin'((S_1,T_1),\dots,(S_k,T_k))$ into a tensor
product of two module transformation trees by adding a 2-valent purple
vertex at the bottom of the right tree. When the inputs are module
diagonal primitives and a module-map primitive, this gives a partial
module diagonal:

\begin{lemma}\label{lem:mm-prim-partial-mm-diag}
  Given an associahedron tree diagonal, compatible module diagonal
  primitives $\TrPMDiag^1$ and $\TrPMDiag^2$, and a compatible
  module-map primitive $\TrPMorDiag$,
  \begin{equation}\label{eq:MapMapPrimExponential}
    \PartTrModMulDiag=
    [\IdTree\otimes \pcorolla{1}]\circ \LRjoin'((\TrPMDiag^1)^{\otimes\bullet}\otimes
    \TrPMorDiag\otimes (\TrPMDiag^2)^{\otimes\bullet})
  \end{equation}
  is a partial module-map diagonal compatible with
  $\TrMDiag^1=\LRjoin'((\TrPMDiag^1)^{\otimes\bullet})$ and
  $\TrMDiag^2=\LRjoin'((\TrPMDiag^2)^{\otimes\bullet})$.
\end{lemma}
\begin{proof}
  We have
  \begin{align*}
    \bdy\Bigl([\IdTree\otimes
    &\pcorolla{1}]\circ \LRjoin'((\TrPMDiag^1)^{\otimes\bullet}\otimes
      \TrPMorDiag\otimes (\TrPMDiag^2)^{\otimes\bullet})\Bigr)\\
    &=[\IdTree\otimes \pcorolla{1}]\circ\Bigl(
      \LRjoin'((\TrPMDiag^1)^{\otimes\bullet}\otimes(\bdy \TrPMDiag^1)
      \otimes (\TrPMDiag^1)^{\otimes\bullet}\otimes \TrPMorDiag\otimes
      (\TrPMDiag^2)^{\otimes\bullet})
      +
      \LRjoin'((\TrPMDiag^1)^{\otimes\bullet}\otimes (\bdy\TrPMorDiag)\otimes
      (\TrPMDiag^2)^{\otimes\bullet})\\
    &\quad +
      \LRjoin'((\TrPMDiag^1)^{\otimes\bullet}\otimes \TrPMorDiag\otimes
      (\TrPMDiag^2)^{\otimes\bullet}\otimes(\bdy \TrPMDiag^2)
      \otimes (\TrPMDiag^2)^{\otimes\bullet})
      +
      \LRjoin'((\TrPMDiag^2)^{\otimes\bullet})\circ_1\LRjoin'((\TrPMDiag^1)^{\otimes\bullet}\otimes\TrPMorDiag\otimes
      (\TrPMDiag^2)^{\otimes\bullet})\\
      &\quad +
      \LRjoin'((\TrPMDiag^1)^{\otimes\bullet}\otimes\TrPMorDiag\otimes
      (\TrPMDiag^2)^{\otimes\bullet})\circ_1\LRjoin'((\TrPMDiag^1)^{\otimes\bullet})
      +
      \LRjoin'((\TrPMDiag^1)^{\otimes\bullet}\otimes \LRjoin((\TrPMDiag^1)^{\otimes\bullet}\otimes
      \TrPMorDiag\otimes (\TrPMDiag^2)^{\otimes\bullet})\otimes
      (\TrPMDiag^2)^{\otimes\bullet})\\
      &\quad +
      \LRjoin'((\TrPMDiag^1)^{\otimes\bullet}\otimes\LRjoin((\TrPMDiag^1)^{\otimes\bullet})\otimes
      (\TrPMDiag^1)^{\otimes\bullet}\otimes
      \TrPMorDiag\otimes (\TrPMDiag^2)^{\otimes\bullet})\\
    &\quad +
      \LRjoin'((\TrPMDiag^1)^{\otimes\bullet}\otimes
      \TrPMorDiag\otimes (\TrPMDiag^2)^{\otimes\bullet}\otimes\LRjoin((\TrPMDiag^2)^{\otimes\bullet})\otimes
      (\TrPMDiag^2)^{\otimes\bullet})
      \Bigr)\\
    &=[\IdTree\otimes \pcorolla{1}]\circ\Bigl(
      \LRjoin'((\TrPMDiag^1)^{\otimes\bullet}\otimes(\TrPMDiag^1\circ\TrDiag)
      \otimes (\TrPMDiag^1)^{\otimes\bullet}\otimes \TrPMorDiag\otimes
      (\TrPMDiag^2)^{\otimes\bullet})
      +
      \LRjoin'((\TrPMDiag^1)^{\otimes\bullet}\otimes (\TrPMorDiag\circ
      \TrDiag)\otimes
      (\TrPMDiag^2)^{\otimes\bullet})\\
    &\quad +
      \LRjoin'((\TrPMDiag^1)^{\otimes\bullet}\otimes \TrPMorDiag\otimes
      (\TrPMDiag^2)^{\otimes\bullet}\otimes(\TrPMDiag^2\circ \TrDiag)
      \otimes (\TrPMDiag^2)^{\otimes\bullet})
      +
      \LRjoin'((\TrPMDiag^2)^{\otimes\bullet})\circ_1\LRjoin'((\TrPMDiag^1)^{\otimes\bullet}\otimes\TrPMorDiag\otimes
      (\TrPMDiag^2)^{\otimes\bullet})\\
      &\quad +
      \LRjoin'((\TrPMDiag^1)^{\otimes\bullet}\otimes\TrPMorDiag\otimes
        (\TrPMDiag^2)^{\otimes\bullet})\circ_1\LRjoin'((\TrPMDiag^1)^{\otimes\bullet})\Bigr)
    \\
    &=
      \left(\sum_{i=2}^\infty \PartTrModMulDiag\circ_i\TrDiag\right)+
      \PartTrModMulDiag\circ_1\TrMDiag^1
      +
      [\IdTree\otimes \pcorolla{1}]\circ\TrMDiag^2\circ_1\LRjoin'((\TrPMDiag^1)^{\otimes\bullet}\otimes\TrPMorDiag\otimes
      (\TrPMDiag^2)^{\otimes\bullet})
  \end{align*}
  where the first equality uses Lemma~\ref{lem:mod-join-diff}, the
  second uses the structure equations, and the third uses the
  definitions. So, we want to show that
  \[
    [\IdTree\otimes \pcorolla{1}]\circ_1\TrMDiag^2\circ_1\LRjoin'((\TrPMDiag^1)^{\otimes\bullet}\otimes\TrPMorDiag\otimes
    (\TrPMDiag^2)^{\otimes\bullet})
    -
    \TrMDiag^2\circ_1 [\IdTree\otimes \pcorolla{1}]\circ_1\LRjoin'((\TrPMDiag^1)^{\otimes\bullet}\otimes\TrPMorDiag\otimes
    (\TrPMDiag^2)^{\otimes\bullet})
  \]
  is of the form $\bdy\mathbf{s}+\mathbf{t}$ for some elements
  $\mathbf{s},\mathbf{t}\in \ModTransTrees{n}\otimes
  \Filt_{\geq 3}\ModTransTrees{n}$.
  This is clear: the element $\mathbf{s}$ is the sum of all ways of
  declaring one of the internal vertices in $\TrMDiag^2$ in
  $\TrMDiag^2\circ_1
  \LRjoin'((\TrPMDiag^1)^{\otimes\bullet}\otimes\TrPMorDiag\otimes
  (\TrPMDiag^2)^{\otimes\bullet})$ to be purple.
\end{proof}

\begin{remark}
  It is not clear whether the partial module-map diagonals constructed
  by Equation~\eqref{eq:MapMapPrimExponential} can be extended to a full
  module-map diagonal. However, by
  Lemma~\ref{lem:part-mod-maps-htpic}, that partial module-map
  diagonal is homotopic to an extendible partial module-map diagonal.
\end{remark}

\subsection{\texorpdfstring{$\mathrm{DADD}$}{DADD} diagonals}\label{sec:DADD-diags}
\begin{definition}\label{def:DADD-diag}
  Fix associahedron diagonals $\TrDiag^1$ and $\TrDiag^2$.  A \emph{DADD diagonal}
  compatible with $\TrDiag^1$ and $\TrDiag^2$ is a collection of chains
  \[
    \TrDADD_n=\sum_{(S,T)}n_{S,T}(S,T)\in \cellC{*}(J_n) \otimes \Trees_{n}
  \]
  in dimension $n-1$, satisfying the following conditions:
  \begin{itemize}
  \item \textbf{Compatibility}:
    \begin{equation}\label{eq:DADD-compat}
      \bdy(\TrDADD_n)=\sum_{k+\ell=n+1} \TrDADD_k\circ\TrDiag^1_\ell
      +\sum_{j_1+\cdots+j_k=n}
      \TrDiag^2_k\circ(\TrDADD_{j_1},\cdots,\TrDADD_{j_k}).
    \end{equation}
    (The second term on the right hand side is the result of feeding
    the outputs of $\TrDADD_{j_1},\dots,\TrDADD_{j_k}$ into
    $\TrDiag_k^2$, analogous to Equation~\eqref{eq:multi-comp-pair}.)
  \item \textbf{Non-degeneracy}:
    $\TrDADD_1=\pcorolla{1}\otimes \IdTree$, the tensor product of a
    1-input purple corolla and the identity tree.
  \end{itemize}
\end{definition}

\begin{figure}
  \centering
  % Font is 12 point
  \includegraphics{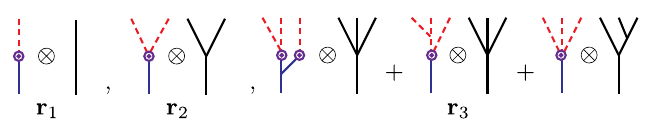}
  \caption[The first few terms in a DADD diagonal]{\textbf{The first few terms in a DADD diagonal.} This
    diagonal is compatible with the associahedron diagonal from
    Figure~\ref{fig:diag-cells}. The cases $\TrDADD_1$, $\TrDADD_2$,
    and $\TrDADD_3$ are shown.}
  \label{fig:DADD-eg}
\end{figure}

Schematically, the compatibility condition is:
\begin{equation}\label{eq:DADD-diag-compat}
  \bdy \left(\mathcenter{
  \tikzsetnextfilename{def-DADD-diag-1}
    \begin{tikzpicture}
      \node at (0,0) (tl) {};
      \node at (0,-1) (DADDL) {$\TrDADD$};
      \node at (0,-2) (bl) {};
      \node at (1,0) (tr) {};
      \node at (1,-1) (DADDR) {$\TrDADD$};
      \node at (1,-2) (br) {};
      \draw[taa, red] (tl) to (DADDL);
      \draw[alga, blue] (DADDL) to (bl);
      \draw[taa] (tr) to (DADDR);
      \draw[alga] (DADDR) to (br);
    \end{tikzpicture}
  }\right)
=
\left(
  \mathcenter{
  \tikzsetnextfilename{def-DADD-diag-2}
    \begin{tikzpicture}
      \node at (-.5,0) (tll) {};
      \node at (0.5,0) (tlr) {};
      \node at (0,0) (tl) {};
      \node at (0,-1) (DeltaL) {$\TrDiag^1$};
      \node at (0,-2) (DADDL) {$\TrDADD$};
      \node at (0,-3) (bl) {};
      \node at (1.5,0) (trl) {};
      \node at (2,0) (tr) {};
      \node at (2.5,0) (trr) {};
      \node at (2,-1) (DeltaR) {$\TrDiag^1$};
      \node at (2,-2) (DADDR) {$\TrDADD$};
      \node at (2,-3) (br) {};
      \draw[taa, red] (tll) to (DADDL);
      \draw[taa, red] (tlr) to (DADDL);
      \draw[taa] (trl) to (DADDR);
      \draw[taa] (trr) to (DADDR);
      \draw[taa, red] (tl) to (DeltaL);
      \draw[taa] (tr) to (DeltaR);
      \draw[alga, red] (DeltaL) to (DADDL);
      \draw[alga] (DeltaR) to (DADDR);
      \draw[alga, blue] (DADDL) to (bl);
      \draw[alga] (DADDR) to (br);
    \end{tikzpicture}}\right)
+
\left(
  \mathcenter{
  \tikzsetnextfilename{def-DADD-diag-3}
    \begin{tikzpicture}
      \node at (-1,0) (tl1) {};
      \node at (0,0) (tl2) {};
      \node at (-1,-1) (DADDL1) {$\TrDADD$};
      \node at (-.5,-1) (dotsl) {$\cdots$};
      \node at (0,-1) (DADDL2) {$\TrDADD$};
      \node at (-.5,-2) (DeltaL) {$\TrDiag^2$};
      \node at (-.5,-3) (bl) {};
      \node at (1,0) (tr1) {};
      \node at (2,0) (tr2) {};      
      \node at (1,-1) (DADDR1) {$\TrDADD$};
      \node at (1.5,-1) (dotsr) {$\cdots$};
      \node at (2,-1) (DADDR2) {$\TrDADD$};
      \node at (1.5,-2) (DeltaR) {$\TrDiag^2$};
      \node at (1.5,-3) (br) {};
      \draw[taa, red] (tl1) to (DADDL1);
      \draw[taa, red] (tl2) to (DADDL2);
      \draw[alga, blue] (DADDL1) to (DeltaL);
      \draw[alga, blue] (DADDL2) to (DeltaL);
      \draw[alga, blue] (DeltaL) to (bl);
      \draw[taa] (tr1) to (DADDR1);
      \draw[taa] (tr2) to (DADDR2);
      \draw[alga] (DADDR1) to (DeltaR);
      \draw[alga] (DADDR2) to (DeltaR);
      \draw[alga] (DeltaR) to (br);
    \end{tikzpicture}
  }
\right).
\end{equation}
The first few terms in a particular DADD diagonal are shown in Figure~\ref{fig:DADD-eg}.

\begin{lemma}\label{lem:DADD-exists}
  Given any pair of associahedron diagonals, there is a compatible DADD
  diagonal.
\end{lemma}
\begin{proof}
  As usual, since the associahedron and multiplihedron are
  contractible, it suffices to verify:
  \begin{enumerate}
  \item The right hand side of the compatibility equation is a cycle.
  \item Solutions to the compatibility equation exist when the right
    side is in dimension $0$ (with respect to the usual grading on the
    tensor product $\cellC{*}(J_n) \otimes \Trees_{n})$.
  \end{enumerate}
  The first statement is clear, particularly from the schematic
  version of the compatibility condition. The second corresponds to
  the cases $\TrDADD_1$ and $\TrDADD_2$, which are shown in
  Figure~\ref{fig:DADD-eg}.
\end{proof}

The reader might notice a similarity between Figure~\ref{fig:DADD-eg}
and the module diagonal primitive in
Figure~\ref{fig:mprim-terms}. This similarity is not coincidental,
though a module diagonal primitive involves only one associahedron
diagonal. Nonetheless, we have:
\begin{proposition}\label{prop:DADD-to-prim}
  Fix a DADD diagonal $\TrDADD_n$ compatible with $\TrDiag^1$ and
  $\TrDiag^2$. Given a pair of trees
  $(S,T)\in\TrDADD_n$, call $S$ \emph{unital} if every blue vertex of
  $S$ has valence $3$. Given $(S,T)\in\TrDADD_n$ with $S$ unital, let
  $S'\in\Trees_{n+1}$ be the result of deleting all of the blue edges
  and vertices from $S$ and then running a new strand on the left
  through all of the purple vertices. Then
  \[
    \TrPMDiag_n\coloneqq \sum_{\substack{(S,T)\in\TrDADD_n\\ S\text{ unital}}}(S',T)
  \]
  is a module diagonal primitive
  compatible with~$\TrDiag^1$.
\end{proposition}

The proof of Proposition~\ref{prop:DADD-to-prim} will use the following
lemma:
\begin{lemma}\label{lem:diag-some-term}
  Let $\TrDiag$ be an associahedron diagonal. Then for each $n\geq 2$ there
  are an odd number of pairs of trees $(S,T)$ in $\TrDiag$ where $S$
  is a binary tree. For each of these pairs of trees, $T$ is the
  corolla $\corolla{n}$.
\end{lemma}
\begin{proof}
  We will give an indirect argument, relying on material from
  Section~\ref{sec:algebra}.

  Note that if $S$ is binary then $T$ must be a corolla for grading
  reasons.  Observe that the diagonal from
  Example~\ref{eg:explicit-diag} has exactly one pair of trees $(S,T)$
  with $S$ binary.

  Now, consider the $\Ainf$-algebra
  $\Alg=\Ring\langle 1,a_1,\dots,a_n,b\rangle$ with the following
  nontrivial products:
  \begin{align*}
   \forall x \in \Alg: \mu_2(1,x)=\mu_2(x,1)&=x\\
   \mu_n(a_1,\dots,a_n)&=b.
  \end{align*}
  (All other products are~$0$.)
  Consider also the ordinary algebra $\Blg$ generated by elements
  $x_1,\dots,x_n$ with the relations $x_ix_j=0$ if $j\neq i+1$. So, a
  basis for $\Blg$ is in bijection with the set of intervals
  $[k,\ell]\subset\{1,\dots,n\}$, via $[k,\ell]\leftrightarrow
  x_{[k,\ell]}=x_kx_{k+1}\dots x_\ell$.

  By Theorem~\ref{thm:Algebras}, up to isomorphism, the
  $\Ainf$-algebra $\Alg\ADtp\Blg$ is independent of the choice of
  diagonal $\AsDiag$. By nondegeneracy of~$\AsDiag$, the algebra
  $\Alg\ADtp\Blg$ has operations
  \[
    \mu_2\bigl((1\otimes x_{[k,\ell]}),(1\otimes
    x_{[\ell,m]})\bigr)=1\otimes x_{[k,m]}.
  \]
  If $\TrDiag$ has an odd number of pairs $(S,T)$ as in the
  statement of the lemma then $\Alg\ADtp\Blg$ has an $\Ainf$ operation
  \[
    \mu_n((a_1\otimes x_1),\dots,(a_n\otimes x_n))=b\otimes(x_1\cdots x_n).
  \]
  Otherwise, $\Alg\ADtp\Blg$ has no higher $\Ainf$ operations at
  all. In either case, these are the only nontrivial $\Ainf$
  operations. It is clear that there is no quasi-isomorphism between
  the two cases. This implies the result.
\end{proof}

\begin{proof}[Proof of Proposition~\ref{prop:DADD-to-prim}]
  Given a pair of trees
  $(S,T)\in \cellC{*}(J_n) \otimes \Trees_{n}$, let
  $\Phi(S,T)=0$ if $S$ is not unital and let $\Phi(S,T)$ be the result
  of deleting the blue vertices of $S$ and running a new strand
  through the purple vertices of the result if $S$ is unital. Our goal
  is to show that the pairs of trees $\Phi(\TrDADD_n)$ form a module
  diagonal primitive.

  From the nondegeneracy condition for a DADD diagonal, we have that
  $\TrPMDiag_2=\corolla{2}\otimes\IdTree$. It remains to verify that
  $\TrPMDiag$ satisfies the compatibility condition for a module
  diagonal primitive, i.e., that
  \[
    \partial \TrPMDiag =
    \LRjoin(\TrPMDiag^{\otimes\bullet})+\TrPMDiag\circ'\TrDiag^1.
  \]
  Of course, we will deduce this from the compatibility condition
  for~$\TrDADD$.

  First, we claim that $\Phi(\bdy(\TrDADD_n))=\bdy(\Phi(\TrDADD_n))$. 
  Consider first the pairs of trees $(S,T)$ in $\bdy(\TrDADD_n)$ where
  $S$ is unital. Such trees can arise either as the boundary of a pair
  $(S',T')$ where $S'$ is unital or from a pair $(S',T')$ where $S'$
  has a single blue vertex with valence $4$. The second kinds of terms
  occur in pairs, and the result of deleting the blue vertices from
  the two terms in the pair is the same. Also, taking the
  differential of a tree $S$ at a purple vertex and then applying
  $\Phi$ is the same as applying $\Phi$ and then taking the
  differential at the corresponding vertex. These two observations
  prove that $\Phi(\bdy(\TrDADD_n))=\bdy(\Phi(\TrDADD_n))$.

  Next, observe that applying $\Phi$ to the first term on the right of
  the DADD compatibility equation (Equation~\eqref{eq:DADD-compat})
  gives the second term on the right of the primitive compatibility
  equation (Equation~\eqref{eq:M-prim-compat}).

  Finally, we claim that applying $\Phi$ to the second term on the
  right of Equation~\eqref{eq:DADD-compat} 
  gives the first term on the
  right of Equation~\eqref{eq:M-prim-compat}.
    Note that for
  $(S,T)\in \TrDiag^2_k$,
  \[
    \Phi\bigl((S,T)\circ(\TrDADD_{j_1},\cdots,\TrDADD_{j_k})\bigr)=0
  \]
  unless $S$ is a binary tree. In this case, by
  Lemma~\ref{lem:diag-some-term}, there are an
  odd number of such pairs of trees and $T$ is a corolla. Given such a
  pair of trees,
  \[
    \Phi\bigl((S,T)\circ(\TrDADD_{j_1},\cdots,\TrDADD_{j_k})\bigr)=\LRjoin(\TrDADD_{j_1},\cdots,\TrDADD_{j_k}).
  \]
  This finishes the proof.
\end{proof}

%%% Local Variables: 
%%% mode: latex
%%% TeX-master: "AbstractDiagonal"
%%% TeX-command-extra-options: "-shell-escape"
%%% End: 

\section{Applications of diagonals}\label{sec:algebra}
\begin{convention}\label{conv:Ring}
  Throughout this section, $\Ring$ is a commutative
  $\FF_2$-algebra (possibly graded) and $\Ground$ (or $\Ground_1$,
  $\Ground_2$, etc.)\ is a free or, more generally, flat commutative
  $\Ring$-algebra (again, possibly graded). All modules are assumed to
  be free or, more generally, projective over $\Ring$.
  Further, the action of $\Ring$ on algebras and bimodules is central: if
  $r\in R$ and $a$ is an element of a $\Ring$-algebra or bimodule, we assume
  that $ra=ar$.

  In this section, undecorated tensor products are over whichever ring $\Ground$ (or $\Ground_i$)
  is relevant, and not over $R$ or $\FF_2$.
\end{convention}

\subsection{The \dg category of chain complexes}
When talking about $\Ainf$-modules or type $D$ structures, we find it
convenient to use the language of \dg categories (or more generally
$\Ainf$-categories), so we spell out some of this terminology
now. (See, for instance, Keller~\cite{Keller06:DGCategories} for a more
extensive introduction.) Briefly, a \dg category is a category
enriched over chain complexes. That is, given objects $M$ and $N$ in a
\dg category $\Cat$, the space of morphisms $\Mor(M,N)$ between $M$
and $N$ is a chain complex, and composition is a chain map
$\Mor(N,P)\otimes \Mor(M,N)\to \Mor(M,P)$.

The first example of a \dg category is the category of chain complexes
itself, where the $i\th$ graded part $\Mor_i(M,N)$ of $\Mor(M,N)$
consists of tuples of maps of abelian groups
\[
\{f_k\co M_k\to N_{k+i}\}_{k\in\ZZ},
\]
and the differential is
\[
d(f)=\bdy\circ f +f\circ\bdy.
\]
In particular, the (degree 0) cycles in $\Mor(M,N)$ are the (degree 0)
chain maps, and the boundaries are the nullhomotopic chain maps.

Below, we will work with the \dg category of $\Ainf$-modules. The
definition of this category is standard; see
\cite[Section~2.2.2]{LOT2} for
a brief review. We will also work with the $\Ainf$-category of type
$D$ structures;
see~\cite[Section~2.2.3]{LOT2} for the
definition of this $\Ainf$-category and, for instance,
Seidel~\cite{Seidel02:FukayaDef} for a
discussion of $\Ainf$-categories in general.

We will typically use $\Mor$ to denote a chain complex of maps of
$\Ground$-modules or bimodules; so, for instance, if $\Alg$ is an
$\Ainf$-algebra over $\Ground$ then $\Mor(\Alg^{\kotimes{\Ground} n},\Alg)$ is
the chain complex of maps of (graded) $\Ground$-bimodules
$\Alg^{\kotimes{\Ground} n}\to\Alg$, with differential induced from $\mu_1$
on~$\Alg$. The complex $\Mor(\Alg^{\kotimes{\Ground} n},\Alg)$ is an $\Ring$-module.

\begin{convention}\label{conv:grading-shift}
  If $C=\{C_d\}$ is a graded module then $C\grs{n}$ is the graded module
  \[
    C\grs{n}_d=C_{d-n}.
  \]
  Thus the \emph{degree $n$} part of $\Mor(M,N)$, i.e., the morphisms
  sending $M_d$ to $N_{d+n}$, is the same as degree 0 morphisms
  $M\to N\grs{-n}$ or $M\grs{n}\to N$.
\end{convention}

\subsection{Tensor products of \texorpdfstring{$\Ainf$}{A-infinity}-algebras}
Let $A$ be a $\Ground$-bimodule, equipped with a differential $\mu_1\co
A\to A\grs{1}$.  A collection of functions $\{\mu_d\co A^{\kotimes{\Ground} d}
\to A\grs{2-d}\}_{d=2}^{\infty}$ can be extended to a collection of
functions $\{\cellC{*}(K_{n})\to\Mor(A^{\kotimes{\Ground} n},A)\}$ by the
following procedure. Given a tree $T$ with $n$ inputs (thought of as a
cell of $K_n$), replace each vertex with $d>0$ inputs by the
operation $\mu_{d}$, and compose these operations as specified by the
edges in $T$. Denote the induced function $\mu(T)\co A^{\kotimes{\Ground} n}\to A\grs{-\dim(T)}$. 
Extend linearly to linear combinations of trees to
get a map $\mu\co \cellC{d}(K_{n})\to\Mor_d(A^{\kotimes{\Ground} n},A)$.

As above, $\mu_1$ induces a differential on $\Mor(A^{\kotimes{\Ground} n},A)$.

We say the operations $\{\mu_d\}$ satisfy the \emph{$\Ainf$-algebra relations} if for any $n\geq 1$ and any $a_1,\dots,a_n\in A$,
\begin{equation}\label{eq:A-inf-alg-rel}
\sum_{k=0}^{n-1}\sum_{i=1}^{n-k}\mu_{n-k}(a_1\otimes\dots\otimes a_{i-1}\otimes\mu_{k+1}(a_i\otimes\dots\otimes a_{i+k})\otimes a_{i+k+1}\otimes\dots\otimes a_n)=0.
\end{equation}
The following reformulation of this condition is well known.
\begin{lemma}\label{lem:A-inf-alg-is}
  The collection $\{\mu_n\}$ satisfies the $\Ainf$-algebra relations~\eqref{eq:A-inf-alg-rel} if
  and only if for each $n$ the map $\mu\co \cellC{*}(K_n)\to \Mor(A^{\kotimes{\Ground} n},A)$ is a
  chain map.
\end{lemma}
\begin{proof}
  The statement that $\bdy \mu(\corolla{n})=\mu(\bdy\corolla{n})$ is
  exactly the $\Ainf$-algebra relation~\eqref{eq:A-inf-alg-rel}. To
  see that this implies that $\mu$ is a chain map in general, observe that if $\bdy[\mu(S)]=\mu(\bdy(S))$ and $\bdy[\mu(T)]=\mu(\bdy(T))$ then 
  \begin{align*}
  \bdy[\mu(T\circ_i S)] &= \bdy[\mu(T)\circ (\Id\otimes \mu(S)\otimes \Id)]\\
  &=(\bdy[\mu(T)])\circ(\Id\otimes \mu(S)\otimes \Id) + \mu(T)\circ (\Id\otimes \bdy[\mu(S)]\otimes\Id)\\
  &=(\mu(\bdy(T)))\circ(\Id\otimes \mu(S)\otimes \Id) + \mu(T)\circ (\Id\otimes \mu(\bdy(S))\otimes\Id)\\
  &=\mu(\bdy(T)\circ_i S)+\mu(T\circ_i \bdy(S)) = \mu(\bdy(T\circ_i S)).
  \end{align*}
  Since any tree can be written as a composition of corollas, the result follows.
\end{proof}
A collection of chain maps
$\mu\co \cellC{*}(K_n)\to \Mor(A^{\kotimes{\Ground} n},A)$ specifies
an $\Ainf$-algebra homomorphism if and only if it is compatible with
composition (stacking of trees). In other words, an $\Ainf$-algebra is
a module over the associahedron operad.

\begin{definition}\label{def:Ainf-alg}
  An \emph{$\Ainf$-algebra} over $\Ground$ is a \dg $\Ground$-bimodule $(A,\mu_1)$ together with
  operations $\mu_n\co A^{\otimes n}\to A$, $n\geq 2$, satisfying the
  $\Ainf$-algebra relations.
\end{definition}

When studying type \DD\ structures below, we will find it convenient to equip
our $\Ainf$-algebras with a unit. An $\Ainf$-algebra $\Alg=(A,\{\mu_n\})$ is
\emph{strictly unital} if there is an element $\unit\in A$ so that
$\mu_2(\unit,a)=\mu_2(a,\unit)=a$ for all $a\in A$ and $\mu_n(a_1,\dots,a_n)=0$
if $n\neq 2$ and some $a_i=\unit$. See also
Definition~\ref{def:strict-unital}. There is further discussion of units in
Section~\ref{subsec:UnitsHomPert}.

The following is a rephrasing of a definition of
Saneblidze-Umble~\cite[Definition 30]{SU04:Diagonals}.
\begin{definition}\label{def:Alg-tp}
  Fix an associahedron diagonal $\AsDiag$. Given $\Ainf$-algebras
  $\Alg_1=(A_1,\{\mu^1_n\})$ and $\Alg_2=(A_2,\{\mu^2_n\})$ over
  $\Ground_1$ and $\Ground_2$, let
  $\Alg_1\ADtp\Alg_2$ be the \dg bimodule
  $A_1\rotimes{\Ring}A_2$ over $\Ground=\Ground_1\rotimes{\Ring}\Ground_2$,
  which we endow with an $\Ainf$-algebra
  structure over $\Ground$ via the maps
  \begin{align*}
  \cellC{*}(K_n)&\xrightarrow{\AsDiag_n} \cellC{*}(K_n)\rotimes{\Ring} \cellC{*}(K_n)\\
  &\xrightarrow{\mu^1_n\otimes \mu^2_n} \Mor(A_1^{\kotimes{\Ground_1} n},A_1)\rotimes{\Ring}\Mor(A_2^{\kotimes{\Ground_2} n},A_2)\hookrightarrow \Mor((A_1\rotimes{\Ring} A_2)^{\kotimes{\Ground} n},A_1\rotimes{\Ring} A_2)
  \end{align*}
  (compare Lemma~\ref{lem:A-inf-alg-is}). The compatibility of the
  diagonal~$\AsDiag$ under stacking implies that this map
  $\cellC{*}(K_n)\to \Mor((A_1\rotimes{\Ring} A_2)^{\kotimes{\Ground} n},A_1\rotimes{\Ring} A_2)$
  is induced by its values on corollas.

  Explicitly, if $\TrDiag=\sum n_i(S_i,T_i)$ is the associahedron tree diagonal 
  corresponding to $\AsDiag$ then the operation
  $(\mu^1 \ADtp \mu^2)_n$ on
  $\Alg_1\ADtp\Alg_2$ is defined by
  \begin{equation}\label{eq:Ainf-tens-trees}
  (\mu^1 \ADtp \mu^2)_n((a_1\otimes b_1)\otimes\dots\otimes(a_n\otimes b_n))\\
   =\sum_i n_i[\mu^1(S_i)(a_1\otimes\cdots\otimes a_n)]\otimes[\mu^2(T_i)(b_1\otimes\cdots\otimes b_n)].
   \end{equation}
 \end{definition}

\begin{lemma}\label{lem:Alg-tp}
  With notation as above, the maps $\mu^1 \ADtp \mu^2$ from an
  $\Ainf$-algebra.
\end{lemma}

\begin{proof}
  The map $\mu^1 \ADtp \mu^2$ is defined as a composition of chain
  maps, and so is a chain map. So, it follows from compatibility of
  $\AsDiag$ under stacking that $\mu^1 \ADtp \mu^2$ specifies an
  $\Ainf$-algebra.
\end{proof}

Later, we will shorthand the right side of
Formula~\eqref{eq:Ainf-tens-trees} as
\begin{equation}\label{eq:Ainf-tens-trees-graph}
  \mathcenter{
  \tikzsetnextfilename{AinfTensTrees}
  \begin{tikzpicture}[smallpic]
    \node at (0,0) (tl) {};
    \node at (0,-1) (gammaL) {$\TrDiag$};
    \node at (0,-2) (bl) {};
    \node at (1,0) (tr) {};
    \node at (1,-1) (gammaR) {$\TrDiag$};
    \node at (1,-2) (br) {};
    \draw[taa] (tl) to (gammaL);
    \draw[alga] (gammaL) to (bl);
    \draw[taa] (tr) to (gammaR);
    \draw[alga] (gammaR) to (br);
  \end{tikzpicture}}
\end{equation}

We turn next to maps of $\Ainf$-algebras. 

\begin{definition}
  \label{def:AlgebraHomomorphism}
Fix $\Ainf$-algebras $\Alg$
and $\Blg$ over $\Ground$. We say that a collection of degree 0 $\Ground$-bimodule maps $f=\{f_n\co A^{\kotimes{\Ground}
  n} \to B\grs{1-n}\}_{n=1}^{\infty}$ is an \emph{$\Ainf$-algebra
  homomorphism} if for each $n$ and each sequence of elements $a_1,\dots,a_n\in A$, the maps $f_k$ satisfy
\begin{multline}
  \label{eq:Ainf-homo}
  \sum_{k=1}^{n-1}\sum_{i=1}^{n-k} f_{n-k}\bigl(a_1\otimes\dots\otimes a_{i-1}\otimes \mu^{\Alg}_{k+1}(a_1\otimes\cdots\otimes a_{k+i})\otimes a_{k+i+1}\otimes\cdots\otimes a_n\bigr)\\
  +\!\!\!\!
  \sum_{n=m_1+\cdots+m_k}\!\!\!\!
  \mu^\Blg_k\bigl(f_{m_1}(a_1\otimes \cdots\otimes a_{m_1})\otimes\cdots\otimes f_{m_k}(a_{m_1+\cdots+m_{k-1}+1}\otimes\cdots \otimes a_n)\bigr)
  =0.
\end{multline}
Given $\Ainf$-algebra homomorphisms $f\co \Alg\to\Blg$ and
$g\co\Blg\to\Clg$, their \emph{composition} is a map $g\circ f$,
defined by
\[
  (g\circ f)_n(a_1,\dots,a_n)=\sum_{k=1}^n \sum_{i_1+\cdots+i_k=n}g_k\bigl(f_{i_1}(a_1,\dots,a_{i_1}),f_{i_2}(a_{i_1+1},\dots,a_{i_1+i_2}),\dots,f_{i_k}(a_{i_1+\dots+i_{k-1}},\dots,a_n)\bigr).
\]
The \emph{identity map} of $\Alg$ is defined by $(\Id_{\Alg})_1=\Id_A$
and $(\Id_{\Alg})_n=0$ for $n\neq 1$. A homomorphism
$f\co \Alg\to\Blg$ is an \emph{isomorphism} if there is a homomorphism
$g\co \Blg\to\Alg$ so that $g\circ f=\Id_{\Alg}$ and
$f\circ g=\Id_{\Blg}$. 

An $\Ainf$-algebra homomorphism $f$ is a \emph{quasi-isomorphism} if
the chain map $f_1$ is a quasi-isomorphism.
\end{definition}

Recall from Section~\ref{sec:multiplihedron} that $J_n$ denotes the
multiplihedron, a CW complex whose cells correspond to transformation
trees.  A collection of maps $\{f_n\co A^{\kotimes{\Ground} n} \to
B\grs{1-n}\}_{n=1}^{\infty}$ can be extended to a give a map
\[
  \phi\co \cellC{d}(J_n) \to \Mor_d(A^{\kotimes{\Ground} n},B)
\]
as follows. Given a transformation tree~$T$, replace
each red vertex of $T$ with $k$ inputs by the operation~$\mu_k^\Alg$,
each blue vertex of $T$ with $k$ inputs by the operation~$\mu_k^\Blg$
and each purple vertex of $T$ with $k$ inputs by the
function~$f_k$. Compose these operations as specified by the
edges of~$T$. We denote the induced function by
$f(T)\co A^{\kotimes{\Ground} n}\to B$. The analogue of
Lemma~\ref{lem:A-inf-alg-is} is:
\begin{lemma}\label{lem:A-inf-alg-map-is}
  The collection $\{f_n\}$ satisfies the $\Ainf$-algebra homomorphism
  relations~\eqref{eq:Ainf-homo} if and only if for each $n$, the map $f\co
  \cellC{*}(J_n)\to \Mor(A^{\kotimes{\Ground} n},B)$ is a chain map.
\end{lemma}
The proof is similar to the proof of Lemma~\ref{lem:A-inf-alg-is}, and
is left to the reader.

The following is the analogue of Definition~\ref{def:Alg-tp}:
\begin{definition}
  \label{def:Alg-map-tp}
  Fix associahedron diagonals $\AsDiag_1$ and $\AsDiag_2$ and a
  multiplihedron diagonal $\MulDiag$ compatible with $\AsDiag_1$ and
  $\AsDiag_2$. Let $\Alg_1=(A_1,\{\mu^{\Alg_1}_n\})$ and
  $\Blg_1=(B_1,\{\mu^{\Blg_1}_n\})$ be $\Ainf$-algebras over
  $\Ground_1$, $\Alg_2=(A_2,\{\mu^{\Alg_2}_n\})$ and
  $\Blg_2=(B_2,\{\mu^{\Blg_2}_n\})$ be $\Ainf$-algebras over
  $\Ground_2$, and $f^1=\{f^1_n\co A_1^{\kotimes{\Ground_1} n}\to B_1\}$ and
  $f^2=\{f^2_n\co A_2^{\kotimes{\Ground_2} n}\to B_2\}$ be $\Ainf$-algebra
  homomorphisms. Define a map
  \[
  (f^1\otimes_{\MulDiag} f^2)\co \cellC{*}(J_n)\to \Mor\bigl((A_1\rotimes{\Ring} A_2)^{\kotimes{\Ground} n}, (B_1\rotimes{\Ring} B_2)\bigr)
  \]
  to be the composition
  \begin{align*}
    \cellC{*}(J_n)&\stackrel{\MulDiag}{\longrightarrow}\cellC{*}(J_n)\rotimes{\Ring}\cellC{*}(J_n)\\
    &\xrightarrow{f^1\otimes_{\MulDiag} f^2}\Mor(A_1^{\kotimes{\Ground_1} n},B_1)\rotimes{\Ring}\Mor(A_2^{\kotimes{\Ground_2} n},B_2)
    \hookrightarrow \Mor(A_1^{\kotimes{\Ground_1} n}\otimes A_2^{\kotimes{\Ground_2} n},B_1\otimes B_2).
\end{align*}

  Explicitly, in terms of trees, if we write $\TrMulDiag=\sum_i n_i(S_i,T_i)$ (where $n_i\in\Ring$ and $S_i$ and $T_i$ are transformation trees) then 
  \[
  (f^1\otimes_{\MulDiag} f^2)_n((a_1\otimes a'_1)\otimes\dots\otimes(a_n\otimes a'_n))=\sum_i n_i[f^1_n(S_i)(a_1\otimes\dots\otimes a_n)]\otimes [f^2_n(T_i)(a'_1\otimes\dots\otimes a'_n)].
  \]
\end{definition}

\begin{lemma}\label{lem:Alg-map-tp}
  With notation as above, the maps $f_1 \otimes_\MulDiag f_2$ form an
  $\Ainf$-algebra homomorphism from
  $\Alg_1 \ADtp[\AsDiag_1] \Alg_2$ to
  $\Blg_1 \ADtp[\AsDiag_2] \Blg_2$.
\end{lemma}

\begin{proof}
  This follows from Lemma~\ref{lem:A-inf-alg-map-is}.
  The map $f_1 \otimes_\MulDiag f_2$ is defined as a composition of
  chain maps, and so is a chain map. Compatibility with stacking
  (composition of trees) is immediate from the definition and the fact
  that $\MulDiag$ is compatible with $\AsDiag_1$ and~$\AsDiag_2$.
\end{proof}

Before proving Theorem~\ref{thm:Algebras}, we recall an equivalent
definition of isomorphisms of $\Ainf$-algebras:
\begin{lemma}\label{lem:Alg-iso-is}\cite[Lemma 2.1.14]{LOT2}
  A homomorphism $f\co \Alg\to\Blg$ of $\Ainf$-algebras is an
  isomorphism (i.e., has an inverse $g\co \Blg\to\Alg$ so that
  $f\circ g=\Id_{\Blg}$ and $g\circ f=\Id_{\Alg}$) if and only if
  $f_1\co A\to B$ is an isomorphism of $\Ground$-bimodules.
\end{lemma}
(The proof is to use invertibility of $f_1$ to inductively construct
the inverse $g$.)

\begin{proof}[Proof of Theorem~\ref{thm:Algebras}]
  Existence of associahedron diagonals was verified in
  Lemma~\ref{lem:ExistAssociahedronDiagonal}. It remains to verify the
  rest of the theorem.

  Lemma~\ref{lem:Alg-tp} states that $\Alg_1\ADtp\Alg_2$ is an
  $\Ainf$-algebra.  By definition, the underlying $\Ring$-module of
  $\Alg_1\ADtp\Alg_2$ is $A_1\otimes A_2$
  (point~\ref{item:Alg-thm-vec-space}).

  For point~\ref{item:Alg-thm-dg}, note that the non-degeneracy
  condition implies that $\mu_2$ is the usual multiplication on
  $\Alg_1\otimes\Alg_2$. For $n>2$, if $(S,T)$ corresponds to an
  $(n-2)$-cell of $K_n\times K_n$ then one of $S$ or $T$ has a vertex
  with $3$ or more inputs. It follows that, if $\Alg_1$ and $\Alg_2$
  are \dg algebras, the higher multiplications on
  $\Alg_1\ADtp\Alg_2$ vanish.

  We will prove points~\ref{item:Alg-thm-qi}
  and~\ref{item:Alg-thm-change-diag} together. So, fix
  $\Ainf$-algebras $\Alg_1$, $\Alg_2$, $\Alg'_1$ and $\Alg'_2$ and
  associahedron diagonals $\AsDiag_1$ and $\AsDiag_2$.  Suppose that
  $f^1\co \Alg_1\to \Alg_1'$ and $f^2\co \Alg_2\to\Alg_2'$ are
  $\Ainf$-quasi-isomorphisms.  Fix a multiplihedron diagonal $\MulDiag$
  compatible with $\AsDiag_1$ and $\AsDiag_2$; such a $\MulDiag$
  exists by Lemma~\ref{lem:mul-diag-exists}.
  Lemma~\ref{lem:Alg-map-tp} then gives an $\Ainf$-algebra
  homomorphism
  \[
  f^1\otimes_{\MulDiag} f^2\co \Alg_1\ADtp[\AsDiag_1]\Alg_2\to \Alg'_1\ADtp[\AsDiag_2]\Alg'_2.
  \]
  To see that this is a quasi-isomorphism, note that the non-degeneracy condition on $\MulDiag$ implies that 
  \[
  (f^1\otimes_{\MulDiag} f^2)_1=(f^1_1\otimes f^2_1)\co A_1\rotimes{\Ring}
  A_2\to A'_1\rotimes{\Ring} A'_2,
  \]
  and since $f^1_1$ and $f^2_1$ induce isomorphisms on homology, and
  the $A_i$ and $A'_i$ are flat over $\Ring$,
  so
  does $f^1_1\otimes f^2_1$. Point~\ref{item:Alg-thm-qi} corresponds
  to the special case that $\AsDiag_1=\AsDiag_2$ and
  Point~\ref{item:Alg-thm-change-diag} corresponds to the special case
  that $\Alg_1=\Alg'_1$ and $\Alg_2=\Alg'_2$ and the $f_i$ are the
  identity maps. For Point~\ref{item:Alg-thm-change-diag}, it follows
  from the non-degeneracy condition for multiplihedron diagonals and
  Lemma~\ref{lem:Alg-iso-is} that $(f^1\otimes_{\MulDiag} f^2)$ is not
  just a quasi-isomorphism but an isomorphism.

  Finally, to prove point~\ref{item:Alg-thm-assoc} (associativity),
  one defines an \emph{associahedron double-diagonal} to be a
  collection of maps 
  \[
  \AsDiag\AsDiag_n\co \cellC{*}(K_n)\to \cellC{*}(K_n)\rotimes{\Ring} \cellC{*}(K_n)\rotimes{\Ring} \cellC{*}(K_n)
  \]
  satisfying
  \[
  \AsDiag\AsDiag_n\circ \phi^n_{i,j} = (\phi^n_{i,j}\otimes \phi^n_{i,j}\otimes \phi^n_{i,j})\circ (\AsDiag\AsDiag_{j-i+1}\otimes \AsDiag\AsDiag_{n+i-j}),
  \]
  and with $\AsDiag\AsDiag_2$ the standard isomorphism from
  $\cellC{*}(\pt)\cong
  \cellC{*}(\pt)\rotimes{\Ring}\cellC{*}(\pt)\rotimes{\Ring}\cellC{*}(\pt)$. (Recall
  that $\phi^n_{i,j}$ is the face inclusion of the associahedron; see
  Formula~\eqref{eq:face-maps}.)
  Then:
  \begin{enumerate}
  \item\label{item:double-diag-1} Both $(\AsDiag\otimes\Id)\circ\AsDiag$ and
    $(\Id\otimes\AsDiag)\circ\AsDiag$ are associahedron
    double-diagonals.
  \item An associahedron double-diagonal is exactly what one needs to
    define the triple tensor product
    $(\Alg_1\otimes\Alg_2\otimes\Alg_3)_{\AsDiag\AsDiag}$. For the two
    associahedron diagonals in point~(\ref{item:double-diag-1}), the
    triple tensor products are
    $(\Alg_1\ADtp\Alg_2)\ADtp\Alg_3$ and
    $\Alg_1\ADtp(\Alg_2\ADtp\Alg_3)$,
    respectively.
  \item Any two associahedron double-diagonals are related by a
    \emph{multiplihedron double-diagonal}. This is the analogue of
    Lemma~\ref{lem:mul-diag-exists}, and uses the obvious analogue of
    Definition~\ref{def:multiplihedron-diag}.
  \end{enumerate}
  Given this, the proof of associativity follows along the same lines
  as the proof of independence from the associahedron diagonals.
\end{proof}

\begin{remark}
  Associativity of the tensor product can hold only up to
  isomorphism, in general, according to a result of Markl and
  Shnider~\cite[Theorem 6.1]{MS06:AssociahedraProdAinf}.
\end{remark}

In Section~\ref{sec:typeD}, some boundedness properties of our algebra
$\Alg$ will be relevant. In previous papers~\cite{LOT1,LOT2}, we
called an $\Ainf$-algebra $(A,\{\mu_i\})$ \emph{operationally bounded}
if there exists an $N$ so that $\mu_i=0$ for all $i>N$. This property
seems not to be preserved by tensor products, but the following
stronger property is.
\begin{definition}\label{def:alg-bonsai}
  An $\Ainf$-algebra $\Alg=(A,\{\mu_i\})$ is \emph{bonsai} if there is
  an integer $N$ so that for any associahedron tree $T$ with
  $\dim(T)>N$, we have
  $\mu(T)=0$. An $\Ainf$-algebra homomorphism $f\co\Alg\to\Blg$ is \emph{bonsai} if there is
  an integer $N$ so that for any multiplihedron tree $T$ with
  $\dim(T)>N$, we have $f(T)=0$. We will sometimes call such an $N$  a
  \emph{bonsai constant} of~$\Alg$.
\end{definition}
Since $\corolla{n}$ has dimension $n-2$, being bonsai implies being
operationally bounded.

\begin{lemma}\label{lem:alg-bounded}
  For any associahedron diagonal $\AsDiag$, if $\Alg_1$ and $\Alg_2$
  are bonsai
  then so is $\Alg_1\ADtp\Alg_2$. Similarly, for any
  associahedron diagonals $\AsDiag_1$ and~$\AsDiag_2$ and compatible
  multiplihedron
  diagonal~$\Theta$, if
  $\Alg_i$, $\Blg_i$, and $f^i\co\Alg_i\to\Blg_i$ are bonsai for $i=1,2$ then $f^1\otimes_{\MulDiag}f^2$ is bonsai.
\end{lemma}
\begin{proof}
  This is immediate from the fact that $\AsDiag$ (respectively
  $\MulDiag$) preserves the grading (dimension).
\end{proof}
(We will relax this and other boundedness assumptions in
Section~\ref{sec:boundedness}.)

\subsection{Tensor products of \texorpdfstring{$\Ainf$}{A-infinity}-modules}
\label{sec:TensorAinftyModules}
We turn next to modules. Throughout this section we will focus on
right modules; analogous statements hold for left modules, using left
module diagonals.

Let $\Alg=(A,\{\mu_d\})$ be an $\Ainf$-algebra
over~$\Ground$, and let $(M,m_1\co M\to M\grs{1})$ be a right \dg
module over~$\Ground$. A
collection of functions $\{m_{d+1}\co M\kotimes{\Ground} A^{\kotimes{\Ground} d}\to
M\grs{1-d}\}_{d=1}^\infty$ can be extended to a collection of functions
$\cellC{d}(K_{n+1})\to \Mor_d(M\kotimes{\Ground} A^{\kotimes{\Ground} n},M)$ by the
following procedure. Given a tree $T$ with $n$ inputs (thought of as a
cell of $K_n$), we replace each vertex on the leftmost strand of $T$
with $d>1$ inputs by $m_d$ and each vertex not on the leftmost strand
of $T$ with $d>1$ inputs by $\mu_d$, and compose these operations as
specified by the edges in $T$. We denote the induced function $m(T)\co
M\kotimes{\Ground} A^{\kotimes{\Ground} n}\to M\grs{-\dim(T)}$. Extending linearly to linear
combinations of trees we get a map $m\co \cellC{*}(K_{n+1})\to
\Mor(M\kotimes{\Ground} A^{\kotimes{\Ground} n},M)$.

If $(M,m_1)$ is a differential $\Ground$-module then there is an
induced differential on $\Mor(M\kotimes{\Ground} A^{\kotimes{\Ground} n},M)$,
incorporating both the differential on $M$ and on $A$.

We say the operations $\{m_d\}$ satisfy the \emph{$\Ainf$-module
  relations} if for any $n\geq 0$, any $x\in M$ and any
$a_1,\dots,a_{n}\in A$,
\begin{multline}\label{eq:A-inf-mod-rel}
  0=\sum_{k=0}^{n} m_{n-k+1}(m_{k+1}(x\otimes a_1\otimes \dots\otimes a_k)\otimes a_{k+1}\otimes \dots\otimes a_n)\\
  +\sum_{k=0}^{n}\sum_{i=1}^{n-k}m_{n-k+1}(x\otimes \dots\otimes a_{i-1}\otimes \mu_{k+1}(a_i\otimes \dots\otimes a_{i+k})\otimes a_{i+k+1}\otimes \dots\otimes a_n).
\end{multline}

\begin{lemma}\label{lem:Ainf-mod-is}
  The collection $\{m_d\}$ satisfies the $\Ainf$-module relations if
  and only if the map $m\co \cellC{*}(K_{n+1})\to
  \Mor(M\kotimes{\Ground} A^{\kotimes{\Ground} n},M)$ is a chain map.
\end{lemma}
\begin{proof}
  The proof is similar to the proof of Lemma~\ref{lem:A-inf-alg-is}
  and is left to the reader.
\end{proof}

\begin{definition}\label{def:Mod-tp}
  Let $\Alg_i$, $i=1,2$, be an $\Ainf$-algebra over $\Ground_i$ and let $\cModule_i=(M_i,\{m_{n+1}^{i})$,
  $i=1,2$, be a right $\Ainf$-module over $\Alg_i$. Fix an associahedron
  diagonal $\AsDiag$ and a module diagonal $\MDiag$ compatible with $\AsDiag$. Then
  we can endow $M_1\rotimes{\Ring} M_2$ with the structure of a right
  $\Ainf$-module over $\Alg_1\ADtp \Alg_2$ via the maps
  \begin{align*}
    \cellC{*}(K_{n+1})&\xrightarrow{\MDiag_{n+1}} \cellC{*}(K_{n+1})\rotimes{\Ring} \cellC{*}(K_{n+1})\\
    &\xrightarrow{m^1_{n+1}\otimes m^2_{n+1}} \Mor(M_1\kotimes{\Ground_1} \Alg_1^{\kotimes{\Ground_1} n},M_1)\rotimes{\Ring}\Mor(M_2\kotimes{\Ground_2}\Alg_2^{\kotimes{\Ground_2} n},M_2)\\
    &\hookrightarrow \Mor((M_1\rotimes{\Ring} M_2)\kotimes{\Ground}(\Alg_1\rotimes{\Ring}\Alg_2)^{\kotimes{\Ground} n},M_1\rotimes{\Ring} M_2)
  \end{align*}
  (compare Definition~\ref{def:Alg-tp}).
  Let $\cModule_1\MDtp \cModule_2$ denote this $\Ainf$-module.

  Explicitly, if $\TrMDiag=\sum n_i(S_i,T_i)$ is the linear
  combination of trees corresponding to $\MDiag$ then the operation
  $m_{n+1}$ on $\cModule_1\MDtp \cModule_2$ is given by
  \begin{multline*}
  m_{n+1}((x\otimes y)\otimes (a_1\otimes b_1)\otimes\dots\otimes(a_n\otimes b_n))\\
    =\sum_i n_i[m^1(S_i)(x\otimes a_1\otimes\cdots\otimes a_n)]\otimes[m^2(T_i)(y\otimes b_1\otimes\cdots\otimes b_n)].
  \end{multline*}
\end{definition}

Similarly, we can use module-map diagonals to define tensor
products of $\Ainf$-module homomorphisms. We start by recalling some
definitions, and an analogue of Lemma~\ref{lem:A-inf-alg-map-is}. Fix
an $\Ainf$-algebra $\Alg$ and $\Ainf$-modules $\cModule$ and
$\cNodule$ over~$\Alg$. An \emph{$\Ainf$-morphism} from $\cModule$ to
$\cNodule$ is a collection of
maps
\[
f_{1+n}\co M\kotimes{\Ground} A^{\kotimes{\Ground} n}\to N\grs{-n}.
\]
Let $\Mor(\cModule,\cNodule)$ denote the set of $\Ainf$-morphisms from $\cModule$ to $\cNodule$. We define a differential on $\Mor(\cModule,\cNodule)$ by 
\begin{multline}\label{eq:d-on-mor}
d(f)(x\otimes a_1\otimes\cdots\otimes a_n)\\
\begin{aligned}
&=\sum_{j=0}^n f_{1+n-j}(m^M_{1+j}(x\otimes a_1\otimes\cdots\otimes a_j)\otimes a_{j+1}\otimes\cdots\otimes a_n)
\\
&\quad+
\sum_{j=0}^n m^N_{1+j}(f_{1+n-j}(x\otimes a_1\otimes\cdots\otimes a_{n-j})\otimes a_{n-j+1}\otimes\cdots\otimes a_n)\\
&\quad+\sum_{1\leq i\leq j\leq n}f_{1+n-j+i}(x\otimes a_1\otimes\cdots\otimes a_{i-1}\otimes \mu_{j-i+1}(a_i\otimes\cdots\otimes a_j)\otimes a_{j+1}\otimes\cdots\otimes a_n)
\end{aligned}
\end{multline}
(where $f=\{f_{1+n}\}$).

Recall that $\ModTransTrees{n}$ denotes the chain complex of module
transformation trees with $1$ module input and $n-1$ algebra
inputs. (The complex $\ModTransTrees{n}$ is isomorphic to
$\cellC{*}(K_{n+1})$.) We get an induced map
$F_f\co \ModTransTrees{n+1,d}\to \Mor_d(M\kotimes{\Ground}
A^{\kotimes{\Ground} n},M)$ by replacing the distinguished vertex in
each module transformation tree by $f$, replacing the other internal vertices
on the left-most strand by $m^{\cModule}$ or $m^{\cNodule}$ depending
on whether they come above or below the distinguished vertex,
replacing all other internal vertices by $\mu$, and composing these maps
according to the tree.
\begin{lemma}\label{lem:F-f}
  The association $f\mapsto F_f$ is a chain map. That is, 
  \begin{equation}\label{eq:mod-map-diff}
  F_{d(f)}=d(F_f)=F_f\circ\bdy_{\ModTransTrees{n+1}}^{\cell}+\bdy^{\Mor}\circ F_f.
  \end{equation}
\end{lemma}
\begin{proof}
  As in the proof of Lemma~\ref{lem:A-inf-alg-is}, it suffices to verify that 
  \begin{equation}\label{eq:mod-chain-map-pf}
  F_{d(f)}(\corolla{n+2})(x\otimes a_1\otimes\cdots\otimes a_n)=d(F_f)(\corolla{n+2})(x\otimes a_1\otimes\cdots\otimes a_n).
  \end{equation}
  The left side of Formula~\eqref{eq:mod-chain-map-pf} is given
  by Formula~(\ref{eq:d-on-mor}).
  After expanding the right side of
  Formula~\eqref{eq:mod-chain-map-pf} by
  Formula~\eqref{eq:mod-map-diff}, the second term contributes
  the terms in Formula~(\ref{eq:d-on-mor}) involving $m_1^M$, $m_1^N$
  or $\mu_1$, and the first term contributes the remaining terms
  in Formula~(\ref{eq:d-on-mor}).
\end{proof}

\begin{definition}\label{def:tensor-mod-maps}
  Fix an associahedron diagonal $\AsDiag$, module diagonals $\MDiag^1$
  and $\MDiag^2$ compatible with $\AsDiag$, and a module-map diagonal
  $\ModMulDiag$ compatible with $\MDiag^1$ and $\MDiag^2$. Suppose
  that $\cModule_1$ and $\cModule_2$ are $\Ainf$-modules over an $\Ainf$-algebra $\Alg$ and $\cNodule_1$ and $\cNodule_2$ are $\Ainf$-modules over an $\Ainf$-algebra $\Blg$, and
  $f\in\Mor(\cModule_1, \cModule_2)$ and $g\in\Mor(\cNodule_1, \cNodule_2)$ are morphisms in the
  \dg category of $\Ainf$-modules. Then there is a corresponding
  morphism $(f\MDtp[\ModMulDiag] g)\in \Mor(\cModule_1\MDtp[{\MDiag^1}]
  \cNodule_1,\cModule_2\MDtp[{\MDiag^2}] \cNodule_2)$ so that
  $F_{f\MDtp[\ModMulDiag] g}$ is the composition
  \begin{equation}\label{eq:define-mod-mor-tp}
  \begin{split} \ModTransTrees{n+1}&\stackrel{\ModMulDiag}{\longrightarrow}\ModTransTrees{n+1}\rotimes{\Ring}\ModTransTrees{n+1}\\
    &\xrightarrow{F_{f}\otimes F_{g}}\,\Mor(M_1\kotimes{\Ground_1} A^{\kotimes{\Ground_1} n},M_2)\rotimes{\Ring} \Mor(N_1\kotimes{\Ground_2} B^{\kotimes{\Ground_2} n},N_2)\\
    &\hookrightarrow \Mor(M_1\rotimes{\Ring} N_1\kotimes{\Ground}(A\rotimes{\Ring} B)^{\kotimes{\Ground} n},M_2\rotimes{\Ring} N_2).
  \end{split}
  \end{equation}

  Explicitly, if $\TrModMulDiag=\sum n_i(S_i,T_i)$ is the linear
  combination of module transformation trees corresponding to
  $\ModMulDiag$ then the operation $(f\MDtp[\ModMulDiag] g)_{1+n}$ is
  given by
  \begin{multline*}
  (f\MDtp[\ModMulDiag] g)_{1+n}\bigl((x\otimes y)\otimes(a_1\otimes b_1)\otimes\cdots\otimes (a_n\otimes b_n)\bigr)\\
  =\sum_i n_i\bigl[F_f(S_i)(x\otimes a_1\otimes\cdots\otimes a_n)\bigr]\otimes \bigl[F_g(T_i)(y\otimes b_1\otimes\cdots\otimes b_n)\bigr].
  \end{multline*}
\end{definition}

\begin{lemma}\label{lem:tens-mod-maps}
  Definition~\ref{def:tensor-mod-maps} defines a chain map
  \[
    \Mor(\cModule_1,\cModule_2)\rotimes{\Ring}\Mor(\cNodule_1,\cNodule_2)\to\Mor(\cModule_1\MDtp[\MDiag^1]\cNodule_1,\cModule_2\MDtp[\MDiag^2]\cNodule_2).
  \]
  Moreover, homotopic module-map diagonals give homotopic chain maps.
\end{lemma}
\begin{proof}
  This is immediate from the definition and Lemma~\ref{lem:F-f}.
\end{proof}

\begin{proposition}\label{prop:dg-bifunctor}
  Fix an associahedron diagonal~$\AsDiag$,
  module diagonals $\MDiag^1$, $\MDiag^2$, and $\MDiag^3$ compatible
  with $\AsDiag$, and for $1 \le i < j \le 3$ module-map diagonals
  $\ModMulDiag^{ij}$ compatible with $\MDiag^i$ and $\MDiag^j$. Then
  for any $\Ainf$-algebras $\Alg$ and $\Blg$, $\Ainf$-modules $\cModule_1$, $\cModule_2$,
  and $\cModule_3$ over $\Alg$, $\Ainf$-modules $\cNodule_1$,
  $\cNodule_2$, and $\cNodule_3$ over $\Blg$, and morphisms
  \[
    \cModule_1\stackrel{f_1}{\longrightarrow}\cModule_2\stackrel{f_2}{\longrightarrow}\cModule_3\qquad\qquad
    \cNodule_1\stackrel{g_1}{\longrightarrow}\cNodule_2\stackrel{g_2}{\longrightarrow}\cNodule_3,
  \]
  we have a chain homotopy
  \begin{equation}
    (f_2\circ f_1)\MDtp[\ModMulDiag^{13}](g_2\circ g_1)\sim(f_2\MDtp[\ModMulDiag^{23}] g_2)\circ (f_1\MDtp[\ModMulDiag^{12}] g_1).
      \label{eq:mod-map-tensor-comp}
  \end{equation}
\end{proposition}
\begin{proof}
  We introduce yet another
  type of trees to encode the two sides of Equation~\eqref{eq:mod-map-tensor-comp}. A \emph{module
    bi-transformation tree} is a tree~$T$ together with \emph{two}
  distinguished internal vertices $v_1$, $v_2$ on the left-most strand. As
  before, we require that all internal vertices except for the $v_i$ have
  valence~$3$ or more. We think of the $v_i$ as dividing the left-most
  strand into three sections, corresponding to three modules involved
  (e.g., $\cModule_1$, $\cModule_2$, and $\cModule_3$). The differential of a module
  bi-transformation tree~$T$ is the sum over all ways of inserting an
  edge in~$T$ so as to get a new module bi-transformation tree.

  The module bi-transformation trees complex with $n$ inputs can be
  viewed as a subcomplex of the associahedron trees complex with $n+2$
  inputs: given a bi-transformation tree, the corresponding
  associahedron tree is obtained by adding two inputs on the far left
  and routing them to the two distinguished vertices. The
  corresponding quotient complex is spanned by trees where the two
  leftmost inputs go to the same internal vertex. This quotient
  complex is clearly acyclic so the module bi-transformation trees
  complex, like the associahedron trees complex, is contractible.
  
  Given a sequence of maps
  $\cModule_1\stackrel{f_1}{\longrightarrow}\cModule_2\stackrel{f_2}{\longrightarrow}\cModule_3$
  between $\Ainf$-modules over~$\Alg$, a module bi-transformation
  tree~$T$ with $n$ inputs can be turned into a map
  $M_1\kotimes{\Ground} A^{\kotimes{\Ground}(n-1)} \to M_3$ by turning
  each internal vertex of
  valence~$d$ into $\mu_d$ if it is not on the leftmost strand or
  $(m_1)_d$, $(f_1)_d$, $(m_2)_d$, $(f_2)_d$, or $(m_3)_d$ if it is on the
  leftmost strand and above $v_1$, at $v_1$, between $v_1$ and~$v_2$,
  at $v_2$, or below $v_2$, respectively.

  \begin{figure}
    \centering
    %Font is 12 point.
    \includegraphics{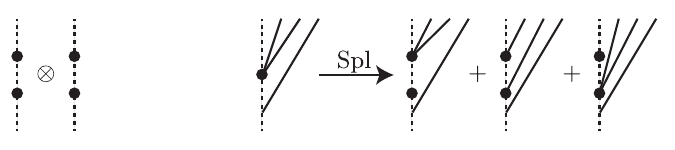}
    \caption[Module bi-transformation trees]{\textbf{Module
        bi-transformation trees.} Left: the pair of trees in
      $\TrModMulDiag_1$. Right: an example of the splitting operation $\Spl$.}
    \label{fig:mod-bi-tree}
  \end{figure}
  
  A \emph{module-bi-map tree diagonal} $\TrModMulDiag$ compatible with
  $\TrMDiag^1$ and $\TrMDiag^3$ is a sum of pairs of module
  bi-transformation trees satisfying the compatibility and
  non-degeneracy conditions of Definition~\ref{def:mod-map-diag}, but
  with module bi-transformation trees in place of module
  transformation trees. (For instance, $\TrModMulDiag_1$ is still the
  unique pair of module bi-transformation trees with one input; see
  Figure~\ref{fig:mod-bi-tree}.)

  Both sides of Equation~\eqref{eq:mod-map-tensor-comp} give
  module-bi-map tree diagonals, as follows. Let $\TrModMulDiag^{ij}$
  be the module-bi-map tree diagonal corresponding to the
  module-bi-map diagonal $\ModMulDiag^{ij}$.  For a module
  transformation tree~$T$, define $\Spl(T)$ to be the sum of module
  bi-transformation trees obtained by splitting the distinguished
  vertex of~$T$ into two distinguished vertices in all possible
  ways; see Figure~\ref{fig:mod-bi-tree}.
  For a pair of trees $T_1 \otimes T_2$, define
  $\Spl(T_1 \otimes T_2) = \Spl(T_1) \otimes \Spl(T_2)$. Then the left side
  of Equation~\eqref{eq:mod-map-tensor-comp} is the map corresponding to
  $\Spl(\TrModMulDiag^{13})$, by the definition of composition of
  $\Ainf$-module maps.

  Given two module transformation trees $T_1$ and~$T_2$, the composition
  $T_1 \circ_1 T_2$ is naturally a module bi-transformation tree.  Similarly,
  given two pairs of module transformation trees $(S_1,T_1)$ and $(S_2,T_2)$,
  $(S_2\circ S_1,T_2\circ T_1)$ is a pair of module bi-transformation
  trees. Then the right side of Equation~\eqref{eq:mod-map-tensor-comp} is the map
  corresponding to $\TrModMulDiag^{23} \circ_1 \TrModMulDiag^{12}$.

  Let us check that $\Spl(\TrModMulDiag^{13})$ and
  $\TrModMulDiag^{12} \circ_1 \TrModMulDiag^{23}$ are both module-bi-map tree diagonals. Normalization is
  immediate in both cases. Compatibility follows easily from
  \begin{align*}
    \bdy\Spl(\TrModMulDiag) &= \Spl(\bdy\TrModMulDiag)\\
    \bdy(\TrModMulDiag^1 \circ_1 \TrModMulDiag^2) &= 
      (\bdy\TrModMulDiag^1) \circ_1 \TrModMulDiag^2 +
      \TrModMulDiag^1 \circ_1 (\bdy\TrModMulDiag^2).
  \end{align*}

  Just as for Proposition~\ref{prop:mod-map-diag-homotopic}, any two
  module-bi-map tree diagonals are homotopic: we can make a
  chain complex out of linear combinations of pairs of module
  bi-transformation trees so that the diagonals are cycles, and then
  inductively build up a homotopy between any two diagonals. The
  homotopy between $\Spl(\TrModMulDiag^{13})$ and
  $\TrModMulDiag^{12}\circ_1\TrModMulDiag^{23}$ gives the desired
  homotopy between the two sides of
  Equation~\eqref{eq:mod-map-tensor-comp}.
\end{proof}

\begin{lemma}\label{lem:mod-iso}
  Let $\cModule$ and $\cNodule$ be $\Ainf$-modules over $\Alg$ and $f\co \cModule\to\cNodule$ a homomorphism. Then $f$ is an isomorphism if and only if $f_1$ is an isomorphism.
\end{lemma}
(See also Lemma~\ref{lem:Alg-iso-is}.)
\begin{proof}
  Certainly if $f$ is an isomorphism then $f_1$ is an isomorphism. To prove the converse, it suffices to construct left and right inverses to $f$. We will construct the left inverse $g$ to $f$; construction of the right inverse is similar. Let $g_1$ be the inverse to $f_1$. Suppose we have defined $g_i$ for $i\leq n$ so that $(g\circ f)_i=0$ for $1<i\leq n$. We have 
  \[
    (g\circ f)_{n+1}(x,a_1,\dots,a_n)=\sum_{k=0}^n g_{n-k+1}(f_{k+1}(x,a_1,\dots,a_k),a_{k+1},\dots,a_n).
  \]
  So, define 
  \[
    g_{n+1}(y,a_1,\dots,a_n)=\sum_{k=1}^n g_{n-k+1}(f_{k+1}(g_1(y),a_1,\dots,a_k),a_{k+1},\dots,a_n).
  \]
  It is clear that the $g$ constructed this way is a left inverse to
  $f$.  As noted above, a similar argument shows that $f$ has a right
  inverse $h$, and it follows that $g=h$. It remains to verify that
  $g$ is an $\Ainf$-homomorphism, rather than just a morphism. But
  \[
    0=d(\Id)=d(g\circ f)=d(g)\circ f + g\circ d(f)=d(g)\circ f.
  \]
  Composing on the right by $g$ gives that $d(g)=0$, as desired.
\end{proof}

\begin{corollary}\label{cor:mod-id-tens-id}
  Given $\Ainf$-modules $\cModule$ and $\cNodule$, module diagonals
  $\MDiag^1$ and $\MDiag^2$, and a module-map diagonal $\ModMulDiag$
  compatible with $\MDiag^1$ and $\MDiag^2$, the map
  \[
    \Id_{\cModule}\MDtp[\ModMulDiag]\Id_{\cNodule}\co \cModule\MDtp[\MDiag^1]\cNodule\to \cModule\MDtp[\MDiag^2]\cNodule
  \]
  is an isomorphism.
\end{corollary}
\begin{proof}
  The fact that $(\Id_{\cModule}\MDtp[\ModMulDiag]\Id_{\cNodule})_1$
  is an isomorphism is immediate from the non-degeneracy condition for
  $\ModMulDiag$. The corollary now follows from
  Lemma~\ref{lem:mod-iso}.
\end{proof}

\begin{corollary}\label{cor:mod-id-tens-id-is-sometimes-id}
  Given a $\Ainf$-modules $\cModule$ and $\cNodule$, a module diagonal
  $\MDiag$, and a module-map diagonal $\ModMulDiag$
  compatible with $\MDiag$ and $\MDiag$, the map
  \[
    \Id_{\cModule}\MDtp[\ModMulDiag]\Id_{\cNodule}\co \cModule\MDtp[\MDiag]\cNodule\to \cModule\MDtp[\MDiag]\cNodule
  \]
  is homotopic to the identity map.
\end{corollary}
\begin{proof}
  By Proposition~\ref{prop:dg-bifunctor},
    \[
    (\Id_{\cModule}\MDtp[\ModMulDiag]\Id_{\cNodule})\circ(\Id_{\cModule}\MDtp[\ModMulDiag]\Id_{\cNodule}) \sim (\Id_{\cModule}\MDtp[\ModMulDiag]\Id_{\cNodule}).
  \]
  By Corollary~\ref{cor:mod-id-tens-id},
  $(\Id_{\cModule}\MDtp[\ModMulDiag]\Id_{\cNodule})$ is an
  isomorphism. Composing both sides with
  $(\Id_{\cModule}\MDtp[\ModMulDiag]\Id_{\cNodule})^{-1}$ gives the
  result.
\end{proof}

\begin{remark}
  Even when $\MDiag^1=\MDiag^2$,
  $\Id_{\cModule}\MDtp[\ModMulDiag]\Id_{\cNodule}$ may not be the
  identity map unless $\ModMulDiag$, and perhaps the $\MDiag^i$, are
  chosen carefully. Specifically, the differential of the pair of
  trees
  %Font is 12 point
  \begin{center}
    \includegraphics{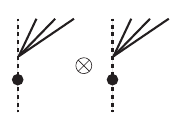}
  \end{center}
  is in the same grading as $\TrModMulDiag_4$. If $\TrModMulDiag$ has
  the property that $\Id_{\cModule}\MDtp[\ModMulDiag]\Id_{\cNodule}$
  is necessarily the identity map then the new module multiplihedron
  diagonal obtained by adding the differential of this pair of trees
  to $\TrModMulDiag_4$ no longer has this property.
\end{remark}

\begin{proof}[Proof of Theorem~\ref{thm:ModuleDiagonalExists}]
  Existence of a compatible module diagonal is trivial: we can take
  $\MDiag=\AsDiag$ (compare Example~\ref{eg:M-diag}). It remains to
  verify the rest of the theorem.

  The fact that the tensor product is an $\Ainf$-module is immediate
  from the definition and Lemma~\ref{lem:Ainf-mod-is}.

  Part~\ref{item:Mod-thm-dg} follows from the non-degeneracy condition
  for module diagonals.  
  Part~\ref{item:Mod-thm-qi} follows from
  Proposition~\ref{prop:dg-bifunctor} and Corollary~\ref{cor:mod-id-tens-id}.\signissue 
  The homotopy equivalence for
  Part~\ref{item:Mod-thm-change-diag} is obtained by tensoring the
  identity maps $\Id_{\cModule_1}$ and $\Id_{\cModule_2}$ using any module-map
  diagonal compatible with $\MDiag_1$ and $\MDiag_2$; existence of
  such a module-map diagonal is guaranteed by
  Lemma~\ref{lem:mod-map-diag-exists}.
\end{proof}

Again, in Section~\ref{sec:typeD}, some boundedness properties of our
modules will be relevant. 
\begin{definition}\label{def:mod-bonsai}
  An $\Ainf$-module $\cModule$ is \emph{bonsai} if there is an integer $N$ so that for any associahedron tree $T$ with dimension $\dim(T)>N$, $m^\cModule(T)=0$. An $\Ainf$-module map $f\co \cModule\to\cNodule$ is \emph{bonsai} if there is an integer $N$ so that for any module transformation tree $T$ with $\dim(T)>N$, $f(T)=0$.
\end{definition}

It is immediate from the definition that the bonsai morphisms form a subcomplex of $\Mor(\cModule,\cNodule)$.

\begin{lemma}\label{lem:mod-bounded}
  Fix an associahedron diagonal $\AsDiag$, module
  diagonals~$\MDiag^1$ and~$\MDiag^2$ compatible with $\AsDiag$, and
  module-map diagonal~$\ModMulDiag$ compatible with $\MDiag^1$ and $\MDiag^2$.
  If $\Alg_1$ and $\Alg_2$ are bonsai
  $\Ainf$-algebras and $\cModule_1$ and $\cModule_2$ are bonsai
  $\Ainf$-modules over $\Alg_1$ and $\Alg_2$ then
  $\cModule_1\otimes_{\MDiag^1}\cModule_2$ is also bonsai. Similarly, if
  $\cModule_i$, $\cNodule_i$, and $f^i\co\cModule_i\to\cNodule_i$ are
  bonsai for $i=1,2$ then so is $f^1\MDtp[\ModMulDiag]f^2$.
\end{lemma}
\begin{proof}
  As was the case for Lemma~\ref{lem:alg-bounded}, this is immediate
  from the fact that $\MDiag$ (respectively $\ModMulDiag$) is
  grading-preserving.
\end{proof}

\subsection{Type \texorpdfstring{$D$}{D} and \texorpdfstring{\DD}{DD} structures}\label{sec:typeD}
Next we turn to type $D$ structures and their bimodule analogue, type
\DD\ structures. We will focus on left type $D$ structures and
left-left type \DD\ structures.

Fix an $\Ainf$-algebra $\Alg$ over $\Ground$. Recall from the
introduction that a (left) type $D$ structure over $\Alg$ consists of a projective
$\Ground$-module $P$ together with a $\Ground$-module map 
$\delta^1_P\co P \to (A\kotimes{\Ground} P)\grs{1}$ satisfying
Equation~\eqref{eq:typeD}. We need to make sure the sums involved are finite.

\begin{definition}
  A type $D$ structure is \emph{bounded} if all sufficiently large
  iterates of the structure map~$\delta^1$ vanish.
\end{definition}

\begin{remark}
  Sometimes (e.g., in~\cite{LOT2}), we refer to our notion of being
  ``bounded'' as being \emph{operationally bounded}, since a bounded
  type $D$ structure need not be supported in finitely many gradings.
\end{remark}

\begin{convention}\label{conv:d-bounded}
  When talking about a type $D$ structure~$P$ over $\cA$, we assume that
  either $P$ is bounded or $\cA$ is bonsai.
\end{convention}
(We will weaken this restriction somewhat in
Section~\ref{sec:boundedness}.)

\begin{definition}\label{def:D-category}
  Given type $D$ structures $\lsup{\Alg}P$ and $\lsup{\Alg}Q$ over $\Alg$,
  the morphism complex $\Mor(\lsup{\Alg}P,\lsup{\Alg}Q)$ consists of
  all $\Ground$-module maps $f^1\co P\to A\kotimes{\Ground} Q$. The differential on this complex
  is defined by
  \[
    d(f^1)=\sum_{m,n=0}^\infty
    (\mu_{m+1+n}\otimes\Id)\circ(\Id_{A^{\kotimes{\Ground}
        m+1}}\otimes\delta^n_Q)\circ(\Id_{A^{\kotimes{\Ground} m}}\otimes f^1)\circ\delta^m_P
  \]
  or, graphically,
  \tikzsetnextfilename{dstreqn}
  \[
    d(f^1)=
    \mathcenter{
    \begin{tikzpicture}[smallpic]
      \node at (0,0) (tc) {};
      \node at (0,-1) (delta1) {$\delta_P$};
      \node at (0,-2) (f) {$f^1$};
      \node at (0,-3) (delta2) {$\delta_Q$};
      \node at (0,-5) (bc) {};
      \node at (-1,-4) (mu) {$\mu$};
      \node at (-1,-5) (bl) {};
      \draw[dmoda] (tc) to (delta1);
      \draw[dmoda] (delta1) to (f);
      \draw[dmoda] (f) to (delta2);
      \draw[dmoda] (delta2) to (bc);
      \draw[taa,bend right=4] (delta1) to (mu);
      \draw[alga,bend right=2] (f) to (mu);
      \draw[taa] (delta2) to (mu);
      \draw[alga] (mu) to (bl);
    \end{tikzpicture}}.
  \]

  We will sometimes write $f\co \lsup{\Alg}P\to\lsup{\Alg}Q$ to
  indicate a morphism $f^1\co P\to A\kotimes{\Ground} Q$.

  These morphism complexes make the category of type $D$ structures
  into a (non-unital) $\Ainf$-category, with (higher) composition
  of a string $\lsup{\Alg}P_0\stackrel{f_1}{\longrightarrow}\lsup{\Alg}P_1\stackrel{f_2}{\longrightarrow}\cdots\stackrel{f_{n}}{\longrightarrow}\lsup{\Alg}P_n$ given by
  \tikzsetnextfilename{dcatainf}
  \[
    \circ_n(f_n,\dots,f_1)=\mathcenter{
    \begin{tikzpicture}[smallpic]
      \node at (0,0) (tc) {};
      \node at (0,-1) (delta1) {$\delta_{P_0}$};
      \node at (0,-2) (f1) {$f^1_1$};
      \node at (0,-3) (delta2) {$\delta_{P_1}$};
      \node at (0,-4) (f2) {$f^1_2$};
      \node at (0,-5) (vdots) {$\vdots$};
      \node at (0,-6) (f3) {$f^1_n$};
      \node at (0,-7) (delta3) {$\delta_{P_n}$};
      \node at (0,-9) (bc) {};
      \node at (-2,-8) (mu) {$\mu$};
      \node at (-2,-9) (bl) {};
      \draw[dmoda] (tc) to (delta1);
      \draw[dmoda] (delta1) to (f1);
      \draw[dmoda] (f1) to (delta2);
      \draw[dmoda] (delta2) to (f2);
      \draw[dmoda] (f2) to (vdots);
      \draw[dmoda] (vdots) to (f3);
      \draw[dmoda] (f3) to (delta3);
      \draw[dmoda] (delta3) to (bc);
      \draw[taa,bend right=6] (delta1) to (mu);
      \draw[alga,bend right=4] (f1) to (mu);
      \draw[taa,bend right=2] (delta2) to (mu);
      \draw[alga] (f2) to (mu);
      \draw[taa] (delta3) to (mu);
      \draw[alga] (f3) to (mu);
      \draw[alga] (mu) to (bl);
    \end{tikzpicture}}
  \]
  \cite[Lemma 2.2.27]{LOT2}.

  In the case that $\Alg$ is strictly unital (see
  Section~\ref{subsec:UnitsHomPert} below), the category of type $D$
  structures over $\Alg$ is also strictly unital, with strict unit
  \[
    \Id_P^1(x)=\unit\otimes x.
  \]
\end{definition}

Recall that a type \DD\ structure over \dg algebras $\Alg$ and $\Blg$
is a type $D$ structure over $\Alg\rotimes{\Ring}\Blg$.
\begin{definition}\label{def:DD}
  Let $\Alg$ and $\Blg$ be $\Ainf$-algebras over $\Ground_1$
  and~$\Ground_2$, respectively. Fix an
  associahedron diagonal $\AsDiag$. Then a \emph{left-left type \DD\
    structure over $\Alg$ and $\Blg$ with respect to $\AsDiag$} is
  a (left) type $D$ structure over $\Alg\ADtp \Blg$. So, the category of
  left-left type \DD\ structures inherits the structure of a non-unital
  $\Ainf$-category from Definition~\ref{def:D-category}.
\end{definition}

Working with left-left modules is awkward, so we will often
think of a left-left module over $\Alg$ and $\Blg$ as a left-right
module over $\Alg$ and $\Blg^\op$, with a structure operation
\[
\delta^1 \co P \to A \otimes P \otimes B^\op\grs{1}
\]
obtained from $\delta^1$ by switching the factors.
Graphically, we can denote the operation $\delta^1\co P\to A\otimes P\otimes B^\op\grs{1}$ on a left-right type $D$ structure as
\[
  \tikzsetnextfilename{def-DD-1}
\mathcenter{\begin{tikzpicture}[smallpic]
  \node at (0,0) (tc) {};
  \node at (0,-1) (delta) {$\delta^1$};
  \node at (0,-2) (bc) {};
  \node at (-1,-2) (bl) {};
  \node at (1,-2) (br) {};
  \draw[dmoda] (tc) to (delta);
  \draw[dmoda] (delta) to (bc);
  \draw[alga] (delta) to (bl);
  \draw[alga] (delta) to (br);
\end{tikzpicture}}.
\]
To get the order of operations correct, the portion of the diagrams to
the right of the dotted line should be folded over to put it on the
left. In particular, the inputs to the algebra operations appear in
the opposite order, and the trees on that side of the
diagonal~$\AsDiag$ will appear reversed. In other words, the diagrams 
respect the notions of closer to / farther from the factor $P$, not the notion of left / right.

Similarly to type $D$ structures, the map $\delta^1$ induces a map 
\[
\delta\co P\to \overline{\Tensor}^*(A\grs{1})\kotimes{\Ground_1} P\kotimes{\Ground_2} \overline{\Tensor}^*(B^\op)
\]
by
\[
\mathcenter{
  \tikzsetnextfilename{def-DD-2}
\begin{tikzpicture}[smallpic]
  \node at (0,0) (tc) {};
  \node at (0,-1) (delta) {$\delta$};
  \node at (0,-2) (bc) {};
  \node at (-1,-2) (bl) {};
  \node at (1,-2) (br) {};
  \draw[dmoda] (tc) to (delta);
  \draw[dmoda] (delta) to (bc);
  \draw[taa] (delta) to (bl);
  \draw[taa] (delta) to (br);
\end{tikzpicture}}
=
\mathcenter{
  \tikzsetnextfilename{def-DD-3}
\begin{tikzpicture}[smallpic]
  \node at (0,0) (tc) {};
  \node at (0,-1) (delta) {$\delta^1$};
  \node at (0,-2) (bc) {};
  \node at (-1,-2) (bl) {};
  \node at (1,-2) (br) {};
  \draw[dmoda] (tc) to (delta);
  \draw[dmoda] (delta) to (bc);
  \draw[alga] (delta) to (bl);
  \draw[alga] (delta) to (br);
\end{tikzpicture}
}
+
\mathcenter{
  \tikzsetnextfilename{def-DD-4}
\begin{tikzpicture}[smallpic]
  \node at (0,0) (tc) {};
  \node at (0,-1) (delta1) {$\delta^1$};
  \node at (0,-2) (delta2) {$\delta^1$};
  \node at (0,-3) (bc) {};
  \node at (-1,-3) (bl) {};
  \node at (1,-3) (br) {};
  \node at (-1.2,-3) (bll) {};
  \node at (1.2,-3) (brr) {};
  \draw[dmoda] (tc) to (delta1);
  \draw[dmoda] (delta1) to (delta2);
  \draw[dmoda] (delta2) to (bc);
  \draw[alga] (delta1) to (bll);
  \draw[alga] (delta1) to (brr);
  \draw[alga] (delta2) to (bl);
  \draw[alga] (delta2) to (br);
\end{tikzpicture}}
+\quad\cdots.
\]
The type \DD\ structure is \emph{bounded} if the image of $\delta$ is contained in 
${\Tensor}^*(A\grs{1})\kotimes{\Ground_1} P\kotimes{\Ground_2} {\Tensor}^*(B^\op)$.
If our associahedron tree diagonal is
\[
\TrDiag_n=\sum_{(S,T)\in\Trees_n\otimes\Trees_n} n_{(S,T)}(S\otimes T)\in \Trees_n\otimes \Trees_n,
\]
then the structure equation for a type \DD\ structure can be encoded
as
\[
\sum_n\sum_{(S,T)\in\Trees_n\otimes\Trees_n}
n_{(S,T)}
\mathcenter{
  \tikzsetnextfilename{def-DD-5}
\begin{tikzpicture}[smallpic]
  \node at (0,0) (tc) {};
  \node at (0,-1) (delta) {$\delta$};
  \node at (-1,-2) (S) {$\mu(S)$};
  \node at (1,-2) (T) {$\mu(T^\op)$};
  \node at (0,-3) (bc) {};
  \node at (-1,-3) (bl) {};
  \node at (1,-3) (br) {};
  \draw[dmoda] (tc) to (delta);
  \draw[dmoda] (delta) to (bc);
  \draw[taa] (delta) to (S);
  \draw[taa] (delta) to (T);
  \draw[alga] (S) to (bl);
  \draw[alga] (T) to (br);
\end{tikzpicture}}
=
d\left(
\mathcenter{
  \tikzsetnextfilename{def-DD-6}
\begin{tikzpicture}[smallpic]
  \node at (0,0) (tc) {};
  \node at (0,-1) (delta) {$\delta^1$};
  \node at (0,-2) (bc) {};
  \node at (-1,-2) (bl) {};
  \node at (1,-2) (br) {};
  \draw[dmoda] (tc) to (delta);
  \draw[dmoda] (delta) to (bc);
  \draw[alga] (delta) to (bl);
  \draw[alga] (delta) to (br);
\end{tikzpicture}
}\right)
\]
or even more succinctly
\begin{equation}\label{eq:DD-struct}
\mathcenter{
  \tikzsetnextfilename{def-DD-7}
\begin{tikzpicture}[smallpic]
  \node at (0,0) (tc) {};
  \node at (0,-1) (delta) {$\delta$};
  \node at (-1,-2) (S) {$\TrDiag$};
  \node at (1,-2) (T) {$\TrDiag^\op$};
  \node at (0,-3) (bc) {};
  \node at (-1,-3) (bl) {};
  \node at (1,-3) (br) {};
  \draw[dmoda] (tc) to (delta);
  \draw[dmoda] (delta) to (bc);
  \draw[taa] (delta) to (S);
  \draw[taa] (delta) to (T);
  \draw[alga] (S) to (bl);
  \draw[alga] (T) to (br);
\end{tikzpicture}}
=
d\left(
\mathcenter{
  \tikzsetnextfilename{def-DD-8}
\begin{tikzpicture}[smallpic]
  \node at (0,0) (tc) {};
  \node at (0,-1) (delta) {$\delta^1$};
  \node at (0,-2) (bc) {};
  \node at (-1,-2) (bl) {};
  \node at (1,-2) (br) {};
  \draw[dmoda] (tc) to (delta);
  \draw[dmoda] (delta) to (bc);
  \draw[alga] (delta) to (bl);
  \draw[alga] (delta) to (br);
\end{tikzpicture}
}\right).
\end{equation}
(See Equation~\eqref{eq:Ainf-tens-trees-graph} for this notation. On
the right-hand side, the operation $d$ is induced by $\mu_1$ on the tensor
product algebra $A\otimes B$.)
As in Convention~\ref{conv:d-bounded}, for these structure equations
to make sense we must require either that $\Alg$ and $\Blg$ are
bonsai, or that $P$ is bounded.

\begin{warning}\label{warn:non-unital-cat}
  Non-unital categories lack a notion of isomorphism, not to mention
  equivalences of categories. This makes it difficult to prove that
  changing the diagonal gives an equivalent category of type \DD\
  structures. One way to handle this issue is to consider homotopy
  unital algebras. See Section~\ref{sec:hu-D} and, in particular,
  Theorem~\ref{thm:DD-indep-diag} for independence of the category of
  type \DD\ structures from the choice of diagonal.
\end{warning}

\subsection{Box products of type \texorpdfstring{$\DD$}{DD} structures}\label{sec:box}
In Section~\ref{sec:TensorAinftyModules}, we considered external
tensor products of
$\Ainf$-modules, i.e., tensor products over the ground
ring~$\Ring$. We turn now to internal tensor products, i.e., tensor
products over an $\Ainf$-algebra~$\Alg$.  Recall that the
$\Ainf$-tensor product of $\Ainf$-modules $\cModule_\Alg$ and
$\lsub{\Alg}\cNodule$ is defined to be the $\Ring$-module\signissue
\[
M\kotimes{\Ground}\Barop(\Alg)\kotimes{\Ground} N,
\]
where $\Barop(\Alg)=\Tensor^*(A\grs{1})$ denotes the bar complex of
$\Alg$ (see, e.g.,~\cite[Section 6.3]{AinftyAlg} or~\cite[Section 3]{Keller:OtherAinfAlg}). This tensor product is equipped with a differential involving (all of) the
operations $m_i$ on $\cModule$ and $\cNodule$ and $\mu_i$ on
$\Alg$. This
$\Ainf$-tensor product is a generalization of the derived tensor product of
ordinary modules (or chain complexes of modules, or \dg modules).

If we replace $\lsub{\Alg}\cNodule$ with a type $D$ structure $\lsup{\Alg}P$,
we can form another kind of tensor product, denoted $\DT$
in~\cite{LOT1,LOT2}, which is closer to the classical
tensor product of ordinary modules. 
If $\cModule$ is bonsai or $P$ is bounded then we can define a differential $\bdy^\DT$ on
$M\kotimes{\Ground} P$ by the composition
\[
M\kotimes{\Ground} P\xrightarrow{\Id\otimes\delta}M\kotimes{\Ground}\overline{\Tensor}^*(A)\kotimes{\Ground} P\xrightarrow{m\otimes\Id}M\kotimes{\Ground} P,
\]
i.e.,
\[
  \bdy^\DT=
  \tikzsetnextfilename{sec-box-1}
  \mathcenter{\begin{tikzpicture}[smallpic]
    \node at (0,0) (tl) {};
    \node at (2,0) (tr) {};
    \node at (2,-1) (delta) {$\delta_P$};
    \node at (0,-2) (m) {$m_M$};
    \node at (2,-3) (br) {};
    \node at (0,-3) (bl) {};
    \draw[moda] (tl) to (m);
    \draw[moda] (m) to (bl);
    \draw[dmoda] (tr) to (delta);
    \draw[dmoda] (delta) to (br);
    \draw[taa] (delta) to (m);
  \end{tikzpicture}}.
\]
Let $\cModule\DT P$ denote the complex $(M\kotimes{\Ground} P,\bdy^\DT)$. Further, 
given $\phi\in\Mor(\cModule_1,\cModule_2)$ and a sequence of
  maps $f_i^1\in \Mor(P_i,P_{i+1})$, $i=1,\dots {n}$ (for any $n\geq 1$) define
  \[
    \DT_{0,n}(;f_1^1,\dots,f_n^1)=
  \tikzsetnextfilename{sec-box-2}
    \mathcenter{\begin{tikzpicture}[smallpic]
        \node at (0,0) (tc) {};
        \node at (2,0) (tr) {};
        \node at (2,-1) (r1) {$\delta_{P_1}$};
        \node at (2,-2) (r2) {$f_1^1$};
        \node at (2,-3) (r3) {$\delta_{P_2}$};
        \node at (2,-4) (r4) {$\vdots$};
        \node at (2,-5) (r5) {$f_{n}^1$};
        \node at (2,-6) (r6) {$\delta_{P_{n+1}}$};
        \node at (0,-7) (phi) {$m_{M_1}$};
        \node at (0,-8) (bc) {};
        \node at (2,-8) (br) {};
        \draw[dmoda] (tr) to (r1);
        \draw[dmoda] (r1) to (r2);
        \draw[dmoda] (r2) to (r3);
        \draw[dmoda] (r3) to (r4);
        \draw[dmoda] (r4) to (r5);
        \draw[dmoda] (r5) to (r6);
        \draw[dmoda] (r6) to (br);
        \draw[taa,bend right=4] (r1) to (phi);
        \draw[alga,bend right=2] (r2) to (phi);
        \draw[taa] (r3) to (phi);
        \draw[alga] (r5) to (phi);
        \draw[taa] (r6) to (phi);
        \draw[moda] (tc) to (phi);
        \draw[moda] (phi) to (bc);
    \end{tikzpicture}}
  \qquad\qquad
    \DT_{1,n}(\phi;f_1^1,\dots,f_n^1)=
  \tikzsetnextfilename{sec-box-3}
    \mathcenter{\begin{tikzpicture}[smallpic]
        \node at (0,0) (tc) {};
        \node at (2,0) (tr) {};
        \node at (2,-1) (r1) {$\delta_{P_1}$};
        \node at (2,-2) (r2) {$f_1^1$};
        \node at (2,-3) (r3) {$\delta_{P_2}$};
        \node at (2,-4) (r4) {$\vdots$};
        \node at (2,-5) (r5) {$f_{n}^1$};
        \node at (2,-6) (r6) {$\delta_{P_{n+1}}$};
        \node at (0,-7) (phi) {$\phi$};
        \node at (0,-8) (bc) {};
        \node at (2,-8) (br) {};
        \draw[dmoda] (tr) to (r1);
        \draw[dmoda] (r1) to (r2);
        \draw[dmoda] (r2) to (r3);
        \draw[dmoda] (r3) to (r4);
        \draw[dmoda] (r4) to (r5);
        \draw[dmoda] (r5) to (r6);
        \draw[dmoda] (r6) to (br);
        \draw[taa] (r1) to (phi);
        \draw[taa] (r3) to (phi);
        \draw[taa] (r6) to (phi);
        \draw[alga] (r2) to (phi);
        \draw[alga] (r5) to (phi);
        \draw[moda] (tc) to (phi);
        \draw[moda] (phi) to (bc);
    \end{tikzpicture}}.
  \]
  Define $\DT_{i,j}$ to vanish if $i>1$.

The following is a slight strengthening of~\cite[Lemma~2.30]{LOT1} and~\cite[Lemma~2.3.3]{LOT2}, the proof of which is left to the reader:
\begin{lemma}\label{lem:DT-mod-defined}
  If $\cModule$ is bonsai or $P$ is bounded then $\cModule\DT P$ is a chain
  complex. Moreover, the operation~$\DT$ is a (non-unital)
  $\Ainf$-bifunctor,
  i.e., the operations $\DT_{i,j}$ satisfy the $\Ainf$-bimodule relations. 
\end{lemma}

\begin{remark}
  If $\Alg$ and $\cModule$ are strictly unital and we restrict to the
  strictly unital morphism complex then $\DT$ is a strictly unital
  bifunctor.
\end{remark}

\begin{remark}
  Note that $\DT_{1,1}$ is a higher composition of $\phi$ and $f^1_1$:
  the ordinary composition is defined to be either
  $\DT_{0,1}(;f^1_1)\circ \DT_{1,0}(\phi;)$ or $\DT_{1,0}(\phi;)\circ
  \DT_{0,1}(;f^1_1)$, a global choice.
\end{remark}

We turn now to the generalizations of $\DT$ to
bimodules. The most obvious generalization is a triple tensor product
of a \DD\ structure with two $\Ainf$-modules.
\begin{definition}\label{def:triple-prod}
  Fix an associahedron diagonal $\AsDiag$ and a module diagonal
  $\MDiag$ compatible with $\AsDiag$. Let $\lsup{\Alg,\Blg}P$ be
  a left-left type \DD\ structure over $\Alg$ and $\Blg$, and let
  $\cModule_{\Alg}$, $\cNodule_\Blg$ be $\Ainf$-modules over $\Alg$ and
  $\Blg$. Assume that either $P$ is bounded or $\cModule$,
  $\cNodule$, $\Alg$, and $\Blg$ are all bonsai. Define
  \[
  [\cModule_{\Alg}\DT \lsup{\Alg}P^{\Blg^\op}\DT \lsub{\Blg^\op}\cNodule]_{\MDiag}
  \coloneqq
  (\cModule\otimes_\MDiag \cNodule)_{\Alg\ADtp\Blg} \,\DT\, \lsup{\Alg\ADtp \Blg}P.
  \]
\end{definition}

Recall that Theorem~\ref{thm:TripleTensorProduct} asserts that the
triple box product is a chain complex and, up to homotopy equivalence,
is independent of the choice of module diagonal.

\begin{proof}[Proof of Theorem~\ref{thm:TripleTensorProduct}]
  The fact that the triple box product is a well-defined chain
  complex is immediate from its definition, together with
  Lemmas~\ref{lem:DT-mod-defined} and~\ref{lem:mod-bounded}. Next, if $\MDiag$ and $\MDiag'$ are
  different module diagonals, it follows from
  Part~\ref{item:Mod-thm-change-diag} of
  Theorem~\ref{thm:ModuleDiagonalExists} that $\cModule\otimes_\MDiag \cNodule$
  and $\cModule\otimes_{\MDiag'} \cNodule$ are homotopy equivalent. If
  $\Alg$, $\Blg$,
  $\cModule$, and $\cNodule$ are bonsai,  
   Lemma~\ref{lem:mod-bounded} implies that the
  homotopy equivalence is also bonsai.
  So, invariance of the box product under
  homotopy equivalence,~\cite[Corollary~2.3.5]{LOT2},
  implies that the triple box products are homotopy equivalent as
  well.
\end{proof}

\begin{remark}\label{remark:ttp-indep-AsDiag}
  There is also a sense in which the triple box product is
  independent of the associahedron diagonal $\AsDiag$, but since the
  structure equation for a type \DD\ structure uses an associahedron
  diagonal, this is rather cumbersome to spell out, and is left to the
  reader.
\end{remark}

Our next goal is to interpret this triple box product as tensoring of
$\cModule$ with~$P$ and then tensoring the result with $\cNodule$. 
This interpretation uses the following notion, which
requires $\Blg$ to be strictly unital. 

\begin{definition}\label{def:one-side-DT}
  Fix a module diagonal primitive
  \[ 
 \TrPMDiag_n=\sum_{(S,T)}n_{S,T}(S,T)\in \Trees_n\rotimes{\Ring}
  \Trees_{n-1}
  \]
  compatible with an associahedron diagonal $\AsDiag$. Let $\Alg$,
  $\Blg$ be $\Ainf$-algebras over $\Ground_1$, $\Ground_2$
  respectively, and assume that $\Blg$ is strictly unital with unit $\unit$. Let
  $\lsup{\Alg,\Blg}P$ be a type \DD\ structure with respect
  to~$\AsDiag$ and let $\cModule_\Alg$ be an $\Ainf$-module
  over~$\Alg$. Assume that either $P$ is bounded or 
  $\cModule$, $\Alg$, and $\Blg$ are bonsai.
  Define $\cModule\DT^{\TrPMDiag}P$ to be the vector space
  $M\kotimes{\Ground_1}P$ equipped with the operation
  \[
  \delta^1\co M\DT^{\TrPMDiag}P\to (M\DT^{\TrPMDiag}P)
  \kotimes{\Ground_2} B^\op
  \]
  given by
  \begin{equation*}
  \delta^1(x\otimes y)=m_1^M(x)\otimes y \otimes \unit+
  \sum_{(S,T)}n_{S,T}(m^M(S) \otimes \Id_{P}\otimes
  \mu^{\Blg^{\op}}(T^\op))\circ (\Id_M \otimes\delta_P).
\end{equation*}
This is shown graphically in Figure~\ref{fig:one-side-DT}.
\end{definition}
\begin{figure}
  \begin{align*}
    \delta^1(x\otimes y)
    &=
    \mathcenter{
  \tikzsetnextfilename{def-one-side-DT-1}
      \begin{tikzpicture}[smallpic]
        \node at (-1,0) (tl) {$x$};
        \node at (0,0) (tc) {$y$};
        \node at (-1,-1) (m) {$m_1$};
        \node at (1,-1) (one) {$1$};
        \node at (-1,-2) (bl) {};
        \node at (0,-2) (bc) {};
        \node at (1,-2) (br) {};
        \draw[dmoda] (tc) to (bc);
        \draw[moda] (tl) to (m);
        \draw[moda] (m) to (bl);
        \draw[alga] (one) to (br);
      \end{tikzpicture}
    }+\,
    \sum_{S,T}
    n_{S,T}
    \mathcenter{
  \tikzsetnextfilename{def-one-side-DT-2}
      \begin{tikzpicture}[smallpic]
        \node at (-1.5,0) (x) {$x$};
        \node at (0,0) (y) {$y$};
        \node at (0, -1) (delta1) {$\delta^1$};
        \node at (0, -2) (delta2) {$\delta^1$};
        \node at (0, -3) (vdots) {$\vdots$};
        \node at (0, -4) (delta3) {$\delta^1$};
        \node at (-1.5, -4.5) (S) {$S$};
        \node at (1.5, -4.5) (T) {$T^\op$};
        \node at (-1.5, -5.5) (bl) {};
        \node at (0, -5.5) (bc) {};
        \node at (1.5, -5.5) (br) {};
        \draw[moda] (x) to (S);
        \draw[moda] (S) to (bl);
        \draw[dmoda] (y) to (delta1);
        \draw[dmoda] (delta1) to (delta2);
        \draw[dmoda] (delta2) to (vdots);
        \draw[dmoda] (vdots) to (delta3);
        \draw[dmoda] (delta3) to (bc);
        \draw[alga] (T) to (br);
        \draw[alga] (delta1) to (S);
        \draw[alga] (delta2) to (S);
        \draw[alga] (delta3) to (S);
        \draw[blga] (delta1) to (T);
        \draw[blga] (delta2) to (T);
        \draw[blga] (delta3) to (T);
      \end{tikzpicture}}\\
    &=
    \mathcenter{
  \tikzsetnextfilename{def-one-side-DT-3}
      \begin{tikzpicture}[smallpic]
        \node at (-1,0) (tl) {$x$};
        \node at (0,0) (tc) {$y$};
        \node at (-1,-1) (m) {$m_1$};
        \node at (1,-1) (one) {$\unit$};
        \node at (-1,-2) (bl) {};
        \node at (0,-2) (bc) {};
        \node at (1,-2) (br) {};
        \draw[dmoda] (tc) to (bc);
        \draw[moda] (tl) to (m);
        \draw[moda] (m) to (bl);
        \draw[alga] (one) to (br);
      \end{tikzpicture}
    }+
    \mathcenter{
  \tikzsetnextfilename{def-one-side-DT-4}
      \begin{tikzpicture}[smallpic]
        \node at (-1,0) (x) {$x$};
        \node at (0,0) (y) {$y$};
        \node at (0, -1) (delta) {$\delta$};
        \node at (-1, -2) (S) {$\TrPMDiag$};
        \node at (1, -2) (T) {$\TrPMDiag^\op$};
        \node at (-1, -3) (bl) {};
        \node at (0, -3) (bc) {};
        \node at (1, -3) (br) {};
        \draw[dmoda] (y) to (delta);
        \draw[dmoda] (delta) to (bc);
        \draw[moda] (x) to (S);
        \draw[moda] (S) to (bl);
        \draw[taa] (delta) to (S);
        \draw[tbb] (delta) to (T);
        \draw[blga] (T) to (br);
      \end{tikzpicture}}  
  \end{align*}
  \caption[The structure operation on the one-sided box tensor product]{The structure operation on $\cModule \DT^{\TrPMDiag} P$, shown
    graphically with two different notations. We will prefer the
    second notation.}
  \label{fig:one-side-DT}
\end{figure}

\begin{lemma}\label{lem:DT-defined}
  The object $(\cModule\DT^{\TrPMDiag}P,\delta^1)$ from
  Definition~\ref{def:one-side-DT} is a right type $D$ structure over
  $\Blg^\op$, or equivalently a left type $D$ structure over $\Blg$.
\end{lemma}
\begin{proof}
  The boundedness assumptions imply that the sum involved in defining
  $\delta^1$ is finite. It remains to verify that $\bdy^2=0$ or, more
  precisely, that Equation~(\ref{eq:typeD}) holds for
  $P\DT^{\TrPMDiag}\cModule$.  We will give the proof in compact graphical
  notation. Using the fact that $\Blg$ is strictly unital,
  Equation~(\ref{eq:typeD}) is equivalent to
\[ 
    \mathcenter{
  \tikzsetnextfilename{lem-DT-defined-1}
      \begin{tikzpicture}[smallpic]
        \node at (0,0) (tc) {};
        \node at (0,-1) (delta1) {$\delta$};
        \node at (0,-2) (vdots) {$\vdots$};
        \node at (0,-3) (delta2) {$\delta$};
        \node at (0,-6) (bc) {};
        \node at (-1,0) (tl) {};
        \node at (-1,-2) (p1) {$\TrPMDiag$};
        \node at (-1,-3) (lvdots) {$\vdots$};
        \node at (-1,-4) (p2) {$\TrPMDiag$};
        \node at (-1,-6) (bl) {};
        \node at (1,-2) (p3) {$\TrPMDiag^\op$};
        \node at (1,-4) (p4) {$\TrPMDiag^\op$};
        \node at (2, -5) (mu) {$\mu$};
        \node at (2,-6) (br) {};
        \draw[dmoda] (tc) to (delta1);
        \draw[dmoda] (delta1) to (vdots);
        \draw[dmoda] (vdots) to (delta2);
        \draw[dmoda] (delta2) to (bc);
        \draw[taa] (delta1) to (p1);
        \draw[tbb] (delta1) to (p3);
        \draw[taa] (delta2) to (p2);
        \draw[tbb] (delta2) to (p4);
        \draw[alga] (p3) to (mu);
        \draw[alga] (p4) to (mu);
        \draw[alga] (mu) to (br);
        \draw[moda] (tl) to (p1);
        \draw[moda] (p1) to (lvdots);
        \draw[moda] (lvdots) to (p2);
        \draw[moda] (p2) to (bl);
      \end{tikzpicture}}
      +
  \tikzsetnextfilename{lem-DT-defined-2}
      \mathcenter{
        \begin{tikzpicture}[smallpic]
          \node at (0,0) (tc) {};
          \node at (0,-2) (delta) {$\delta$};
          \node at (0,-4) (bc) {};
          \node at (-1,0) (tl) {};
          \node at (-1,-1) (m) {$m_1$};
          \node at (-1,-3) (pl) {$\TrPMDiag$};
          \node at (-1, -4) (bl) {};
          \node at (1,-3) (pr) {$\TrPMDiag^\op$};
          \node at (1,-4) (br) {};
          \draw[dmoda] (tc) to (delta);
          \draw[dmoda] (delta) to (bc);
          \draw[tbb] (delta) to (pl);
          \draw[moda] (tl) to (m);
          \draw[moda] (m) to (pl);
          \draw[moda] (pl) to (bl);
          \draw[taa] (delta) to (pr);
          \draw[alga] (pr) to (br);
        \end{tikzpicture}
      }
      +
  \tikzsetnextfilename{lem-DT-defined-3}
      \mathcenter{
        \begin{tikzpicture}[smallpic]
          \node at (0,0) (tc) {};
          \node at (0,-1) (delta) {$\delta$};
          \node at (0,-4) (bc) {};
          \node at (-1,0) (tl) {};
          \node at (-1,-3) (m) {$m_1$};
          \node at (-1,-2) (pl) {$\TrPMDiag$};
          \node at (-1, -4) (bl) {};
          \node at (1,-2) (pr) {$\TrPMDiag^\op$};
          \node at (1,-4) (br) {};
          \draw[dmoda] (tc) to (delta);
          \draw[dmoda] (delta) to (bc);
          \draw[tbb] (delta) to (pl);
          \draw[moda] (tl) to (pl);
          \draw[moda] (pl) to (m);
          \draw[moda] (m) to (bl);
          \draw[taa] (delta) to (pr);
          \draw[alga] (pr) to (br);
        \end{tikzpicture}
      }
      =0.
      \]
  The module diagonal primitive structure equation,
  Equation~(\ref{eq:M-prim-compat}), together
  Lemmas~\ref{lem:A-inf-alg-is} and~\ref{lem:Ainf-mod-is}, implies
  that the sum of these three terms is equal to
  \[
  \mathcenter{
  \tikzsetnextfilename{lem-DT-defined-4}
    \begin{tikzpicture}[smallpic]
      \node at (0,0) (tc) {};
      \node at (0,-1) (delta1) {$\delta$};
      \node at (0,-2) (delta2) {$\delta$};
      \node at (0,-3) (delta3) {$\delta$};
      \node at (0,-5) (bc) {};
      \node at (-2,-4) (pl) {$\TrPMDiag$};
      \node at (2,-4) (pr) {$\TrPMDiag^\op$};
      \node at (-2,-5) (bl) {};
      \node at (2,-5) (br) {};     
      \node at (-2,0) (tl) {};
      \node at (-1,-3) (lasdiag) {$\TrDiag$};
      \node at (1,-3) (rasdiag) {$\TrDiag$};
      \draw[dmoda] (tc) to (delta1);
      \draw[dmoda] (delta1) to (delta2);
      \draw[dmoda] (delta2) to (delta3);
      \draw[dmoda] (delta3) to (bc);
      \draw[moda] (tl) to (pl);
      \draw[taa] (delta1) to (pl);
      \draw[taa] (delta2) to (lasdiag);
      \draw[alga] (lasdiag) to (pl);
      \draw[taa] (delta3) to (pl);
      \draw[moda] (pl) to (bl);
      \draw[tbb] (delta1) to (pr);
      \draw[tbb] (delta2) to (rasdiag);
      \draw[blga] (rasdiag) to (pr);
      \draw[tbb] (delta3) to (pr);
      \draw[alga] (pr) to (br);
    \end{tikzpicture}
  }
  +
  \mathcenter{
  \tikzsetnextfilename{lem-DT-defined-5}
    \begin{tikzpicture}[smallpic]
      \node at (0,0) (tc) {};
      \node at (0, -1) (delta1) {$\delta$};
      \node at (0, -2) (delta2) {$\delta^1$};
      \node at (0,-3) (delta3) {$\delta$};
      \node at (0,-5) (bc) {};
      \node at (-2,0) (tl) {};
      \node at (-1,-3) (mu) {$\mu_1$};
      \node at (-2,-4) (pl) {$\TrPMDiag$};
      \node at (-2,-5) (bl) {};
      \node at (1,-4) (pr) {$\TrPMDiag^\op$};
      \node at (1,-5) (br) {};
      \draw[dmoda] (tc) to (delta1);
      \draw[dmoda] (delta1) to (delta2);
      \draw[dmoda] (delta2) to (delta3);
      \draw[dmoda] (delta3) to (bc);
      \draw[moda] (tl) to (pl);
      \draw[moda] (pl) to (bl);
      \draw[taa] (delta1) to (pl);
      \draw[alga] (delta2) to (mu);
      \draw[alga] (mu) to (pl);
      \draw[taa] (delta3) to (pl);
      \draw[alga] (pr) to (br);
      \draw[tbb] (delta1) to (pr);
      \draw[blga] (delta2) to (pr);
      \draw[tbb] (delta3) to (pr);
    \end{tikzpicture}
  }
  +
  \mathcenter{
  \tikzsetnextfilename{lem-DT-defined-6}
    \begin{tikzpicture}[smallpic]
      \node at (0,0) (tc) {};
      \node at (0, -1) (delta1) {$\delta$};
      \node at (0, -2) (delta2) {$\delta^1$};
      \node at (0,-3) (delta3) {$\delta$};
      \node at (0,-5) (bc) {};
      \node at (1,-3) (mu) {$\mu_1$};
      \node at (-1,0) (tl) {};
      \node at (-1,-4) (pl) {$\TrPMDiag$};
      \node at (-1,-5) (bl) {};
      \node at (2,-4) (pr) {$\TrPMDiag^\op$};
      \node at (2,-5) (br) {};
      \draw[dmoda] (tc) to (delta1);
      \draw[dmoda] (delta1) to (delta2);
      \draw[dmoda] (delta2) to (delta3);
      \draw[dmoda] (delta3) to (bc);
      \draw[moda] (tl) to (pl);
      \draw[moda] (pl) to (bl);
      \draw[taa] (delta1) to (pl);
      \draw[blga] (delta2) to (mu);
      \draw[blga] (mu) to (pr);
      \draw[taa] (delta3) to (pl);
      \draw[alga] (pr) to (br);
      \draw[tbb] (delta1) to (pr);
      \draw[alga] (delta2) to (pl);
      \draw[tbb] (delta3) to (pr);
    \end{tikzpicture}.
  }
  \]
  These terms cancel by the structure equation
  for a type \DD\ structure (Equation~\eqref{eq:DD-struct}).
\end{proof}

\begin{lemma}\label{lem:tp-is-tp}
  Fix an associahedron diagonal $\AsDiag$ and a module diagonal
  primitive $\TrPMDiag$ compatible with~$\AsDiag$.  Let $\MDiag$ be
  the module diagonal associated to $\TrPMDiag$
  (Definition~\ref{def:assoc-mod-diag}).  Let
  $\lsup{\Alg,\Blg}P$ be a type \DD\ structure over
  $\Alg$ and $\Blg$ and let $\cModule_{\Alg}$, $\cNodule_{\Blg}$ be
  $\Ainf$-modules over $\Alg$ and $\Blg$, respectively. Assume that either $P$
  is bounded or $\cModule$, $\cNodule$, $\Alg$, and $\Blg$ are
  bonsai. Assume also that $\Blg$ and $\cNodule$ are strictly unital. Then 
  \begin{equation}\label{eq:tp-is-tp}
    [\cModule_{\Alg}\DT \lsup{\Alg}P^{\Blg^\op}\DT \lsub{\Blg^\op}\cNodule]_{\MDiag}
    \cong
    (\cModule\DT^{\TrPMDiag}_{\Alg}P)\DT_{\Blg^\op}\cNodule.
  \end{equation}
\end{lemma}
\begin{proof}
  Graphically, both sides of Equation~\eqref{eq:tp-is-tp} are given by
  \[
  \left(
  \mathcenter{
  \tikzsetnextfilename{lem-tp-is-tp-1}
    \begin{tikzpicture}[smallpic]
      \node at (0,0) (tc) {};
      \node at (0,-1) (delta1) {$\delta$};
      \node at (0,-2) (vdots) {$\vdots$};
      \node at (0,-3) (delta2) {$\delta$};
      \node at (0,-6) (bc) {};
      \node at (-1,-2) (pl1) {$\TrPMDiag$};
      \node at (-1,-3) (lvdots) {$\vdots$};
      \node at (-1,-4) (pl2) {$\TrPMDiag$};
      \node at (1,-2) (pr1) {$\TrPMDiag^\op$};
      \node at (1,-3) (rvdots) {$\vdots$};
      \node at (1,-4) (pr2) {$\TrPMDiag^\op$};
      \node at (-1,0) (tl) {};
      \node at (2,0) (tr) {};
      \node at (2,-6) (br) {};
      \node at (-1,-6) (bl) {};
      \node at (2,-5) (m) {$m$};
      \draw[dmoda] (tc) to (delta1);
      \draw[dmoda] (delta1) to (vdots);
      \draw[dmoda] (vdots) to (delta2);
      \draw[dmoda] (delta2) to (bc);
      \draw[taa] (delta1) to (pl1);
      \draw[taa] (delta2) to (pl2);
      \draw[alga] (pr1) to (m);
      \draw[alga] (pr2) to (m);
      \draw[moda] (tr) to (m);
      \draw[moda] (m) to (br);
      \draw[moda] (tl) to (pl1);
      \draw[moda] (pl1) to (lvdots);
      \draw[moda] (lvdots) to (pl2);
      \draw[moda] (pl2) to (bl);
      \draw[tbb] (delta1) to (pr1);
      \draw[tbb] (delta2) to (pr2);
    \end{tikzpicture}
  }
  \right)
  +
  \left(
  \mathcenter{
  \tikzsetnextfilename{lem-tp-is-tp-2}
    \begin{tikzpicture}[smallpic]
      \node at (-1,0) (tl) {};
      \node at (0,0) (tc) {};
      \node at (1,0) (tr) {};
      \node at (-1,-2) (bl) {};
      \node at (0,-2) (bc) {};
      \node at (1,-2) (br) {};
      \node at (-1,-1) (m) {$m_1$};
      \draw[moda] (tr) to (br);
      \draw[dmoda] (tc) to (bc);
      \draw[moda] (tl) to (m);
      \draw[moda] (m) to (bl);
    \end{tikzpicture}
  }\right).
  \]
\end{proof}

\begin{proof}[Proof of Theorem~\ref{thm:prim-DT}]
  This is a restatement of Lemmas~\ref{lem:DT-defined} and~\ref{lem:tp-is-tp}.
\end{proof}

\subsubsection{Functoriality of one-sided box tensor products}
We turn next to functoriality of
$\DT^{\TrPMDiag}$. Let $P$ and $Q$ be type \DD\ structures over
$\Alg$ and $\Blg$, with $\Blg$ strictly unital, and let $f\co P\to Q$
be a morphism. That is, $f$ is a 
map $P\to A\kotimes{\Ground_1} Q\kotimes{\Ground_2} B^\op$. Given an $\Ainf$-module $\cModule$
over $\Alg$, we can consider the morphism of right type $D$ structures over $\Blg$
\[
\Id_{\cModule}\DT^{\TrPMDiag}f\co \cModule\DT^{\TrPMDiag}{P}\to \cModule\DT^{\TrPMDiag}{Q}
\]
defined graphically by
\begin{equation}\label{eq:DD-map-tens-id}
\Id_{\cModule}\DT^{\TrPMDiag}f=
\mathcenter{
  \tikzsetnextfilename{eq-DD-map-tens-id-1}
  \begin{tikzpicture}[smallpic]
    \node at (0,0) (tc) {};
    \node at (-1,0) (tl) {};
    \node at (0, -1) (delta1) {$\delta$};
    \node at (0, -2) (f) {$f$};
    \node at (0, -3) (delta2) {$\delta$};
    \node at (-1, -4) (pl) {$\TrPMDiag$};
    \node at (1, -4) (pr) {$\TrPMDiag^\op$};
    \node at (-1, -5) (bl) {};
    \node at (0, -5) (bc) {};
    \node at (1, -5) (br) {};
    \draw[dmoda] (tc) to (delta1);
    \draw[dmoda] (delta1) to (f);
    \draw[dmoda] (f) to (delta2);
    \draw[dmoda] (delta2) to (bc);
    \draw[moda] (tl) to (pl);
    \draw[moda] (pl) to (bl);
    \draw[taa] (delta1) to (pl);
    \draw[taa] (delta1) to (pr);
    \draw[alga] (f) to (pl);
    \draw[blga] (f) to (pr);
    \draw[tbb] (delta2) to (pl);
    \draw[tbb] (delta2) to (pr);
    \draw[alga] (pr) to (br);
  \end{tikzpicture}}
\end{equation}

Similarly, fix a module-map primitive $\TrPMorDiag$ compatible with
module diagonal primitives $\TrPMDiag^1$ and $\TrPMDiag^2$. Given
another $\Ainf$-module $\cNodule$ over $\Alg$ and a morphism $g\co \cModule\to \cNodule$
we can consider the map
\[
  g \DT^{\TrPMorDiag}\Id_P
  \co \cModule\DT^{\TrPMDiag_1}P\to \cNodule\DT^{\TrPMDiag_2}P
\]
defined graphically by
\begin{equation}\label{eq:DD-id-tens-map}
  \mathcenter{
  \tikzsetnextfilename{eq-DD-id-tens-map-1}
  \begin{tikzpicture}[smallpic]
    \node at (0,0) (tc) {};
    \node at (-1,0) (tl) {};
    \node at (0,-1) (delta) {$\delta$};
    \node at (-1,-2) (pl) {$\TrPMorDiag g$};
    \node at (1,-2) (pr) {$\TrPMorDiag^\op$};
    \node at (-1,-3) (bl) {};
    \node at (0,-3) (bc) {};
    \node at (1,-3) (br) {};
    \draw[dmoda] (tc) to (delta);
    \draw[dmoda] (delta) to (bc);
    \draw[moda] (tl) to (pl);
    \draw[moda] (pl) to (bl);
    \draw[taa] (delta) to (pl);
    \draw[taa] (delta) to (pr);
    \draw[alga] (pr) to (br);
  \end{tikzpicture}}
\end{equation}
where $\TrPMorDiag g$ means that we apply $g$ at the distinguished
vertex of $\TrPMorDiag$. There is a special case corresponding to
$\delta^0$, which is the term
\[
  \mathcenter{
  \tikzsetnextfilename{eq-DD-id-tens-map-2}
  \begin{tikzpicture}[smallpic]
    \node at (0,0) (tc) {};
    \node at (-1,0) (tl) {};
    \node at (-1,-1) (pl) {$g_1$};
    \node at (1,-1) (pr) {$\unit$};
    \node at (-1,-2) (bl) {};
    \node at (0,-2) (bc) {};
    \node at (1,-2) (br) {};
    \draw[dmoda] (tc) to (bc);
    \draw[moda] (tl) to (pl);
    \draw[moda] (pl) to (bl);
    \draw[alga] (pr) to (br);
  \end{tikzpicture}},
\]
making use of the fact that $\Blg$ is strictly unital.
As usual, for these operations to be defined requires
some boundedness: either $P$ and $Q$ should be bounded or
$\Alg$, $\Blg$, $\cModule$, $\cNodule$, and~$g$ should be bonsai.

\begin{remark}
  A module-map primitive has a first term
  $\TrPMorDiag_1=\pcorolla{1}\otimes\stump$, giving the special case
  above, while a module diagonal primitive does not have a
  corresponding term, instead starting with
  $\TrPMDiag_2=\corolla{2}\otimes \IdTree$. As a result there is an
  extra term in Figure~\eqref{fig:one-side-DT} compared to
  Equation~\eqref{eq:DD-id-tens-map}.
\end{remark}

\begin{lemma}\label{lem:DT-Id-dg-func}
  Suppose that either $P$ and $Q$ are bounded or $\Alg$, $\Blg$, $\cModule$, $\cNodule$, and $g$ are bonsai. Assume also that $\Blg$ is strictly unital. Then
  the maps $\Id_{\cModule}\DT^{\TrPMDiag}\cdot$ and $\cdot \DT^{\TrPMorDiag}\Id_P$ are chain maps
  \begin{align*}
    \Id_{\cModule}\DT^{\TrPMDiag}\cdot\co \Mor(P,Q) &\to \Mor(\cModule\DT^{\TrPMDiag}{P}, \cModule\DT^{\TrPMDiag}{Q})\\
    \cdot \DT^{\TrPMorDiag}\Id_P\co \Mor(\cModule,\cNodule)&\to \Mor(\cModule\DT^{\TrPMDiag_1}P,\cNodule\DT^{\TrPMDiag_2}P).                                                          
  \end{align*}
  Moreover, for $\cdot \DT^{\TrPMorDiag}\Id_P$, homotopic module-map
  primitives give homotopic chain maps.
\end{lemma}
\begin{proof}
  For the first statement,
  \begin{align*}
    d(\Id\DT f)&=
  \tikzsetnextfilename{lem-DT-Id-dg-func-1}
                 \mathcenter{
}
      \\
               &=\Id\DT (d(f)).
  \end{align*}

  Here, the second equality uses the module diagonal primitive relation.
  The third uses the type \DD\ structure relation and an analogue of
  Lemma~\ref{lem:F-f}.
  
  Turning to $\cdot \DT^{\TrPMorDiag}\Id_P$, we have
  \begin{align*}
    d(g \DT^{\TrPMorDiag}\Id_P)&=
    \mathcenter{
  \tikzsetnextfilename{lem-DT-Id-dg-func-15}
    %
}.
    \end{align*}
    Here, the second equality uses the structure relation for module-map primitives, the third uses Lemma~\ref{lem:F-f}, and the fourth uses the structure relation for a type \DD\ structure. This proves that $\cdot \DT^{\TrPMorDiag}\Id_P$ is a chain map. The proof that homotopic module-map primitives give homotopic chain maps is similar, and is left to the reader.
\end{proof}

\begin{warning}\label{warn:tensor-id}
  In the expression $g \DT^{\TrPMorDiag}\Id_P$, $\Id_P$ is just
  notation. Indeed, to define the identity map of $P$ requires that
  $\Alg\ADtp\Blg$ be strictly unital (Definition~\ref{def:strict-unital}), which
  is not guaranteed, even if $\Alg$ and~$\Blg$ are.
\end{warning}

For the sake of simplicity, we only give a partial analogue of Lemma~\ref{lem:DT-mod-defined}:
\begin{proposition}\label{prop:DT-funct}
  Fix:
  \begin{itemize}
  \item An associahedron tree diagonal $\TrDiag$.
  \item Module diagonal primitives $\TrPMDiag^1$, $\TrPMDiag^2$, and
    $\TrPMDiag^3$ compatible with $\TrDiag$.
  \item Module-map primitives $\TrPMorDiag^{12}$ compatible with
    $\TrPMDiag^1$ and $\TrPMDiag^2$; $\TrPMorDiag^{23}$ compatible
    with $\TrPMDiag^2$ and $\TrPMDiag^3$; and $\TrPMorDiag^{13}$ compatible
    with $\TrPMDiag^1$ and $\TrPMDiag^3$.
  \item $\Ainf$-algebras $\Alg$ and $\Blg$ with $\Blg$ strictly unital.
  \item $\Ainf$-modules $\cModule_1$, $\cModule_2$, and $\cModule_3$ over $\Alg$.
  \item Type \DD\ structures $P_1$, $P_2$, and $P_3$ over $\Alg$ and $\Blg$.
  \end{itemize}
  Then the following diagrams commute up to homotopy:
  \begin{equation}
    \label{eq:mm-com-1}
    \tikzsetnextfilename{ugly1}
    \begin{tikzpicture}[smallpic]
      \node at (0,0) (tlcorn) {$\Mor_{\Alg}(\cModule_1,\cModule_2) \rotimes{\Ring}
        \Mor_{\Alg}(\cModule_2,\cModule_3)$};
      \node at (13,0) (trcorn) {${\begin{array}{l}
                                    \Mor^{\Blg}(\cModule_1\DT^{\TrPMDiag^1} P_1,\cModule_2 \DT^{\TrPMDiag^2} P_1) \\
                                  \qquad\rotimes{\Ring} \Mor^{\Blg}(\cModule_2\DT^{\TrPMDiag^2} P_1,\cModule_3\DT^{\TrPMDiag^3} P_1)
                                \end{array}}$};
                              \node at (0,-2) (blcorn) {$\Mor_{\Alg}(\cModule_1,\cModule_3)$};
                              \node at (13,-2) (brcorn) {$\Mor^{\Blg}(\cModule_1\DT^{\TrPMDiag^1} P_1, \cModule_3\DT^{\TrPMDiag^3} P_1);$};
                              \draw[->] (tlcorn) to node[above]{\lab{f\otimes g \mapsto (f\DT^{\TrPMDiag^{12}}\Id_{P_1})\otimes(g\DT^{\TrPMorDiag^{23}} \Id_{P_1})}} (trcorn);
                              \draw[->] (tlcorn) to node[left]{\lab{(f\otimes g)\mapsto
                                  g\circ f}} (blcorn);
                              \draw[->] (trcorn) to node[right]{\lab{(k^1\otimes l^1)\mapsto l^1\circ k^1}} (brcorn);
                              \draw[->] (blcorn) to node[below]{\lab{h\mapsto (h\DT^{\TrPMorDiag^{13}} \Id_{P_1})}} (brcorn);
                            \end{tikzpicture}
      \end{equation}
      \begin{equation}
        \label{eq:mm-com-2}
    \tikzsetnextfilename{ugly2}
                  \begin{tikzpicture}[smallpic]
            \node at (0,0) (tlcorn) {$\Mor^{\Alg}(P_1,P_2)\rotimes{\Ring} \Mor^{\Alg}(P_2,P_3)$};
            \node at (13,0) (trcorn) {${\begin{array}{l}
              \Mor_{\Ground}(M\DT^{\TrPMDiag^1} P_1,M \DT^{\TrPMDiag^1} P_2) \\
              \qquad \rotimes{\Ring} \Mor_{\Ground}(M \DT^{\TrPMDiag^1} P_2,M\DT^{\TrPMDiag^1} P_3)
            \end{array}}$};
            \node at (0,-2) (blcorn) {$\Mor^{\Alg}(P_1,P_3)$};
            \node at (13,-2) (brcorn) {$\Mor_{\Ground}(M\DT^{\TrPMDiag^1} P_1, M\DT^{\TrPMDiag^1} P_3).$};
            \draw[->] (tlcorn) to node[above]{\lab{f^1\otimes g^1 \mapsto (\Id_M\DT^{\TrPMDiag^1} f^1)\otimes(\Id_M\DT^{\TrPMDiag^1}
            g^1)}} (trcorn);
            \draw[->] (tlcorn) to node[left]{\lab{(f^1\otimes g^1)\mapsto g^1\circ f^1}} (blcorn);
            \draw[->] (trcorn) to node[right]{\lab{(k^1\otimes l^1)\mapsto
                l\circ k}} (brcorn);
            \draw[->] (blcorn) to node[below]{\lab{h^1\mapsto (\Id_M\DT^{\TrPMDiag^1} h^1)}} (brcorn);
          \end{tikzpicture}
      \end{equation}
      \begin{equation}
        \label{eq:mm-com-3}
    \tikzsetnextfilename{ugly3}
        \begin{tikzpicture}[smallpic]
            \node at (0,0) (tlcorn) {$\Mor_\Alg(\cModule_1,\cModule_2)\rotimes{\Ring}
          \Mor^{\Alg}(P_1,P_2)$};
            \node at (13,0) (trcorn) {${\begin{array}{l}\Mor_\Ground(\cModule_1\DT^{\TrPMDiag^1} P_1, \cModule_1\DT^{\TrPMDiag^1} P_2)\\ \qquad\rotimes{\Ring} \Mor_\Ground(\cModule_1\DT^{\TrPMDiag^1} P_2, \cModule_2\DT^{\TrPMDiag^2} P_2)\end{array}}$};
            \node at (0,-3) (blcorn) {$\begin{array}{l}\Mor_\Ground(\cModule_1\DT^{\TrPMDiag^1} P_1, \cModule_2\DT^{\TrPMDiag^2} P_1) \\
              \qquad \rotimes{\Ring} \Mor_\Ground(\cModule_2\DT^{\TrPMDiag^2} P_1, \cModule_2\DT^{\TrPMDiag^2} P_2)
            \end{array}$};
            \node at (13,-3) (brcorn) {$\Mor_\Ground(\cModule_1\DT^{\TrPMDiag^1} P_1, \cModule_2\DT^{\TrPMDiag^2} P_2)$};
            \draw[->] (tlcorn) to node[above]{\lab{f \otimes g^1\mapsto (\Id_{\cModule_1}\DT^{\TrPMDiag^1}
            g^1)\otimes (f\DT^{\TrPMorDiag^{12}} \Id_{P_2})}} (trcorn);
            \draw[->] (tlcorn) to node[left]{
              \lab{\substack{(f\otimes g^1)\mapsto(f\DT^{\TrPMorDiag^{12}} \Id_{P_1})\hfill\\
                  \qquad\otimes (\Id_{\cModule_2}\DT^{\TrPMDiag^2} g^1)}}} (blcorn);
            \draw[->] (trcorn) to node[right]{\lab{(k^1\otimes
                l^1)\mapsto l^1\circ k^1}} (brcorn);
            \draw[->] (blcorn) to node[below]{\lab{(k^1\otimes
                l^1)\mapsto l^1\circ k^1}} (brcorn);
          \end{tikzpicture}
      \end{equation}
      (Here, we either assume that the $P_i$ are all bounded, or else that
      $\Alg$, $\Blg$, and the $\cModule_i$ are bonsai and that the morphism complexes
      denote the complexes of bonsai morphisms.)
\end{proposition}
\begin{proof}
  To construct the homotopy for the square~\eqref{eq:mm-com-1} one uses the
  notion of a module bi-map diagonal primitive (the definition of
  which is left to the reader). The square~\eqref{eq:mm-com-2}
  follows from the definitions and the module primitive structure equation. The
  homotopy for square~\eqref{eq:mm-com-2} is the map
  \[
  \tikzsetnextfilename{prop-DT-funct-1}
  \begin{tikzpicture}[smallpic]
    \node at (0,0) (tc) {};
    \node at (-1,0) (tl) {};
    \node at (0, -1) (delta1) {$\delta$};
    \node at (0, -2) (f) {$g$};
    \node at (0, -3) (delta2) {$\delta$};
    \node at (-1, -4) (pl) {$\TrPMorDiag^{12}f$};
    \node at (1, -4) (pr) {$\TrPMorDiag^{12,\op}$};
    \node at (-1, -5) (bl) {};
    \node at (0, -5) (bc) {};
    \node at (1, -5) (br) {};
    \draw[dmoda] (tc) to (delta1);
    \draw[dmoda] (delta1) to (f);
    \draw[dmoda] (f) to (delta2);
    \draw[dmoda] (delta2) to (bc);
    \draw[moda] (tl) to (pl);
    \draw[moda] (pl) to (bl);
    \draw[taa] (delta1) to (pl);
    \draw[taa] (delta1) to (pr);
    \draw[alga] (f) to (pl);
    \draw[blga] (f) to (pr);
    \draw[tbb] (delta2) to (pl);
    \draw[tbb] (delta2) to (pr);
    \draw[alga] (pr) to (br);
  \end{tikzpicture}.    
  \]
  Further details are left to the reader.
\end{proof}

Recall that for $\Blg$ a strictly unital $\Ainf$-algebra and $\lsup{\Blg}Q$ a type $D$
structure over $\Blg$, the identity map of $Q$ is defined by
\[
  \Id_Q^1(x)=\unit\otimes x.
\]
This identity map gives rise to the notion of isomorphic or homotopy equivalent
type $D$ structures.  (See also Warnings~\ref{warn:non-unital-cat} and~\ref{warn:tensor-id}.)
\begin{lemma}\label{lem:Id-DT-Id}
  Fix module diagonal primitives $\TrPMDiag_1$ and $\TrPMDiag_2$ and a
  module-map primitive $\TrPMorDiag$ compatible with $\TrPMDiag_1$ and
  $\TrPMDiag_2$.
  Given $\Ainf$-algebras $\Alg$ and $\Blg$ with $\Blg$ strictly unital, an
  $\Ainf$-module $\cModule_{\Alg}$ and a type \DD\ structure
  $\lsup{\Alg,\Blg}P$, with $P$ bounded or $\cA$, $\cB$, and $\cM$ bonsai, the map
  \[
    \Id_{\cModule} \DT^{\TrPMorDiag}\Id_P\co \cModule\DT^{\TrPMDiag_1}P\to \cModule\DT^{\TrPMDiag_2}P
  \]
  is a homotopy equivalence.
\end{lemma}
\begin{proof}
  First, consider the special case that $\TrPMDiag_1=\TrPMDiag_2$ and
  $\TrPMorDiag$ is the module-map primitive from
  Lemma~\ref{lem:mor-prim-example}. The term
  \[
  \tikzsetnextfilename{lem-Id-DT-Id-1}
    \begin{tikzpicture}[smallpic]
      \node at (0,0) (tl) {};
      \node at (0,-1) (Id) {$\Id$};
      \node at (0,-2) (bl) {};
      \node at (1,0) (tc) {};
      \node at (1,-2) (bc) {};
      \node at (2,-1) (one) {$\unit$};
      \node at (2,-2) (br) {};
      \draw[moda] (tl) to (Id);
      \draw[moda] (Id) to (bl);
      \draw[dmoda] (tc) to (bc);
      \draw[alga] (one) to (br);
    \end{tikzpicture}
  \]
  contributes the identity map to
  $\Id_{\cModule} \DT^{\TrPMorDiag}\Id_P$, so it suffices to show that
  no other terms contribute. This is immediate from the facts that
  $(\Id_{\cModule})_n=0$ for $n>1$ and the distinguished vertex in
  $\TrPMorDiag_n$ has valence $\geq 3$ if $n>1$.

  Now, for any other module-map primitive $\TrPMorDiag'$ compatible
  with $\TrPMDiag_1=\TrPMDiag_2$, the result follows from the fact
  that $\TrPMorDiag$ and $\TrPMorDiag'$ are homotopic
  (Proposition~\ref{prop:mod-map-prim-homotopic}) and
  Lemma~\ref{lem:DT-Id-dg-func}. Finally, the case of general module
  diagonal primitives follows from Equation~(\ref{eq:mm-com-1}) in Proposition~\ref{prop:DT-funct}, with $f=g=\Id_\cModule$.
\end{proof}

\begin{corollary}\label{cor:Id-DT-Id}
  Fix a module diagonal primitive $\TrPMDiag$ and a
  module-map primitive $\TrPMorDiag$ compatible with $\TrPMDiag$ and
  $\TrPMDiag$.
  Given $\Ainf$-algebras $\Alg$ and $\Blg$ with $\Blg$ strictly unital, an
  $\Ainf$-module $\cModule_{\Alg}$ and a type \DD\ structure
  $\lsup{\Alg,\Blg}P$, the map
  \[
    \Id_{\cModule} \DT^{\TrPMorDiag}\Id_P\co \cModule\DT^{\TrPMDiag}P\to \cModule\DT^{\TrPMDiag}P
  \]
  is homotopic to the identity map.
\end{corollary}
\begin{proof}
  By Equation~(\ref{eq:mm-com-1}) in Proposition~\ref{prop:DT-funct},
  \[
    (\Id_{\cModule} \DT^{\TrPMorDiag}\Id_P)\circ (\Id_{\cModule} \DT^{\TrPMorDiag}\Id_P)\sim (\Id_{\cModule}\circ \Id_{\cModule})\DT^{\TrPMorDiag}\Id_P=\Id_{\cModule} \DT^{\TrPMorDiag}\Id_P.
  \]
  Composing with the homotopy inverse of $\Id_{\cModule} \DT^{\TrPMorDiag}\Id_P$
  (from Lemma~\ref{lem:Id-DT-Id}) gives the result.
\end{proof}

\begin{corollary}\label{cor:DT-preserve-hequiv}
  Fix $\Ainf$-modules $\cModule$ and $\cNodule$, a homotopy equivalence $f\co\cModule\to\cNodule$, and a type \DD\ structure $P$ such that either $P$ is bounded or $\cModule$, $\cNodule$, and $f$ are bonsai.
  Then for any choice of module diagonal primitives
  and module-map primitive,  
  \[
    f \DT^{\TrPMorDiag}\Id_P
    \co \cModule\DT^{\TrPMDiag_1}P\to \cNodule\DT^{\TrPMDiag_2}P
  \]
  is a homotopy equivalence.
\end{corollary}
\begin{proof}
  This is immediate from Lemma~\ref{lem:Id-DT-Id} and
  Proposition~\ref{prop:DT-funct}.
\end{proof}

\subsubsection{Associativity of the box tensor product of morphisms}
Given a type \DD\ structure $P$, modules $\cModule_1$, $\cModule_2$,
$\cNodule_1$, and $\cNodule_2$, and morphisms $f\co
\cModule_1\to\cModule_2$ and $g\co\cNodule_1\to\cNodule_2$, as well as
an associahedron diagonal $\AsDiag$, module diagonals $\MDiag^1$ and
$\MDiag^2$ compatible with $\AsDiag$, and a module-map diagonal
$\ModMulDiag$ compatible with $\MDiag^1$ and $\MDiag^2$,
Definition~\ref{def:tensor-mod-maps} gives a triple box product
\[
  [f\DT\Id_P\DT g]_{\ModMulDiag}=(f\MDtp[\ModMulDiag]g)\DT \Id_P
\]
(under appropriate boundedness assumptions). This is a little trickier
than it looks, since $\Alg\ADtp\Blg$ is a non-unital algebra, so the
category of type $D$ structures over it is a non-unital
category. Nonetheless, for $f\co \cModule\to \cNodule$ a map of modules
over a non-unital algebra $\Clg$ and $\lsup{\Clg}P$ a type $D$
structure over $\Clg$, it still makes sense to define
\tikzsetnextfilename{fDTidNonunital}
\[
  f\DT\Id_P=
  \mathcenter{\begin{tikzpicture}[smallpic]
    \node at (0,0) (tl) {};
    \node at (2,0) (tr) {};
    \node at (2,-1) (delta) {$\delta_P$};
    \node at (0,-2) (m) {$f$};
    \node at (2,-3) (br) {};
    \node at (0,-3) (bl) {};
    \draw[moda] (tl) to (m);
    \draw[moda] (m) to (bl);
    \draw[dmoda] (tr) to (delta);
    \draw[dmoda] (delta) to (br);
    \draw[taa] (delta) to (m);
  \end{tikzpicture}}.
\]
In the case of interest, $\Clg=\Alg\ADtp\Blg$ and this reduces to 
\tikzsetnextfilename{fDTidDD}
\[
  [f\DT\Id_P\DT g]_{\ModMulDiag}=
    \mathcenter{\begin{tikzpicture}[smallpic]
    \node at (0,0) (tc) {};
    \node at (-1,0) (tl) {};
    \node at (1,0) (tr) {};
    \node at (0,-1) (delta) {$\delta$};
    \node at (-1,-2) (pl) {$\TrModMulDiag f$};
    \node at (1,-2) (pr) {$\TrModMulDiag^\op g$};
    \node at (-1,-3) (bl) {};
    \node at (0,-3) (bc) {};
    \node at (1,-3) (br) {};
    \draw[dmoda] (tc) to (delta);
    \draw[dmoda] (delta) to (bc);
    \draw[moda] (tl) to (pl);
    \draw[moda] (pl) to (bl);
    \draw[moda] (tr) to (pr);
    \draw[taa] (delta) to (pl);
    \draw[taa] (delta) to (pr);
    \draw[moda] (pr) to (br);
  \end{tikzpicture}},
\]
where $\TrModMulDiag f$ means we apply $f$ at the distinguished vertex
on the left tree in $\TrModMulDiag$ and $\TrModMulDiag^\op g$ means we apply $g$ at the
distinguished vertex of the right tree in $\TrModMulDiag$.

In the special case that $g$ is the identity map, $[f\DT \Id_P\DT g]_{\ModMulDiag}$ depends only on the partial module-map diagonal induced by $\ModMulDiag$. That is, given a partial module-map diagonal $\PartTrModMulDiag$ compatible with $\MDiag^1$ and $\MDiag^2$ we can define
\begin{equation}
  \label{eq:part-mor-tens}
\mathcenter{  [f\DT\Id_P\DT \Id_{\cNodule}]_{\PartTrModMulDiag}=}
\tikzsetnextfilename{eq-part-mor-tens-1}
  \mathcenter{\begin{tikzpicture}[smallpic]
    \node at (0,0) (tc) {};
    \node at (-1,0) (tl) {};
    \node at (1,0) (tr) {};
    \node at (0,-1) (delta) {$\delta$};
    \node at (-1,-2) (pl) {$\PartTrModMulDiag f$};
    \node at (1,-2) (pr) {$\PartTrModMulDiag^\op\Id$};
    \node at (-1,-3) (bl) {};
    \node at (0,-3) (bc) {};
    \node at (1,-3) (br) {};
    \draw[dmoda] (tc) to (delta);
    \draw[dmoda] (delta) to (bc);
    \draw[moda] (tl) to (pl);
    \draw[moda] (pl) to (bl);
    \draw[moda] (tr) to (pr);
    \draw[taa] (delta) to (pl);
    \draw[taa] (delta) to (pr);
    \draw[moda] (pr) to (br);
  \end{tikzpicture}}
\end{equation}
where $\PartTrModMulDiag f$ means that we apply $f$ at the distinguished vertex of (the left tree in) $\PartTrModMulDiag$ and $\PartTrModMulDiag \Id$ means we apply $\Id$ at the distinguished vertex of (the right tree in) $\PartTrModMulDiag$.
\begin{lemma}\label{lem:part-mod-map-diag-DT}
  Fix module diagonals $\MDiag^1$ and $\MDiag^2$, a partial
  module-map diagonal $\PartTrModMulDiag$, $\Ainf$-algebras $\Alg$ and
  $\Blg$, a type \DD\ structure $P$ over $\Alg$ and $\Blg$,
  $\Ainf$-modules $\cModule_1$ and $\cModule_2$ over $\Alg$, and an $\Ainf$-module
  $\cNodule$ over $\Blg$.  
  Assume that either $P$ is bounded or $\cModule_1$, $\cModule_2$, 
  $\cNodule$, and the morphisms $f\co \cModule_1\to\cModule_2$ under consideration are bonsai. Then
  Formula~\eqref{eq:part-mor-tens} defines a chain map
  \[
    \Mor(\cModule_1,\cModule_2)\to \Mor([\cModule_1\DT P\DT \cNodule]_{\MDiag^1},[\cModule_2\DT P\DT\cNodule]_{\MDiag^2}).
  \]
  Moreover, homotopic partial module-map diagonals induce homotopic chain maps. Finally, if $\PartTrModMulDiag$ is induced by a module-map primitive $\TrPMorDiag$ then $[f\DT\Id_P\DT \Id_{\cNodule}]_{\PartTrModMulDiag}$ is homotopic to $(f\DT^{\TrPMorDiag}\Id_P)\DT\Id_{\cNodule}$.
\end{lemma}
\begin{proof}
  For the first statement, $d[f\DT\Id_P\DT
  \Id_{\cNodule}]_{\PartTrModMulDiag}$ is given by
  \[
\tikzsetnextfilename{lem-part-mod-map-diag-DT-1}
  \mathcenter{
}.
  \]
  From the definition of a partial module-map diagonal, this sum is equal to
  \begin{multline*}
\tikzsetnextfilename{lem-part-mod-map-diag-DT-7}
    \mathcenter{%
}
  \end{multline*}
  where $\alphas,\betas\in \ModTransTrees{}\rotimes{\Ring} \Filt_{\geq 3}\ModTransTrees{}$.

  Applying Lemma~\ref{lem:F-f}, this sum is equal to
  \begin{multline*}
\tikzsetnextfilename{lem-part-mod-map-diag-DT-16}
    \mathcenter{%
}.
  \end{multline*}
  The second term vanishes because $d(\Id)=0$. The third, fourth, and fifth
  terms cancel by the type \DD\ structure relation. The last term vanishes
  because the identity map $\Id$ has $\Id_n=0$ for $n>1$. So, applying
  Lemma~\ref{lem:F-f} to the second-to-last term we are left with
  \[
\tikzsetnextfilename{lem-part-mod-map-diag-DT-23}
    \mathcenter{%
}.
  \]
  Again, the second term vanishes because $d(\Id)=0$ and the third,
  fourth, and fifth terms vanish
  because $\Id_n=0$ for $n\geq 2$. Thus, we are left with $[(df)\DT\Id_P\DT \Id_{\cNodule}]_{\PartTrModMulDiag}$, as desired.

  For the second statement, applying the argument above but with a homotopy $\zeta$ of partial module-map diagonals in place of $\PartTrModMulDiag$ gives
  \[
\tikzsetnextfilename{lem-part-mod-map-diag-DT-26}
    d\left(  \mathcenter{%
}
    \right)
    +[f\DT\Id_P\DT \Id_{\cNodule}]_{\PartTrModMulDiag^1}-[f\DT\Id_P\DT \Id_{\cNodule}]_{\PartTrModMulDiag^2},
  \]
  as desired.

  Finally, if $\PartTrModMulDiag$ comes from a primitive $\TrPMorDiag$ (and
  $\TrPMDiag^i$ comes from a primitive $\TrPMDiag^i$) then
  \[
    [f\DT\Id_P\DT \Id_{\cNodule}]_{\PartTrModMulDiag}=
    \mathcenter{
\tikzsetnextfilename{lem-part-mod-map-diag-DT-28}
      %
    }=(f\DT_{\TrPMorDiag}\Id_P)\DT \Id_N,
  \]
  as claimed.

  The boundedness assumptions imply that all the sums considered above
  are finite.
\end{proof}

\begin{corollary}
  Fix a module-map primitive $\TrPMorDiag$ compatible with module
  diagonal primitives $\TrPMDiag^1$ and~$\TrPMDiag^2$. Let
  $\TrMDiag^i$ be the module diagonal induced by $\TrPMDiag^i$ and let
  $\TrModMulDiag$ be any module-map diagonal compatible with
  $\TrMDiag^1$ and $\TrMDiag^2$. Then for any $\Ainf$-modules
  $\cModule_1$, $\cModule_2$, and $\cNodule$, type \DD\ structure $P$,
  and $\Ainf$-module homomorphism $f\co \cModule_1\to\cModule_2$ (with
  either $P$ bounded or all the other objects bonsai),
  the chain maps
  \[
    [f\DT \Id_P\DT \Id_{\cNodule}]_{\TrModMulDiag}
  \]
  and
  \[
    (f\DT_{\TrPMorDiag}\Id_P)\DT \Id_{\cNodule}
  \]
  are chain homotopic.
\end{corollary}
\begin{proof}
  This is immediate from Lemma~\ref{lem:part-mod-map-diag-DT} and the fact that all partial module-map diagonals are homotopic (Lemma~\ref{lem:part-mod-maps-htpic}).
\end{proof}

\subsection{Further tensor products of modules and bimodules}
\label{sec:box2}

We collect here a few other tensor products of bimodules, generalizing
the earlier constructions.  We subdivide these into two types:
internal tensor products and external tensor products.

Examples of internal tensor products that we have already
defined are:
\begin{itemize}
\item The box product of an $\Ainf$-module and a type $D$ structure,
  as recalled in Section~\ref{sec:box}.
\item Given \DA\ structures $\lsup{\Alg}P_{\Blg}$ and
  $\lsup{\Blg}Q_{\Clg}$ (see Definition~\ref{def:DA-str}
  or~\cite[Definition 2.2.43]{LOT2}), the \DA\ structure
  $\lsup{\Alg}P_{\Blg}\DT \lsup{\Blg}Q_{\Clg}$~\cite[Figure~4]{LOT2}
  (see also Figure~\ref{fig:AADADADA}).
\item Given an $\Ainf$-bimodule $\lsub{\Alg}\cModule_{\Blg}$ and a \DA\
  structure $\lsup{\Blg}Q_{\Clg}$, the $\Ainf$-bimodule
  $\lsub{\Alg}\cModule_{\Blg}\DT \lsup{\Blg}Q_{\Clg}$~\cite[Figure~4]{LOT2}
  (see also Figure~\ref{fig:AADADADA}).
\item The triple box product of a type \DD\ structure and two
  $\Ainf$-modules (Section~\ref{sec:box}), using a module diagonal.
\item The one-sided box tensor product of a type \DD\ structure and an
  $\Ainf$-module (Section~\ref{sec:box}), using a module diagonal
  primitive.
\end{itemize}

\begin{figure}
  \centering
  \[
\tikzsetnextfilename{fig-AADADADA-1}
    \mathcenter{\begin{tikzpicture}[smallpic]
        \node at (0,0) (tl) {};
        \node at (2,0) (tr) {};
        \node at (0,-2) (m) {$\delta^1_P$};
        \node at (2,-1) (delta) {$\delta_Q$};
        \node at (2,-3) (br) {};
        \node at (0,-3) (bl) {};
        \node at (4,0) (trr) {};
        \node at (-2,-3) (bll) {};
        \draw[moda] (tl) to (m);
        \draw[moda] (m) to (bl);
        \draw[dmoda] (tr) to (delta);
        \draw[dmoda] (delta) to (br);
        \draw[taa] (delta) to (m);
        \draw[taa] (trr) to (delta);
        \draw[alga] (m) to (bll);
      \end{tikzpicture}}
    \qquad\qquad\qquad
\tikzsetnextfilename{fig-AADADADA-2}
    \mathcenter{\begin{tikzpicture}[smallpic]
        \node at (-2,0) (tll) {};
        \node at (0,0) (tl) {};
        \node at (2,0) (tr) {};
        \node at (0,-2) (m) {$m_M$};
        \node at (2,-1) (delta) {$\delta_Q$};
        \node at (2,-3) (br) {};
        \node at (0,-3) (bl) {};
        \node at (4,0) (trr) {};
        \draw[moda] (tl) to (m);
        \draw[moda] (m) to (bl);
        \draw[dmoda] (tr) to (delta);
        \draw[dmoda] (delta) to (br);
        \draw[taa] (delta) to (m);
        \draw[taa] (trr) to (delta);
        \draw[taa] (tll) to (m);
      \end{tikzpicture}}.
  \]
  \caption[Box products of two \DA\ structure, and of an
      A-infinity bimodule and a \DA\ structure.]{\textbf{Box products of two \DA\ structure, and of an
      $\Ainf$-bimodule and a \DA\ structure.} Left: the structure map for
    a box tensor product $\lsup{\Alg}P_{\Blg}\DT\lsup{\Blg}Q_\Clg$ of
    two \DA\ structure. Right: the structure map of a box tensor
    product $\lsub{\Alg}\cModule_{\Blg}\DT\lsup{\Blg}Q_\Clg$ of an
    $\Ainf$-bimodule and a \DA\ structure.  In each case, boundedness
    criteria must be satisfied for the definition to make
    sense. See~\cite{LOT2} for more details.}
  \label{fig:AADADADA}
\end{figure}

Given $\Ainf$-modules $\cModule_\Alg$ and $\cNodule_\Blg$, the tensor product
$(\cModule\MDtp \cNodule)_{\Alg\ADtp \Blg}$ from
Section~\ref{sec:TensorAinftyModules} is an example of an external
tensor product.  An even simpler example of an external tensor
product occurs if $\Alg$ and $\Blg$ are $\Ainf$-algebras over $\Ground$
with $\Blg$ strictly unital, and $\lsup{\Blg}P$ and $\cModule_{\Alg}$ are a
type $D$ structure and an $\Ainf$-module, respectively. In this
case, $P\kotimes{\Ground} M$ can be endowed with the structure of a
type \DA\ structure (in the sense recalled below), with operations
\begin{align*}
  \delta^1_1(x\otimes y)&=\One\otimes x\otimes m_1(y) 
                          + \delta^1(x)\otimes y \\
  \delta^1_{\ell+1}(x\otimes y,a_1,\dots,a_{\ell})
                        &=\One\otimes x\otimes m_{\ell+1}(y,a_1,\dots,a_{\ell})).
\end{align*}

In the next sections we handle a few more cases.

\subsubsection{Box product of an \texorpdfstring{$\Ainf$-}{A-infinity-}bimodule and a type \texorpdfstring{\DD}{DD} structure.}

Fix $\Ainf$-algebras $\Alg$, $\Blg$, and $\Clg$, over $\Ground_1$,
$\Ground_2$, and $\Ground_3$, respectively, with $\Clg$ strictly
unital with unit $\unit$. Fix also an $\Ainf$-bimodule
$\lsub{\Alg}\cModule_{\Blg}$ and a type \DD\ structure
$\lsup{\Blg,\Clg}P$. The goal of this section is to use a primitive
to construct a type \DA\ structure
$\lsub{\Alg}\cModule_{\Blg}\DT\lsup{\Blg}P^{\Clg^\op}$.

We start by recalling the definition of a type \DA\ structure.
\begin{definition}\label{def:DA-str}
Let $\Alg$ and $\Blg$ be algebras over $\Ground_1$ and $\Ground_2$
respectively, and let $Q=\lsub{\Ground_1}Q_{\Ground_2}$ be a
$(\Ground_1,\Ground_2)$-bimodule. Given a collection of bimodule homomorphisms
\[
  \{\delta^1_{1+n}\co Q\kotimes{\Ground_2} B^{\kotimes{\Ground_2} n}\to A\kotimes{\Ground_1} Q\grs{1-n}\}_{n=1}^{\infty},
\]
we say that
$(Q,\{\delta^1_{1+n}\})$ is a \emph{type \DA\ structure} if the
following condition holds. Define
$\delta^1=\sum_j\delta^1_j\co Q\kotimes{\Ground_2} \Tensor^*(B\grs{1})\to A\grs{1}\kotimes{\Ground_1}
Q$, and inductively define
\[
  \delta^n=(\Id_{A^{\kotimes{\Ground_1} (n-1)}}\otimes
  \delta^1)\circ(\delta^{n-1}\otimes\Id_{\Tensor^*B})\circ(\Id_Q\otimes \Delta)
\]
where $\Delta\co \Tensor^*B\to \Tensor^*B\otimes \Tensor^*B$ is the comultiplication.
Let
\[
  \delta=\sum_n\delta^n\co Q\kotimes{\Ground_2} \Tensor^*(B\grs{1})\to
  \overline{\Tensor}^*(A\grs{1})\kotimes{\Ground_1} Q.
\]
Then the structure equation is:
\[
  \tikzsetnextfilename{DAstreq1}
  \mathcenter{\begin{tikzpicture}[smallpic]
    \node at (0,0) (tc) {};
    \node at (0,-1) (delta) {$\delta$};
    \node at (0,-3) (bc) {};
    \node at (-1,-2) (mu) {$\mu^\Alg$};
    \node at (-1,-3) (bl) {};
    \node at (1,0) (tr) {};
    \draw[damoda] (tc) to (delta);
    \draw[dmoda] (delta) to (bc);
    \draw[taa] (tr) to (delta);
    \draw[taa] (delta) to (mu);
    \draw[alga] (mu) to (bl);
  \end{tikzpicture}}
  +
  \tikzsetnextfilename{DAstreq2}
  \mathcenter{
  \begin{tikzpicture}[smallpic]
    \node at (0,0) (tc) {};
    \node at (0,-2) (delta) {$\delta$};
    \node at (0,-3) (bc) {};
    \node at (.75,-1) (mu) {$\mu^\Blg$};
    \node at (-1,-3) (bl) {};
    \node at (.75,0) (tr1) {};
    \node at (1.25,0) (tr2) {};
    \node at (2.25,0) (tr3) {};
    \draw[damoda] (tc) to (delta);
    \draw[dmoda] (delta) to (bc);
    \draw[alga] (delta) to (bl);
    \draw[alga] (mu) to (delta);
    \draw[taa] (tr1) to (delta);
    \draw[taa] (tr2) to (mu);
    \draw[taa] (tr3) to (delta);
  \end{tikzpicture}}
  =0.
\]
This makes sense if either $\Alg$ is bonsai, or for each $k$ there is
an $N$ so that $\delta^n_k=0$ for all $n>N$.
(Here, $\delta^n_k$ is the component of $\delta^n$ taking $k$ total inputs.)
When talking about a type \DA\ structure, we always assume one of
these conditions holds.  If the second condition holds, we say that
$P$ is \emph{left bounded}. We say that $P$ is \emph{bounded} if there
is an $N$ so that $\delta^n_k=0$ whenever $n+k>N$.
\end{definition}
Note that left boundedness is automatic if $\delta^1_1=0$ (cf. Section~\ref{sec:DA-tens-DD}).

To define $\lsub{\Alg}\cModule_{\Blg}\DT\lsup{\Blg}P^{\Clg^\op}$ we need one
further kind of tree. A \emph{bimodule tree} is a pair of an element
$T\in \Trees_n$ and an integer $1\leq i\leq n$. The integer $i$ corresponds to a
distinguished input of $T$, which we think of as the module input. Call vertices
between the $i\th$ input of $T$ and the output of $T$ \emph{bimodule vertices},
and call vertices to the left (respectively right) of the bimodule vertices
\emph{left algebra} (respectively \emph{right algebra}) vertices. Given an
$\Ainf$-bimodule $\lsub{\Alg}\cModule_{\Blg}$ and a bimodule tree $(T,i)$ there
is an induced map
$m(T,i)\co A^{\otimes i-1}\otimes M\otimes B^{\otimes n-i}\to M$, by applying
the operation $m$ on $M$ at the bimodule vertices of $T$ and the operation
$\mu^\Alg$ (respectively $\mu^\Blg$) at the left (respectively right) algebra
vertices.  We will typically suppress the integer $i$ and refer to a bimodule
tree $T$.

The operation $\LeftJoin$ of left joining extends to an operation on bimodule trees by composing
at the distinguished inputs. The original case of $\LeftJoin$ is the case that
all of the distinguished inputs were left-most.

Let $S$ be a tree with $1+\ell$ inputs. Given an integer $k\geq 0$, a {\em
  $k$-left enlargement} of $S$ is a tree $S'$ with $k+1+\ell$ inputs
and the same number of internal vertices, with the property that if we
remove the edges adjacent to the first $k$ inputs of $S'$ the
resulting tree is $S$.  We view $S'$ as a bimodule tree where the $(k+1)\st$ vertex is distinguished.
There is a corresponding map
\[
  \lsub{k}E\co \cellC{*}(K_{\ell+1})\to\cellC{*+k}(K_{k+\ell+1})
\]
that sends a tree $S$ to the sum of all of its $k$-enlargements.
The operation $\lsub{k}E$ satisfies
\begin{align}
  \bdy(\lsub{k}E(T))+\lsub{k}E(\bdy(T))
  &=\sum_{\ell=0}^{k-1}\lsub{\ell}E(T)\circ_{\ell+1}\corolla{k-\ell+1}
    +\sum_{\ell=0}^{k-1}\corolla{k-\ell+1}\circ_{k-\ell+1}\lsub{\ell}E(T)
    +\sum_{\ell=1}^{k-1}\sum_{j=1}^\ell \lsub{\ell}E(T)\circ_j \corolla{k-\ell+1}
    \label{eq:left-enlarge-1}
  \\
  \lsub{*}E(\LeftJoin(T_1,\dots,T_n))&=\LeftJoin(\lsub{*}E(T_1),\dots,\lsub{*}E(T_n))
    \label{eq:left-enlarge-2}
\end{align}
where, in Equation~\eqref{eq:left-enlarge-2}, $\lsub{*}E(T)=\sum_k \lsub{k}E(T)$
and $\LeftJoin$ is the usual left joining on the left side of the equation, and
the operation $\LeftJoin$ on bimodule trees on the right side of the
equation. The terms on the right side Equation~\eqref{eq:left-enlarge-1}
correspond to gluing a corolla at the top of $T$, and at the bottom of $T$, and
to the left of $T$, respectively.

Fix a module diagonal primitive
\[ 
  \TrPMDiag_n=\sum_{(S,T)}n_{S,T}(S,T)\in \Trees_n\rotimes{\Ring}
  \Trees_{n-1}.
\]

Generalizing Lemma~\ref{lem:Ainf-mod-is}, the bimodule $\lsub{\Alg}\cModule_\Blg$ gives a
chain map $m\co \cellC{*}(K_{k+1+\ell})\to\Mor(A^{\kotimes{\Ground_1}
  k}\kotimes{\Ground_1} M \kotimes{\Ground_2} B^{\otimes \ell},M)$.
We define operations
\[
  \delta^1\co\Tensor^*(A)\kotimes{\Ground_1} M\kotimes{\Ground_2}P
  \to M\kotimes{\Ground_2} P\kotimes{\Ground_3}C^\op
\]
by
\[
  \delta^1_{k+1}((a_1\otimes\dots\otimes a_k)\otimes x\otimes y)=m^{\cModule}({\underline a}\otimes x)\otimes y \otimes \unit+
  \sum_{(S,T)}n_{S,T}(m^M(\lsub{k}E(S)) \otimes \Id_{P}\otimes
  \mu^{\Clg^{\op}}(T^\op))\circ (\Id_{A^{\otimes k}}\otimes \Id_M \otimes\delta_P)
\]
or, graphically:
\[ \delta^1({\underline a}\otimes x\otimes y)
  =
  \mathcenter{
\tikzsetnextfilename{sec-AADADADA-1}
    \begin{tikzpicture}[smallpic]
      \node at (-2,0) (tll) {${\underline a}$};
      \node at (-1,0) (tl) {$x$};
      \node at (0,0) (tc) {$y$};
      \node at (-1,-1) (m) {$m_{k,1,0}$};
      \node at (1,-1) (one) {$\unit$};
      \node at (-1,-2) (bl) {};
      \node at (0,-2) (bc) {};
      \node at (1,-2) (br) {};
      \draw[taa] (tll) to (m);
      \draw[dmoda] (tc) to (bc);
      \draw[moda] (tl) to (m);
      \draw[moda] (m) to (bl);
      \draw[alga] (one) to (br);
    \end{tikzpicture}
  }+
  \mathcenter{
\tikzsetnextfilename{sec-AADADADA-2}
    \begin{tikzpicture}[smallpic]
      \node at (-2,0) (a) {${\underline a}$};
      \node at (-1,0) (x) {$x$};
      \node at (0,0) (y) {$y$};
      \node at (0, -1) (delta) {$\delta$};
      \node at (-1, -2) (S) {$\lsub{k}E(\TrPMDiag)$};
      \node at (1, -2) (T) {$\TrPMDiag^\op$};
      \node at (-1, -3) (bl) {};
      \node at (0, -3) (bc) {};
      \node at (1, -3) (br) {};
      \draw[dmoda] (y) to (delta);
      \draw[dmoda] (delta) to (bc);
      \draw[moda] (x) to (S);
      \draw[moda] (S) to (bl);
      \draw[taa] (delta) to (S);
      \draw[taa] (a) to (S);
      \draw[tbb] (delta) to (T);
      \draw[blga] (T) to (br);
    \end{tikzpicture}
  },
\]
where ${\underline a}=a_1\otimes\dots\otimes a_k$.

Let
\[
  \lsub{\Alg} \cModule_{\Blg}\DT_{\TrPMDiag} \lsup{\Blg}P^{\Clg^\op}=(M\kotimes{\Ground_2} P,\delta^1_{1+n}).
\]

\begin{proposition}\label{prop:DT-AA-DD-defined}
  If $\lsub{\Alg}\cModule_\Blg$ is bonsai or $\lsup{\Blg,\Clg}P$ is bounded,
  then the operation $\delta^1_{1+n}$ on
  $\lsub{\Alg} \cModule_{\Blg}\DT_{\TrPMDiag} \lsup{\Blg}P^{\Clg^\op}$ satisfies
  the type \DA\ structure relation.
\end{proposition}

\begin{proof}
  This is a straightforward generalization of the proof of
  Lemma~\ref{lem:DT-defined}, using
  Equations~(\ref{eq:left-enlarge-1}) and~(\ref{eq:left-enlarge-2})
  when needed.
\end{proof}

\subsubsection{Box product of a \texorpdfstring{\DA}{DA}
  structure and a \texorpdfstring{\DD}{DD} structure}\label{sec:DA-tens-DD}

We can use DADD diagonals to define a tensor product of a type \DA\
structure with a type \DD\ structure, under a simplifying
assumption. Specifically, suppose that the \DA\ structure
$\lsup{\Alg}Q_{\Blg}$ satisfies the condition that the operation
$\delta^1_1$ vanishes identically.  (For instance, this assumption
holds for the type \DA\ structure associated to an $\Ainf$-algebra
homomorphism~\cite[Definition 2.2.48]{LOT2} or the type \DA\ structures associated to
arcslides in bordered Floer theory~\cite{LOT4}.) In this case, the operations
$\delta^1_{1+n}$ give rise to an action
\[
  \lambda^Q\co \cellC{*}(J_n)\to \Mor(Q\kotimes{\Ground_2} B^{\kotimes{\Ground_2} n},
  A\kotimes{\Ground_1} Q)
\]
as follows.
Given a multiplihedron tree $S$ representing a generator of
$\cellC{*}(J_n)$, the action of $\lambda^Q(S)$ is obtained by running an
additional strand labelled by $Q$ through all purple vertices, and
thinking of the purple vertices as corresponding to $\delta^1_{1+n}$
operations, red vertices correspond to $\mu^{\Alg}_n$ operations, and
blue vertices correspond to $\mu^{\Blg}_n$ operations.
(Note that this agrees with the standard multiplihedron action in the
case the DA structure is associated to an $\Ainf$-algebra
homomorphism.)
\begin{lemma}\label{lem:DA-lambda-chain-map}
  For each $n$, the map
  $\lambda^Q\co \cellC{*}(J_n)\to \Mor(Q\kotimes{\Ground_2}
  B^{\kotimes{\Ground_2} n}, A\kotimes{\Ground_1} Q)$
  defined above is a chain map.
\end{lemma}
\begin{proof}
  This follows from the \DA\ structure relation for
  $\lsup{\Alg}Q_{\Blg}$ and the fact that the operation $\delta^1_1$
  vanishes identically.
\end{proof}

We can now specify the \DD\ operations on
$\lsup{\Alg}Q_{\Blg}\otimes \lsup{\Blg}P^{\Clg}$, as follows.  Fix a
DADD diagonal
\[
  \TrDADD_n=\sum_{(S,T)}n_{S,T}(S,T)\in \cellC{*}(J_n) \rotimes{\Ring} \Trees_{n}
\]
compatible with associahedron diagonals $\TrDiag^1$ and $\TrDiag^2$
(Definition~\ref{def:DADD-diag}) and a type \DD\ structure
$\lsup{\Blg,\Clg}P$ with respect to $\TrDiag^1$.
Define
\[
  \delta^1\co Q\kotimes{\Ground_2}P\to A\kotimes{\Ground_1}
  Q\kotimes{\Ground_2} P\kotimes{\Ground_3} C^\op
\]
by
\[
  \delta^1(q\otimes p)=
  \sum_{(S,T)}n_{S,T}\bigl(\lambda^Q(S)\otimes\Id_P \otimes \mu^{\Clg}(T^{\op})\bigr)
  \circ (\Id_Q\otimes \delta_P),
\]
or, graphically:
\[
  \delta^1=
  \mathcenter{
\tikzsetnextfilename{sec-DA-tens-DD-1}
    \begin{tikzpicture}[smallpic]
      \node at (-2,-3) (bll) {};
      \node at (-1,0) (x) {};
      \node at (0,0) (y) {};
      \node at (0, -1) (delta) {$\delta$};
      \node at (-1, -2) (S) {$\TrDADD$};
      \node at (1, -2) (T) {$\TrDADD^\op$};
      \node at (-1, -3) (bl) {};
      \node at (0, -3) (bc) {};
      \node at (1, -3) (br) {};
      \draw[dmoda] (y) to (delta);
      \draw[dmoda] (delta) to (bc);
      \draw[moda] (x) to (S);
      \draw[moda] (S) to (bl);
      \draw[taa] (delta) to (S);
      \draw[alga] (S) to (bll);
      \draw[tbb] (delta) to (T);
      \draw[blga] (T) to (br);
    \end{tikzpicture}}.
\]
This formula and the corresponding type \DD\ structure relation make sense so long as either:
\begin{itemize}
\item $\lsup{\Blg,\Clg}P$ is bounded or %and $\Alg$ is bonsai or
\item $\lsup{\Alg}Q_\Blg$ is bounded and $\Blg$ and $\Clg$ are bonsai.
\end{itemize}

\begin{proposition}\label{prop:DADD-defined}
  Let $\lsup{\Alg}Q_\Blg$ be a type \DA\ structure with $\delta^1_1=0$
  and $\lsup{\Blg,\Clg}P$ a type \DD\ structure with respect to $\TrDiag^1$.  If
  $\lsup{\Alg}Q_\Blg$ is bounded and $\Blg$ and $\Clg$ bonsai, or
  $\lsup{\Blg,\Clg}P$ is bounded and $\Alg$ is bonsai, then the
  operation
  \[
    \delta^1\co Q\kotimes{\Ground_2} P \to A\kotimes{\Ground_1}
    Q\kotimes{\Ground_2}P \kotimes{\Ground_3} C^{\op}
  \]
  defined above gives $Q\kotimes{\Ground_2} P$ the structure of a type
  \DD\ structure
  $\lsup{\Alg}Q_{\Blg}\DT_{\TrDADD}\lsup{\Blg}P^{\Clg^\op}$ over
  $\Alg$ and $\Clg$ with
  respect to $\TrDiag^2$.
\end{proposition}
\begin{proof}
  We have:
  \begin{align*}
    &\phantom{=}
    \tikzsetnextfilename{DADD-prop-1}
    \mathcenter{
    \begin{tikzpicture}[smallpic]
      \node at (0,0) (tl) {};
      \node at (1,0) (tr) {};
      \node at (1,-1) (delta1) {$\delta$};
      \node at (1,-2) (rdots) {$\vdots$};
      \node at (1,-3) (delta2) {$\delta$};
      \node at (1,-6) (br) {};
      \node at (0,-2) (r1) {$\TrDADD$};
      \node at (0,-3) (ldots) {$\vdots$};
      \node at (0,-4) (r2) {$\TrDADD$};
      \node at (0,-6) (bl) {};
      \node at (-1,-5) (gammal) {$\TrDiag$};
      \node at (3,-5) (gammar) {$\TrDiag^\op$};
      \node at (2,-2) (rop1) {$\TrDADD^\op$};
      \node at (2,-3) (rrdots) {$\vdots$};
      \node at (2,-4) (rop2) {$\TrDADD^\op$};
      \node at (3,-6) (brr) {};
      \node at (-1,-6) (bll) {};
      \draw[damoda] (tl) to (r1);
      \draw[damoda] (r1) to (ldots);
      \draw[damoda] (ldots) to (r2);
      \draw[damoda] (r2) to (bl);
      \draw[dmoda] (tr) to (delta1);
      \draw[dmoda] (delta1) to (rdots);
      \draw[dmoda] (rdots) to (delta2);
      \draw[dmoda] (delta2) to (br);
      \draw[taa] (delta1) to (r1);
      \draw[taa] (delta1) to (rop1);
      \draw[taa] (delta2) to (r2);
      \draw[taa] (delta2) to (rop2);
      \draw[alga] (r1) to (gammal);
      \draw[alga] (r2) to (gammal);
      \draw[alga] (rop1) to (gammar);
      \draw[alga] (rop2) to (gammar);
      \draw[alga] (gammal) to (bll);
      \draw[alga] (gammar) to (brr);
    \end{tikzpicture}}
    +
    \tikzsetnextfilename{DADD-prop-2}
    \mathcenter{
    \begin{tikzpicture}[smallpic]
      \node at (0,0) (tl) {};
      \node at (1,0) (tr) {};
      \node at (1,-1) (delta) {$\delta$};
      \node at (0,-2) (r) {$\TrDADD$};
      \node at (2,-2) (rop) {$\TrDADD^\op$};
      \node at (-1,-3) (mul) {$\mu_1^\Alg$};
      \node at (-1,-4) (bll) {};
      \node at (2,-4) (brr) {};
      \node at (0,-4) (bl) {};
      \node at (1,-4) (br) {};
      \draw[dmoda] (tr) to (delta);
      \draw[dmoda] (delta) to (br);
      \draw[damoda] (tl) to (r);
      \draw[damoda] (r) to (bl);
      \draw[taa] (delta) to (r);
      \draw[taa] (delta) to (rop);
      \draw[alga] (r) to (mul);
      \draw[alga] (mul) to (bll);
      \draw[alga] (rop) to (brr);
    \end{tikzpicture}
    }
    +
    \tikzsetnextfilename{DADD-prop-3}
    \mathcenter{
    \begin{tikzpicture}[smallpic]
      \node at (0,0) (tl) {};
      \node at (1,0) (tr) {};
      \node at (1,-1) (delta) {$\delta$};
      \node at (0,-2) (r) {$\TrDADD$};
      \node at (2,-2) (rop) {$\TrDADD^\op$};
      \node at (2,-3) (mur) {$\mu_1^\Clg$};
      \node at (-1,-4) (bll) {};
      \node at (2,-4) (brr) {};
      \node at (0,-4) (bl) {};
      \node at (1,-4) (br) {};
      \draw[dmoda] (tr) to (delta);
      \draw[dmoda] (delta) to (br);
      \draw[damoda] (tl) to (r);
      \draw[damoda] (r) to (bl);
      \draw[taa] (delta) to (r);
      \draw[taa] (delta) to (rop);
      \draw[alga] (rop) to (mur);
      \draw[alga] (mur) to (brr);
      \draw[alga] (r) to (bll);
    \end{tikzpicture}
      }\\
    &=
      \tikzsetnextfilename{DADD-prop-4}
      \mathcenter{
      \begin{tikzpicture}[smallpic]
        \node at (-1,0) (tl) {};
        \node at (1,0) (tr) {};
        \node at (1,-1) (delta1) {$\delta$};
        \node at (1,-2) (delta2) {$\delta$};
        \node at (1,-3) (delta3) {$\delta$};
        \node at (0,-3) (gammal) {$\TrDiag$};
        \node at (2,-3) (gammar) {$\TrDiag$};
        \node at (-1,-4) (r) {$\TrDADD$};
        \node at (3,-4) (rop) {$\TrDADD$};
        \node at (3,-5) (brr) {};
        \node at (-1,-5) (bl) {};
        \node at (1,-5) (br) {};
        \node at (-2,-5) (bll) {};
        \draw[dmoda] (tr) to (delta1);
        \draw[dmoda] (delta1) to (delta2);
        \draw[dmoda] (delta2) to (delta3);
        \draw[dmoda] (delta3) to (br);
        \draw[damoda] (tl) to (r);
        \draw[damoda] (r) to (bl);
        \draw[taa] (delta1) to (r);
        \draw[taa] (delta2) to (gammal);
        \draw[taa] (delta3) to (r);
        \draw[taa] (delta1) to (rop);
        \draw[taa] (delta2) to (gammar);
        \draw[taa] (delta3) to (rop);
        \draw[alga] (gammal) to (r);
        \draw[alga] (gammar) to (rop);
        \draw[alga] (rop) to (brr);
        \draw[alga] (r) to (bll);
      \end{tikzpicture}
      }
      +
      \tikzsetnextfilename{DADD-prop-5}
    \mathcenter{
    \begin{tikzpicture}[smallpic]
      \node at (-1,0) (tl) {};
      \node at (1,0) (tr) {};
      \node at (1,-1) (delta1) {$\delta$};
      \node at (1,-2) (delta2) {$\delta^1$};
      \node at (1,-3) (delta3) {$\delta$};
      \node at (-1,-4) (r) {$\TrDADD$};
      \node at (2,-4) (rop) {$\TrDADD^\op$};
      \node at (0,-3) (mul) {$\mu_1^\Blg$};
      \node at (-2,-5) (bll) {};
      \node at (2,-5) (brr) {};
      \node at (-1,-5) (bl) {};
      \node at (1,-5) (br) {};
      \draw[dmoda] (tr) to (delta1);
      \draw[dmoda] (delta1) to (delta2);
      \draw[dmoda] (delta2) to (delta3);
      \draw[dmoda] (delta3) to (br);
      \draw[damoda] (tl) to (r);
      \draw[damoda] (r) to (bl);
      \draw[alga] (delta2) to (mul);
      \draw[alga] (delta2) to (rop);
      \draw[taa] (delta1) to (rop);
      \draw[taa] (delta3) to (rop);
      \draw[taa] (delta1) to (r);
      \draw[taa] (delta3) to (r);
      \draw[alga] (mul) to (r);
      \draw[alga] (r) to (bll);
      \draw[alga] (rop) to (brr);
    \end{tikzpicture}
    }
    +
    \tikzsetnextfilename{DADD-prop-6}
    \mathcenter{
    \begin{tikzpicture}[smallpic]
      \node at (0,0) (tl) {};
      \node at (1,0) (tr) {};
      \node at (1,-1) (delta1) {$\delta$};
      \node at (1,-2) (delta2) {$\delta^1$};
      \node at (1,-3) (delta3) {$\delta$};
      \node at (0,-4) (r) {$\TrDADD$};
      \node at (3,-4) (rop) {$\TrDADD^\op$};
      \node at (2,-3) (mur) {$\mu_1^\Clg$};
      \node at (-1,-5) (bll) {};
      \node at (3,-5) (brr) {};
      \node at (0,-5) (bl) {};
      \node at (1,-5) (br) {};
      \draw[dmoda] (tr) to (delta1);
      \draw[dmoda] (delta1) to (delta2);
      \draw[dmoda] (delta2) to (delta3);
      \draw[dmoda] (delta3) to (br);
      \draw[damoda] (tl) to (r);
      \draw[damoda] (r) to (bl);
      \draw[alga] (delta2) to (mur);
      \draw[alga] (delta2) to (r);
      \draw[taa, bend left=10] (delta1) to (rop);
      \draw[taa] (delta3) to (rop);
      \draw[taa] (delta1) to (r);
      \draw[taa] (delta3) to (r);
      \draw[alga] (mur) to (rop);
      \draw[alga] (r) to (bll);
      \draw[alga] (rop) to (brr);
    \end{tikzpicture}
      }
      =0.
  \end{align*}
  Here, the first equation follows from the definition of a DADD diagonal
  and Lemma~\ref{lem:DA-lambda-chain-map}, while the second follows
  from the \DD\ structure relation.
\end{proof}

Finally, we relate these tensor products with the one-sided box tensor
product of a $\Ainf$-module and a type \DD\ structure. Let $\Alg$ be a
strictly unital $\Ainf$-algebra and $\cModule_\Alg$ an $\Ainf$-module over
$\Alg$. We can view $\cModule_\Alg$ as a type \DA\ structure $\lsup{\Ring}\cModule_\Alg$ over $\Ring$ and
$\Alg$ by declaring that
\[
  \delta^1_{1+n}(x,a_1,\dots,a_n)=\unit\otimes m_{1+n}(x,a_1,\dots,a_n).
\]
So, a DADD diagonal also induces a tensor product of an $\Ainf$-module
with a type \DD\ structure. We relate this tensor product to the one
coming from a module diagonal primitive:
\begin{proposition}\label{prop:DADD-prim-box}
  Fix a strictly unital $\Ainf$-algebras $\Alg$ and $\Blg$.  Fix also
  a DADD diagonal $\TrDADD$ and let $\TrPMDiag$ be the induced module
  diagonal primitive, as in Proposition~\ref{prop:DADD-to-prim}. Let
  $\cModule_\Alg$ be an $\Ainf$-module over $\Alg$ and $\lsup{\Alg}P^\Blg$ a
  type \DD\ structure over $\Alg$ and $\Blg$. Assume that $P$ is
  bounded or $\cModule$, $\Alg$, and $\Blg$ are bonsai. Then there is an
  isomorphism of type \DA\ structures over $\Ring$ and $\Blg$
  \[
    (\lsup{\Ring}\cModule_\Alg)\DT_{\TrDADD}\lsup{\Alg}P^\Blg\cong \lsup{\Ring}\bigl(\cModule_\Alg\DT_{\TrPMDiag}\lsup{\Alg}P^\Blg\bigr).
  \]
\end{proposition}
\begin{proof}
  This is immediate from the definitions.
\end{proof}

\subsubsection{External tensor product of a DA structure and an
  \texorpdfstring{$\Ainf$-}{A-infinity }module}

Given a type \DA\ structure $\lsub{\Blg}X^{\Clg}$, an $\Ainf$-module
$\lsub{\Alg}\cModule$, and a module diagonal primitive, we can form the \DA\
structure $\lsub{\Alg^{\op}\otimes\Blg}(\cModule^{\op}\otimes X)^{\Clg}$.
Diagrams depicting the structure map on $\cModule^{\op}\otimes X$
are reminiscent of the proof of Lemma~\ref{eq:tp-is-tp}. We draw the
algebra elements
${\underline a}=(a_1\otimes b_1)\otimes\dots(a_\ell\otimes b_\ell)$ in
the middle:
\begin{equation}
  \label{eq:AxDA}
  \delta^1(x,{\underline a},y)=
  \left(
    \mathcenter{
\tikzsetnextfilename{eq-AxDA-1}
      \begin{tikzpicture}[smallpic]
        \node at (-1,0) (tl) {x};
        \node at (0,0) (tc) {$\emptyset$};
        \node at (1,0) (tr) {y};
        \node at (2,-2) (brr) {};
        \node at (-1,-2) (bl) {};
        \node at (0,-2) (bc) {};
        \node at (1,-2) (br) {};
        \node at (2,-1) (one) {$\One$};
        \node at (-1,-1) (m) {$m_1$};
        \draw[moda] (tr) to (br);
        \draw[moda] (tl) to (m);
        \draw[moda] (m) to (bl);
        \draw[alga] (one) to (brr);
      \end{tikzpicture}
    }\right)+
  \left(
    \mathcenter{
\tikzsetnextfilename{eq-AxDA-2}
      \begin{tikzpicture}[smallpic]
        \node at (0,0) (tc) {${\underline a}$};
        \node at (0,-1) (delta1) {$\Delta$};
        \node at (0,-2) (vdots) {$\vdots$};
        \node at (0,-3) (delta2) {$\Delta$};
        \node at (0,-6) (bc) {};
        \node at (-1,-2) (pl1) {$\TrPMDiag$};
        \node at (-1,-3) (lvdots) {$\vdots$};
        \node at (-1,-4) (pl2) {$\TrPMDiag$};
        \node at (1,-2) (pr1) {$\TrPMDiag^\op$};
        \node at (1,-3) (rvdots) {$\vdots$};
        \node at (1,-4) (pr2) {$\TrPMDiag^\op$};
        \node at (-1,0) (tl) {x};
        \node at (2,0) (tr) {y};
        \node at (2,-6) (br) {};
        \node at (-1,-6) (bl) {};
        \node at (2,-5) (m) {$\delta^1$};
        \node at (3,-6) (brr) {};
        \draw[alga](m) to (brr);
        \draw[taa] (tc) to (delta1);
        \draw[taa] (delta1) to (vdots);
        \draw[taa] (vdots) to (delta2);
        \draw[taa] (delta1) to (pl1);
        \draw[taa] (delta2) to (pl2);
        \draw[alga] (pr1) to (m);
        \draw[alga] (pr2) to (m);
        \draw[moda] (tr) to (m);
        \draw[moda] (m) to (br);
        \draw[moda] (tl) to (pl1);
        \draw[moda] (pl1) to (lvdots);
        \draw[moda] (lvdots) to (pl2);
        \draw[moda] (pl2) to (bl);
        \draw[tbb] (delta1) to (pr1);
        \draw[tbb] (delta2) to (pr2);
      \end{tikzpicture}
    }
  \right),
\end{equation}
Here,
\[
  \Delta\co \Tensor^*(A\otimes B)\to \Tensor^*(A\otimes B)\otimes
  \Tensor^*(A\otimes B)\otimes \Tensor^*(A\otimes B)
\]
denotes all ways of splitting up the sequence $\underline{a}$ into three (or, at
the bottom, two) sequences (without breaking up any individual factor
$a_i\otimes b_i$). The first term contributes only when the incoming sequence of
algebra elements is empty (i.e., $\ell=0$).

\begin{proposition}
  \label{prop:ExternalAwithDA}
  Let $\Alg$ and $\Blg$ be $\Ainf$-algebras over $\Ground_1$ and
  $\Ground_2$ and $\Clg$
  a strictly unital $\Ainf$-algebra over $\Ground_3$, and assume
  that $\cModule_\Alg$ and $\lsub{\Blg}N^{\Clg}$ are a bonsai module and
  bounded \DA\ 
  structure, respectively.  Then $M\rotimes{\Ring} N$ can be given the
  structure of a \DA\ structure, with input algebra
  $\Alg\rotimes{\Ring}\Blg^{\op}$ on the right and output algebra $\Clg$ (on
  the left), with operations specified in Equation~\eqref{eq:AxDA}
\end{proposition}
\begin{proof}
  The proof is a minor adaptation of the proof of Lemma~\ref{lem:DT-defined}.
\end{proof}

\subsection{Units in \texorpdfstring{$\Ainf$-}{A-infinity-}algebras}
\label{subsec:UnitsHomPert}
In this section we recall that a weaker notion of a unit an $\Ainf$-algebra can be
strictified to a strict unit, and discuss corresponding notions for modules and
morphisms.

\begin{definition}
  \label{def:WeakUnitalAinf}
  An $\Ainf$-algebra $\Alg$ is \emph{weakly unital} if there is an
  element $\unit\in A$ so that $\mu_1(\unit)=0$ and for all $a\in A$,
  $\mu_2(a,\unit)=\mu_2(\unit,a)=a$.
\end{definition}

\begin{definition}\label{def:strict-unital}
  A weakly unital $\Ainf$-algebra $\Alg$ is \emph{strictly unital} if
  for all $n\geq 3$ and $a_1,\dots,a_n\in \Alg$, if some $a_i=\unit$
  then $\mu_n(a_1,\dots,a_n)=0$.
\end{definition}

If $\Alg$ and $\Blg$ are weakly unital $\Ainf$-algebras, it follows
from the non-degeneracy condition of an associahedron diagonal that
$\Alg\ADtp\Blg$ is weakly unital. By contrast, even if
$\Alg$ and $\Blg$ are strictly unital, $\Alg\ADtp\Blg$ may
not be (unless $\AsDiag$ satisfies an extra condition). Thus, if one
prefers to work with strictly unital algebras, one needs a procedure
to turn a weakly unital one into a strictly unital one (or to work
with special diagonals).

There is an analogous notion for homomorphisms:
\begin{definition}
  Let $\Alg$ and $\Blg$ be two weakly unital $\Ainf$-algebras over $\Ground$.  A
  homomorphism $f\co \Alg\to \Blg$ (in the sense of
  Definition~\ref{def:AlgebraHomomorphism}) is {\em weakly
    unital} if $f_1$ carries the unit $\unit_\Alg$ in $\Alg$ to the
  unit $\unit_\Blg$ in $\Blg$. A homomorphism $f$ is \emph{strictly
    unital} if $f_1(\unit_\Alg)=\unit_\Blg$ and for each $n>1$,
  $f_n(a_1,\dots,a_n)=0$ if some $a_i=\unit_\Alg$.
\end{definition}

The following is a special case of the well-known result that
homologically unital $\Ainf$-algebras (a notion we will not use) are
quasi-isomorphic to strictly unital
ones~\cite{Lazarev,LefevreAInfinity,SeidelBook}; our proof follows Seidel~\cite{SeidelBook}:
\begin{theorem}\label{thm:UnitalIsUnital}
  Every weakly unital $\Ainf$-algebra $\Alg$ over $\Ground$ is isomorphic to 
  a strictly unital $\Ainf$-algebra~$\Alg'$, and the isomorphism can
  be chosen to be weakly unital. If $\Alg$ is bonsai then
  $\Alg'$ and the isomorphism can both be chosen to be bonsai as well.
\end{theorem}

We sketch a proof of Theorem~\ref{thm:UnitalIsUnital}, after some preliminary
discussion.
To set it up, let 
\[ \phi_n\co \overbrace{A\kotimes{\Ground}\dots\kotimes{\Ground} A}^n\to B \]
be a map of $\Ground$-bimodules. Define
\[ d \phi_n = \mu_1\circ \phi_n + \sum_{i=1}^n \phi_n\circ (\Id_{A^{\kotimes{\Ground} i-1}} \otimes \mu_1 \otimes \Id_{A^{\kotimes{\Ground} {n-i}}}).\]
(This is the differential of $\phi_n$ thought of as an element of
$\Mor(A^{\kotimes{\Ground} n},B)$, the morphism complex of chain complexes.)

\begin{lemma}\label{lem:ModifyAction}\cite[Lemme 3.2.2.2]{LefevreAInfinity}
Suppose that $\Alg$ is an $\Ainf$-algebra over $\Ground$, and 
$\phi_n\co A ^{\kotimes{\Ground} n}\to A$ is any $\Ground$-bimodule map
of degree $n-1$ 
with $n>1$ (where the tensor products are taken over $\Ground$).
Then there is an $\Ainf$-algebra ${\overline\Alg}$ isomorphic to $\Alg$ so that 
${\overline\mu}_i=\mu_i$ for $i<n$ and
\[ {\overline \mu}_n = \mu_n + d \phi_n.\]
Moreover, if $n>2$ or $n=2$ and $d\phi_2=0$ then
\[ 
{\overline\mu}_{n+1}=\mu_{n+1}+\mu_2\circ (\Id_{A}\otimes \phi_n) + \mu_2\circ (\phi_n\otimes \Id_{A}) + 
\sum_{i=1}^{n-1} \phi_n(\Id_{A^{\kotimes{\Ground} i-1}}\otimes \mu_2\otimes
\Id_{A^{\kotimes{\Ground} n-i-1}}).
\]
\end{lemma}

\begin{proof}
  By the \emph{trunk vertex} of a rooted, planar tree we mean the internal
  vertex adjacent to the root.
  
  We will express the operations on ${\overline\Alg}$ graphically, as follows.  An {\em
    ${\overline\Alg}$-operation tree} is a planar, rooted tree $T$
  with internal vertices labelled as follows:
  \begin{itemize}
    \item There is one distinguished vertex, labelled by the
      operation~$\mu_k$ with $k \ge 1$, with $k$ inputs.
    \item All other vertices have $n$ inputs, and are labelled by $\phi$.
    \item The distinguished vertex is either the trunk vertex, or it
      is one of the parents of the trunk vertex.
  \end{itemize}
  Let ${\mathcal T}'$ denote the set of ${\overline\Alg}$-operation trees.
  For example, if $j < n$, then there is a unique ${\overline\Alg}$-operation
  tree with $j$ inputs: the corolla $\corolla{j}$ marked by $\mu_j$.
  There are $n+2$
  $\overline{\Alg}$-operation trees with $n$ inputs
  (a distinguished corolla, or a distinguished $2$-valent vertex
  attached to an undistinguished vertex), and for $n > 2$,
  there are $n+3$ different ${\overline\Alg}$-operation trees with $n+1$ inputs,
  as in
  Figure~\ref{fig:AlgPrimeOperationTrees}.
  \begin{figure}
    \centering
    \input{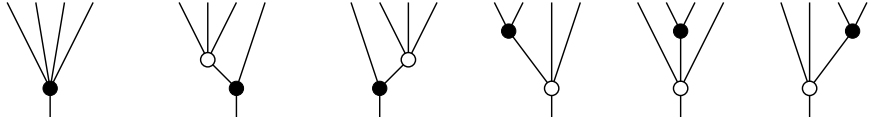}
    \caption[Operation trees used to make a weakly unital algebra strictly unital]{{\bf The ${\overline\Alg}$-operation trees with $n=3$ and $j=4$ inputs.}
      The white vertices are labelled by $\phi_3$, and the black
      vertex by a $\mu$-operation.}
    \label{fig:AlgPrimeOperationTrees}
  \end{figure}

  Summing over all ${\overline\Alg}$-operation trees with $j$ inputs induces a map
  \[
    {\overline\mu}_j\co A^{\kotimes{\Ground} j}\to A.
  \]

  We assert that the operations ${\overline\mu}_j$ satisfy the $\Ainf$ relations
  \[ 
  \sum_{k=0}^{m-1}\sum_{i=1}^{m-k}{\overline\mu}_{m-k}\circ 
  (\Id_{A^{\kotimes{\Ground}(i-1)}}\otimes {\overline\mu}_{k+1}\otimes \Id_{A^{\kotimes{\Ground}(m-i-k)}})=0.
  \]
  This amounts to proving that the sum over all trees $T_1, T_2\in
  {\mathcal T}'$ of the operation obtained by composing $T_1$
  and $T_2$ vanishes.

  To this end, let ${\mathcal S}'$ be the set of trees $S$ that can be
  obtained by composing $T_1,T_2\in {\mathcal T}'$. This
  means that $S$ has exactly two distinguished vertices that are
  marked by operations in $\Alg$, at least one of them has distance
  $\leq 1$ from the trunk vertex, and all other vertices are marked by~$\phi_n$.

  Suppose that $S\in{\mathcal S}'$ has two consecutive distinguished
  vertices. Such trees $S$ can be uniquely decomposed as the
  juxtaposition of $T_1$ and $T_2$: the trees $T_1$ and $T_2$ are
  obtained by cutting along the edge connecting the two distinguished
  vertices in $S$. The set of such trees can be alternatively
  obtained from the set ${\mathcal T}'$ by inserting an edge at the
  distinguished vertex. From this description, it is clear that the
  $\Ainf$ relation on $\Alg$ ensures that the sum of the
  contributions of $S\in {\mathcal S'}$ with two consecutive
  distinguished vertices vanishes.

  To complete the argument we claim that if $S\in {\mathcal S}'$ can
  be written as a juxtaposition of $T_1,T_2\in {\mathcal T}'$ so that
  the distinguished vertices are not consecutive, then $S$ has exactly
  two such decompositions.  So, suppose that $S$ is of this
  form, and let $d_1$ and $d_2$ denote the distances of the two
  distinguished vertices $v_1$ and $v_2$ to the root, chosen so that
  $d_1\leq d_2$. Since we assumed that the distinguished vertices are
  not consecutive, there are two cases: either $d_1=d_2=1$ or $d_2>1$.
  If $d_1=d_2=1$, the two decompositions are obtained by cutting the
  edge immediately below $v_1$ or $v_2$.  If $d_2>1$, then we can
  either cut the edge immediately below $v_2$ or two below $v_2$ to
  obtain the two decompositions of $S$ into trees in ${\mathcal T}'$.
  See Figure~\ref{fig:ProveItsAinfty}.
\begin{figure}
  \centering
  \input{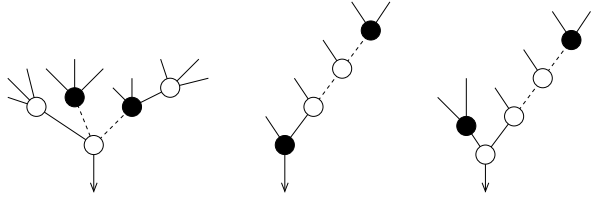}
  \caption[Cancellations of tree juxtapositions]{{\bf Cancellations of tree juxtapositions.}
    We have drawn trees $S\in{\mathcal S}'$, with distinguished vertices
    black. (In the first case, the valence $n$ of the undistinguished vertices is $3$;
    while in the second two $n=2$.)
    In each case, cutting either of the two dotted edges gives juxtapositions
    of pairs of trees in ${\mathcal T}'$.}
  \label{fig:ProveItsAinfty}
\end{figure}

Next, we argue that $\Alg$ and ${\overline\Alg}$ are isomorphic.
We construct a morphism $\phi\co \Alg\to {\overline\Alg}$ that has $\phi_1$ the
identity map, $\phi_n$ as specified in the lemma, and $\phi_i=0$ for
$i\neq 1,n$. Since any $\Ainf$-algebra homomorphism $\phi$ with
$\phi_1$ an isomorphism is an $\Ainf$-algebra isomorphism~\cite[Lemma
2.1.14]{LOT2}, all that remains is to verify the $\Ainf$ relations for $\phi$.

To this end, consider one of the $\Ainf$ relations and all the trees
that appear in it, once the purple vertices have been replaced by
either the identity (and erased) or $\phi$, the red vertices have been
replaced by a $\mu_k$, and the blue vertices have been replaced by an
$\overline{\Alg}$-operation tree as above. These trees all have a
single distinguished vertex, and, in fact, are all in ${\mathcal
  T}'$. We must verify that each tree $T\in {\mathcal T}'$ appears an
even number of times in the $\Ainf$ relation.  Call a vertex $v$ of
$T$ \emph{highest} if all parents of $v$ are leaves of $T$. Suppose
that $T$ has $k$ highest, undistinguished vertices. Then $T$ appears
with multiplicity $2^k$ in part of the $\Ainf$ relation corresponding
to ${\overline\mu}\circ \phi$. (The case $k=2$ is shown in
Figure~\ref{fig:ProveItsAinftyMap}.) When $k=0$, the tree also appears
once in $\phi\circ\mu$; and indeed, those are precisely the trees that
appear in $\phi\circ\mu$.
\begin{figure}
  \centering
  \input{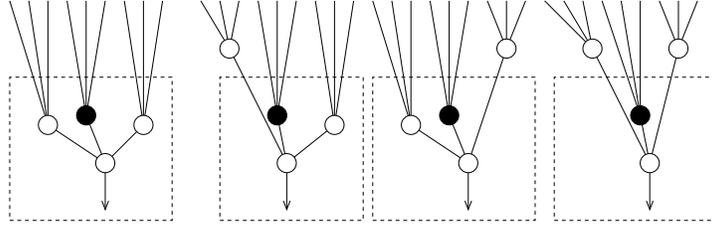}
  \caption[Verifying the A-infinity relation for a particular map]{{\bf Verifying the $\Ainf$ relation for the map.}
    The tree pictured here has $k=2$, and it has $4$ contributions
    to the $\Ainf$ relation for the morphism determined by $\phi$.
    The portions of the trees corresponding to the differential $\mu'$
    are contained in the dashed boxes.}
  \label{fig:ProveItsAinftyMap}
\end{figure}

The statement about the form of $\overline{\mu}_{n+1}$ follows from our more
general description of ${\overline\mu}_j$.
\end{proof}

\begin{remark}
  \label{rem:ModifyActionSimple}
  A shorter, but less explicit proof of Lemma~\ref{lem:ModifyAction}
  is the following (cf.~\cite[Lemme 3.2.2.2]{LefevreAInfinity},~\cite[Section (1c)]{SeidelBook}). Given an
  $\Ainf$-algebra $\Alg$ and a collection of maps $f_n\co A^{\otimes
    n}\to B$ with $f_1$ an isomorphism, one can inductively solve the
  $\Ainf$-homomorphism relations to define operations $\mu_n^\Blg\co
  B^{\otimes n}\to B$ satisfying Equation~\eqref{eq:Ainf-homo}. It
  follows from the fact that $\Blg=(B,\{\mu_n^\Blg\})$ is isomorphic
  to $\Alg$ (via $\{f_n\}$) that the $\mu_n^\Blg$ satisfy the
  $\Ainf$-algebra relations. Lemma~\ref{lem:ModifyAction} is the
  special case $f_1=\Id$, $f_n=\phi_n$ given, and $f_i=0$ for
  $i\not\in\{1,n\}$.
\end{remark}

We will also use the following property of bonsai-ness. Suppose that
$T$ is a planar, rooted tree with $n$ input leaves, but which is
allowed to have $2$-valent vertices, and $\Alg=(A,\mu)$ is an
$\Ainf$-algebra. There is still an induced map
$\mu(T)\co A^{\kotimes{\Ground} n}\to A$. Define the dimension of $T$ as in
Definition~\ref{def:dim}; note that each $2$-valent vertex decreases
the dimension of $T$ by $1$. 
\begin{lemma}\label{lem:bonsai-2-valent}
  If $\Alg$ is a bonsai $\Ainf$-algebra with bonsai constant $N$ and
  $T$ is a planar, rooted tree possibly with $2$-valent vertices with
  $\dim(T)>N$ then $\mu(T)=0$. Similar statements hold for modules,
  algebra homomorphisms, and module morphisms.
\end{lemma}
\begin{proof}
  Fix a tree $T$ with $\dim(T)>N$.
  Inductively applying the $\Ainf$ relation at each $n>2$ valent
  vertex which is followed by a $2$-valent vertex, $\mu(T)$ is equal
  to a linear combination of operations $\mu(T')$ where all of the
  $2$-valent vertices of $T'$ are adjacent to the input leaves, with
  $\dim(T')=\dim(T)$. So, it suffices to prove the result for such
  trees $T'$. Suppose that $T'$ has $k$ $2$-valent vertices. Let $S$
  be the result of forgetting the $2$-valent vertices of $T'$. Then
  $\dim(S)=\dim(T')+k>N+k$. So, $\mu(S)$ vanishes; but this implies
  that $\mu(T')$ vanishes, as well.
\end{proof}

\begin{proof}[Proof of Theorem~\ref{thm:UnitalIsUnital}]
  We wish to find an isomorphic model for $\Alg$ satisfying
  \begin{equation}
    \label{eq:PartialUnitality}
    \mu_{m+1}(a_1,\dots,a_{j-1},\unit,a_{j},\dots,a_{m})=0 
  \end{equation}
  for all $m>1$ and $j=1,\dots,m+1$. We will inductively modify the
  actions $\mu_n$ on $\Alg$, for $n=3,4,\dots$, to obtain
  Equation~\eqref{eq:PartialUnitality}.

  So, suppose there are $n\geq 2$ and $k$ so that
  Equation~\eqref{eq:PartialUnitality} holds for all $1<m<n$ and for
  $m=n$ and $j<k$. We will construct an isomorphic algebra $\Alg$
  which has the same underlying $\Ground$-bimodule and $\mu_i$ actions
  for $i<n+1$, and a new $\mu_{n+1}$ action satisfying
  Equation~\eqref{eq:PartialUnitality} for all $j<k+1$, provided that
  $k+1<n+1$; the case $k=n$ is dealt with separately at the end.

  Let
  \[
    \phi_n(a_1,\dots,a_n)=\mu_{n+1}(a_1,\dots,a_{k},\unit,a_{k+1},\dots,a_n).
  \]
  The $\Ainf$ relation with input
  $(a_1,\dots,a_k,\unit,a_{k+1},\dots,a_n)$ and the inductive hypothesis shows
  that $d \phi_n=0$. Thus, we can apply Lemma~\ref{lem:ModifyAction}, to construct a new action
  with 
  \[
    {\overline \mu}_{n+1}=\mu_{n+1}+\mu_2\circ (\Id_{A}\otimes \phi_n) + \mu_2\circ (\phi_n\otimes \Id_{A}) + 
    \sum_{i=1}^{n-1} \phi_n(\Id_{\Alg^{\kotimes{\Ground} i-1}}\otimes \mu_2\otimes \Id_{\Alg^{\kotimes{\Ground} n-i}}).
  \]
  The $\Ainf$ relation with input $(a_1,\dots,a_k,\unit,\unit,a_{k+1},\dots,a_n)$
  proves that
  \[{\overline\mu}_{n+1}(a_1,\dots,a_k,\unit,a_{k+1},\dots,a_n)+(d \psi_{n+1})(a_1,\dots,a_k,\unit,a_{k+1},\dots,a_n)=0,\]
  where $\psi_{n+1}=\psi^k_{n+1}\co A^{\kotimes{\Ground} n+1}\to A$ is the function defined by
  \[
    \psi_{n+1}(a_1,\dots,a_{n+1})=\mu_{n+2}(a_1,\dots,a_k,\unit,a_{k+1},\dots,a_{n+1}).\]
  Further, it follows from the $\Ainf$ relation with input
  $(a_1,\dots,a_{i-1},\unit,a_{i+1},\dots,a_k,\unit,a_{k+1},\dots,a_{n+1})$
  that
  \[
    (d\psi_{n+1})(a_1,\dots,a_{i-1},\unit,a_{i+1},\dots,a_{n+1})=0
  \]
  if $i\leq k$.  Applying Lemma~\ref{lem:ModifyAction} once again,
  this time using $\psi_{n+1}$, we find a new algebra
  ${\overline{\overline\Alg}}$ isomorphic to $\Alg$ whose operations
  satisfy Equation~\eqref{eq:PartialUnitality} with the following
  possible pairs $(m,j)$:
  \begin{itemize}
  \item for all $1<m<n$ and $j=1,\dots,m+1$ and
  \item for $m=n$ and $j<k+1$.
  \end{itemize}
  Thus, by induction on $k$, we can arrange for
  \[ \mu_{m+1}(a_1,\dots,a_{j},\unit,a_{j+1},\dots,a_m)=0\]
  if $m<n$ or $m=n$ and $j<n$.
  We can extend to the case where $m=n$ and $j=n$
  by another application of Lemma~\ref{lem:ModifyAction}.
  Specifically, 
  for $\psi_{n+1}=\mu_{n+1}(a_1,\dots,a_{n},a_{n+1},\unit)$
  Lemma~\ref{lem:ModifyAction} gives a new $\Ainf$ structure
  with
  ${\overline \mu}_{n+1}=\mu_{n+1}+d\psi_{n+1}$.
  The $\Ainf$ relation with input
  $(a_1,\dots,a_n,\unit,\unit)$ ensures that
  \[   \mu_{n+1}(a_1,\dots,a_{n},\unit)+d \psi_{n+1}(a_1,\dots,a_{n},\unit)=0. \]

  The above inductive procedure constructs a sequence of isomorphic
  algebras ${\Alg}_{n,k}$, where the algebra operations $\mu_i$ on
  ${\Alg}_{n,k}$ and $\Alg_{m,\ell}$ coincide as long as
  $i\leq\min\{m,n\}$. The desired algebra $\Alg_{\infty,\infty}$ has
  the same $\mu_i$ as ${\Alg}_{n,k}$ for any $n\geq i$.

  Observe that at each stage of the induction, the isomorphism $f$
  given by Lemma~\ref{lem:ModifyAction} has $f_1=\Id$. In
  particular, $f$ is weakly unital.
  
  Finally, we must show that if $\Alg$ is bonsai then so is
  $\Alg_{\infty,\infty}$. Observe that, from the proof of
  Lemma~\ref{lem:ModifyAction} applied to the $\psi_{n+1}^k$ defined
  above, each $\Ainf$ operation $\mu_n$
  on $\Alg_{\infty,\infty}$ is a finite sum of operations of the form
  $\mu^{\Alg}(T)$, with $\unit$ fed into some of the inputs and some
  2-valent vertices allowed. In
  particular, any $\Ainf$ operation tree for $\Alg_{\infty,\infty}$ is
  a linear combination of $\Ainf$ operation trees $\mu(T)$ for $\Alg$,
  of the same grading (dimension), with $\unit$ as some of the
  inputs.
  Thus, if
  the $\mu^{\Alg}(T)$ vanish for $T$ of sufficiently large dimension,
  by Lemma~\ref{lem:bonsai-2-valent},
  the same is true for $\Alg_{\infty,\infty}$. The isomorphism
  from $\Alg$ to $\Alg_{\infty,\infty}$ is also bonsai, by the same
  reasoning.
\end{proof}

This discussion has the following adaptation to modules.

\begin{definition}
  An $\Ainf$-module $\cModule$ over a weakly unital $\Ainf$-algebra is 
  called {\em weakly unital} if $m_2(x,\unit)=x$ for all $x\in M$.
  A weakly unital $\Ainf$-module over a strictly unital $\Ainf$-algebra
  is called {\em strictly unital} if for all $n\geq 2$,
  $x\in M$, and $a_1,\dots,a_{n}\in \Alg$, we have that
  $m_{1+n}(x,a_1,\dots,a_{n})=0$ if some $a_i=\unit$.
\end{definition}

We have the following analogue of Theorem~\ref{thm:UnitalIsUnital}:

\begin{theorem}
  \label{thm:UnitalIsUnitalM}
  Every weakly unital $\Ainf$-module $\cModule$ over a strictly unital
  $\Ainf$-algebra is isomorphic to a strictly unital module
  $\cModule'$. If $\cModule$ is bonsai we may arrange that $\cModule'$
  and the isomorphism from $\cModule$ to $\cModule'$ are
  also bonsai.
\end{theorem}

The proof rests on the following analogue of Lemma~\ref{lem:ModifyAction}:

\begin{lemma}
  \label{lem:ModifyActionM}
  Suppose that $\cModule$ is an $\Ainf$-module and
  \[ \phi_n \co M \kotimes{\Ground} A \kotimes{\Ground}\cdots\kotimes{\Ground} A \to M\]
  is any $\Ground$-module map of degree $n-1$. 
  There is an $\Ainf$-module ${\overline \cModule}$ isomorphic to 
  $\cModule$ so that
  ${\overline m}_i=m_i$ for $i<n$ and 
  \[ {\overline m}_n= m_n + d \phi_n.\]
\end{lemma}

\begin{proof}
  The proof is similar to the proof of Lemma~\ref{lem:ModifyAction},
  with the following modifications. The ${\overline\cModule}$-operation
  trees are planar, rooted trees with the following properties:
  \begin{itemize}
  \item There is a distinguished internal vertex,
  which is labelled $m_k$ if it is contained along the leftmost
  path, and $\mu_k$ otherwise.
  \item All other vertices occur along the leftmost path, and are
    labelled $\phi_n$.
  \item The distinguished vertex is either the trunk vertex (in which 
    case it is labelled $m_k$), or it is a parent of the trunk vertex (in which case it
    can be either $m_k$ or $\mu_k$).
  \end{itemize}
  These trees induce the operations
  \[ {\overline m}_k \co M\kotimes{\Ground} A^{\kotimes{\Ground} (k-1)}\to M.\]
  The proof that they satisfy the $\Ainf$ relations is as before.
  The $\Ainf$-map $\phi\co \cModule \to {\overline\cModule}$ is once again
  given by $\phi_1=\Id$ and $\phi_n$ as specified. It is straightforward to verify that $\phi$ is a homomorphism, and it then follows from Lemma~\ref{lem:mod-iso} that $\phi$ is an isomorphism.
\end{proof}

\begin{proof}[Proof of Theorem~\ref{thm:UnitalIsUnitalM}]
  This is similar to the proof of Theorem~\ref{thm:UnitalIsUnital}. Again, we assume that
  $m_{m+1}(x,a_1,\dots,\unit,a_{i+1},\dots,a_{m})=0$ if $1<m<n$ or $m=n$ and $i<k$.
  Applying the $\Ainf$ relation to 
  $(x,a_1,\dots,a_{k-1},\unit,\unit,a_{k+1},\dots,a_n)$, we find that if we choose
  \[
    \phi_{n+1}(x,a_1,\dots,a_n)=m_{n+1}(x,a_1,\dots,a_{k-1},\unit,a_k,a_{k+1},\dots,a_n),
  \]
  then
  \[ 
    m_{n+1}(x,a_1,\dots,a_{k-1},\unit,a_{k+1},\dots,a_n)+
    d\phi_{n+1}(x,a_1,\dots,a_{k-1},\unit,a_{k+1},\dots,a_n)=0.
  \]
  Further, applying the $\Ainf$ relation with input
  $(a_1,\dots,\unit,a_i,\dots,\unit,a_k,\dots,a_n)$, $i<k$, implies that 
  \[
    (d \phi_{n+1})(x,a_1,\dots,a_{i-1},\unit,a_i,\dots,a_n)=0
  \]
  if $i<k$.  Applying Lemma~\ref{lem:ModifyActionM} and induction
  gives the desired module $\cModule'$. The fact that $\cModule'$ and
  the isomorphism are bonsai if $\cModule$ is bonsai follows by
  the same reasoning as in Theorem~\ref{thm:UnitalIsUnital}.
\end{proof}

\begin{definition}
  Let $\Alg$ be a strictly unital $\Ainf$-algebra over $\Ground$, and fix strictly unital
  $\Ainf$-modules $\cModule$ and $\cNodule$ over $\Alg$.
  The \emph{strictly unital morphism complex} $\suMor(\cModule,\cNodule)$ is 
  the subcomplex of 
  \[ \Mor(\cModule,\cNodule)=\prod_{i=0}^{\infty} \Mor(M\kotimes{\Ground} A^{\kotimes{\Ground} i},N)\]
  (equipped with the differential specified in
  Equation~\eqref{eq:d-on-mor}) consisting of those $f_n$ with the
  property that
  $f_{1+n}(x,a_1,\dots,a_{i-1},\One,a_{i+1},\dots,a_{n-1})=0$ for all
  $n\ge2$ and $1\leq i<n$.
\end{definition}
(To see that the strictly unital morphisms form a subcomplex note that
the two terms in the differential there $\One$ is fed into a $\mu_2$
or $m_2$ cancel.)

\begin{proposition}\label{prop:mod-mor-cx-unital}
  Given strictly unital $\Ainf$-modules $\cModule$ and $\cNodule$ over a
  strictly unital $\Ainf$-algebra $\Alg$, the inclusion of the 
  strictly unital morphism complex $\suMor(\cModule,\cNodule)$ into the full
  morphism complex $\Mor(\cModule,\cNodule)$ is a chain homotopy equivalence. If $\cModule$ and $\cNodule$ are bonsai then the statement also holds if we let $\Mor(\cModule,\cNodule)$ denote the complex of bonsai morphisms and $\suMor(\cModule,\cNodule)$ the subcomplex of strictly unital bonsai morphisms.
\end{proposition}

\begin{proof}
  Let ${\mathcal F}_i \subset
  \Mor(\cModule,\cNodule)$ consist of those morphisms
  so that for all  $n$ and all $1\leq j\leq i$,
  \[f_n(x,a_1,\dots,a_{j-1},\One,a_{j+1},\dots,a_{n-1})=0.\]
  (This condition is vacuous when $j\geq n-1$, of course.)
  It is straightforward to check that ${\mathcal F}_i$
  is indeed a subcomplex.
  
  Consider the map
  \begin{equation}
    \label{eq:DefHi}
    H_i\co \Mor(\cModule,\cNodule) \to \Mor(\cModule,\cNodule)
  \end{equation}
  defined by
  \[ H_i(f)(x,a_1,\dots,a_{n-2})=
    \begin{cases}
      0 &{\text{if $n-1<i$}} \\
      f_n(x,a_1,\dots,a_{i-1},\One,a_{i},\dots,a_{n-2}) & \text{otherwise}.
    \end{cases}
  \]
  It is easy to see that $H_i$ maps the subcomplex ${\mathcal F}_i$ to
  itself.  
  
  Let $\Pi_i \co \Mor(\cModule,\cNodule)\to\Mor(\cModule,\cNodule)$  be the map $\Pi_i=\Id + \partial \circ
  H_i + H_i\circ \partial$. The map $\Pi_i$ clearly preserves the
  subcomplex of bonsai morphisms and, in fact, the bonsai constant is
  unchanged by~$\Pi_i$.
  Most of the rest of the proof is to check that $\Pi_i({\mathcal F}_i)\subset  {\mathcal F}_{i+1}$.  
  
  Note that $\partial \circ H_i (f)$ consists of operations associated
  to trees with two internal vertices, of the following kinds:
  \begin{enumerate}[label=(t-\arabic*),ref=(t-\arabic*)]
  \item 
    \label{Nmult} A vertex labelled $f_{k}$ that feeds into another vertex labelled $m_\ell^{\cNodule}$,
    and the $i^{th}$ algebra input to the $f_{k}$-vertex is $\One$.
  \item 
    \label{Cmult}
    A vertex labelled $\mu_k^{\Alg}$ that feeds into another vertex
    labelled by $f_{\ell}$, and the $i^{th}$ input of the $f_\ell$ is $\One$.
  \item 
    \label{Mmult}
    A vertex labelled by $m_k^{\cModule}$ that feeds into another
    vertex labelled by $f_\ell$, and the $i^{th}$ input of the $f_\ell$ is $\One$.
  \end{enumerate}
  We call these contributions {\em $DH$-trees}.  Terms in
  $H_i\circ \partial$ have the same three kinds of vertices, except that
  instead of the $i\th$ algebra input of the $f_k$ vertex being $\One$, the
  $i\th$ algebra input to the whole configuration is
  $\One$.  We call
  these trees {\em $HD$-trees}.  See Figures~\ref{fig:MorHDTrees}
  and~\ref{fig:MorDHTrees}.
  
  In view of the strict unitality of $\cModule$,
  the $HD$ trees of Type~\ref{Nmult} with more than $i+1$ inputs and
  where the $f_k$ vertex has $k<i$ have zero contribution.
  In the special case where the total number of inputs is exactly $i+1$, there is an $HD$ tree
  of Type~\ref{Nmult}, where the vertex labelled $m^{\cNodule}_\ell$ has $\ell=2$; this corresponds to 
  the operation
  \[ (x,a_1,\dots,a_{i-1})\mapsto m_2^{\cNodule}(f_{i}(x,a_1,\dots,a_{i-1}),\One).\]
  This map is vacuously in ${\mathcal F}_{i+1}$, since $f\in{\mathcal
    F}_i$ and this tree has only $i$ inputs.
  All the other $DH$-trees and
  $HD$-trees of Type~\ref{Nmult} cancel against each other.
  
  Similarly, the $DH$-trees of Type~\ref{Cmult} cancel against the corresponding $HD$-trees, provided that the $\mu_k^{\Alg}$
  feeds into some input after the $i^{th}$ input. We call the remaining terms {\em Type~\ref{Cmult} leftovers}.

  \begin{figure}
    \centering
    \input{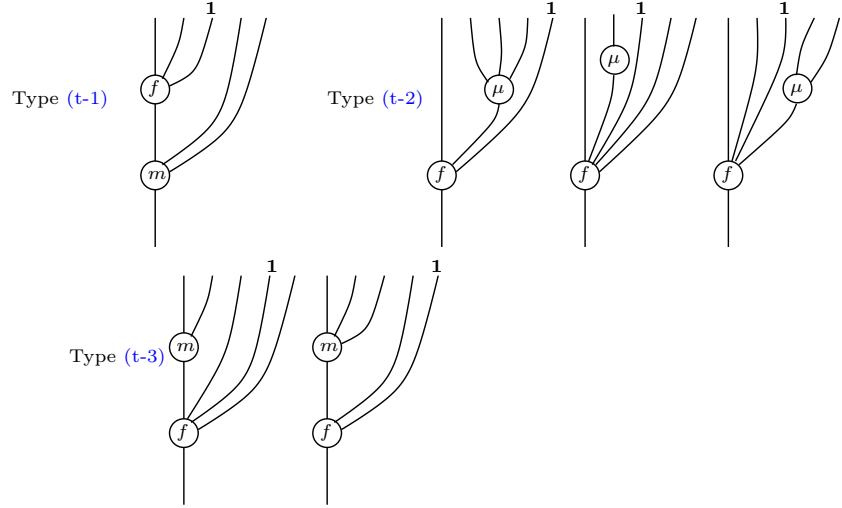}
    \caption[Trees of type \textit{DH}, for relating the unital and non-unital morphism complexes]{{\bf Trees of type $DH$.}
      We have drawn various trees of type $DH$ (as in the proof of Proposition~\ref{prop:mod-mor-cx-unital}),
      with $i=2$, that correspond to trees in Figure~\ref{fig:MorHDTrees}.}
    \label{fig:MorDHTrees}
  \end{figure}
  
  \begin{figure}
  \centering
  \input{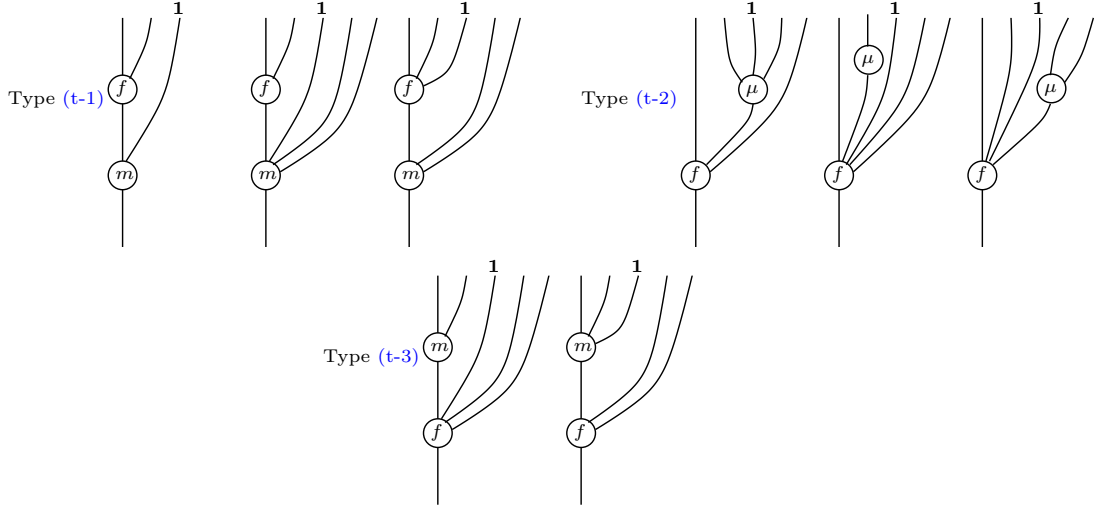}
  \caption[Trees of type \textit{HD}, for relating the unital and non-unital morphism complexes]{{\bf Trees of type $HD$.}
    We have drawn various trees of type $HD$ (as in the proof of Proposition~\ref{prop:mod-mor-cx-unital}),
    with $i=2$.
  The first Type~\ref{Nmult} tree does not have a corresponding $DH$ pair; but it is clearly
  in ${\mathcal F}_2$. The first of the type~\ref{Cmult} trees is a leftover (it does not cancel with the corresponding
  $HD$ tree).}
  \label{fig:MorHDTrees}
  \end{figure}
  
  The $DH$-trees of Type~\ref{Mmult} where the $m_k^{\cModule}$-labelled vertex has
  $k=1$ cancel against the corresponding $HD$-trees of Type~\ref{Mmult}.
  The remaining trees of Type~\ref{Mmult} are in ${\mathcal
    F}_{i+1}$: for some trees, this follows from the fact that
  $\cModule$ is strictly unital; for others, it follows from the fact
  that $f\in {\mathcal F}_i$. Similarly, 
  the Type~\ref{Cmult}-leftovers are typically in ${\mathcal F}_{i+1}$,
  again either because $\Alg$ is strictly unital 
  (provided that the $\mu_k$-labelled vertex has $k>2$)
  or because $f\in {\mathcal F}_i$
  (when the $i^{th}$ input is channelled directly into the $f_\ell$ vertex).
  There are three remaining Type~\ref{Cmult}-leftovers which have been unaccounted for:
  the two $HD$ trees where the $i^{th}$ vertex is channeled into a $\mu_2$-operation,
  and the $DH$ tree where the $i^{th}$ vertex is channeled into a $\mu_2$-operation.
  The first two terms correspond to 
  \begin{align*} (x,a_1,\dots,a_{n-1})&\to f_n(x,a_1,\dots,a_{i-1},\mu_2(\One,a_{i}),a_{i+1},\dots,a_n) \\
   (x,a_1,\dots,a_{n-1})&\to f_n(x,a_1,\dots,a_{i-2},\mu_2(a_{i-1},\One),a_{i},\dots,a_n) 
  \end{align*}
  and their contributions cancel.
  The $DH$ tree where the $i^{th}$ input is fed into a $\mu_2$ 
  corresponds to the operation
  \[ (x,a_1,\dots,a_{n-1})\mapsto f_{n}(x,a_1,\dots,\mu_2(a_{i-1},a_{i}),\One,a_{i+1},\dots,a_{n-1}).\]
  When we set $a_i=\One$, though, the contribution of this tree is the same as the contribution of
  the identity map, i.e.,
  \[ (x,a_1,\dots,a_{n-1})\mapsto f_{n}(x,a_1,\dots,a_{n-1}),\]
  when we set $a_i=1$. This completes the verification that $(\Id+\partial\circ H_i+H_i\circ \partial)({\mathcal F}_i)\subset {\mathcal F}_{i+1}$.

  Clearly, $\suMor(\cModule,\cNodule) =\bigcap_{i=1}^{\infty}
  {\mathcal F}_i$.  Observe that $\partial\circ H_i +
  H_i\circ \partial$ annihilates all $f_k$ with $k<i$.  (This once
  again uses strict unitality of the modules and the algebras.)  Thus,
  we can form the infinite composite 
  \begin{equation}
    \label{eq:DefOfPi}
    \Pi = \dots\circ \Pi_{i}\circ
    \Pi_{i-1}\circ\dots\circ \Pi_1.
  \end{equation}
  The map
  $\Pi\co \Mor(\cModule,\cNodule)\to
  \suMor(\cModule,\cNodule)$ is the homotopy inverse to the inclusion
  map. Finally, the fact that $\Pi_i$ respects bonsai morphisms and does
  not change the bonsai constant implies that $\Pi$ takes bonsai
  morphisms to bonsai morphisms.
\end{proof}

Proposition~\ref{prop:mod-mor-cx-unital} implies, for example, that
every $\Ainf$-homomorphism between strictly unital $\Ainf$-modules $\cModule$ and
$\cNodule$ is homotopic to a strictly unital $\Ainf$-homomorphism.

We now turn to strictly unital algebra homomorphisms.

\begin{proposition}\label{prop:qi-unital}
  If $\Alg$ and $\Blg$ are strictly unital $\Ainf$-algebras and $f\co\Alg\to\Blg$ is a weakly unital 
  quasi-isomorphism then there is a strictly unital quasi-isomorphism
  $g\co\Alg\to\Blg$. If $\Alg$, $\Blg$, and $f$ are bonsai then $g$ may be chosen to be bonsai as well.
\end{proposition}

The proof rests on the following:

\begin{lemma}
  \label{lem:ModifyHomomorphism}
  Let $\Alg$ and $\Blg$ be $\Ainf$-algebras over $\Ground$, and let $f\co
  \Alg\to \Blg$ be an $\Ainf$-algebra homomorphism.
  Let $h\co A^{\kotimes{\Ground} n}\to B$ be a map of degree $n-1$. 
  Then, there is a new homomorphism $F\co \Alg\to \Blg$
  with $F_\ell=f_\ell$ for all $\ell<n$, and
  \begin{equation}
    \label{eq:dhn0}
    F_n = f_n+d h,
  \end{equation}
  where $d h \co A^{\kotimes{\Ground} n}\to B$ is given by
  \[
    d h(a_1,\dots,a_n)=\mu_1^B \circ h(a_1,\dots,a_n) + 
  \sum_{i=1}^n h(a_1,\dots,a_{i-1},\mu_1^A(a_i),a_{i+1},\dots a_n).
  \]
  If moreover $d h=0$, so $F_n=f_n$, then
  \begin{equation} \label{eq:dh1}
  \begin{aligned}
    F_{n+1}&(a_1,\dots,a_{n+1})=f_{n+1}(a_1,\dots,a_{n+1})\\
    &+\sum_{i=1}^n h(a_1,\dots,a_{i-1},\mu_2^{\Alg}(a_i,a_{i+1}),a_{i+2},\dots,a_{n+1}) \\
  &+ \mu_2^{\Blg}(f_1(a_1),h(a_2,\dots,a_{n+1}))+
  \mu_2^{\Blg}(h(a_1,\dots,a_{n}),f_1(a_{n+1})).
\end{aligned}
\end{equation}
  Note that in all cases, $f_1$ and $F_1$ induce the same maps on homology
  $H(A)\to H(B)$.
\end{lemma}

\begin{proof}
  The map $F_{\ell}\co A^{\kotimes{\Ground} \ell}\to B$ is a sum of contributions of trees of three types.
  (Note that the third type is described inductively.)
  \begin{enumerate}[label=(Ft-\arabic*),ref=(Ft-\arabic*)]
  \item \label{type:f}
    A corolla with operation $f$.
  \item \label{type:hd}
    A tree with exactly two internal vertices, where the trunk vertex is labelled by $h$
    and the other vertex is labelled by some $\mu^A_{\ell-n+1}$. 
  \item 
    \label{type:dh}
    A tree with a trunk vertex labelled by $\mu^B$, with a
    distinguished parent vertex labelled by $h$; all $x\geq 0$ parent vertices
    to the left of the distinguished parent  are
    labelled by $f$; and all $y\geq 0$ parent vertices to
    the right of the distinguished parent, are labelled by $F$.
  \end{enumerate}
  See Figure~\ref{fig:Ftrees} for some examples.

\begin{figure}
  \centering
  \input{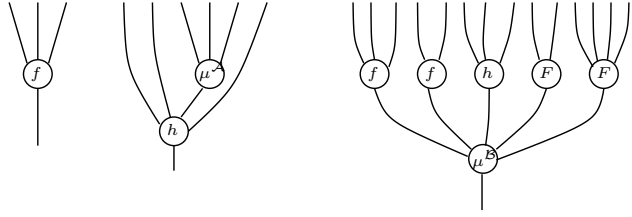}
  \caption[Trees appearing in a particular homomorphism $F$]{{\bf Trees appearing in the homomorphism $F$.}
    We have drawn a tree of Type~\ref{type:f},~\ref{type:hd}, and~\ref{type:dh},
    respectively.}
  \label{fig:Ftrees}
\end{figure}

We claim that $F$ satisfies the $\Ainf$ relation.
This is proved by induction on the number $n$ of inputs in the $\Ainf$
relation. The base case where this number is $1$ is straightforward.

For the inductive step, consider the $\Ainf$ relation
with $n$ inputs. This involves two kinds of trees,
  which we call $\Alg$-trees and  $\Blg$-trees. The $\Alg$-trees have two internal vertices, one of which
  is labelled $\mu^{\Alg}$, which feeds into
  a vertex labelled $F$; the type $\Blg$ trees consist of a number of
  vertices labelled $F$ 
  channeled into a $\mu^{\Blg}$.

  The $\Alg$-trees are classified into the three above types, according to the type of the $F$-labelled vertex.
  Similarly, we classify the $\Blg$-trees into three types, as follows. If
  all the $F$-labelled vertices are of Type~\ref{type:f}, we say that the $\Blg$-tree is of
  Type~\ref{type:f}. Otherwise, the $\Blg$-tree has the type of the leftmost
  $F$-labelled parent of the $\mu^\Blg$ that is not of Type~\ref{type:f}.

  The $\Alg$-trees of Type~\ref{type:f} cancel against the $\Blg$-trees of Type~\ref{type:f} by the $\Ainf$ relation
  for $f$.

  The $\Alg$-trees of Type~\ref{type:hd} cancel each other.  (When one
  $\mu_{\Alg}$ vertex feeds into the other, the
  cancellation follows from the $\Ainf$ relation on $\Alg$; the
  other kinds of terms come in canceling pairs.)

  The $\Alg$-trees of type~\ref{type:dh} where the $\mu^\Alg$-labelled vertex feeds in to the $h$-labelled vertex
  cancel against the $\Blg$-trees of Type~\ref{type:hd}.

  For the remaining cancellations, we apply the $\Ainf$ relation on $\Blg$, to find a relation between
  the various $\Blg$-trees of Type~\ref{type:dh} and two further types of trees: (1) ones  that cancel with the $\Alg$-trees
  of Type~\ref{type:dh} where the $\mu^{\Alg}$ vertex feeds into an $f$-labelled vertex 
  (this cancellation follows from the $\Ainf$ relation
  for $f$) and (2) ones that cancel with the $\Alg$-trees of Type~\ref{type:dh} where the $\mu^{\Alg}$ vertex feeds into
  an $F$-labelled vertex (this cancellation follows from the $\Ainf$ relation
  for $F$, which holds by the inductive hypothesis).
\end{proof}

\begin{proof}[Proof of Proposition~\ref{prop:qi-unital}]
  This is an inductive argument as in the proof of
  Theorem~\ref{thm:UnitalIsUnital}.  Suppose that there are $n$ and
  $k$ so that $f\co \Alg\to \Blg$ satisfies
  \begin{equation}
    \label{eq:PartialUnitalMap}
    f_{m+1}(a_1,\dots,a_{j-1},\unit,a_{j},\dots,a_m)=0 
  \end{equation}
  for all $1< m \leq n$ and for $m=n+1$ and $j\leq k$.
  When $k=n+1$, we will modify $f$ to construct a new
  $\Ainf$-homomorphism $F$ with $F_m=f_m$ for $m\leq n+1$, and so that 
  Equation~\eqref{eq:PartialUnitalMap} holds for $m=n+2$ and $j=1$.
  When $k<n+1$, 
  we will modify $f$ to construct a new $\Ainf$-homomorphism $F$ with 
  $F_m=f_m$ for $m<n+1$, and so that Equation~\eqref{eq:PartialUnitalMap} holds
  for $n+1$ and $k+1$. 

  Start with the first step.  Let
  $h(a_1,\dots,a_n)=f_{n+1}(\unit,a_1,\dots,a_{n})$.  Apply the
  $\Ainf$-homomorphism relation for $f$ with inputs
  $(\unit,a_1,\dots,a_{n})$.  By the inductive hypothesis, the only
  non-vanishing terms are
  $\mu_2^{\Blg}(f_1(\unit),f_n(a_1,\dots,a_n))$ and
  $f_n(\mu_2(\unit,a_1),a_2,\dots,a_n)$, which cancel against each
  other, and the terms $d h$; i.e., we have verified that
  $d h=0$. Thus, we can modify $f$ as in
  Lemma~\ref{lem:ModifyHomomorphism}.  The $\Ainf$-homomorphism
  relation with inputs $(\One,\One,a_2,\dots,a_{n+1})$, together with
  the inductive hypothesis on $n$ implies that, for $F_{n+1}$ as given
  in Equation~\eqref{eq:dh1},
  \[ F_{n+1}(\One,a_2,\dots,a_{n+1})=0.\]

  Similarly, for  the second step, assume that
  Equation~\eqref{eq:PartialUnitalMap} holds for all $m\leq n$ and for
  $m=n+1$ for all $2\leq j\leq k$. Define
  \[ h(a_1,\dots,a_{n})=f(a_1,\dots,a_{k-1},\unit,a_{k},\dots,a_{n}).\]
  The $\Ainf$ relation with inputs
  $a_1,\dots,a_{k-1},\unit,a_{k},\dots,a_{n}$ implies that $d
  h=0$.  Modifying $f$ as in Lemma~\ref{lem:ModifyHomomorphism}, it is
  straightforward to check that the new map $F_{n+1}$ as specified in
  Equation~\eqref{eq:dh1} satisfies
  $F_{n+1}(a_1,\dots,a_{k-1},\unit,a_{k+1},\dots,a_{n+1})=0$.

  The statement about the bonsai case follows by the same reasoning as
  in the proof of Theorem~\ref{thm:UnitalIsUnital}.
\end{proof}

\begin{remark}
  There is a notion in between weak and strict unitality, called
  \emph{homotopy unitality.} We will need this notion only in the
  weighted case, so we spell it out there, in Section~\ref{sec:hu}.
\end{remark}

%%% Local Variables: 
%%% mode: latex
%%% TeX-master: "AbstractDiagonal"
%%% TeX-command-extra-options: "-shell-escape"
%%% End: 

\section{Weighted \texorpdfstring{$\Ainf$}{A-infinity} algebras and weighted trees}
\label{sec:wAinfty}
\subsection{Weighted \texorpdfstring{$\Ainf$}{A-infinity} algebras and the weighted trees complex}\label{sec:wAlgs}
Convention~\ref{conv:Ring} regarding ground rings continues to be in force in this section.

\begin{definition}
  \label{def:wAinfty}
  A {\em weighted $\Ainf$-algebra} or {\em $w$-algebra} $\wAlg$
  over $\Ground$ consists of:
  \begin{itemize}
  \item a graded $\Ground$-bimodule $A$, and
  \item for each pair of integers $n,w\geq 0$ a linear map
    \[
      \mu^w_n\co \overbrace{A\kotimes{\Ground}\dots\kotimes{\Ground} A}^{n}\to  A\grs{2-n-w\kappa}
    \]
    of $\Ground$-bimodules
  \end{itemize}
  satisfying the structure equation
  \begin{equation}
    \label{eq:wAlg-StructuralEquation}
    \sum_{\substack{1\leq i\leq j\leq n+1\\ u+v=w}}\mu^u_{n-j+i+1}(a_1,\dots,a_{i-1},\mu_{j-i}^v(a_i,\dots,a_{j-1}),a_{j},\dots,a_n)
    = 0
  \end{equation}
  for each pair of non-negative integers $(n,w)$.

  The index $w$ is called the \emph{weight} of the operation
  $\mu^w_n$. The integer $\kappa$ is called the \emph{weight grading}.

  The $w$-algebra is \emph{uncurved} if $\mu^0_0=0$.
\end{definition}

Recall that the grading shift in Definition~\ref{def:wAinfty} means
that $\mu_n^w$ decreases the grading by $2-n-w\kappa$.

The maps $\mu_n^0$ give $A$ the structure of an $\Ainf$-algebra,
which we call the \emph{undeformed $\Ainf$-algebra of $\wAlg$}.  The
undeformed $\Ainf$-algebra of $\wAlg$ is uncurved exactly when $\wAlg$
is uncurved.

\begin{convention}\label{conv:uncurved}
  We will henceforth consider only uncurved $w$-algebras.
\end{convention}

A $w$-algebra gives rise to a curved $\Ainf$-algebra structure on
$A\llbracket t\rrbracket=A\rotimes{\Ring}\Ring\llbracket t\rrbracket$,
the formal power series in $t$ with coefficients in~$A$,
with operations
\begin{equation}
  \label{eq:AssociatedCurvedAlgebra}
\mu_n(a_1,\dots,a_n) = \sum_{w=0}^\infty \mu^w_n(a_1,\dots,a_n)t^w.
\end{equation}
(The variable~$t$ has grading~$-\kappa$.)
The condition that a $w$-algebra is uncurved is the condition that the 
associated curved algebra has curvature 
in $t\cdot A\llbracket t\rrbracket$.

\begin{figure}
  \centering
  %Font is 18 point.
  \includegraphics[scale=.667]{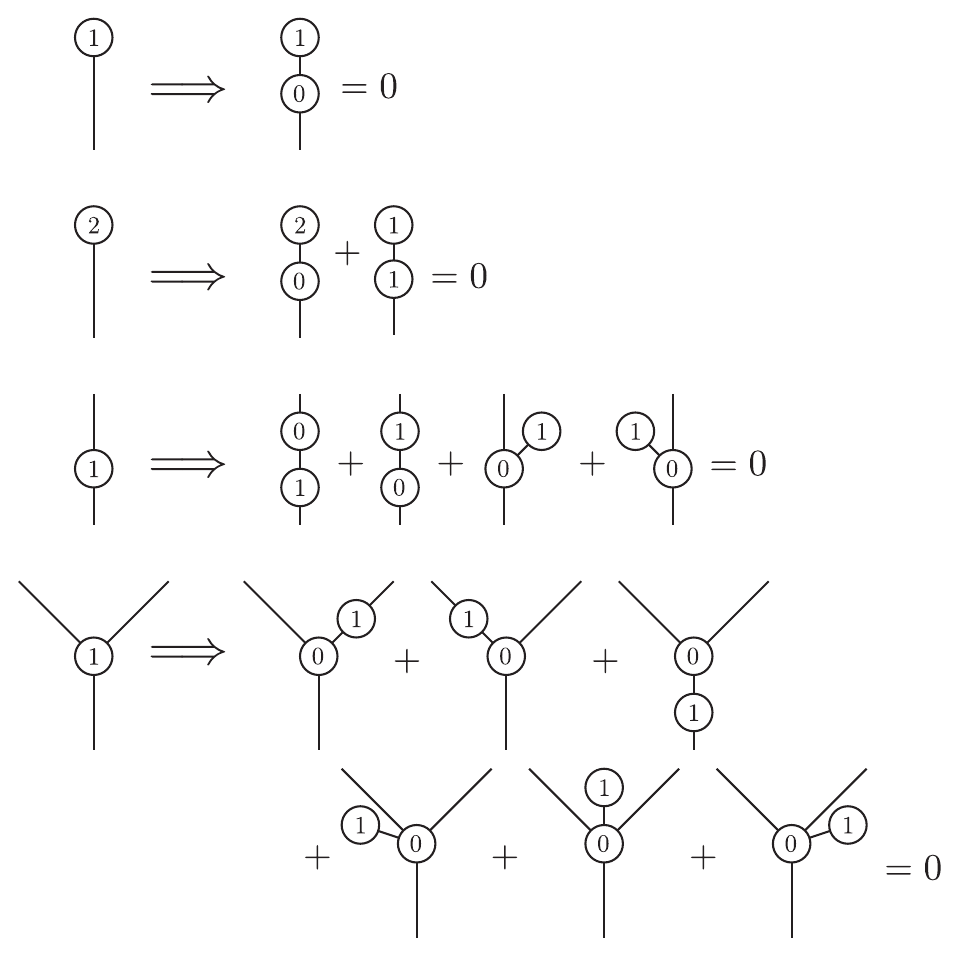}
  \caption[Some weighted A-infinity relations]{\textbf{Uncurved weighted $\Ainf$ relations.} The cases
    $(n,w)=(0,1)$, $(0,2)$, $(1,1)$, and $(2,1)$ are shown. This is the
    uncurved case, i.e., $\mu^0_0=0$.}
  \label{fig:weighted-Ainf-rel}
\end{figure}

The operations $\mu^w_n$ in a $w$-algebra are conveniently expressed
as planar trees with one internal vertex, $n$ inputs, one output, and a label
by the integer $w$. The structure equations can be thought of as
indexed by such trees; their terms consist of all the ways of
inserting an edge in tree and splitting the weight. See
Figure~\ref{fig:weighted-Ainf-rel}.

\begin{remark}\label{remark:weighted-is-deformation}
  A weighted $\Ainf$-algebra $\wAlg=(A,\{\mu_n^w\})$ is the same as a
  one-parameter deformation of $(A,\{\mu_n^0\})$, in the sense, for
  instance, of Seidel~\cite[Section 3b]{Seidel15:quartic}, except that
  Seidel focuses on uncurved deformations ($\mu_0^w=0$) and only
  considers the case $\kappa\neq 0$.

  Weighted $\Ainf$-algebras are also a special case of the weakly curved
  $\Ainf$-algebras studied by Positselski~\cite{Positselski:curved} and of the
  gapped, filtered $\Ainf$-algebras introduced by
  Fukaya-Oh-Ohta-Ono~\cite[Chapter 3]{FOOO1}.
\end{remark}

By a \emph{marked tree} we
mean a planar, rooted tree $T$ together with a subset of the leaves of
$T$, which we call the \emph{inputs} of $T$ and another leaf of $T$,
the \emph{output} of $T$. Leaves of $T$ which are not inputs or
outputs are \emph{popsicles}.
We call a vertex $v$ of $T$ {\em internal} if
$v$ is not an input or output (i.e., is a popsicle or has valence $>1$).

Given a vertex $v$ in a tree $T$ and a non-negative integer $w(v)$, we
call $(v,w(v))$ \emph{stable} if either $v$ has valence $3$ or larger
or $w(v)>0$ (or both). A \emph{weighted tree} is a marked tree $T$
together with a weight function $w$ from the internal vertices of $T$
to the non-negative integers. A {\em stably-weighted tree} is a
weighted tree $T$ with the property that for all internal vertices $v$
of $T$, $(v,w(v))$ is stable. The {\em total weight} of a
weighted tree is the sum of the weights $w(v)$ of all the
internal vertices $v$ of $T$. Let $\wTrees{n}{w}$ denote the set of
stably-weighted trees with $n$ inputs and total weight $w$.

\begin{figure}
  \centering
  % Font is 18 point.
  \includegraphics[scale=.667]{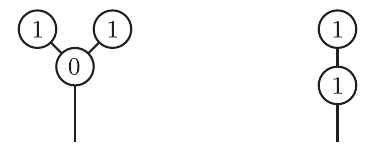}
  \caption[A stably-weighted tree and two contractions of it]{\textbf{A stably-weighted tree and two contractions of it.} Left: a stably-weighted tree. Right: the result of contracting either of the internal edges in the tree.}
  \label{fig:notSimpCx}
\end{figure}

Let $(S,e)$ be a pair consisting of a stably-weighted tree $S$ and an
edge $e$ connecting two internal vertices in $S$.  We can contract the
edge to form a new stably-weighted tree $T$ with one fewer internal
vertex than $S$. The contracted edge gives rise to a new vertex in
$T$ whose weight is the sum of the weights of the two vertices of
$e$. In this case, we say that $S$ is obtained from $T$ by \emph{inserting
an edge}, and we call the pair $(S,e)$ an {\em edge expansion} of
$T$.
Note that if $(S,e)$ is an edge expansion of $T$, then $e$ is typically,
but not always, uniquely determined by $(S,T)$; see, for example,
Figure~\ref{fig:notSimpCx}.

Fix a stably-weighted tree $T$ with $n$ inputs, total weight $w$, and
$v$ internal vertices.  The \emph{dimension} of $T$ is given by
\begin{equation}\label{eq:w-tree-dim}
\dim(T)= n+2w-v-1.
\end{equation}
(See Section~\ref{sec:Associaplex} for justification of the term
``dimension''.) Note that 
\[ 
0 \leq \dim(T)\leq n+2w-2. 
\]

\begin{definition}
  Fix non-negative integers $n$ and $w$.  The {\em $n$-input, weight
    $w$ weighted trees complex $\wTreesCx{n}{w}$} is the chain
  complex with basis $\wTrees{n}{w}$ the set of stably-weighted trees
  with $n$ inputs and total weight $w$, equipped with a grading
  induced by the dimension $\dim(T)$, i.e.,
  \[\wTreesCx[k]{n}{w}=\Ring\langle\{T\mid \dim(T)=k\}\rangle.\]  The
  differential is specified by declaring that for a stably-weighted tree $T$, 
  \[ 
  \partial T = \sum_{(S,e)\text{ an edge expansion of } T} S.
  \]
\end{definition}

\begin{lemma}\label{lem:wTrees-is-cx}
  The $n$-input, weight $w$ weighted trees complex $\wTreesCx{n}{w}$
  is a chain complex, and the differential decreases the grading by
  one. Further, the composition map $\circ_i\co
  \wTreesCx{n_1}{w_1}\rotimes{\Ring}
  \wTreesCx{n_2}{w_2}\to\wTreesCx{n_1+n_2-1}{w_1+w_2}$ is a chain map.
\end{lemma}
\begin{proof}
  To see that $\bdy^2=0$, it is easier to consider the dual complex
  $(\wTreesCx{n}{w})^*$, whose basis is the stably-weighted trees
  and whose (co)differential $\delta(T)$ is the sum of all ways of contracting
  an internal edge (edge not adjacent to an input or the output) in~$T$.
  Then $\delta^2(T)$ is the sum of all ways of contracting two
  edges, and since the contractions can be performed in either order,
  $\delta^2=0$. Next, it is immediate that the differential decrease
  the grading by one, since it preserves the total weight, and
  increases the number of internal vertices by one. The fact that the
  composition maps are chain maps is also immediate from the
  definitions.
\end{proof}

\begin{example}
  The $n$-input weight $0$ weighted trees complex is $\cellC{*}(K_n)$,
  the cellular chain complex of the associahedron $K_{n}$.
\end{example}

\begin{example}
  The chain complex $\wTreesCx{0}{1}$ is one-dimensional, with
  vanishing differential.
\end{example}

\begin{example}
  The chain complex $\wTreesCx{0}{2}$ is generated by three trees, in
  dimensions $2$, $1$, and $0$. The differential of the
  $2$-dimensional tree is the one-dimensional tree, and the
  differential of the one-dimensional tree is zero, because the
  $0$-dimensional tree occurs twice in the boundary of the
  one-dimensional tree. (See Figure~\ref{fig:notSimpCx}.)  It follows
  that the stably-weighted trees complex is not the chain complex
  associated to any polyhedral complex.
\end{example}

Just as with unweighted trees, there are stacking operations for
weighted trees: given $S\in\wTrees{m}{v}$ and $T\in \wTrees{n}{w}$ and
an integer $1\leq i\leq n$ let $T\circ_i S$ be the result of gluing
the output (root) of $S$ to the $i\th$ input of $T$ (and
forgetting the resulting $2$-valent vertex). The operation $\circ_i$
induces a chain map
\begin{align*}
\phi_{i,j,n;v,w}\co \wTreesCx{j-i+1}{v}\rotimes{\Ring}\wTreesCx{n+i-j}{w} &\to \wTreesCx{n}{v+w}\\
S\otimes T&\mapsto T\circ_i S
\end{align*}
for any $1\leq i\leq n$ and $i-1\leq j\leq n$.  (These are the same as
the composition operations from Section~\ref{sec:associahedron}, and
in particular Figure~\ref{fig:SomeStackingTrees} except that $S$ and
$T\circ S$ are allowed to be $0$ or $1$ input trees.)

A \emph{weighted corolla} is a tree with $n$ inputs and one internal
vertex, of some weight $w$, denoted $\wcorolla{w}{n}$. (The corolla $\corolla{n}$ from
Section~\ref{sec:associahedron} is $\wcorolla{n}{0}$.)  All weighted
trees can be built by stacking weighted corollas.

Turning to the connection with weighted $\Ainf$-algebras, suppose
$A$ is a \dg $\Ground$-bimodule with differential $\mu_1^0$ and that we are
also given maps $\mu_{n}^w\co A^{\kotimes{\Ground} n}\to A\grs{2-n-w\kappa}$ for $n\geq 0$,
$w\geq 0$, $(n,w)\not\in\{(0,0),(1,0)\}$.
Given a stably-weighted tree
$T$ with $n$ inputs and total weight $W$ there is an associated map
$\mu(T)\co A^{\kotimes{\Ground} n}\to A\grs{-\dim(T)+(2-\kappa)W}$ obtained by replacing each vertex
with weight $w$ and valence $n+1$ (including the popsicles) with
$\mu_n^w$ and
composing according to the edges of the tree. Extending linearly, we
obtain an $\Ring$-module map $\mu\co \wTreesCx{n}{W}\to \Mor(A^{\kotimes{\Ground} n},A\grs{(2-\kappa)W})$.
\begin{lemma}\label{lem:wAlg-reinterp}
  The collection $\{\mu_n^w\}$ satisfies the weighted $\Ainf$ relations
  if and only if for each $(n,w)$ the map
  $\mu\co\wTreesCx{n}{w}\to \Mor(A^{\kotimes{\Ground} n},A\grs{(2-\kappa)w})$ is a chain
  map.
\end{lemma}
\begin{proof}
  First, suppose the $\mu_n^w$ satisfy the weighted $\Ainf$ relations;
  we will verify that $\mu$ is a chain map.  Since
  $\mu(T\circ_i S)=\mu(T)\circ_i \mu(S)$ and composition of trees and
  composition in the morphism complex are both chain maps, it suffices
  to verify that for any corolla $\wcorolla{n}{w}$,
  \[
    \bdy\circ (\mu(\wcorolla{n}{w}))+\mu(\wcorolla{n}{w})\circ\bdy+\mu(\bdy(\wcorolla{n}{w}))=0.
  \]
  The first term corresponds to the term in
  Equation~(\ref{eq:wAlg-StructuralEquation}) with $u=0$, $j=n+1$, and
  $i=1$. The second term corresponds to the terms in
  Equation~(\ref{eq:wAlg-StructuralEquation}) with $v=0$ and
  $j=i+1$. The third term corresponds to the terms in
  Equation~(\ref{eq:wAlg-StructuralEquation}) in which both
  multiplications are stable, i.e., with ($r>0$ or $j-i>1$) and ($v>0$
  or $n-j+i+1>1$); these are all the remaining terms.

  Conversely, if $\mu$ is a chain map then considering
  $\bdy(\mu(\wcorolla{n}{w}))$ and $\mu(\bdy(\wcorolla{n}{w}))$ and
  making the same identification of terms gives
  Equation~(\ref{eq:wAlg-StructuralEquation}).
\end{proof}

We also have a notion of boundedness for $w$-algebras.

\begin{definition}\label{def:w-bonsai}
  A $w$-algebra $\wAlg=(A,\{\mu_n^w\})$ is \emph{bonsai} if there is an
  integer $N$ so that for any stably-weighted tree $T$ with
  $\dim(T)>N$, we have $\mu(T)=0$.
\end{definition}

\subsection{Weighted algebra maps and weighted transformation trees}
\begin{definition}\label{def:w-TransformationTree}
  A \emph{weighted transformation tree} consists of a weighted tree
  $T$ and a coloring of each edge of $T$ as either \emph{red} or
  \emph{blue}, so that:
  \begin{enumerate}
  \item The edges adjacent to input leaves are red.
  \item The edge adjacent to the output leaf is blue.
  \item For each vertex $v$, all of the inputs of $v$ have the same color (red or blue).
  \item If a vertex $v$ has a red output, then all of the inputs of $v$ are red;
    if a vertex has blue inputs, then its output is also blue.
  \item Call a vertex \emph{red} (respectively \emph{blue}) if all of
    its inputs and its output are red (respectively blue), and
    \emph{purple} if its inputs are red but its output is blue. We
    require that tree be \emph{stable} in the sense that for each
    vertex $v$, either:
    \begin{enumerate}
    \item the weight of $v$ is strictly positive, or
    \item $v$ has valence at least $3$ (i.e., at least two inputs), or
    \item $v$ is purple and has valence $2$.
    \end{enumerate}
  \end{enumerate}
  (Cf.\ Definition~\ref{def:TransformationTree}.)

  Given a (red) weighted tree $S$, an $n$-input weighted transformation
  tree $T$, and an integer $i$, $1\leq i\leq n$, we can form a new
  weighted transformation tree $T\circ_i S$ by connecting the output of
  $S$ to the $i\th$ input of $T$. Also, given an $n$-input (blue)
  weighted tree $S$ and weighted transformation trees $T_1,\dots,T_n$,
  we can form a new weighted transformation tree $S\circ(T_1,\dots,T_n)$
  by connecting the output of each $T_i$ to the $i\th$ input of $S$.
  These operations induce composition maps
  \begin{align}
    \circ_i\co \wTransCx{n_1}{w_1}\rotimes{\Ring}
    \wTreesCx{n_2}{w_2}&\to\wTransCx{n_1+n_2-1}{w_1+w_2},\label{eq:w-trans-compo-1}\\
    \circ\co \wTreesCx{n}{w}\rotimes{\Ring}\wTransCx{m_1}{v_1}\rotimes{\Ring}\cdots\rotimes{\Ring} \wTransCx{m_n}{v_n}
                       &\to\wTransCx{m_1+\cdots+m_n}{v_1+\cdots+v_n+w}\label{eq:w-trans-compo-2}.
  \end{align}
  As usual, given $T\in \wTransCx{n_1}{w_1}$ and
  $S\in \wTreesCx{n_2}{w_2}$, we define
  \[
    T\circ S=\sum_{i=1}^{n_1} T\circ_i S.
  \]
  
  The \emph{weighted multiplex} $\wTransCx{n}{w}$ is the chain complex
  spanned by weighted transformation trees with $n$ inputs and total
  weight $w$, with the following differential. As usual, the
  differential satisfies the Leibniz rule with respect to composition,
  so it suffices to define the differential of a red, purple, or blue
  corolla. The differential of a red corolla $\wrcorolla{n}{w}$ or
  blue corolla $\wbcorolla{n}{w}$ is as in the weighted trees
  complex. The differential of a purple corolla $\wpcorolla{n}{w}$ is
  given by
  \[
    \bdy(\wpcorolla{n}{w})=\sum_{\substack{k+\ell=n+1\\w_1+w_2=w}}\wpcorolla{k}{w_1}\circ\wrcorolla{\ell}{w_2}
    + \hspace{-1em}\sum_{\substack{k_1+\cdots+k_\ell=n\\v+w_1+\cdots+w_\ell=w}}\hspace{-1em}
    \wbcorolla{\ell}{v}\circ(\wpcorolla{k_1}{w_1},\dots,\wpcorolla{k_\ell}{w_\ell}).
  \]
  Note that in the special case $n=0$, there is a term in the second
  sum of the form $\wbcorolla{0}{w}$, corresponding to $\ell=0$.

  Define the dimension to be $\dim(T)= n+2w-r-b-1$, where $r$ and $b$
  are the number of red and blue vertices respectively.
\end{definition}

\begin{lemma}\label{lem:wTrans-is-cx}
  The $n$-input, weight $w$ multiplex $\wTransCx{n}{w}$ is a chain
  complex, and the differential decreases the grading by one. Further,
  the composition maps~\eqref{eq:w-trans-compo-1}
  and~\eqref{eq:w-trans-compo-2}  are chain maps.
\end{lemma}
\begin{proof}
  As in Lemma~\ref{lem:wTrees-is-cx}, to see that $\bdy^2=0$ it is
  easier to consider the dual complex, which has the same basis as
  $\wTransCx{n}{w}$ but where the differential is the sum of all ways of:
  \begin{enumerate}[label=(\arabic*)]
  \item\label{item:cont-red-red} contracting an edge between two red vertices,
  \item\label{item:cont-blue-blue} contracting an edge between two blue vertices,
  \item\label{item:cont-purp-red} contracting an edge between a purple vertex and a red vertex,
    or
  \item\label{item:cont-blue-purp} contracting all the edges into a blue vertex $v$ if all of the
    vertices above $v$ are purple.
  \end{enumerate}
  Terms in $\delta^2$ coming from two contractions of the first three
  kinds cancel in pairs, corresponding to the two possible
  orders. It remains to understand terms where one or both are of the
  fourth type. Several of these cases commute, so cancel in pairs:
  \begin{itemize}
  \item Contractions of type~\ref{item:cont-red-red}
    and~\ref{item:cont-blue-purp}.
  \item Contractions of type~\ref{item:cont-purp-red}
    and~\ref{item:cont-blue-purp}.
  \item Pairs of a contraction of type~\ref{item:cont-blue-blue} and a
    contraction of type~\ref{item:cont-blue-purp} corresponding
    to three distinct blue vertices.
  \item Pairs of contractions of type~\ref{item:cont-blue-purp}
    corresponding to a pair of blue vertices $v_1$, $v_2$ where $v_1$
    is neither above nor below $v_2$.
  \end{itemize}
  The remaining terms are:
  \begin{itemize}
  \item Contraction of an edge connecting two blue vertices $v_1$ and
    $v_2$, giving a new vertex $v$, and then all of the edges above
    $v$ (which must lead to purple vertices). For this to be sensible,
    all of the edges into $v_1$ come from purple vertices, as do all
    of the edges into $v_2$ except the edge from $v_1$.
  \item Contraction of all the edges above a vertex $v_1$ and then all
    the edges above a vertex $v_2$ below $v_1$. For this to be
    sensible, again, all of the edges into $v_1$ come from purple
    vertices, as do all of the edges into $v_2$ except the edge from
    $v_1$.
  \end{itemize}
  These terms, again, cancel in pairs.

  The facts that the differential decreases the grading by $1$ and the
  composition maps are chain maps are immediate from the definitions.
\end{proof}

\begin{definition}\label{def:wAlg-homo}
  Fix weighted $\Ainf$-algebras $\wAlg$ and $\wBlg$ over $\Ground$ with the same
  weight grading $\kappa$.
  A \emph{homomorphism} of weighted $\Ainf$-algebras
  $f\co \wAlg\to\wBlg$ is a homomorphism of curved $\Ainf$-algebras
  whose $0$-input component~$f_0$ lies in
  $t B\llbracket t\rrbracket\subset B\llbracket
  t\rrbracket$.
  Equivalently, a homomorphism is a sequence of maps
  $f_n^w\co A^{\kotimes{\Ground} n}\to B\grs{1-n-\kappa w}$, $n,w=0,\dots,\infty$,
  $(n,w)\neq(0,0)$, satisfying the structure equation
  \begin{multline*}
    \sum_{\substack{1\leq i\leq j\leq n+1\\u+v=w\\2v+j-i\ge 2}}f_{n-j+i+1}^u(a_1,\dots,a_{i-1},\mu_{j-i}^v(a_i,\dots,a_{j-1}),a_{j},\dots,a_n)
    \\
    +\!\!\!\!\sum_{\substack{0=i_0\leq i_1\leq \dots\leq i_{k+1}=n\\
        u_1+\dots+u_k+v=w\\
        2u_j+i_j-i_{j-1} \ge 1,\,\, 2v+k \ge 2
    }}
    \!\!\!\!\!\mu_k^v\bigl(f^{u_1}_{i_1-i_0}(a_1,\dots,a_{i_1}),f^{u_2}_{i_2-i_1}(a_{i_1+1},\dots,a_{i_2}),\dots,f^{u_k}_{n-i_k}(a_{i_k+1},\dots,a_n)\bigr)=0
  \end{multline*}
  for each $(n,w)$.
  
Homomorphisms compose in the obvious way, generalizing the case of
$\Ainf$-algebra homomorphisms. It is clear that the composition of two
weighted $\Ainf$-algebra homomorphisms is again a homomorphism.
The \emph{identity homomorphism} of $\wAlg$, $\Id$, is defined by
\[
  \Id_n^w=
  \begin{cases}
    \Id_{A} & (n,w)=(1,0)\\
    0 & \text{else.}
  \end{cases}
\]
An invertible homomorphism is an \emph{isomorphism}.

Given a weighted $\Ainf$-algebra $\wAlg$, the homology $H_*(\wAlg)$ is
the homology of $A$ with respect to $\mu_1^0$. A homomorphism is
a \emph{quasi-isomorphism} if the induced map $H_*(\wAlg)\to H_*(\wBlg)$
is an isomorphism or, equivalently, if the induced map of undeformed
$\Ainf$-algebras is a quasi-isomorphism.
\end{definition}

\begin{lemma}\label{lem:wAlg-iso-is}
  A weighted algebra homomorphism $f\co \wAlg\to\wBlg$ is an
  isomorphism if and only if $f^0_1$ is an isomorphism.
\end{lemma}
\begin{proof}
  From the same proof as in the unweighted case \cite[Lemma 2.1.14,
  \textit{inter alia}]{LOT2}, but inducting first on the weight and
  then on the number of inputs, we get a left inverse $g_L$ to $f$,
  and a right inverse $g_R$ to $f$. \emph{A priori}, these maps may
  not be homomorphisms; however, associativity of the formula for
  composition shows that $g_L=g_L\circ f\circ g_R=g_R$, so $g=g_L=g_R$
  is a two-sided inverse.

  The map $g$ is a homomorphism if
  \begin{equation}\label{eq:wAlg-iso-is-1}
    \tikzsetnextfilename{wAlg-iso-is-1}
    \sum_{u_1+\cdots+u_n+v=w}
    \mathcenter{
      \begin{tikzpicture}[smallpic]
        \node at (0,0) (tl) {};
        \node at (2,0) (tr) {};
        \node at (0,-1) (g1) {$g^{u_1}$};
        \node at (2,-1) (g2) {$g^{u_k}$};
        \node at (1,-1) (dots) {$\cdots$};
        \node at (1,-2) (mu) {$\mu^v$};
        \node at (1,-3) (bc) {};
        \draw[tbb] (tl) to (g1);
        \draw[tbb] (tr) to (g2);
        \draw[alga] (g1) to (mu);
        \draw[alga] (g2) to (mu);
        \draw[alga] (mu) to (bc);
      \end{tikzpicture}}
    +
    \tikzsetnextfilename{wAlg-iso-is-2}
    \sum_{u+v=w}
    \mathcenter{
      \begin{tikzpicture}[smallpic]
        \node at (0,0) (tl) {};
        \node at (1,0) (tc) {};
        \node at (2,0) (tr) {};
        \node at (1,-1) (mu) {$\mu^u$};
        \node at (1,-2) (g) {$g^{v}$};
        \node at (1,-3) (bc) {};
        \draw[tbb] (tl) to (g);
        \draw[tbb] (tc) to (mu);
        \draw[tbb] (tr) to (g);
        \draw[blga] (mu) to (g);
        \draw[alga] (g) to (bc);
      \end{tikzpicture}}
    =0
  \end{equation}
  (as maps from the tensor algebra on $A$ to $B$). To keep notation
  simple, we will suppress the weights in the rest of the
  proof. Pre-composing the left side by $f$ gives
  \[
    \tikzsetnextfilename{wAlg-iso-is-3}
    \mathcenter{
      \begin{tikzpicture}[smallpic]
        \node at (0,0) (t1) {};
        \node at (2,0) (t2) {};
        \node at (4,0) (t3) {};
        \node at (6,0) (t4) {};
        \node at (0,-1) (f1) {$f$};
        \node at (1,-1) (dots1) {$\cdots$};
        \node at (2,-1) (f2) {$f$};
        \node at (4,-1) (f3) {$f$};
        \node at (5,-1) (dots2) {$\cdots$};
        \node at (6,-1) (f4) {$f$};
        \node at (1,-2) (g12) {$g$};
        \node at (3,-2) (dots3) {$\cdots$};
        \node at (5,-2) (g34) {$g$};
        \node at (3,-3) (mu) {$\mu$};
        \node at (3,-4) (bc) {};
        \draw[taa] (t1) to (f1);
        \draw[taa] (t2) to (f2);
        \draw[taa] (t3) to (f3);
        \draw[taa] (t4) to (f4);
        \draw[blga] (f1) to (g12);
        \draw[blga] (f2) to (g12);
        \draw[blga] (f3) to (g34);
        \draw[blga] (f4) to (g34);
        \draw[alga] (g12) to (mu);
        \draw[alga] (g34) to (mu);
        \draw[alga] (mu) to (bc);
      \end{tikzpicture}
    }
    +
    \tikzsetnextfilename{wAlg-iso-is-4}
    \mathcenter{
      \begin{tikzpicture}[smallpic]
        \node at (0,0) (t1) {};
        \node at (2,0) (t2) {};
        \node at (4,0) (t3) {};
        \node at (6,0) (t4) {};
        \node at (0,-1) (f1) {$f$};
        \node at (1,-1) (dots1) {$\cdots$};
        \node at (2,-1) (f2) {$f$};
        \node at (3,-1) (dots3) {$\cdots$};
        \node at (4,-1) (f3) {$f$};
        \node at (5,-1) (dots2) {$\cdots$};
        \node at (6,-1) (f4) {$f$};
        \node at (3,-2) (mu) {$\mu$};
        \node at (3,-3) (g) {$g$};
        \node at (3,-4) (bc) {};
        \draw[taa] (t1) to (f1);
        \draw[taa] (t2) to (f2);
        \draw[taa] (t3) to (f3);
        \draw[taa] (t4) to (f4);
        \draw[blga] (f1) to (g);
        \draw[blga] (f2) to (mu);
        \draw[blga] (f3) to (mu);
        \draw[blga] (f4) to (g);
        \draw[blga] (mu) to (g);
        \draw[alga] (g) to (bc);
      \end{tikzpicture}
    }
    =
    \mathcenter{
      \tikzsetnextfilename{wAlg-iso-is-5}
      \begin{tikzpicture}[smallpic]
        \node at (0,0) (tc) {};
        \node at (0,-1) (mu) {$\mu$};
        \node at (0,-2) (bc) {};
        \draw[taa] (tc) to (mu);
        \draw[alga] (mu) to (bc);
    \end{tikzpicture}
  }
  +
    \mathcenter{
      \tikzsetnextfilename{wAlg-iso-is-6}
      \begin{tikzpicture}[smallpic]
        \node at (0,0) (tl) {};
        \node at (.5,0) (tcl) {};
        \node at (1,0) (tc) {};
        \node at (1.5,0) (tcr) {};
        \node at (2,0) (tr) {};
        \node at (1,-1) (mu) {$\mu$};
        \node at (0,-2) (fl) {$f$};
        \node at (1,-2) (fc) {$f$};
        \node at (2,-2) (fr) {$f$};
        \node at (1,-3) (g) {$g$};
        \node at (1,-4) (bc) {};
        \draw[taa] (tl) to (fl);
        \draw[taa] (tr) to (fr);
        \draw[taa] (tc) to (mu);
        \draw[taa] (tcl) to (fc);
        \draw[taa] (tcr) to (fc);
        \draw[alga] (mu) to (fc);
        \draw[blga] (fl) to (g);
        \draw[blga] (fc) to (g);
        \draw[blga] (fr) to (g);
        \draw[alga] (g) to (bc);
      \end{tikzpicture}
  }
    =
    \mathcenter{
      \tikzsetnextfilename{wAlg-iso-is-7}
      \begin{tikzpicture}[smallpic]
        \node at (0,0) (tc) {};
        \node at (0,-1) (mu) {$\mu$};
        \node at (0,-2) (bc) {};
        \draw[taa] (tc) to (mu);
        \draw[alga] (mu) to (bc);
    \end{tikzpicture}
  }+
    \mathcenter{
      \tikzsetnextfilename{wAlg-iso-is-8}
      \begin{tikzpicture}[smallpic]
        \node at (0,0) (tc) {};
        \node at (0,-1) (mu) {$\mu$};
        \node at (0,-2) (bc) {};
        \draw[taa] (tc) to (mu);
        \draw[alga] (mu) to (bc);
    \end{tikzpicture}
  }=0.
  \]
  Since pre-composition with $f$ is an isomorphism with inverse pre-composition by $g$, this implies Equation~\eqref{eq:wAlg-iso-is-1}, as desired.
\end{proof}

Given a homomorphism of weighted $\Ainf$-algebras $f\co\wAlg\to\wBlg$
and a weighted transformation tree $T$ with $n$ inputs and total
weight $W$ there is an induced map
$F(T)\co A^{\kotimes{\Ground} n}\to B\grs{-\dim(T)+(2-\kappa)W}$.  Thus, we
have a map
\[
  F\co \wTransCx{n}{w}\to \Mor(A^{\kotimes{\Ground} n},B\grs{(2-\kappa)w}).
\]
\begin{lemma}
  The collection $\{f_{n}^w\}$ satisfy the weighted $\Ainf$-homomorphism
  relations if and only if for each $(n,w)$ the induced map $F$ is a
  chain map.
\end{lemma}
\begin{proof}
  The proof is similar to the proof of
  Lemma~\ref{lem:wAlg-reinterp} or~\ref{lem:A-inf-alg-map-is} and is left to the reader.
\end{proof}

\begin{definition}
  A homomorphism $f \co \wAlg \to \wBlg$ between $w$-algebras is
  \emph{bonsai} if there is a number~$N$ so that for any weighted
  transformation tree~$T$ with $\dim(T) > N$, we have $F(T) = 0$.
\end{definition}

\subsection{Units in weighted \texorpdfstring{$\Ainf$-}{A-infinity-}algebras}
\label{sec:w-units}

Just as in the unweighted case, there are several notions of unitality
one can consider for weighted algebras.  We will discuss two.

\begin{definition}\label{def:w-weakly-unital}
  A $w$-algebra $\wAlg$ is \emph{weakly unital} if there is an
  element $\unit\in A$ so that:
  \begin{enumerate}
  \item $\mu_1^0(\unit)=0$.
  \item For all $a\in A$, $\mu_2^0(a,\unit)=m_2^0(\unit,a)=a$.
  \end{enumerate}
\end{definition}

\begin{definition}\label{def:w-alg-strict-unital}
  A $w$-algebra is called {\em strictly unital}
  if there is an element $\unit\in A$ with the following properties:
  \begin{enumerate}
  \item For any $a\in A$, 
    $\mu^0_2(\unit,a)=\mu^0_2(a,\unit)=a$.
  \item For any sequence of elements 
    $a_1,\dots,a_n\in A$ with the property
    that $a_i=\unit$ for some $i$, and any weight $w$, 
    $\mu^w_n(a_1,\dots,a_n)=0$, except in the special case
    where $w=0$ and $n=2$ described above.
  \end{enumerate}
  A strictly unital $w$-algebra $A$ is \emph{split unital} if there is a
  $\Ring$-module splitting $A=\Ring\langle\unit\rangle\oplus A'$.
\end{definition}

When taking tensor products in Section~\ref{sec:wDiagApps}, we will
require that the weighted $\Ainf$-algebras and modules be strictly
unital. (There was no analogous requirement in the unweighted case.)
The results of these tensor products, however, are in general only
weakly unital. Fortunately, Theorem~\ref{thm:UnitalIsUnital} has the
following adaptation to the weighted case:

\begin{theorem}
  \label{thm:UnitalIsUnitalW}
  Every weakly unital $w$-algebra $\wAlg$ is isomorphic to a
  strictly unital one $\wAlg'$. Further, if $\wAlg$ is bonsai then
  $\wAlg'$ and the isomorphism can be chosen to be bonsai as well.
\end{theorem}

The proof rests on the following generalization of Lemma~\ref{lem:ModifyAction}.
To set it up, let
\[ 
  \phi_n^w\co \overbrace{A\kotimes{\Ground}\dots\kotimes{\Ground} A}^n\to A 
\]
be a map of $\Ground$-bimodules. Define
\[ 
  d \phi^w_n = \mu_1^0\circ \phi^w_n + \sum_{i=1}^n \phi^w_n\circ (\Id_{A^{\kotimes{\Ground} i-1}} \otimes \mu_1^0 \otimes \Id_{A^{\kotimes{\Ground} {n-i+1}}}).
\]
(This is the differential of $\phi_n^w$ thought of as an element of $\Mor(A^{\kotimes{\Ground} n},A)$.)

\begin{lemma}
  \label{lem:ModifyActionw}
  Suppose that $\wAlg$ is a $w$-algebra. Fix integers $n,v\geq 0$ with $n+2v\geq 2$, and let
  $\phi_n^v\co A^{\kotimes{\Ground} n}\to A$ be any map
  of degree $(2-\kappa)v+n-1$.
  There is a $w$-algebra $\wAlg'=(A,\overline{\mu}_n^w)$ isomorphic to $\wAlg$, so that 
  ${\overline\mu}^w_i=\mu^w_i$ for $i< n$ or $w<v$ and
  \[ {\overline \mu}^v_n = \mu^v_n + d \phi^v_n.\]
  Moreover if $d \phi^v_n=0$ or $n>2$ or $v>0$ then 
  \[ {\overline\mu}^v_{n+1}=\mu^v_{n+1}+\mu_2^{0}\circ (\Id_{A}\otimes \phi^v_n) + \mu_2^{0}\circ (\phi^v_n\otimes \Id_{A}) + 
    \sum_{i=1}^{n-1} \phi_n^v\circ(\Id_{A^{\kotimes{\Ground} i-1}}\otimes \mu_2^{0}\otimes \Id_{A^{\kotimes{\Ground} n+1-i}}).\]
\end{lemma}

\begin{proof}
  The proof of Lemma~\ref{lem:ModifyAction} applies (with
  Lemma~\ref{lem:wAlg-iso-is} in place of~\cite[Lemma 2.1.14]{LOT2}).
\end{proof}

The following is an easy modification of the proof of Theorem~\ref{thm:UnitalIsUnital}:

\begin{proof}[Proof of Theorem~\ref{thm:UnitalIsUnitalW}]
  The proof is a triple induction, on the weight $v$, the number of inputs $n$,
  and the first input $k$ at which the operation $\mu_n^v$ is not strictly
  unital. For the case $v=0$, the same proof as for
  Theorem~\ref{thm:UnitalIsUnital}, with Lemma~\ref{lem:ModifyActionw} in place
  of Lemma~\ref{lem:ModifyAction}, gives a weighted $\Ainf$-algebra whose
  undeformed honest $\Ainf$-algebra is strictly unital.

  Next, if $\mu_n^w$ is strictly unital for all $w<v$ (i.e.,
  $\mu_n^w(a_1,\dots,\unit,\dots,a_n)=0$ if $w<v$ and $(n,w)\neq(2,0)$) then it
  follows from the weight $v$ $\Ainf$ relation with input $(\unit,\unit)$ that
  $\mu_1^v(\unit)=\mu_1^0(\mu_2^v(\unit,\unit))$. So, applying
  Lemma~\ref{lem:ModifyActionw} with $\phi_1^v(a)=\mu_2^v(\unit,a)$
  gives an isomorphic weighted algebra with $\overline{\mu}_1^v(\unit)=0$.

  Now, suppose that there is a constant $v\geq 0$ so that $\wAlg$
  is strictly unital for all operations with weight less than $v$.
  Suppose also that there are integers $n\geq1$ and $k\geq 1$  so that
  $\wAlg$ is strictly unital for all weight $v$
  operations fewer than $n+1$ inputs, and moreover
  \begin{equation}
    \label{eq:PartialUnitalityw}
    \mu^v_{n+1}(a_1,\dots,a_{j-1},\unit,a_{j},\dots,a_{n})=0 
  \end{equation}
  holds for all $j<k$. 

  Suppose that $k+1<n+1$.
  We construct an isomorphic model for $\wAlg$
  by modifying only those actions with weight at least $v$ and at least $n+1$
  inputs, so that Equation~\eqref{eq:PartialUnitalityw} holds for
  all $j<k+1$. 
  When $k<n$, let
  $\phi^v_n(a_1,\dots,a_n)=\mu^v_{n+1}(a_1,\dots,a_i,\unit,a_{i+1},\dots,a_n)$.
  Apply Lemma~\ref{lem:ModifyActionw} to construct a new action with 
  \[
    {\overline\mu}^v_{n+1}=\mu^v_{n+1}+\mu_2^0
    \circ (\Id_{A}\otimes \phi^v_n) + \mu_2^0\circ (\phi^v_n\otimes \Id_{A}) + 
    \sum_{i=1}^{n-1} \phi^v_n\circ (\Id_{A^{\kotimes{\Ground} i-1}}\otimes \mu_2^0\otimes \Id_{A^{\kotimes{\Ground} i+1,n}}).
  \]
  The weighted $\Ainf$ relation on $\wAlg$ with input $(a_1,\dots,a_k,\unit,\unit,a_{k+1},\dots,a_n)$ proves that
  \[
    \overline{\mu}^v_{n+1}(a_1,\dots,a_k,\unit,a_{k+1},\dots,a_n)+{d
      \psi}^v_{n+1}(a_1,\dots,a_k,\unit,a_{k+1},\dots,a_n)=0,
  \]
  where
  \[ \psi^v_{n+1}(a_1,\dots,a_{n+1})=\mu^v_{n+2}(a_1,\dots,a_{k},\unit,a_{k+1},\dots,a_{n+1}).\]
  This uses the hypothesis that lower weight actions with 
  an arbitrary number of inputs are strictly unital.
  Applying Lemma~\ref{lem:ModifyActionw} once again, this time using $\psi^v_{n+1}$, we obtain the improved algebra 
  for which
  Equation~\eqref{eq:PartialUnitalityw} holds with $k+1$ in place of
  $k$. 

  By induction, we can now assume that Equation~\eqref{eq:PartialUnitalityw} holds for all $k<n+1$.
  We arrange for Equation~\eqref{eq:PartialUnitalityw} to hold for $k=n+1$ by applying Lemma~\ref{lem:ModifyAction} 
  once again, this time using $\psi^v_{n+1}$ with $k=n$. 

  Successively increasing $n$ and $v$, we obtain the claimed
  result. The fact that $\wAlg'$ and the isomorphism are bonsai
  if $\wAlg$ is bonsai follows by the same argument as in the
  unweighted case.
\end{proof}

\begin{definition}
  A homomorphism $f$ between weakly unital $w$-algebras $\wAlg$ and
  $\wBlg$ is {\em weakly unital} if
  $f_1^0(\unit_{\wAlg})=\unit_{\wBlg}$. A homomorphism between
  strictly unital $w$-algebras $\wAlg$ and $\wBlg$ is {\em strictly
    unital} if it is weakly unital and $f_n^w(a_1,\dots,a_n)=0$ if
  some $a_i=\unit_\wAlg$ and $(n,w)\neq(1,0)$.
\end{definition}

\begin{proposition}\label{prop:w-qi-unital}
  If $\wAlg$ and $\wBlg$ are strictly unital $w$-algebras and
  $f\co \wAlg\to\wBlg$ is a weakly unital quasi-isomorphism then there
  is a strictly unital quasi-isomorphism $g\co \wAlg\to\wBlg$. If
  $\wAlg$, $\wBlg$, and $f$ are bonsai then $g$ can be chosen to be
  bonsai, as well.
\end{proposition}

The proof rests on the following weighted analogue of Lemma~\ref{lem:ModifyHomomorphism}
(in the spirit of Lemma~\ref{lem:ModifyActionMw}):

\begin{lemma}
  \label{lem:ModifyHomomorphismw}
  Let $\wAlg$ and $\wBlg$ be weighted $\Ainf$-algebras, 
  and $f\co \wAlg\to \wBlg$ be a weighted $\Ainf$-algebra homomorphism.
  Let $h^v\co A^{\kotimes{\Ground} n}\to B$ be a map of degree
  $(2-\kappa)v+n-1$, for some $(n,v)\neq(0,0)$. 
  Then, there is a new weighted $\Ainf$-algebra homomorphism 
  $F\co \wAlg\to \wBlg$
  with $F^w_\ell=f^w_\ell$ for all $\ell<n$ or $w<v$, and
  \begin{equation}
    \label{eq:dhn0w}
    F^v_n = f^v_n+d h^v,
  \end{equation}
  where
  \[
    d h^v(a_1,\dots,a_n)=\mu_1^{0} \circ h^v(a_1,\dots,a_n) + 
    \sum_{i=1}^n h(a_1,\dots,a_{i-1},\mu_1^{0}(a_i),a_{i+1},\dots a_n).
  \]
  (Here, the two terms labelled
  $\mu_1^0$ denote the weight zero differential on $\wBlg$ and
  $\wAlg$ respectively.)
  If moreover $d h^v=0$, so $F^v_n=f^v_n$, then
  \begin{align}
    F^v_{n+1}&(a_1,\dots,a_{n+1})=f^v_{n+1}(a_1,\dots,a_{n+1}) \label{eq:dh1w}\\
    &+\sum_{i=1}^n h^v(a_1,\dots,a_{i-1},\mu_2^{0}(a_i,a_{i+1}),a_{i+2},\dots,a_{n+1}) \nonumber \\
    &+ \mu_2^{0}(f_1^{0}(a_1),h^v(a_2,\dots,a_{n+1}))+
    \mu_2^{0}(h^v(a_1,\dots,a_{n}),f_1^{0}(a_{n+1})). \nonumber
  \end{align}
\end{lemma}

\begin{proof}
  We adapt the proof of Lemma~\ref{lem:ModifyHomomorphism}, with the
  following changes.  In the present case, we are attempting to define
  components of a weighted homomorphism
  $F^w_\ell\co A^{\kotimes{\Ground} \ell}\to B$.  Again, we allow trees of the three
  Types~\ref{type:f}-\ref{type:dh} as in the proof of that lemma, with
  the understanding that each vertex is now also decorated by some
  weight, and the sum of the weights in each tree must be $w$.
  The same inductive argument then shows that the maps $F^w_\ell$ described
  above give a weighted homomorphism.
\end{proof}

\begin{proof}[Proof of Proposition~\ref{prop:w-qi-unital}]
  This is an inductive argument, like the proof of
  Proposition~\ref{prop:qi-unital}. 
  Suppose there are $w$, $n$, and $k$ so that the map $f\co \wAlg\to\wBlg$
  satisfies
  \begin{equation} 
    \label{eq:PrettyUnitalWeightedMap}
    f_{m+1}^v(a_1,\dots,a_{j-1},\One,a_j,\dots,a_m)=0 
  \end{equation}
  in the following cases:
  \begin{itemize}
    \item when $v<w$
    \item when $v=w$ and  $1\leq m\leq n$
    \item when $v=w$, $w=n+1$, and  $j\leq k$.
  \end{itemize}We construct the modification of $f$ by
  an outer induction on $w$, and then induction on $k$ and then $n$.
  
  The base case $(w,n)=(0,1)$ is the definition of weak unitality.
  For fixed $w$, the inductive steps on $k$ and $n$ follow as in the
  proof of Proposition~\ref{prop:qi-unital}, using
  Lemma~\ref{lem:ModifyHomomorphismw} in place of
  Lemma~\ref{lem:ModifyHomomorphism}. The base case for this induction
  is $f_1^w(\One)=0$ if $w>0$, and is $f_1^0(\One)=\One$.

  For the inductive step on the weight, suppose that
  Equation~\eqref{eq:PrettyUnitalWeightedMap} holds for all $v<w$.
  Then, the weighted $\Ainf$ relation with inputs $(\One,\One)$
  implies that $f_{1}^{w}(\One)=\mu_1^0(f_2^w(\One,\One))$. So,
  applying Lemma~\ref{lem:ModifyHomomorphismw} with
  $h^w_0=f_2^w(\One,\One)$ gives a new homomorphism $F$ with
  $F_1^w(\One)=0$, as desired. This completes the
  induction.

  The bonsai statement follows as in the unweighted case.
\end{proof}

\subsection{Weighted modules and weighted module trees}\label{sec:w-mods}

\begin{figure}
  \centering
  %Font is 18 point.
  \includegraphics[scale=.66667]{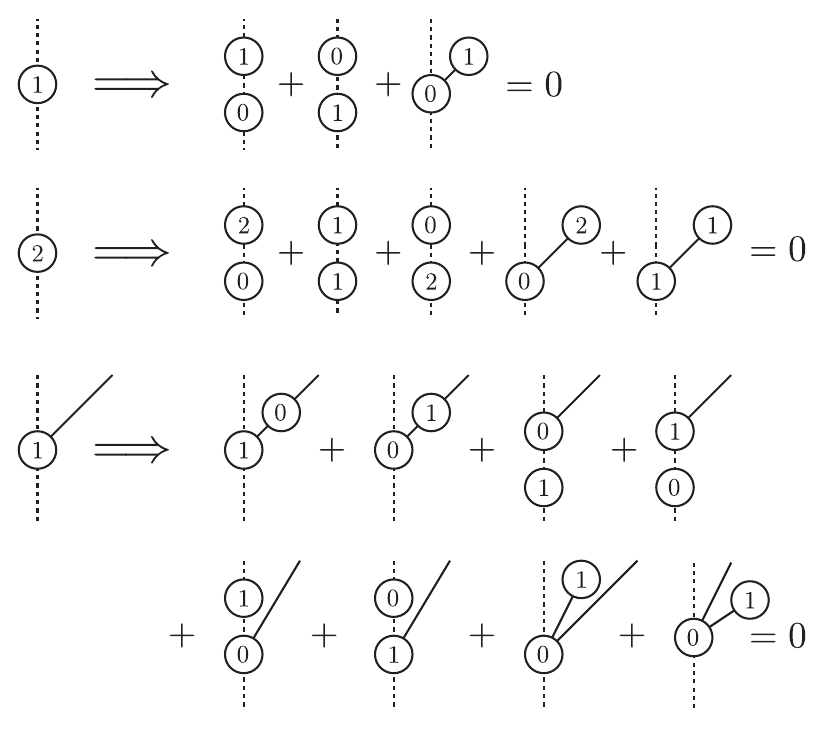}
  \caption[Some weighted A-infinity module relations]{\textbf{Weighted $\Ainf$-module relations.}}
  \label{fig:wModRel}
\end{figure}

\begin{definition}\label{def:wmod}
  Fix a weighted $\Ainf$-algebra $\wAlg$ over $\Ground$ with weight grading
  $\kappa$. A (right) \emph{weighted $\Ainf$-module} or
  \emph{$w$-module} $\wMod$ over $\wAlg$ consists of:
  \begin{itemize}
  \item a right $\Ground$-module $M$ and
  \item for each pair of integers $n,w\geq 0$ a linear map
    \[
      m^w_{1+n}\co
      M\kotimes{\Ground}\overbrace{A\kotimes{\Ground}\dots\kotimes{\Ground} A}^n\to M\grs{1-n-\kappa w}
    \]
    of right $\Ground$-modules 
  \end{itemize}
  satisfying the structure equation
    \begin{equation}
    \label{eq:wMod-StructuralEquation}
    \begin{split}
      \sum_{\substack{1\leq i\leq j\leq n\\
          u+v=w}}&m^u_{2+n-j+i}(x,a_1,\dots,a_{i-1},\mu_{j-i}^v(a_i,\dots,a_{j-1}),a_{j},\dots,a_n)
      \\
      &+\sum_{\substack{0\leq i\leq n\\
          u+v=w}}m^u_{1+n-i}(m_{1+i}^v(x,a_1,\dots,a_{i}),a_{i+1},\dots,a_n)
      = 0
    \end{split}
  \end{equation}
  for each pair of non-negative integers $(n,w)$.

  Again, $w$ is called the \emph{weight} of $m^w_n$.
\end{definition}
The first few terms in the structure
equation~\eqref{eq:wMod-StructuralEquation} are shown in
Figure~\ref{fig:wModRel}.

A weighted $\Ainf$-module $\wMod$ over $\wAlg$ gives rise to an $\Ainf$-module
$M\llbracket t\rrbracket=M\rotimes{\Ring}\Ring\llbracket t\rrbracket$ over $A\llbracket t\rrbracket$ with operations
\[
  m_{1+n}(x,a_1,\dots,a_n)=\sum_{w=0}^\infty m_{1+n}^w(x,a_1,\dots,a_n)t^w.
\]
Equation~\eqref{eq:wMod-StructuralEquation} is equivalent to the
$\Ainf$-module structure equation for $M\llbracket t\rrbracket$.

Next we turn to the analogue of stably-weighted trees. A \emph{stably-weighted module tree} is a stably-weighted tree for which the
left-most leaf is an input. Let $\wMTrees{1+n}{w}$ denote the set of
stably-weighted module trees with total weight $w$ and $1+n$
inputs.
Call a stably-weighted tree which is not a
stably-weighted module tree (i.e., for which the left-most leaf is
a popsicle) a \emph{left-popsicle tree}. Note that if $T$ is a
left-popsicle tree then so is every tree $S$ in $\partial T$, so the
left-popsicle trees form a subcomplex $L^{n,w}_*$ of the weighted trees complex
$\wTreesCx{n}{w}$. The \emph{weighted module trees complex} with $1+n$
inputs and weight $w$, $\wMTreesCx{1+n}{w}$ is the quotient complex of
$\wTreesCx{1+n}{w}$ by the subcomplex of left-popsicle trees. 

The weighted module trees complex inherits a grading from the weighted
trees complex.

As in the case of stably-weighted trees, there are stacking
operations for stably-weighted module trees: given
$S\in\wMTrees{1+m}{v}$ and $T\in \wMTrees{1+n}{w}$ we can glue the
output of $S$ to the first input of $T$ to obtain
$T\circ_1 S\in\wMTrees{1+m+n}{v+w}$, and given $S\in\wTrees{m}{v}$ and
$T\in\wMTrees{1+n}{w}$ and $1<i\leq n+1$ we can glue the output of $S$
to the $i\th$ input of $T$ to obtain
$T\circ_i S\in\wMTrees{m+n}{v+w}$.
The operation $\circ_i$ induces a chain map
\begin{align*}
\phi_{i,j,n;v,w}&\co 
                  \begin{cases}
                    \wTreesCx{j-i+1}{v}\rotimes{\Ring}\wMTreesCx{n+i-j}{w} \to
                    \wMTreesCx{n}{v+w} & i>1\\
                    \wMTreesCx{j}{v}\rotimes{\Ring}\wMTreesCx{n+1-j}{w} \to
                    \wMTreesCx{n}{v+w}
                    & i=1
                  \end{cases}\\
\phi_{i,j,n;v,w}(S,T)&=T\circ_i S
\end{align*}
for any $1\leq i\leq n$ and $i-1\leq j\leq n$.

The connection with weighted modules is as follows. Given a weighted
$\Ainf$-algebra $\wAlg$, a chain complex $(M,m_1^0)$, and
maps $m_{1+n}^w\co M\kotimes{\Ground} A^{\kotimes{\Ground} n}\to M\grs{1-n-\kappa w}$, $n\geq0$,
$w\geq 0$, $(n,w)\neq(0,0)$, for any stably-weighted module tree $T$
with $1+n$ inputs and total weight $W$ there is an induced map
$m(T)\co M\kotimes{\Ground} A^{\kotimes{\Ground} n}\to M\grs{-\dim(T)+(2-\kappa)W}$ as follows: replace each
internal vertex $v$ on the left-most strand of $T$ by $m_{1+n}^w$
(where $w$ is the weight of $v$ and $n+2$ is the valence of $v$) and
each internal vertex $v$ not on the left-most strand by $\mu_n^w$
(where $w$ is the weight of $v$ and $n+1$ is the valence of $v$) and
compose these maps along the edges of $T$. Extending linearly we
obtain a map $\wMTreesCx{1+n}{w}\to \Mor(M\kotimes{\Ground} A^{\kotimes{\Ground} n},M\grs{(2-\kappa)w})$.
\begin{lemma}\label{lem:wMod-reinterp}
  The collection $\{m_{1+n}^w\}$ satisfy the weighted $\Ainf$-module
  relations if and only if for each $(n,w)$ the map
  $m\co\wMTreesCx{1+n}{w}\to \Mor(M\kotimes{\Ground} A^{\kotimes{\Ground} n},M\grs{(2-\kappa)w})$ is
  a chain map.
\end{lemma}
\begin{proof}
  The proof is similar to the proof of Lemma~\ref{lem:wAlg-reinterp},
  and is left to the reader.
\end{proof}

\begin{definition}
  A $w$-module $\wMod$ is \emph{bonsai} if there is an integer~$N$ so
  that for any stably-weighted module tree~$T$ with
  $\dim(T) > N$, we have $m(T) = 0$.
\end{definition}

\subsection{Weighted module transformation trees and weighted module maps}
\label{sec:wmod-morph}
\begin{definition}
  A \emph{weighted module transformation tree} is a weighted tree $T$
  whose leftmost leaf is an input, together with a distinguished
  internal vertex $v$ on the leftmost strand of $T$, subject to the
  condition that every internal vertex of $T$ except $v$ must be
  stable.  (Compare Section~\ref{sec:multiplihedron}; we will
  sometimes call the distinguished vertex \emph{purple}.) Let
  $\wMTransCx{n}{w}$ be the free $\Ring$-module with basis the weighted module
  transformation trees with $n$ inputs and total weight $w$. There is
  a grading on $\wMTransCx{n}{w}$ defined by
  \[
    \dim(T)=n+2w-v
  \]
  where $v$ is the number of internal vertices.
  The differential of a weighted module transformation $T$ is the sum
  of all ways of inserting an edge into $T$ to obtain a new weighted
  module transformation tree.
\end{definition}

There are obvious composition maps
\begin{align}
  \circ_1&\co \wMTransCx{n}{w}\rotimes{\Ring} \wMTreesCx{m}{v}\to \wMTransCx{m+n-1}{v+w}\label{eq:w-mod-trans-compo-1}\\
  \circ_1&\co \wMTreesCx{n}{w}\rotimes{\Ring} \wMTransCx{m}{v}\to \wMTransCx{m+n-1}{v+w}\label{eq:w-mod-trans-compo-2}\\
  \qquad\qquad\circ_i&\co \wMTransCx{n}{w}\rotimes{\Ring} \wTreesCx{m}{v}\to
           \wMTransCx{m+n-1}{v+w} \label{eq:w-mod-trans-compo-3}\qquad\qquad 1<i\leq n.
\end{align}

\begin{lemma}\label{lem:wMTrans-is-cx}
  The $n$-input, weight $w$ module transformation trees complex
  $\wMTransCx{n}{w}$ is a chain complex, and the differential
  decreases the grading by one. Further, the composition
  maps~\eqref{eq:w-mod-trans-compo-1},~\eqref{eq:w-mod-trans-compo-2},
  and~\eqref{eq:w-mod-trans-compo-3} are chain maps.
\end{lemma}
\begin{proof}
  The proof that the differential satisfies $\bdy^2=0$ is similar to
  Lemma~\ref{lem:wTrees-is-cx}.  The facts that the differential
  decreases the grading by $1$ and the composition maps are chain maps
  are immediate from the definitions.
\end{proof}

Given weighted $\Ainf$-modules $\wMod$ and $\wNod$ over $\wAlg$ 
and integers $w\geq 0$ and $d\in\ZZ$, a
{\em morphism of degree $d$ and weight $w$} is a collection of $\Ground$-module maps 
$f_{1+n}^w\co M\kotimes{\Ground} A^{\kotimes{\Ground} n}\to N\grs{-n-d-\kappa w}$.
Thus, the space of morphisms of weight $w$ is the direct product
$\Mor^w(\wMod,\wNod)=\prod_{n=0}^{\infty} \Mor(M\kotimes{\Ground} A^{\kotimes{\Ground} n},N)$.
The {\em morphism space} from $\wMod$ to $\wNod$ is the direct product
\[
  \wMor(\wMod,\wNod)=\prod_{w=0}^{\infty}\Mor^w(\wMod,\wNod).
\]
This is a chain complex 
with 
\[
  d \co \Mor(\wMod,\wNod)\to \Mor(\wMod,\wNod)
\]
given by
\begin{equation}\label{eq:d-of-w-mod-map}
\begin{split}
  (d f)_{1+n}^{v}(x,a_1,\dots,a_n)=
  &\sum_{w=0}^v\sum_{k=0}^{n}\bigl(f_{1+n-k}^{v-w}(m_{\wMod,1+k}^{w}(x,a_1,\dots,a_k),a_{k+1},\dots,a_n)\\
  &+m_{\wNod,1+n-k}^{w}(f_{1+k}^{v-w}(x,a_1,\dots,a_k),a_{k+1},\dots,a_n)\bigr)\\
  &+\sum_{w=0}^v\sum_{2\leq i\leq j\leq n}f_{1+n+i-j}^{v-w}(x,a_1,\dots,\mu_{j-i}^{w}(a_i,\dots,a_{j-1}),\dots,a_n).
\end{split}
\end{equation}
(We adhere to the following notational convention: if $g\in\wMor(\wMod,\wNod)$, then $g^w$ denotes weight $w$ component, in
$\Mor^w(\wMod,\wNod)$.)
The kernel of $d$ is the \emph{$w$-module
  homomorphisms}, and if two morphisms differ by a boundary then they
are \emph{chain homotopic}.

\begin{example}
  The \emph{identity morphism} $\Id_{\wMod}$ of $\wMod$ is defined by
  \[
    \Id_{1+n}^w=
    \begin{cases}
      \Id_M & n=0~\text{and}~w=0\\
      0 & \text{otherwise.}
    \end{cases}
  \]
  It is straightforward to check that $\Id_{\wMod}$ is a homomorphism.
\end{example}

Given a degree $0$ morphism $f\co \wMod \to \wNod$ and a weighted module
transformation tree $T$ with $n+1$ inputs and total weight $w$ there is an induced map
\[
  F(T)\co M\kotimes{\Ground} A^{\kotimes{\Ground} n}\to N\grs{-\dim(T)+(2-\kappa)w}
\]
This gives a map
\[
  F\co \wMTransCx{n+1}{w}\to \Mor(M\kotimes{\Ground} A^{\kotimes{\Ground} n},N\grs{(2-\kappa)w}).
\]
\begin{lemma}\label{lem:w-F-f}
  A morphism $\{f_{n}^w\}$ is a homomorphism if and only if for each
  $(n,w)$ the induced map $F$ is a chain map. More generally, the map
  \[
    \Mor(\wMod,\wNod)\to \prod_{n,w}\Mor(M\kotimes{\Ground} A^{\kotimes{\Ground} n},N\grs{-\dim(T)+(2-\kappa)w})
  \]
  given by $f\mapsto F$ is a chain map.
\end{lemma}
(Compare Lemma~\ref{lem:F-f}.)
\begin{proof}
  Again, the proof is straightforward and is left to the reader.
\end{proof}

\begin{lemma}\label{lem:w-mod-iso}
  Let $\wMod$ and $\wNod$ be $w$-modules over $\wAlg$ and let $f\co
  \wMod\to\wNod$ be a homomorphism. Then $f$ is an isomorphism if and only if $f_1^0$ is an isomorphism.
\end{lemma}
\begin{proof}
  As in the proof of Lemma~\ref{lem:w-mod-iso}, it suffices to
  construct left and right inverses to $f$, and we will focus on the
  left inverse. We construct the left inverse $g_n^w$ inductively in
  $w$ and then $n$. Suppose $g_m^v$ has been defined for all $v<w$ and
  for $v=w$ and $m\leq n$, satisfying the condition that $(g\circ
  f)_m^v=0$ for $v<w$ or $v=w$ and $m\leq n$ (except for $(g\circ
  f)_1^0=\Id$). We have
  \[
    (g\circ f)_{n+1}^w(x,a_1,\dots,a_n)=\sum_{\substack{m\leq n\\v\leq w}}g_{n-m+1}^{w-v}(f_{m+1}^v(x,a_1,\dots,a_m),a_{m+1},\dots,a_n).
  \]
  All of these terms have been defined except
  \[
    g_{n+1}^w(f_1^0(x),a_1,\dots,a_n).
  \]
  So, define
  \[
    g_{n+1}^w(f_1^0(x),a_1,\dots,a_n)=\sum_{\substack{m\leq n\\v\leq w\\(m,v)\neq(0,0)}}g_{n-m+1}^{w-v}(f_{m+1}^v(x,a_1,\dots,a_m),a_{m+1},\dots,a_n),
  \]
  using the fact that $f_1^0$ is an isomorphism, and continue the induction.
\end{proof}

\begin{definition}\label{def:w-mod-map-bonsai}
  Let $\wMod$ and $\wNod$ be bonsai $w$-modules over a bonsai
  $w$-algebra $\wAlg$.  A morphism $f\in \Mor(\wMod,\wNod)$ is
  \emph{bonsai} if there is an integer $N$ so that for any module
  transformation tree $T$ with $\dim(T)>N$, $f(T)=0$.
\end{definition}
As in the unweighted case, it is clear that the bonsai morphisms form
a subcomplex of $\Mor(\wMod,\wNod)$.

\subsection{Units in weighted \texorpdfstring{$\Ainf$-}{A-infinity-}modules}\label{sec:w-mod-units}

The discussion from Section~\ref{sec:w-units} adapts to weighted modules as follows.

\begin{definition}\label{def:w-mod-strict-unital}
  A $w$-module $\wMod$ over a weakly unital $w$-algebra $\wAlg$
  is {\em weakly unital} if $m_2^0(x,\unit)=x$ for all $x\in M$.
  A weakly unital $w$-module $\wMod$ over a strictly unital $w$-algebra
  $\wAlg$ is {\em strictly unital} if for any sequence of 
  elements $a_1,\dots,a_n\in A$ with the property that $a_i=\unit$ for
  some $i$, any $x\in M$, and any weight $w$,
  $m^w_{1+n}(x,a_1,\dots,a_n)=0$, except in the special case where
  $w=0$ and $n=1$ (in which case $m^0_2(x,\unit)=x$).
\end{definition}

The proofs of the following lemma and theorem are essentially the same as the
proofs of Lemma~\ref{lem:ModifyActionM} and Theorem~\ref{thm:UnitalIsUnitalM}
respectively, with an extra induction on the weight.

\begin{lemma}
  \label{lem:ModifyActionMw}
  Suppose that $\wMod$ is a weakly unital $\Ainf$-module and
  \[ \phi^v_n \co M \kotimes{\Ground} A \kotimes{\Ground}\dots\kotimes{\Ground} A \to M\]
  is any $\Ground$-module map of degree $(2-\kappa)v+n-1$.
  There is an $\Ainf$-module ${\overline \wMod}$ isomorphic to $\wMod$ so that
  ${\overline m}^w_i=m^w_i$ for $i<n$ or $w<v$; and 
  \[ {\overline m}^v_n= m^v_n + d \phi^v_n\]
  (where $d$ is the differential on the morphism complex $\Mor(M\kotimes{\Ground} A^{\kotimes{\Ground} (n-1)},A)$).
\end{lemma}

\begin{theorem}
  \label{thm:UnitalIsUnitalMw}
  Every weakly unital $w$-module $\wMod$ over a strictly unital
  $w$-algebra is isomorphic to a strictly unital $w$-module $\wNod$. If
  $\wMod$ is bonsai then $\wNod$ and the isomorphism can be taken to
  be bonsai, as well.
\end{theorem}

\begin{convention}\label{conv:su-mods}
  Unless otherwise noted, when we are over strictly unital $w$-algebras we will
  only consider strictly unital $w$-modules.
\end{convention}

If $\wMod$ and $\wNod$ are strictly unital then a morphism
$f\co\wMod\to\wNod$ is \emph{strictly unital} if
$f_{n+1}^w(m,a_1,\dots,a_n)=0$ if some $a_i=\unit$. The strictly unital
morphisms form a subcomplex $\uwMor(\wMod,\wNod)$ of the complex
$\Mor(\wMod,\wNod)$ of all $w$-module morphisms.

\begin{proposition}\label{prop:w-mod-mor-cx-unital}
  Given strictly unital $w$-modules $\wMod$ and $\wNod$ over a
  strictly unital $w$-algebra $\wAlg$, the inclusion of chain
  complexes $\uwMor(\wMod,\wNod)\to\Mor(\wMod,\wNod)$ is a chain
  homotopy equivalence. If $\wAlg$, $\wMod$, and $\wNod$ are bonsai,
  the same statement holds for the subcomplexes of bonsai strictly
  unital and non-unital morphisms.
\end{proposition}
\begin{proof}
  As a notational preliminary, recall that in the proof of Proposition~\ref{prop:mod-mor-cx-unital}, 
  there was an inclusion map $\suMor(\cModule,\cNodule)\to \Mor(\cModule,\cNodule)$.
  Let 
  \begin{equation}
    \label{eq:RecallH}
    \begin{split}
    H = &\sum_{i_1<\dots<i_k} H_{i_{k}}\circ \partial \circ H_{i_{k-1}}\circ \dots\circ\partial \circ H_{i_1} \\
    &+ \sum_{i_1<\dots<i_k, 1\leq j<k}
     H_{i_k}\circ \partial\circ H_{i_{k-1}}\circ \partial \circ \dots\circ \partial\circ  
    H_{i_{j+1}}\circ H_{i_j}\circ \partial\circ\dots\circ H_{i_1}\circ \partial.
    \end{split}
  \end{equation}
  where $H_i$ are as in Equation~\eqref{eq:DefHi}.  It is
  straightforward to see that the homotopy inverse $\Pi$ (as in
  Equation~\ref{eq:DefOfPi}) to the inclusion map
  $\suMor(\cModule,\cNodule)\to \Mor(\cModule,\cNodule)$ is given by
  $\Id + \partial \circ H + H\circ \partial$.  We now upgrade these maps to the weighted context.

  Given $f\in \wMor(\wMod,\wNod)$, 
  let ${\mathcal H}_i(f)\in \wMor(\wMod,\wNod)$ be the map with 
  \[ {\mathcal H}_i(f)^w(x,a_1,\dots,a_{n-2})=
  \begin{cases}
       0 &{\text{if $n-1<i$}} \\
       f^w_n(x,a_1,\dots,a_{i-1},\One,a_{i+1},\dots,a_{n-2}) & \text{otherwise.}
     \end{cases}
   \]
   Define $\mathcal{H}\co \wMor(\wMod,\wNod)\to \wMor(\wMod,\wNod)$ by 
   \[
    \begin{split}
    \mathcal{H} = 
    &\sum_{i_1<\dots<i_k} \mathcal{H}_{i_{k}}\circ \partial \circ \mathcal{H}_{i_{k-1}}\circ \dots\circ\partial \circ \mathcal{H}_{i_1} \\
    &+ \sum_{i_1<\dots<i_k, 1\leq j<k}
     \mathcal{H}_{i_k}\circ \partial\circ \mathcal{H}_{i_{k-1}}\circ \partial \circ \dots\circ \partial\circ  
    \mathcal{H}_{i_{j+1}}\circ \mathcal{H}_{i_j}\circ \partial\circ\dots\circ \mathcal{H}_{i_1}\circ \partial.
    \end{split}
  \]
  
  For a fixed integer $v\geq 0$, let
  $\uwMor_{v}(\wMod,\wNod)\subset \Mor(\wMod,\wNod)$ consist of those
  morphisms $f$ with
  \[
    f^w(x,a_1,\dots,a_{i-1},\One,a_{i+1},\dots,a_n)=0
  \]
  for all $w<v$ and $1\leq i<n+1$. This specifies a filtration of the chain complex
  $\Mor(\wMod,\wNod)$, with the following properties:
  \begin{itemize}
    \item $\uwMor_{v_1}(\wMod,\wNod)\supset \uwMor_{v_2}(\wMod,\wNod)$
      if $v_1<v_2$,
    \item
      $\uwMor_0(\wMod,\wNod)=\Mor(\wMod,\wNod)$, and
    \item
      $\uwMor(\wMod,\wNod)=\bigcap_{v=0}^{\infty}\uwMor_v(\wMod,\wNod)$.
  \end{itemize}

  Let
  \begin{equation}
    \label{eq:DefOfF}
    F=\Id + \partial\circ{\mathcal H}+{\mathcal H}\circ \partial.
    \end{equation}
  We claim that  
  \begin{equation}
    \label{eq:Funitizes}
    F(\uwMor_{v}(\wMod,\wNod))\subseteq\uwMor_{v+1}(\wMod,\wNod).
  \end{equation}
  To see this, observe that for $f\in\uwMor_{v}(\wMod,\wNod)$ , the terms in the weight $v+1$ part of
  $(\Id+\partial\circ{\mathcal H}+{\mathcal H}\circ \partial)(f)$ involving components of $f^w$ with weight $w<v+1$ vanish after we specialize to $a_i=\One$.
  Thus, 
  \begin{align*} 
    [(\Id+\partial\circ{\mathcal H}+&{\mathcal H}\circ \partial)(f)]^{v+1}(a_1,\dots,a_{i-1},\One,a_{i+1},\dots,a_n) = \\
                                                 &[(\Id + \partial \circ H +  H \circ \partial)(f^{v+1})](a_1,\dots,a_{i-1},\One,a_{i+1},\dots,a_n),
  \end{align*}
  which vanishes by Proposition~\ref{prop:mod-mor-cx-unital},
  verifying Equation~\eqref{eq:Funitizes}.
  (We are actually applying the proof 
  of Proposition~\ref{prop:mod-mor-cx-unital},
  where it is shown that
  $\Pi=\partial \circ H + H\circ \partial$ is a
  projection to the strictly unital morphism space.)

  Let $P_{\leq v}\co \Mor(\wMod,\wNod)\to \bigoplus_{w\leq v}
  \Mor^w(\wMod,\wNod)$ denote the projection to the weight less than
  $v$ portion.  If $\xi\in\uwMor_{v}(\wMod,\wNod)$, then
  $P_{\leq v}\circ F(\xi)=P_{\leq v}(\xi)$.
  Let $F^{(\circ k)}(\xi)$ denote the $k^{th}$ iterate of $F$
  applied to $\xi$. It follows that for any pair of 
  integers $0\leq w \leq i$, 
  \begin{equation}
    \label{eq:PhiIsInfiniteIterate}
    P_{\leq w}\circ F^{(\circ i)}=P_{\leq w} \circ F^{(\circ w)}.
  \end{equation}
   Consider the function $\Phi \co
  \Mor(\wMod,\wNod)\to \suMor(\wMod,\wNod)$, characterized by $P_{\leq
    w} \circ \Phi=P_{\leq w} \circ F^{(\circ w)}$.  We think of this
  as the infinite composition of $F$ with itself, since $F\circ \Phi =
  \Phi$. Moreover, $\Phi$ is the identity map on
  $\uwMor(\wMod,\wNod)\subset\Mor(\wMod,\wNod)$.

  The map $F$ is chain homotopic to the identity. This 
  homotopy can be used to construct a chain homotopy of $\Phi$ to the identity.
  In detail, the chain homotopy is given by
  \[ {\mathcal K}=\sum_{i=0}^{\infty} {\mathcal H}\circ F^{(\circ i)}.\]
  To see that ${\mathcal K}$ is well-defined,
  we show that for any $w$,
  $P_{\leq w}\circ {\mathcal K}$ is a finite
  sum. Specifically,
  \begin{equation}P_{\leq w} (\sum_{i=0}^{\infty} {\mathcal H}\circ F^{(\circ i)})
    = \sum_{i=0}^{\infty} P_{\leq w}\circ {\mathcal H}\circ  F^{(\circ i)} 
    =\sum_{i=0}^{w} P_{\leq w}\circ {\mathcal H} \circ F^{(\circ i)},
    \label{eq:FiniteSum}
  \end{equation}
  since for all $i\geq w$, the image of $F^{(\circ i)}$ is in $\uwMor_{w+1}(\wMod,\wNod)$
  (by Equation~\eqref{eq:Funitizes}) and 
  \[ P_{\leq w}\circ {\mathcal H}|_{\uwMor_{w+1}(\wMod,\wNod)}=0 \]
  (which in turn is clear from the definition of ${\mathcal H}$).

  To see that $\Phi=\Id+\partial \circ {\mathcal K}+{\mathcal
    K}\circ \partial$, 
  it suffices to verify that 
  \begin{equation}
    \label{eq:PhiIsAHomotopyEquivalence}
    P_{\leq w} \circ
    (\Id+\partial \circ {\mathcal K}+{\mathcal K}\circ \partial)=P\circ  F^{\circ w}.
  \end{equation}
  To this end,  abbreviate, for any fixed $w$,
  $P=P_{\leq w}$.
  Using the fact that $P$ is a chain map,
  the definition of ${\mathcal K}$,
  Equation~\eqref{eq:FiniteSum},
  $P$ is a chain map (again),
  $F$ is a chain map,
  Equation~\eqref{eq:DefOfF}, and Equation~\eqref{eq:PhiIsInfiniteIterate}, in succession, we see that
  \begin{align*}
    P\circ (\Id + \partial\circ {\mathcal K} + {\mathcal K}\circ \partial)
    &= P + \partial \circ P \circ {\mathcal K}
    + P\circ {\mathcal K}\circ \partial 
    \\
    %1
    &=P +  \partial \circ P \circ {\mathcal H} \circ \left(\sum_{i=0}^\infty F^{(\circ i)}\right)
    + P \circ {\mathcal H} \circ \left(\sum_{i=0}^\infty F^{(\circ i)}\right)\circ \partial 
    \\
    %2
    &=P + \partial \circ P \circ {\mathcal H} \circ \left(\sum_{i=0}^w F^{(\circ i)}\right)
    + P \circ {\mathcal H} \circ \left(\sum_{i=0}^w F^{(\circ i)}\right)\circ \partial 
    \\
    %3
    &=P + P \circ \partial  \circ {\mathcal H} \circ \left(\sum_{i=0}^w F^{(\circ i)}\right)
    + P \circ {\mathcal H} \circ \left(\sum_{i=0}^w F^{(\circ i)}\right)\circ \partial 
    \\
    %4
    &=P + P \circ \left(\partial  \circ {\mathcal H}+{\mathcal H}\circ \partial\right) \circ\left( \sum_{i=0}^w F^{(\circ i)}\right)
    \\
    \\
    %5
    &=P + P \circ \left(\Id + F\right) \circ\left( \sum_{i=0}^w F^{(\circ i)}\right)
    \\
    %6
    &=P + P \circ (\Id + F^{\circ (w+1)}) = P\circ F^{\circ (w+1)}=P\circ F^{\circ w},  \end{align*}
 verifying Equation~\eqref{eq:PhiIsAHomotopyEquivalence}.

  Clearly, now, $\Phi=\Id+\partial \circ {\mathcal K}+{\mathcal K}\circ \partial$
  is the homotopy inverse to the inclusion
  $\uwMor(\wMod,\wNod)\subset \Mor(\wMod,\wNod)$. The fact
  that $\Phi$ preserves the property of being bonsai follows as in the
  unweighted case.
\end{proof}

\subsection{Weighted type \texorpdfstring{$D$}{D} structures}
\label{sec:wD}

Let $\wAlg$ by a $w$-algebra over $\Ground$ and fix an element
$X\in\Ground$ with grading $-\kappa$. Assume that $X$ acts centrally
on $A$, in the sense that $Xa=aX$ for all $a\in A$.
\begin{definition}\label{def:wD-bounded}
  Let $P$ be a graded left $\Ground$-module and
  $\delta^1\co P\to A\kotimes{\Ground} P\grs{1}$ a left
  $\Ground$-module map. Define
  $\delta^n\co P\to A^{\kotimes{\Ground} n}\kotimes{\Ground} P\grs{n}$
  to be the result of iterating
  $\delta^1$ $n$ times. (In the special case $n=0$, $\delta^0$ is the
  identity map $P\to P$.) We say that $\delta^1$ (or $(P,\delta^1)$)
  is \emph{operationally bounded} (or just \emph{bounded}) if
  \begin{itemize}
  \item $\delta^n$ vanishes identically for $n$ sufficiently large and
  \item the action of $X$ on $A$ is nilpotent, i.e., for $n$
    sufficiently large and all $a\in A$, $X^n\cdot a=0$.
  \end{itemize}
\end{definition}

\begin{definition}\label{def:wD}
  With notation as in Definition~\ref{def:wD-bounded}, suppose that
  either $\delta^1$ is bounded or $\wAlg$ is bonsai. We
  say that $(P,\delta^1)$ is a (left) \emph{type $D$ structure
  over $\wAlg$ with charge $X$}  if 
\begin{equation}
  \label{eq:w-D-str}
  \sum_{w=0}^\infty\sum_{n=0}^\infty X^w(\mu_n^w\otimes\Id)\circ \delta^n=0\in A\kotimes{\Ground} P.
\end{equation}
\end{definition}
Sometimes we abbreviate a type $D$ structure as
$\lsup{\wAlg}\wP=(P,\delta\co P \to A \kotimes{\Ground} P\grs{1})$.

Letting
\[
  \delta=\sum_{i=0}^\infty \delta^i\co P\to \overline{\Tensor}^*(A\grs{1})\kotimes{\Ground} P
  \qquad
  \mu_\bullet^\bullet =\sum_{n,w}X^w\cdot \mu_n^w \co \Tensor^*(A\grs{1})\to A\grs{2}
\]
then we can represent Equation~\eqref{eq:w-D-str} graphically by
\[
  \tikzsetnextfilename{def-wD-1}
  \mathcenter{\begin{tikzpicture}
    \node at (1,0) (tr) {};
    \node at (1,-1) (delta) {$\delta$};
    \node at (0,-2) (mu) {$\mu_\bullet^\bullet$};
    \node at (0,-3) (bl) {};
    \node at (1,-3) (br) {};
    \draw[dmoda] (tr) to (delta);
    \draw[dmoda] (delta) to (br);
    \draw[taa] (delta) to (mu);
    \draw[alga] (mu) to (bl);
  \end{tikzpicture}}
  =0.
\]
The structure equation~\eqref{eq:w-D-str} includes the
terms $\mu_0^w\otimes\Id$.

If $\kappa=0$, we can consider cases where $X=\One$, the
unit in $\Ground$, when the type $D$ structure relation takes the
simpler form
\begin{equation}
  \sum_{w=0}^\infty\sum_{n=0}^\infty
  (\mu_n^w\otimes\Id)\circ \delta^n=0\in A\kotimes{\Ground} P.
\end{equation}

Let $\wAlg$ be a weighted algebra, and 
suppose that $\lsup{\wAlg}\wP$ and $\lsup{\wAlg}\wQ$ are weighted type
$D$ structures (with either $\wAlg$ bonsai or $\lsup{\wAlg}\wP$ and
$\lsup{\wAlg}\wQ$ operationally bounded).
Endow the space $\Hom_{\Ground}(P, A\kotimes{\Ground} Q)$ with the endomorphism
$d$ 
given by 
\[
  d(\phi)=(\mu_\bullet^\bullet\otimes \Id_Q)\circ
  (\Id_{{\overline\Tensor}^*A}\otimes \delta_Q)\circ
  (\Id_{{\overline\Tensor}^*A}\otimes \phi)\circ \delta_P
\]
or, graphically,
\[
  \tikzsetnextfilename{d-D-phi-weighted-1}
  d(\phi)=
  \mathcenter{
  \begin{tikzpicture}
    \node at (0,0) (tl) {};
    \node at (1.5,0) (tr) {};
    \node at (1.5,-1) (delta1) {$\delta_P$};
    \node at (1.5,-2)(phi) {$\phi$};
    \node at (1.5,-3)(delta2){$\delta_Q$};
    \node at (0,-4.5) (mu) {$\mu_\bullet^\bullet$};
    \node at (0,-5.5) (bl) {};
    \node at (1.5,-5.5) (br) {};
    \draw[dmoda] (tr) to (delta1);
    \draw[dmoda] (delta1) to (phi);
    \draw[dmoda] (phi) to (delta2);
    \draw[dmoda] (delta2) to (br);
    \draw[taa] (delta1) to (mu);
    \draw[taa] (delta2) to (mu);
    \draw[alga] (phi) to (mu);
    \draw[alga] (mu) to (bl);
  \end{tikzpicture}}
=\sum_{i,j,w\geq 0} X^w
  \tikzsetnextfilename{d-D-phi-weighted-2}
  \mathcenter{
  \begin{tikzpicture}
    \node at (0,0) (tl) {};
    \node at (1.5,0) (tr) {};
    \node at (1.5,-1) (delta1) {$\delta_P^i$};
    \node at (1.5,-2)(phi) {$\phi$};
    \node at (1.5,-3)(delta2){$\delta_Q^j$};
    \node at (0,-4.5) (mu) {$\mu_{i+j+1}^w$};
    \node at (0,-5.5) (bl) {};
    \node at (1.5,-5.5) (br) {};
    \draw[dmoda] (tr) to (delta1);
    \draw[dmoda] (delta1) to (phi);
    \draw[dmoda] (phi) to (delta2);
    \draw[dmoda] (delta2) to (br);
    \draw[taa] (delta1) to (mu);
    \draw[taa] (delta2) to (mu);
    \draw[alga] (phi) to (mu);
    \draw[alga] (mu) to (bl);
  \end{tikzpicture}}
.
\]

Let $\Mor^{\wAlg}(\lsup{\wAlg}\wP,\lsup{\wAlg}\wQ)$ denote the graded $\Ring$-module
$\Hom_{\Ground}(P,A\kotimes{\Ground} Q)$, endowed with the endomorphism $d$.
We call this the space of {\em weighted morphisms from $\lsup{\wAlg}\wP$ to $\lsup{\wAlg}\wQ$}.

\begin{lemma}
  \label{lem:MorDiff}
  Fix a weighted algebra $\wAlg$ and type $D$ structures
  $\lsup{\wAlg}\wP$ and $\lsup{\wAlg}\wQ$ over $\wAlg$ with charge $X$. 
  If either $\wAlg$ is bonsai or $\lsup{\wAlg}\wP$ and
  $\lsup{\wAlg}\wQ$ are both operationally bounded then the morphism space
  $\Mor(\lsup{\wAlg}\wP,\lsup{\wAlg}\wQ)$ is a chain complex.
\end{lemma}

\begin{proof}
  This is an easy consequence of the weighted $\Ainf$ relation
  on $\wAlg$, and the type $D$ structure relations on
  $\lsup{\wAlg}\wP$ and $\lsup{\wAlg}\wQ$. That is,
  \begin{align*}
    d^2(\phi)&=
    \tikzsetnextfilename{d-mor-d-sq-1}
    \mathcenter{
      \begin{tikzpicture}[smallpic]
        \node at (0,0) (tl) {};
        \node at (1.5,0) (tr) {};
        \node at (1.5,-1) (delta0) {$\delta_P$};
        \node at (1.5,-2) (delta1) {$\delta_P$};
        \node at (1.5,-3)(phi) {$\phi$};
        \node at (1.5,-4)(delta2){$\delta_Q$};
        \node at (1.5,-5)(delta3){$\delta_Q$};
        \node at (0,-5) (mu1) {$\mu_\bullet^\bullet$};
        \node at (-1,-6) (mu2) {$\mu_\bullet^\bullet$};
        \node at (-1,-7) (bl) {};
        \node at (1.5,-7) (br) {};
        \draw[dmoda] (tr) to (delta0);
        \draw[dmoda] (delta0) to (delta1);
        \draw[dmoda] (delta1) to (phi);
        \draw[dmoda] (phi) to (delta2);
        \draw[dmoda] (delta2) to (delta3);
        \draw[dmoda] (delta3) to (br);
        \draw[taa] (delta1) to (mu1);
        \draw[taa] (delta2) to (mu1);
        \draw[alga] (phi) to (mu1);
        \draw[taa] (delta0) to (mu2);
        \draw[taa] (delta3) to (mu2);
        \draw[alga] (mu1) to (mu2);
        \draw[alga] (mu2) to (bl);
      \end{tikzpicture}}
               =
    \tikzsetnextfilename{d-mor-d-sq-2}
    \mathcenter{
      \begin{tikzpicture}[smallpic]
        \node at (0,0) (tl) {};
        \node at (1.5,0) (tr) {};
        \node at (1.5,-1) (delta0) {$\delta_P$};
        \node at (1.5,-2) (delta1) {$\delta_P$};
        \node at (1.5,-3)(delta2){$\delta_P$};
        \node at (1.5,-4)(phi) {$\phi$};
        \node at (1.5,-5)(delta3){$\delta_Q$};
        \node at (.5,-3.75) (mu1) {$\mu_\bullet^\bullet$};
        \node at (-1,-6) (mu2) {$\mu_\bullet^\bullet$};
        \node at (-1,-7) (bl) {};
        \node at (1.5,-7) (br) {};
        \draw[dmoda] (tr) to (delta0);
        \draw[dmoda] (delta0) to (delta1);
        \draw[dmoda] (delta1) to (delta2);
        \draw[dmoda] (delta2) to (phi);
        \draw[dmoda] (phi) to (delta3);
        \draw[dmoda] (delta3) to (br);
        \draw[taa] (delta1) to (mu1);
        \draw[taa] (delta2) to (mu2);
        \draw[alga] (phi) to (mu2);
        \draw[taa] (delta0) to (mu2);
        \draw[taa] (delta3) to (mu2);
        \draw[alga] (mu1) to (mu2);
        \draw[alga] (mu2) to (bl);
      \end{tikzpicture}}
+    \tikzsetnextfilename{d-mor-d-sq-3}
    \mathcenter{
      \begin{tikzpicture}[smallpic]
        \node at (0,0) (tl) {};
        \node at (1.5,0) (tr) {};
        \node at (1.5,-1) (delta0) {$\delta_P$};
        \node at (1.5,-2)(phi) {$\phi$};
        \node at (1.5,-3) (delta1) {$\delta_Q$};
        \node at (1.5,-4)(delta2){$\delta_Q$};
        \node at (1.5,-5)(delta3){$\delta_Q$};
        \node at (.5,-5) (mu1) {$\mu_\bullet^\bullet$};
        \node at (-1,-6) (mu2) {$\mu_\bullet^\bullet$};
        \node at (-1,-7) (bl) {};
        \node at (1.5,-7) (br) {};
        \draw[dmoda] (tr) to (delta0);
        \draw[dmoda] (delta0) to (phi);
        \draw[dmoda] (phi) to (delta1);
        \draw[dmoda] (delta1) to (delta2);
        \draw[dmoda] (delta2) to (delta3);
        \draw[dmoda] (delta3) to (br);
        \draw[taa] (delta1) to (mu2);
        \draw[taa] (delta2) to (mu1);
        \draw[alga] (phi) to (mu2);
        \draw[taa] (delta0) to (mu2);
        \draw[taa] (delta3) to (mu2);
        \draw[alga] (mu1) to (mu2);
        \draw[alga] (mu2) to (bl);
      \end{tikzpicture}}
               =0.
  \end{align*}

  Note that each term
  \[
  X^w(\mu_{i+j+1}^w\otimes \Id_Q)\circ
  (\Id\otimes \delta_Q^j)\circ
  (\Id\otimes \phi)\circ \delta_P^i  
  \]
  in $d(\phi)$ has grading
  \[
    \gr(\phi)-i-j+(-1+i+j+\kappa w)-\kappa w=\gr(\phi)-1.\qedhere
  \]
\end{proof}

Given a sequence of weighted type $D$ structures
$\lsup{\wAlg}\wP_0,\dots,\lsup{\wAlg}\wP_\ell$ (with $\ell\geq 1$),
there are composition maps
\[
  \underline{\mu}_\ell\co 
  \Mor(\lsup{\wAlg}\wP_{\ell-1},\lsup{\wAlg}\wP_\ell)\rotimes{\Ring}\dots\rotimes{\Ring}
  \Mor(\lsup{\wAlg}\wP_0,\lsup{\wAlg}\wP_1)
  \to\Mor(\lsup{\wAlg}\wP_0,\lsup{\wAlg}\wP_\ell)\grs{2-\ell},
\]
defined by
\[
  \underline{\mu}_\ell(\phi_\ell\otimes\dots\otimes \phi_1)=
  \sum_w X^w(\mu^{w}\otimes \Id_{M_{\ell}})\circ(\Id_{\Tensor^*(A)}\otimes \delta_{P_\ell}) 
  \circ(\Id_{\Tensor^*(A)}\otimes \phi_{\ell})\circ
  \dots
  \circ (\Id_{\Tensor^*(A)}\otimes \delta_{P_2})\circ (\Id_{\Tensor^*(A)}\otimes\phi_1)\circ
  (\delta_{P_1}).
\]
That is, graphically:
\[
  {\underline\mu}_\ell(\phi_\ell\otimes\dots\otimes \phi_1)=
  \tikzsetnextfilename{lem-MorDiff-2}
  \mathcenter{\begin{tikzpicture}
    \node at (-2,0) (tl) {};
    \node at (1.5,3) (tr) {};
    \node at (1.5,2) (delta1) {$\delta_{P_0}$};
    \node at (1.5,1)(phi1) {$\phi_1$};
    \node at (1.5,-3)(phiL) {$\phi_{\ell}$};
    \node at (1.5,-4)(deltaL) {$\delta_{P_\ell}$};
    \node at (1.5,0)(delta2){$\delta_{P_1}$};
    \node at (1.5,-2) (dots){$\vdots$};
    \node at (-2,-4.5) (mu) {$\mu_\bullet^\bullet$};
    \node at (-2,-5.5) (bl) {};
    \node at (1.5,-5.5) (br) {};
    \draw[dmoda] (tr) to (delta1);
    \draw[dmoda] (delta1) to (phi1);
    \draw[dmoda] (phi1) to (delta2);
    \draw[dmoda] (delta2) to (dots);
    \draw[dmoda] (dots) to (phiL);
    \draw[dmoda] (phiL) to (deltaL);
    \draw[dmoda] (deltaL) to (br);
    \draw[taa] (delta1) to (mu);
    \draw[taa] (delta2) to (mu);
    \draw[taa] (deltaL) to (mu);
    \draw[alga] (phi1) to (mu);
    \draw[alga] (phiL) to (mu);
    \draw[alga] (mu) to (bl);
  \end{tikzpicture}}\qquad=\qquad
\sum_w\ 
  X^w
  \tikzsetnextfilename{lem-MorDiff-3}
  \mathcenter{\begin{tikzpicture}
    \node at (-2,0) (tl) {};
    \node at (1.5,3) (tr) {};
    \node at (1.5,2) (delta1) {$\delta_{P_0}$};
    \node at (1.5,1)(phi1) {$\phi_1$};
    \node at (1.5,-3)(phiL) {$\phi_{\ell}$};
    \node at (1.5,-4)(deltaL) {$\delta_{P_\ell}$};
    \node at (1.5,0)(delta2){$\delta_{P_1}$};
    \node at (1.5,-2) (dots){$\vdots$};
    \node at (-2,-4.5) (mu) {$\mu_\bullet^w$};
    \node at (-2,-5.5) (bl) {};
    \node at (1.5,-5.5) (br) {};
    \draw[dmoda] (tr) to (delta1);
    \draw[dmoda] (delta1) to (phi1);
    \draw[dmoda] (phi1) to (delta2);
    \draw[dmoda] (delta2) to (dots);
    \draw[dmoda] (dots) to (phiL);
    \draw[dmoda] (phiL) to (deltaL);
    \draw[dmoda] (deltaL) to (br);
    \draw[taa] (delta1) to (mu);
    \draw[taa] (delta2) to (mu);
    \draw[taa] (deltaL) to (mu);
    \draw[alga] (phi1) to (mu);
    \draw[alga] (phiL) to (mu);
    \draw[alga] (mu) to (bl);
  \end{tikzpicture}}.
\]

If $\wAlg$ is strictly unital with unit $\unit$, we can define the
identity map of the type $D$ structure $\lsup{\wAlg}\wP$ by
$\Id(x)=\unit\otimes x$. 

\begin{proposition}
  \label{prop:TypeDStructureswAnfCat}
  Given a bonsai weighted $\Ainf$-algebra $\wAlg$, the set of
  weighted type $D$ structures, with morphism spaces given by $\Mor$,
  and composition maps defined by ${\underline\mu}_\ell$, forms a (nonunital)
  $\Ainf$-category. That is, the composition maps
  ${\underline\mu}_\ell$ have grading $2-\ell$ and satisfy the
  $\Ainf$ relations. If $\wAlg$ is strictly unital then the category
  of type $D$ structures is also unital.
  The statements also hold for the category of
  operationally bounded type $D$ structures over a non-bonsai algebra.
\end{proposition}

\begin{proof}
  The proof is essentially the same as in the unweighted
  case~\cite[Lemma 2.2.27]{LOT2}.
  The statement about gradings is arithmetic.
  The $\Ainf$ relation counts contributions of trees
  which consist of a sequence of operations on the $P_i$, for
  $i=1,\dots,\ell$ starting with $\delta_{P_i}$, alternating between
  $\delta_{P_i}$, $\phi_i$, and ending with $\delta_{P_\ell}$; some
  consecutive sequence of the algebra outputs are channeled into a
  vertex labelled by $\mu^\bullet_\bullet$, and then that output,
  and all others, are channeled into a second vertex labelled by
  $\mu^\bullet_\bullet$.  For the $\Ainf$ relation on the category,
  we consider only those trees for which the first $\mu^\bullet_\bullet$-vertex
  has an input from at least one $\phi_i$. Trees which do not have
  this property have all the inputs channeled into the first
  $\mu^\bullet_\bullet$ which come from a single $\delta_{P_i}$; their
  contribution vanishes by the type $D$ structure relation on $P_i$.
  The sum of all the trees is zero by the weighted $\Ainf$ relation
  on the algebra. It is immediate from the definition of strict
  unitality of a weighted algebra that, in the strictly unital case,
  the identity maps are strict units for composition.
\end{proof}

\begin{remark}
  There is also a notion of homotopy unital weighted $\Ainf$-algebras,
  and the category of type $D$ structures over a homotopy unital
  weighted $\Ainf$-algebra is a homotopy unital $\Ainf$-category. See
  Sections~\ref{sec:hu} and~\ref{sec:hu-D}.
\end{remark}

\subsection{Box products}
\label{sec:wD1}

Given $X\in\Ground$, a $w$-algebra $\wAlg$ over $\Ground$ so that $X$
acts centrally on $A$, a type $D$
structure $P$ over $\wAlg$ with charge $X$, and a $w$-module $\wMod$
over $\wAlg$, we define $\wMod\DT P$ to be the graded $\Ring$-module
$M\kotimes{\Ground} P$ with differential
\begin{equation}\label{eq:w-DT-def}
  \bdy = \sum_{n=0}^\infty\sum_{w=0}^\infty X^w(m_{1+n}^w\otimes\Id)\circ(\Id\otimes\delta^n).
\end{equation}
This sum makes sense if either:
\begin{itemize}
\item $\wMod$ and $\wAlg$ are bonsai or
\item $P$ is operationally bounded.
\end{itemize}

We may represent Equation~\eqref{eq:w-DT-def} graphically as
\[
  \tikzsetnextfilename{D-weighted-box-prod-1}
  \mathcenter{\begin{tikzpicture}
    \node at (0,0) (tl) {};
    \node at (1,0) (tr) {};
    \node at (1,-1) (delta) {$\delta$};
    \node at (0,-2) (mu) {$m_\bullet^\bullet$};
    \node at (0,-3) (bl) {};
    \node at (1,-3) (br) {};
    \draw[dmoda] (tr) to (delta);
    \draw[dmoda] (delta) to (br);
    \draw[taa] (delta) to (mu);
    \draw[moda] (tl) to (mu);
    \draw[moda] (mu) to (bl);
  \end{tikzpicture}}
  =
  \sum_{n,w\geq 0}
  X^w
  \tikzsetnextfilename{D-weighted-box-prod-2}
    \mathcenter{\begin{tikzpicture}
    \node at (0,0) (tl) {};
    \node at (1,0) (tr) {};
    \node at (1,-1) (delta) {$\delta^n$};
    \node at (0,-2) (mu) {$m_{n+1}^w$};
    \node at (0,-3) (bl) {};
    \node at (1,-3) (br) {};
    \draw[dmoda] (tr) to (delta);
    \draw[dmoda] (delta) to (br);
    \draw[taa] (delta) to (mu);
    \draw[moda] (tl) to (mu);
    \draw[moda] (mu) to (bl);
  \end{tikzpicture}}.
\]

\begin{lemma}\label{lem:w-DT-sq-0}
  Let $\wAlg$ be a weighted $\Ainf$-algebra, $\wMod_{\wAlg}$ a
  weighted $\Ainf$-module, and $\lsup{\wAlg}P$ a weighted type $D$
  structure with some charge $X$. Assume that either $\wMod$ and
  $\wAlg$ are bonsai or $P$ is operationally bounded. Then the map
  $\bdy$ on $\wMod\DT P$ has degree $-1$ and satisfies $\bdy^2=0$.
\end{lemma}
\begin{proof}
  The proof is the same as in the unweighted case~\cite[Lemma 2.30]{LOT1}.
\end{proof}

We have the following analogue of the functoriality part of Lemma~\ref{lem:DT-mod-defined}:
\begin{lemma}\label{lem:w-DT-bifunc}
  The operation $\DT$ extends to a (non-unital) $\Ainf$-bifunctor from the
  categories of weighted $\Ainf$-modules and operationally bounded
  type $D$ structures with charge $X$ to the category of chain complexes. It also defines an $\Ainf$-bifunctor from the categories of
  bonsai weighted $\Ainf$-modules and bonsai maps and arbitrary type
  $D$ structures with charge $X$ to the category of chain complexes.
  For strictly unital weighted algebras and modules, $\DT$ is a strictly
  unital functor, i.e., $\Id\DT\Id=\Id$.
  In particular, in the strictly unital case, the operation $\DT$
  respects homotopy equivalences of $\wMod$ and $P$. 
\end{lemma}
\begin{proof}
  The proof of the first statement is the same as in the unweighted
  case~\cite[Lemmas 2.3.3 and 2.3.7]{LOT2}. The second statement is
  immediate from the definitions. The third statement follows
  immediately from the other two.
\end{proof}

\subsection{Homotopy unital weighted algebras and modules}
\label{sec:hu}
Unlike the unweighted case, the tensor product of weighted algebras makes use of
a unit. In Sections~\ref{sec:wDiags} and~\ref{sec:wDiagApps} we will focus
mainly on the case that the weighted algebras are strictly unital, but in order
to prove associativity of the tensor product we need to relax this to allow
homotopy unital weighted algebras. This is a weighted extension of (a simple
variant) of a well-known notion for $\Ainf$-algebras~\cite[Section 3.3]{FOOO1}. As
usual, we will define the notion using a complex of decorated trees; in the
unweighted case, this is essentially the cellular chains on Muro-Tonks's
unital associahedra~\cite{MuroTonks14:unital-assoc} (see
also~\cite{Lyubashenko11:htpy-unital-operad}).

\begin{definition}\label{def:wu-trees-cx}
  A \emph{thorn tree} consists of:
  \begin{itemize}
  \item A planar, rooted tree $T$.
  \item A partition of the leaves of $T$ into the following five sets:
    The \emph{output}, which is the root, the \emph{inputs}, the
    \emph{thorns}, the \emph{stumps}, and the \emph{popsicles}. Of
    these, we regard only popsicles as internal vertices.
  \item A weight function from the internal vertices of $T$ to $\ZZ_{\geq 0}$.
  \end{itemize}
  These are required to satisfy the conditions that:
  \begin{itemize}
  \item Each internal vertex either has valence at least $3$ or weight $>0$.
  \item The tree has at least one vertex which is not a thorn or the
    output. (This disallows exactly one tree.)
  \end{itemize}

  Composition of thorn trees is defined in the obvious way.

  The differential of a thorn tree $T$ is the sum of the following terms:
  \begin{itemize}
  \item all ways of splitting a single internal vertex in $T$ into two vertices
    to obtain a new thorn tree,
  \item all ways of replacing a single thorn in $T$ with a stump, and
  \item all ways of deleting a single thorn in $T$ which is connected
    to a $2$-input
    ($3$-valent), weight $0$ vertex.
\end{itemize}

  The \emph{dimension} of a thorn tree $T$ is
  \[
    \dim(T)=n+2w+t-v-1
  \]
  where $n$ is the number of inputs, $w$ is the total weight, $t$ is
  the number of thorns, and $v$ is the number of internal vertices.

  The \emph{homotopy unital weighted trees complex} $\uwTreesCx{*}{*}$ is the
  free chain complex over $\Ring$ generated by the thorn trees. There is a
  tri-grading on the homotopy unital weighted trees complex by the grading, the number of inputs,
  and the total weight. The summand with $n$ inputs and weight $w$ is denoted
  $\uwTreesCx{n}{w}$
\end{definition}

\begin{figure}
  \centering
  %Font is 12 point.
  \includegraphics{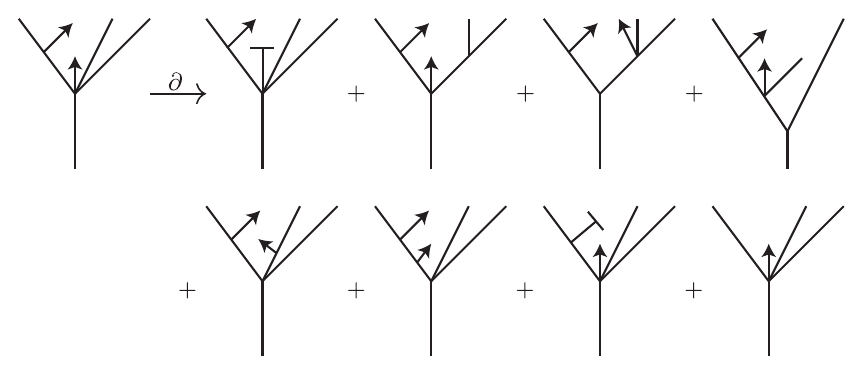}
  \caption[A thorn tree and its differential]{\textbf{A thorn tree and its differential.} The thorn tree on the
    left has 3 inputs, no stumps, and two thorns. We have drawn an example in
    which all vertices have weight $0$, so we have suppressed the weights.}
  \label{fig:thorn-tree}
\end{figure}

See Figure~\ref{fig:thorn-tree}. Note that there is a weight $0$ thorn
tree with no one input and no internal vertices, the \emph{identity
  tree} $\IdTree$, which is the identity for composition of
trees. Abusing terminology, we will often call the $0$-input tree with
no internal vertices or thorns and a single stump the
\emph{stump}. Composing a tree with the stump adds a stump to the tree.

\begin{lemma}\label{lem:uwTrees-is-cx}
  The $n$-input, weight $w$ homotopy unital weighted trees complex
  $\uwTreesCx{n}{w}$ is a chain complex, and the differential
  decreases the grading (by the dimension $\dim(T)$) by one. Further, the composition map
  $\circ_i\co \uwTreesCx{n_1}{w_1}\times
  \uwTreesCx{n_2}{w_2}\to\uwTreesCx{n_1+n_2-1}{w_1+w_2}$ is a chain
  map.
\end{lemma}
\begin{proof}
  As in Lemma~\ref{lem:wTrees-is-cx}, to prove that $\bdy^2=0$ we
  consider the dual complex. In the dual complex, the differential is
  the sum of all ways of contracting an edge, turning a stump into a
  thorn, or adding a thorn to an edge. Each of the following kinds of
  terms in $\delta^2$ cancels in pairs:
  \begin{itemize}
  \item Contracting two edges.
  \item Turning a stump into a thorn and contracting an edge.
  \item Turning two stumps into thorns.
  \item Turning a stump into a thorn and adding a thorn to an edge.
  \item Adding two thorns.
  \item Contracting an edge and adding a thorn to a different edge.
  \end{itemize}
  The remaining terms correspond to adding a thorn to an edge $e$,
  which splits $e$ into two edges $e_1$ and $e_2$ and then contracting
  one of $e_1$ or $e_2$. This has the same effect as adding a thorn at
  one of the vertices adjacent to $e$. Further, each way of adding a
  thorn at a vertex occurs twice, corresponding to the two edges
  adjacent to the thorn. Thus, terms in $\delta^2$ of this form also
  cancel in pairs.
  
  The facts that the differential decreases the grading by $1$ and
  that composition is a chain map are immediate from the definitions.
\end{proof}

\begin{definition}\label{def:hu-w-alg}
  A \emph{homotopy unital weighted algebra} $\wAlg$ over $\Ground$ consists of a
  \dg $\Ground$-bimodule $A$ and a chain map 
  $\mu\co \uwTreesCx{n}{w}\to \Mor(A^{\kotimes{\Ground} n},A\grs{(2-\kappa)w})$
  for each $n,w$, so that $\mu(S\circ_i T)=\mu(S)\circ_i\mu(T)$.

  Given a homotopy unital weighted algebra $\wAlg$, the image of the
  stump $\stump$ is the \emph{homotopy unit} in $\wAlg$, and is denoted
  $\One$.

  A homotopy unital weighted algebra $\wAlg$ is \emph{split homotopy
    unital} if there is a $\Ring$-module splitting
  $A=\Ring\langle\unit\rangle\oplus A'$.
\end{definition}

There is an inclusion from the weighted trees complex into the homotopy
unital weighted trees complex, so we can regard any homotopy unital
weighted algebra as an ordinary weighted algebra. In other words, a
homotopy unital weighted algebra has an \emph{underlying weighted
  algebra}. Conversely, given a
strictly unital weighted algebra $\wAlg$ there is an induced homotopy
unital weighted algebra as follows. Given a thorn tree $T$, with $>0$
thorns, $\mu(T)=0$. For a thorn tree $T$ with no thorns, $\mu(T)$ is
obtained by inserting the unit $\One$ at each stump of $T$, and then
applying the $k$-input, weight $w$ operation from $\wAlg$ at each
$k$-input, weight $w$ vertex of $T$. It is straightforward to verify
that this defines a chain map.

\begin{figure}
  \centering
  \includegraphics{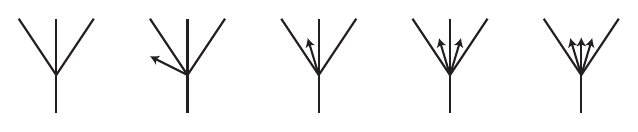}
  \caption[Notation for weighted corollas with thorns]{\textbf{Notation for weighted corollas with thorns.} From
    left to right, $\wcorolla{3}{0}$, $\wcorolla{0\uparrow 3}{0}$,
    $\wcorolla{1\uparrow 2}{0}$, $\wcorolla{1\uparrow 1\uparrow
      1}{0}$, $\wcorolla{1\uparrow 0\uparrow0\uparrow1}{0}$.}
  \label{fig:thorn-corolla-not}
\end{figure}

A homotopy unital weighted algebra is determined by the unit $\One$
and the operations corresponding to a corolla with some thorns coming
out. That is, given a sequence of non-negative integers $n_1,\dots,n_k$ there is a
tree $\wcorolla{n_1\uparrow n_2\uparrow\cdots\uparrow n_k}{w}$ with
one internal vertex with weight $w$, with $n_1$ inputs followed by a
thorn followed by $n_2$ inputs, another thorn, and so on. See
Figure~\ref{fig:thorn-corolla-not}. (The operations on the underlying
weighted algebra correspond to the case $k=1$.) The operations
\[
  \mu_{n_1\uparrow\cdots \uparrow n_k}^w=\mu(\wcorolla{n_1\uparrow
    n_2\uparrow\cdots\uparrow n_k}{w})\co A^{\kotimes{\Ground}
    (n_1+\cdots+n_k)}\to A,
\]
together with the unit $\One$, determine the homotopy unital weighted
algebra. (We have suppressed the grading shift in this formula.)
It is not hard to rephrase the structure relations for a
homotopy unital weighted algebra in terms of the unit $\One$ and the
operations $\mu_{n_1\uparrow\cdots\uparrow n_k}^w$.

Given a corolla
$\wcorolla{n_1\uparrow n_2\uparrow\cdots\uparrow n_k}{w}$, with
$n=n_1+\cdots+n_k$, there is a monotone injection
$\sigma_{n_1\uparrow\cdots\uparrow n_k}\co \{1,\dots,n\}\into
\{1,\dots,n+k-1\}$ with image
$\{1,\dots,n+k-1\}\setminus
\{n_1+1,n_1+n_2+2,\dots,n_1+\cdots+n_{k-1}+k-1\}$. For the operation
$\mu_{n_1\uparrow\cdots\uparrow n_k}^w(a_1,\dots,a_n)$ we call
$\sigma_{n_1\uparrow\cdots\uparrow n_k}(i)$ the \emph{apparent
  position} of $a_i$.

\begin{definition}\label{def:hu-w-mod}
  The \emph{homotopy unital weighted module trees complex} is the
  subquotient of $\uwTreesCx{*}{*}$ obtained by requiring that there
  are no stumps, popsicles, or thorns to the left of the leftmost
  input.

  A \emph{homotopy unital weighted module} $\wMod$ over a homotopy unital
  weighted algebra $\wAlg$ consists of a chain complex $M$ over
  $\Ground$ and grading-preserving chain maps
  $m\co \uwMTreesCx{n}{w}\to \Mor(M\kotimes{\Ground} A^{\kotimes{\Ground} n-1},M\grs{(2-\kappa)w}) $ which
  respect composition in the obvious sense.
\end{definition}

Similarly to algebras, homotopy unital weighted modules are determined
by the operations
\[
  \mu_{1+n_1\uparrow\cdots\uparrow n_k}^w\co M\kotimes{\Ground} A^{\kotimes{\Ground} (n_1+\cdots+n_k)}\to A
\]
corresponding to a corolla with valence $n_1+\cdots+n_k+k$ with a mix
of thorn and input leaves, where the thorns are at the positions
$n_1+2,n_1+n_2+3,\dots$. (Again, we have suppressed the grading in this formula.)

As in the case of algebras, there is a forgetful map from homotopy
unital weighted modules to weighted modules. Further, given a strictly
unital weighted module over a strictly unital weighted algebra, there
is an induced homotopy unital weighted module, by feeding the unit
$\One$ in at the stumps and declaring trees with thorns to act by
zero.

\begin{figure}
  \centering
  %Font is 12 point
  \includegraphics{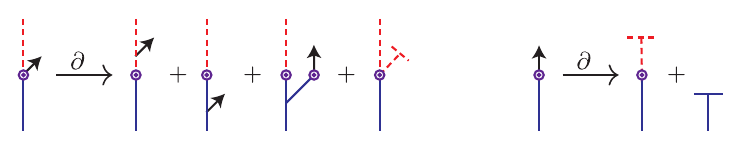}
  \caption[Transformation thorn trees and their differentials]{\textbf{Transformation thorn trees and their
      differentials.} Two examples are shown.}
  \label{fig:trans-thorn-diff}
\end{figure}

\begin{definition}
  A \emph{transformation thorn tree} consists of:
  \begin{itemize}
  \item A planar, rooted tree $T$.
  \item A partition of the leaves of $T$ into the following five sets:
    The \emph{output}, which is the root, the \emph{inputs}, the
    \emph{thorns}, the \emph{stumps}, and the \emph{popsicles}. Of
    these, we regard only popsicles as internal vertices.
  \item A coloring of each edge of $T$ not incident to a thorn as
    either red or blue, subject to the following conditions:
    \begin{enumerate}
    \item The edges adjacent to input leaves are red.
    \item The edge adjacent to the output leaf is blue.
    \item For each vertex $v$, all of the colored inputs of $v$ have
      the same color (red or blue).
    \item If a vertex $v$ has a red output, then all of the colored
      inputs of $v$ are red; if a vertex has blue inputs, then its
      output is also blue.
    \end{enumerate}
    (Compare Definition~\ref{def:TransformationTree}.)
  \item A coloring of the non-thorn vertices of $T$ as \emph{red}, \emph{blue}, or \emph{purple}, subject to the conditions that
    \begin{enumerate}
    \item if a vertex has any red inputs then it is red or purple,
    \item if a vertex has a blue output then it is blue or purple,
    \item if a vertex has a red output then it is red,
    \item if a vertex has a blue input then it is blue, and
    \item no purple vertex is above another purple vertex. 
    \end{enumerate}
    (These conditions determine the color of every vertex whose inputs
    are not all thorns.)
  \item A weight function from the internal vertices of $T$ to
    $\ZZ_{\geq 0}$.
  \end{itemize}
  These are required to satisfy the condition that each internal
  vertex either has valence at least $3$ or has weight $>0$ or is
  colored purple.

  The differential of a transformation thorn tree is the sum of all
  ways of obtaining a new transformation thorn tree by doing one of the following:
  \begin{itemize}
  \item splitting a blue (respectively red) vertex into two blue
    (respectively red) vertices,
  \item splitting a purple vertex into a red vertex feeding into a
    purple vertex,
  \item splitting a purple vertex into a layer of purple vertices
    feeding into a blue vertex, with each thorn ending up coming from
    either the blue vertex or one of the purple vertices,
  \item replacing a thorn by a stump,
  \item deleting a thorn adjacent to a $2$-input ($3$-valent), weight $0$ blue or red vertex,
    or
  \item replacing a $2$-valent, weight $0$ purple vertex whose input
    is a thorn by a blue stump.
  \end{itemize}
  See Figure~\ref{fig:trans-thorn-diff}.
  
  The dimension (grading) of a weighted thorn tree $T$ with $n$ inputs, total
  weight $w$, $t$ thorns, $r$ internal red vertices, and $b$
  internal blue vertices is
  \[
    \dim(T)=n+2w+t-r-b-1.
  \]
  
  The \emph{homotopy unital weighted transformation trees complex}
  $\uwTransCx{n}{w}$ is the complex generated by all transformation
  thorn trees with $n$ inputs and total weight $w$.
\end{definition}

As in the non-unital setting, there are composition maps
\begin{align}
  \circ_i\co \uwTransCx{n_1}{w_1}\rotimes{\Ring}
  \uwTreesCx{n_2}{w_2}&\to\uwTransCx{n_1+n_2-1}{w_1+w_2},\label{eq:uw-trans-compo-1}\\
  \circ\co \uwTreesCx{n}{w}\rotimes{\Ring}\uwTransCx{m_1}{v_1}\rotimes{\Ring}\cdots\rotimes{\Ring} \uwTransCx{m_n}{v_n}
                     &\to\uwTransCx{m_1+\cdots+m_n}{v_1+\cdots+v_n+w}\label{eq:uw-trans-compo-2}.
\end{align}

\begin{lemma}\label{lem:wTrans-is-cx-hu}
  The $n$-input, weight $w$ homotopy unital weighted transformation
  trees complex $\uwTransCx{n}{w}$ is a chain complex, and the
  differential decreases the grading by one. Further, the composition
  maps~\eqref{eq:uw-trans-compo-1} and~\eqref{eq:uw-trans-compo-2} are
  chain maps. 
\end{lemma}
\begin{proof}
  The proof is similar to, but more complicated than, the proof of
  Lemma~\ref{lem:wTrans-is-cx}. In the dual complex to
  $\uwTransCx{n}{w}$, the differential~$\delta$ is the sum of all ways of:
  \begin{enumerate}[label=(\arabic*)]
  \item\label{item:cont-red-red-hu} contracting an edge between two red vertices,
  \item\label{item:cont-blue-blue-hu} contracting an edge between two blue vertices,
  \item\label{item:cont-purp-red-hu} contracting an edge between a purple vertex and a red vertex,
  \item\label{item:cont-blue-purp-hu} contracting all the colored
    edges into a blue vertex $v$ if all of the vertices above $v$ are
    purple or un-colored (thorns),
  \item\label{item:add-thorn-hu} inserting a new $3$-valent, weight
    $0$ vertex on an edge, with one of its inputs a new thorn,
  \item\label{item:stump-to-thorn-hu} replacing a (red or blue) stump with a thorn, or
  \item\label{item:stump-to-vert-and-thorn-hu} replacing a blue stump
    with a thorn feeding into a $2$-valent, weight $0$ purple vertex.
  \end{enumerate}
  Terms in $\delta^2$ involving only
  operations~\ref{item:cont-red-red-hu},~\ref{item:cont-blue-blue-hu},~\ref{item:cont-purp-red-hu},
  and~\ref{item:cont-blue-purp-hu} cancel as in the proof of
  Lemma~\ref{lem:wTrans-is-cx}. Terms coming from two operations of
  types~\ref{item:add-thorn-hu},~\ref{item:stump-to-thorn-hu},
  or~\ref{item:stump-to-vert-and-thorn-hu} cancel in pairs in an
  obvious way. Similarly, terms in $\delta^2$ coming from a pair of
  operations
  in $\{$\ref{item:cont-red-red-hu},\ref{item:cont-blue-blue-hu},\ref{item:cont-purp-red-hu},\ref{item:stump-to-thorn-hu},\ref{item:stump-to-vert-and-thorn-hu}$\}$
  cancel in pairs in an obvious way. This leaves terms coming from the
  following pairs:
  \begin{itemize}
  \item (\ref{item:cont-red-red-hu},\ref{item:add-thorn-hu}) and
    (\ref{item:cont-blue-blue-hu},\ref{item:add-thorn-hu}). These
    cancel in pairs as in the proof of Lemma~\ref{lem:uwTrees-is-cx}.
  \item (\ref{item:cont-purp-red-hu},\ref{item:add-thorn-hu}). When
    the red vertex in the type~\ref{item:cont-purp-red-hu} operation is the
    3-valent vertex created by the type~\ref{item:add-thorn-hu}
    operation, these cancel against terms of type
    (\ref{item:cont-blue-purp-hu},\ref{item:add-thorn-hu}). The rest
    of these terms cancel in pairs.
  \item (\ref{item:cont-blue-purp-hu},\ref{item:add-thorn-hu}). When the
    blue vertex in the type~\ref{item:cont-blue-purp-hu} operation is the
    $3$-valent vertex created by the type~\ref{item:add-thorn-hu}
    operation, these cancel against terms of
    type (\ref{item:cont-purp-red-hu},\ref{item:add-thorn-hu}). The
    rest of these terms cancel in pairs.
  \item (\ref{item:cont-blue-purp-hu},\ref{item:stump-to-thorn-hu}).
    When the stump in the type~\ref{item:stump-to-thorn-hu} operation
    is incident to the blue vertex in the
    type~\ref{item:cont-blue-purp-hu} operation, these cancel
    against terms of type
    (\ref{item:cont-blue-purp-hu},\ref{item:stump-to-vert-and-thorn-hu}).
    The rest of these terms cancel in pairs.
  \item (\ref{item:cont-blue-purp-hu},\ref{item:stump-to-vert-and-thorn-hu}).
    When the purple vertex in the
    type~\ref{item:stump-to-vert-and-thorn-hu} operation
    feeds into the blue vertex in the
    type~\ref{item:cont-blue-purp-hu} operation, these cancel
    against terms of type
    (\ref{item:cont-blue-purp-hu},\ref{item:stump-to-thorn-hu}).
    The rest of these terms cancel in pairs.
  \end{itemize}
  This completes the proof that $\delta^2=0$.  The facts that the
  differential decreases the grading by $1$ and that the composition
  maps are chain maps are immediate from the definitions.
\end{proof}

\begin{definition}
  Given homotopy unital weighted algebras $\wAlg$ and $\wBlg$, a
  \emph{homotopy unital homomorphism} from $\wAlg$ to $\wBlg$ consists of
  chain maps
  \[
    \uwTransCx{n}{w}\to \Mor(A^{\kotimes{\Ground} n},B\grs{(2-\kappa)w})
  \]
  for each $n$ and $w$ which are compatible with composition in the
  obvious sense.
\end{definition}

\begin{figure}
  \centering
  %Font is 12 point.
  \includegraphics{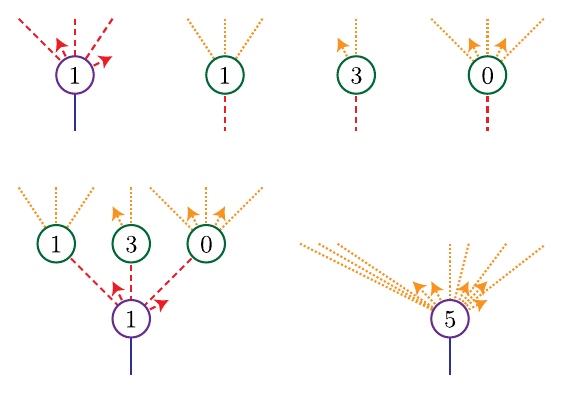}
  \caption[The collapsed composition of homotopy unital purple
      corollas]{\textbf{The collapsed composition of homotopy unital purple
      corollas.} Top: four homotopy unital purple corollas. In order,
    these are $\wpcorolla{1\uparrow2\uparrow0}{1}$, $\wpcorolla{3}{1}$,
    $\wpcorolla{0\uparrow1}{3}$, and
    $\wpcorolla{1\uparrow1\uparrow1}{0}$. Bottom: their composition
    (left) and collapsed composition (right), which is $\wpcorolla{3\uparrow0\uparrow2\uparrow1\uparrow1\uparrow0}{5}$.}
  \label{fig:collapsed-comp}
\end{figure}

A homotopy unital homomorphism $f$ is determined by maps
\[
  f_{n_1\uparrow\cdots\uparrow n_k}^w\co A^{\kotimes{\Ground} n_1+\cdots+n_k}\to B
\]
corresponding to a purple corolla
$\wpcorolla{n_1\uparrow\cdots\uparrow n_k}{w}$ with valence $n+t+1$ with a mix of thorn and
input leaves; see Figure~\ref{fig:collapsed-comp}. (We have suppressed the grading.)
Again, there is a forgetful map from homotopy unital homomorphisms to
homomorphisms of underlying non-unital weighted algebras. A homomorphism is a
\emph{quasi-isomorphism} if the induced map of non-unital weighted algebras is a
quasi-isomorphism.

\begin{example}
  The \emph{identity homomorphism} of a homotopy unital algebra $\wAlg$ is the map
  $\Id_{\wAlg}$ with $(\Id_{\wAlg})_{1}^0=\Id_A$ and $(\Id_{\wAlg})_{n_1\uparrow\cdots\uparrow n_k}^w=0$ otherwise (i.e., if $w>0$ or $k>1$ or $n_1>1$).
  It is immediate from the definitions that this defines a homotopy unital
  homomorphism.
\end{example}

Given non-negative integers $n_1,\dots,n_k$,
$m_{i,1},\dots,m_{i,\ell_i}$ (for $i=1,\dots,n\coloneqq n_1+\cdots+n_k$), $w$, and $v_1,\dots,v_n$, there is a \emph{collapsed
  composition} of $\wpcorolla{n_1\uparrow\cdots\uparrow n_k}{w}$ with
$(\wpcorolla{m_{1,1}\uparrow\cdots\uparrow
  m_{i,\ell_i}}{v_1},\dots,\wpcorolla{m_{n,1}\uparrow\cdots\uparrow
  m_{n,\ell_n}}{v_n})$ by composing $\wpcorolla{n_1\uparrow\cdots\uparrow n_k}{w}$ with
$(\wpcorolla{m_{1,1}\uparrow\cdots\uparrow
  m_{i,\ell_i}}{v_1},\dots,\wpcorolla{m_{n,1}\uparrow\cdots\uparrow
  m_{n,\ell_n}}{v_n})$ and then collapsing all of the internal edges
to obtain a new corolla. In formulas, the collapsed composition is
\[
  \wpcorolla{m_{1,1}\uparrow m_{1,2}\uparrow \cdots\uparrow
    m_{1,\ell_1}+m_{2,1}\uparrow m_{2,2}\uparrow\cdots\uparrow
    m_{2,\ell_2}+m_{3,1}\uparrow\cdots\uparrow
    m_{n_1,\ell_{n_1}}\uparrow m_{n_1+1,1}\uparrow m_{n_1+1,2}\uparrow\cdots}{v_1+\cdots+v_n+w}.
\]
See Figure~\ref{fig:collapsed-comp}.

\begin{definition}
  Given homotopy unital homomorphisms $f$ and $g$, their composition is defined by
  \[
    (g\circ f)_{p_1\uparrow\cdots\uparrow p_\ell}^w=\sum
    g_{n_1\uparrow\cdots\uparrow n_k}^v\circ(f_{m_{1,1},\uparrow
      \cdots\uparrow m_{1,j_1}}^{u_1}\otimes\cdots\otimes f_{m_{n,1},\uparrow
      \cdots\uparrow m_{n,j_n}}^{u_n}),
  \]
  where the sum is over data
  $n_1,\dots,n_k,m_{1,1},\dots,m_{n,j_n},v,u_1,\dots,u_n$
  whose collapsed composition is $\wpcorolla{p_1\uparrow\cdots\uparrow
  p_\ell}{w}$. (Here, $n=n_1+\cdots+n_k$.)
\end{definition}
Compare Definition~\ref{def:AlgebraHomomorphism}.
\begin{lemma}
  The composition of two homotopy unital algebra homomorphisms is a
  homotopy unital algebra homomorphism.
\end{lemma}
\begin{proof}
  The proof is left to the reader.
\end{proof}

A homotopy unital homomorphism $f\co \wAlg\to\wBlg$ is an
\emph{isomorphism} if there is a homotopy unital homomorphism $g\co
\wBlg\to\wAlg$ so that $g\circ f=\Id_{\wAlg}$ and $f\circ g=\Id_{\wBlg}$.
We have the following analogue of Lemma~\ref{lem:wAlg-iso-is}:
\begin{lemma}\label{lem:hu-Alg-iso-is}
  A homotopy unital algebra homomorphism $f\co \wAlg\to\wBlg$ is an
  isomorphism if and only if $f^0_1$ is an isomorphism.
\end{lemma}
\begin{proof}
  As in the proof of Lemma~\ref{lem:wAlg-iso-is}, we construct left
  and right inverses of $f$ inductively, and then verify that they are
  homomorphisms afterwards. So, we first construct maps
  $g^w_{n_1\uparrow\dots\uparrow n_k}$ formally satisfying
  $g \circ f = \Id$. Let $g^0_1 = \bigl(f^0_1\bigr)^{-1}$, which
  exists by hypothesis.  Construct the rest of $g$ by induction on the
  weight~$w$, number of thorns~$(k-1)$, and the number of non-thorn
  inputs $\sum n_i$, in that order. Inspecting the terms that
  contribute to $(g \circ f)^w_{n_1\uparrow \dots \uparrow n_k}$,
  there is one maximal term, namely
  $g^w_{n_1 \uparrow \dots \uparrow n_k} \circ (f^0_1 \otimes\dots f^0_1)$,
  and other terms that involve previously-defined values of~$g$. Since
  $f^0_1$ is invertible, we thus determine the new term uniquely.

  We can similarly construct maps $h^w_{n_1\uparrow\dots\uparrow n_k}$
  so that $f \circ h = \Id$. Then the usual argument (reducing $g \circ f
  \circ h$ in two ways) shows that $g = h$.

  To see $g$ is a homomorphism, we modify the proof of
  Lemma~\ref{lem:wAlg-iso-is}. The relation in the case that the input
  is a single thorn is easily seen to be satisfied. To verify
  the remaining cases, we adopt the following notation.  Given a map
  $h\co \bigoplus_n A^{\otimes n}\to B$, let
  \[
    \underline{h}=\sum_{m=1}^\infty
    (\overbrace{h\otimes\cdots\otimes h}^m)\co \bigoplus_n A^{\otimes n}\to
    B^{\otimes m}.
  \]
  Let
  \[
    \tikzsetnextfilename{hu-Alg-iso-is-1}
    \mathcenter{
      \begin{tikzpicture}[smallpic]
        \node at (0,0) (tc) {};
        \node at (0,-1) (F) {${f}^*$};
        \node at (0,-2) (bc) {};
        \draw[taa] (tc) to (F);
        \draw[blga] (F) to (bc);
      \end{tikzpicture}
    }
  \]
  denote feeding an element of the tensor algebra on $A$ into
  $f(\wpcorolla{k_1\uparrow k_2\uparrow\cdots \uparrow k_\ell}{w})$ for
  some $\ell$, $w$, and sequence $k_1,\dots,k_\ell$. Extend the notation
  $\underline{f}^*$ and $\mu^*$ similarly, by again
  interspersing thorns in an arbitrary way. Then the homotopy unital
  homomorphism relation for $g$ (except for the case with a single
  thorn as input) is
  \[
    \tikzsetnextfilename{hu-Alg-iso-is-2}
    \mathcenter{
      \begin{tikzpicture}
        \node at (0,0) (tc) {};
        \node at (0,-1) (g) {$\underline{g}^*$};
        \node at (0,-2) (mu) {$\mu^*$};
        \node at (0,-3) (bc) {};
        \draw[tbb] (tc) to (g);
        \draw[taa] (g) to (mu);
        \draw[alga] (mu) to (bc);
      \end{tikzpicture}
    }
    +
    \tikzsetnextfilename{hu-Alg-iso-is-3}
    \mathcenter{
      \begin{tikzpicture}
        \node at (-1,0) (tl) {};
        \node at (0,0) (tc) {};
        \node at (1,0) (tr) {};
        \node at (0,-1) (mu) {$\mu^*$};
        \node at (0,-2) (g) {$g^*$};
        \node at (0,-3) (bc) {};
        \draw[tbb] (tl) to (g);
        \draw[tbb] (tc) to (mu);
        \draw[tbb] (tr) to (g);
        \draw[blga] (mu) to (g);
        \draw[alga] (g) to (bc);
      \end{tikzpicture}
    }
    +
    \tikzsetnextfilename{hu-Alg-iso-is-4}
    \mathcenter{
      \begin{tikzpicture}
        \node at (-1,0) (tl) {};
        \node at (1,0) (tr) {};
        \node at (0,-1) (g) {$g^*$};
        \node at (0,-2) (bc) {};
        \node at (0,-.1) (one) {};
        \draw[tbb] (tl) to (g);
        \draw[tbb] (tr) to (g);
        \draw[alga] (g) to (bc);
        \draw[stumpa] (one) to (g);
      \end{tikzpicture}
    }=0.
  \]
  Since composition with $f$ is an isomorphism (with inverse given by
  composition by $g$, though $g$ is not known yet to be a
  homomorphism), it suffices to show that pre-composing the left side
  by $f$ vanishes. This is
  \begin{align*}
    \tikzsetnextfilename{hu-Alg-iso-is-5}
    \mathcenter{
      \begin{tikzpicture}
        \node at (0,1) (tc) {};
        \node at (0,0) (f) {$\underline{f}^*$};
        \node at (0,-1) (g) {$\underline{g}^*$};
        \node at (0,-2) (mu) {$\mu^*$};
        \node at (0,-3) (bc) {};
        \draw[taa] (tc) to (f);
        \draw[tbb] (f) to (g);
        \draw[taa] (g) to (mu);
        \draw[alga] (mu) to (bc);
      \end{tikzpicture}
    }
    +
    \tikzsetnextfilename{hu-Alg-iso-is-6}
    \mathcenter{
      \begin{tikzpicture}
        \node at (-1,0) (tl) {};
        \node at (0,0) (tc) {};
        \node at (1,0) (tr) {};
        \node at (-1,-1) (fl) {$\underline{f}^*$};
        \node at (0,-1) (fc) {$\underline{f}^*$};
        \node at (1,-1) (fr) {$\underline{f}^*$};
        \node at (0,-2) (mu) {$\mu^*$};
        \node at (0,-3) (g) {$g^*$};
        \node at (0,-4) (bc) {};
        \draw[taa] (tl) to (fl);
        \draw[taa] (tc) to (fc);
        \draw[taa] (tr) to (fr);
        \draw[tbb] (fl) to (g);        
        \draw[tbb] (fc) to (mu);        
        \draw[tbb] (fr) to (g);
        \draw[blga] (mu) to (g);
        \draw[alga] (g) to (bc);
      \end{tikzpicture}
    }
    +
    \tikzsetnextfilename{hu-Alg-iso-is-7}
    \mathcenter{
      \begin{tikzpicture}
        \node at (-1,0) (tl) {};
        \node at (1,0) (tr) {};
        \node at (-1,-1) (fl) {$\underline{f}^*$};
        \node at (1,-1) (fr) {$\underline{f}^*$};
        \node at (0,-1) (one) {};
        \node at (0,-2) (g) {$g^*$};
        \node at (0,-3) (bc) {};
        \draw[taa] (tl) to (fl);
        \draw[taa] (tr) to (fr);
        \draw[tbb] (fl) to (g);
        \draw[tbb] (fr) to (g);
        \draw[alga] (g) to (bc);
        \draw[stumpa] (one) to (g);
      \end{tikzpicture}
    }
    &=
    \tikzsetnextfilename{hu-Alg-iso-is-8}
    \mathcenter{
      \begin{tikzpicture}
        \node at (0,0) (tc) {};
        \node at (0,-1) (mu) {$\mu^*$};
        \node at (0,-2) (bc) {};
        \draw[taa] (tc) to (mu);
        \draw[alga] (mu) to (bc);
      \end{tikzpicture}
    }
    +
    \tikzsetnextfilename{hu-Alg-iso-is-9}
    \mathcenter{
      \begin{tikzpicture}
        \node at (-1,0) (tl) {};
        \node at (-.5,0) (tcl) {};
        \node at (0,0) (tc) {};
        \node at (.5,0) (tcr) {};
        \node at (1,0) (tr) {};
        \node at (-1,-1) (fl) {$\underline{f}^*$};
        \node at (0,-1) (mu) {$\mu$};
        \node at (1,-1) (fr) {$\underline{f}^*$};
        \node at (0,-2) (fc) {$f^*$};
        \node at (0,-3) (g) {$g^*$};
        \node at (0,-4) (bc) {};
        \draw[taa] (tl) to (fl);
        \draw[taa] (tc) to (mu);
        \draw[taa] (tcl) to (fc);
        \draw[taa] (tcr) to (fc);
        \draw[taa] (tr) to (fr);
        \draw[tbb] (fl) to (g);        
        \draw[tbb] (fc) to (g);        
        \draw[tbb] (fr) to (g);
        \draw[alga] (mu) to (fc);
        \draw[alga] (g) to (bc);
      \end{tikzpicture}
    }
    +
    \tikzsetnextfilename{hu-Alg-iso-is-10}
    \mathcenter{
      \begin{tikzpicture}
        \node at (-1,0) (tl) {};
        \node at (-.5,0) (tcl) {};
        \node at (0,0) (tc) {};
        \node at (.5,0) (tcr) {};
        \node at (1,0) (tr) {};
        \node at (-1,-1) (fl) {$\underline{f}^*$};
        \node at (0,-1) (one) {};
        \node at (1,-1) (fr) {$\underline{f}^*$};
        \node at (0,-2) (fc) {$f^*$};
        \node at (0,-3) (g) {$g^*$};
        \node at (0,-4) (bc) {};
        \draw[taa] (tl) to (fl);
        \draw[taa] (tcl) to (fc);
        \draw[taa] (tcr) to (fc);
        \draw[taa] (tr) to (fr);
        \draw[tbb] (fl) to (g);        
        \draw[tbb] (fc) to (g);        
        \draw[tbb] (fr) to (g);
        \draw[stumpa] (one) to (fc);
        \draw[alga] (g) to (bc);
      \end{tikzpicture}
      }\\
    =
    \tikzsetnextfilename{hu-Alg-iso-is-11}
    \mathcenter{
      \begin{tikzpicture}
        \node at (0,0) (tc) {};
        \node at (0,-1) (mu) {$\mu^*$};
        \node at (0,-2) (bc) {};
        \draw[taa] (tc) to (mu);
        \draw[alga] (mu) to (bc);
      \end{tikzpicture}
    }
    +
    \tikzsetnextfilename{hu-Alg-iso-is-12}
    \mathcenter{
      \begin{tikzpicture}
        \node at (0,0) (tc) {};
        \node at (0,-1) (mu) {$\mu^*$};
        \node at (0,-2) (bc) {};
        \draw[taa] (tc) to (mu);
        \draw[alga] (mu) to (bc);
      \end{tikzpicture}
    }
    +
    \tikzsetnextfilename{hu-Alg-iso-is-13}
    \mathcenter{
      \begin{tikzpicture}
        \node at (-1,0) (tl) {};
        \node at (0,-.25) (one) {};
        \node at (1,0) (tr) {};
        \node at (0,-1) (Id) {$\Id$};
        \node at (0,-2) (bc) {};
        \draw[taa] (tl) to (Id);
        \draw[taa] (tr) to (Id);
        \draw[alga] (Id) to (bc);
        \draw[stumpa] (one) to (Id);
      \end{tikzpicture}
    }=0.
  \end{align*}
  Here, the last term vanishes because $\Id$ vanishes if there is more
  than one input, and we checked the homotopy unital $\Ainf$ relation
  in the case that the input is a single thorn separately. This completes the proof.
\end{proof}

\begin{definition}
  A \emph{module transformation thorn tree} is the same as a module
  thorn tree except that there is a distinguished (\emph{purple}) internal vertex which
  is allowed to have one input and weight~$0$.
    
  If $S$ is a module thorn tree and $T$ is a module transformation
  thorn tree then $T\circ_1 S$ and $S\circ_1 T$ are defined in the
  obvious way. Similarly, if $S$ is an (algebra) thorn tree then
  $T\circ_i S$, $i>1$, is defined in the obvious way.
  
  The differential of a module transformation thorn tree is defined in
  the same way as for a module thorn tree, with the understanding that
  the purple vertex is allowed to split into a 2-valent, weight 0
  purple vertex and another (allowed) vertex.
  Also, for a $3$-valent, weight $0$ purple
  vertex with one input a thorn, there is \emph{not} a term in the
  differential corresponding to erasing that thorn.

  Let $\uwMTransCx{n}{w}$ be the complex of module transformation
  thorn trees with $n$ inputs and total weight $w$. We call
  $\uwMTransCx{n}{w}$ the \emph{homotopy unital weighted module
    transformation trees complex}.
\end{definition}

As in the non-unital case, there are obvious composition maps
\begin{align}
  \circ_1&\co \uwMTransCx{n}{w}\rotimes{\Ring} \uwMTreesCx{m}{v}\to \uwMTransCx{m+n-1}{v+w}\label{eq:uw-mod-trans-compo-1}\\
  \circ_1&\co \uwMTreesCx{n}{w}\rotimes{\Ring} \uwMTransCx{m}{v}\to \uwMTransCx{m+n-1}{v+w}\label{eq:uw-mod-trans-compo-2}\\
  \qquad\qquad\circ_i&\co \uwMTransCx{n}{w}\rotimes{\Ring} \uwTreesCx{m}{v}\to
           \uwMTransCx{m+n-1}{v+w} \label{eq:uw-mod-trans-compo-3}\qquad\qquad 1<i\leq n.
\end{align}

\begin{lemma}
  The $n$-input, weight $w$ homotopy unital weighted module
  transformation trees complex $\uwMTransCx{n}{w}$ is a chain complex,
  and the differential decreases the grading by one. Further, the
  composition
  maps~\eqref{eq:uw-mod-trans-compo-1},~\eqref{eq:uw-mod-trans-compo-2},
  and~\eqref{eq:uw-mod-trans-compo-3} are chain maps.
\end{lemma}
\begin{proof}
  The proof is essentially the same as the proof of Lemma~\ref{lem:uwTrees-is-cx}, and is left to the reader.  
\end{proof}

\begin{figure}
  \centering
  %Font is 12 point
  \includegraphics{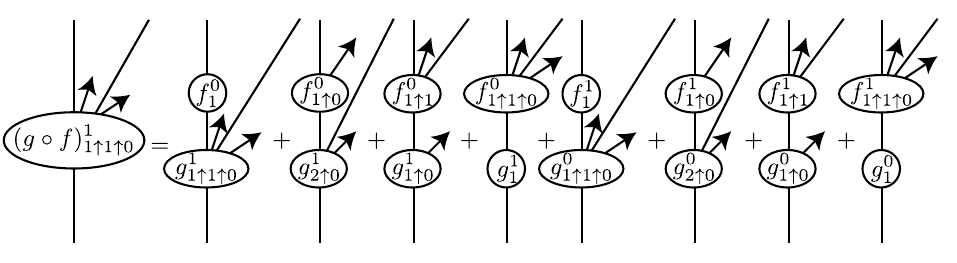}
  \caption[Composing homotopy unital weighted module
      morphisms]{\textbf{Composing homotopy unital weighted module
      morphisms.} The arrows are the thorns.}
  \label{fig:compose-hu-w-mod-maps}
\end{figure}

\begin{definition}
  Given homotopy unital modules $\wMod$ and $\wNod\!$, a homotopy unital
  weighted module map from $\wMod$ to $\wNod$ consists of maps
  \[
    \uwMTransCx{n}{w}\to \Mor(M\kotimes{\Ground} A^{\kotimes{\Ground} n-1},N\grs{(2-\kappa)w})
  \]
  which are compatible with composition in the obvious sense.
\end{definition}
  
A homotopy unital module map $f$ is determined by maps
\[
  f_{1+n_1\uparrow\cdots\uparrow n_k}^w\co M\otimes A^{\kotimes{\Ground} n}\to N
\]
corresponding to a purple corolla with weight $w$ and valence $n_1+\cdots+n_k+k-1$
with a mix of thorn and input leaves, where the thorns are at the
positions $n_1+2,n_1+n_2+3,\dots$. (The grading has been suppressed
here.)

\begin{definition}
  Given two homotopy unital weighted module maps
  $f\co \wMod_1\to \wMod_2$ and $g\co \wMod_2\to\wMod_3$, the
  \emph{composition} $g\circ f$ is defined by setting $(g\circ f)(S)$
  to be the
  sum of all ways of splitting the purple vertex in $S$ into two
  adjacent purple vertices and applying $f$ at the first of these
  purple vertices and $g$ at the second. Equivalently,
  \begin{multline*}
    (g\circ f)_{1+n_1\uparrow n_2\uparrow\cdots\uparrow n_k}(x,a_1,\dots,a_n)\\
    =\sum_{\ell=1}^k\sum_{n'_\ell=0}^{n_\ell}
    g_{1+n_\ell-n'_\ell\uparrow n_{\ell+1}\uparrow\cdots\uparrow
      n_k}(f_{1+n_1\uparrow\cdots\uparrow n_{\ell-1}\uparrow n'_\ell}(x,a_1,\dots,a_{n_1+\cdots+n_{\ell-1}+n'_\ell}),a_{n_1+\cdots+n'_\ell+1},\dots,a_n).
  \end{multline*}  
\end{definition}

\begin{example}
  The \emph{identity map} of $\wMod$ is the map $\Id_{\wMod}$ with
  \[
    (\Id_{\wMod})_{1+n_1\uparrow\cdots\uparrow n_k}^w=
    \begin{cases}
      \Id_M & k=1,\ w=0,\text{ and }n_1=0\\
      0 & \text{otherwise.}
    \end{cases}
  \]
\end{example}

\begin{lemma}
  The homotopy unital module maps from $\wMod$ to $\wNod$ form a
  subcomplex of the morphism complex
  \[
    \prod_{n,w} \Mor(\uwMTransCx{n}{w}\rotimes{\Ring} M\kotimes{\Ground} A^{\kotimes{\Ground} n-1},N\grs{(2-\kappa)w}).
  \]
  Further, the identity map is a cycle,
  $f\circ \Id_{\wMod}=\Id_{\wNod}\circ f=f$, and composition is a
  chain map.  Thus, homotopy unital modules and homotopy unital module
  maps forms a (strictly unital) \dg category.
\end{lemma}
\begin{proof}
  For the first statement, we need to check that if $\phi$ respects
  composition then $d(\phi)$ also respects composition. We spell out
  the case that $T=T_1\circ_i T_2$, with the purple vertex
  on $T_1$ and $i>1$; the other cases are similar. Let $k$ be the
  number of inputs of $T_2$. Then
  \begin{align*}
    (d\phi)&(T_1\circ_i T_2)(m,a_1,\dots,a_n)\\
    &=\phi((\bdy T_1)\circ_i T_2)(m,a_1,\dots,a_n)+\phi(T_1\circ_i(\bdy T_2))(m,a_1,\dots,a_n)
      +\bdy(\phi(T_1\circ_i T_2)(m,a_1,\dots,a_n))\\
    &\qquad+\phi(T_1\circ_i T_2)(\bdy m,a_1,\dots,a_n)
      +\sum_{j=1}^n \phi(T_1\circ_i T_2)(m,a_1,\dots,\bdy a_j,\dots,a_n)\\
           &=(\phi(\bdy T_1)\circ_i \mu(T_2))(m,a_1,\dots,a_n)+(\phi(T_1)\circ_i\mu(\bdy T_2))(m,a_1,\dots,a_n)\\
    &\qquad+\bdy(\phi(T_1)(m,a_1,\dots,a_{i-1},\mu(T_2)(a_i,\dots,a_{i+k-1}),a_{i+k},\dots,a_n))\\
      &\qquad+\phi(T_1)(\bdy m,a_1,\dots,a_{i-1},\mu(T_2)(a_i,\dots,a_{i+k-1}),a_{i+k},\dots,a_n)\\
      &\qquad+\sum_{j=1}^{i-1}\phi(T_1)(m,a_1,\dots,\bdy(a_j),\dots,a_{i-1},\mu(T_2)(a_i,\dots,a_{i+k-1}),a_{i+k},\dots,a_n)\\
      &\qquad+\sum_{j=i}^{i+k-1}\phi(T_1)(m,a_1,\dots,a_{i-1},\mu(T_2)(a_i,\dots,\bdy(a_j),\dots,a_{i+k-1}),a_{i+k},\dots,a_n)\\
      &\qquad +\sum_{j=i+k}^{n}\phi(T_1)(m,a_1,\dots,a_{i-1},\mu(T_2)(a_i,\dots,a_{i+k-1}),a_{i+k},\dots,\bdy(a_j),\dots,a_n)\\
    &=(\phi(\bdy T_1)\circ_i \mu(T_2))(m,a_1,\dots,a_n)+\bdy(\phi(T_1)(m,a_1,\dots,a_{i-1},\mu(T_2)(a_i,\dots,a_{i+k-1}),a_{i+k},\dots,a_n))\\
      &\qquad+\phi(T_1)(\bdy m,a_1,\dots,a_{i-1},\mu(T_2)(a_i,\dots,a_{i+k-1}),a_{i+k},\dots,a_n)\\
      &\qquad+\sum_{j=1}^{i-1}\phi(T_1)(m,a_1,\dots,\bdy(a_j),\dots,a_{i-1},\mu(T_2)(a_i,\dots,a_{i+k-1}),a_{i+k},\dots,a_n)\\
      &\qquad+\phi(T_1)(m,a_1,\dots,a_{i-1},\bdy(\mu(T_2)(a_i,\dots,,a_{i+k-1})),a_{i+k},\dots,a_n)\\
      &\qquad+\sum_{j=i+k}^{n}\phi(T_1)(m,a_1,\dots,a_{i-1},\mu(T_2)(a_i,\dots,a_{i+k-1}),a_{i+k},\dots,\bdy(a_j),\dots,a_n)\\
    &=(d\phi)(T_1)(m,a_1,\dots,a_{i-1},\mu(T_2)(a_{i},\dots,a_{i+k-1}),a_{i+k},\dots,a_n),
  \end{align*}
  as desired.

  Turning to the remaining statements, it is clear that $\Id_{\wMod}$
  is a cycle and that $f\circ \Id_{\wMod}=\Id_{\wNod}\circ f=f$. To
  verify that composition is a chain map, it suffices to check that
  $d(g\circ f)$ and $(dg\circ f)+(g\circ df)$ take the same value on 
  each corolla (with thorns).
  This is straightforward from the definitions, and is left to the
  reader.
\end{proof}

\begin{definition}\label{def:hu-bonsai}
  A homotopy unital weighted algebra $\wAlg$ is \emph{bonsai} if
  there is an integer $N$ so that $\mu(T)=0$ whenever $T$ is a thorn
  tree with $\dim(T)>N$. Bonsai weighted algebra homomorphisms,
  weighted modules, and weighted module morphisms are defined
  similarly.
\end{definition}

Finally, we have a rectification result:
\begin{theorem}\label{thm:hu-rectify}
  Given a split homotopy unital weighted algebra $\wAlg$ there is a strictly
  unital weighted algebra $\wBlg$ and a homotopy unital
  isomorphism $f\co\wAlg\to\wBlg$. Further, if $\wAlg$ is bonsai
  then $\wBlg$ and $f$ can be chosen to be bonsai as well.
\end{theorem}

The proof uses a version of Lemma~\ref{lem:ModifyActionw}:
\begin{lemma}
  \label{lem:ModifyActionHU}
  Suppose that $\wAlg$ is a homotopy unital $w$-algebra. Fix integers
  $v,j,m_1,\dots,m_j\geq 0$, with $m_1+\cdots+m_j+j-1+2v>1$. Let
  $m=m_1+\cdots+m_j$, and let
  $\phi_{m_1\uparrow\cdots\uparrow m_j}^v\co A^{\kotimes{\Ground} m}\to A$ be any map of
  degree $(2-\kappa)v+m+j-2$. 
  There is a homotopy unital $w$-algebra
  $\wAlg'=(A,\overline{\mu}_{n_1\uparrow \cdots\uparrow n_k}^w)$ isomorphic to
  $\wAlg$, so that
  $\overline{\mu}^w_{n_1\uparrow\cdots\uparrow
    n_k}=\mu^w_{n_1\uparrow\cdots\uparrow n_k}$ if
  \begin{itemize}
  \item $w<v$, or
  \item $w=v$ and $n+k<m+j$.
  \end{itemize}
  Further, 
  the only cases with $w=v$ and $n+k=m+j$ for which
  $\overline{\mu}^w_{n_1\uparrow\cdots\uparrow n_k}\neq
  \mu^w_{n_1\uparrow\cdots\uparrow n_k}$ are
  \begin{multline*}
    \begin{aligned}
      \overline{\mu}^v_{m_1\uparrow\cdots\uparrow m_j} &= \mu^v_{m_1\uparrow\cdots\uparrow m_j} + d \phi^v_{m_1\uparrow\cdots\uparrow m_j}\\
      \overline{\mu}^v_{m_1\uparrow\cdots\uparrow m_{i-1}\uparrow m_i'\uparrow m_i''\uparrow m_{i+1}\uparrow\cdots\uparrow m_j}(a_1,\dots,a_{m-1})
      &=\mu^v_{m_1\uparrow\cdots\uparrow m_{i-1}\uparrow m_i'\uparrow
        m_i''\uparrow m_{i+1}\uparrow\cdots\uparrow
        m_j}(a_1,\dots,a_{m-1})
    \end{aligned}\\
    +
    \phi^v_{m_1\uparrow\cdots\uparrow m_j}(a_1,\dots,a_{m_1+\cdots+m_{i-1}+m'_i},\unit,a_{m_1+\cdots+m_{i-1}+m'_i+1},\dots,a_{m-1}),
  \end{multline*}
  where $m_i'+m_i''=m_i-1$.
\end{lemma}
\begin{proof}[Proof sketch]
  As in Lemma~\ref{lem:ModifyAction}, we construct the modified algebra
  $\wAlg'$, with operations coming from trees as follows:
  \begin{itemize}
  \item There is one distinguished vertex labelled by either a stump
    or some $\mu_{n_1\uparrow\cdots\uparrow n_k}^w$.
  \item All other vertices have $m$ inputs and weight~$v$ and are
    labelled by~$\phi$.
  \item The distinguished vertex is either the trunk vertex or a
    parent of the trunk vertex. If the distinguished vertex is a
    stump, it must be a parent of the trunk vertex.
  \end{itemize}
  The tree contributes to the operation on $\wAlg'$ with the given
  inputs, with arrows interpolated at the positions of the
  thorns. Additionally, if the distinguished vertex is a stump, it
  contributes an extra arrow at its position.

  The proof of Lemma~\ref{lem:ModifyAction} applies almost without
  change to show that these operations form a homotopy unital
  $w$-algebra, except that some decompositions are disallowed because
  the stump cannot be the root of the modified action trees. These
  disallowed terms correspond instead to terms in the homotopy unital
  $w$-algebra relations that turn a thorn into a stump.

  We can similarly see that $\wAlg'$ is isomorphic to~$\wAlg$ by
  constructing a map that is a perturbation of the identity by~$\phi$.
  
  The terms in the last part of the lemma statement come from trees
  with one undistinguished vertex and either a distinguished
  $2$-valent vertex or a distinguished stump.
\end{proof}

\begin{remark}
  One can also prove Lemma~\ref{lem:ModifyActionHU} by
  following the strategy of Remark~\ref{rem:ModifyActionSimple}.
\end{remark}

\begin{proof}[Proof of Theorem~\ref{thm:hu-rectify}]
  We prove the unweighted case; the weighted case is obtained by
  wrapping the argument below in an induction on the weight $w$.

  We will modify the algebra so that
  $\mu_{n_1\uparrow\cdots\uparrow n_k}=0$ for $k>1$.  The construction is
  inductive on:
  \begin{enumerate}
  \item $n_1+\cdots+n_k+k$ (outer induction),
  \item $n_1$ (middle induction), and
  \item $k$ (inner induction).
  \end{enumerate}
  We will use the following lemma, which gives a precise relationship
  between vanishing of operations with thorns and strict unitality.
  \begin{lemma}\label{lem:mu-pu}
    Suppose, in a homotopy unital algebra, that $\mu_{n_1\uparrow\cdots\uparrow n_k}=0$ whenever
    $n_1+\cdots+n_k+k\leq T$ and $n_1<N$. Then, if $m = \sum m_i > 2$,
    $m+\ell\leq T$, and $i<N$, we have
    \[
      \mu_{m_1\uparrow\cdots\uparrow m_\ell}(a_1,\dots,a_{i-1},\unit,a_{i+1},\dots,a_m)=0.
    \]
  \end{lemma}
  \begin{proof}
    If $m_1<N$ then $\mu_{m_1\uparrow\cdots\uparrow m_\ell}$ vanishes
    by assumption. Otherwise, consider the $\Ainf$ relation
    corresponding to
    $\corolla{i-1\uparrow m_1-i\uparrow m_2\uparrow\cdots\uparrow
      m_\ell}$
    with inputs $(a_1,\dots,a_{i-1},a_{i+1},\dots,a_m)$. One term is
    \[
      \mu_{m_1\uparrow\cdots\uparrow m_\ell}(a_1,\dots,a_{i-1},\unit,a_{i+1},\dots,a_m)
    \]
    and the other terms vanish by hypothesis.
  \end{proof}

  Now, fix an $\Ring$-module splitting
  $A=\Field\langle \unit\rangle \oplus A_0$. Assume that
  $\mu_{m_1\uparrow\cdots\uparrow m_\ell}=0$ if $\ell>1$ and either
  \begin{itemize}
  \item $m_1+\cdots+m_\ell+\ell<T$; or
  \item $m_1+\cdots+m_\ell+\ell=T$ and $m_1<N$; or
  \item $m_1+\cdots+m_\ell+\ell=T$ and $m_1=N$ and $\ell<K$.
  \end{itemize}
  Suppose that the first non-vanishing homotopy-unital operation is
  $\mu_{N\uparrow n_2\uparrow\cdots\uparrow n_t}\neq 0$.
  Define
  \[
    \phi_{N+1+n_2\uparrow n_3\uparrow\cdots\uparrow n_t}(a_1,\dots,a_{n+1})=
    \begin{cases}
     \mu_{N\uparrow\cdots\uparrow n_t}(a_1,\dots,a_{N},a_{N+2},\dots,a_{n+1}) & a_{N+1}=\unit\\
     0 & a_{N+1}\in A_0.
    \end{cases}
  \]
  Let $\overline{\mu}$ be the operations on the induced homotopy unital
  $\Ainf$-algebra from Lemma~\ref{lem:ModifyActionHU}. By that lemma,
  $\overline{\mu}_{m_1\uparrow\cdots\uparrow m_j}=\mu_{m_1\uparrow\cdots\uparrow m_j}$
  if $m_1+\cdots+m_j+j<N+n_2+\cdots+n_t+t$, and
  \[
    \overline{\mu}_{N\uparrow\cdots\uparrow n_t}(a_1,\dots,a_{n})=
    \mu_{N\uparrow\cdots\uparrow n_t}(a_1,\dots,a_{n})+
    \phi_{N+1+n_2\uparrow\cdots\uparrow n_t}(a_1,\dots,a_{N},\unit,a_{N+1},\dots,a_{n})=0.
  \]

  Although not needed for the induction, note that
  the operation
  $\mu_{N+1+n_2\uparrow n_3\uparrow\cdots\uparrow n_t}$ with $(N+1)\st$
  input the unit changes to
  \begin{multline*}
    \overline{\mu}_{N+1+n_2\uparrow n_3\uparrow\cdots\uparrow n_t}(a_1,\dots,a_{N},\unit,a_{N+2},\dots,a_{n+1})=
    {\mu}_{N+1+n_2\uparrow n_3\uparrow\cdots\uparrow n_t}(a_1,\dots,a_{N},\unit,a_{N+2},\dots,a_{n+1})\\
    +(d\mu_{N\uparrow n_2\uparrow\cdots\uparrow n_t})(a_1,\dots,a_{N},a_{N+2},\dots,a_{n+1}). 
  \end{multline*}
  Applying the $\Ainf$ relation to the $d\mu$ term on the right
  side, the inductive
  hypothesis implies that almost all terms vanish, except for one term that
  cancels the first term on the right side.
  So,
  if $t = 2$, the operation $\mu_{n+1}$ has
  just become strictly unital in the $(N+1)\st$ place.

  According to
  Lemma~\ref{lem:ModifyActionHU},
  the remaining cases with $m_1+\cdots+m_\ell+\ell=T$ for which we may
  have $\overline{\mu} \ne \mu$ are operations
  obtained from $N+1+n_2\uparrow n_3\uparrow\cdots\uparrow n_t$ by
  replacing $N$ or some $n_i$ by $n'\uparrow n''$.  The cases
  corresponding to some $n_i$ with $i>1$ are later in the induction,
  hence irrelevant. For the case of replacing $N$, by
  Lemma~\ref{lem:mu-pu}, this operation is, in fact, unchanged.

  Thus, we have eliminated the operation
  $\mu_{N\uparrow n_2\uparrow\cdots\uparrow n_t}$ without changing any operation
  earlier or in the same step in the induction. Repeat this for each
  operation at this stage of the induction and continue with the
  induction.

  The bonsai statement follows from the fact that the new
  multiplications and maps are defined via trees of the old
  multiplications and maps (possibly with $2$-valent, weight $0$
  vertices) with the same dimension.
\end{proof}

%%% Local Variables: 
%%% mode: latex
%%% TeX-master: "AbstractDiagonal.tex"
%%% TeX-command-extra-options: "--shell-escape"
%%% End: 

\section{The weighted trees complexes are contractible}\label{sec:wAcyclic}
In this section, we compute the homology of the various complexes of weighted trees introduced in Section~\ref{sec:wAinfty}.

\subsection{Signs on the weighted trees complex}\label{sec:w-signs}
In Section~\ref{sec:Associaplex} we construct a CW complex whose
cellular chain complex agrees with the weighted trees complex. Even
though the rest of this paper is in characteristic 2, it seems natural
to prove this identification of chain complexes over the integers. So,
in this section we define a signed refinement of the weighted trees
complex; we prove this signed version is contractible in
Section~\ref{sec:w-alg-cx-acyclic}. Of course, this also implies
contractibility in characteristic 2. In
Section~\ref{sec:w-mod-cx-acyclic} we revert to working in
characteristic 2. Our construction of signs is a trivial adaptation of
Markl-Schneider~\cite[Section 3]{MS06:AssociahedraProdAinf}.

Fix a commutative ring $\Ring$.

An \emph{orientation} for a weighted tree $T$ is an equivalence class
of total orderings of the internal edges of $T$ (i.e., the edges not
adjacent to the inputs or output), where total orderings $\omega$ and
$\omega'$ are equivalent if they differ by an even permutation. We use
Markl-Schneider's suggestive notation
$e_{i_1}\wedge \cdots\wedge e_{i_k}$ for the ordering
$e_{i_1}<\cdots<e_{i_k}$. Let $\wTreesCx{n}{w}$ be the free
$\Ring$-module generated by pairs $(T,\omega)$ of stably-weighted
trees and orientations for them, modulo the relation
$(T,\omega)\sim -(T,\omega')$ if $\omega$ and $\omega'$ differ by an
odd permutation. The boundary map is defined by
\begin{equation}\label{eq:signed-diff}
  \bdy(T,\omega)=\sum (T',e\wedge \omega)
\end{equation}
where $T$ is obtained from $T'$ by collapsing the edge $e$
(cf.~\cite[Equation (3.2)]{MS06:AssociahedraProdAinf}). It is
immediate from the definition and the relation
$e\wedge f\wedge \omega=-f\wedge e\wedge \omega$ that $\bdy^2=0$. If
$\Ring$ has characteristic $2$, this complex is exactly the weighted
trees complex from Section~\ref{sec:wAlgs}.

Composition with signs can also be defined following
Markl-Schneider~\cite[Equation (3.1)]{MS06:AssociahedraProdAinf}: for $(S,\eta)\in \wTreesCx{m}{v}$ and $(T,\omega)\in \wTreesCx{n}{w}$ let
\[
  (T,\omega)\circ_i(S,\eta)=(-1)^{n\dim(S)+mi}(T\circ_i S,\omega\wedge \eta\wedge e)
\]
where $e$ is the new internal edge of $T\circ_i S$.  It is
straightforward to verify that composition is a chain map and
satisfies the usual operad associativity law.

\subsection{Contractibility of the weighted trees complex}\label{sec:w-alg-cx-acyclic}
\begin{theorem}
  \label{thm:AssociaplexAcyclic}
  The homology of the $n$-input, weight $w$ weighted trees complex $\wTreesCx{n}{w}$ is
  isomorphic to $\Ring$, supported in dimension zero. A generator for the homology is
  represented by the right-associated planar weighted tree, such that all of
  its inputs are to the left of all its popsicles,
  all of its valence $>1$ vertices have weight $0$, and all of its
  popsicles have weight $1$.
\end{theorem}

Define a map $K\co \wTreesCx{n}{w}\to \wTreesCx{n}{w}$ as follows.
Call a vertex $v$ of a weighted tree $T$ {\em potentially movable} if
\begin{itemize}
   \item its weight $w(v)>1$ or
   \item its weight $w(v)=1$ and its valence $d(v)>2$ or
   \item its weight $w(v)=1$, its valence $d(v)=2$, and the parent of $v$ is not an input.
\end{itemize}
Find the first potentially movable vertex $v$ in $T$ using a
depth-first search.  If $w(v)=1$, $d(v)=2$, and the parent of $v$ is
not an input, let $e$ be the edge pointing into $v$, let $T'=T/e$ be
the result of collapsing $e$, and define
$K(T,e\wedge \omega)=-(T',\omega)$. Otherwise (if $w(v)>1$ or
$d(v)>2$), let $K(T)=0$.

Let $Y$ be the subcomplex of $\wTreesCx{n}{w}$ with the property that all
vertices with positive weight have $w(v)=1$ and either $d(v)=1$ or $d(v)=2$ and the parent of $v$ is an input.

\begin{lemma}\label{lem:associa-htpy}
  For any tree $T$ and any sufficiently large $m$,
  \[
    (\Id + \partial \circ K + K \circ \partial)^m(T)\in Y.
  \]
  Further, if $T$ is in $Y$ then
  \[
    (\Id + \partial \circ K + K \circ \partial)(T)=T.
  \]
\end{lemma}

\begin{proof}
  Recall depth first search, as defined in Equation~\eqref{eq:DFS}; and note that
  $\depth$ satisfies the following monotonicity property under
  edge insertions.  Suppose $S$ is obtained from $T$ by inserting an edge
  at $v$, and let $u$ be any vertex in $T$ other than $v$. Then the
  vertex $u'$ in $S$ corresponding to $u$ has $\depth(u')\geq \depth(u)$.

  Let $\dpm(T)$ denote the depth $\depth(v)$ of
  first potentially movable vertex $v$ in $T$ with respect to the
  depth-first search ordering. We set $\dpm(T)=\infty$ if there are no
  such vertices, i.e., if the tree is a generator of $Y$. We also
  abuse notation and write $\dpm(T,\omega)=\dpm(T)$.
  By the monotonicity property of $\depth$, it
  follows that $\dpm$ induces a filtration on the chain
  complex. Moreover, it is clear from the construction of $K$ that if
  $K(T,\omega)\neq 0$, then $\dpm(K(T,\omega))=\dpm(T)$.

  It follows that for any tree $S$ that appears with non-zero
  multiplicity in $(\Id + \partial \circ K + K \circ \partial)(T,\omega)$,
  $\dpm(T)\leq \dpm(S)$. We claim in fact that $\dpm(T)<\dpm(S)$,
  unless $\dpm(T)=\infty$. This follows from a case analysis:

  \begin{figure}
    \centering
    \includegraphics{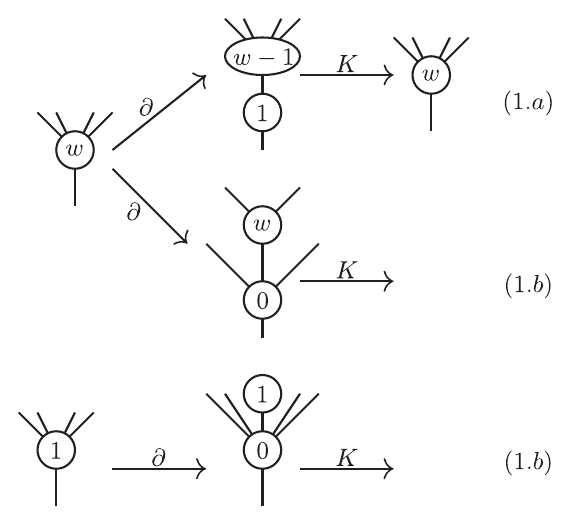}
    \caption[A case in the proof that the weighted trees complex is contractible]{\textbf{Case~(\ref{item:assoHtpy1}) of the proof of Lemma~\ref{lem:associa-htpy}.}}
    \label{fig:AssociaplexHtpy1}
  \end{figure}

  \begin{figure}
    \centering
    \includegraphics{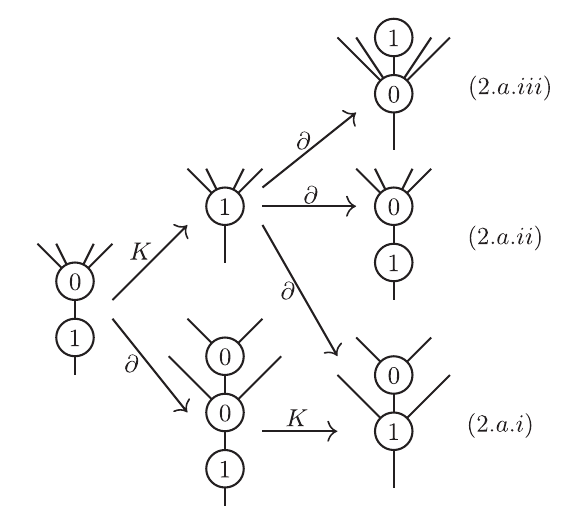}
    \caption[Another case in the proof that the weighted trees complex is contractible]{\textbf{Case~(\ref{item:assoHtpy2}) of the proof of Lemma~\ref{lem:associa-htpy}.}}
    \label{fig:AssociaplexHtpy2}
  \end{figure}
    
  \begin{enumerate}
  \item\label{item:assoHtpy1} Suppose that the first potentially
    movable vertex $v$ in $T$ has $w(v)>1$ or that $w(v)=1$ and
    $d(v)>2$; so that (in either case) $K(T,\omega)=0$.  The
    differential $\partial (T,\omega)$ is a sum over trees $S$
    obtained by inserting edges at the various vertices of $T$. If
    $K(S,\omega')\neq 0$, then $S$ must be obtained by inserting an
    edge $e$ in $T$ at $v$, and $\omega'=e\wedge \omega$. There are
    the following two subcases:
    \begin{enumerate}
    \item The vertex $v$ is replaced by two vertices $v_1$ and $v_2$,
      so that $d(v_1)=2$, $w(v_1)=1$, and $w(v_2)=w(v)-1$; and $v_2$
      is the parent of $v_1$ (so that $v_1$ is the first potentially
      movable vertex in $S$).  In this case,
      $K(S,e\wedge \omega)=-(T,\omega)$ (which cancels with
      $\Id(T,\omega)=(T,\omega)$).
    \item The vertex $v$ is replaced by two vertices $v_1$ and $v_2$
      so that $w(v_1)=0$ and $w(v_2)=w(v)$, and $v_2$ is a parent of
      $v_1$.  Clearly, $v_1$ is not potentially movable, so it follows
      that $\dpm(S)>\dpm(T)$.
    \end{enumerate}
  \item\label{item:assoHtpy2} Suppose that the first potentially
    movable vertex $v$ in $T$ has $w(v)=1$ and $d(v)=2$, so that the
    parent is not one of the inputs, and $K(T,\omega)\neq 0$. In this
    case, $K(T,\omega)$ is obtained by contracting the parent edge $e$
    of $v$, to form a new vertex $v'$.  Consider the terms in
    $\partial(K(T,\omega))$, which come in the following types:
    \begin{enumerate}
    \item Trees obtained by inserting an edge at $v'$, connecting
      vertices $v_1$ and $v_2$ (so that $v_2$ is the parent of $v_1$).
      Such trees $S$ can be of three basic types:
      \begin{enumerate}
      \item\label{item:assoHtpy2ai} Trees with $w(v_1)>1$ or $w(v_1)=1$ and $d(v_1)>2$; in
        such cases, there are corresponding canceling trees contained
        in $K(\partial (T,\omega))$, obtained by inserting edges $e'$ in the parent
        of $v$. Given an orientation $e\wedge \omega$ of $T$, this
        term in $\bdy K(T,e\wedge \omega)$ has orientation $-e'\wedge
        \omega$, while the corresponding term in
        $K(\bdy(T,e\wedge\omega))$ is $K(T',e'\wedge e\wedge
        \omega)=(T'',e'\wedge\omega)$, so the terms indeed appear with
        opposite signs.
      \item Trees with $w(v_1)=1$ and $d(v_1)=2$, in which case
        $S=T$. Keeping track of signs, given an orientation
        $e\wedge\omega$ for $T$, this term in $\bdy(K(T,e\wedge
        \omega))=-\bdy(T',\omega)$ is $-(T,e\wedge\omega)$, so this
        cancels with $\Id(T,e\wedge\omega)$.
      \item Trees with $w(v_1)=0$.  Since $v_1$ is not potentially
        movable, it follows that $\dpm(S)>\dpm(T)$.
      \end{enumerate}
    \item Trees obtained by inserting an edge at some vertex in $K(T)$
      other than $v'$.  These trees cancel with trees appearing in
      $K(\partial T)$. (As in case~(\ref{item:assoHtpy2ai}),
      skew-commutativity of the orientation implies that the signs
      work out in this cancellation.)
    \end{enumerate}
    There are potentially two remaining trees in $K(\partial T)$ not canceling with trees
    in $\partial K(T)$: the trees corresponding to
    inserting an edge at $v$.  Both of those trees $S$ have
    $\dpm(S)>\dpm(T)$, for the following reason. 
    The vertex $v$ is the first potentially movable vertex in $T$.
    Inserting an edge at $T$ gives a new tree $T'$ (with $K(T')=S$),
    whose first potentially movable vertex $u'$ corresponds to the
    second potentially movable vertex of $T$. (If there is no such
    second potentially movable vertex, then $K(T')=0$.) Thus,
    \[ \dpm(T)<\depth(u)\leq \depth(u')=\dpm(T')=\dpm(K(T')),\]
    as claimed.
    \end{enumerate}
    Since there are finitely many trees of weight $w$ and $n$ inputs,
    $\dpm$ takes values in a finite set.  Thus, for $m$
    sufficiently large,
    $(\Id + \partial \circ K + K \circ \partial)^m(T)=\infty$, i.e.,
    $(\Id + \partial \circ K + K \circ \partial)^m(T)\in Y$.

    The last part of the statement follows from the fact that $K$
    vanishes on the subcomplex $Y$.
\end{proof}

Since there are finitely many trees in $\wTreesCx{n}{w}$, we may
choose a single $m$ so that for any $T$,
$(\Id + \partial \circ K + K \circ \partial)^m(T)\in Y$. Then,
$(\Id + \partial \circ K + K \circ \partial)^m$ is a homotopy
equivalence from $\wTreesCx{n}{w}$ to $Y$.

\begin{proof}[Proof of Theorem~\ref{thm:AssociaplexAcyclic}.]
  There is a filtration of $Y$ by 
  \[
    \Filt(T)=\#\{v\in\Vertices(T)\mid v=\wcorolla{0}{1}\},
  \]
  i.e., by the number of (weight $1$) popsicles in $T$. Consider the spectral
  sequence associated to the filtration $\Filt$. The $E^0$-page is
  identified with many copies of the associahedron. More precisely,
  consider the set of sequences of symbols $i,p,s$ (for ``input'',
  ``popsicle'', and ``stick''), so that the number of occurrences of
  $p$ plus the number of occurrences of $s$ is $w$ and the number of
  occurrences of $i$ plus the number of occurrences of $s$ is
  $n$. Call such a sequence a \emph{leaf sequence}. The $E^0$-page
  consists of a copy of the unweighted trees complex
  $\Trees_{\#i+\#p+\#s}$
  for each leaf sequence (with $\#i$ $i$'s, etc.). (See Figure~\ref{fig:leaf-seq}.)

  \begin{figure}
    \centering
    %Font is 12 point.
    \includegraphics{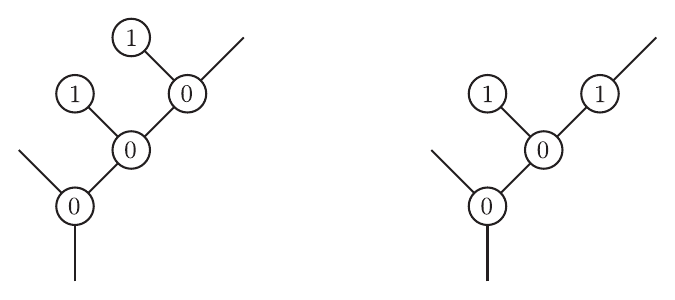}
    \caption[Examples of leaf sequences, for contractibility of the weighted trees complex]{\textbf{Leaf sequences.} Left: the sequence
      $(i,p,p,i)$. Right: the sequence $(i,p,s)$.}
    \label{fig:leaf-seq}
  \end{figure}

  Thus, the $E^1$-page has a single copy of $\Ring$ for each leaf
  sequence. Of course, the isomorphism with $\Ring$ depends on a
  choice of orientation. Given a leaf sequence $L$ and an orientation
  $\omega$ for $L$, the $d^1$ differential of a leaf sequence is the sum of
  all ways of replacing a copy of $s$ with $(p,i)+(i,p)$ and replacing
  $\omega$ with $e\wedge \omega$, where $e$ is the new edge.

  Define a map $H\co E^1 \to E^1$, the \emph{popsicle-stick homotopy},
  as follows. Given a leaf sequence $L$, consider the left-most
  instance of $p$ or $s$ in the sequence. If it is $s$ then
  $H(L,\omega)=0$. If $p$ is the first term in the sequence then let $L'$ be
  the rest of the sequence and define $H(L,\omega)=((p,H(L')),\omega)$. Finally, if
  the first $p$ or $s$ is $p$ and the previous term is $i$, and $e$ is
  the edge into this popsicle $p$, let $L'$ be the result of replacing
  this pair $(p,i)$ by $s$ and define
  $H(L,e\wedge \omega)=-(L',\omega)$. Then for $m$ sufficiently large, 
  \[
    (\Id + \partial \circ H + H \circ \partial)^m
  \]
  is a chain homotopy equivalence to the copy of $\Ring$ spanned by
  the leaf sequence $(p,p,\dots,p,i,i,\dots,i)$.
\end{proof}

\subsection{The weighted module trees complex is acyclic}\label{sec:w-mod-cx-acyclic}
In this section and subsequent ones, we again fix a ring $\Ring$ of
characteristic 2.  For the weighted module trees complex we have the
following analogue of Theorem~\ref{thm:AssociaplexAcyclic}:

\begin{theorem}
  \label{thm:ModAssociaplexAcyclic}
  The homology of the $n$-input, weight $w$ weighted module trees
  complex $\wMTreesCx{n}{w}$ is:
  \begin{itemize}
  \item Isomorphic to $R$ if $w=0$.
  \item Trivial if $w>0$.
  \end{itemize}
\end{theorem}

\begin{proof}
The weight-zero part of $\wMTreesCx{n}{0}$ is isomorphic to the
cellular chain complex of the associahedron, $\cellC{*}(K_n)$, so the
weight-zero case of Theorem~\ref{thm:ModAssociaplexAcyclic} follows
from the fact that the associahedron is a polyhedron. So, we will
focus on the positive-weight case.

Recall that $L^{n,w}_*\subset\wTreesCx{n}{w}$ denotes the subcomplex
of left-unmarked trees. We will show that the inclusion map
$L^{n,w}_*\into\wTreesCx{n}{w}$ induces an isomorphism on homology.
Consider the three steps in the proof that $\wTreesCx{n}{w}$ is
contractible:
\begin{enumerate}
\item First, we applied a homotopy $K$ to retract to a subcomplex
  $Y\subset \wTreesCx{n}{w}$. Note that $K$ preserves $L^{n,w}_*$, and
  hence gives a contraction of $L^{n,w}_*$ to $L^{n,w}_*\cap Y$.
\item We then filtered the complex $Y$ and observed that the
  $E^0$-page of the associated spectral sequence is a direct sum of
  copies of unweighted trees complexes. Each of these complexes is
  contractible. Further, each of these summands either lies in
  $L^{n,w}_*$ or intersects $L^{n,w}_*$ trivially, so these
  contractions again preserve $L^{n,w}_*$.
\item Finally, we define a retraction of the $E^1$-page to a
  particular tree $T_0$ in $L^{n,w}_*$. The homotopy $H$ used to
  define this retraction again preserves $L^{n,w}_*$.
\end{enumerate}
Thus, as claimed, the inclusion $L^{n,w}_*\into \wTreesCx{n}{w}$ induces an
isomorphism in homology, so
$\wMTreesCx{n}{w}=\wTreesCx{n}{w}/L^{n,w}_*$ is acyclic.
\end{proof}

\subsection{The weighted transformation trees complex is contractible}

In this section we show that the weighted transformation trees complex
is contractible:
\begin{proposition}
  \label{prop:wTransCxContract}
  The weighted transformation trees complex $\wTransCx{n}{w}$ has homology
  $H_0(\wTransCx{n}{w})=\Ring$ and $H_i(\wTransCx{n}{w})=0$ for all $i>0$.
\end{proposition}

\begin{proof}
  Define a homotopy operator $H$ as follows. 
  Given a weighted transformation tree $T$, find via depth first search
  the first purple vertex which
  is not a $2$-valent weight $0$ vertex $v$ whose predecessor is an input.
  If $v$ is $2$-valent with weight $0$, let $H(T)$ be the weighted
  transformation tree obtained by contracting the edge into $v$;
  otherwise, let $H(T)=0$.
  
  Consider the function
  \[
    F(T)=\sum_{\{u\in\Vertices(T)\mid u~\text{is not blue}\}} (\valence(u)+2w(u)-2).
  \]  
  Clearly,
  $F(T)\geq 0$ for all $T$, and the differential respects the
  filtration determined by $F$. Finally, if $F(T)>0$, then 
  all terms in 
  \[ (\partial\circ H + H \circ \partial + \Id)(T)\] 
  are in filtration level strictly less than $F(T)$;
  while if $F(T)=0$, then $(\partial\circ H + H \circ \partial)(T) =0$.
  It follows that $\wTransCx{n}{w}$ is chain homotopy equivalent to the subcomplex
  generated by $T$ where $F(T)=0$. If $F(T)=0$ then $T$ consists of a valence $2$, 
  weight $0$ purple vertex and a tree of blue vertices. Thus, this subcomplex
  is identified with the
  associaplex with $n$ inputs and weight $w$. The result now follows from
  the corresponding fact for the associaplex.
\end{proof}

\subsection{The weighted module transformation trees complex is acyclic}
Finally we show that the weighted module transformation trees complex 
is acyclic:
\begin{proposition}\label{prop:wModTransAcyclic}
  The weighted module transformation trees complex $\wMTransCx{n}{w}$
  has homology $H_0(\wMTransCx{n}{0})=\Ring$,
  $H_i(\wMTransCx{n}{0})=0$ for $i>0$, and $H_*(\wMTransCx{n}{w})=0$
  if $w>0$.
\end{proposition}
\begin{proof}
  Define a homotopy $K$ of the module transformation trees complex by
  declaring that $K(T)=0$ if the distinguished vertex $v$ of $T$ has weight
  $w(v)>0$ or valence $\valence(v)>2$ or if $w(v)=0$ and $\valence(v)=2$ but the predecessor to $v$ is an
  input of $T$; and if $w(v)=0$ and $\valence(v)=2$ and the
  predecessor to $v$ is not an input of $T$ then $K(T)$ is the result
  of contracting the edge into $v$. Then for $N$ sufficiently large,
  $(K\circ\bdy+\bdy\circ K+\Id)^N$ sends any tree $T$ to a
  transformation tree $T'$ where $v$ is $2$-valent, weight $0$, and at
  the top of $T'$. The resulting complex is identified with the
  weighted module trees complex. The proposition follows.
\end{proof}

\subsection{The homotopy unital complexes are contractible}
\begin{definition}\label{def:proj-hu-to-u}
For $n+2w\geq 2$ there is a chain map
$\pi\co \uwTreesCx{n}{w}\to \wTreesCx{n}{w}$ which erases all stumps
adjacent to $2$-input, weight $0$ vertices and forgets the vertex
adjacent to the stump, and then sends the tree to zero if there are
any thorns or any remaining stumps.  There are similar projections for
the module trees complexes, the transformation trees complexes, and
the module transformation trees complexes.
\end{definition}

\begin{theorem}\label{thm:hu-contractible}
  For $n+2w\geq 2$ the projection
  $\pi\co \uwTreesCx{n}{w}\to \wTreesCx{n}{w}$ induces an isomorphism on
  homology. The homology of $\uwTreesCx{0}{0}$ is generated by $\stump$, while the homology of $\uwTreesCx{1}{0}$ is generated by $\IdTree$.  
  In particular, for each $n,w$, the homology of $\uwTreesCx{n}{w}$ is
  one-dimensional, supported in dimension $0$.
\end{theorem}

(In the unweighted case, this is immediate from~\cite[Corollary 4.5]{MuroTonks14:unital-assoc}.)

\begin{proof}
  There is a filtration on
  $\uwTreesCx{n}{w}$ by the number of thorns. The associated graded
  complex is a direct sum of copies of $\wTreesCx{n+s+t}{w}$,
  corresponding to trees with $s$ stumps and $t$ thorns. There is one
  copy of $\wTreesCx{n+s+t}{w}$ for each pair of subsets
  $S,T\subset \{1,\dots, n+s+t\}$ with $|S|=s$, $|T|=t$, and
  $S\cap T=\emptyset$ (where $S$ and $T$ correspond to the locations
  of the stumps and thorns, respectively). Thus, by
  Theorem~\ref{thm:AssociaplexAcyclic}, the $E^1$-page is
  $1$-dimensional for each pair $S,T$ as above, represented by the
  left-associated, binary trees with all weight in weight-1 popsicles,
  where some inputs are thorns or stumps, and the popsicles are to the
  right of all other inputs.

  The $d^1$-differential corresponds to the sum of all ways of
  deleting a thorn or replacing a thorn by a stump. Filter the
  $E^1$-page by the number of internal vertices. The associated graded
  complex is acyclic unless $s=t=0$. At this point, the homology is
  entirely supported in dimension $0$, so the spectral sequences
  collapse.

  We have shown that, for each $n,w$, the homology of
  $\uwTreesCx{n}{w}$ is isomorphic to $\Field$ and, moreover, the
  homology is represented by the left-associated binary tree with all
  the weight in weight-1 popsicles. The map $\pi$ takes this binary tree to a
  generator for the homology of $\wTreesCx{n}{w}$. It follows that
  $\pi$ is a quasi-isomorphism.
\end{proof}

\begin{theorem}\label{thm:hu-mod-acyclic}
  For $n+2w\geq 2$, the projection
  $\uwMTreesCx{n}{w}\to \wMTreesCx{n}{w}$ induces an isomorphism on
  homology. For $(n,w)=(1,0)$ the homology of $\uwMTreesCx{n}{w}$ is
  generated by the identity tree. In particular, the homology of
  $\uwMTreesCx{n}{w}$ is one-dimensional, supported in dimension $0$, if
  $w=0$ and is trivial otherwise.
\end{theorem}
\begin{proof}
  The proof is essentially the same as the proof of Theorem~\ref{thm:hu-contractible}.
\end{proof}

\begin{theorem}\label{thm:hu-trans-contract}
  The projection of the homotopy unital weighted transformation trees
  complex to the weighted transformation trees complex induces an
  isomorphism on homology.
\end{theorem}
\begin{proof}
  Again, the proof is essentially the same as the proof of
  Theorem~\ref{thm:hu-contractible}.
\end{proof}

\begin{theorem}\label{thm:hu-mod-trans-acyclic}
  The projection of the homotopy unital weighted module transformation
  trees complex to the weighted module transformation trees complex
  induces an isomorphism on homology.
\end{theorem}
\begin{proof}
  Again, the proof is essentially the same as the proof of
  Theorem~\ref{thm:hu-contractible}.
\end{proof}

%%% Local Variables: 
%%% mode: latex
%%% TeX-master: "AbstractDiagonal.tex"
%%% TeX-command-extra-options: "-shell-escape"
%%% End: 

\section{Weighted diagonals}\label{sec:wDiags}

\begin{convention}\label{conv:wDiags-tens-Ring}
  Fix a commutative $\FF_2$-algebra $\Ring$. In this section, undecorated tensor
  products are over $\Ring$. (The rings $\Ground$ will not appear.)
\end{convention}

\subsection{Dimension and weight}
In this section, we work over the polynomial ring $\Ring[Y_1,Y_2]$. We extend
the notion of the dimension
\[
  \dim(T)=n+2w-v-1
\]
(Equation~\eqref{eq:w-tree-dim}) of a weighted tree to
$\wTreesCx{n}{w}\otimes\Ring[Y_1,Y_2]$ by declaring that the dimension
$\dim(Y_i)=0$; and similarly for weighted transformation trees, weighted module
trees, and so on. (Another grading will be used in Section~\ref{sec:wDiagApps},
so we use the word dimension here to avoid a conflict of terminology.) Extend
$\dim$ additively to tensor products.

Given a weighted tree $T\in \wTreesCx{n}{w}$, let $\wgr(T)=w$ be the total weight of $T$, i.e., the
sum of the weights of the vertices of $T$. Given an element
\[
  Y_1^aY_2^b(S\otimes T)\in \Ring[Y_1,Y_2]\otimes \wTreesCx{m}{v}\otimes\wTreesCx{n}{w}
\]
let
\begin{align*}
  \wgr_1(Y_1^aY_2^b(S\otimes T))&=a+v=a+\wgr(S)\\
  \wgr_2(Y_1^aY_2^b(S\otimes T))&=b+w=b+\wgr(T).
\end{align*}
Again, we make the analogous definitions for weighted transformation trees, weighted module
trees, and so on.

\subsection{Weighted algebra diagonals}
\label{sec:walg-diag}

In order to define a notion of weighted diagonals general enough for
our applications, it is convenient to consider a slightly larger
complex of trees. Specifically, we include the following
\emph{generalized weighted trees}:
\begin{itemize}
\item The tree $\IdTree$ with one input and no internal vertices. We
  define $\xwTrees{1}{0}$ to consist of this
  tree. Formula~\eqref{eq:w-tree-dim} gives $\dim(\IdTree)=0$.
\item A \emph{stump} $\stump$ with no inputs and no internal
  vertices. We define $\xwTrees{0}{0}$ to consist of this stump. 
  (The tree $\stump$ represents feeding in the identity in a
  weighted $\Ainf$-algebra.) We think of $\stump$ as having $-1$
  internal vertices, so Formula~\eqref{eq:w-tree-dim} gives
  $\dim(\stump)=0$ for the stump.
\end{itemize}
Neither of these is stable in the sense of
Section~\ref{sec:wAlgs}. Let $\xwTreesCx{n}{w}$ denote
$\wTreesCx{n}{w}$ if $w>0$ or $w=0$ and $n>1$,
$\xwTreesCx{1}{0}=\Ring\langle \xwTrees{1}{0}\rangle$, and
$\xwTreesCx{0}{0}=\Ring\langle\xwTrees{0}{0}\rangle$.

We extend the composition map $\circ_i$ to $\IdTree$ and $\stump$ as
follows:
\begin{itemize}
\item Composing with $\IdTree$ (in any way) is the identity map
  $\xwTreesCx{n}{w}\to\xwTreesCx{n}{w}$.
\item If $T$ is a stable weighted tree then the composition $T\circ_i
  \stump=0$ unless the $i\th$ input of $T$ feeds immediately into a
  valence $3$ (2-input, 1-output) vertex with weight $0$. If the
  successor of the $i\th$ input of $T$ has valence $3$ and weight $0$
  then $T\circ_i\stump$ is the result of erasing the $i\th$ input of
  $T$ and forgetting the successor of the $i\th$ input of $T$, so that
  $T\circ_i\stump$ is an element of $\xwTrees{n-1}{w}$ with one fewer
  internal vertex than $T$.
\end{itemize}
The differential on $\xwTreesCx{n}{w}$ is induced from the
differential on $\wTreesCx{n}{w}$ and, in particular, vanishes on
$\xwTreesCx{1}{0}$ and $\xwTreesCx{0}{0}$.

\begin{lemma}\label{lem:extend-circ}
  This extension of $\circ_i$ is a chain map.
\end{lemma}
\begin{proof}
  The fact that composing with $\IdTree$ (which is the identity map)
  is a chain map is obvious. To see that $\circ_i\stump$ is a chain
  map, suppose that the $i\th$ input of $T$ feeds into a vertex $v$ of
  valence $n$ and weight $w$. We have the following cases:
  \begin{enumerate}
  \item If $n>3$ then $T\circ_i\stump=0$. There are exactly two terms $S_1$, $S_2$ in $\bdy T$ so that $S_j\circ_i\stump\neq 0$: $S_1$ and $S_2$ are the two ways of splitting $v$ so that the $i\th$ input feeds into a valence $3$, weight $0$ input. Further, $S_1\circ_i\stump=S_2\circ_i\stump$, so $(\bdy T)\circ_i\stump=0$.
  \item If $n=2$ or $n=3$ and $w>0$ then again $T\circ_i\stump=0$ and
    there are exactly two terms $S_1,S_2$ in $\bdy T$ so that
    $S_j\circ_i\stump\neq 0$. For $n=3$, $S_1$ and $S_2$ are the
    result of pushing all of the weight of $v$ onto either of the two
    edges incident to $v$ not coming from the $i\th$ input. For $n=2$,
    $S_1$ and $S_2$ come from pushing all of $v$'s weight into a
    popsicle feeding into $v$ from either the left or the right.
  \item If $n=3$ and $w=0$ then $T\circ_i\stump$ is the result of
    deleting the $i\th$ input and forgetting the vertex $v$. Further
    for every term $S$ in $\bdy T$, the $i\th$ input of $S$ feeds into
    a valence $3$ vertex, so $S\circ_i\stump$ is the result of
    deleting the $i\th$ input of $S$ and forgetting this valence $3$
    vertex. In particular, $\bdy(T\circ_i\stump)=(\bdy T)\circ_i\stump$.
  \end{enumerate}
  This completes the proof.
\end{proof}

\begin{definition}
  A \emph{seed for a weighted algebra diagonal} is
  an $\Ring[Y_1,Y_2]$-linear combination $\wSeed$
  of the following terms:
  \[
  \wcorolla{0}{1}\otimes\wcorolla{0}{1},\
  Y_1\stump\otimes\wcorolla{0}{1},\ Y_2\wcorolla{0}{1}\otimes \stump.
  \]
  (This linear combination is to be interpreted as an element of
  $[(\xwTreesCx{0}{1}\otimes\xwTreesCx{0}{1})\oplus(\xwTreesCx{0}{0}\otimes\xwTreesCx{0}{1})\oplus(\xwTreesCx{0}{1}\otimes\xwTreesCx{0}{0})]\otimes\Ring[Y_1,Y_2]$.)
\end{definition}

\begin{definition}\label{def:w-alg-diag}
  A \emph{weighted algebra diagonal with seed $\wSeed$} is a collection of
  chain maps
  \begin{equation}\label{eq:wdiag-form}
    \wDiag{n}{w}\co \wTreesCx{n}{w}\to \bigoplus_{w_1,w_2\leq
      w}\xwTreesCx{n}{w_1}\otimes\xwTreesCx{n}{w_2}\otimes\Ring[Y_1,Y_2]
  \end{equation}
  with the following properties:
  \begin{itemize}
  \item \textbf{Dimension homogeneity}: The map
    $\wDiag{n}{w}$ is dimension-preserving, i.e.,
    \[
      \dim(\wDiag{n}{w}(T))=\dim(T).
    \]
    (Recall that $\dim(Y_i)=0$.)
  \item \textbf{Weight homogeneity}: The projection of the image of $\wDiag{n}{w}$ to
    $\xwTreesCx{n}{w_1}\otimes\xwTreesCx{n}{w_2}\otimes\Ring[Y_1,Y_2]$ is
    contained in
    $Y_1^{w-w_1} Y_2^{w-w_2}\xwTreesCx{n}{w_1}\otimes\xwTreesCx{n}{w_2}$. In other words,
    \[
      \wgr_1(\wDiag{n}{w}(T))=\wgr_2(\wDiag{n}{w}(T))=\wgr(T)=w.
    \]
  \item \textbf{Compatibility under stacking}: 
    \[
    \wDiag{n}{v+w}\circ \phi_{i,j,n;v,w} =
    \sum_{\substack{v_1+v_2\leq v\\w_1+w_2\leq w}}(\phi_{i,j,n;v_1,w_1}\otimes\phi_{i,j,n;v_2,w_2})\circ(\wDiag{j-i+1}{v}\otimes\wDiag{n+i-j}{w}).
    \]
  \item \textbf{Non-degeneracy}:
    \begin{enumerate}
    \item By the weight homogeneity condition, $\wDiag{2}{0}$ is a map
      $\Ring=\wTreesCx{2}{0}\to\Ring\otimes_\Ring\Ring=\wTreesCx{2}{0}\otimes\wTreesCx{2}{0}$.
      We require that this map is the canonical isomorphism. (This is the
      same as non-degeneracy in
      Definition~\ref{def:AssociahedronDiagonal}.)
    \item The image under $\wDiag{0}{1}$ of
      $\wcorolla{0}{1}\in \wTreesCx{0}{1}$ is the seed $\wSeed$.
    \item The image of $\wDiag{n}{w}$ is contained in 
      \[
      \bigoplus_{w_1,w_2\leq
        w}(\xwTreesCx{n}{w_1}\otimes\wTreesCx{n}{w_2}\otimes\Ring[Y_1,Y_2]+\wTreesCx{n}{w_1}\otimes\xwTreesCx{n}{w_2}\otimes\Ring[Y_1,Y_2]),
      \]
      i.e., at most one of each pair of trees in the diagonal is a
      generalized tree.
    \end{enumerate}
  \end{itemize}
\end{definition}

In particular, restricting to the diagonals $\wDiag{n}{0}$ gives an
associahedron diagonal.

\begin{figure}
  \centering
  %Font is 18 point
  \includegraphics[scale=.667]{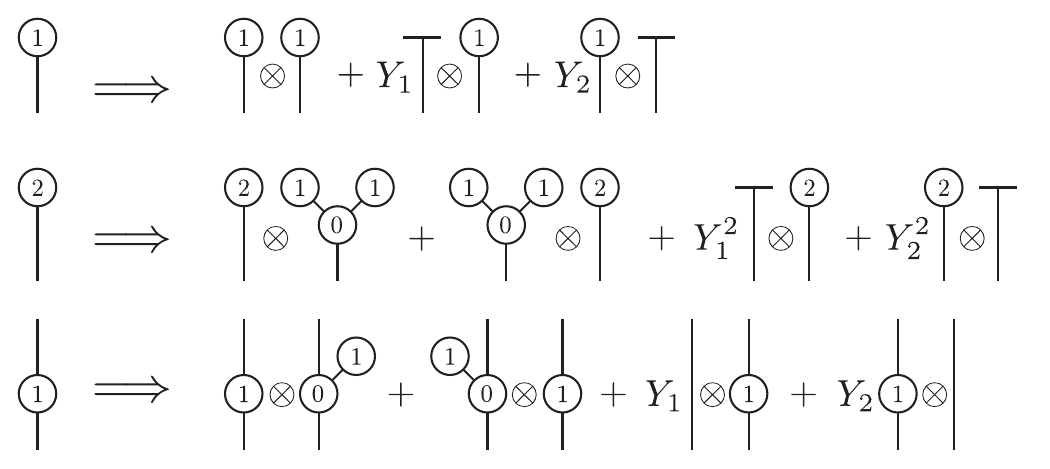}
  \caption[Terms in a maximal weighted algebra diagonal]{\textbf{Terms in a maximal weighted algebra diagonal.} The images of $\wcorolla{0}{1}$, $\wcorolla{0}{2}$, and
    $\wcorolla{1}{1}$ are shown. The tree $\stump$ is denoted as a vertical line
    with horizontal line on top. For the image of weight $0$ corollas
    see Figure~\ref{fig:diag-cells}.  We have drawn one of the two
    possible choices for $\wDiag{1}{1}$.}
  \label{fig:wdiag-terms}
\end{figure}

In the applications, we will be most interested in the \emph{maximal seed}
\[
\wSeed=\wcorolla{0}{1}\otimes\wcorolla{0}{1}+
Y_1\stump\otimes\wcorolla{0}{1}+Y_2\wcorolla{0}{1}\otimes \stump.
\]
The first few terms in a diagonal with this seed are shown in
Figure~\ref{fig:wdiag-terms}.

\begin{definition}\label{def:wDiagCells}
  A \emph{collection of weighted diagonal cells} consists of a chain 
  \[
  \wDiagCell{n}{w}\in \bigoplus_{w_1,w_2\leq
    w}Y_1^{w-w_1} Y_2^{w-w_2}\xwTreesCx{n}{w_1}\otimes\xwTreesCx{n}{w_2}
  \subset \bigoplus_{w_1,w_2\leq
    w}\xwTreesCx{n}{w_1}\otimes\xwTreesCx{n}{w_2}\otimes\Ring[Y_1,Y_2]
  \]
  of dimension $\dim(\wcorolla{n}{w})=n+2w-2$, for each $n,m\geq 0$ with
  $n+2w\geq 2$, satisfying the following properties:
  \begin{itemize}
  \item \textbf{Compatibility}: 
    \begin{equation}\label{eq:wCell-compat}
    \bdy\wDiagCell{n}{x}=\sum_{\substack{v+w=x\\1\leq i\leq n\\i-1\leq j\leq
        n}}
    \wDiagCell{n+i-j}{v}\circ_i\wDiagCell{j-i+1}{w}.
  \end{equation}
  \item \textbf{Non-degeneracy}:
    \begin{itemize}
    \item $\wDiagCell{2}{0}=\wcorolla{2}{0}\otimes\wcorolla{2}{0}$.
    \item $\wDiagCell{n}{w}$ is contained in 
      \[
      \bigoplus_{w_1,w_2\leq
        w}(\xwTreesCx{n}{w_1}\otimes\wTreesCx{n}{w_2}\otimes\Ring[Y_1,Y_2]+\wTreesCx{n}{w_1}\otimes\xwTreesCx{n}{w_2}\otimes\Ring[Y_1,Y_2]),
      \]
      i.e., at most one of each pair of trees in the diagonal is a
      generalized tree.
    \end{itemize}
  \end{itemize}
  The \emph{seed} of the collection of weighted diagonal cells is $\wDiagCell{0}{1}$.
\end{definition}

The compatibility condition can also be written
\begin{equation}\label{eq:wCell-compat-redux}
\bdy\wDiagCell{n}{x}=\sum_{\substack{v+w=x\\0\leq j\leq n\\1\leq i\leq
  n-j+1}}
\wDiagCell{n-j+1}{v}\circ_i \wDiagCell{j}{w}=\sum_{\substack{v+w=x\\j+k=n+1}}
\wDiagCell{k}{v}\circ \wDiagCell{j}{w}.
\end{equation}
with the understanding that
$\wDiagCell{1}{0}=\wDiagCell{0}{0}=0$. (For the abuse of notation
$\circ$ used here, see Section~\ref{sec:associahedron}.)

\begin{construction}\label{const:wcell-to-wdiag}
  Given a collection of diagonal cells $\{\wDiagCell{n}{w}\}$ define a
  collection of maps $\wDiag{n}{w}$ of the form~\eqref{eq:wdiag-form}
  inductively by the following two rules:
  \begin{itemize}
  \item $\wDiag{n}{w}(\wcorolla{n}{w})=\wDiagCell{n}{w}$.
  \item $\wDiag{n}{w}(T_1\circ_i
    T_2)=\wDiag{n_1}{w_1}(T_1)\circ_i\wDiag{n_2}{w_2}(T_2)$. (Here,
    $T_i$ has weight $w_i$ and $n_i$ inputs.)
  \end{itemize}
\end{construction}

\begin{lemma}
  Given a weighted algebra diagonal $\{\wDiag{n}{w}\}$ with seed
  $\wSeed=\wDiag{0}{1}(\wcorolla{0}{1})$, the chains
  $\{\wDiagCell{n}{w}=\wDiag{n}{w}(\wcorolla{n}{w})\}$ form a collection
  of weighted diagonal cells with seed $\wSeed$.

  Conversely, given a collection of weighted diagonal cells
  $\{\wDiagCell{n}{w}\}$, Construction~\ref{const:wcell-to-wdiag}
  defines a weighted algebra diagonal with the same seed.
\end{lemma}
\begin{proof}
  This is immediate from the definitions.
\end{proof}

\begin{theorem}\label{thm:wDiag-exists}
  Given any seed $\wSeed$ there is a weighted
  algebra diagonal with seed $\wSeed$.
\end{theorem}
\begin{proof}
  We must check that the right hand side of
  Equation~\eqref{eq:wCell-compat-redux} is a cycle. Then, for
  $(n,w)\not\in\{(2,0),(0,1),(3,0),(1,1)\}$ the dimension of the right
  hand side of Equation~\eqref{eq:wCell-compat} is greater than $0$,
  so it follows from Theorem~\ref{thm:AssociaplexAcyclic} that the
  right hand side is also a boundary. For the base cases, $(2,0)$ and
  $(0,1)$ are specified by the non-degeneracy condition, a solution
  for $\wDiagCell{3}{0}$ is shown in Figure~\ref{fig:diag-cells} and a
  solution for $\wDiagCell{1}{1}$ is shown in
  Figure~\ref{fig:wdiag-terms} (third line) for the maximal seed. Given another
  seed
  \[
    r_1\wcorolla{0}{1}\otimes\wcorolla{0}{1}+
    r_2 Y_1\stump\otimes\wcorolla{0}{1}+r_3 Y_2\wcorolla{0}{1}\otimes
    \stump,
  \]
  this base case is obtained by multiplying the first two pairs of
  trees in Figure~\ref{fig:wdiag-terms} by $r_1$, the third pair by
  $r_2$, and the fourth pair by $r_3$.

  To see that the right hand side of
  Equation~\eqref{eq:wCell-compat-redux} is a cycle, note that
  inductively
  \begin{align*}
    \bdy\sum_{\substack{v+w=x\\j+k=n+1}}
    \wDiagCell{k}{v}\circ \wDiagCell{j}{w}&=
    \sum_{\substack{v+w=x\\j+k=n+1}}
    (\bdy\wDiagCell{k}{v})\circ \wDiagCell{j}{w} +
    \wDiagCell{k}{v}\circ (\bdy\wDiagCell{j}{w})\\
    &=2\sum_{\substack{u+v+w=x\\j+k+\ell=n+2}}
    \wDiagCell{k}{u}\circ \wDiagCell{j}{v}\circ\wDiagCell{\ell}{w}=0.\qedhere
  \end{align*}
\end{proof}

\begin{remark}\label{rmk:extend-wDiag}
  The proof of Theorem~\ref{thm:wDiag-exists} shows that for any
  $W_0$, if we are given chains $\wDiagCell{n}{w}$ for all $w<W_0$,
  satisfying Condition~(\ref{eq:wCell-compat}) for all $x<W_0$, we can
  extend $\wDiagCell{n}{w}$ to a weighted algebra diagonal. In
  particular, any associahedron diagonal can be extended to a weighted
  algebra diagonal.
\end{remark}

\subsection{Weighted map diagonals}
For any $1\leq i\leq j\leq n$, $v, w\geq 0$;
$0=i_1\leq i_2<\dots<i_\ell=n$, and $w_1,\dots,w_\ell\geq 0$, there are chain maps
\begin{align*}
  \psi_{i,j;v,w}^n\co \wTreesCx{j-i}{v}\otimes \wTransCx{n+i-j+1}{w}&\to \wTransCx{n}{v+w}\\
  \xi^n_{i_1,\dots,i_\ell}\co \wTransCx{i_2-i_1}{w_1}\otimes\cdots\otimes \wTransCx{i_{\ell}-i_{\ell-1}}{w_{\ell-1}}
  \otimes\wTreesCx{\ell-1}{w_\ell}&\to \wTransCx{n}{w_1+\dots+w_\ell}.
\end{align*}
In terms of trees, the map $\psi^n_{i,j}$ corresponds to
$(S,T)\mapsto T\circ_iS$ (where $T$ is a weighted transformation tree
and $S$ is a (red) weighted algebra tree). The map
$\xi^n_{i_1,\dots,i_\ell}$ corresponds to
\[
(T_1,\dots,T_\ell,S)\mapsto S\circ(T_1,\dots,T_\ell)=((S\circ_{\ell}T_\ell)\circ_{\ell-1}\cdots)\circ_1 T_1
\]
(where the $T_i$ are weighted transformation trees and $S$ is a (blue)
weighted algebra tree).

We extend these maps to the generalized weighted trees complex, by declaring that:
\begin{itemize}
\item Composing a purple vertex $v$ with a red or blue identity $\IdTree$, in any valid
  way, gives the same purple vertex $v$.
\item Composing a purple vertex $v$ with a red stump \textcolor{red}{$\stump$} gives a blue stump
  \textcolor{blue}{$\stump$} if the purple vertex has $1$ input and weight $0$, and $0$ otherwise.
\end{itemize}
In particular, the second composition lands in the \emph{extended
  weighted transformation trees complex} $\xwTransCx{\ell}{w}$, which
agrees with $\wTransCx{\ell}{w}$ if $(\ell,w)\neq (0,0)$ and is
generated by the blue stump \textcolor{blue}{$\stump$} if
$(\ell,w)=(0,0)$.

\begin{figure}
  \centering
  %Font is 12 point
  \includegraphics[scale=.75]{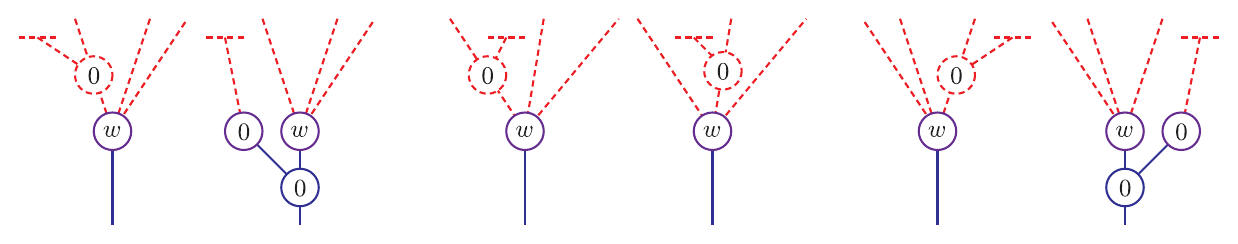}
  \caption[Differential of a purple corolla composed with a stump]{\textbf{Differential of a purple corolla composed with a stump.} The two non-zero terms in each of the three cases are shown.}
  \label{fig:purp-comp-stump}
\end{figure}

\begin{lemma}
  These extensions make $(S,T_1,\dots,T_\ell)\mapsto
  S\circ(T_1,\dots,T_\ell)$ and $(S,T)\mapsto S\circ_i T$ into chain maps
  \[
    \xwTreesCx{k}{v}\otimes \xwTransCx{\ell_1}{w_1}\otimes\cdots\otimes\xwTransCx{\ell_k}{w_k}\to
    \xwTransCx{\ell_1+\cdots+\ell_k}{v+w_1+\cdots+w_k}
  \]
  and
  \[
    \xwTransCx{k}{v}\otimes \xwTreesCx{\ell}{w}\to
    \xwTransCx{k+\ell-1}{v+w}
  \]
  and hence make the extensions of $\psi$ and $\xi$ to generalized
  weighted trees into chain maps.
\end{lemma}
\begin{proof}
  The fact that composing with $\IdTree$ is a chain map is obvious.  For
  $\stump$, there is one interesting case: the composition of an $n$-input,
  weight $w$ purple corolla, $n+2w\geq 2$, and a stump. There are three sub-cases, depending on
  whether $i=1$, $1<i<n$, or $i=n$; these are shown in
  Figure~\ref{fig:purp-comp-stump}. In all cases, both nonzero terms give an
  $(n-1)$-input, weight $w$ purple corolla.
\end{proof}

\begin{definition}\label{def:w-map-diag}
  Fix a weighted algebra diagonals $\wDiag{*}{*}_1$ and $\wDiag{*}{*}_2$.  A
  \emph{weighted map diagonal} compatible with $\wDiag{*}{*}_1$ and $\wDiag{*}{*}_2$
  consists of a sequence of 
  chain maps
  \begin{multline}
    \wMulDiag{n}{w}\co \wTransCx{n}{w}\to
    \bigoplus_{w_1,w_2\leq w}\left(\bigl(\wTransCx{n}{w_1}\otimes\xwTransCx{n}{w_2}\bigr)+\bigl(\xwTransCx{n}{w_1}\otimes\wTransCx{n}{w_2}\bigr)\right)\otimes\bigl(\Ring[Y_1,Y_2]\bigr)\\
    \subset\bigoplus_{w_1,w_2\leq w}\xwTransCx{n}{w_1}\otimes\xwTransCx{n}{w_2}\otimes\Ring[Y_1,Y_2]
  \end{multline}
  for $w\geq 0$ and $n\geq 0$, $(w,n)\neq (0,0)$,
  satisfying the following conditions:
  \begin{itemize}
    \item \textbf{Dimension homogeneity}: 
      The map $\wMulDiag{n}{w}$ is dimension-preserving, i.e.,
      \[
        \dim(\wMulDiag{n}{w}(T))=\dim(T).
      \]
    \item \textbf{Weight homogeneity:} The image of $\wMulDiag{n}{w}$
      is contained in
      $Y_1^{w-w_1}Y_2^{w-w_2}\xwTransCx{n}{w}\otimes\xwTransCx{n}{w}$
      or, equivalently,
      \[
        \wgr_1(\wMulDiag{n}{w}(T))=\wgr_2(\wMulDiag{n}{w}(T))=\wgr(T)=w.
      \]
    \item \textbf{Compatibility under stacking}:
    \begin{align*}
      \wMulDiag{n}{v+w}\circ \psi^n_{i,j;v,w}&=
      \sum_{\substack{v_1+v_2\leq v\\ w_1+w_2\leq w}}
        (\psi_{i,j;v_1,w_1}^{n}\otimes\psi_{i,j;v_2,w_2}^{n})\circ(\wDiag{j-i}{v}_1\otimes\wMulDiag{n+i-j+1}{w}) \\
      \wMulDiag{n}{w}\circ \xi_{i_1,\dots,i_\ell}^n&=
      \sum_{w_1+\cdots+w_\ell=w}(\xi_{i_1,\dots,i_\ell}^n\otimes \xi_{i_1,\dots,i_\ell}^n)\circ(\wMulDiag{i_2-i_1}{w_1}\otimes \cdots\otimes \wMulDiag{i_\ell-i_{\ell-1}}{w_{\ell-1}}\otimes \wDiag{\ell}{w_\ell}_2),
    \end{align*}
    with the understanding that the compositions on the right-hand
    side involve shuffling of factors (compare
    Formula~\eqref{eq:as-diag-shuf}).
  \item \textbf{Non-degeneracy}:
    $\wMulDiag{1}{0}(\wpcorolla{1}{0})=\wpcorolla{1}{0}\otimes\wpcorolla{1}{0}$.
    (Note that $\wTransCx{1}{0}\cong \Ring\langle\wpcorolla{1}{0}\rangle$.)
  \end{itemize}
\end{definition}

Weighted map diagonals can be phrased in terms of trees as follows:

\begin{definition}\label{def:w-map-diag-cell}
  In terms of trees, if $\wDiag{*}{*}_1$ and $\wDiag{*}{*}_2$ correspond to 
  collections of weighted diagonal cells $\wDiagCell{*}{*}_1$ and $\wDiagCell{*}{*}_2$, a
  weighted map diagonal is specified by a sequence of {\em weighted map cells}, which are elements
  \begin{multline*}
    \wTrMulDiag{n}{w}=\wMulDiag{n}{w}(\wcorolla{n}{w})\in  \bigoplus_{w_1,w_2\leq w}Y_1^{w-w_1}Y_2^{w-w_2}\left(\xwTransCx{n}{w_1}\otimes\wTransCx{n}{w_2}+\wTransCx{n}{w_1}\otimes\xwTransCx{n}{w_2}\right)\\
    \subset \xwTransCx{n}{w_1}\otimes\xwTransCx{n}{w_2}\otimes\Ring[Y_1,Y_2]
  \end{multline*}
  in dimension $n+2w-1$ satisfying:
  \begin{itemize}
  \item \textbf{Compatibility}: 
  \begin{equation}\label{eq:wMulDiagCompatTree}
    \bdy(\wTrMulDiag{n}{w})=\sum_{\substack
        {i+j=n+1 \\
      u+v=w}}\wTrMulDiag{i}{u}\circ \wDiagCell{j}{v}_1+ 
      \sum_{k=0}^\infty\sum_{\substack {m_1+\cdots+m_k=n\\ v+w_1+\dots+w_k=w}
      }
      \wDiagCell{k}{v}_2\circ(\wTrMulDiag{m_1}{w_1},\dots,\wTrMulDiag{m_k}{w_k}).
  \end{equation}
  (Note that in the second sum, $k=0$ only occurs if $n=0$.)
\item \textbf{Non-degeneracy}:
  $\wTrMulDiag{1}{0}=\wpcorolla{1}{0}\otimes\wpcorolla{1}{0}\in\wTreesCx{1}{0}\otimes\wTreesCx{1}{0}$ is the (unique)
  pair of $1$-input, weight $0$ transformation trees (with $1$ internal vertex
  each).
  \end{itemize}
\end{definition}

\begin{figure}
  \centering
  %Font is 12 point.
  \includegraphics{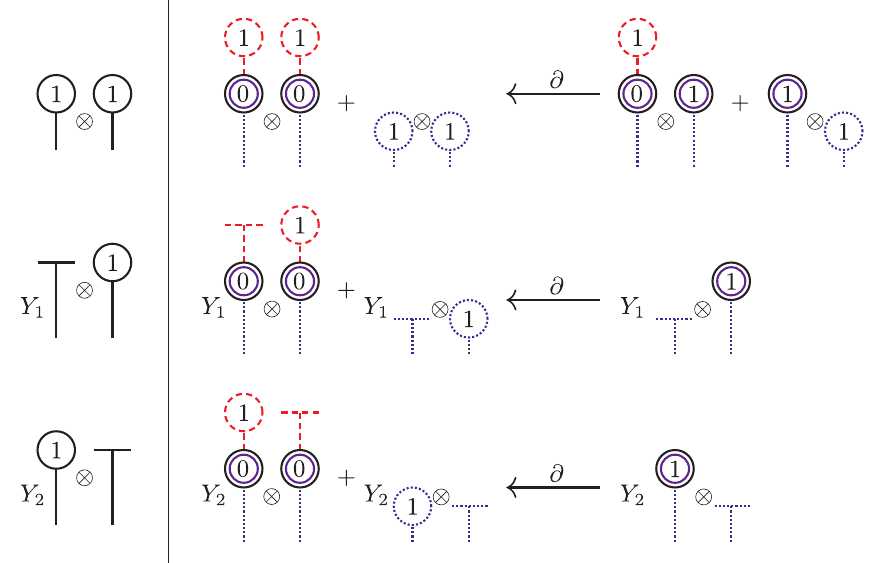}
  \caption[Existence of weighted map diagonals]{\textbf{Existence of weighted map diagonals.} The base case
    for each possible term in the seed is shown. On the left is the
    term in the seed, in the center are the corresponding terms on the
    right side of Equation~\eqref{eq:wMulDiagCompatTree}, and on the
    right is an option for the weighted map diagonal.}
  \label{fig:w-map-diag-base-case}
\end{figure}

\begin{lemma}\label{lem:wmul-diag-exists}
  Given any two weighted algebra diagonals $\wDiag{*}{*}_1$ and
  $\wDiag{*}{*}_2$ with the same seeds, there exists a weighted map diagonal
  $\wMulDiag{*}{*}$ compatible with $\wDiag{*}{*}_1$ and
  $\wDiag{*}{*}_2$.
\end{lemma}
\begin{proof}
  By Proposition~\ref{prop:wTransCxContract}, the weighted
  transformation trees complex is contractible. So, as usual, it
  suffices to verify that the right-hand side of the weighted
  map diagonal equation~\eqref{eq:wMulDiagCompatTree} is a
  cycle and that when the right-hand side lies in dimension 0, it is a
  boundary. We leave the first statement, that the right-hand side is
  a cycle, to the reader. For the second, note that
  \begin{align*}
    \dim\bigl(\wTrMulDiag{i}{u}\circ \wDiagCell{j}{v}_1\bigr)&=i+j+2u+2v-3\\
    \dim\bigl(\wDiagCell{k}{v}_2\circ(\wTrMulDiag{m_1}{w_1},\dots,\wTrMulDiag{m_k}{w_k})\bigr)&=2(v+w_1+\cdots+w_k)-2+m_1+\cdots+m_k.
  \end{align*}
  Thus, the dimension $0$ part consists of
  $\wTrMulDiag{1}{0}\circ\wDiagCell{2}{0}_1$ and
  $(\wDiagCell{2}{0}_2\circ_2\wTrMulDiag{1}{0})\circ_1\wTrMulDiag{1}{0}$,
  which lie in the classical multiplihedron, and
  $\wTrMulDiag{1}{0}\circ\wDiagCell{0}{1}_1$ and
  $\wDiagCell{0}{1}_2$. This reduces to a case check,
  depending on the seed of the weighted algebra diagonals; see
  Figure~\ref{fig:w-map-diag-base-case}.
\end{proof}

\subsection{Weighted module diagonals}
Like we did with $\wTreesCx{n}{w}$, we extend $\wMTreesCx{n}{w}$ by
allowing the $1$-input, $0$-internal vertex tree
$\IdTree$. Specifically, let $\xwMTreesCx{n}{w}=\wMTreesCx{n}{w}$
except that $\xwMTreesCx{1}{0}$ is generated by $\IdTree$ (rather than
being trivial).
\begin{definition}\label{def:w-mod-diag}
  Fix a seed $\wSeed$ and a weighted algebra diagonal $\wDiag{n}{w}$
  with seed $\wSeed$. A \emph{weighted (right) module diagonal}
  compatible with $\wDiag{n}{w}$ is a
  collection of chain maps
  \begin{multline}\label{eq:wmdiag-form}
  \wModDiag{n}{w}\co \wMTreesCx{n}{w}\to \bigoplus_{w_1,w_2\leq
    w}\left(\wMTreesCx{n}{w_1}\otimes\xwMTreesCx{n}{w_2}+\xwMTreesCx{n}{w_1}\otimes\wMTreesCx{n}{w_2}\right)\otimes\Ring[Y_1,Y_2]\\
  \subset\bigoplus_{w_1,w_2\leq
    w}\xwMTreesCx{n}{w_1}\otimes\xwMTreesCx{n}{w_2}\otimes\Ring[Y_1,Y_2]
  \end{multline}
  with the following properties:
  \begin{itemize}
  \item \textbf{Dimension homogeneity}:
    The map $\wModDiag{n}{w}$ is grading-preserving, i.e.,
    \[
      \dim(\wModDiag{n}{w}(T))=\dim(T).
    \]
  \item \textbf{Weight homogeneity}: The projection of the image of $\wModDiag{n}{w}$ to
    $\xwMTreesCx{n}{w_1}\otimes\xwMTreesCx{n}{w_2}\otimes\Ring[Y_1,Y_2]$ is
    contained in
    $Y_1^{w-w_1}Y_2^{w-w_2}\xwMTreesCx{n}{w_1}\otimes\xwMTreesCx{n}{w_2}$. In
    other words,
    \[
      \wgr_1(\wModDiag{n}{w}(T))=\wgr_2(\wModDiag{n}{w}(T))=\wgr(T)=w.
    \]
  \item \textbf{Compatibility under stacking}: 
    \[
    \wModDiag{n}{v+w}\circ \phi_{i,j,n;v,w} =
    \begin{cases}
    \sum_{\substack{v_1+v_2\leq v\\w_1+w_2\leq
        w}}(\phi_{1,j,n;v_1,w_1}\otimes\phi_{1,j,n;v_2,w_2})\circ(\wModDiag{j}{v}\otimes\wModDiag{n+1-j}{w})
      & i=1\\
    \sum_{\substack{v_1+v_2\leq v\\w_1+w_2\leq
        w}}(\phi_{i,j,n;v_1,w_1}\otimes\phi_{i,j,n;v_2,w_2})\circ(\wDiag{j-i+1}{v}\otimes\wModDiag{n+i-j}{w})
    & i>1
    \end{cases}
    \]
  \item \textbf{Non-degeneracy}:
    \begin{enumerate}
    \item By the weight homogeneity condition, the map
      $\wModDiag{2}{0}$ is a map
      $\Ring=\wMTreesCx{2}{0}\to\Ring\otimes_\Ring\Ring=\wMTreesCx{2}{0}\otimes\wMTreesCx{2}{0}$.
      We require that $\wModDiag{2}{0}$ is the canonical isomorphism,
      i.e., $\wModDiag{2}{0}(\wcorolla{2}{0})=\wcorolla{2}{0}\otimes\wcorolla{2}{0}$. (This
      is the same as non-degeneracy in Definition~\ref{def:ModuleDiagonal}.)
    \end{enumerate}
  \end{itemize}

  Equivalently, we can view $\wModDiag{n}{w}$ as a formal linear
  combination of pairs of trees 
  \begin{multline}
    \wModDiagCell{n}{w}\in \bigoplus_{w_1,w_2\leq w}
    Y_1^{w-w_1}Y_2^{w-w_2}\left(\xwMTreesCx{n}{w_1}\otimes\wMTreesCx{n}{w_2}+\wMTreesCx{n}{w_1}\otimes\xwMTreesCx{n}{w_2}\right)
    \\
    \subset \bigoplus_{w_1,w_2\leq w}
    \xwMTreesCx{n}{w_1}\otimes\xwMTreesCx{n}{w_2}\otimes\Ring[Y_1,Y_2]
  \end{multline}
  of grading $\dim(\wcorolla{n}{w})=n+2w-2$, satisfying the following
  conditions:
  \begin{itemize}
  \item \textbf{Compatibility}: 
    \begin{equation}\label{eq:wModDiagCell-compat}
    \bdy\wModDiagCell{n}{x}=\sum_{v+w=x}\left(
      \sum_{j=1}^n\wModDiagCell{n-j+1}{v}\circ_1\wModDiagCell{j}{w}+\sum_{2\leq
      i\leq j\leq n}\wModDiagCell{n+1+i-j}{v}\circ_i\wDiagCell{j-i}{w}
    \right).
    \end{equation}
  \item \textbf{Non-degeneracy}:
    $\wModDiagCell{2}{0}=\wcorolla{2}{0}\otimes\wcorolla{2}{0}$.
  \end{itemize}
\end{definition}

\begin{figure}
  \centering
  %Font is 18 point
  \includegraphics[scale=.667]{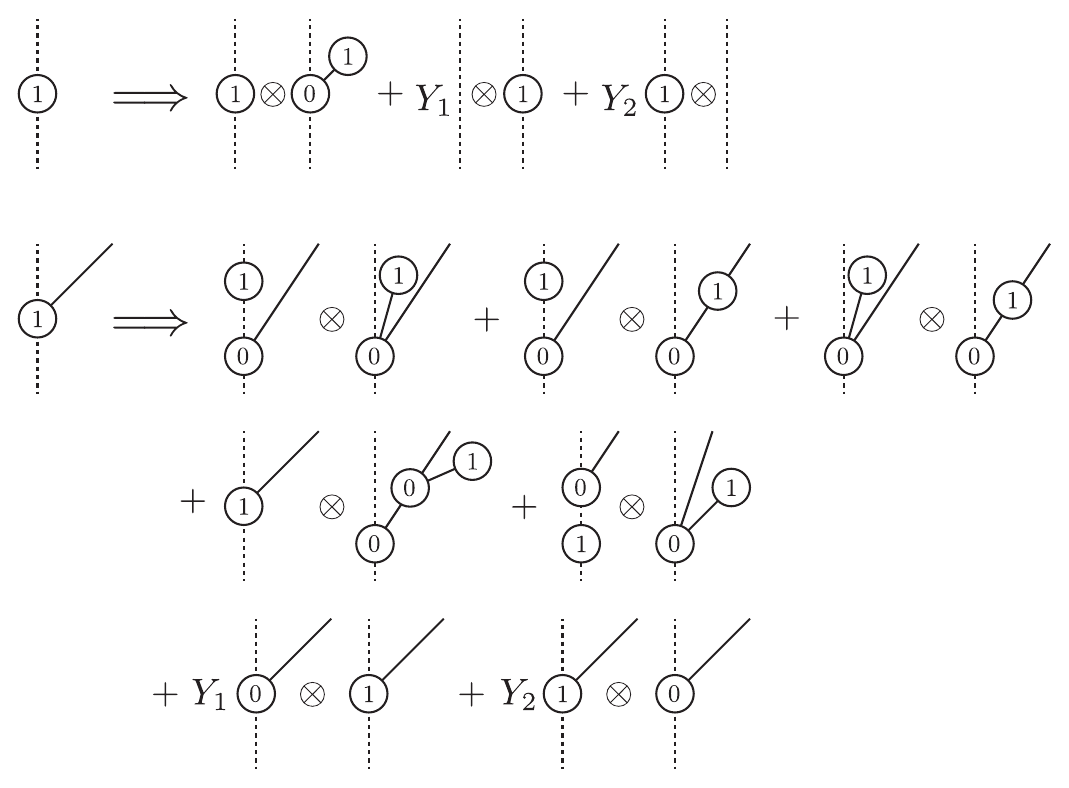}
  \caption[Terms in a weighted module diagonal compatible with
      Figure~\ref{fig:wdiag-terms}]{\textbf{Terms in a weighted module diagonal compatible with
      Figure~\ref{fig:wdiag-terms}.} We have shown $\wModDiagCell{1}{1}$ and $\wModDiagCell{2}{1}$.}
  \label{fig:wmod-diag-terms}
\end{figure}

The first few terms in a particular weighted module diagonal
are shown in Figure~\ref{fig:wmod-diag-terms}.

A module diagonal can be expressed  alternatively as a formal power series, in the following sense.
Let
$\wModDiag{*}{*}\llbracket
Y_1,Y_2\rrbracket=\wModDiag{*}{*}\otimes\Ring\llbracket
Y_1,Y_2\rrbracket$, a module over the ring of formal power series $\Ring\llbracket Y_1,Y_2\rrbracket$.
A module diagonal forms an element 
\[ \wModDiagCell{*}{*}=\prod_{n,w=0}^{\infty}\wModDiagCell{n}{w}\in\wModDiag{*}{*}\llbracket Y_1,Y_2\rrbracket.\]

\begin{theorem}\label{thm:wModDiag-exists}
  Given any weighted algebra diagonal $\wDiagCell{n}{w}$ there is a weighted
  module diagonal $\wModDiagCell{n}{w}$ compatible with $\wDiagCell{n}{w}$.
\end{theorem}
\begin{proof}
  The image of $\wDiagCell{n}{w}$ in 
  \[
    \bigoplus_{w_1,w_2\leq w}
    \xwMTreesCx{n}{w_1}\otimes\xwMTreesCx{n}{w_2}\otimes\Ring[Y_1,Y_2]
    =
    \bigoplus_{w_1,w_2\leq w}
    (\xwTreesCx{n}{w_1}/L_*)\otimes(\xwTreesCx{n}{w_2}/L_*)\otimes\Ring[Y_1,Y_2]
  \]
  is a weighted module diagonal.
\end{proof}

\begin{remark}
  As in the case of weighted algebra diagonals
  (Remark~\ref{rmk:extend-wDiag}), given terms $\wModDiagCell{1+n}{w}$
  for all $w<W_0$ satisfying Equation~\ref{eq:wModDiagCell-compat} for
  all $x<W_0$, these terms can be extended to a weighted module diagonal.
\end{remark}

\subsection{Weighted primitives}\label{sec:w-prim}
In this section, we define primitives for weighted module diagonals. We work over a weighted algebra diagonal with maximal seed. Module diagonal primitives with respect to certain other seeds can be obtained (and defined) by setting $Y_1=0$ and/or $Y_2=0$, and the more general case is left to the reader.

Before defining weighted module diagonal primitives, we extend the
definitions of root joining and left-root joining to the weighted
case. Given weighted trees $S_1,\dots,S_n$ the
\emph{root joining} of $S_1,\dots,S_n$ is the sum over all
non-negative weights $w$ of the result of joining the outputs of
$S_1,\dots,S_n$ into a single new vertex of weight $w$, multiplied by $Y_1^w$. In the case
$n=1$, we require that the new vertex have positive weight. In the
case $n=0$ we define the root joining of zero trees to be
$\sum_{w\geq 1}Y_1^w \wcorolla{0}{w}\in\wTreesCx{0}{w}\llbracket Y_1\rrbracket$.
In formulas, given $n\geq 0$ and weighted trees $S_1,\dots,S_n$, define
\begin{equation}
  \label{eq:wRootJoin}
  \wRootJoin^{w}(S_1,\dots,S_n)
  =Y_1^w\wcorolla{n}{w}\circ(S_1,\dots,S_n)
\end{equation}
and
\begin{equation}
  \wRootJoin^{*}(S_1,\dots,S_n)=\sum_{w=0}^{\infty}\wRootJoin^w(S_1,\dots,S_n).
\end{equation}

For weighted trees, left joining is defined exactly as it was in the
unweighted case. Define the left joining of the empty list of trees
(the case $n=0$) to be the generalized tree $\IdTree$.

Given a sequence
$(T_1,S_1),(T_2,S_2),\dots,(T_n,S_n)$ of pairs of trees, with each
$T_i$ having one more input than the corresponding~$S_i$, define the
\emph{left-root joining} of the sequence to be
\begin{align}
  \LRjoinW^w((T_1,S_1),\dots,(T_n,S_n))&=\LeftJoin(T_1,\dots,T_n)\otimes\wRootJoin^{w}(S_1,\dots,S_n)\\
  \label{eq:LRjoinW}
  \LRjoinW^*((T_1,S_1),\dots,(T_n,S_n))&=\LeftJoin(T_1,\dots,T_n)\otimes\wRootJoin^{*}(S_1,\dots,S_n)\\
                                         &=\sum_w \LRjoinW^w((T_1,S_1),\dots,(T_n,S_n))\nonumber.
\end{align}
Extend $\LRjoinW^*$ multi-linearly to a function
\[
  \LRjoinW^*\co (\xwMTreesCx{*}{*}\otimes\xwTreesCx{*}{*})^{\otimes n}\llbracket Y_1,Y_2\rrbracket\to (\xwMTreesCx{*}{*}\otimes\xwTreesCx{*}{*})\llbracket Y_1,Y_2\rrbracket.
\]
for $n \ge 0$.

It follows from the definitions that
\[ \LRjoinW^*()= \sum_{w=1}^{\infty} Y_1^{w} (\IdTree\otimes \wcorolla{w}{0})\in (\xwMTreesCx{*}{*}\otimes \wTreesCx{*}{*})\llbracket Y_1 \rrbracket.\]

\begin{definition}\label{def:wM-prim}
  Fix a weighted algebra diagonal $\wDiagCell{n}{w}$.  A \emph{(right)
  weighted module diagonal primitive compatible with
  $\wDiagCell{*}{*}$} consists of a linear combination of trees
  \begin{equation}\label{eq:wM-prim-1}
  \wTrPMDiag{n}{w}=\sum_{(S,T)}(S,T)\in \bigoplus_{w_1,w_2\leq w}Y_1^{w-w_1} Y_2^{w-w_2}\bigl(\wMTreesCx{n}{w_1}\otimes \xwTreesCx{n-1}{w_2}\bigr),
  \end{equation}
  of dimension $n+2w-2$, for each $n \geq 1$ and $w\geq 0$, $(n,w)\neq(1,0)$, satisfying the following conditions:
  \begin{itemize}
  \item \textbf{Compatibility}:
    \begin{equation}
      \label{eq:wM-prim-compat-long}
      \partial \wTrPMDiag{n}{w} = \hspace{-1.5em}
      \sum_{\substack{v+w_1+\cdots+w_k=w\\n_1+\cdots+n_k=n+k-1}}\hspace{-1.5em}
      \LRjoinW^v(\wTrPMDiag{n_1}{w_1}\otimes\cdots\otimes\wTrPMDiag{n_k}{w_k})
      +
      \hspace{-.5em}\sum_{\substack{w_1+w_2=w\\n_1+n_2=n+1}}\hspace{-.5em}\wTrPMDiag{n_1}{w_1}\circ'\wDiagCell{n_2}{w_2}
    \end{equation}
    or, more succinctly,
    \begin{equation}
      \label{eq:wM-prim-compat}
      \partial \wTrPMDiag{*}{*} = 
      \LRjoinW^*((\wTrPMDiag{*}{*})^{\otimes\bullet})+\wTrPMDiag{*}{*}\circ'\wDiagCell{*}{*}.
    \end{equation}
    Here,
    $\wTrPMDiag{*}{*}=\sum_{n,w}\wTrPMDiag{n}{w}$
    and
    $(\wTrPMDiag{*}{*})^{\otimes\bullet} = \sum_{n \ge 0}(\wTrPMDiag{*}{*})^{\otimes n}$.
  \item \textbf{Non-degeneracy}: $\wTrPMDiag{2}{0}=\wcorolla{2}{0}\otimes \IdTree$.
  \end{itemize}
\end{definition}

The weight-zero piece of a weighted module diagonal primitive
$\wTrPMDiag{n}{0}$ is a
module diagonal primitive (Definition~\ref{def:M-prim}).  Note also
that while the definition of $\LRjoinW^*$ uses a power series ring, the
coefficients of $\wTrPMDiag{n}{w}$ are polynomials in $Y_1$ and $Y_2$,
not power series.

See Figure~\ref{fig:wdiag-terms} for some positive-weight terms in a
weighted module diagonal primitive (and Figure~\ref{fig:mprim-terms}
for some weight-zero terms).

\begin{figure}
  \centering
  \includegraphics[scale=.666667]{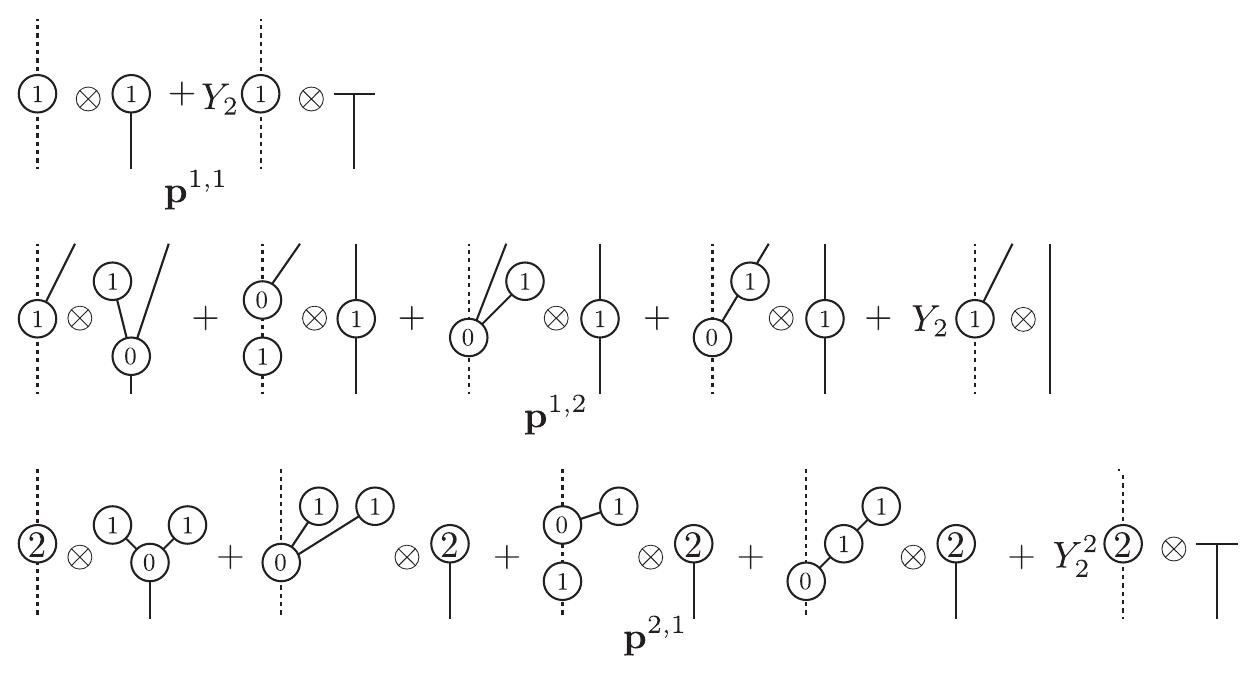}
  %Font is 18 point
  \caption[Terms in a weighted module diagonal primitive]{\textbf{Terms in a weighted module diagonal primitive.} The module diagonal primitive shown is compatible with $\wDiagCell{n}{w}$ from Figure~\ref{fig:wdiag-terms} and the weight-0 module diagonal primitive in Figure~\ref{fig:mprim-terms}.  Top: $\wTrPMDiag{1}{1}$. Middle: $\wTrPMDiag{2}{1}$. Bottom: $\wTrPMDiag{1}{2}$.}
  \label{fig:wmod-prim-diag-terms}
\end{figure}
%Check of third line in figure is in "calc" folder.

The construction of primitives hinges on the following identities
(cf. Lemma~\ref{lem:join-diff}):

\begin{lemma}\label{lem:wjoin-diff}
  Given $x_1\otimes\cdots\otimes x_n\in (\Trees\otimes\Trees)^{\otimes n}$,
  the operation $\LRjoinW^*$ satisfies
  \begin{multline*}
  \bdy(\LRjoinW^*(x_1\otimes\cdots\otimes x_n))=\sum_{i=1}^n \LRjoinW^*(x_1\otimes\cdots\otimes \bdy(x_i)\otimes\cdots\otimes x_n) \\
  + \sum_{j=0}^{n-1}\sum_{i=0}^{n-j}\LRjoinW^*(x_1\otimes\cdots\otimes x_i\otimes \LRjoinW^*(x_{i+1}\otimes\cdots\otimes x_{i+j})\otimes x_{i+j+1}\otimes\cdots\otimes x_n).
  \end{multline*}
  The operation $\circ'$ satisfies $\bdy(x\circ'y)=\bdy(x)\circ' y+x\circ'\bdy(y)$.
  Finally,
  \[ \LRjoinW^*(x_1\otimes\dots\otimes x_n)\circ'y=\sum_{i=1}^n \LRjoinW^*(x_1,\dots,x_{i-1}, x_i\circ' y,x_{i+1},\dots,x_n).\]
\end{lemma}
\begin{proof}
  The proof is straightforward.
\end{proof}

Note that, by definition, $\LRjoinW^*()\circ' y = 0$.

Proposition~\ref{prop:prim-exist} has the following generalization to the weighted case.

\begin{proposition}\label{prop:wprim-exist}
  For any weighted algebra diagonal, $\wDiagCell{*}{*}$, there exists a
  compatible weighted module diagonal primitive $\wTrPMDiag{*}{*}$.
\end{proposition}
\begin{proof}
  For fixed $w$ and $n$, Equation~\eqref{eq:wM-prim-compat}
  expresses $\partial \wTrPMDiag{n}{w}$ in terms of elements $\wTrPMDiag{n'}{w'}$ 
  with $w'<w$ or $w=w'$ and $n'<n$.
  Moreover, if each $\wTrPMDiag{n'}{w'}$ with $w'<w$ or $w=w'$ and $n'<n$ is contained in
  grading $n'+2w'-2$, then the right hand side of Equation~\eqref{eq:wM-prim-compat}
  is contained in grading $n+2w-3$. This is true because if
  $x_1,\dots,x_m$ is a sequence of pairs of trees with
  $x_i\in Y_1^{w_i-u_i}Y_2^{w_i-v_i}\wMTreesCx{n_i}{u_i}\otimes \wTreesCx{n_i-1}{v_i}$
  with grading $n_i+2w_i-2$, then 
  \[
    \LRjoinW^*(x_1,\dots,x_m)\in
    \prod_{s=0}^\infty Y_1^{w-u+s}Y_2^{w-v}\wMTreesCx{n}{u}\otimes
    \wTreesCx{n-1}{v+s},
  \]
  where $n=1-m+\sum_{i=1}^m n_i$, $u=\sum_{i=1}^m u_i$, and
  $v=\sum_{i=1}^mv_i$, and the grading of $\LRjoinW^*(x_1,\dots,x_m)$ is
  \[
    2s+m-2+\sum_{i=1}^m (n_i+2w_i-2)=n-3+2s+2\sum_{i=1}^mw_i=\gr(\wTrPMDiag{n}{w})-1
  \]
  (where $w=s+\sum w_i$).
  
  Thus, we can prove the existence of $\wTrPMDiag{n}{w}$ by induction first on $w$ and then, for each $w$, on $n$.

  The base case $w=0$ is an unweighted module diagonal
  primitive, which was shown to exist in Proposition~\ref{prop:prim-exist}. 
  Now, fix $w>0$ and
  suppose that $\wTrPMDiag{n}{w}$ has been constructed for all $w'<w$.
  Consider first  the base cases when
  $n+2w-2\in\{0,1\}$, i.e., $(n,w)\in\{(2,0),(3,0),(1,1)\}$; the first
  is specified in the non-degeneracy condition, the second in
  Figure~\ref{fig:mprim-terms}, and the third in
  Figure~\ref{fig:wmod-prim-diag-terms}.
  
  For the inductive step on $n$,
  the structure equation for $\wTrPMDiag{n}{w}$ and Lemma~\ref{lem:wjoin-diff} ensure that
  \[
    \partial(\LRjoinW^*((\wTrPMDiag{*}{*})^{\otimes\bullet})+\wTrPMDiag{*}{*}\circ'\wDiagCell{*}{*})^{n,w}=0.
  \]
  Since this cycle has dimension $n+2w-2$ with $n,w\geq 1$, 
  Theorems~\ref{thm:AssociaplexAcyclic} and~\ref{thm:ModAssociaplexAcyclic} guarantees that there is an element $\wTrPMDiag{n}{w}$ with 
  \[
    \partial \wTrPMDiag{n}{w}=(\LRjoinW^*((\wTrPMDiag{*}{*})^{\otimes\bullet})+\wTrPMDiag{*}{*}\circ'\wDiagCell{*}{*})^{w,n},
  \]
  as needed.
\end{proof}

To construct weighted module diagonals from primitives, we use the following generalization of $\LRjoinW^{\prime}$ to the weighted case.  For all $k\geq 0$, define
\begin{align}\label{eq:LRjoinWprime}
    \LRjoinW^{\prime,w}\bigl((S_1,T_1),\dots,(S_k,T_k)\bigr) &=
    \LeftJoin(S_1,\dots,S_k) \otimes
    \wRootJoin^{w}(\IdTree,T_1,\dots,T_k)\\    
  \LRjoinW^{\prime,*}\bigl((S_1,T_1),\dots,(S_k,T_k)\bigr) &=
    \LeftJoin(S_1,\dots,S_k) \otimes
    \wRootJoin^{*}(\IdTree,T_1,\dots,T_k)\\
    &=\sum_w\LRjoinW^{\prime,w}\bigl((S_1,T_1),\dots,(S_k,T_k)\bigr)
    \in
    \xwMTreesCx{*}{*}\otimes\xwMTreesCx{*}{*}\otimes\Ring\llbracket
    Y_1,Y_2\rrbracket,\nonumber
\end{align}
with the understanding that, for the case $k=0$,
\[ \wRootJoin^{*}(\IdTree) = \sum_{w=1}^{\infty} Y_1^w\wcorolla{1}{w}\in \wMTreesCx{*}{*}\otimes\Ring\llbracket Y_1\rrbracket,\]
and therefore
\[\LRjoinW^{\prime,*}() = \sum_{w=1}^{\infty} Y_1^w(\IdTree\otimes \wcorolla{1}{w})\in\xwMTreesCx{*}{*}\otimes\wMTreesCx{*}{*}\otimes\Ring\llbracket Y_2\rrbracket.\]

We generalize  Lemma~\ref{lem:mod-join-diff}. 

\begin{lemma}\label{lem:mod-join-diffw}
  Given $x_1\otimes\cdots\otimes x_n\in (\wMTreesCx{*}{*}\otimes\wTreesCx{*}{*})^{\otimes n}$,
  the operation $\LRjoinW^{\prime,*}$ satisfies
  \begin{align*}
  \bdy(\LRjoinW^{\prime,*}(x_1\otimes\cdots\otimes x_n))
   &=\sum_{i=1}^n \LRjoinW^{\prime,*}(x_1\otimes\cdots\otimes \bdy(x_i)\otimes\cdots\otimes x_n) \\
  &+ \sum_{i=0}^{n} \LRjoinW^{\prime,*}(x_1\otimes\dots\otimes x_i) \circ_1 \LRjoinW^{\prime,*}(x_{i+1}\otimes\dots\otimes x_n)\\
  &+ \sum_{j=0}^n\sum_{i=0}^{n-j}\LRjoinW^{\prime,*}(x_1\otimes\cdots\otimes \LRjoinW^*(x_{i+1}\otimes\cdots\otimes x_{i+j})\otimes\cdots\otimes x_n).
  \end{align*}
\end{lemma}

\begin{proof}
  The proof is straightforward. 
\end{proof}

The above formulas give the recipe for constructing weighted module diagrams from their primitives:

\begin{lemma}\label{lem:wprim-gives-diag}
  Given a collection of weighted diagonal cells $\wDiagCell{n}{w}$, a compatible module
  diagonal primitive $\wTrPMDiag{*}{*}$ gives rise to a module diagonal compatible with $\wDiag{*}{*}$ by the
  formula
  \[
    \wModDiagCell{n}{w}=\sum_{\substack{v+w_1+\cdots+w_k=w\\n_1+\cdots+n_k=n+k-1}}\LRjoinW^{\prime,v}(\wTrPMDiag{n_1}{w_1}\otimes\cdots\otimes\wTrPMDiag{n_k}{w_k})
  \]
  or, more succinctly,
  \[
    \wModDiagCell{*}{*}=\LRjoinW^{\prime,*}((\wTrPMDiag{*}{*})^{\otimes\bullet}).
  \]
\end{lemma}
\begin{proof}
  This is an immediate consequence of the structure equations and
  Lemma~\ref{lem:mod-join-diffw}.
\end{proof}

\subsection{Weighted module-map diagonals}
Consider the weighted module transformation trees complex $\wMTransCx{n}{w}$. There are chain maps
\begin{align*}
 \chi_{i,j,n;v,w}\co\wTreesCx{n+i-j}{w}\otimes\wMTransCx{j-i+1}{v}  &\to\wMTransCx{n}{v+w}  & i>1\\
 \chi_{i,j,n;v,w}\co\wMTreesCx{n+i-j}{w}\otimes\wMTransCx{j-i+1}{v} &\to\wMTransCx{n}{v+w}  & i=1
\end{align*}
given in both cases by $(S,T)\mapsto T\circ_i S$, and a map
\[
 \zeta_{j,n;v,w}\co \wMTransCx{j}{v}\otimes\wMTreesCx{n+1-j}{w}\to\wMTransCx{n}{v+w}
\]
given by $(S,T)\mapsto T\circ_1 S$. We extend these maps to the identity tree in
$\xwTreesCx{1}{0}$ by declaring that the composition of a purple vertex with the
identity tree (in any way) is the original purple vertex, and to stumps in
$\xwTreesCx{0}{0}$ by declaring that the composition of a purple vertex with a
stump is $0$.

\begin{definition}\label{def:w-mod-map-diag}
  Fix a weighted algebra diagonal $\wDiag{*}{*}$ and weighted module diagonals
  $\wModDiag{*}{*}_1$ and $\wModDiag{*}{*}_2$ compatible with $\wDiag{*}{*}$. A
  \emph{weighted module-map diagonal} compatible with $\wModDiag{*}{*}_1$ and $\wModDiag{*}{*}_2$
  is a collection of chain maps
  \[
    \wModMulDiag{n}{w}\co
    \wMTransCx{n}{w}\to \bigoplus_{w_1,w_2\leq w}\wMTransCx{n}{w_1}\otimes\wMTransCx{n}{w_2}\otimes\Ring[Y_1,Y_2]
  \]
  satisfying the following conditions:
  \begin{itemize}
  \item \textbf{Dimension homogeneity}:
    The map $\wModMulDiag{n}{w}$ is
    dimension-preserving, i.e.,
    \[
      \dim(\wModMulDiag{n}{w}(T))=\dim(T).
    \]
  \item \textbf{Weight homogeneity}: The image of $\wModMulDiag{n}{w}$
    is contained in
    $Y_1^{w-w_1}Y_2^{w-w_2}\wMTransCx{n}{w_1}\otimes\wMTransCx{n}{w_2}$. In
    other words,
    \[
      \wgr_1(\wModMulDiag{n}{w}(T))=\wgr_2(\wModMulDiag{n}{w}(T))=\wgr(T)=w.
    \]
  \item \textbf{Compatibility}:
    \[
      \wModMulDiag{n}{v+w}\circ \chi_{i,j,n;v,w}=
      \begin{cases}
        \sum_{\substack{v_1+v_2\leq v\\w_1+w_2\leq w}}(\chi_{1,j,n;v_1,w_1}\otimes \chi_{1,j,n;v_2,w_2})\circ (\wModDiag{j}{v}_1\otimes \wModMulDiag{n+1-j}{w}) & i=1\\
        \sum_{\substack{v_1+v_2\leq v\\w_1+w_2\leq w}}(\chi_{i,j,n;v_1,w_1}\otimes \chi_{i,j,n;v_2,w_2})\circ (\wDiag{j-i+1}{v}\otimes \wModMulDiag{n+i-j}{w}) & i>1
      \end{cases}
    \]
    and
    \[
      \wModMulDiag{n}{v+w}\circ \zeta_{j,n;v,w}=
        \sum_{\substack{v_1+v_2\leq v\\w_1+w_2\leq w}}(\zeta_{j,n;v_1,w_1}\otimes \zeta_{j,n;v_2,w_2})\circ ( \wModMulDiag{n+1-j}{w}\otimes\wModDiag{j}{v}_2).
    \]
  \item \textbf{Non-degeneracy}:
    $\wModMulDiag{1}{0}(\wpcorolla{1}{0})=\wpcorolla{1}{0}\otimes\wpcorolla{1}{0}$.
  \end{itemize}

  Equivalently, we can view $\wModMulDiag{n}{w}$ as a formal linear
  combination of pairs of weighted module transformation trees
  $\wTrModMulDiag{n}{w}$ of dimension $n-1+2w$ (i.e., a \emph{weighted
    module-map tree diagonal})
  \[
    \wTrModMulDiag{n}{w}=\wModMulDiag{n}{w}(\wcorolla{n}{w})\in \bigoplus_{w_1,w_2\leq w}Y_1^{w-w_1}Y_2^{w-w_2}\wMTransCx{n}{w_1}\otimes\wMTransCx{n}{w_2}.
  \] 
  These trees must satisfy:
  \begin{itemize}
  \item \textbf{Compatibility}:
    \[
    \bdy(\wTrModMulDiag{n}{w})=\sum_{w_1+w_2=w}\sum_{k+\ell=n+1}\Bigl(\wTrModMulDiag{k}{w_1}\circ_1\wModDiagCell{\ell}{w_2}_1+\wModDiagCell{k}{w_1}_2\circ_1\wTrModMulDiag{\ell}{w_2} + \sum_{i=2}^{k+1}\wTrModMulDiag{k}{w_1}\circ_i\wDiagCell{\ell}{w_2}\Bigr).
    \]
  \item \textbf{Non-degeneracy}: $\wTrModMulDiag{1}{0}$ is the (unique)
    pair of module transformation trees with one input.
  \end{itemize}
\end{definition}

\begin{figure}
  \centering
  \includegraphics[scale=.66667]{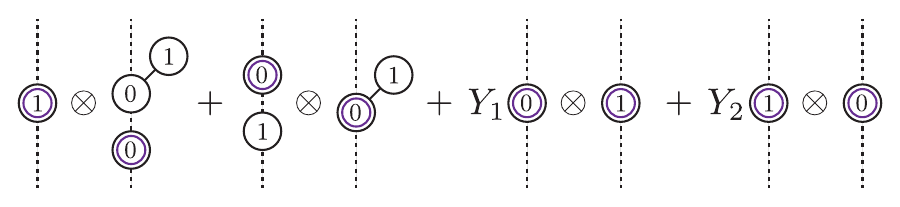}
  \caption[A weighted module-map diagonal]{\textbf{A weighted module-map diagonal.} The figure shows a choice for $\wTrModMulDiag{1}{1}$ compatible with the weighted module diagonal from Figure~\ref{fig:wmod-diag-terms}, and extending the weight-0 case from Figure~\ref{fig:m-tree-diag}.}
  \label{fig:w-mod-map-diag-eg}
\end{figure}

See Figure~\ref{fig:w-mod-map-diag-eg} for  
the first weight-1 term in a particular weighted module-map diagonal.

\begin{definition}\label{def:w-homotopy-module-map-diag}
  Given weighted module-map tree diagonals $\wTrModMulDiag{n}{w}_1$ and
  $\wTrModMulDiag{n}{w}_2$, a \emph{homotopy} from $\wTrModMulDiag{n}{w}_1$ to
  $\wTrModMulDiag{n}{w}_2$ is a collection of chains
  \[
    \eta^{n,w}\in \bigoplus_{w_1,w_2\leq w}Y_1^{w-w_1}Y_2^{w-w_2}\wMTransCx{n}{w_1}\otimes\wMTransCx{n}{w_2},
  \]
  $n\geq 1$, $w\geq 0$, satisfying 
  \[
    \wTrModMulDiag{n}{w}_2-\wTrModMulDiag{n}{w}_1=\bdy(\eta^{n,w})+\sum_{\substack{n_1+n_2=n+1\\w_1+w_2=w}}\Bigl(\eta^{n_1,w_1}\circ_1\wModDiagCell{n_2}{w_2}_1+\wModDiagCell{n_1}{w_1}_2\circ_1\eta^{n_2,w_2} + \sum_{i=2}^{n_1} \eta^{n_1,w_1}\circ_i\wDiagCell{n_2}{w_2}\Bigr)
  \]
  for each $n,w$.
\end{definition}

We have the following analogue of Lemma~\ref{lem:mod-map-diag-exists} and Proposition~\ref{prop:mod-map-diag-homotopic}:
\begin{lemma}\label{lem:w-mod-map-diag-exists-unique}
  Given any weighted algebra diagonal $\wDiagCell{*}{*}$ and weighted
  module diagonals $\wModDiagCell{*}{*}_1$ and $\wModDiagCell{*}{*}_2$
  compatible with $\wDiagCell{*}{*}$ there is a weighted module-map
  diagonal $\wModMulDiag{*}{*}$ compatible with $\wModDiagCell{*}{*}_1$
  and $\wModDiagCell{*}{*}_2$. Further, all weighted module-map diagonals
  compatible with $\wModDiagCell{*}{*}_1$ and $\wModDiagCell{*}{*}_2$ are
  homotopic.
\end{lemma}
\begin{proof}
  The inductive step in the proof is essentially the same as in the proofs of Lemma~\ref{lem:mod-map-diag-exists} and Proposition~\ref{prop:mod-map-diag-homotopic}, using Proposition~\ref{prop:wModTransAcyclic}. Further, the base case is the same as in the unweighted case, since $\wMTransCx{n}{w}$ has trivial homology if $w>0$.
\end{proof}

Note that the maps $\LRjoinW^*$ and $\LRjoinW^{\prime,*}$ from Formulas~\eqref{eq:LRjoinW}
and~\eqref{eq:LRjoinWprime} extend by the same formulas to give maps
\begin{align*}
  \LRjoinW^*\co (\xwMTreesCx{*}{*}\otimes\xwTreesCx{*}{*})^{\otimes \bullet}\otimes (\wMTransCx{*}{*}\otimes \xwTreesCx{*}{*}) \otimes (\xwMTreesCx{*}{*}\otimes\xwTreesCx{*}{*})^{\otimes \bullet}\llbracket Y_1,Y_2\rrbracket
  &\to (\wMTransCx{*}{*}\otimes \xwTreesCx{*}{*})\llbracket Y_1,Y_2\rrbracket\\
  \LRjoinW^{\prime,*}\co (\xwMTreesCx{*}{*}\otimes\xwTreesCx{*}{*})^{\otimes \bullet}\otimes (\wMTransCx{*}{*}\otimes \xwTreesCx{*}{*}) \otimes (\xwMTreesCx{*}{*}\otimes\xwTreesCx{*}{*})^{\otimes \bullet}\llbracket Y_1,Y_2\rrbracket
  &\to (\wMTransCx{*}{*}\otimes\xwMTreesCx{*}{*})\llbracket Y_1,Y_2\rrbracket.
\end{align*}
We can compose $\LRjoinW^{\prime,*}((S_1,T_1),\dots,(S_k,T_k))$ with 
$\IdTree\otimes \wpcorolla{1}{0}$ to obtain a tensor product of
two weighted module transformation trees; see Lemma~\ref{lem:w-mmprim-to-pmm}.

\begin{figure}
  \centering
  %Font is 18 point.
  \includegraphics[scale=.66667]{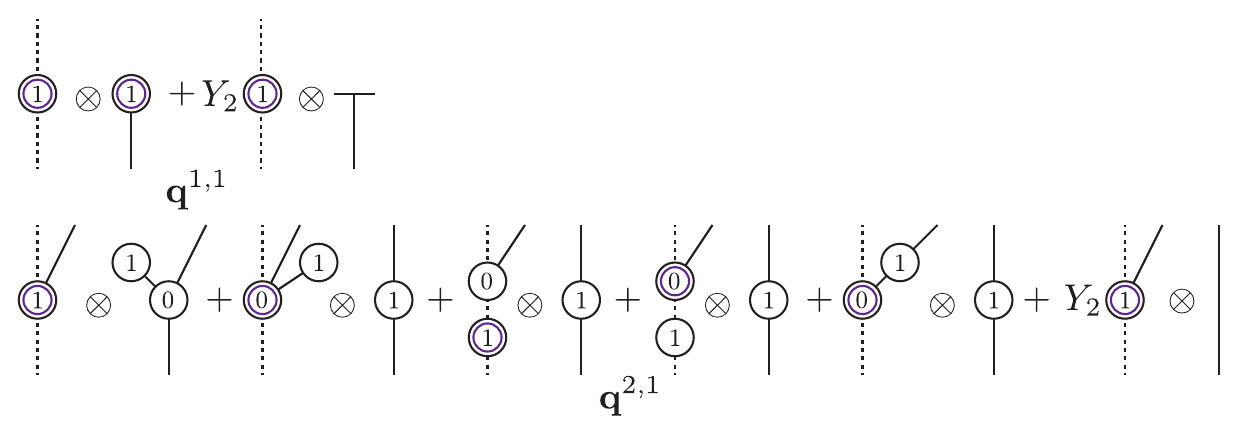}
  \caption[A weighted module-map primitive]{\textbf{A weighted module-map primitive.} The terms $\wTrPMorDiag{1}{1}$ and $\wTrPMorDiag{2}{1}$ are shown; the weight-0 terms are as in Figure~\ref{fig:m-map-prim}. This weighted module-map primitive is compatible with the weighted module primitive in Figure~\ref{fig:wmod-prim-diag-terms} (on both sides).}
  \label{fig:w-map-prim}
\end{figure}

\begin{definition}\label{def:w-mod-map-prim}
  Fix a weighted algebra diagonal $\wDiagCell{*}{*}$ and two weighted module diagonal
  primitives $\wTrPMDiag{*}{*}_1$ and $\wTrPMDiag{*}{*}_2$ compatible with
  $\wDiagCell{*}{*}$. A \emph{weighted module-map primitive} compatible with
  $\wTrPMDiag{*}{*}_1$ and $\wTrPMDiag{*}{*}_2$ consists of chains
  \begin{equation}\label{eq:w-ModMap-prim-1}
    \wTrPMorDiag{n}{w}\in \bigoplus_{w_1,w_2\leq w}
    Y_1^{w-w_1} Y_2^{w-w_2}\bigl(\wMTransCx{n}{w_1}\otimes \xwTreesCx{n}{w_2}\bigr)
  \end{equation}
  of dimension $n+2w-1$, for each $n\geq 1$, $w\geq 0$.  These chains are required
  to satisfy the following conditions:
  \begin{itemize}
  \item \textbf{Compatibility}:
    \begin{equation}\label{eq:w-Mor-prim-compat-long}
      \partial \wTrPMorDiag{n}{w} =
      \hspace{-1.5em}\sum_{\substack{v+w_1+\cdots+w_k=w\\n_1+\cdots+n_k=n+k-1\\i,k}}\hspace{-1.5em}
      \LRjoinW^*(\wTrPMDiag{n_1}{w_1}\otimes\cdots\otimes \wTrPMDiag{n_{i-1}}{w_{i-1}}
      \otimes \wTrPMorDiag{n_i}{w_i}\otimes
      \wTrPMDiag{n_{i+1}}{w_{i+1}}\otimes\cdots\otimes \wTrPMDiag{n_{k}}{w_{k}})+\hspace{-.5em}\sum_{\substack{w_1+w_2=w\\n_1+n_2=n+1}}\hspace{-.5em}\wTrPMorDiag{n_1}{w_1}\circ'\wDiagCell{n_2}{w_2}
    \end{equation}
    or, more succinctly, 
    \begin{equation}
      \label{eq:w-Mor-prim-compat}
      \partial \wTrPMorDiag{*}{*} = 
      \LRjoinW^*((\wTrPMDiag{*}{*}_1)^{\otimes\bullet}
      \otimes {\wTrPMorDiag{*}{*}}\otimes
      (\wTrPMDiag{*}{*}_2)^{\otimes\bullet})+\wTrPMorDiag{*}{*}\circ'\wDiagCell{*}{*}.
    \end{equation}
  \item \textbf{Non-degeneracy}: $\wTrPMorDiag{1}{0}$ is the tensor
    product of the $1$-input, weight $0$ module transformation tree
    (with one internal vertex) and the stump, i.e.,
    \[
      \wTrPMorDiag{1}{0}=\wpcorolla{1}{0}\otimes\stump.
    \]
  \end{itemize}
\end{definition}
The first few terms in a particular weighted module-map primitive are shown in Figure~\ref{fig:w-map-prim}.

\begin{definition}
  Weighted module-map primitives $\wTrPMorDiag{n}{w}_a$ and
  $\wTrPMorDiag{n}{w}_b$ compatible with weighted module diagonal primitives
  $\wTrPMDiag{*}{*}_1$ and $\wTrPMDiag{*}{*}_2$ are \emph{homotopic} if there is
  a collection of pairs of trees
  \[
    \xi^{n,w}\in\bigoplus_{w_1,w_2\leq w}
    Y_1^{w-w_1} Y_2^{w-w_2}\bigl(\wMTransCx{n}{w_1}\otimes \xwTreesCx{n}{w_2}\bigr)
  \]
  so that 
  \[
    \partial \xi^{*,*} + \LRjoinW^*((\wTrPMDiag{*}{*}_1)^{\otimes\bullet}\otimes \xi^{*,*}\otimes (\wTrPMDiag{*}{*}_2)^{\otimes\bullet})+\xi^{*,*}\circ'\wDiagCell{*}{*}
    = \wTrPMorDiag{*}{*}_1-\wTrPMorDiag{*}{*}_2.
  \]
\end{definition}

We have the following analogue of Lemmas~\ref{lem:m-map-prim-exist}
Proposition~\ref{prop:mod-map-prim-homotopic}:
\begin{lemma}\label{lem:w-m-map-prim}
  Given any weighted algebra diagonal $\wDiag{n}{w}$ and module
  diagonal primitives $\wTrPMDiag{n}{w}_1$ and $\wTrPMDiag{n}{w}_2$
  compatible with $\wDiag{n}{w}$ there is a weighted module-map
  primitive $\wTrPMorDiag{n}{w}$ compatible with $\wTrPMDiag{n}{w}_1$
  and $\wTrPMDiag{n}{w}_2$, and all weighted module-map primitives
  compatible with $\wTrPMDiag{n}{w}_1$ and $\wTrPMDiag{n}{w}_2$ are
  homotopic.
\end{lemma}
\begin{proof}
  The proof is similar to the proofs of
  Lemma~\ref{lem:m-map-prim-exist} and
  Proposition~\ref{prop:mod-map-prim-homotopic} and is left to the
  reader.
\end{proof}

As in the unweighted case, it will be useful to describe one
particular class of weighted module-map primitives:
\begin{lemma}\label{lem:w-mor-prim-example}
  Let $\wTrPMDiag{*}{*}$ be a weighted module diagonal primitive. Let
  $\wTrPMorDiag{1}{0}=\wcorolla{1}{0}\otimes\stump$ and for $n\geq 2$ or
  $w\geq 1$, let
  $\wTrPMorDiag{n}{w}$ be the sum over $(S,T)$ in $\wTrPMDiag{n}{w}$ of all
  pairs $(S',T)$ where $S'$ is obtained from $S$ by making one vertex
  on the leftmost strand of $S$ distinguished (purple). Then $\wTrPMorDiag{*}{*}$ is a
  module-map primitive compatible with $\wTrPMDiag{*}{*}$ and $\wTrPMDiag{*}{*}$.
\end{lemma}
\begin{proof}
  This is the same as the proof of Lemma~\ref{lem:mor-prim-example},
  with the amendment that in cases~\ref{item:mor-prim-eg-term-2}
  and~\ref{item:mor-prim-eg-term-3} in that proof, the distinguished
  vertex is 2-valent and has weight $0$.
\end{proof}

Next we discuss the relationship between weighted module-map
primitives and weighted module-map diagonals. As in the unweighted
case, this relationship is through partial module-map diagonals.
\begin{definition}\label{def:w-partial-mod-map-diag}
  Let $\Filt_{\neq(2,0)}\wMTransCx{n}{w}\subset\wMTransCx{n}{w}$
  denote the subspace spanned by trees where the distinguished
  (purple) vertex has either weight $>0$ or valence $>2$.  Let 
  \[
    \gwMTransCx{n}{w}=\wMTransCx{n}{w}/\bigl(\Filt_{\neq(2,0)}\wMTransCx{n}{w}+\bdy (\Filt_{\neq(2,0)}\wMTransCx[*+1]{n}{w})\bigr).
  \]
  
  A \emph{partial weighted module-map diagonal} consists of elements
  \[
    \wPartTrModMulDiag{n}{w}\in\bigoplus_{w_1,w_2\leq w}Y_1^{w-w_1}Y_2^{w-w_2}\wMTransCx{n}{w_1}\otimes\gwMTransCx{n}{w_2}
  \]
  of dimension $n-1+2w$ satisfying
  \begin{itemize}
  \item \textbf{Compatibility}:
    \begin{multline*}
    \bdy(\wPartTrModMulDiag{n}{w})-\sum_{\substack{w_1+w_2=w\\k+\ell=n+1}}\Bigl(\wPartTrModMulDiag{k}{w_1}\circ_1\wModDiagCell{\ell}{w_2}_1+\wModDiagCell{k}{w_1}_2\circ_1\wPartTrModMulDiag{\ell}{w_2}
    +
    \sum_{i=2}^{k+1}\wPartTrModMulDiag{k}{w_1}\circ_i\wDiagCell{\ell}{w_2}\Bigr)=0\\
    \in
    \bigoplus_{w_1,w_2\leq w}Y_1^{w-w_1}Y_2^{w-w_2}\wMTransCx{n}{w_1}\otimes\gwMTransCx{n}{w_2}.
    \end{multline*}
  \item \textbf{Non-degeneracy}: $\wPartTrModMulDiag{1}{0}$ is the (unique)
    pair of module transformation trees with one input and weight $0$.
  \end{itemize}

  Partial weighted module-map diagonals $\wPartTrModMulDiag{*}{*}_1$ and
  $\wPartTrModMulDiag{*}{*}_2$ are \emph{homotopic} if there is a collection
  of elements
  $\zeta^{n,w}\in \bigoplus_{i+j=n-1}\wMTransCx{n}{w}\otimes\Filt_{\leq 2}\wMTransCx{n}{w}$ satisfying
    \begin{multline}\label{eq:w-part-mod-map-htpy}
    \bdy(\zeta^{n,w})-\sum_{\substack{w_1+w_2=w\\k+\ell=n+1}}\Bigl(\zeta^{k,w_1}\circ_1\wModDiagCell{\ell}{w_2}_1+\wModDiagCell{k}{w_1}_2\circ_1\zeta^{\ell,w_2}
    +
    \sum_{i=2}^{k+1}\zeta^{k,w_1}\circ_i\wDiagCell{\ell}{w_2}\Bigr)+\wPartTrModMulDiag{n}{w}_1-\wPartTrModMulDiag{n}{w}_2
    =0\\
    \in\bigoplus_{w_1,w_2\leq w}Y_1^{w-w_1}Y_2^{w-w_2}\wMTransCx{n}{w_1}\otimes\gwMTransCx{n}{w_2}.
    \end{multline}
\end{definition}
Explicitly, $\gwMTransCx{n}{w}$ is spanned by weighted module
transformation trees where the purple vertex is $2$-valent and has
weight $0$, and if two trees differ only in the location of the purple
vertex they are equivalent in $\gwMTransCx{n}{w}$.

As in the unweighted case, given a module-map diagonal $\wTrModMulDiag{*}{*}$, the image
$\wPartTrModMulDiag{*}{*}$ of $\wTrModMulDiag{*}{*}$ under the quotient map
\[
  \wMTransCx{*}{*}\otimes\wMTransCx{*}{*}\to \wMTransCx{*}{*}\otimes
  \gwMTransCx{*}{*}
\]
is a partial weighted module-map diagonal. In this case, we say that
$\wTrModMulDiag{*}{*}$ is a module-map diagonal \emph{extending}
$\wPartTrModMulDiag{*}{*}$.

\begin{lemma}\label{lem:w-mod-trans-subquot}
  The complex $\gwMTransCx{n}{w}$ is contractible if $w=0$ (i.e., has
  homology $\Ring$ in dimension $0$ and $0$ in all other dimensions) and is
  acyclic if $w>0$ (i.e., has trivial homology).
\end{lemma}
\begin{proof}
  Forgetting the purple vertex gives an isomorphism between
  $\gwMTransCx{n}{w}$ and the module trees complex $\wMTreesCx{n}{w}$,
  so the result follows from
  Theorem~\ref{thm:ModAssociaplexAcyclic}.
\end{proof}

\begin{lemma}\label{lem:w-part-mod-maps-htpic}
  All partial weighted module-map diagonals compatible with
  $\wModDiagCell{*}{*}_1$ and $\wModDiagCell{*}{*}_2$ are homotopic.
\end{lemma}
\begin{proof}
  The proof is essentially the same as the proof of
  Lemma~\ref{lem:part-mod-maps-htpic}.  As there, we build the maps
  $\zeta^{n,w}$ inductively, now inducting on $w$ first and then $n$. By
  the non-degeneracy condition, we can take $\zeta^{1,0}=0$. For the
  inductive step one verifies that 
  \[
\sum_{\substack{w_1+w_2=w\\k+\ell=n+1}}\Bigl(\zeta^{k,w_1}\circ_1\wModDiagCell{\ell}{w_2}_1+\wModDiagCell{k}{w_1}_2\circ_1\zeta^{\ell,w_2}
    +
    \sum_{i=2}^{k+1}\zeta^{k,w_1}\circ_i\wDiagCell{\ell}{w_2}\Bigr)+\wPartTrModMulDiag{n}{w}_1-\wPartTrModMulDiag{n}{w}_2    
  \]
  is a cycle in $\bigoplus_{w_1,w_2\leq
    w}Y_1^{w-w_1}Y_2^{w-w_2}\wMTransCx{n}{w_1}\otimes\gwMTransCx{n}{w_2}$
  so, by Lemma~\ref{lem:w-mod-trans-subquot}, is also the boundary of
  some $\zeta^{n,w}$.  
\end{proof}

\begin{lemma}\label{lem:w-mmprim-to-pmm}
  Given an associahedron tree diagonal, compatible module diagonal
  primitives $\wTrPMDiag{*}{*}_1$ and $\wTrPMDiag{*}{*}_2$, and a compatible
  module-map primitive $\wTrPMorDiag{*}{*}$,
  \begin{equation}\label{eq:w-MapMapPrimExponential}
    \wPartTrModMulDiag{*}{*}=
    [\IdTree\otimes \wpcorolla{1}{0}]\circ \LRjoin'((\wTrPMDiag{*}{*}_1)^{\otimes\bullet}\otimes
    \wTrPMorDiag{*}{*}\otimes (\wTrPMDiag{*}{*}_2)^{\otimes\bullet})
  \end{equation}
  is a partial weighted module-map diagonal compatible with
  $\wModDiagCell{*}{*}_1=\LRjoin'((\wTrPMDiag{*}{*}_1)^{\otimes\bullet})$ and
  $\wModDiagCell{*}{*}_2=\LRjoin'((\wTrPMDiag{*}{*}_2)^{\otimes\bullet})$.
\end{lemma}
\begin{proof}
  The proof is the same as the proof of
  Lemma~\ref{lem:mm-prim-partial-mm-diag}.
\end{proof}

\subsection{Weighted \texorpdfstring{$\mathrm{DADD}$}{DADD} diagonals}\label{sec:w-DADD-diags}
In this section, we define the weighted analogue of DADD diagonals
(Section~\ref{sec:DADD-diags}). In this paper, we do not discuss
weighted type \DA\ structures, so this section does not have an
algebraic application in this paper, but we include it for completeness.
\begin{definition}\label{def:w-DADD-diag}
  Fix weighted algebras diagonals $\wDiagCell{*}{*}_1$ and $\wDiagCell{*}{*}_2$.  A \emph{DADD diagonal}
  compatible with $\wDiagCell{*}{*}_1$ and $\wDiagCell{*}{*}_2$ is a collection of elements
  \[
    \wTrDADD{n}{w}\in \bigoplus_{w_1+w_2\leq w}Y_1^{w-w_1}Y_2^{w-w_2}\bigl(\wTransCx{n}{w_1} \otimes \xwTreesCx{n}{w_2}\bigr)
  \]
  in dimension $n-1+2w$, satisfying the following conditions:
  \begin{itemize}
  \item \textbf{Compatibility}:
    \[
      \bdy(\wTrDADD{n}{w})=\sum_{\substack{k+\ell=n+1\\u+v=w}} \wTrDADD{k}{u}\circ\wDiagCell{\ell}{v}_1
      +\sum_{\substack{j_1+\cdots+j_k=n\\u+v_1+\cdots+v_k=w}}\wDiagCell{u}{k}_2\circ(\wTrDADD{j_1}{v_1},\dots,\wTrDADD{j_k}{v_k}).
    \]
  \item \textbf{Non-degeneracy}:
    $\wTrDADD{1}{0}=\wpcorolla{1}{0}\otimes \IdTree$, the tensor product of a
    1-input purple corolla and the identity tree.
  \end{itemize}
\end{definition}

Schematically, the compatibility condition is the same as Equation~\eqref{eq:DADD-diag-compat}, except that the weight is distributed among all the vertices.

\begin{figure}
  \centering
  %Font is 18 point
  \includegraphics[scale=.666667]{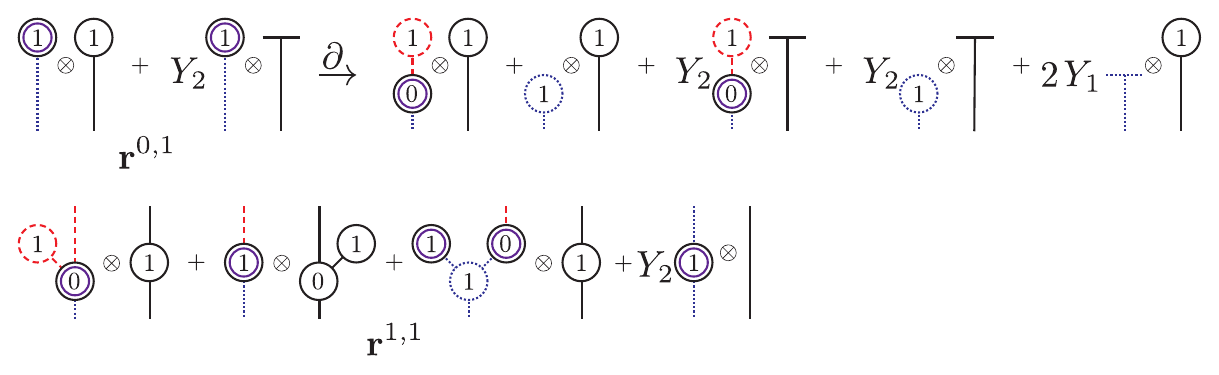}
  \caption[Weighted DADD diagonals]{\textbf{Weighted DADD diagonals.} Top: the element $\wTrDADD{0}{1}$
    and its differential, which is the base case in
    Lemma~\ref{lem:wDADD-exists}. (Recall that we are working in characteristic
    2.) Bottom: one option for the element $\wTrDADD{1}{1}$, compatible with the
    weighted diagonal from Figure~\ref{fig:wdiag-terms}.}
  \label{fig:wDADD}
\end{figure}

\begin{lemma}\label{lem:wDADD-exists}
  Given any weighted algebra diagonals there is a compatible weighted DADD
  diagonal.
\end{lemma}
\begin{proof}
  As in the unweighted case (Lemma~\ref{lem:DADD-exists}), since the
  associahedron and multiplihedron are contractible, it suffices to verify:
  \begin{enumerate}
  \item The right hand side of the compatibility equation is a cycle.
  \item Solutions to the compatibility equation exist when the right
    side is in dimension $0$.
  \end{enumerate}
  The first statement is clear. For the second, we need to check the cases
  $(n,w)\in\{(1,0),(2,0),(0,1)\}$. The first two cases are the same as the
  unweighted case. The third case is shown in Figure~\ref{fig:wDADD} (along with
  the next term in a particular weighted DADD diagonal).
\end{proof}

\subsection{Homotopy unital diagonals}
\begin{definition}\label{def:uwDiag}
  A \emph{homotopy
    unital algebra diagonal} consists of 
  chain maps
  \[
    \uwDiag{n}{w}\co\uwTreesCx{n}{w}\to \bigoplus_{w_1,w_2\leq
      w}\uwTreesCx{n}{w_1}\otimes\uwTreesCx{n}{w_2}\otimes\Ring[Y_1,Y_2]
  \]
  satisfying the dimension homogeneity, weight homogeneity, and
  compatibility under stacking conditions from
  Definition~\ref{def:w-alg-diag} and the non-degeneracy conditions:
  \begin{align*}
    \uwDiag{0}{0}(\stump)&=\stump\otimes\stump\\
    \uwDiag{1}{0}(\IdTree)&=\IdTree\otimes\IdTree\\
    \uwDiag{2}{0}(\wcorolla{2}{0})&=\wcorolla{2}{0}\otimes \wcorolla{2}{0}.
  \end{align*}
  The \emph{seed} of $\uwDiag{n}{w}$ is
  $\uwDiag{0}{1}(\wcorolla{0}{1})$, the image of the weight $1$,
  $0$-input corolla.
  
  Let $\wDiag{n}{w}$ be a weighted algebra diagonal.
  We say that $\uwDiag{n}{w}$ \emph{extends $\wDiag{n}{w}$} if
  $\pi(\uwDiag{n}{w}(\wcorolla{n}{w}))=\wDiag{n}{w}(\wcorolla{n}{w})$
  for each $n,w$ with $n+2w\geq 2$. (The map
  $\pi\co \uwTreesCx{n}{w}\to \wTreesCx{n}{w}$ is from
  Definition~\ref{def:proj-hu-to-u}.) 
  
  Define a \emph{homotopy unital module diagonal} $\uwMDiag{n}{w}$ extending a weighted
  module diagonal similarly, using $\uwMTreesCx{n}{w}$ in place of
  $\uwTreesCx{n}{w}$ and Definition~\ref{def:w-mod-diag} in place of
  Definition~\ref{def:w-alg-diag}, and the non-degeneracy conditions
  \begin{align*}
    \uwDiag{1}{0}(\IdTree)&=\IdTree\otimes\IdTree\\
    \uwDiag{2}{0}(\wcorolla{2}{0})&=\wcorolla{2}{0}\otimes \wcorolla{2}{0}.
  \end{align*}
  We say that $\uwMDiag{n}{w}$ extends $\wModDiag{n}{w}$ if 
  $\pi(\uwMDiag{n}{w}(\wcorolla{n}{w}))=\wModDiag{n}{w}(\wcorolla{n}{w})$
  for each $n,w$ with $n+2w\geq 2$.
  
  Define a \emph{homotopy unital map diagonal} $\uwMulDiag{n}{w}$ using
  $\uwTransCx{n}{w}$ in place of $\uwTreesCx{n}{w}$,
  Definition~\ref{def:w-map-diag} in place of
  Definition~\ref{def:w-alg-diag}, and the non-degeneracy condition
  \begin{align*}
    \uwMulDiag{1}{0}(\wpcorolla{1}{0})=\wpcorolla{1}{0}\otimes \wpcorolla{1}{0}.
  \end{align*}
  The notion of a homotopy unital map diagonal extending an algebra
  map diagonal is defined similarly to the previous cases.

  Define a \emph{homotopy unital module-map diagonal}
  $\uwModMulDiag{n}{w}$ similarly, using $\uwMTransCx{n}{w}$ in place
  of $\uwTreesCx{n}{w}$, Definition~\ref{def:w-mod-map-diag} in place
  of Definition~\ref{def:w-alg-diag}, and the non-degeneracy
  condition
  \begin{align*}
    \uwModMulDiag{1}{0}(\wpcorolla{1}{0})=\wpcorolla{1}{0}\otimes \wpcorolla{1}{0}.
  \end{align*}
  The notion of a homotopy unital module-map diagonal extending a
  module-map diagonal is defined similarly to the previous cases.
\end{definition}

\begin{figure}
  \centering
  \includegraphics[scale=.65]{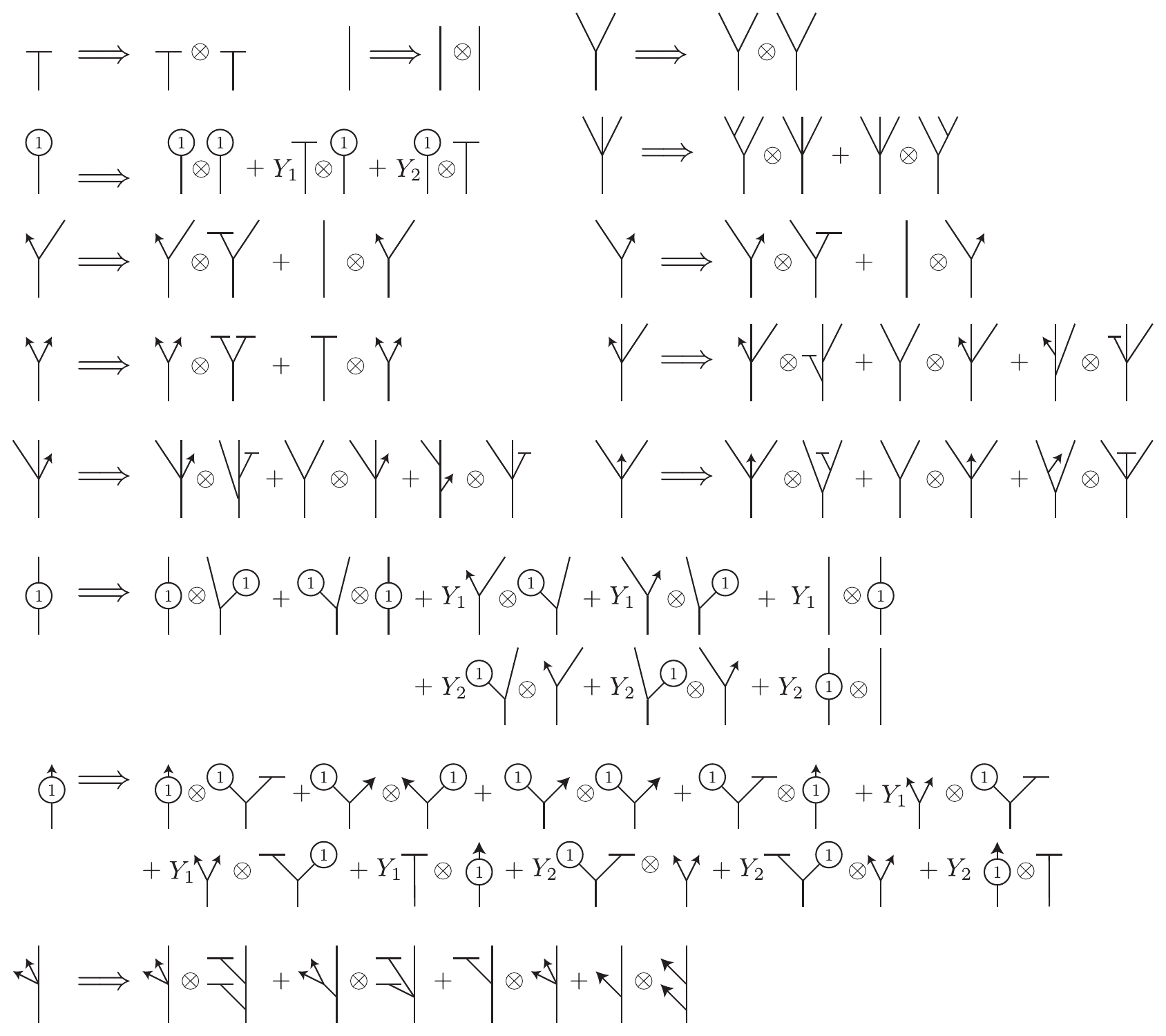}
  \caption[Terms in a homotopy unital algebra diagonal]{\textbf{Terms in a homotopy unital algebra diagonal.} Top
    row: the non-degeneracy conditions. Second row: (maximal) seed for
    the diagonal and base case for the inductive construction.
    Remaining rows: a few other terms in the diagonal. Throughout,
    unlabelled vertices have weight $0$. This diagonal extends the
    diagonal $\wDiag{n}{w}$ from Figures~\ref{fig:diag-cells}
    and~\ref{fig:wdiag-terms}.}
  \label{fig:hu-alg-diag}
\end{figure}

\begin{figure}
  \centering
  \includegraphics[scale=.65]{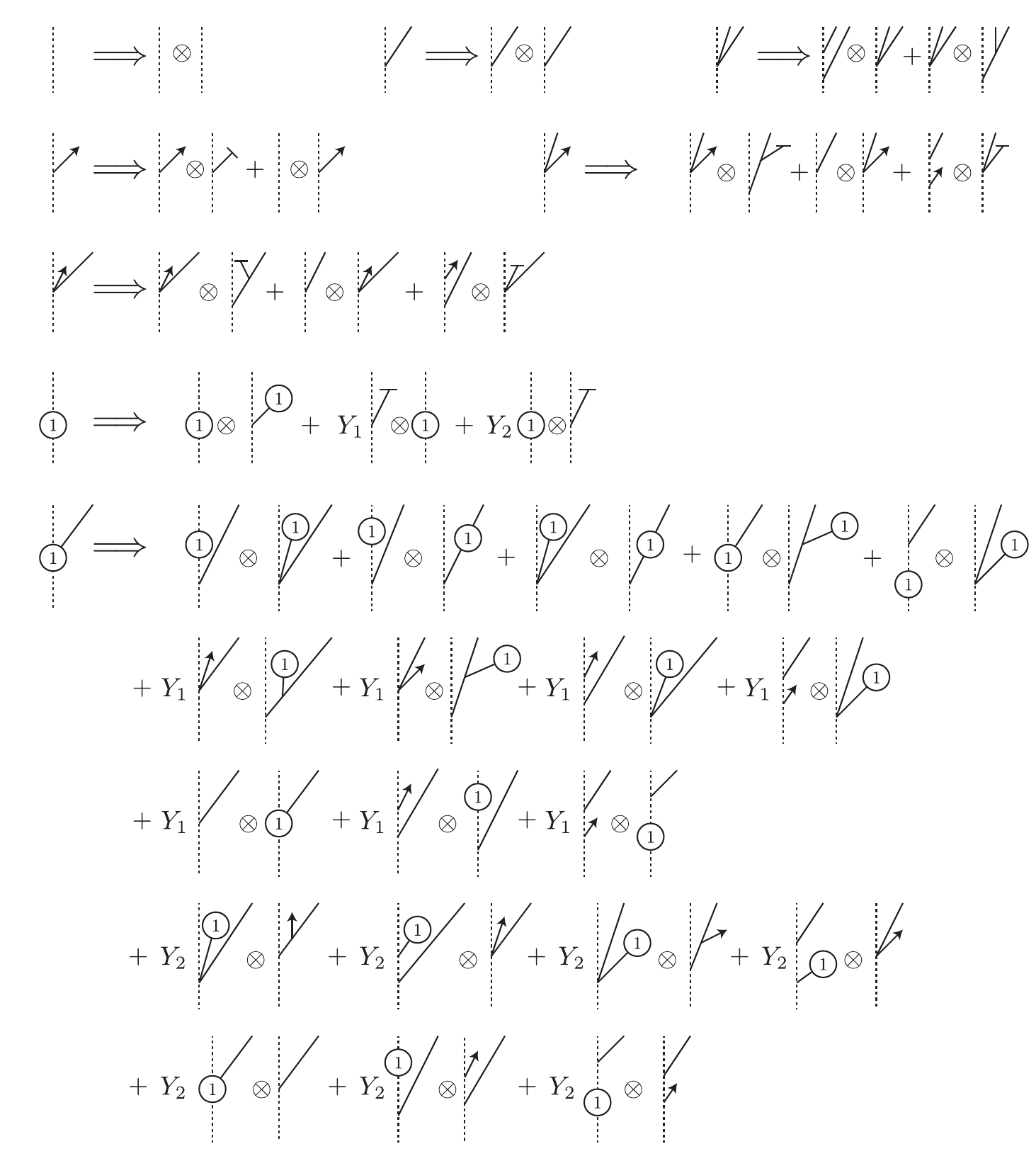}
  \caption[Terms in a homotopy unital module diagonal]{\textbf{Terms in a homotopy unital module diagonal.}
    Unlabelled vertices have weight $0$. This diagonal extends the
    diagonal in Figures~\ref{fig:diag-cells}
    and~\ref{fig:wmod-diag-terms}.}
  \label{fig:hu-mod-diag}
\end{figure}

\begin{figure}
  \centering
  \includegraphics{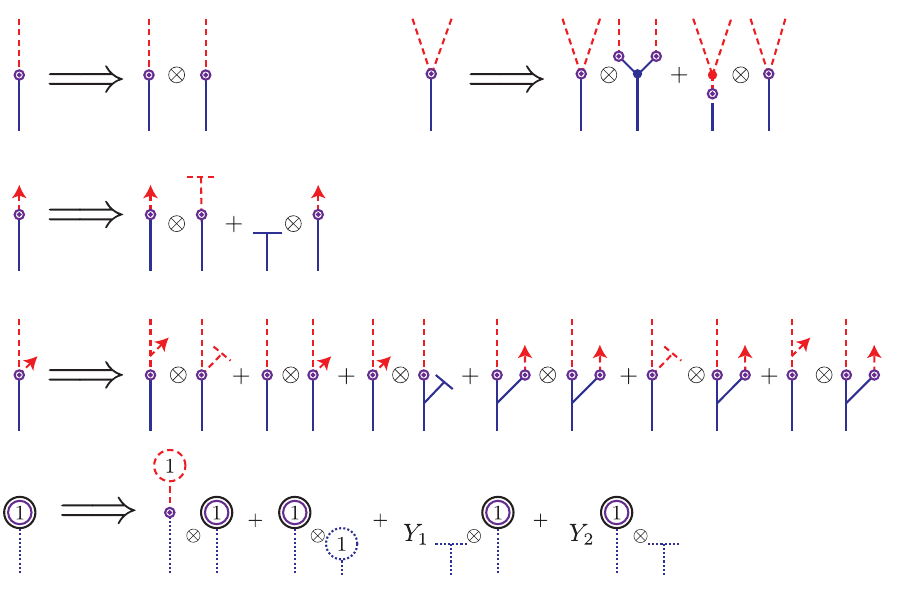}
  \caption[Terms in a homotopy unital map diagonal]{\textbf{Terms in a
      homotopy unital map diagonal.} Vertices not labelled by integers have
    weight $0$. This diagonal extends the diagonal in
    Figures~\ref{fig:mult-diag} and~\ref{fig:w-map-diag-base-case}.}
  \label{fig:hu-map-diag}
\end{figure}

\begin{figure}
  \centering
  \includegraphics{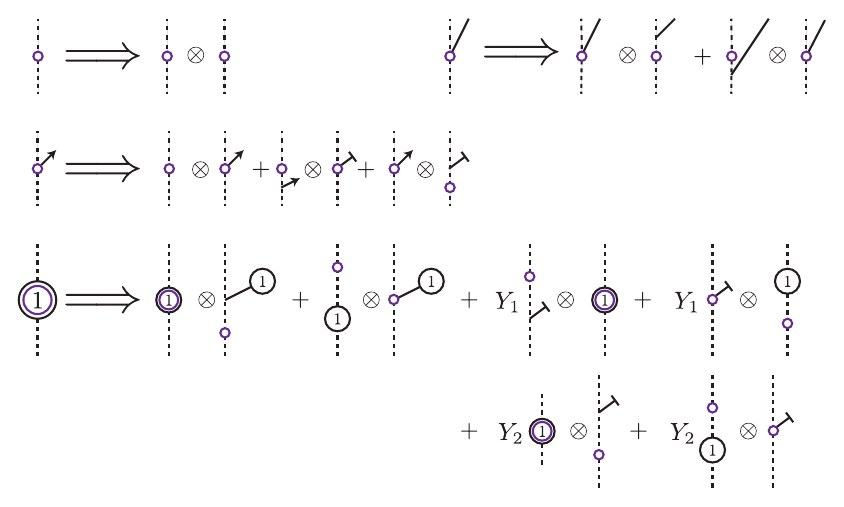}  
  \caption[Terms in a homotopy module-map diagonal]{\textbf{Terms in a homotopy module-map diagonal.} The
    diagonal shown extends the diagonal from
    Figure~\ref{fig:m-tree-diag}. Vertices not labelled by integers have
    weight $0$.}
  \label{fig:hu-mod-map-diag}
\end{figure}

See
Figures~\ref{fig:hu-alg-diag},~\ref{fig:hu-mod-diag},~\ref{fig:hu-map-diag},
and~\ref{fig:hu-mod-map-diag} for some terms in a homotopy unital
algebra diagonal, homotopy unital module diagonal, homotopy unital map
diagonal, and homotopy unital module map diagonal, respectively.

\begin{remark}
  Notice that $\uwDiag{1}{1}(\wcorolla{1}{1})$ in
  Figure~\ref{fig:hu-alg-diag} includes trees with thorns, even though
  $\wcorolla{1}{1}$ has no thorns. In particular, this
  suggests that, unlike the unweighted case, there is no sensible
  tensor product of non-unital weighted algebras (with maximal seed).
\end{remark}

\begin{lemma}\label{lem:hu-diag-exists}
  Given a weighted algebra diagonal $\wDiag{n}{w}$ there is a homotopy
  unital algebra diagonal $\uwDiag{n}{w}$ extending
  $\wDiag{n}{w}$. Similar statements hold for module diagonals, map
  diagonals, and module-map diagonals.
\end{lemma}
\begin{proof}
  To define $\uwDiag{n}{w}$ we need to define the elements
  $\uwDiag{n}{w}(\wcorolla{n_1\uparrow\cdots\uparrow n_k}{w})$.  
  The non-degeneracy condition dictates the values
  $\uwDiag{0}{0}(\stump)$, $\uwDiag{1}{0}(\IdTree)$, and
  $\uwDiag{2}{0}(\wcorolla{2}{0})$.
  Define the image of the weight $1$, $0$-input corolla to be 
  \[
     \uwDiag{0}{1}(\wcorolla{0}{1})=\wDiag{0}{1}(\wcorolla{0}{1}),
  \]
  viewed as a linear combination of pairs of thorn trees in the obvious way.
  
  Next, we define $\uwDiag{n}{w}(\wcorolla{n}{w})$ by
  induction on $(n,w)$, as in the proof of
  Theorem~\ref{thm:wDiag-exists}.  The base cases
  $\uwDiag{0}{0}(\wcorolla{0}{0})$,
  $\uwDiag{1}{0}(\wcorolla{1}{1})$,
  $\uwDiag{2}{0}(\wcorolla{2}{0})$ and
  $\uwDiag{0}{1}(\wcorolla{0}{1})$ have already been
  defined. The base cases $\uwDiag{3}{0}(\wcorolla{3}{0})$ and
  $\uwDiag{1}{1}(\wcorolla{1}{1})$ are shown in
  Figure~\ref{fig:wdiag-terms}. 
  Now, suppose
  $\uwDiag{n}{w}(\wcorolla{n}{w})$ for $w<W$ and for $w=W$
  and $n<N$, satisfying the compatibility condition
  $\bdy \uwDiag{n}{w}(\wcorolla{n}{w})=\uwDiag{n}{w}(\bdy
  \wcorolla{n,\emptyset}{w})$
  and the condition that
  $\pi(\uwDiag{n}{w}(\wcorolla{n}{w}))=\wDiag{n}{w}(\wcorolla{n}{w})$
  for $w<W$ and for $w=W$ and $n<N$. Observe that
  $\uwDiag{N}{W}(\bdy \wcorolla{N}{W})$ is a cycle in
  dimension $N+2W-3>0$, and 
  $\pi(\uwDiag{N}{W}(\bdy \wcorolla{N}{W}))=\wDiag{N}{W}(\bdy \wcorolla{N}{W})$.
  Since $\pi$ is surjective, we can choose an element %
  $x\in \bigoplus_{w_1,w_2\leq W}\uwTreesCx{N}{w_1}\otimes\uwTreesCx{N}{w_2}\otimes\Ring[Y_1,Y_2]$ %
  with $\pi(x)=\wDiag{N}{W}(\wcorolla{N}{W})$. Thus,
  \[
    \pi(\bdy x + \uwDiag{N}{W}(\bdy \wcorolla{N}{W}))=0
  \]
  so, since $\pi$ induces an isomorphism on homology,
  $\bdy x + \uwDiag{N}{W}(\bdy \wcorolla{N}{W})=\bdy y$ where
  $\pi(y)=0$. Let $\uwDiag{N}{W}(\wcorolla{N}{W})=x+y$ and
  continue the induction.
  
  Finally, we perform induction on the number $k-1$ of thorns to define
  $\uwDiag{n}{w}(\wcorolla{n_1\uparrow\cdots\uparrow n_k}{w})$ for $k>1$. 
  Suppose we have
  defined $\uwDiag{n}{w}(\wcorolla{n_1\uparrow\cdots\uparrow n_k}{w})$ for $k<K$ and for $k=K$ if
  $w<W$ or $w=W$ and $n<N$.
  Then we have already defined
  $\uwDiag{n}{w}(\bdy \wcorolla{n_1\uparrow\cdots\uparrow n_k}{w})$.  It is clear that
  $\uwDiag{n}{w}(\bdy \wcorolla{n_1\uparrow\cdots\uparrow n_k}{w})$ is a cycle. There are several cases:
  \begin{itemize}
    \item $N=2$, $W=0$, and $K=2$.
  \item $N=0$, $W=1$, and $K=2$.
  \item $N+2W=0$ and $K\geq 3$.
  \item $N+2W>3$.
  \end{itemize}
  In the first two cases, by inspection the class
  $\uwDiag{n}{w}(\bdy \wcorolla{n_1\uparrow\cdots\uparrow n_K}{w})$ is a boundary; see Figure~\ref{fig:hu-alg-diag}. In the last two cases, it
  follows from Theorem~\ref{thm:hu-contractible} that
  $\uwDiag{n}{w}(\bdy \wcorolla{n_1\uparrow\cdots\uparrow n_K}{w})$ is a boundary. So, in any case, choose
  $\uwDiag{n}{w}(\wcorolla{n_1\uparrow\cdots\uparrow n_K}{w})$ so that
  $\bdy\uwDiag{n}{w}(\wcorolla{n_1\uparrow\cdots\uparrow
    n_K}{w})=\uwDiag{n}{w}(\bdy\wcorolla{n_1\uparrow\cdots\uparrow n_K}{w})$
  and continue the induction.
  
  This completes the construction of $\uwDiag{n}{w}$.  The
  constructions of module diagonals, map diagonals, and module-map
  diagonals are similar. See
  Figures~\ref{fig:hu-mod-diag},~\ref{fig:hu-map-diag},
  and~\ref{fig:hu-mod-map-diag} for some low-dimensional cases
  (including the base cases for the inductions).
\end{proof}

\begin{lemma}\label{lem:hu-mod-map-diag-unique}
  Given a homotopy unital algebra diagonal $\uwDiag{*}{*}$ and two homotopy
  unital module diagonals $\uwMDiag{*}{*}_1$ and $\uwMDiag{*}{*}_2$ compatible
  with $\uwDiag{*}{*}$, all homotopy unital module-map diagonals with
  $\uwMDiag{*}{*}_1$ and $\uwMDiag{*}{*}_2$ are homotopic.
\end{lemma}
\begin{proof}
  The proof is the usual inductive argument, using
  Theorem~\ref{thm:hu-mod-trans-acyclic}, and is left to the reader.
\end{proof}

%%% Local Variables: 
%%% mode: latex
%%% TeX-master: "AbstractDiagonal.tex"
%%% End: 

\section{Algebraic applications of weighted diagonals}\label{sec:wDiagApps}
Weighted diagonals can be used to construct various tensor products.

\begin{convention}\label{conv:wDiagApps-ground}
  As in Sections~\ref{sec:algebra} and~\ref{sec:wAinfty}
  (Convention~\ref{conv:Ring}), fix a
  commutative $\FF_2$-algebra $\Ring$ and commutative $\Ring$-algebras
  $\Ground_1$ and $\Ground_2$ (and occasionally
  $\Ground_3$). Throughout, we assume that the action of $\Ring$ on
  bimodules is central, i.e., if $r\in R$ and $a$ is an element of an
  $\Ring$-bimodule then $ra=ar$.

  Fix also elements $Y_i\in \Ground_i$. The elements $Y_i$ make
  $\Ground_i$ into a $\Ring[Y_i]$-algebra. We assume that the elements
  $Y_i$ also act centrally on $\Ground_i$-bimodules, i.e., if $a$ is
  an element of a $\Ground_i$-bimodule then $Y_ia=aY_i$.
  
  We will typically let $\Ground=\Ground_1\otimes_\Ring\Ground_2$.
  Undecorated tensor products are over the appropriate ring
  $\Ground_i$ or $\Ground$.
\end{convention}
The convention that the $Y_i$ acts centrally implies that for a
$w$-algebra $\wAlg$ over $\Ground_1$, say,
\[
  \mu_n^w(a_1,\dots,Y_1a_k,\dots,a_n)=Y_1\mu_n^w(a_1,\dots,a_k,\dots,a_n),
\]
by repeatedly using the relation $a_jY_1=Y_1a_j$ and multi-linearity
of the operation over $\Ground_1$. It does not imply that
\[
  \mu_2^0(a,Y_1)=\mu_2^0(Y_1,a)
\]
unless $\wAlg$ is weakly unital.

\subsection{Tensor products of weighted algebras}\label{sec:w-alg-tens}
Fix ground rings $\Ground_1$ and $\Ground_2$ and let
$\Ground=\Ground_1\rotimes{\Ring}\Ground_2$. Fix also integers $\kappa_1$,
$\kappa_2$. Assume that the $Y_i\in\Ground_i$ satisfy
\[
  \gr(Y_i)=\kappa_i-2.
\]
(The grading $\gr(Y_i)$ should not be confused with the dimension
$\dim(Y_i)=0$ or the weight $\wgr_j(Y_i)=\delta_{i,j}$ from
Section~\ref{sec:wDiags}.)

Fix strictly unital weighted $\Ainf$-algebras $\wAlg=(A,\{\mu_n^w\})$ and
$\wBlg=(B,\{\nu_n^w\})$ over $\Ground_1$ and $\Ground_2$,
with weight gradings $\kappa_1$ and $\kappa_2$, respectively.
Given a weighted algebra diagonal $\wDiag{*}{*}$, we will associate a weighted
$\Ainf$-algebra $\wAlg\otimes_{\wDiagNS}\wBlg$ over $\Ground$ with
weight grading $\kappa_1+\kappa_2-2$.

By Lemma~\ref{lem:wAlg-reinterp}, the operations on $\wAlg$
(respectively $\wBlg$) induce chain maps
$\mu\co \wTreesCx{n}{w}\to \Mor(A^{\kotimes{\Ground_1} n},A\grs{(2-\kappa_1)w})$ (respectively
$\nu\co \wTreesCx{n}{w}\to \Mor(B^{\kotimes{\Ground_2} n},B\grs{(2-\kappa_2)w})$). We can
extend these maps to $\xwTreesCx{n}{w}$ by declaring that
$\mu(\IdTree)=\Id\co A\to A$ and $\mu(\stump)=1\in A$ (and
similarly for $\nu$). 

\begin{lemma}\label{lem:w-extended-mu}
  If the algebra $\wAlg$ is strictly unital then the extension of
  $\mu$ to $\xwTreesCx{n}{w}$ is a chain map and intertwines stacking
  of trees and composition of operations.
\end{lemma}
\begin{proof}
  The fact that the extension is a chain map follows from the fact
  that $d(1)=0$. 
  The fact that the extension to $\IdTree$ is compatible under stacking
  is trivial. The fact that the extension to $\stump$ is compatible
  with stacking follows from the fact that $\mu_2^w(\unit,a)=\mu_2^w(a,\unit)=a$
  and $\mu_n^w(a_1,\dots,a_i,\unit,a_{i+2},\dots,a_n)=0$ if $(n,w)\neq(2,0)$.
\end{proof}

\begin{definition}\label{def:wAlg-tp}
  Fix a weighted algebra diagonal
  \[    
    \wDiag{n}{w}\co \wTreesCx{n}{w}\to \bigoplus_{w_1+w_2\leq
      2w}\xwTreesCx{n}{w_1}\otimes_\Ring\xwTreesCx{n}{w_2}\otimes_\Ring\Ring[Y_1,Y_2]. 
  \]
  Fix also strictly unital weighted $\Ainf$-algebras $\wAlg$ and
  $\wBlg$ over $\Ground_1$ and $\Ground_2$, respectively.  Define maps
  \begin{equation}\label{eq:wAlg-tensor}
    (\mu\otimes_{\wDiagNS}\nu)\co \wTreesCx{n}{w}\to
    \Mor((A\rotimes{\Ring} B)^{\kotimes{\Ground} n},A\rotimes{\Ring} B\grs{(4-\kappa_1-\kappa_2)w})
  \end{equation}
  to be the composition $(\mu\otimes\nu)\circ\wDiag{*}{*}$.
  Since
  $\wDiag{*}{*}$ is a collection of chain maps and is compatible under stacking, by Lemma~\ref{lem:wAlg-reinterp} the operations induce
  a weighted $\Ainf$-algebra structure on
  $\wAlg\wADtp\wBlg$ over $\Ground$.
  The induced weighted algebra is called the {\em weighted tensor
    product via the diagonal $\wDiag{*}{*}$} of $\wAlg$ and $\wBlg$.
  It is denoted $\wAlg\wADtp\wBlg$.
\end{definition}

\begin{lemma}\label{lem:wADtp-grading}
  The map $\mu\otimes_{\wDiagNS}\nu$ in Definition~\ref{def:wAlg-tp}
  does, in fact, shift the grading by $(4-\kappa_1-\kappa_2)w$, so
  $\wAlg\wADtp\wBlg$ is a $w$-algebra with weight grading
  $\kappa_{12}\coloneqq \kappa_1+\kappa_2-2$. Further, if $\wAlg$ and
  $\wBlg$ are bonsai, so is $\wAlg\wADtp\wBlg$.
\end{lemma}
\begin{proof}
  Define three new gradings on $\wTreesCx{n}{w}$ by:
  \begin{align*}
    \gr_1(T)&=\dim(T)+(\kappa_1-2)w\\
    \gr_2(T)&=\dim(T)+(\kappa_2-2)w\\
    \gr_{12}(T)&=\dim(T)+(\kappa_1+\kappa_2-4)w=\dim(T)+(\kappa_{12}-2)w.
  \end{align*}
  If we grade $\wTreesCx{n}{w}$ by $\gr_1$ then
  $\mu\co \wTreesCx{n}{w}\to \Mor(A^{\kotimes{\Ground_1} n},A)$ is
  grading-preserving.  Similarly, grading $\wTreesCx{n}{w}$ by
  $\gr_2$, $\nu$ is grading-preserving. We must check that, with
  respect to $\gr_{12}$, $(\mu\otimes_{\wDiagNS}\nu)$ is
  grading-preserving.

  Since $\mu\otimes_{\wDiagNS}\nu$ is grading-preserving with respect
  to the grading $\gr_1\otimes \gr_2$ on
  $\xwTreesCx{*}{*}\rotimes{\Ring}\xwTreesCx{*}{*}$, it suffices to show that
  $\wDiag{*}{*}$ is grading-preserving with respect to the grading
  $\gr_{12}$ on the source and $\gr_1\otimes\gr_2$ on the target.
  This is straightforward: given a term $Y_1^aY_2^b(T_1,T_2)$ in
  $\wDiag{n}{w}(S)$ we have
  \begin{align*}
    \gr_{12}(\wDiag{n}{w}(S))&=\dim(\wDiag{n}{w}(S))+(\kappa_1+\kappa_2-4)w\\
                             &=\dim(T_1)+\dim(T_2)+(\kappa_1-2)\wgr_1(Y_1^aT_1)+(\kappa_2-2)\wgr_2(Y_2^bT_2)\\
                             &=\gr_1(Y_1^aT_1)+\gr_2(Y_2^bT_2),
  \end{align*}
  as desired.

  The statement about bonsai-ness follows from dimension homogeneity
  of the diagonal.
\end{proof}

The weight $1$, $0$-input operation on $\wAlg\wADtp\wBlg$ depends on
the seed $\wSeed$. A term $\wcorolla{0}{1}\otimes\wcorolla{0}{1}$ in
$\wSeed$ contributes $\mu_0^1\otimes\mu_0^1$ to $(\mu\wADtp\nu)_0^1$,
and a term $Y_1\stump\otimes\wcorolla{0}{1}$ (respectively
$Y_2\wcorolla{0}{1}\otimes \stump$) contributes a copy of
$Y_1\otimes \nu_0^1$ (respectively $\mu_0^1\otimes Y_2$).

\begin{definition}\label{def:w-alg-map-tensor}
  Fix:
  \begin{itemize}
  \item Weighted algebra diagonals $\wDiag{*}{*}_1$ and
    $\wDiag{*}{*}_2$,
  \item A weighted map diagonal $\wMulDiag{*}{*}$ compatible with
    $\wDiag{*}{*}_1$ and $\wDiag{*}{*}_2$,
  \item Strictly unital weighted $\Ainf$-algebras $\wAlg_1$, $\wAlg_2$
    over $\Ground_1$ with weight gradings $\kappa_1$,
  \item Strictly unital weighted $\Ainf$-algebras $\wBlg_1$, $\wBlg_2$
    over $\Ground_2$ with weight gradings $\kappa_2$, and
  \item Strictly unital weighted algebra homomorphisms
    $f\co\wAlg_1\to\wAlg_2$ and $g\co\wBlg_1\to\wBlg_2$.
  \end{itemize}
  The composition
  \begin{align*}
    \wTransCx{n}{w}&\xrightarrow{\wMulDiag{n}{w}}
                     \bigoplus_{w_1,w_2\leq w} \Ring[Y_1,Y_2]\rotimes{\Ring} \xwTransCx{n}{w_1}\rotimes{\Ring}
                     \xwTransCx{n}{w_2}\\
                   &\xrightarrow{f\otimes g}
                     \Mor(A_1^{\kotimes{\Ground_1} n},A_2\grs{(2-\kappa_1)w})\rotimes{\Ring}\Mor(B_1^{\kotimes{\Ground_2}
                     n},B_2\grs{(2-\kappa_2)w})\\
    &\hookrightarrow\Mor((A_1\rotimes{\Ring} B_1)^{\kotimes{\Ground} n},A_2\rotimes{\Ring} B_2\grs{(4-\kappa_1-\kappa_2)w})
  \end{align*}
  specifies a weighted algebra homomorphism $f\wADtp[\wMulDiagNS]g$, the
  \emph{tensor product via the diagonal $\wMulDiag{*}{*}$} of $f$ and $g$.
\end{definition}
(The proof that the grading shifts are as specified is similar to the
proof of Lemma~\ref{lem:wADtp-grading}, and it is clear that if $f$
and $g$ are bonsai then so is $f\wADtp[\wMulDiagNS]g$.)

\begin{proof}[Proof of Theorem~\ref{thm:wAlgebras}]
  Existence of weighted algebra diagonals is Theorem~\ref{thm:wDiag-exists}.
  Part~\ref{item:wAlg-thm-Ainf} is clear from the definition.  For
  Parts~\ref{item:wAlg-thm-qi} and~\ref{item:wAlg-thm-change-diag}, we
  use weighted map diagonals, whose existence is ensured by
  Lemma~\ref{lem:wmul-diag-exists}.  With these remarks in place, the
  proof of Theorem~\ref{thm:Algebras} applies with little change to
  establish the weighted version.
\end{proof}

\begin{warning}
  As in the unweighted case, even though we are assuming $\wAlg$ and
  $\wBlg$ are strictly unital, $\wAlg\wADtp\wBlg$ is not necessarily
  strictly unital. See Sections~\ref{sec:w-htpy-unital} for further
  discussion.
\end{warning}

\begin{remark}\label{rem:Amorim}
  By using homological perturbation theory, Amorim constructed tensor products
  of gapped, filtered $\Ainf$-algebras, which includes weighted $\Ainf$-algebras
  as a special case~\cite{Amorim16:tensor}. In our language, his construction
  does not include the term $\mu_0^1\otimes\mu_0^1$ in the seed. It seems
  reasonable to expect that the tensor product he constructs is quasi-isomorphic
  to the tensor product arising from any weighted algebra diagonal with seed
  $\stump\otimes\wcorolla{0}{1}+\wcorolla{0}{1}\otimes \stump$, but we have not
  verified this. His techniques do not seem to adapt easily to seeds including
  the term $\mu_0^1\otimes\mu_0^1$; perhaps correspondingly, note that many of
  the terms in Figure~\ref{fig:wdiag-terms} do not arise for seeds without this
  term. (In particular, with the terms coming from $\mu_0^1\otimes\mu_0^1$ there
  is not a naive tensor product of what Amorim calls filtered \dg algebras.)
\end{remark}

\subsection{Tensor products of weighted modules}
Let $\wAlg=(A,\{\mu_n^w\})$ be a strictly unital $w$-algebra over $\Ground_1$, say,
with weight grading $\kappa$ and let
$\wMod_{\wAlg}=(M,\{m_n^w\})$ be a strictly unital $w$-module over $\wAlg$.
Lemma~\ref{lem:wMod-reinterp} gives a corresponding map
\[ 
  m\co \wMTreesCx{1+n}{w}\to \Mor(M\kotimes{\Ground_1} A^{\kotimes{\Ground_1} n},M\grs{(2-\kappa)w}) 
\]
These maps can be extended to 
\[ 
  m\co \xwMTreesCx{1+n}{w}\to \Mor(M\kotimes{\Ground_1} A^{\kotimes{\Ground_1} n},M\grs{(2-\kappa)w}) 
\]
by the convention that $m(\IdTree)=\Id\co M\to M$. As in
Lemma~\ref{lem:w-extended-mu}, it follows from strict unitality of
$\wMod$ and $\wAlg$ that this extension is compatible with stacking.

\begin{definition}\label{def:wMod-tp}
  Let $\wAlg$ and $\wBlg$ be strictly unital $w$-algebras with weight
  gradings $\kappa_1$ and $\kappa_2$, respectively, and let
  $\wMod_{\wAlg}$ and $\wNod_{\wBlg}$ be strictly unital $w$-modules
  over $\wAlg$ and $\wBlg$, respectively. Fix a weighted algebra
  diagonal $\wDiag{*}{*}$ and a weighted module diagonal
  $\wModDiag{*}{*}$ compatible with $\wDiag{*}{*}$. Define chain maps
  \begin{align*}
    \wMTreesCx{n+1}{w}&\xrightarrow{\wModDiag{n+1}{w}} \bigoplus_{w_1,w_2\leq w}\Ring[Y_1,Y_2]\rotimes{\Ring} \xwMTreesCx{n+1}{w_1}\rotimes{\Ring} \xwMTreesCx{n+1}{w_2}\\
    &\xrightarrow{m^M\otimes m^{N}} \Mor(M\kotimes{\Ground_1} A^{\kotimes{\Ground_1} n},M\grs{(2-\kappa_1)w})\rotimes{\Ring}
    \Mor(N\kotimes{\Ground_2} B^{\kotimes{\Ground_2} n},N\grs{(2-\kappa_2)w})\\
    &\hookrightarrow \Mor((M\rotimes{\Ring} N)\kotimes{\Ground}(A\rotimes{\Ring} B)^{\kotimes{\Ground} n},M\rotimes{\Ring} N\grs{(4-\kappa_1-\kappa_2)w}).
  \end{align*}
  By Lemma~\ref{lem:wMod-reinterp}, these operations endow
  $M\rotimes{\Ring} N$ with the structure of a 
  $w$-module, the {\em tensor product} of $\wMod_{\wBlg}$ and
  $\wNod_{\wAlg}$. We will denote this tensor product
  $\wMod_{\wAlg}\wMDtp \wNod_{\wBlg}$.
\end{definition}

\begin{lemma}\label{lem:wMDtp-grading}
  The grading shifts in Definition~\ref{def:wMod-tp} are as stated, so
  Definition~\ref{def:wMod-tp} does define a $w$-module over
  $\wAlg\wADtp\wBlg$. Further, if $\wAlg$ and $\wBlg$ are bonsai, so
  is $\wAlg\wADtp\wBlg$.
\end{lemma}
\begin{proof}
  The proof is the same as the proof of Lemma~\ref{lem:wADtp-grading}.
\end{proof}

Next we define tensor products of morphisms of $w$-modules.

\begin{definition}\label{def:w-tensor-mod-maps}
  Fix:
  \begin{itemize}
  \item A weighted algebra diagonal $\wDiag{*}{*}$,
  \item Weighted module diagonals $\wModDiag{*}{*}_1$ and
    $\wModDiag{*}{*}_2$ compatible with $\wDiag{*}{*}$,
  \item A weighted module-map diagonal $\wModMulDiag{*}{*}$ compatible
    with $\wModDiag{*}{*}_1$ and $\wModDiag{*}{*}_2$,
  \item a strictly unital $w$-algebra $\wAlg$ with weight grading $\kappa_1$,
  \item a strictly unital $w$-algebra $\wBlg$ with weight grading $\kappa_2$,
  \item strictly unital $w$-modules $\wMod^1$ and $\wMod^2$ over $\wAlg$,
  \item strictly unital $w$-modules $\wNod^1$ and $\wNod^2$ over
    $\wBlg$,
  \item and morphisms $f\in\Mor(\wMod^1, \wMod^2)$ and
    $g\in\Mor(\wNod^1, \wNod^2)$ in the \dg category of strictly
    unital $w$-modules.
  \end{itemize}
  Then there is a corresponding
  morphism $(f\MDtp[\wModMulDiagNS] g)\in \Mor(\wMod^1\MDtp[{\wModDiagNS_1}]
  \wNod^1,\wMod^2\MDtp[{\wModDiagNS_2}] \wNod^2)$ so that
  $F_{f\MDtp[\ModMulDiag] g}$ is the composition
  \begin{equation}\label{eq:w-define-mod-mor-tp}
    \begin{split} \wMTransCx{n}{w}&\stackrel{\ModMulDiag}{\longrightarrow}\bigoplus_{w_1+w_2\leq w}\Ring[Y_1,Y_2]\rotimes{\Ring} \wMTransCx{n}{w_1}\rotimes{\Ring}\wMTransCx{n}{w_2}\\
      &\xrightarrow{F_{f}\otimes F_{g}}\Mor(M^1\kotimes{\Ground_1} A^{\kotimes{\Ground_1} n},M^2\grs{(2-\kappa_1)w})\rotimes{\Ring} \Mor(N^1\kotimes{\Ground_2} B^{\kotimes{\Ground_2} n},N^2\grs{(2-\kappa_2)w})\\
      &\hookrightarrow \Mor\bigl((M^1\rotimes{\Ring} N^1)\kotimes{\Ground}(A\rotimes{\Ring} B)^{\kotimes{\Ground} n},M^2\rotimes{\Ring} N^2\grs{(2-\kappa)w}\bigr).
    \end{split}
  \end{equation}
\end{definition}

(As usual, the proof that the grading shifts are as specified is
similar to the proof of Lemma~\ref{lem:wADtp-grading}, and it is clear
that this operation takes a pair of bonsai morphisms to a bonsai
morphism.)

We have the following analogue of Proposition~\ref{prop:dg-bifunctor}:
\begin{proposition}\label{prop:w-dg-bifunctor}
  Fix a weighted algebra diagonal $\wDiag{*}{*}$, weighted module diagonals
  $\wModDiag{*}{*}_1$, $\wModDiag{*}{*}_2$, and $\wModDiag{*}{*}_3$
  compatible with $\wDiag{*}{*}$, and for $1 \le i < j \le 3$ weighted
  module-map diagonals $\wModMulDiag{*}{*}_{ij}$ compatible with
  $\wModDiag{*}{*}_i$ and $\wModDiag{*}{*}_j$. Then for any $w$-algebras $\wAlg$
  and $\wBlg$, $w$-modules $\wMod^1$, $\wMod^2$, and $\wMod^3$ over
  $\wAlg$, $w$-modules $\wNod^1$, $\wNod^2$, and $\wNod^3$ over
  $\wBlg$, and morphisms
  \[
    \wMod^1\stackrel{f_1}{\longrightarrow}\wMod^2\stackrel{f_2}{\longrightarrow}\wMod^3\qquad\qquad
    \wNod^1\stackrel{g_1}{\longrightarrow}\wNod^2\stackrel{g_2}{\longrightarrow}\wNod^3,
  \]
  we have
  \begin{align}
    (f_2\circ f_1)\MDtp[\wModMulDiagNS_{13}](g_2\circ g_1)&\sim(f_2\MDtp[\wModMulDiagNS_{23}] g_2)\circ (f_1\MDtp[\wModMulDiagNS_{12}] g_1)
       \label{eq:w-mod-map-tensor-comp}\\
    d(f_1\MDtp[\wModMulDiagNS_{12}]g_1)&=(df_1)\MDtp[\wModMulDiagNS_{12}]g_1+f_1\MDtp[\wModMulDiagNS_{12}](dg_1)
      \label{eq:w-mod-map-tensor-diff}
  \end{align}
  where in \eqref{eq:w-mod-map-tensor-comp} the notation $\sim$
  indicates that the two sides are homotopic morphisms.
\end{proposition}
\begin{proof}
  The proof is similar to the proof of Lemma~\ref{lem:tens-mod-maps} and
  Proposition~\ref{prop:dg-bifunctor}, and is left to the reader.
\end{proof}

\begin{lemma}\label{lem:w-mod-id-tens-id}
  Given $w$-modules $\wMod$ and $\wNod$, weighted module diagonals
  $\wModDiag{*}{*}_1$ and $\wModDiag{*}{*}_2$, and a weighted module-map diagonal $\wModMulDiag{*}{*}$
  compatible with $\wModDiag{*}{*}_1$ and $\wModDiag{*}{*}_2$, the map
  \[
    \Id_{\wMod}\MDtp[\wModMulDiagNS]\Id_{\wNod}\co \wMod\MDtp[\wModDiagNS_1]\wNod\to \wMod\MDtp[\wModDiagNS_2]\wNod
  \]
  is an isomorphism.
\end{lemma}
\begin{proof}
  The fact that $(\Id_{\wMod}\MDtp[\ModMulDiag]\Id_{\wNod})_1^0$ is an
  isomorphism is immediate from the non-degeneracy condition for
  $\wModMulDiag{*}{*}$, so this follows from Lemma~\ref{lem:w-mod-iso}.
\end{proof}

\begin{corollary}\label{cor:w-mod-id-tens-id-is-sometimes-id}
  Given $w$-modules $\wMod$ and $\wNod$, a weighted module diagonal
  $\wModDiag{*}{*}$, and a weighted module-map diagonal $\wModMulDiag{*}{*}$
  compatible with $\wModDiag{*}{*}$ and $\wModDiag{*}{*}$, the map
  \[
    \Id_{\wMod}\MDtp[\wModMulDiagNS]\Id_{\wNod}\co \wMod\MDtp[\wModDiagNS]\wNod\to \wMod\MDtp[\wModDiagNS]\wNod
  \]
  is homotopic to the identity map.
\end{corollary}
\begin{proof}
  This follows from Proposition~\ref{prop:w-dg-bifunctor} and Lemma~\ref{lem:w-mod-id-tens-id} (cf.\ Corollary~\ref{cor:mod-id-tens-id-is-sometimes-id}).
\end{proof}

\begin{corollary}\label{cor:w-tens-he}
  Fix a weighted algebra diagonal $\wDiag{*}{*}$, weighted module
  diagonals $\wModDiag{*}{*}_1$, $\wModDiag{*}{*}_2$ compatible with
  $\wDiag{*}{*}$, and $w$-modules $\wMod$ and $\wNod$ over
  $w$-algebras $\wAlg$ and $\wBlg$. Then
  $\wMod\MDtp[\wModDiagNS_1]\wNod$ and
  $\wMod\MDtp[\wModDiagNS_2]\wNod$ are homotopy equivalent
  $w$-modules.
\end{corollary}
\begin{proof}
  By Lemma~\ref{lem:w-mod-map-diag-exists-unique}, there is a weighted
  module-map diagonal $\wModMulDiag{*}{*}_{12}$ compatible with
  $\wModDiag{*}{*}_1$ and $\wModDiag{*}{*}_2$. So, the result is
  immediate from Proposition~\ref{prop:w-dg-bifunctor} and
  Lemma~\ref{lem:w-mod-id-tens-id}.
\end{proof}

\begin{proof}[Proof of Theorem~\ref{thm:wMods}]
  Both the fact that $\wMod_1\wMDtp\wMod_2$ is a weighted module and
  the fact that this tensor product reduces to the tensor product of
  $\Ainf$-modules in the trivially weighted case
  (part~\ref{item:w-Mod-thm-unweighted} of the theorem) are immediate
  from the definitions. Invariance under homotopy equivalences
  (part~\ref{item:w-Mod-thm-qi}) is
  Corollary~\ref{cor:w-tens-he}. Invariance under changing the diagonal
  (part~\ref{item:w-Mod-thm-change-diag}) is Lemma~\ref{lem:w-mod-id-tens-id}.
\end{proof}

\subsection{Type \textalt{\DD}{DD} structures and box products}
\subsubsection{The category of weighted type \textalt{\DD}{DD} bimodules}
The following is the weighted analogue of Definition~\ref{def:DD}:
\begin{definition}\label{def:wDD}
  Let $\wAlg$ and $\wBlg$ be $w$-algebras over $\Ground_1$
  and~$\Ground_2$, respectively. Let $\kappa=\kappa_1+\kappa_2-2$ and
  fix an element $X_{12}\in\Ground$ of grading $-\kappa$, such that
  $X_{12}$ acts centrally on all $\Ground$-bimodules under
  consideration.
  Fix also a weighted algebra diagonal $\wDiag{n}{w}$ with some
  seed~$\wSeed$. A
  \emph{left-left type \DD\ structure over $\wAlg$ and $\wBlg$ with
    respect to $\wDiag{n}{w}$ with charge $X_{12}$} is a (left) type
  $D$ structure over $\wAlg\wADtp \wBlg$ with charge
  $X_{12}$ (Definition~\ref{def:wD}). Similarly, the morphism complex
  $\Mor(\lsup{\wAlg,\wBlg}P,\lsup{\wAlg,\wBlg}Q)$ between weighted
  type \DD\ structures is defined to be the morphism complex between
  the corresponding type $D$ structure over $\wAlg\wADtp \wBlg$ and,
  more generally, the $\Ainf$-category of type \DD\ structures over
  $\wAlg$ and $\wBlg$ with charge $X_{12}$ is the category of type $D$
  structures over $\wAlg\wADtp \wBlg$.
\end{definition}

Explicitly, the structure map of a weighted type \DD\ structure
$\lsup{\wAlg,\wBlg}P$ is a map
\[
  \delta^1\co P\to (A\rotimes{\Ring} B)\otimes_{\Ground} P\grs{1}\cong A\otimes_{\Ground_1} P\otimes_{\Ground_2} B^\op\grs{1}.
\]
(Here, for a change, we have not suppressed the ground rings from the tensor product notation.)
We can iterate $\delta^1$ to get maps
\begin{align*}
  \delta^i &\co P\to (A\rotimes{\Ring} B)^{\kotimes{\Ground} i}\kotimes{\Ground} P\grs{i}\\
  \delta &\co P\to \prod_{i\geq 0} (A\rotimes{\Ring} B)^{\kotimes{\Ground} i}\kotimes{\Ground} P\grs{i}.
\end{align*}
The structure equation for a weighted type \DD\ structure is given
graphically by
\begin{equation}\label{eq:w-DD-str-eq}
  \mathcenter{
    \displaystyle\sum_{w=0}^\infty
    \displaystyle\sum_{i=0}^\infty
    X_{12}^w}
  \mathcenter{
  \tikzsetnextfilename{def-wDD-1}
    \begin{tikzpicture}[smallpic]
      \node at (0,0) (tc) {};
      \node at (0,-1) (delta) {$\delta^i$};
      \node at (0,-3) (bc) {};
      \node at (-1,-2) (muL) {$\wDiagCell{i}{w}$};
      \node at (1,-2) (muR) {$(\wDiagCell{i}{w})^\op$};
      \node at (-1,-3) (bl) {};
      \node at (1,-3) (br) {};
      \draw[dmoda] (tc) to (delta);
      \draw[dmoda] (delta) to (bc);
      \draw[taa] (delta) to (muL);
      \draw[taa] (delta) to (muR);
      \draw[alga] (muL) to (bl);
      \draw[alga] (muR) to (br);
    \end{tikzpicture}
  }
  \mathcenter{=0.}
\end{equation}
This makes sense as long as $\wAlg$ and $\wBlg$ are bonsai (so
$\wAlg\wADtp\wBlg$ is) or $\lsup{\wAlg,\wBlg}P$ is \emph{operationally
  bounded}, in the sense that:
\begin{enumerate}
\item $\delta^n=0$ for $n$ sufficiently large and
\item $X_{12}$ is nilpotent, i.e., for sufficiently large $n$,
  $X_{12}^n\cdot (a\otimes b)=0$ for all $a\otimes b\in A\otimes B$.
\end{enumerate}

Similarly, given type \DD\ structures $\lsup{\wAlg,\wBlg}P$,
$\lsup{\wAlg,\wBlg}Q$ over $\wAlg$ and $\wBlg$ with the same charge $X_{12}$,
the differential of a morphism $f^1\co P\to (A\rotimes{\Ring} B)\kotimes{\Ground} Q$ is given by
\[
  \mathcenter{d(f^1)}=
  \mathcenter{
    \displaystyle\sum_{w=0}^\infty
    \displaystyle\sum_{i,j=0}^\infty
    X_{12}^{w}}
  \mathcenter{
  \tikzsetnextfilename{def-wDD-2}
    \begin{tikzpicture}[smallpic]
      \node at (0,0) (tc) {};
      \node at (0,-1) (delta1) {$\delta^i$};
      \node at (0,-2) (f) {$f^1$};
      \node at (0,-3) (delta2) {$\delta^j$};      
      \node at (0,-5) (bc) {};
      \node at (-1,-4) (muL) {$\wDiagCell{i}{w}$};
      \node at (1,-4) (muR) {$(\wDiagCell{i}{w})^\op.$};
      \node at (-1,-5) (bl) {};
      \node at (1,-5) (br) {};
      \draw[dmoda] (tc) to (delta1);
      \draw[dmoda] (delta1) to (f);
      \draw[dmoda] (f) to (delta2);
      \draw[dmoda] (delta2) to (bc);
      \draw[taa] (delta1) to (muL);
      \draw[taa] (delta1) to (muR);
      \draw[taa] (delta2) to (muL);
      \draw[taa] (delta2) to (muR);
      \draw[alga] (f) to (muL);
      \draw[alga] (f) to (muR);
      \draw[alga] (muL) to (bl);
      \draw[alga] (muR) to (br);
    \end{tikzpicture}
  }
\]

\begin{warning}
  As in the unweighted case, the category of type \DD\ structures over
  $\wAlg$ and $\wBlg$ is a non-unital category, since
  $\wAlg\wADtp \wBlg$ is a non-unital weighted algebra. This can be
  rectified by considering homotopy unital weighted algebras. See
  Section~\ref{sec:hu-D} and, in particular,
  Theorem~\ref{thm:DD-indep-diag} for independence of the category of
  type \DD\ structures from the choice of diagonal.
\end{warning}

\subsubsection{The triple box product}
The following is the weighted analogue of Definition~\ref{def:triple-prod}:
\begin{definition}
  \label{def:wtriple-prod}
  Fix a weighted algebra diagonal $\wDiag{*}{*}$ and a compatible
  weighted module diagonal $\wModDiag{*}{*}$.  Let $\wAlg$ and $\wBlg$
  be strictly unital $w$-algebras, $\wMod_{\wAlg}$ and $\wNod_{\wBlg}$
  strictly unital $w$-modules, and $\lsup{\wAlg,\wBlg}P$ a type \DD\
  structure with charge $X_{12}$. Assume that either
  \begin{itemize}
  \item $\wMod$, $\wNod$, $\wAlg$, and $\wBlg$ are bonsai or
  \item $P$ is operationally bounded.
  \end{itemize}
  Then define 
  \[
    [\wMod_{\wAlg}\DT \lsup{\wAlg}P^{\wBlgop}\DT \lsub{\wBlgop}\wNod]_{\wModDiagNS}
    \coloneqq
    (\wMod\otimes_\wModDiagNS \wNod)_{\wAlg\otimes_{\wDiagNS}\wBlg} \,\DT\, \lsup{\wAlg\otimes_{\wDiagNS} \wBlg}P.
  \]
\end{definition}
Explicitly, the differential on the triple box product is
\[
  \tikzsetnextfilename{w-trip-box-prod-1}
  \sum_{n,w\geq 0} X_{12}^w\mathcenter{
    \begin{tikzpicture}[smallpic]
      \node at (0,0) (tl) {};
      \node at (2,0) (tc) {};
      \node at (4,0) (tr) {};
      \node at (2,-1) (delta) {$\delta^n$};
      \node at (0,-2) (ml) {$\wModDiagCell{n+1}{w}$};
      \node at (4,-2) (mr) {$(\wModDiagCell{n+1}{w})^\op$};
      \node at (0,-3) (bl) {};
      \node at (2,-3) (bc) {};
      \node at (4,-3) (br) {};
      \draw[dmoda] (tc) to (delta);
      \draw[dmoda] (delta) to (bc);
      \draw[taa] (delta) to (ml);
      \draw[taa] (delta) to (mr);
      \draw[moda] (tr) to (mr);
      \draw[moda] (tl) to (ml);
      \draw[moda] (mr) to (br);
      \draw[moda] (ml) to (bl);
    \end{tikzpicture}
  }
\]

Recall that Theorem~\ref{thm:wTripleTensorProduct} asserts that the
weighted triple box product is well-defined and, up to homotopy equivalence,
independent of $\wModDiag{*}{*}$:
\begin{proof}[Proof of Theorem~\ref{thm:wTripleTensorProduct}]
  The first statement follows from the facts that tensor product of two weighted
  modules is a weighted module (see Definition~\ref{def:wMod-tp}) and the box product of a weighted
  module and a weighted type $D$ structure is a chain complex (Lemma~\ref{lem:w-DT-sq-0}). The
  second statement follows from the facts that different module
  diagonals give homotopy equivalent tensor products of weighted modules (Corollary~\ref{cor:w-tens-he})
  and the box tensor product respects
  homotopy equivalence (Lemma~\ref{lem:w-DT-bifunc}). 
\end{proof}

\subsubsection{The one-sided box tensor product}
Next we turn to the weighted analogue of the one-sided box product
(Definition~\ref{def:one-side-DT}). Recall from Section~\ref{sec:wAlgs} that given a $w$-algebra
$(\wAlg,\mu^{A})$ and a weighted (algebra) tree $T$ with $n$ inputs and weight $w$
there is a corresponding map 
\[
  \mu^{A}(T)\co A^{\kotimes{\Ground_1} n}\to A\grs{-\dim(T)+(2-\kappa)w}.
\]
Similarly, as described in Section~\ref{sec:w-mods}, given a $w$-module
$(\wMod,m^M)$ over $\wAlg$ and a weighted module tree $T$ with $1+n$
inputs there is a corresponding map
\[
  m^{M}(T)\co M\kotimes{\Ground_1} A^{\kotimes{\Ground_1} n}\to M\grs{-\dim(T)+(2-\kappa)w}.
\]

\begin{definition}\label{def:w-one-sided-DT}
  Fix a weighted algebra diagonal $\wDiagCell{*}{*}$ and a weighted
  module diagonal primitive $\wTrPMDiag{*}{*}$ compatible with
  $\wDiagCell{*}{*}$.  Let $\wAlg$, $\wBlg$ be strictly unital
  $w$-algebras over $\Ground_1$ and $\Ground_2$,
  respectively, $\lsup{\wAlg,\wBlg}P$ a type \DD\ structure with
  respect to $\wDiagCell{*}{*}$, with charge $X_{12}$, and $\wMod_{\wAlg}$ a strictly unital
  $w$-module over $\wAlg$. Assume that either $\lsup{\wAlg,\wBlg}P$ is operationally
  bounded or $\wMod$, $\wAlg$, and $\wBlg$ are bonsai. Define $\wMod\DT^{\wTrPMDiagNS}P$ to be the type $D$
  structure $M\otimes_{\Ground_1} P$ with structure map
  \[ 
    \delta^1\co M\otimes_{\Ground_1} P
    \to B\otimes_{\Ground_2}(M\otimes_{\Ground_1} P)\grs{1} 
    =(M\otimes_{\Ground_1} P)\otimes_{\Ground_2}B^\op\grs{1}
  \]
  defined by
  \[
    \delta^1(x\otimes y)=m_1^0(x)\otimes y \otimes 1+
    \sum_{n,w\geq 0}\sum_{(S,T)\in \wTrPMDiag{n+1}{w}}X_{12}^w(m^M(S) \otimes \Id_{P}\otimes
    \mu^{\wBlg^{\op}}(T^\op))\circ (\Id_M \otimes\delta^n_P),
  \]
  and charge
  \[
    X_2=X_{12}Y_1.
  \]
\end{definition}
(Unlike our usual convention, we have not suppressed the ground rings for the
tensor products here.)

Graphically, the differential on $\wMod\DT^{\wTrPMDiagNS}P$ is
\[
  \tikzsetnextfilename{w-1-sided-DT-1}
  \mathcenter{\begin{tikzpicture}[smallpic]
      \node at (0,0) (tl) {};
      \node at (.75,0) (tc) {};
      \node at (0,-1) (m) {$m_1^0$};
      \node at (1.5,-1) (one) {$\unit$};
      \node at (0,-2) (bl) {};
      \node at (.75,-2) (bc) {};
      \node at (1.5,-2) (br) {};
      \draw[dmoda] (tc) to (bc);
      \draw[moda] (tl) to (m);
      \draw[moda] (m) to (bl);
      \draw[alga] (one) to (br);    
    \end{tikzpicture}}
  \qquad+\qquad\sum_{n,w\neq (0,0)}
  X_{12}^w
  \tikzsetnextfilename{w-1-sided-DT-2}
  \mathcenter{
    \begin{tikzpicture}[smallpic]
      \node at (0,0) (tl) {};
      \node at (1.5,0) (tc) {};
      \node at (0,-2) (pl) {$\wTrPMDiag{n+1}{w}$};
      \node at (3,-2) (pr) {$(\wTrPMDiag{n+1}{w})^\op$};
      \node at (1.5,-1) (delta) {$\delta^n$};
      \node at (0,-3) (bl) {};
      \node at (1.5,-3) (bc) {};
      \node at (3,-3) (br) {};
      \draw[dmoda] (tc) to (delta);
      \draw[dmoda] (delta) to (bc);
      \draw[taa] (delta) to (pl);
      \draw[taa] (delta) to (pr);
      \draw[moda] (tl) to (pl);
      \draw[moda] (pl) to (bl);
      \draw[alga] (pr) to (br);    
    \end{tikzpicture}}.
\]

\begin{proposition}\label{prop:one-sided-DT-works}
  Assuming that either $\lsup{\wAlg,\wBlg}P$ is operationally bounded
  or $\wMod$, $\wAlg$, and $\wBlg$ are bonsai, the operation
  $\delta^1$ from Definition~\ref{def:w-one-sided-DT} respects the
  gradings and makes $\wMod\DT^{\wTrPMDiagNS}P$ into a weighted type
  $D$ structure with charge $X_{12} Y_1$. If $\lsup{\wAlg,\wBlg}P$ is
  operationally bounded then so is $\wMod\DT^{\wTrPMDiagNS}P$.

  Further, if $N_{\wBlg}$ is a $w$-module (bonsai if
  $\lsup{\wAlg,\wBlg}P$ is not operationally
  bounded)
  and $\wModDiagCell{*}{*}$ is the weighted module diagonal induced by
  $\wTrPMDiag{*}{*}$ then there is an isomorphism
  \[
  (\wMod\DT^{\wTrPMDiagNS} P)\DT \wNod \cong
  [\wMod\DT P \DT \wNod]_{\wModDiagNS},
  \] 
  where the right-hand-side denotes the weighted triple box product of
  Definition~\ref{def:wtriple-prod}.
\end{proposition}
\begin{proof}
  First, we check the gradings. The first term in the definition of
  $\delta^1$ decreases the grading by $1$. For the second term, we
  have
  \begin{align*}
    \gr(X_{12}^w)&=-(\kappa_1+\kappa_2-2)w\\
    \gr(\delta_P^n)&=-n\\
    \gr(m^M(S)\otimes\Id_P\otimes \mu^{\wBlg^\op}(T^\op))
                 &=\dim(S)-(2-\kappa_1)w+\dim(T)-(2-\kappa_2)w\\
                 &=(n+1+2w-2)-(4-\kappa_1-\kappa_2)w\\
                 &=n-1+(\kappa_1+\kappa_2-2)w.
  \end{align*}
  Summing these, the operation $\delta^1$ has grading $-1$, as
  claimed.
  
  Next, we check the type $D$ structure relation
  (Equation~\eqref{eq:w-D-str}). The left side of the type $D$
  structure relation is
  \[
    \sum_{\substack{k\geq 1\\n_1,\dots,n_k,w_1,\dots,w_k,w\geq 0}} (X_{12}Y_1)^w(X_{12}^{w_1+\dots+w_k})
    \tikzsetnextfilename{w-DT-behaves-1}
    \mathcenter{\begin{tikzpicture}[smallpicwide]
        \node at (-1,0) (tl) {};
        \node at (0,0) (tc) {};
        \node at (0,-1) (delta1) {$\delta^{n_1}$};
        \node at (0,-2) (vdotsc) {$\vdots$};
        \node at (0,-3) (delta2) {$\delta^{n_k}$};
        \node at (0,-6) (bc) {};
        \node at (-1,-2) (p1l) {$\wTrPMDiag{n_1+1}{w_1}$};
        \node at (-1,-3) (vdotsl) {$\vdots$};
        \node at (-1,-4) (p2l) {$\wTrPMDiag{n_k+1}{w_k}$};
        \node at (-1,-6) (bl) {};
        \node at (1,-4) (p2r) {$(\wTrPMDiag{n_k}{w_k})^\op$};
        \node at (3,-4) (p1r) {$(\wTrPMDiag{n_1}{w_1})^\op$};
        \node at (2,-4) (cdots) {$\cdots$};
        \node at (2,-5) (mu) {$\mu_k^w$};
        \node at (2,-6) (br) {};
        \draw[dmoda] (tc) to (delta1);
        \draw[dmoda] (delta1) to (vdotsc);
        \draw[dmoda] (vdotsc) to (delta2);
        \draw[dmoda] (delta2) to (bc);
        \draw[moda] (tl) to (p1l);
        \draw[moda] (p1l) to (vdotsl);
        \draw[moda] (vdotsl) to (p2l);
        \draw[moda] (p2l) to (bl);
        \draw[taa] (delta1) to (p1l);
        \draw[taa] (delta2) to (p2l);
        \draw[taa] (delta1) to (p1r);
        \draw[taa] (delta2) to (p2r);
        \draw[alga] (p1r) to (mu);
        \draw[alga] (p2r) to (mu);
        \draw[alga] (mu) to (br);
      \end{tikzpicture}}
  \]
  (where in the case $(n_i,w_i)=(0,0)$ we
  let $(\wTrPMDiag{n_1+1}{w_1})\otimes (\wTrPMDiag{n_1+1}{w_1})^\op$ denote
  $m_1^0\otimes\Id$). From the primitive structure equation, this is
  \[
    \sum_{\substack{u+v=w+w_1+\cdots+w_k\\ m_1,m_3,u,v\geq 0}} X_{12}^{w+w_1+\cdots+w_k}
    \mathcenter{
    \tikzsetnextfilename{w-DT-behaves-2}
      \begin{tikzpicture}[smallpicwide]
        \node at (-2,0) (tl) {};
        \node at (-2,-5) (pl)  {$\wTrPMDiag{m_1+m_3+2}{v}$};
        \node at (-2,-6) (bl) {};
        \node at (0,0) (tc) {};
        \node at (0,-1) (delta1) {$\delta^{m_1}$};
        \node at (0,-2) (delta2) {$\delta^{m_2}$};
        \node at (0,-3) (delta3) {$\delta^{m_3}$};
        \node at (0,-6) (bc) {};
        \node at (-1,-3) (gammal) {$\wDiagCell{m_2}{u}$};
        \node at (1,-3) (gammar) {$\wDiagCell{m_2}{u}$};
        \node at (2,-5) (pr) {$(\wTrPMDiag{m_1+m_3+2}{v})^\op$};
        \node at (2,-6) (br) {};
        \draw[dmoda] (tc) to (delta1);
        \draw[dmoda] (delta1) to (delta2);
        \draw[dmoda] (delta2) to (delta3);
        \draw[dmoda] (delta3) to (bc);
        \draw[moda] (tl) to (pl);
        \draw[moda] (pl) to (bl);
        \draw[taa, bend right=20] (delta1) to (pl);
        \draw[taa, bend left=20] (delta1) to (pr);
        \draw[taa] (delta2) to (gammal);
        \draw[taa] (delta2) to (gammar);
        \draw[taa] (delta3) to (pl);
        \draw[taa] (delta3) to (pr);
        \draw[alga] (gammal) to (pl);
        \draw[alga] (gammar) to (pr);
        \draw[alga] (pr) to (br);
      \end{tikzpicture}
    }
  \]
  (where in the case $(m_2,u)=(0,0)$ we
  let $(\wDiagCell{m_2}{u})\otimes(\wDiagCell{m_2}{u})^\op$ denote
  $\mu_1^0\otimes\Id$).
  So, the result follows from the type \DD\ structure equation
  (Equation~\eqref{eq:w-DD-str-eq}). Note, in particular, that this
  computation explains the (perhaps surprising) charge $X_{12} Y_1$.

  For the statement about boundedness, if $X_{12}$ is nilpotent then
  so is $X_{12}Y_1$. Since $(m_1^0)^2=0$, every term in the operation $\delta^{2n+2}$ on
  $\wMod\DT^{\wTrPMDiagNS}P$ involves $\delta^u$ and a coefficient of $X_{12}^v$ where $u+v\geq n$.
  Consequently, $\delta^m$ on $\wMod\DT^{\wTrPMDiagNS}P$
  vanishes for $m$ sufficiently large.

  The proof of associativity of $\DT$ is similar but easier, and is
  left to the reader.
\end{proof}

\subsubsection{Functoriality of the one-sided box tensor product}
\label{sec:w-one-sided-DT-func}
Next, we turn to functoriality of the one-sided box product. Given
weighted type \DD\ structures $P$ and $Q$ over $\wAlg_1$ and $\wAlg_2$
with respect to a weighted algebra diagonal $\wDiag{*}{*}$, with the same charge $X_{12}$; a
morphism $f\co P\to Q$; a $w$-module $\wMod$ over $\wAlg_1$; and a
weighted module diagonal primitive $\wTrPMDiag{*}{*}$, under suitable
boundedness hypotheses (see
Lemma~\ref{lem:w-one-side-DT-id-chain-map}) there is a type $D$
structure morphism
\[
\Id_{\wMod}\DT^{\wTrPMDiagNS}f^1\co \wMod\DT^{\wTrPMDiagNS}{P}\to \wMod\DT^{\wTrPMDiagNS}{Q}
\]
defined by
\begin{equation}\label{eq:w-DD-map-tens-id}
  \Id_{\cModule}\DT^{\TrPMDiag}f^1=
  \sum_{n,w\geq 0}X_{12}^w
  \mathcenter{
    \tikzsetnextfilename{w-Id-DT-g}
  \begin{tikzpicture}[smallpic]
    \node at (0,0) (tc) {};
    \node at (-1.5,0) (tl) {};
    \node at (0, -1) (delta1) {$\delta$};
    \node at (0, -2) (f) {$f^1$};
    \node at (0, -3) (delta2) {$\delta$};
    \node at (-1.5, -4) (pl) {$\wTrPMDiag{n}{w}$};
    \node at (1.5, -4) (pr) {$(\wTrPMDiag{n}{w})^\op$};
    \node at (-1.5, -5) (bl) {};
    \node at (0, -5) (bc) {};
    \node at (1.5, -5) (br) {};
    \draw[dmoda] (tc) to (delta1);
    \draw[dmoda] (delta1) to (f);
    \draw[dmoda] (f) to (delta2);
    \draw[dmoda] (delta2) to (bc);
    \draw[moda] (tl) to (pl);
    \draw[moda] (pl) to (bl);
    \draw[taa] (delta1) to (pl);
    \draw[taa] (delta1) to (pr);
    \draw[alga] (f) to (pl);
    \draw[blga] (f) to (pr);
    \draw[tbb] (delta2) to (pl);
    \draw[tbb] (delta2) to (pr);
    \draw[alga] (pr) to (br);
  \end{tikzpicture}}
\end{equation}
(cf.~Formula~\eqref{eq:DD-map-tens-id}).

Similarly, given a weighted module-map primitive $\wTrPMorDiag{*}{*}$
compatible with weighted module diagonal primitives
$\wTrPMDiag{*}{*}_1$ and $\wTrPMDiag{*}{*}_2$, another $w$-module
$\wNod$ over $\wAlg_1$ and a morphism $g\co \wMod\to \wNod$, if either
$P$ is operationally bounded or $\wAlg_1$, $\wAlg_2$, $\wMod$, and
$\wNod$ are bonsai, there is a type $D$ structure morphism
\[
  g \DT^{\wTrPMorDiagNS}\Id_P
  \co \wMod\DT^{\wTrPMDiagNS_1}P\to \wNod\DT^{\wTrPMDiagNS_2}P
\]
defined by
\begin{equation}\label{eq:w-DD-id-tens-map}
  \sum_{n,w\geq 0}X_{12}^w
  \mathcenter{
    \tikzsetnextfilename{w-g-DT-Id}
  \begin{tikzpicture}[smallpic]
    \node at (0,0) (tc) {};
    \node at (-1.5,0) (tl) {};
    \node at (0,-1) (delta) {$\delta$};
    \node at (-1.5,-2) (pl) {$\wTrPMorDiag{n}{w} g$};
    \node at (1.5,-2) (pr) {$(\wTrPMorDiag{n}{w})^\op$};
    \node at (-1.5,-3) (bl) {};
    \node at (0,-3) (bc) {};
    \node at (1.5,-3) (br) {};
    \draw[dmoda] (tc) to (delta);
    \draw[dmoda] (delta) to (bc);
    \draw[moda] (tl) to (pl);
    \draw[moda] (pl) to (bl);
    \draw[taa] (delta) to (pl);
    \draw[taa] (delta) to (pr);
    \draw[alga] (pr) to (br);
  \end{tikzpicture}}
\end{equation}
(cf.~Formula~\eqref{eq:DD-id-tens-map}).

\begin{lemma}\label{lem:w-one-side-DT-id-chain-map}
  Suppose that either $P$ is operationally bounded or $\wAlg$,
  $\wBlg$, $\wMod$, $\wNod$, and $g$ are bonsai. Then
  the maps $\Id_{\wMod}\DT^{\wTrPMDiagNS}\cdot$ and $\cdot
  \DT^{\wTrPMorDiagNS}\Id_P$ are grading-preserving chain maps
  \begin{align*}
    \Id_{\wMod}\DT^{\wTrPMDiagNS}\cdot\co \Mor(P,Q) &\to \Mor(\wMod\DT^{\wTrPMDiagNS}{P}, \wMod\DT^{\wTrPMDiagNS}{Q})\\
    \cdot \DT^{\wTrPMorDiagNS}\Id_P\co \Mor(\wMod,\wNod)&\to \Mor(\wMod\DT^{\wTrPMDiagNS_1}P,\wNod\DT^{\wTrPMDiagNS_2}P).                                                          
  \end{align*}
  Moreover, for $\cdot \DT^{\wTrPMorDiagNS}\Id_P$, homotopic module-map
  primitives give homotopic chain maps.
\end{lemma}
\begin{proof}
  This follows from the same argument as Lemma~\ref{lem:DT-Id-dg-func}.
\end{proof}

\begin{proposition}\label{prop:w-one-side-DT-bifunc}
  Fix a weighted algebra diagonal $\wDiag{*}{*}$; weighted module
  diagonal primitives $\wTrPMDiag{*}{*}_1$, $\wTrPMDiag{*}{*}_2$, and
  $\wTrPMDiag{*}{*}_3$ compatible with $\wDiagCell{*}{*}$; and weighted
  module-map primitives $\wTrPMorDiag{*}{*}_{12}$ compatible with
  $\wTrPMDiag{*}{*}_1$ and $\wTrPMDiag{*}{*}_2$,
  $\wTrPMorDiag{*}{*}_{23}$ compatible with $\wTrPMDiag{*}{*}_2$ and
  $\wTrPMDiag{*}{*}_3$, and $\wTrPMorDiag{*}{*}_{13}$ compatible with
  $\wTrPMDiag{*}{*}_1$ and $\wTrPMDiag{*}{*}_3$. Let $\wAlg$ and
  $\wBlg$ be $w$-algebras, $\wMod_1$, $\wMod_2$, and
  $\wMod_3$ $w$-modules over $\wAlg$, and $P_1$, $P_2$, and $P_3$
  type \DD\ structures over $\Alg$ and $\Blg$ with the same charge $X_{12}$. Assume that either the
  $P_i$ are operationally bounded or $\wAlg$, $\wBlg$, and the $\wMod_i$ are bonsai.
  Then the analogues of
  Diagrams~\eqref{eq:mm-com-1},~\eqref{eq:mm-com-2}, and~\eqref{eq:mm-com-3}
  commute up to homotopy where, in the case that the $\wMod_i$ are
  bonsai, we consider the bonsai morphism complexes.
\end{proposition}
\begin{proof}
  The proof is similar to the proof of
  Proposition~\ref{prop:w-one-side-DT-bifunc}, and is left to the reader.
\end{proof}

\begin{lemma}\label{lem:w-Id-DT-Id}
  Fix weighted module diagonal primitives $\wTrPMDiag{*}{*}_1$ and $\wTrPMDiag{*}{*}_2$ and a
  weighted module-map primitive $\wTrPMorDiag{*}{*}$ compatible with $\wTrPMDiag{*}{*}_1$ and
  $\wTrPMDiag{*}{*}_2$.
  Given $w$-algebras $\wAlg$ and $\wBlg$, a
  $w$-module $\wMod_{\wAlg}$ and a type \DD\ structure
  $\lsup{\wAlg,\wBlg}P$, the map
  \[
    \Id_{\wMod} \DT^{\wTrPMorDiagNS}\Id_P\co \wMod\DT^{\wTrPMDiagNS_1}P\to \wMod\DT^{\wTrPMDiagNS_2}P
  \]
  is a homotopy equivalence. Further, if
  $\wTrPMDiag{*}{*}_1=\wTrPMDiag{*}{*}_2$ then
  $\Id_{\wMod} \DT^{\wTrPMorDiagNS}\Id_P$
  is homotopic to the identity map.
\end{lemma}
\begin{proof}
  The proof is essentially the same as the proofs
  Lemma~\ref{lem:Id-DT-Id} and Corollary~\ref{cor:Id-DT-Id}. First, in
  the special case that $\wTrPMDiag{*}{*}_1=\wTrPMDiag{*}{*}_2$ and
  $\wTrPMorDiag{*}{*}$ is the weighted module-map primitive from
  Lemma~\ref{lem:w-mor-prim-example},
  $\Id_{\wMod} \DT^{\wTrPMorDiagNS}\Id_P$ is the identity map. The
  fact that $\Id_{\wMod} \DT^{\wTrPMorDiagNS}\Id_P$ is homotopic to
  the identity map for any weighted module-map primitive compatible
  with $\wTrPMDiag{*}{*}_1$ and $\wTrPMDiag{*}{*}_1$ follows from the
  facts that all weighted module-map primitives are homotopic
  (Lemma~\ref{lem:w-mod-map-diag-exists-unique}) and homotopic
  primitives give homotopic maps
  (Lemma~\ref{lem:w-one-side-DT-id-chain-map}). The case that
  $\wTrPMDiag{*}{*}_1\neq \wTrPMDiag{*}{*}_2$ then follows from
  Proposition~\ref{prop:w-one-side-DT-bifunc}.
\end{proof}

\begin{corollary}\label{cor:w-DT-preserve-hequiv}
  Fix strictly unital $w$-modules $\wMod$ and $\wNod$, a homomorphism
  $f\co\wMod\to\wNod$, and a type \DD\ structure $P$ such that either
  $P$ is operationally bounded or $\wMod$, $\wNod$, and $f$ are
  bonsai.  If $f$ is a homotopy equivalence then for
  any choice of module diagonal primitives and module-map primitive,
  \[
    f \DT^{\wTrPMorDiagNS}\Id_P
    \co \wMod\DT^{\wTrPMDiagNS_1}P\to \wNod\DT^{\wTrPMDiagNS_2}P
  \]
  is a homotopy equivalence.
\end{corollary}
\begin{proof}
  This is immediate from Lemma~\ref{lem:w-Id-DT-Id} and
  Proposition~\ref{prop:w-one-side-DT-bifunc}.
\end{proof}

\subsubsection{Associativity of the box tensor product of weighted morphisms}
Similarly to the unweighted case, given a weighted type \DD\ structure
$P$, strictly unital $w$-modules $\wMod_1$, $\wMod_2$, $\wNod_1$, and
$\wNod_2$, and morphisms $f\co \wMod_1\to\wMod_2$ and
$g\co\wNod_1\to\wNod_2$, as well as a weighted algebra diagonal
$\wDiag{*}{*}$ , weighted module diagonals $\wModDiag{*}{*}_1$ and
$\wModDiag{*}{*}_2$ compatible with $\wDiag{*}{*}$, and a weighted
module-map diagonal $\wModMulDiag{*}{*}$ compatible with
$\wModDiag{*}{*}_1$ and $\wModDiag{*}{*}_2$,
Definition~\ref{def:w-tensor-mod-maps} gives a triple box product
\[
  [f\DT\Id_P\DT g]_{\wModMulDiagNS}=(f\wMDtp[\wModMulDiagNS]g)\DT \Id_P
\]
(under appropriate boundedness assumptions).
In the special case that $g$ is the identity map, $[f\DT \Id_P\DT
g]_{\wModMulDiagNS}$ depends only on the partial weighted module-map diagonal
induced by $\wModMulDiag{*}{*}$. That is, given a partial weighted module-map diagonal $\wPartTrModMulDiag{*}{*}$ compatible with $\wModDiag{*}{*}_1$ and $\wModDiag{*}{*}_2$ we can define
\begin{equation}
  \label{eq:w-part-mor-tens}
\mathcenter{  [f\DT\Id_P\DT \Id_{\wNod}]_{\wPartTrModMulDiagNS}=}
  \tikzsetnextfilename{eq-w-part-mor-tens-1}
  \mathcenter{\begin{tikzpicture}[smallpic]
    \node at (0,0) (tc) {};
    \node at (-1,0) (tl) {};
    \node at (1,0) (tr) {};
    \node at (0,-1) (delta) {$\delta$};
    \node at (-1,-2) (pl) {$\wPartTrModMulDiagNS f$};
    \node at (1,-2) (pr) {$\wPartTrModMulDiagNS^\op\Id$};
    \node at (-1,-3) (bl) {};
    \node at (0,-3) (bc) {};
    \node at (1,-3) (br) {};
    \draw[dmoda] (tc) to (delta);
    \draw[dmoda] (delta) to (bc);
    \draw[moda] (tl) to (pl);
    \draw[moda] (pl) to (bl);
    \draw[moda] (tr) to (pr);
    \draw[taa] (delta) to (pl);
    \draw[taa] (delta) to (pr);
    \draw[moda] (pr) to (br);
  \end{tikzpicture}}=
\sum_{n,w\geq 0}X_{12}^w
  \tikzsetnextfilename{eq-w-part-mor-tens-2}
  \mathcenter{\begin{tikzpicture}[smallpic]
    \node at (0,0) (tc) {};
    \node at (-1.5,0) (tl) {};
    \node at (1.5,0) (tr) {};
    \node at (0,-1) (delta) {$\delta^n$};
    \node at (-1.5,-2) (pl) {$\wPartTrModMulDiag{n+1}{w} f$};
    \node at (1.5,-2) (pr) {$(\wPartTrModMulDiag{n+1}{w})^\op\Id$};
    \node at (-1.5,-3) (bl) {};
    \node at (0,-3) (bc) {};
    \node at (1.5,-3) (br) {};
    \draw[dmoda] (tc) to (delta);
    \draw[dmoda] (delta) to (bc);
    \draw[moda] (tl) to (pl);
    \draw[moda] (pl) to (bl);
    \draw[moda] (tr) to (pr);
    \draw[taa] (delta) to (pl);
    \draw[taa] (delta) to (pr);
    \draw[moda] (pr) to (br);
  \end{tikzpicture}}
\end{equation}
where $\wPartTrModMulDiagNS f$ means that we apply $f$ at the
distinguished vertex of (the left tree in) $\wPartTrModMulDiag{*}{*}$
and $\wPartTrModMulDiagNS \Id$ means we apply $\Id$ at the
distinguished vertex of (the right tree in)
$\wPartTrModMulDiag{*}{*}$, and we multiply by $X_{12}$ to the weight.
\begin{lemma}\label{lem:w-part-mod-map-diag-DT}
  If either $P$ is operationally bounded or $\wAlg$, $\wBlg$, $\wMod_i$, $\wNod_i$, $i=1,2$, and $f$ are bonsai then
  Formula~\eqref{eq:w-part-mor-tens} defines a chain map
  \[
    \Mor(\wMod_1,\wMod_2)\to \Mor([\wMod_1\DT P\DT \wNod]_{\wModDiagNS_1},[\wMod_2\DT P\DT\wNod]_{\wModDiagNS_2}).
  \]
  Moreover, homotopic partial weighted module-map diagonals induce homotopic
  chain maps. Finally, if $\wPartTrModMulDiag{*}{*}$ is induced by a
  weighted module-map primitive $\TrPMorDiag$ then $[f\DT\Id_P\DT
  \Id_{\wNod}]_{\wPartTrModMulDiagNS}$ is homotopic to
  $(f\DT^{\wTrPMorDiagNS}\Id_P)\DT\Id_{\wNod}$.
\end{lemma}
\begin{proof}
  This is a straightforward adaptation of the proof of
  Lemma~\ref{lem:part-mod-map-diag-DT} to the weighted case, but we repeat the
  proof keeping track of the weights and variables. (Nothing surprising happens
  with them, either.)

  To keep the notation from getting out of hand, we will not indicate
  which module an $m^{1,0}$ or $\wModDiagCell{n}{w}$ comes from, as it is
  determined by whether the operation comes before or after the copy
  of $\wPartTrModMulDiag{*}{*}$. We will also suppress the $\op$s from the notation.

  For the first statement, $d[f\DT\Id_P\DT
  \Id_{\cNodule}]_{\PartTrModMulDiag}$ is given by
  \begin{align*}
    \sum_{u,v}X_{12}^{u+v}
  \tikzsetnextfilename{wpmm-1-1}
  &\mathcenter{\begin{tikzpicture}[smallpic]
    \node at (0,2) (tc) {};
    \node at (-1,2) (tl) {};
    \node at (1,2) (tr) {};
    \node at (0,1) (tdelta) {$\delta$};
    \node at (-1,0) (ml) {$\wModDiagCell{*}{u}$};
    \node at (1,0) (mr) {$\wModDiagCell{*}{u}$};
    \node at (0,-1) (delta) {$\delta$};
    \node at (-1,-2) (pl) {$\wPartTrModMulDiag{*}{v} f$};
    \node at (1,-2) (pr) {$\wPartTrModMulDiag{*}{v}\Id$};
    \node at (-1,-3) (bl) {};
    \node at (0,-3) (bc) {};
    \node at (1,-3) (br) {};
    \draw[dmoda] (tc) to (tdelta);
    \draw[dmoda] (tdelta) to (delta);
    \draw[dmoda] (delta) to (bc);
    \draw[moda] (tl) to (ml);
    \draw[moda] (ml) to (pl);
    \draw[moda] (pl) to (bl);
    \draw[moda] (tr) to (mr);
    \draw[moda] (mr) to (pr);
    \draw[taa] (tdelta) to (ml);
    \draw[taa] (tdelta) to (mr);
    \draw[taa] (delta) to (pl);
    \draw[taa] (delta) to (pr);
    \draw[moda] (pr) to (br);
  \end{tikzpicture}}
    +
  \tikzsetnextfilename{wpmm-1-4}
  \mathcenter{\begin{tikzpicture}[smallpic]
    \node at (0,2) (tc) {};
    \node at (-1,2) (tl) {};
    \node at (1,2) (tr) {};
    \node at (0,1) (delta) {$\delta$};
    \node at (-1,-2) (ml) {$\wModDiagCell{*}{u}$};
    \node at (1,-2) (mr) {$\wModDiagCell{*}{u}$};
    \node at (0,-1) (tdelta) {$\delta$};
    \node at (-1,0) (pl) {$\wPartTrModMulDiag{*}{v} f$};
    \node at (1,0) (pr) {$\wPartTrModMulDiag{*}{v}\Id$};
    \node at (-1,-3) (bl) {};
    \node at (0,-3) (bc) {};
    \node at (1,-3) (br) {};
    \draw[dmoda] (tc) to (delta);
    \draw[dmoda] (delta) to (tdelta);
    \draw[dmoda] (tdelta) to (bc);
    \draw[moda] (tl) to (pl);
    \draw[moda] (pl) to (ml);
    \draw[moda] (ml) to (bl);
    \draw[moda] (tr) to (pr);
    \draw[moda] (pr) to (mr);
    \draw[taa] (tdelta) to (ml);
    \draw[taa] (tdelta) to (mr);
    \draw[taa] (delta) to (pl);
    \draw[taa] (delta) to (pr);
    \draw[moda] (mr) to (br);
  \end{tikzpicture}}\\
    +
    \sum_{w}X_{12}^{w}
  \tikzsetnextfilename{wpmm-1-2}
  &\mathcenter{\begin{tikzpicture}[smallpic]
    \node at (0,1) (tc) {};
    \node at (-1,1) (tl) {};
    \node at (1,1) (tr) {};
    \node at (-1,0) (ml) {$m^{1,0}$};
    \node at (0,-1) (delta) {$\delta$};
    \node at (-1,-2) (pl) {$\wPartTrModMulDiag{*}{w} f$};
    \node at (1,-2) (pr) {$\wPartTrModMulDiag{*}{w}\Id$};
    \node at (-1,-3) (bl) {};
    \node at (0,-3) (bc) {};
    \node at (1,-3) (br) {};
    \draw[dmoda] (tc) to (delta);
    \draw[dmoda] (delta) to (bc);
    \draw[moda] (tl) to (ml);
    \draw[moda] (ml) to (pl);
    \draw[moda] (pl) to (bl);
    \draw[moda] (tr) to (pr);
    \draw[taa] (delta) to (pl);
    \draw[taa] (delta) to (pr);
    \draw[moda] (pr) to (br);
  \end{tikzpicture}}+
  \tikzsetnextfilename{wpmm-1-3}
  \mathcenter{\begin{tikzpicture}[smallpic]
    \node at (0,1) (tc) {};
    \node at (-1,1) (tl) {};
    \node at (1,1) (tr) {};
    \node at (1,0) (mr) {$m^{1,0}$};
    \node at (0,-1) (delta) {$\delta$};
    \node at (-1,-2) (pl) {$\wPartTrModMulDiag{*}{w} f$};
    \node at (1,-2) (pr) {$\wPartTrModMulDiag{*}{w}\Id$};
    \node at (-1,-3) (bl) {};
    \node at (0,-3) (bc) {};
    \node at (1,-3) (br) {};
    \draw[dmoda] (tc) to (delta);
    \draw[dmoda] (delta) to (bc);
    \draw[moda] (tl) to (pl);
    \draw[moda] (pl) to (bl);
    \draw[moda] (tr) to (mr);
    \draw[moda] (mr) to (pr);
    \draw[taa] (delta) to (pl);
    \draw[taa] (delta) to (pr);
    \draw[moda] (pr) to (br);
  \end{tikzpicture}}
+
  \tikzsetnextfilename{wpmm-1-5}
  \mathcenter{\begin{tikzpicture}[smallpic]
    \node at (0,0) (tc) {};
    \node at (-1,0) (tl) {};
    \node at (1,0) (tr) {};
    \node at (-1,-3) (ml) {$m^{1,0}$};
    \node at (0,-1) (delta) {$\delta$};
    \node at (-1,-2) (pl) {$\wPartTrModMulDiag{*}{w} f$};
    \node at (1,-2) (pr) {$\wPartTrModMulDiag{*}{w}\Id$};
    \node at (-1,-4) (bl) {};
    \node at (0,-4) (bc) {};
    \node at (1,-4) (br) {};
    \draw[dmoda] (tc) to (delta);
    \draw[dmoda] (delta) to (bc);
    \draw[moda] (tl) to (pl);
    \draw[moda] (pl) to (ml);
    \draw[moda] (ml) to (bl);
    \draw[moda] (tr) to (pr);
    \draw[taa] (delta) to (pl);
    \draw[taa] (delta) to (pr);
    \draw[moda] (pr) to (br);
  \end{tikzpicture}}+
  \tikzsetnextfilename{wpmm-1-6}
  \mathcenter{\begin{tikzpicture}[smallpic]
    \node at (0,0) (tc) {};
    \node at (-1,0) (tl) {};
    \node at (1,0) (tr) {};
    \node at (1,-3) (mr) {$m^{1,0}$};
    \node at (0,-1) (delta) {$\delta$};
    \node at (-1,-2) (pl) {$\wPartTrModMulDiag{*}{w} f$};
    \node at (1,-2) (pr) {$\wPartTrModMulDiag{*}{w}\Id$};
    \node at (-1,-4) (bl) {};
    \node at (0,-4) (bc) {};
    \node at (1,-4) (br) {};
    \draw[dmoda] (tc) to (delta);
    \draw[dmoda] (delta) to (bc);
    \draw[moda] (tl) to (pl);
    \draw[moda] (pl) to (bl);
    \draw[moda] (tr) to (pr);
    \draw[moda] (pr) to (mr);
    \draw[taa] (delta) to (pl);
    \draw[taa] (delta) to (pr);
    \draw[moda] (mr) to (br);
  \end{tikzpicture}}.
  \end{align*}
  From the definition of a partial module-map diagonal, this sum is equal to
  \begin{align*}
    \sum_w X_{12}^w
  \tikzsetnextfilename{wpmm-2-1}
    &\mathcenter{\begin{tikzpicture}[smallpic]
        \node at (0,0) (tc) {};
        \node at (-1,0) (tl) {};
        \node at (1,0) (tr) {};
        \node at (0,-1) (delta) {$\delta$};
        \node at (-1,-2) (pl) {$(\bdy\wPartTrModMulDiag{*}{w}) f$};
        \node at (1,-2) (pr) {$\wPartTrModMulDiag{*}{w}\Id$};
        \node at (-1,-3) (bl) {};
        \node at (0,-3) (bc) {};
        \node at (1,-3) (br) {};
        \draw[dmoda] (tc) to (delta);
        \draw[dmoda] (delta) to (bc);
        \draw[moda] (tl) to (pl);
        \draw[moda] (pl) to (bl);
        \draw[moda] (tr) to (pr);
        \draw[taa] (delta) to (pl);
        \draw[taa] (delta) to (pr);
        \draw[moda] (pr) to (br);
      \end{tikzpicture}}
    +
  \tikzsetnextfilename{wpmm-2-2}
        \mathcenter{\begin{tikzpicture}[smallpic]
        \node at (0,0) (tc) {};
        \node at (-1,0) (tl) {};
        \node at (1,0) (tr) {};
        \node at (0,-1) (delta) {$\delta$};
        \node at (-1,-2) (pl) {$\wPartTrModMulDiag{*}{w}f$};
        \node at (1,-2) (pr) {$(\bdy\wPartTrModMulDiag{*}{w})\Id$};
        \node at (-1,-3) (bl) {};
        \node at (0,-3) (bc) {};
        \node at (1,-3) (br) {};
        \draw[dmoda] (tc) to (delta);
        \draw[dmoda] (delta) to (bc);
        \draw[moda] (tl) to (pl);
        \draw[moda] (pl) to (bl);
        \draw[moda] (tr) to (pr);
        \draw[taa] (delta) to (pl);
        \draw[taa] (delta) to (pr);
        \draw[moda] (pr) to (br);
      \end{tikzpicture}}
      +
  \tikzsetnextfilename{wpmm-2-3}
  \mathcenter{\begin{tikzpicture}[smallpic]
    \node at (0,1) (tc) {};
    \node at (-1,1) (tl) {};
    \node at (1,1) (tr) {};
    \node at (-1,0) (ml) {$m^{1,0}$};
    \node at (0,-1) (delta) {$\delta$};
    \node at (-1,-2) (pl) {$\wPartTrModMulDiag{*}{w} f$};
    \node at (1,-2) (pr) {$\wPartTrModMulDiag{*}{w}\Id$};
    \node at (-1,-3) (bl) {};
    \node at (0,-3) (bc) {};
    \node at (1,-3) (br) {};
    \draw[dmoda] (tc) to (delta);
    \draw[dmoda] (delta) to (bc);
    \draw[moda] (tl) to (ml);
    \draw[moda] (ml) to (pl);
    \draw[moda] (pl) to (bl);
    \draw[moda] (tr) to (pr);
    \draw[taa] (delta) to (pl);
    \draw[taa] (delta) to (pr);
    \draw[moda] (pr) to (br);
  \end{tikzpicture}}+
  \tikzsetnextfilename{wpmm-2-4}
  \mathcenter{\begin{tikzpicture}[smallpic]
    \node at (0,1) (tc) {};
    \node at (-1,1) (tl) {};
    \node at (1,1) (tr) {};
    \node at (1,0) (mr) {$m^{1,0}$};
    \node at (0,-1) (delta) {$\delta$};
    \node at (-1,-2) (pl) {$\wPartTrModMulDiag{*}{w} f$};
    \node at (1,-2) (pr) {$\wPartTrModMulDiag{*}{w}\Id$};
    \node at (-1,-3) (bl) {};
    \node at (0,-3) (bc) {};
    \node at (1,-3) (br) {};
    \draw[dmoda] (tc) to (delta);
    \draw[dmoda] (delta) to (bc);
    \draw[moda] (tl) to (pl);
    \draw[moda] (pl) to (bl);
    \draw[moda] (tr) to (mr);
    \draw[moda] (mr) to (pr);
    \draw[taa] (delta) to (pl);
    \draw[taa] (delta) to (pr);
    \draw[moda] (pr) to (br);
  \end{tikzpicture}}+
  \tikzsetnextfilename{wpmm-2-5}
  \mathcenter{\begin{tikzpicture}[smallpic]
    \node at (0,0) (tc) {};
    \node at (-1,0) (tl) {};
    \node at (1,0) (tr) {};
    \node at (-1,-3) (ml) {$m^{1,0}$};
    \node at (0,-1) (delta) {$\delta$};
    \node at (-1,-2) (pl) {$\wPartTrModMulDiag{*}{w} f$};
    \node at (1,-2) (pr) {$\wPartTrModMulDiag{*}{w}\Id$};
    \node at (-1,-4) (bl) {};
    \node at (0,-4) (bc) {};
    \node at (1,-4) (br) {};
    \draw[dmoda] (tc) to (delta);
    \draw[dmoda] (delta) to (bc);
    \draw[moda] (tl) to (pl);
    \draw[moda] (pl) to (ml);
    \draw[moda] (ml) to (bl);
    \draw[moda] (tr) to (pr);
    \draw[taa] (delta) to (pl);
    \draw[taa] (delta) to (pr);
    \draw[moda] (pr) to (br);
  \end{tikzpicture}}\\&+
  \tikzsetnextfilename{wpmm-2-6}
  \mathcenter{\begin{tikzpicture}[smallpic]
    \node at (0,0) (tc) {};
    \node at (-1,0) (tl) {};
    \node at (1,0) (tr) {};
    \node at (1,-3) (mr) {$m^{1,0}$};
    \node at (0,-1) (delta) {$\delta$};
    \node at (-1,-2) (pl) {$\wPartTrModMulDiag{*}{w} f$};
    \node at (1,-2) (pr) {$\wPartTrModMulDiag{*}{w}\Id$};
    \node at (-1,-4) (bl) {};
    \node at (0,-4) (bc) {};
    \node at (1,-4) (br) {};
    \draw[dmoda] (tc) to (delta);
    \draw[dmoda] (delta) to (bc);
    \draw[moda] (tl) to (pl);
    \draw[moda] (pl) to (bl);
    \draw[moda] (tr) to (pr);
    \draw[moda] (pr) to (mr);
    \draw[taa] (delta) to (pl);
    \draw[taa] (delta) to (pr);
    \draw[moda] (mr) to (br);
  \end{tikzpicture}}+      
  \tikzsetnextfilename{wpmm-2-7}
    \sum_{u,v}X_{12}^{u+v}
    \mathcenter{
      \begin{tikzpicture}[smallpic]
        \node at (0,0) (tc) {};
        \node at (-2,0) (tl) {};
        \node at (2,0) (tr) {};
        \node at (0,-1) (delta1) {$\delta$};
        \node at (0,-2) (delta2) {$\delta$};
        \node at (0,-3) (delta3) {$\delta$};
        \node at (-2,-4) (kl) {$\wPartTrModMulDiag{*}{v} f$};
        \node at (2,-4) (kr) {$\wPartTrModMulDiag{*}{v} \Id$};
        \node at (-1,-3) (gammal) {$\wTrDiag{*}{u}$};
        \node at (1,-3) (gammar) {$\wTrDiag{*}{u}$};
        \node at (-2,-5) (bl) {};
        \node at (2,-5) (br) {};
        \node at (0,-5) (bc) {};
        \draw[dmoda] (tc) to (delta1);
        \draw[dmoda] (delta1) to (delta2);
        \draw[dmoda] (delta2) to (delta3);
        \draw[dmoda] (delta3) to (bc);
        \draw[taa, bend right=10] (delta1) to (kl);
        \draw[taa, bend left=10] (delta1) to (kr);
        \draw[taa] (delta2) to (gammal);
        \draw[taa] (delta2) to (gammar);
        \draw[taa, bend left=5] (delta3) to (kl);
        \draw[taa, bend right=5] (delta3) to (kr);
        \draw[alga] (gammal) to (kl);
        \draw[alga] (gammar) to (kr);
        \draw[moda] (tl) to (kl);
        \draw[moda] (tr) to (kr);
        \draw[moda] (kl) to (bl);
        \draw[moda] (kr) to (br);
      \end{tikzpicture}}
    +
  \tikzsetnextfilename{wpmm-2-8}
        \mathcenter{\begin{tikzpicture}[smallpic]
        \node at (0,0) (tc) {};
        \node at (-1,0) (tl) {};
        \node at (1,0) (tr) {};
        \node at (0,-1) (delta) {$\delta$};
        \node at (-1,-2) (pl) {$\alphas f$};
        \node at (1,-2) (pr) {$(\bdy\alphas)\Id$};
        \node at (-1,-3) (bl) {};
        \node at (0,-3) (bc) {};
        \node at (1,-3) (br) {};
        \draw[dmoda] (tc) to (delta);
        \draw[dmoda] (delta) to (bc);
        \draw[moda] (tl) to (pl);
        \draw[moda] (pl) to (bl);
        \draw[moda] (tr) to (pr);
        \draw[taa] (delta) to (pl);
        \draw[taa] (delta) to (pr);
        \draw[moda] (pr) to (br);
      \end{tikzpicture}}
    +
  \tikzsetnextfilename{wpmm-2-9}
    \mathcenter{\begin{tikzpicture}[smallpic]
        \node at (0,0) (tc) {};
        \node at (-1,0) (tl) {};
        \node at (1,0) (tr) {};
        \node at (0,-1) (delta) {$\delta$};
        \node at (-1,-2) (pl) {$\betas f$};
        \node at (1,-2) (pr) {$\betas\Id$};
        \node at (-1,-3) (bl) {};
        \node at (0,-3) (bc) {};
        \node at (1,-3) (br) {};
        \draw[dmoda] (tc) to (delta);
        \draw[dmoda] (delta) to (bc);
        \draw[moda] (tl) to (pl);
        \draw[moda] (pl) to (bl);
        \draw[moda] (tr) to (pr);
        \draw[taa] (delta) to (pl);
        \draw[taa] (delta) to (pr);
        \draw[moda] (pr) to (br);
      \end{tikzpicture}}
  \end{align*}
  where $\alphas,\betas\in \wMTransCx{*}{*}\rotimes{\Ring} \Filt_{\neq2}\wMTransCx{*}{*}\rotimes{\Ring}\Ring[Y_1,Y_2]$.

  Applying Lemma~\ref{lem:w-F-f}, this sum is equal to
  \begin{align*}
    \sum_w
    X_{12}^w
    \tikzsetnextfilename{wpmm-3-1}
    &\mathcenter{\begin{tikzpicture}[smallpic]
        \node at (0,0) (tc) {};
        \node at (-1,0) (tl) {};
        \node at (1,0) (tr) {};
        \node at (0,-1) (delta) {$\delta$};
        \node at (-1,-2) (pl) {$\wPartTrModMulDiag{*}{w}(df)$};
        \node at (1,-2) (pr) {$\wPartTrModMulDiag{*}{w}\Id$};
        \node at (-1,-3) (bl) {};
        \node at (0,-3) (bc) {};
        \node at (1,-3) (br) {};
        \draw[dmoda] (tc) to (delta);
        \draw[dmoda] (delta) to (bc);
        \draw[moda] (tl) to (pl);
        \draw[moda] (pl) to (bl);
        \draw[moda] (tr) to (pr);
        \draw[taa] (delta) to (pl);
        \draw[taa] (delta) to (pr);
        \draw[moda] (pr) to (br);
      \end{tikzpicture}}
    +
  \tikzsetnextfilename{wpmm-3-2}
        \mathcenter{\begin{tikzpicture}[smallpic]
        \node at (0,0) (tc) {};
        \node at (-1,0) (tl) {};
        \node at (1,0) (tr) {};
        \node at (0,-1) (delta) {$\delta$};
        \node at (-1,-2) (pl) {$\wPartTrModMulDiag{*}{w} f$};
        \node at (1,-2) (pr) {$\wPartTrModMulDiag{*}{w}(d\Id)$};
        \node at (-1,-3) (bl) {};
        \node at (0,-3) (bc) {};
        \node at (1,-3) (br) {};
        \draw[dmoda] (tc) to (delta);
        \draw[dmoda] (delta) to (bc);
        \draw[moda] (tl) to (pl);
        \draw[moda] (pl) to (bl);
        \draw[moda] (tr) to (pr);
        \draw[taa] (delta) to (pl);
        \draw[taa] (delta) to (pr);
        \draw[moda] (pr) to (br);
      \end{tikzpicture}}
      +    \mathcenter{
  \tikzsetnextfilename{wpmm-3-3}
      \begin{tikzpicture}[smallpic]
        \node at (0,0) (tc) {};
        \node at (-2,0) (tl) {};
        \node at (2,0) (tr) {};
        \node at (0,-1) (delta1) {$\delta$};
        \node at (0,-2) (delta2) {$\delta^1$};
        \node at (0,-3) (delta3) {$\delta$};
        \node at (-2,-4) (kl) {$\wPartTrModMulDiag{*}{w} f$};
        \node at (2,-4) (kr) {$\wPartTrModMulDiag{*}{w} \Id$};
        \node at (-1,-3) (gammal) {$\mu^{1,0}$};
        \node at (-2,-5) (bl) {};
        \node at (2,-5) (br) {};
        \node at (0,-5) (bc) {};
        \draw[dmoda] (tc) to (delta1);
        \draw[dmoda] (delta1) to (delta2);
        \draw[dmoda] (delta2) to (delta3);
        \draw[dmoda] (delta3) to (bc);
        \draw[taa, bend right=10] (delta1) to (kl);
        \draw[taa, bend left=10] (delta1) to (kr);
        \draw[alga] (delta2) to (gammal);
        \draw[alga] (delta2) to (kr);
        \draw[taa, bend left=5] (delta3) to (kl);
        \draw[taa, bend right=5] (delta3) to (kr);
        \draw[alga] (gammal) to (kl);
        \draw[moda] (tl) to (kl);
        \draw[moda] (tr) to (kr);
        \draw[moda] (kl) to (bl);
        \draw[moda] (kr) to (br);
      \end{tikzpicture}}
      +
      \mathcenter{
  \tikzsetnextfilename{wpmm-3-4}
      \begin{tikzpicture}[smallpic]
        \node at (0,0) (tc) {};
        \node at (-2,0) (tl) {};
        \node at (2,0) (tr) {};
        \node at (0,-1) (delta1) {$\delta$};
        \node at (0,-2) (delta2) {$\delta^1$};
        \node at (0,-3) (delta3) {$\delta$};
        \node at (-2,-4) (kl) {$\wPartTrModMulDiag{*}{w} f$};
        \node at (2,-4) (kr) {$\wPartTrModMulDiag{*}{w} \Id$};
        \node at (1,-3) (gammar) {$\mu^{1,0}$};
        \node at (-2,-5) (bl) {};
        \node at (2,-5) (br) {};
        \node at (0,-5) (bc) {};
        \draw[dmoda] (tc) to (delta1);
        \draw[dmoda] (delta1) to (delta2);
        \draw[dmoda] (delta2) to (delta3);
        \draw[dmoda] (delta3) to (bc);
        \draw[taa, bend right=10] (delta1) to (kl);
        \draw[taa, bend left=10] (delta1) to (kr);
        \draw[alga] (delta2) to (kl);
        \draw[alga] (delta2) to (gammar);
        \draw[taa, bend left=5] (delta3) to (kl);
        \draw[taa, bend right=5] (delta3) to (kr);
        \draw[alga] (gammar) to (kr);
        \draw[moda] (tl) to (kl);
        \draw[moda] (tr) to (kr);
        \draw[moda] (kl) to (bl);
        \draw[moda] (kr) to (br);
      \end{tikzpicture}}
    \\
    + &
        \sum_{u,v}
        X_{12}^{u+v}
    \mathcenter{
  \tikzsetnextfilename{wpmm-3-5}
      \begin{tikzpicture}[smallpic]
        \node at (0,0) (tc) {};
        \node at (-2,0) (tl) {};
        \node at (2,0) (tr) {};
        \node at (0,-1) (delta1) {$\delta$};
        \node at (0,-2) (delta2) {$\delta$};
        \node at (0,-3) (delta3) {$\delta$};
        \node at (-2,-4) (kl) {$\wPartTrModMulDiag{*}{v} f$};
        \node at (2,-4) (kr) {$\wPartTrModMulDiag{*}{v} \Id$};
        \node at (-1,-3) (gammal) {$\wTrDiag{*}{u}$};
        \node at (1,-3) (gammar) {$\wTrDiag{*}{u}$};
        \node at (-2,-5) (bl) {};
        \node at (2,-5) (br) {};
        \node at (0,-5) (bc) {};
        \draw[dmoda] (tc) to (delta1);
        \draw[dmoda] (delta1) to (delta2);
        \draw[dmoda] (delta2) to (delta3);
        \draw[dmoda] (delta3) to (bc);
        \draw[taa, bend right=10] (delta1) to (kl);
        \draw[taa, bend left=10] (delta1) to (kr);
        \draw[taa] (delta2) to (gammal);
        \draw[taa] (delta2) to (gammar);
        \draw[taa, bend left=5] (delta3) to (kl);
        \draw[taa, bend right=5] (delta3) to (kr);
        \draw[alga] (gammal) to (kl);
        \draw[alga] (gammar) to (kr);
        \draw[moda] (tl) to (kl);
        \draw[moda] (tr) to (kr);
        \draw[moda] (kl) to (bl);
        \draw[moda] (kr) to (br);
      \end{tikzpicture}}
    +
  \tikzsetnextfilename{wpmm-3-6}
        \mathcenter{\begin{tikzpicture}[smallpic]
        \node at (0,0) (tc) {};
        \node at (-1,0) (tl) {};
        \node at (1,0) (tr) {};
        \node at (0,-1) (delta) {$\delta$};
        \node at (-1,-2) (pl) {$\alphas f$};
        \node at (1,-2) (pr) {$(\bdy\alphas)\Id$};
        \node at (-1,-3) (bl) {};
        \node at (0,-3) (bc) {};
        \node at (1,-3) (br) {};
        \draw[dmoda] (tc) to (delta);
        \draw[dmoda] (delta) to (bc);
        \draw[moda] (tl) to (pl);
        \draw[moda] (pl) to (bl);
        \draw[moda] (tr) to (pr);
        \draw[taa] (delta) to (pl);
        \draw[taa] (delta) to (pr);
        \draw[moda] (pr) to (br);
      \end{tikzpicture}}
    +
  \tikzsetnextfilename{wpmm-3-7}
    \mathcenter{\begin{tikzpicture}[smallpic]
        \node at (0,0) (tc) {};
        \node at (-1,0) (tl) {};
        \node at (1,0) (tr) {};
        \node at (0,-1) (delta) {$\delta$};
        \node at (-1,-2) (pl) {$\betas f$};
        \node at (1,-2) (pr) {$\betas\Id$};
        \node at (-1,-3) (bl) {};
        \node at (0,-3) (bc) {};
        \node at (1,-3) (br) {};
        \draw[dmoda] (tc) to (delta);
        \draw[dmoda] (delta) to (bc);
        \draw[moda] (tl) to (pl);
        \draw[moda] (pl) to (bl);
        \draw[moda] (tr) to (pr);
        \draw[taa] (delta) to (pl);
        \draw[taa] (delta) to (pr);
        \draw[moda] (pr) to (br);
      \end{tikzpicture}}.
  \end{align*}
  The second term vanishes because $d(\Id)=0$. The third, fourth, and fifth
  terms cancel by the type \DD\ structure relation. (The cancellation happens
  after factoring out $X_{12}^v$ from the fifth term. The third and fourth terms
  are like special cases where $v=0$, if one views $\wTrDiag{1}{0}$ as
  $\mu^{1,0}\otimes\Id+\Id\otimes \mu^{1,0}$.)  The last term vanishes because
  the identity map $\Id$ has $\Id_n=0$ for $n>1$. So, applying
  Lemma~\ref{lem:w-F-f} to the second-to-last term we are left with
  \[
    \sum_w
    X_{12}^w
  \tikzsetnextfilename{wpmm-4-1}
    \mathcenter{\begin{tikzpicture}[smallpic]
        \node at (0,0) (tc) {};
        \node at (-1,0) (tl) {};
        \node at (1,0) (tr) {};
        \node at (0,-1) (delta) {$\delta$};
        \node at (-1,-2) (pl) {$\wPartTrModMulDiag{*}{w}(df)$};
        \node at (1,-2) (pr) {$\wPartTrModMulDiag{*}{w}\Id$};
        \node at (-1,-3) (bl) {};
        \node at (0,-3) (bc) {};
        \node at (1,-3) (br) {};
        \draw[dmoda] (tc) to (delta);
        \draw[dmoda] (delta) to (bc);
        \draw[moda] (tl) to (pl);
        \draw[moda] (pl) to (bl);
        \draw[moda] (tr) to (pr);
        \draw[taa] (delta) to (pl);
        \draw[taa] (delta) to (pr);
        \draw[moda] (pr) to (br);
      \end{tikzpicture}}+
  \tikzsetnextfilename{wpmm-4-2}
    \mathcenter{\begin{tikzpicture}[smallpic]
        \node at (0,0) (tc) {};
        \node at (-1,0) (tl) {};
        \node at (1,0) (tr) {};
        \node at (0,-1) (delta) {$\delta$};
        \node at (-1,-2) (pl) {$\alphas f$};
        \node at (1,-2) (pr) {$\alphas(d\Id)$};
        \node at (-1,-3) (bl) {};
        \node at (0,-3) (bc) {};
        \node at (1,-3) (br) {};
        \draw[dmoda] (tc) to (delta);
        \draw[dmoda] (delta) to (bc);
        \draw[moda] (tl) to (pl);
        \draw[moda] (pl) to (bl);
        \draw[moda] (tr) to (pr);
        \draw[taa] (delta) to (pl);
        \draw[taa] (delta) to (pr);
        \draw[moda] (pr) to (br);
      \end{tikzpicture}}+
  \tikzsetnextfilename{wpmm-4-3}
    \mathcenter{
      \begin{tikzpicture}[smallpic]
        \node at (0,0) (tc) {};
        \node at (-2,0) (tl) {};
        \node at (2,0) (tr) {};
        \node at (0,-1) (delta1) {$\delta$};
        \node at (0,-2) (delta2) {$\delta^1$};
        \node at (0,-3) (delta3) {$\delta$};
        \node at (-2,-4) (kl) {$\alphas f$};
        \node at (2,-4) (kr) {$\alphas \Id$};
        \node at (1,-3) (gammar) {$\mu^{1,0}$};
        \node at (-2,-5) (bl) {};
        \node at (2,-5) (br) {};
        \node at (0,-5) (bc) {};
        \draw[dmoda] (tc) to (delta1);
        \draw[dmoda] (delta1) to (delta2);
        \draw[dmoda] (delta2) to (delta3);
        \draw[dmoda] (delta3) to (bc);
        \draw[taa, bend right=10] (delta1) to (kl);
        \draw[taa, bend left=10] (delta1) to (kr);
        \draw[alga] (delta2) to (kl);
        \draw[alga] (delta2) to (gammar);
        \draw[taa, bend left=5] (delta3) to (kl);
        \draw[taa, bend right=5] (delta3) to (kr);
        \draw[alga] (gammar) to (kr);
        \draw[moda] (tl) to (kl);
        \draw[moda] (tr) to (kr);
        \draw[moda] (kl) to (bl);
        \draw[moda] (kr) to (br);
      \end{tikzpicture}}.
  \]
  Again, the second term vanishes because $d(\Id)=0$ and the third term vanishes
  because $\Id_n=0$ for $n\geq 2$. Thus, we are left $[(df)\DT\Id_P\DT \Id_{\cNodule}]_{\PartTrModMulDiag}$, as desired.
  
  For the second statement, applying the argument above but with a homotopy $\zeta$ of partial module-map diagonals in place of $\PartTrModMulDiag$ gives
  \[
    \sum_w
    X_{12}^w
    d\left(
  \tikzsetnextfilename{wpmm-5-1}      
      \mathcenter{\begin{tikzpicture}[smallpic]
          \node at (0,0) (tc) {};
          \node at (-1,0) (tl) {};
          \node at (1,0) (tr) {};
          \node at (0,-1) (delta) {$\delta$};
          \node at (-1,-2) (pl) {$\zeta^{*,w} f$};
          \node at (1,-2) (pr) {$\zeta^{*,w}\Id$};
          \node at (-1,-3) (bl) {};
          \node at (0,-3) (bc) {};
          \node at (1,-3) (br) {};
          \draw[dmoda] (tc) to (delta);
          \draw[dmoda] (delta) to (bc);
          \draw[moda] (tl) to (pl);
          \draw[moda] (pl) to (bl);
          \draw[moda] (tr) to (pr);
          \draw[taa] (delta) to (pl);
          \draw[taa] (delta) to (pr);
          \draw[moda] (pr) to (br);
        \end{tikzpicture}}
    \right)
    =
    \sum_w
    X_{12}^w
  \tikzsetnextfilename{wpmm-5-2}      
    \left(  \mathcenter{\begin{tikzpicture}[smallpic]
          \node at (0,0) (tc) {};
          \node at (-1,0) (tl) {};
          \node at (1,0) (tr) {};
          \node at (0,-1) (delta) {$\delta$};
          \node at (-1,-2) (pl) {$\zeta^{*,w} (df)$};
          \node at (1,-2) (pr) {$\zeta^{*,w}\Id$};
          \node at (-1,-3) (bl) {};
          \node at (0,-3) (bc) {};
          \node at (1,-3) (br) {};
          \draw[dmoda] (tc) to (delta);
          \draw[dmoda] (delta) to (bc);
          \draw[moda] (tl) to (pl);
          \draw[moda] (pl) to (bl);
          \draw[moda] (tr) to (pr);
          \draw[taa] (delta) to (pl);
          \draw[taa] (delta) to (pr);
          \draw[moda] (pr) to (br);
        \end{tikzpicture}}
    \right)
    +[f\DT\Id_P\DT \Id_{\cNodule}]_{\PartTrModMulDiag^1}-[f\DT\Id_P\DT \Id_{\cNodule}]_{\PartTrModMulDiag^2},
  \]
  as desired.

  Finally, if $\PartTrModMulDiag$ comes from a primitive $\TrPMorDiag$ (and
  $\TrMDiag^i$ comes from a primitive $\TrPMDiag^i$) then
  \[
    [f\DT\Id_P\DT \Id_{\cNodule}]_{\PartTrModMulDiag}=
    \sum_{k,\ell,u_1,\dots,u_\ell,v}
    X_{12}^{u_1+\cdots+u_\ell+v}Y_1^v
    \mathcenter{
  \tikzsetnextfilename{wpmm-6-1}
      \begin{tikzpicture}[smallpic]
        \node at (-1,1) (tl) {};
        \node at (-1,-1) (p1l) {$\wTrPMDiag{*}{u_1}$};
        \node at (-1,-2) (p1ldots) {$\vdots$};
        \node at (-1,-3) (p2l) {$\wTrPMDiag{*}{u_{k-1}}$};
        \node at (-1,-4) (ql) {$\wTrPMorDiag{*}{u_{k}} f$};
        \node at (-1,-5) (p3l) {$\wTrPMDiag{*}{u_{k+1}}$};
        \node at (-1,-6) (p2ldots) {$\vdots$};
        \node at (-1,-7) (p4l) {$\wTrPMDiag{*}{u_\ell}$};
        \node at (-1,-10) (bl) {};
        \node at (0,1) (tc) {};
        \node at (0,0) (delta1) {$\delta$};
        \node at (0,-1) (deltadots1) {$\vdots$};
        \node at (0,-2) (delta2) {$\delta$};
        \node at (0,-3) (delta3) {$\delta$};
        \node at (0,-4) (delta4) {$\delta$};
        \node at (0,-5) (deltadots2) {$\vdots$};
        \node at (0,-6) (delta5) {$\delta$};
        \node at (0,-10) (bc) {};
        \node at (1,-1) (p1r) {$\wTrPMDiag{*}{u_1}$};
        \node at (1,-2) (p1rdots) {$\vdots$};
        \node at (1,-3) (p2r) {$\wTrPMDiag{*}{u_{k-1}}$};
        \node at (1,-4) (qr) {$\wTrPMorDiag{*}{u_{k}}$};
        \node at (1,-5) (p3r) {$\wTrPMDiag{*}{u_{k+1}}$};
        \node at (1,-6) (p2rdots) {$\vdots$};
        \node at (1,-7) (p4r) {$\wTrPMDiag{*}{u_\ell}$};
        \node at (3,1) (tr) {};
        \node at (3,-8) (mr) {$m^{\ell,v}$};
        \node at (3,-9) (Idr) {$\Id$};
        \node at (3,-10) (br) {};
        \draw[moda] (tl) to (p1l);
        \draw[moda] (p1l) to (p1ldots);
        \draw[moda] (p1ldots) to (p2l);
        \draw[moda] (p2l) to (ql);
        \draw[moda] (ql) to (p3l);
        \draw[moda] (p3l) to (p2ldots);
        \draw[moda] (p2ldots) to (p4l);
        \draw[moda] (p4l) to (bl);
        \draw[dmoda] (tc) to (delta1);
        \draw[dmoda] (delta1) to (deltadots1);
        \draw[dmoda] (deltadots1) to (delta2);
        \draw[dmoda] (delta2) to (delta3);
        \draw[dmoda] (delta3) to (delta4);
        \draw[dmoda] (delta4) to (deltadots2);
        \draw[dmoda] (deltadots2) to (delta5);
        \draw[dmoda] (delta5) to (bc);
        \draw[taa] (delta1) to (p1l);
        \draw[taa] (delta2) to (p2l);
        \draw[taa] (delta3) to (ql);
        \draw[taa] (delta4) to (p3l);
        \draw[taa] (delta1) to (p1r);
        \draw[taa] (delta2) to (p2r);
        \draw[taa] (delta3) to (qr);
        \draw[taa] (delta4) to (p3r);
        \draw[taa] (delta5) to (p4l);
        \draw[taa] (delta5) to (p4r);
        \draw[alga, bend left=15] (p1r) to (mr);
        \draw[alga, bend left=12] (p2r) to (mr);
        \draw[alga, bend left=6] (p3r) to (mr);
        \draw[alga, bend left=3] (p4r) to (mr);
        \draw[alga, bend left=9] (qr) to (mr);
        \draw[moda] (tr) to (mr);
        \draw[moda] (mr) to (Idr);
        \draw[moda] (Idr) to (br);
      \end{tikzpicture}
    }=(f\DT_{\TrPMorDiag}\Id_P)\DT \Id_N,
  \]
  as claimed. (This uses the fact that $X_2=X_{12}Y_1$.)

  The boundedness assumptions imply that all the sums considered above
  are finite. 
\end{proof}

\begin{corollary}\label{cor:w-morph-DT-assoc}
  Fix a module-map primitive $\TrPMorDiag$ compatible with module diagonal primitives $\TrPMDiag^1$ and $\TrPMDiag^2$. Let $\TrMDiag^i$ be the module diagonal induced by $\TrPMDiag^i$ and let $\TrModMulDiag$ be any module-map diagonal compatible with $\TrMDiag^1$ and $\TrMDiag^2$. Then for any $\Ainf$-modules $\wMod_1$, $\wMod_2$, and $\wNod$, type \DD\ structure $P$, and $\Ainf$-module homomorphism $f\co \wMod_1\to\wMod_2$ (with either $P$ operationally bounded or all the other objects bonsai), the chain maps
  \[
    [f\DT \Id_P\DT \Id_{\wNod}]_{\TrModMulDiag}
  \]
  and
  \[
    (f\DT_{\TrPMorDiag}\Id_P)\DT \Id_{\wNod}
  \]
  are chain homotopic.
\end{corollary}
\begin{proof}
  As in the unweighted case, this is immediate from
  Lemma~\ref{lem:w-part-mod-map-diag-DT} and the fact that all partial
  weighted module-map diagonals are homotopic
  (Lemma~\ref{lem:w-part-mod-maps-htpic}).
\end{proof}

\subsection{Tensor products of homotopy unital weighted algebras and associativity of the tensor product}\label{sec:w-htpy-unital}

Homotopy unital diagonals allow us to define various tensor products:
\begin{itemize}
\item Given a homotopy unital algebra diagonal $\uwDiag{*}{*}$ and
  homotopy unital $w$-algebras $\wAlg$ and $\wBlg$ over
  $\Ground_1$ and $\Ground_2$, respectively, the tensor
  product $A\rotimes{\Ring} B$ inherits the structure of a homotopy unital
  $w$-algebra $\wAlg\wADtp[\uwDiagNS]\wBlg$ over $\Ground$
  via the composition
  \begin{align*}
    \uwTreesCx{n}{w}&\xrightarrow{\uwDiag{*}{*}}\bigoplus_{w_1+w_2\leq w}\uwTreesCx{n}{w_1}\rotimes{\Ring}\uwTreesCx{n}{w_2}\rotimes{\Ring}\Ring[Y_1,Y_2]\\
                    &\xrightarrow{\mu_{\wAlg}\otimes\mu_{\wBlg}}
                      \Mor(A^{\kotimes{\Ground_1} n},A)\rotimes{\Ring}\Mor(B^{\kotimes{\Ground_2} n},B)\to \Mor((A\kotimes{\Ground} B)^{\kotimes{\Ground} n},A\rotimes{\Ring} B\grs{(4-\kappa_1-\kappa_2)w}).
  \end{align*}
\item Similarly, given a homotopy unital module diagonal
  $\uwMDiag{*}{*}$ compatible with $\uwDiag{*}{*}$, and homotopy
  unital $w$-modules $\wMod$ and $\wNod$ over $\wAlg$ and $\wBlg$, the
  tensor product $M\rotimes{\Ring} N$ inherits the structure of a homotopy
  unital $w$-module $\wMod\wADtp[\uwMDiagNS]\wNod$ over
  $\wAlg\wADtp[\uwDiagNS]\wBlg$ via the composition
    \begin{align*}
    \uwMTreesCx{1+n}{w}&\xrightarrow{\uwMDiag{*}{*}}\bigoplus_{w_1+w_2\leq w}\uwMTreesCx{1+n}{w_1}\rotimes{\Ring}\uwMTreesCx{1+n}{w_2}\rotimes{\Ring}\Ring[Y_1,Y_2]\\
                    &\xrightarrow{m_{\wMod}\otimes m_{\wNod}}
                      \Mor((M\rotimes{\Ring} N)\kotimes{\Ground} (A\rotimes{\Ring} B)^{\kotimes{\Ground} n},M\rotimes N\grs{(4-\kappa_1-\kappa_2)w}).
  \end{align*}
\item Similarly, given a homotopy unital map diagonal
  $\uwMulDiag{n}{w}$ compatible with $\uwDiag{*}{*}_1$ and
  $\uwDiag{*}{*}_2$, as well as homotopy unital algebra maps
  $f\co \wAlg_1\to\wAlg_2$ and $g\co \wBlg_1\to\wBlg_2$ there is an
  induced tensor product
  \[
    f\wADtp[\uwMulDiagNS]g\co \wAlg_1\wADtp[\uwDiagNS_1]\wBlg_1\to\wAlg_2\wADtp[\uwDiagNS_2]\wBlg_2.
  \]
\item Similarly, given a homotopy unital module-map diagonal
  $\uwModMulDiag{n}{w}$ compatible with $\uwMDiag{*}{*}_1$ and
  $\uwMDiag{*}{*}_2$, as well as homotopy unital module maps
  $f\co \wMod_1\to\wMod_2$ and $g\co \wNod_1\to\wNod_2$ there is an
  induced tensor product
  \[
    f\wMDtp[\uwModMulDiagNS]g\co \wMod_1\wADtp[\uwMDiagNS_1]\wNod_1\to\wMod_2\wMDtp[\uwMDiagNS_2]\wNod_2.
  \]
\end{itemize}

These homotopy unital tensor product have similar properties to the
cases spelled out above:
\begin{theorem}\label{thm:hu-tens-prod-props}
  These tensor products have the following properties:
  \begin{enumerate}[label=(\arabic*)]
  \item\label{item:hu-tens-alg} If $f\co\wAlg_1\to\wAlg_2$ and $g\co\wBlg_1\to\wBlg_2$ are
    quasi-isomorphisms then
    \[
      f\wADtp[\uwMulDiagNS]g\co \wAlg_1\wADtp[\uwDiagNS_1]\wBlg_1\to\wAlg_2\wADtp[\uwDiagNS_2]\wBlg_2.
    \]
    is a quasi-isomorphism.  In particular:
    \begin{enumerate}
    \item If $\wAlg_1$ is quasi-isomorphic to $\wBlg_1$ and $\wAlg_2$ is
      quasi-isomorphic to $\wBlg_2$, then
      $\wAlg_1\wADtp[\uwDiagNS]\wAlg_2$ is quasi-isomorphic to
      $\wBlg_1\wADtp[\uwDiagNS]\wBlg_2$.
    \item If $\uwDiag{*}{*}_1$ and $\uwDiag{*}{*}_2$ are homotopy unital
      weighted algebra diagonals with the same seed, then
      $\wAlg_1\wADtp[\uwDiagNS_1]\wAlg_2$ is isomorphic to
      $\wAlg_1\wADtp[\uwDiagNS_2]\wAlg_2$.
    \end{enumerate}
  \item\label{item:hu-tens-forget} If $\wAlg$ and $\wBlg$ are strictly unital weighted algebras
    and $\uwDiag{*}{*}$ is a homotopy unital algebra diagonal
    extending $\wDiag{*}{*}$ then the (weakly unital) weighted algebra
    $\wAlg\wADtp[\wDiagNS]\wBlg$ is the image of the homotopy unital
    weighted algebra $\wAlg\wADtp[\uwDiagNS]\wBlg$ under the forgetful map.
    Similar statements hold for modules, algebra maps, and module maps.
  \item\label{item:hu-tens-mod-morph} Given homotopy unital module morphisms
    \[
      \wMod^1\stackrel{f_1}{\longrightarrow}\wMod^2\stackrel{f_2}{\longrightarrow}\wMod^3\qquad\qquad
      \wNod^1\stackrel{g_1}{\longrightarrow}\wNod^2\stackrel{g_2}{\longrightarrow}\wNod^3,
    \]
    we have
    \begin{align}
      (f_2\circ f_1)\MDtp[\uwModMulDiagNS_{13}](g_2\circ g_1)&\sim(f_2\MDtp[\uwModMulDiagNS_{23}] g_2)\circ (f_1\MDtp[\uwModMulDiagNS_{12}] g_1)\label{eq:u-w-mod-map-tensor-comp}\\
      d(f_1\MDtp[\uwModMulDiagNS_{12}]g_1)&=(df_1)\MDtp[\uwModMulDiagNS_{12}]g_1+f_1\MDtp[\uwModMulDiagNS_{12}](dg_1).\label{eq:u-w-mod-map-tensor-diff}
    \end{align}
  \item\label{item:hu-tens-htpic-mod-morph} If $\uwModMulDiag{*}{*}_1$ and $\uwModMulDiag{*}{*}_2$ are
    homotopic homotopy unital module-map diagonals then
    $f\MDtp[\uwModMulDiagNS_{12}]g$ and
    $f\MDtp[\uwModMulDiagNS_{12}]g$ are homotopic homotopy unital
    module-maps.
  \item\label{item:hu-tens-id} If $\wMod$ and $\wNod$ are homotopy unital weighted modules then the tensor
    product of their identity maps
    \[
      \Id_{\wMod}\wMDtp[\uwModMulDiagNS]\Id_{\wNod}\co
      \wMod\wADtp[\uwMDiagNS_1]\wNod\to\wMod\wMDtp[\uwMDiagNS_2]\wNod
    \]
    is an isomorphism. Consequently, the tensor product of two
    homotopy equivalences is a homotopy equivalence.
  \end{enumerate}
\end{theorem}
\begin{proof}
  The proof Part~\ref{item:hu-tens-alg} is analogous to the proof of
  Theorem~\ref{thm:Algebras}
  or~\ref{thm:wAlgebras}. Part~\ref{item:hu-tens-forget} is immediate
  from the definitions. The proofs of
  Parts~\ref{item:hu-tens-mod-morph}
  and~\ref{item:hu-tens-htpic-mod-morph} are similar to the proof of
  Proposition~\ref{prop:w-dg-bifunctor} (or
  Proposition~\ref{prop:dg-bifunctor}) and
  Lemma~\ref{lem:tens-mod-maps}.  The proof of
  Part~\ref{item:hu-tens-id} is similar to the proof of
  Lemma~\ref{lem:w-mod-id-tens-id}.
\end{proof}

\begin{theorem}\label{thm:hu-tens-assoc}
  Let $\wAlg_1$, $\wAlg_2$, and $\wAlg_3$ be homotopy unital
  $w$-algebras over $\Ground_1$, $\Ground_2$ and
  $\Ground_3$, respectively. Fix homotopy unital algebra
  diagonals $\uwDiag{*}{*}_1$ and $\uwDiag{*}{*}_2$. View
  $\wAlg_1\wADtp[\uwDiagNS_1]\wAlg_2$ as a homotopy unital weighted
  algebra over $\Ground_{12}\coloneqq\Ground_1\rotimes{\Ring}\Ground_2$, and
  view $\wAlg_2\wADtp[\uwDiagNS_2]\wAlg_3$ as a homotopy unital
  weighted algebra over
  $\Ground_{23}\coloneqq\Ground_2\rotimes{\Ring}\Ground_3$. View
  $\Ground_{12}$ as an
  $\Ring[Y_1]$-algebra by sending $Y_1\mapsto Y_{12}\coloneqq Y_1Y_2$, and
  view $\Ground_{23}$ as an $\Ring[Y_2]$-algebra by sending $Y_2\mapsto
  Y_{23}\coloneqq Y_2Y_3$.
  Then there is a homotopy unital isomorphism
  \[
    (\wAlg_1\wADtp[\uwDiagNS_1]\wAlg_2)\wADtp[\uwDiagNS_2]\wAlg_3
    \cong
    \wAlg_1\wADtp[\uwDiagNS_1](\wAlg_2\wADtp[\uwDiagNS_2]\wAlg_3)
  \]
  over $\Ground_1\rotimes{\Ring}\Ground_2\rotimes{\Ring}\Ground_3$.
\end{theorem}
\begin{proof}
  The proof is similar to the unweighted case,
  point~\ref{item:Alg-thm-assoc} in Theorem~\ref{thm:Algebras}. We
  define a \emph{homotopy unital algebra double diagonal} to be a
  collection of chain maps
  \[
    \uwDiag{n}{w}\co\uwTreesCx{n}{w}\to \bigoplus_{w_1,w_2,w_3\leq
      w}\uwTreesCx{n}{w_1}\rotimes{\Ring}\uwTreesCx{n}{w_2}\rotimes{\Ring}\uwTreesCx{n}{w_1}\rotimes{\Ring}\Ring[Y_1,Y_2,Y_3]
  \]
  which respect composition and so that
  \begin{align*}
    \uwDiag{0}{0}(\stump)&=\stump\otimes\stump\otimes\stump\\
    \uwDiag{1}{0}(\IdTree)&=\IdTree\otimes\IdTree\otimes\IdTree\\
    \uwDiag{2}{0}(\wcorolla{2}{0})&=\wcorolla{2}{0}\otimes
                                              \wcorolla{2}{0}\otimes
                                              \wcorolla{2}{0}.
  \end{align*}
  The \emph{seed} of $\uwDiag{*}{*}$ is $\uwDiag{0}{1}$. 
  Define a homotopy unital algebra map double diagonal similarly,
  using the transformation thorn trees complex. Then:
  \begin{enumerate}[label=(\arabic*)]
  \item\label{item:hu-doub-diag-1} Both
    $(\Id\otimes
    \uwDiag{*}{*}_2|_{(Y_1,Y_2)\to(Y_2,Y_3)})\circ(\uwDiag{*}{*}_1|_{Y_2\to
      Y_2Y_3})$ and
    $(\uwDiag{*}{*}_1\otimes \Id)\circ
    (\uwDiag{*}{*}_2|_{(Y_1,Y_2)\to(Y_1Y_2,Y_3)})$ are homotopy unital
    algebra double diagonals.  Here, the
    subscripts indicate variable substitutions; for example,
    $(Y_1,Y_2)\to (Y_2,Y_3)$ indicates the substitution replacing
    $Y_1$ by $Y_2$ and $Y_2$ by $Y_3$.
  \item A homotopy unital algebra double diagonal allows one to define
    a triple tensor product $\wAlg_1\rotimes{\Ring}\wAlg_2\rotimes{\Ring}\wAlg_3$ and
    a homotopy unital algebra map double diagonal allows one to define
    a triple tensor product of maps $f_1\otimes f_2\otimes f_3$.
  \item A triple tensor product of quasi-isomorphisms is a
    quasi-isomorphism and a triple tensor product of isomorphisms is
    an isomorphism. (Both statements use the nondegeneracy condition
    for the diagonals, and the second statement uses Lemma~\ref{lem:hu-Alg-iso-is}.)
  \item Any two homotopy algebra double diagonals with the same seed are related by a
    homotopy unital algebra map double diagonal.
  \item The two homotopy unital double diagonals in point~\ref{item:hu-doub-diag-1} have the same seed.
  \end{enumerate}
  The result follows.
\end{proof}

\begin{corollary}\label{cor:w-tens-assoc}
  Let $\wAlg$, $\wBlg$, and $\wClg$ be split unital weighted
  algebras over $\Ground_1$, $\Ground_2$, and
  $\Ground_3$, and
  $\wDiag{*}{*}_1$ and $\wDiag{*}{*}_2$ weighted algebra diagonals. Let
  $\overline{\wAlg\wADtp[\wDiagNS_1]\wBlg}$ be the result of tensoring together
  $\wAlg$ and $\wBlg$ using $\wDiag{*}{*}_1$ as in Definition~\ref{def:wAlg-tp}
  and then making the result strictly unital as in
  Theorem~\ref{thm:UnitalIsUnitalW}. Similarly, let
  $\overline{\wBlg\wADtp[\wDiagNS_2]\wClg}$ be the result of tensoring together
  $\wBlg$ and $\wClg$ using $\wDiag{*}{*}_2$
  and then making the result strictly unital. View
  $\Ground_{12}$ as an
  $\Ring[Y_1]$-algebra by sending $Y_1\mapsto Y_{12}\coloneqq Y_1Y_2$, and
  view $\Ground_{23}$ as an $\Ring[Y_2]$-algebra by sending $Y_2\mapsto
  Y_{23}\coloneqq Y_2Y_3$.
  Then there is an isomorphism
  \[
    \overline{\wAlg\wADtp[\wDiagNS_1]\wBlg}\wADtp[\wDiagNS_2]\wClg\cong
    \wAlg\wADtp[\wDiagNS_1]\overline{\wBlg\wADtp[\wDiagNS_2]\wClg}.
  \]
\end{corollary}
\begin{proof}
  Fix homotopy unital weighted algebra diagonals $\uwDiag{*}{*}_i$
  extending $\wDiag{*}{*}_i$. Let $\Forget$ be the forgetful functor
  from homotopy unital algebras to ordinary algebras. By
  Theorem~\ref{thm:hu-tens-prod-props},
  \begin{align*}
    \Forget(\wAlg\wADtp[\uwDiagNS_1]\wBlg)&\cong \wAlg\wADtp[\wDiagNS_1]\wBlg\\
    \Forget(\wBlg\wADtp[\uwDiagNS_2]\wClg)&\cong\wBlg\wADtp[\wDiagNS_2]\wClg.
  \end{align*}
  Let $\overline{\wAlg\wADtp[\uwDiagNS_1]\wBlg}$ and
  $\overline{\wBlg\wADtp[\uwDiagNS_2]\wClg}$ be strictly unital
  weighted algebras which are homotopy unitally isomorphic to
  $\wAlg\wADtp[\uwDiagNS_1]\wBlg$ and $\wBlg\wADtp[\uwDiagNS_2]\wClg$,
  respectively, the existence of which is guaranteed by
  Theorem~\ref{thm:hu-rectify}. Then
  \begin{align*}
    \Forget(\overline{\wAlg\wADtp[\uwDiagNS_1]\wBlg})&\cong \overline{\wAlg\wADtp[\wDiagNS_1]\wBlg}\\
    \Forget(\overline{\wBlg\wADtp[\uwDiagNS_2]\wClg})&\cong\overline{\wBlg\wADtp[\wDiagNS_2]\wClg}.
  \end{align*}
  By Theorems~\ref{thm:hu-tens-prod-props} and~\ref{thm:hu-tens-assoc},
  \begin{align*}
    (\overline{\wAlg\wADtp[\uwDiagNS_1]\wBlg})\wADtp[\uwDiagNS_2]\wClg
    &\cong
      (\wAlg\wADtp[\uwDiagNS_1]\wBlg)\wADtp[\uwDiagNS_2]\wClg\\
    &\cong \wAlg\wADtp[\uwDiagNS_1](\wBlg\wADtp[\uwDiagNS_2]\wClg)\\
    &\cong
    \wAlg\wADtp[\uwDiagNS_1]\overline{\wBlg\wADtp[\uwDiagNS_2]\wClg}.
  \end{align*}
  Applying the forgetful functor and
  Theorem~\ref{thm:hu-tens-prod-props} again gives the result.
\end{proof}

\begin{remark}
  If we think of $Y_i$ as associated to the algebra $\wAlg_i$, then
  Corollary~\ref{cor:w-tens-assoc} says that the $Y$-variable associated to
  $\wAlg_1\wADtp\wAlg_2$ should be $Y_{12}=Y_1Y_2$. Similarly, a category of
  weighted type $D$ structures is associated to a pair $(\wAlg_i,X_i)$, and
  Definition~\ref{def:w-one-sided-DT} implies that these charges should satisfy
  $X_2=X_{12}Y_1$ (and similarly $X_1=X_{12}Y_2$). (That is, the $1$-sided box
  tensor product relates the category of type \DD\ structures associated to
  $(\wAlg_1\wADtp\wAlg_2,X_{12})$, the category of $w$-modules over
  $(\wAlg_1,Y_1)$, and the category of type $D$ structures over
  $(\wAlg_2,X_{12}Y_1)$.) In particular, this suggests considering variables satisfying
  \begin{equation}\label{eq:XY-relations}
    X_1Y_1 = X_2Y_2 = X_{12}Y_{12}.
  \end{equation}
  Recall that the grading of $X_i$ is $-\kappa_i$, and of $Y_i$ is
  $\kappa_i-2$. So, all terms in Equation~\eqref{eq:XY-relations} have grading
  $-2$. (In bordered Floer theory, $X_iY_i$ will be the variable $U$, which
  corresponds in monopole Floer theory to a generator of $H^2(\CC P^\infty)$.)

  On a related point, if $\lsup{\wAlg_1}P$ and $\lsup{\wAlg_2}Q$ are type $D$
  structures then one might try to define an operation $\delta^1$ on
  $P\rotimes{\Ring} Q$, making $P\rotimes{\Ring} Q$ into a type \DD\ structure, by
  $\delta^1=\delta^1_P\otimes \unit + \unit\otimes\delta^1_Q$. It is not hard to
  see that the type \DD\ structure relation for this type \DD\ structure would
  require the charge to be $X_2/Y_1=X_1/Y_2$. (Even then, it is not apparent
  that this operation $\delta^1$ satisfies the type \DD\ structure relation for
  a general diagonal; a correct, general definition of the external tensor
  product of two type $D$ structures, when variables are chosen as above, can be
  given using a DADD diagonal, as in Section~\ref{sec:DA-tens-DD}.)
\end{remark}

\subsection{The category of type \textalt{$D$}{D} structures  over a homotopy
  unital algebra and independence of the category of type \textalt{\DD}{DD}
  structures from the algebra diagonal}\label{sec:hu-D}

We now show that the $\Ainf$-category of type $D$ structures over a
homotopy unital weighted $\Ainf$-algebra is itself an (unweighted)
homotopy unital category. We start with some definitions.

\begin{definition}\label{def:hu-D}
  Let $\wAlg$ be a homotopy unital weighted $\Ainf$-algebra over
  $\Ground$ and let $X\in\Ground$ be an element acting centrally on
  every module under consideration. Recall
  that there is an underlying weighted $\Ainf$-algebra $\wAlg'$ of
  $\wAlg$. A \emph{type $D$ structure over $\wAlg$ with charge $X$} is simply a type
  $D$ structure over $\wAlg'$ with charge $X$  (Definition~\ref{def:wD}).
\end{definition}

\begin{remark}
  Implicitly, in Definition~\ref{def:hu-D}, we are assuming that either
  $\wAlg'$ is bonsai or the type $D$ structure is bounded.
\end{remark}

\begin{definition}
  The \emph{unweighted thorn trees complex} $\uTreesCx{n}$ is the
  subcomplex $\uwTreesCx{n}{0}$ of the thorn trees complex
  (Definition~\ref{def:wu-trees-cx}).
\end{definition}

\begin{definition}
  A \emph{homotopy unital $\Ainf$-category} $\Cat$ consists of
  \begin{itemize}
  \item a collection of objects $\ob(\Cat)$;
  \item for each pair of objects $C_1$,
    $C_2$, a chain complex $\Mor_\Cat(C_0,C_1)$; and
  \item for each $n\ge 0$ and each sequence $C_0,\dots,C_n$ of
    objects of $\Cat$, a chain map
    \[
      \mu_\Cat \co \uTreesCx{n} \to
      \Mor\bigl(\Mor_\Cat(C_0,C_1)\otimes\cdots\otimes\Mor_\Cat(C_{n-1},C_n),
         \Mor_\Cat(C_0,C_n)\bigr)
    \]
  \end{itemize}
  so that $\mu_\Cat(S \circ_i T) = \mu_\Cat(S) \circ_i
  \mu_\Cat(T)$.
  (In the special case $n=0$, $\mu_\Cat\co \uTreesCx{n}\to \Mor(C_0,C_0)$.)
\end{definition}

To define a homotopy unital structure on the category of type~$D$
structures, it suffices to define the unit $\Id^1_P$ and the corollas
$\mu_{n_1\uparrow\cdots\uparrow n_k}$. We distinguish between the maps $\mu_\wAlg$ on
the underlying homotopy-unital $\Ainf$-algebra and $\mu_D$ for the
category of type~$D$ structures.

\begin{definition}\label{def:hu-D-id}
  Given a homotopy unital weighted algebra $\wAlg$ and a type $D$
  structure $\lsup{\wAlg}P$ over $\wAlg$ with charge $X$, with either
  $\wAlg$ bonsai or $\lsup{\wAlg}P$ operationally bounded, the
  \emph{identity map of $\lsup{\wAlg}P$} is defined by
  \[
    \Id^1(x)=\One\otimes x+\sum_{n_1,n_2,w}X^w(\mu_{\wAlg,n_1\uparrow n_2}^w\otimes\Id_P)\circ\delta^{n_1+n_2}(x)
  \]
  or, graphically,
  \[
    \Id^1=
  \tikzsetnextfilename{def-hu-D-id-1}
    \mathcenter{\begin{tikzpicture}[smallpic]
        \node at (.5,1) (tr) {};
        \node at (0,-1) (One) {}; 
        \node at (0,-3) (bl) {};
        \node at (.5,-3) (br) {};
        \draw[dmoda] (tr) to (br);
        \draw[stumpa] (One) to (bl);
      \end{tikzpicture}}
    \mathcenter{\qquad+\qquad}
    \mathcenter{
  \tikzsetnextfilename{def-hu-D-id-2}
      \begin{tikzpicture}[smallpic]
        \node at (1,0) (tr) {};
        \node at (1,-1) (delta1) {$\delta$};
        \node at (1,-2) (delta2) {$\delta$};
        \node at (1,-4) (br) {};
        \node at (0,-3) (mu) {$\mu_\wAlg^\bullet$};
        \node at (0,-4) (bl) {};
        \node at ($ (delta1)!.5!(delta2) $) (thead) {};
        \node at ($ (thead)!.35!(mu) $) (thead3) {};
        \draw[dmoda] (tr) to (delta1);
        \draw[dmoda] (delta1) to (delta2);
        \draw[dmoda] (delta2) to (br);
        \draw[taa,bend right=10] (delta1) to (mu);
        \draw[taa, bend left=10] (delta2) to (mu);
        \draw[alga] (mu) to (bl);
        \draw[->,thick,>=stealth] (mu) to (thead3);
      \end{tikzpicture}
    }.
  \]
  (In the formula, $\Id_P$ means the identity map of the vector space
  $P$. Also, recall that $\mu^\bullet=\sum_{n,w}X^w\mu_n^w$.)
\end{definition}

\begin{convention}\label{conv:boundedness-suppressed}
  To keep statements simpler, in the rest of this section we will
  suppress the boundedness hypotheses. In each case, either the
  weighted $\Ainf$-algebras and weighted $\Ainf$-algebra homomorphisms
  involved are be required to be bonsai or the weighted type $D$ and
  \DD\ structures are required to be operationally bounded.
\end{convention}

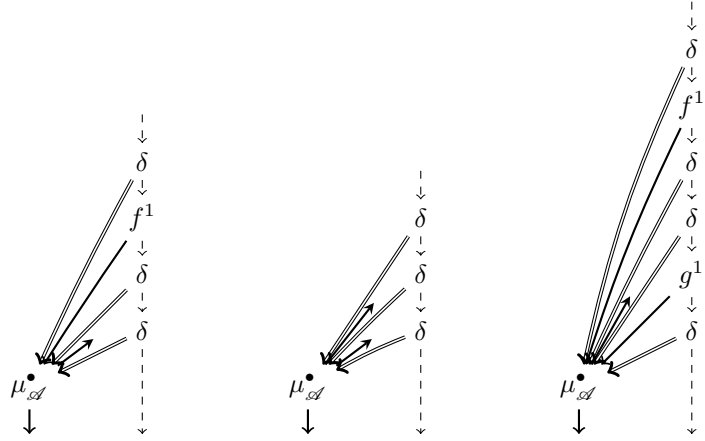
\begin{figure}
  \centering
  \[
    \tikzsetnextfilename{huD-example-1}
    \begin{tikzpicture}[smallpic]
      \node at (0,0) (tr) {};
      \node at (0,-1) (delta1) {$\delta$};
      \node at (0,-2) (f) {$f^1$};
      \node at (0,-3) (delta2) {$\delta$};
      \node at (0,-4) (delta3) {$\delta$};
      \node at (0,-6) (br) {};
      \node at (-2,-5) (mu) {$\mu_\wAlg^\bullet$};
      \node at (-2,-6) (bl) {};
      \node at ($ (delta2)!.5!(delta3) $) (d23avg) {};
      \node at ($ (d23avg)!.35!(mu) $) (thead1) {};
      \draw[dmoda] (tr) to (delta1);
      \draw[dmoda] (delta1) to (f);
      \draw[dmoda] (f) to (delta2);
      \draw[dmoda] (delta2) to (delta3);
      \draw[dmoda] (delta3) to (br);
      \draw[taa,bend right=2] (delta1) to (mu);
      \draw[alga,bend right=1] (f) to (mu);
      \draw[taa] (delta2) to (mu);
      \draw[taa] (delta3) to (mu);
      \draw[alga] (mu) to (bl);
      \draw[->,thick,>=stealth] (mu) to (thead1);
    \end{tikzpicture}\qquad\qquad
    \tikzsetnextfilename{huD-example-2}
    \begin{tikzpicture}[smallpic]
      \node at (0,0) (tr) {};
      \node at (0,-1) (delta1) {$\delta$};
      \node at (0,-2) (delta2) {$\delta$};
      \node at (0,-3) (delta3) {$\delta$};
      \node at (0,-5) (br) {};
      \node at (-2,-4) (mu) {$\mu_\wAlg^\bullet$};
      \node at (-2,-5) (bl) {};
      \node at ($ (delta2)!.5!(delta3) $) (d23avg) {};
      \node at ($ (d23avg)!.35!(mu) $) (thead2) {};
      \node at ($ (delta1)!.45!(delta2) $) (d12avg) {};
      \node at ($ (d12avg)!.35!(mu) $) (thead1) {};
      \draw[dmoda] (tr) to (delta1);
      \draw[dmoda] (delta1) to (delta2);
      \draw[dmoda] (delta2) to (delta3);
      \draw[dmoda] (delta3) to (br);
      \draw[taa] (delta1) to (mu);
      \draw[taa,bend right=2] (delta2) to (mu);
      \draw[taa,bend right=4] (delta3) to (mu);
      \draw[alga] (mu) to (bl);
      \draw[->,thick,>=stealth] (mu) to (thead1);
      \draw[->,thick,>=stealth] (mu) to (thead2);
    \end{tikzpicture}\qquad\qquad    
    \tikzsetnextfilename{huD-example-3}
    \begin{tikzpicture}[smallpic]
      \node at (0,0) (tr) {};
      \node at (0,-1) (delta1) {$\delta$};
      \node at (0,-2) (f) {$f^1$};
      \node at (0,-3) (delta2) {$\delta$};
      \node at (0,-4) (delta3) {$\delta$};
      \node at (0,-5) (g) {$g^1$};
      \node at (0,-6) (delta4) {$\delta$};
      \node at (0,-8) (br) {};
      \node at (-2,-7) (mu) {$\mu_\wAlg^\bullet$};
      \node at (-2,-8) (bl) {};
      \node at ($ (delta2)!.45!(delta3) $) (d23avg) {};
      \node at ($ (d23avg)!.5!(mu) $) (thead1) {};
      \draw[dmoda] (tr) to (delta1);
      \draw[dmoda] (delta1) to (f);
      \draw[dmoda] (f) to (delta2);
      \draw[dmoda] (delta2) to (delta3);
      \draw[dmoda] (delta3) to (g);
      \draw[dmoda] (g) to (delta4);
      \draw[dmoda] (delta4) to (br);
      \draw[taa,bend right=6] (delta1) to (mu);
      \draw[alga,bend right=4] (f) to (mu);
      \draw[taa,bend right=2] (delta2) to (mu);
      \draw[->,thick,>=stealth] (mu) to (thead1);
      \draw[taa] (delta3) to (mu);
      \draw[alga] (g) to (mu);
      \draw[taa] (delta4) to (mu);
      \draw[alga] (mu) to (bl);
    \end{tikzpicture}
  \]
  \caption[Examples of homotopy unital compositions on the category of type $D$ structures]{\textbf{Examples of homotopy unital compositions on the category of type $D$ structures.} Left: $\mu_D(\wcorolla{1\uparrow}{0})$. Center:
    $\mu_D(\wcorolla{\uparrow\uparrow}{0})$. Right: $\mu_D(\wcorolla{1\uparrow1}{0})$}
  \label{fig:hu-D-comps}
\end{figure}

\begin{definition}\label{def:hu-D-str-cat}
  Given a homotopy unital weight algebra $\wAlg$,
  a corolla $\wcorolla{n_1\uparrow \cdots \uparrow n_k}{w}$ with a
  total of $N = \sum_{i=1}^k n_i$ inputs, a sequence of $N+1$ type~$D$
  structures $\lsup{\wAlg}P_j$ with structure map~$\delta^1_j$ and
  charge~$X$, and type $D$ morphisms~$\phi_j$
  between them, define the interpretation of the corolla by
  interspersing an arbitrary
  number of $\delta_j^1$ inputs among the inputs to the corolla, as
  shown by example in Figure~\ref{fig:hu-D-comps}
  and defined formally by
  \begin{multline}\label{eq:huD-mu}
    \mu_D(\corolla{n_1\uparrow\cdots\uparrow n_k})(\phi_1,\dots,\phi_N)
    = \sum_w X^w \sum_{m_0,\dots,m_N=0}^\infty \,\,\sideset{}{'}\sum_{M_1,\dots,M_k=0}^\infty
    \Bigl(\bigl(\mu^w_{\wAlg,n_1+M_1\uparrow n_2+M_2\uparrow\cdots\uparrow n_k+M_k}\otimes\Id_{P_N}\bigr) \\
    \circ
    (\Id \otimes \delta^{m_N}_N) \circ (\Id \otimes \phi_N) \circ\dots\circ
    (\Id \otimes \delta^{m_1}_1) \circ (\Id \otimes \phi_1) \circ
    \delta^{m_0}_0\Bigr)
  \end{multline}
  where the sum over the $M_i$ is constrained to make the arrows
  properly interspersed: we require
  \[
    \sum_{i=1}^k M_i = \sum_{i=0}^N m_i
  \]
  and, for $1 \le j \le k$,
  \[
    \sum_{i=0}^{n_1+\dots+n_{j-1}} m_i \le \sum_{i=1}^j M_i
      \le \sum_{i=0}^{n_1+\dots+n_j} m_i.
  \]
\end{definition}

\begin{proposition}\label{prop:hu-D-is-hu-Ainf-cat}
  Fix a weighted $\Ainf$-algebra $\wAlg$ over $\Ground$ and an element $X\in\Ring$.
  The structure maps in Definition~\ref{def:hu-D-str-cat} make the
  category of type $D$ structures $\lsup{\wAlg,X}\ModCat$ over 
  $\wAlg$ with charge $X$, where $\Alg$ is a
  homotopy unital weighted algebra, into a homotopy unital
  $\Ainf$-category.
\end{proposition}
\begin{proof}
  We have to check that the map $\mu_D$ defined above is a chain map. We check
  two special cases before the general case. First, we check that the identity
  map of $P$ is a cycle in the morphism complex, i.e.,
  $\partial(\mu_D(\stump)) = \mu_D(\partial \stump) = 0$. We have
  \[
    \partial(\mu_D(\stump))=
    \mathcenter{
      \tikzsetnextfilename{prop-huD-id-1}
      \begin{tikzpicture}[smallpic]
        \node at (0,0) (tr) {};
        \node at (0,-1) (delta1) {$\delta$};
       \node at (-.25,-2.25) (one) {};
        \node at (0,-3) (delta2) {$\delta$};      
        \node at (0,-5) (br) {};
        \node at (-1.5,-4) (mu) {$\mu_A^\bullet$};
        \node at (-1.5,-5) (bl) {};
        \draw[dmoda] (tr) to (delta1);
        \draw[dmoda] (delta1) to (delta2);
        \draw[dmoda] (delta2) to (br);
        \draw[alga] (mu) to (bl);
        \draw[taa,bend right=2] (delta1) to (mu);
        \draw[taa] (delta2) to (mu);
        \draw[stumpa] (one) to (mu);
      \end{tikzpicture}}
    +
    \mathcenter{
      \tikzsetnextfilename{prop-huD-id-2}
      \begin{tikzpicture}[smallpic]
        \node at (0,0) (tr) {};
        \node at (0,-1) (delta1) {$\delta$};
        \node at (0,-2) (delta2) {$\delta$};
        \node at (0,-3.5) (delta3) {$\delta$};
        \node at (0,-4.5) (delta4) {$\delta$};
        \node at (-1,-4.25) (muin) {$\mu_A^\bullet$};
        \node at (-2,-5.5) (muout) {$\mu_A^\bullet$};
        \node at ($ (delta2)!.5!(delta3) $) (d23avg) {};
        \node at ($ (d23avg)!.35!(muin) $) (thead) {};
        \node at (0,-6.5) (br) {};
        \node at (-2,-6.5) (bl) {};
        \draw[dmoda] (tr) to (delta1);
        \draw[dmoda] (delta1) to (delta2);
        \draw[dmoda] (delta2) to (delta3);
        \draw[dmoda] (delta3) to (delta4);
        \draw[dmoda] (delta4) to (br);
        \draw[taa,bend right=4] (delta1) to (muout);
        \draw[taa,bend right=2] (delta2) to (muin);
        \draw[taa] (delta3) to (muin);
        \draw[->,thick,>=stealth] (muin) to (thead);
        \draw[taa] (delta4) to (muout);
        \draw[alga] (muin) to (muout);
        \draw[alga] (muout) to (bl);
      \end{tikzpicture}
    }
    =
    \mathcenter{
      \tikzsetnextfilename{prop-huD-id-3}
      \begin{tikzpicture}[smallpic]
        \node at (0,0) (tr) {};
        \node at (0,-1) (delta1) {$\delta$};
        \node at (0,-2.5) (delta2) {$\delta$};
        \node at (0,-3.5) (delta3) {$\delta$};
        \node at (0,-4.5) (delta4) {$\delta$};
        \node at (-1,-4.5) (muin) {$\mu_A^\bullet$};
        \node at (-2,-5.5) (muout) {$\mu_A^\bullet$};
        \node at ($ (delta1)!.3!(delta2) $) (d12avg) {};
        \node at ($ (d12avg)!.6!(muout) $) (thead) {};
        \node at (0,-6.5) (br) {};
        \node at (-2,-6.5) (bl) {};
        \draw[dmoda] (tr) to (delta1);
        \draw[dmoda] (delta1) to (delta2);
        \draw[dmoda] (delta2) to (delta3);
        \draw[dmoda] (delta3) to (delta4);
        \draw[dmoda] (delta4) to (br);
        \draw[taa,bend right=4] (delta1) to (muout);
        \draw[->,thick,>=stealth] (muout) to (thead);
        \draw[taa,bend right=2] (delta2) to (muout);
        \draw[taa] (delta3) to (muin);
        \draw[taa] (delta4) to (muout);
        \draw[alga] (muin) to (muout);
        \draw[alga] (muout) to (bl);
      \end{tikzpicture}
    } +
    \mathcenter{
      \tikzsetnextfilename{prop-huD-id-4}
      \begin{tikzpicture}[smallpic]
        \node at (0,0) (tr) {};
        \node at (0,-1) (delta1) {$\delta$};
        \node at (0,-2) (delta2) {$\delta$};
        \node at (0,-3) (delta3) {$\delta$};
        \node at (0,-4.5) (delta4) {$\delta$};
        \node at (-1,-3.75) (muin) {$\mu_A^\bullet$};
        \node at (-2,-5.5) (muout) {$\mu_A^\bullet$};
        \node at ($ (delta3)!.5!(delta4) $) (d34avg) {};
        \node at ($ (d34avg)!.5!(muout) $) (thead) {};
        \node at (0,-6.5) (br) {};
        \node at (-2,-6.5) (bl) {};
        \draw[dmoda] (tr) to (delta1);
        \draw[dmoda] (delta1) to (delta2);
        \draw[dmoda] (delta2) to (delta3);
        \draw[dmoda] (delta3) to (delta4);
        \draw[dmoda] (delta4) to (br);
        \draw[taa,bend right=4] (delta1) to (muout);
        \draw[taa] (delta2) to (muin);
        \draw[taa,bend right=2] (delta3) to (muout);
        \draw[->,thick,>=stealth] (muout) to (thead);
        \draw[taa] (delta4) to (muout);
        \draw[alga] (muin) to (muout);
        \draw[alga] (muout) to (bl);
      \end{tikzpicture}
    }
    =0,
  \]
  The second equality comes from the homotopy unital
  $\Ainf$-algebra relations.
  The third equality uses the structure relation for~$P$.

  It remains to prove that
  \[
    \partial(\mu_D(\corolla{n_1\uparrow\cdots\uparrow n_k})) =
    \mu_D(\partial\corolla{n_1\uparrow\cdots\uparrow n_k}).
  \]
  (Actually, we must prove this identity for any tree, but as usual it suffices
  to prove it for the basic building blocks, namely stumps (above) and corollas
  (here).) The cases $\corolla{\uparrow1}$ and $\corolla{1\uparrow}$ are
  slightly special, since the differential of a 3-valent vertex with one input a thorn has a term
  corresponding to erasing that thorn. We check that $\corolla{\uparrow1}$ is a
  chain map; the case of $\mu_D(\corolla{1\uparrow})$ is similar.  We have
  \begin{align*}
    (\partial(\mu_D(\corolla{\uparrow1})))(f)
    &=\partial((\mu_D(\corolla{\uparrow1}))(f))+(\mu_D(\corolla{\uparrow1}))(\partial f)\\
    &=
      \tikzsetnextfilename{prop-huD-id-new-1}
      \mathcenter{\begin{tikzpicture}[smallpic]
          \node at (0,0) (tc) {};
          \node at (0,-1) (delta1) {$\delta$};
          \node at (0,-2) (delta2) {$\delta$};
          \node at (0,-3) (delta3) {$\delta$};
          \node at (0,-4) (f) {$f^1$};
          \node at (0,-5) (delta4) {$\delta$};
          \node at (0,-6) (delta5) {$\delta$};
          \node at (0,-8) (bc) {};
          \node at (-1.5,-6) (mu1) {$\mu_\wAlg^\bullet$};
          \node at (-3,-7) (mu2) {$\mu_\wAlg^\bullet$};
          \node at (-3,-8) (bl) {};
          \node at ($ (delta2)!.3!(delta3) $) (d23avg) {};
          \node at ($ (d23avg)!.6!(mu1) $) (thead) {};
          \draw[dmoda] (tc) to (delta1);
          \draw[dmoda] (delta1) to (delta2);
          \draw[dmoda] (delta2) to (delta3);
          \draw[dmoda] (delta3) to (f);
          \draw[dmoda] (f) to (delta4);
          \draw[dmoda] (delta4) to (delta5);
          \draw[dmoda] (delta5) to (bc);
          \draw[taa, bend right=10] (delta1) to (mu2);
          \draw[taa, bend right=10] (delta2) to (mu1);
          \draw[taa, bend left=5] (delta3) to (mu1);
          \draw[taa] (delta4) to (mu1);
          \draw[taa] (delta5) to (mu2);
          \draw[alga] (f) to (mu1);
          \draw[alga] (mu1) to (mu2);
          \draw[alga] (mu2) to (bl);
        \draw[->,thick,>=stealth] (mu1) to (thead);
      \end{tikzpicture}}
      +
      \tikzsetnextfilename{prop-huD-id-new-2}
      \mathcenter{\begin{tikzpicture}[smallpic]
          \node at (0,0) (tc) {};
          \node at (0,-1) (delta1) {$\delta$};
          \node at (0,-2) (delta2) {$\delta$};
          \node at (0,-3) (delta3) {$\delta$};
          \node at (0,-4) (f) {$f^1$};
          \node at (0,-5) (delta4) {$\delta$};
          \node at (0,-6) (delta5) {$\delta$};
          \node at (0,-8) (bc) {};
          \node at (-1.5,-6) (mu1) {$\mu_\wAlg^\bullet$};
          \node at (-3,-7) (mu2) {$\mu_\wAlg^\bullet$};
          \node at (-3,-8) (bl) {};
          \node at ($ (delta1)!.3!(delta2) $) (d12avg) {};
          \node at ($ (d12avg)!.6!(mu2) $) (thead) {};
          \draw[dmoda] (tc) to (delta1);
          \draw[dmoda] (delta1) to (delta2);
          \draw[dmoda] (delta2) to (delta3);
          \draw[dmoda] (delta3) to (f);
          \draw[dmoda] (f) to (delta4);
          \draw[dmoda] (delta4) to (delta5);
          \draw[dmoda] (delta5) to (bc);
          \draw[taa, bend right=10] (delta1) to (mu2);
          \draw[taa, bend left=5] (delta2) to (mu2);
          \draw[taa] (delta3) to (mu1);
          \draw[taa] (delta4) to (mu1);
          \draw[taa] (delta5) to (mu2);
          \draw[alga] (f) to (mu1);
          \draw[alga] (mu1) to (mu2);
          \draw[alga] (mu2) to (bl);
        \draw[->,thick,>=stealth] (mu2) to (thead);
      \end{tikzpicture}}\\
      (\mu_D(\bdy\corolla{\uparrow1}))(f)
    &=\mathcenter{
      \tikzsetnextfilename{prop-huD-id-new-3}
      \begin{tikzpicture}
        \node at (0,0) (tc) {};
        \node at (0,-1) (f) {$f^1$};
        \node at (0,-2) (bc) {};
        \node at (-1,-2) (bl) {};
        \draw[dmoda] (tc) to (f);
        \draw[dmoda] (f) to (bc);
        \draw[alga] (f) to (bl);
      \end{tikzpicture}}
      +
      \tikzsetnextfilename{prop-huD-id-new-4}
      \mathcenter{
      \begin{tikzpicture}[smallpic]
        \node at (0,0) (tr) {};
        \node at (0,-1) (delta1) {$\delta$};
       \node at (-.6,-3.25) (one) {};
        \node at (0,-3) (delta2) {$\delta$};      
        \node at (0,-4.5) (f) {$f^1$};      
        \node at (0,-5.5) (delta3) {$\delta$};      
        \node at (0,-7) (br) {};
        \node at (-1.5,-6) (mu) {$\mu_\wAlg^\bullet$};
        \node at (-1.5,-7) (bl) {};
        \draw[dmoda] (tr) to (delta1);
        \draw[dmoda] (delta1) to (delta2);
        \draw[dmoda] (delta2) to (f);
        \draw[dmoda] (f) to (delta3);
        \draw[dmoda] (delta3) to (br);
        \draw[alga] (mu) to (bl);
        \draw[taa,bend right=10] (delta1) to (mu);
        \draw[taa] (delta2) to (mu);
        \draw[taa] (delta3) to (mu);
        \draw[alga] (f) to (mu);
        \draw[stumpa] (one) to (mu);
      \end{tikzpicture}}
      +
      \tikzsetnextfilename{prop-huD-id-new-5}
      \mathcenter{\begin{tikzpicture}[smallpic]
          \node at (0,0) (tc) {};
          \node at (0,-1) (delta1) {$\delta$};
          \node at (0,-2) (delta2) {$\delta$};
          \node at (0,-3) (delta3) {$\delta$};
          \node at (0,-4) (delta4) {$\delta$};
          \node at (0,-5) (f) {$f^1$};
          \node at (0,-6) (delta5) {$\delta$};
          \node at (0,-8) (bc) {};
          \node at (-1.5,-4) (mu1) {$\mu_\wAlg^\bullet$};
          \node at (-3,-7) (mu2) {$\mu_\wAlg^\bullet$};
          \node at (-3,-8) (bl) {};
          \node at ($ (delta2)!.3!(delta3) $) (d23avg) {};
          \node at ($ (d23avg)!.4!(mu1) $) (thead) {};
          \draw[dmoda] (tc) to (delta1);
          \draw[dmoda] (delta1) to (delta2);
          \draw[dmoda] (delta2) to (delta3);
          \draw[dmoda] (delta3) to (delta4);
          \draw[dmoda] (delta4) to (f);
          \draw[dmoda] (f) to (delta5);
          \draw[dmoda] (delta5) to (bc);
          \draw[taa, bend right=10] (delta1) to (mu2);
          \draw[taa, bend right=10] (delta2) to (mu1);
          \draw[taa, bend left=5] (delta3) to (mu1);
          \draw[taa] (delta4) to (mu2);
          \draw[taa] (delta5) to (mu2);
          \draw[alga] (f) to (mu2);
          \draw[alga] (mu1) to (mu2);
          \draw[alga] (mu2) to (bl);
        \draw[->,thick,>=stealth] (mu1) to (thead);
      \end{tikzpicture}}.
  \end{align*}
  These are most of the terms coming from the homotopy unital $\Ainf$-relation
  at the $\mu_\wAlg^\bullet$-vertex of
  \[
      \mathcenter{\begin{tikzpicture}[smallpic]
          \node at (0,-1) (tc) {};
          \node at (0,-2) (delta2) {$\delta$};
          \node at (0,-3) (delta3) {$\delta$};
          \node at (0,-4) (f) {$f^1$};
          \node at (0,-5) (delta4) {$\delta$};
          \node at (0,-7) (bc) {};
          \node at (-1.5,-6) (mu1) {$\mu_\wAlg^\bullet$};
          \node at (-1.5,-7) (bl) {};
          \node at ($ (delta2)!.3!(delta3) $) (d23avg) {};
          \node at ($ (d23avg)!.6!(mu1) $) (thead) {};
          \draw[dmoda] (tc) to (delta2);
          \draw[dmoda] (delta2) to (delta3);
          \draw[dmoda] (delta3) to (f);
          \draw[dmoda] (f) to (delta4);
          \draw[dmoda] (delta4) to (bc);
          \draw[taa, bend right=10] (delta2) to (mu1);
          \draw[taa, bend left=5] (delta3) to (mu1);
          \draw[taa] (delta4) to (mu1);
          \draw[alga] (f) to (mu1);
          \draw[alga] (mu1) to (bl);
        \draw[->,thick,>=stealth] (mu1) to (thead);
      \end{tikzpicture}}.
  \]
  (The term with just the map $f^1$ corresponds to the case that all of the
  instances $\delta$ are the map $\delta^0=\Id$.)  The remaining terms coming
  from this $\Ainf$-relation, like
  \[
      \mathcenter{\begin{tikzpicture}[smallpic]
          \node at (0,1) (tc) {};
          \node at (0,0) (delta0) {$\delta$};
          \node at (0,-1) (delta1) {$\delta$};
          \node at (0,-2) (delta2) {$\delta$};
          \node at (0,-3) (delta3) {$\delta$};
          \node at (0,-4) (f) {$f^1$};
          \node at (0,-5) (delta4) {$\delta$};
          \node at (0,-7) (bc) {};
          \node at (-2,-6) (mu1) {$\mu_\wAlg^\bullet$};
          \node at (-1,-2) (mu2) {$\mu_\wAlg^\bullet$};
          \node at (-2,-7) (bl) {};
          \node at ($ (delta2)!.3!(delta3) $) (d23avg) {};
          \node at ($ (d23avg)!.6!(mu1) $) (thead) {};
          \draw[dmoda] (tc) to (delta0);
          \draw[dmoda] (delta0) to (delta1);
          \draw[dmoda] (delta1) to (delta2);
          \draw[dmoda] (delta2) to (delta3);
          \draw[dmoda] (delta3) to (f);
          \draw[dmoda] (f) to (delta4);
          \draw[dmoda] (delta4) to (bc);
          \draw[taa, bend right=10] (delta2) to (mu1);
          \draw[taa, bend left=5] (delta3) to (mu1);
          \draw[taa, bend right=20] (delta0) to (mu1);
          \draw[taa] (delta4) to (mu1);
          \draw[taa] (delta1) to (mu2);
          \draw[alga, bend right=5] (mu2) to (mu1);
          \draw[alga] (f) to (mu1);
          \draw[alga] (mu1) to (bl);
        \draw[->,thick,>=stealth] (mu1) to (thead);
      \end{tikzpicture}}
  \]
  cancel by the type $D$ structure relation.

  The strategy for the remaining cases is similar. We will again
  show that on the level of operation trees,
  \begin{multline}\label{eq:huD-big-check}
    \partial(\mu_D(\wcorolla{n_1\uparrow\cdots\uparrow n_k}{w})) +
    \mu_D(\partial\wcorolla{n_1\uparrow\cdots\uparrow n_k}{w})\\
    = \sum_{m_i,M_i} \Bigr(\bigl(\bigl(\mu_{\wAlg}^\bullet(\partial \wcorolla{n_1+M_1\uparrow\dots\uparrow n_k+M_k}{w})+\partial\mu_{\wAlg}^\bullet(\wcorolla{n_1+M_1\uparrow\dots\uparrow n_k+M_k}{w})\bigr) \otimes \Id_{P_N}\bigr)\\
    \circ
    (\Id \otimes \delta^{m_N}_N) \circ (\Id \otimes \phi_N) \circ\dots\circ
    (\Id \otimes \delta^{m_1}_1) \circ (\Id \otimes \phi_1) \circ
    \delta^{m_0}_0\Bigr)\\
    + \bigl(\text{terms that vanish by the type $D$ structure relations}\bigr)
  \end{multline}
  where the sum over the $m_i$ and $M_i$ is constrained as in
  Equation~\eqref{eq:huD-mu}. We will describe the various terms and how they
  cancel in pairs; Figure~\ref{fig:prop-gen-case} shows an example where terms
  of all of the kinds described appear.

  \begin{figure}
    \centering
    \[
      \mu_D(\bdy(\corolla{1\uparrow1}))(f^1)=
      \mathcenter{
        \tikzsetnextfilename{prop-huD-id-fig-1}
        \begin{tikzpicture}[smallpic]
          \node at (0,0) (tc) {};
          \node at (0,-1) (delta1) {$\delta$};
          \node at (0,-2) (f1) {$f_1^1$};
          \node at (0,-3) (delta2) {$\delta$};
          \node at (0,-4) (delta3) {$\delta$};
          \node at (0,-5) (f2) {$f_2^1$};
          \node at (0,-6) (delta4) {$\delta$};
          \node at (0,-8) (bc) {};
          \node at (-2,-7) (mu) {$\mu_\wAlg^\bullet$};
          \node at (-2,-8) (bl) {};
          \node at ($ (delta2)!.3!(delta3) $) (d23avg) {};
          \node at ($ (d23avg)!.5!(mu) $) (thead) {};
          \draw[dmoda] (tc) to (delta1);
          \draw[dmoda] (delta1) to (f1);
          \draw[dmoda] (f1) to (delta2);
          \draw[dmoda] (delta2) to (delta3);
          \draw[dmoda] (delta3) to (f2);
          \draw[dmoda] (f2) to (delta4);
          \draw[dmoda] (delta4) to (bc);
          \draw[taa] (delta1) to (mu);
          \draw[taa, bend right=5] (delta2) to (mu);
          \draw[taa, bend left=5] (delta3) to (mu);
          \draw[alga, bend right=2] (f1) to (mu);
          \draw[alga, bend left=2] (f2) to (mu);
          \draw[taa] (delta4) to (mu);
          \draw[alga] (mu) to (bl);
          \draw[stumpa] (thead) to (mu);
        \end{tikzpicture}
      }
      +
      \mathcenter{
        \tikzsetnextfilename{prop-huD-id-fig-2}
        \begin{tikzpicture}[smallpic]
          \node at (0,0) (tc) {};
          \node at (0,-1) (delta1) {$\delta$};
          \node at (0,-2) (f1) {$f_1^1$};
          \node at (0,-3) (delta2) {$\delta$};
          \node at (0,-4) (delta2a) {$\delta$};
          \node at (0,-5) (delta2b) {$\delta$};
          \node at (0,-6) (delta3) {$\delta$};
          \node at (0,-7) (f2) {$f_2^1$};
          \node at (0,-8) (delta4) {$\delta$};
          \node at (0,-10) (bc) {};
          \node at (-1,-6.5) (mub) {$\mu_\wAlg^\bullet$};
          \node at (-2,-9) (mu) {$\mu_\wAlg^\bullet$};
          \node at (-2,-10) (bl) {};
          \node at ($ (delta2a)!.3!(delta2b) $) (davg) {};
          \node at ($ (davg)!.5!(mub) $) (thead) {};
          \draw[dmoda] (tc) to (delta1);
          \draw[dmoda] (delta1) to (f1);
          \draw[dmoda] (f1) to (delta2);
          \draw[dmoda] (delta2) to (delta2a);
          \draw[dmoda] (delta2a) to (delta2b);
          \draw[dmoda] (delta2b) to (delta3);
          \draw[dmoda] (delta3) to (f2);
          \draw[dmoda] (f2) to (delta4);
          \draw[dmoda] (delta4) to (bc);
          \draw[taa] (delta1) to (mu);
          \draw[taa, bend right=5] (delta2) to (mu);
          \draw[taa, bend left=5] (delta3) to (mu);
          \draw[alga, bend right=2] (f1) to (mu);
          \draw[alga, bend left=2] (f2) to (mu);
          \draw[taa, bend right=5] (delta2a) to (mub);
          \draw[taa, bend left=5] (delta2b) to (mub);
          \draw[taa] (delta4) to (mu);
          \draw[alga] (mu) to (bl);
          \draw[alga] (mub) to (mu);
          \draw[->,thick,>=stealth] (mub) to (thead);
        \end{tikzpicture}
      }
      +
      \mathcenter{
        \tikzsetnextfilename{prop-huD-id-fig-3}
        \begin{tikzpicture}[smallpic]
          \node at (0,0) (tc) {};
          \node at (0,-1) (delta1) {$\delta$};
          \node at (0,-2) (delta2) {$\delta$};
          \node at (0,-3) (f1) {$f_1^1$};
          \node at (0,-4) (delta2a) {$\delta$};
          \node at (0,-5) (delta2b) {$\delta$};
          \node at (0,-6) (delta3) {$\delta$};
          \node at (0,-7) (f2) {$f_2^1$};
          \node at (0,-8) (delta4) {$\delta$};
          \node at (0,-10) (bc) {};
          \node at (-1,-6.5) (mub) {$\mu_\wAlg^\bullet$};
          \node at (-2,-9) (mu) {$\mu_\wAlg^\bullet$};
          \node at (-2,-10) (bl) {};
          \node at ($ (delta2a)!.3!(delta2b) $) (davg) {};
          \node at ($ (davg)!.5!(mub) $) (thead) {};
          \draw[dmoda] (tc) to (delta1);
          \draw[dmoda] (delta1) to (delta2);
          \draw[dmoda] (delta2) to (f1);
          \draw[dmoda] (f1) to (delta2a);
          \draw[dmoda] (delta2a) to (delta2b);
          \draw[dmoda] (delta2b) to (delta3);
          \draw[dmoda] (delta3) to (f2);
          \draw[dmoda] (f2) to (delta4);
          \draw[dmoda] (delta4) to (bc);
          \draw[taa] (delta1) to (mu);
          \draw[taa, bend right=5] (delta2) to (mub);
          \draw[taa, bend left=5] (delta3) to (mu);
          \draw[alga, bend right=5] (f1) to (mub);
          \draw[alga, bend left=2] (f2) to (mu);
          \draw[taa, bend right=5] (delta2a) to (mub);
          \draw[taa, bend left=5] (delta2b) to (mub);
          \draw[taa] (delta4) to (mu);
          \draw[alga] (mu) to (bl);
          \draw[alga] (mub) to (mu);
          \draw[->,thick,>=stealth] (mub) to (thead);
        \end{tikzpicture}
      }
      +
      \mathcenter{
        \tikzsetnextfilename{prop-huD-id-fig-4}
        \begin{tikzpicture}[smallpic]
          \node at (0,0) (tc) {};
          \node at (0,-1) (delta1) {$\delta$};
          \node at (0,-2) (f1) {$f_1^1$};
          \node at (0,-3) (delta2) {$\delta$};
          \node at (0,-4) (delta2a) {$\delta$};
          \node at (0,-5) (delta2b) {$\delta$};
          \node at (0,-6) (f2) {$f_2^1$};
          \node at (0,-7) (delta3) {$\delta$};
          \node at (0,-8) (delta4) {$\delta$};
          \node at (0,-10) (bc) {};
          \node at (-1,-7.5) (mub) {$\mu_\wAlg^\bullet$};
          \node at (-2,-9) (mu) {$\mu_\wAlg^\bullet$};
          \node at (-2,-10) (bl) {};
          \node at ($ (delta2a)!.3!(delta2b) $) (davg) {};
          \node at ($ (davg)!.5!(mub) $) (thead) {};
          \draw[dmoda] (tc) to (delta1);
          \draw[dmoda] (delta1) to (f1);
          \draw[dmoda] (f1) to (delta2);
          \draw[dmoda] (delta2) to (delta2a);
          \draw[dmoda] (delta2a) to (delta2b);
          \draw[dmoda] (delta2b) to (f2);
          \draw[dmoda] (f2) to (delta3);
          \draw[dmoda] (delta3) to (delta4);
          \draw[dmoda] (delta4) to (bc);
          \draw[taa] (delta1) to (mu);
          \draw[taa] (delta2) to (mu);
          \draw[taa, bend left=5] (delta3) to (mub);
          \draw[alga] (f1) to (mu);
          \draw[alga, bend left=5] (f2) to (mub);
          \draw[taa, bend right=5] (delta2a) to (mub);
          \draw[taa, bend left=5] (delta2b) to (mub);
          \draw[taa] (delta4) to (mu);
          \draw[alga] (mu) to (bl);
          \draw[alga] (mub) to (mu);
          \draw[->,thick,>=stealth] (mub) to (thead);
        \end{tikzpicture}
      }
    \]
    \[
      \bdy(\mu_D(\corolla{1\uparrow1})(f^1))=
      \mathcenter{
        \tikzsetnextfilename{prop-huD-id-fig-5}
        \begin{tikzpicture}[smallpic]
          \node at (0,0) (tc) {};
          \node at (0,-1) (delta1) {$\delta$};
          \node at (0,-2) (delta2) {$\delta$};
          \node at (0,-3) (f1) {$f_1^1$};
          \node at (0,-4) (delta2a) {$\delta$};
          \node at (0,-5) (delta2b) {$\delta$};
          \node at (0,-6) (f2) {$f_2^1$};
          \node at (0,-7) (delta3) {$\delta$};
          \node at (0,-8) (delta4) {$\delta$};
          \node at (0,-10) (bc) {};
          \node at (-1,-7.5) (mub) {$\mu_\wAlg^\bullet$};
          \node at (-2,-9) (mu) {$\mu_\wAlg^\bullet$};
          \node at (-2,-10) (bl) {};
          \node at ($ (delta2a)!.3!(delta2b) $) (davg) {};
          \node at ($ (davg)!.5!(mub) $) (thead) {};
          \draw[dmoda] (tc) to (delta1);
          \draw[dmoda] (delta1) to (delta2);
          \draw[dmoda] (delta2) (f1);
          \draw[dmoda] (f1) to (delta2a);
          \draw[dmoda] (delta2a) to (delta2b);
          \draw[dmoda] (delta2b) to (f2);
          \draw[dmoda] (f2) to (delta3);
          \draw[dmoda] (delta3) to (delta4);
          \draw[dmoda] (delta4) to (bc);
          \draw[taa, bend right=5] (delta1) to (mu);
          \draw[taa, bend right=5] (delta2) to (mub);
          \draw[taa, bend left=5] (delta3) to (mub);
          \draw[alga, bend right=5] (f1) to (mub);
          \draw[alga, bend left=5] (f2) to (mub);
          \draw[taa, bend right=5] (delta2a) to (mub);
          \draw[taa, bend left=5] (delta2b) to (mub);
          \draw[taa, bend left=5] (delta4) to (mu);
          \draw[alga] (mu) to (bl);
          \draw[alga] (mub) to (mu);
          \draw[->,thick,>=stealth] (mub) to (thead);
        \end{tikzpicture}
      }
      \qquad
      \mu_D(\corolla{1\uparrow1})(\bdy f^1)=
      \mathcenter{
        \tikzsetnextfilename{prop-huD-id-fig-6}
        \begin{tikzpicture}[smallpic]
          \node at (0,0) (tc) {};
          \node at (0,-1) (delta1) {$\delta$};
          \node at (0,-2) (delta1a) {$\delta$};
          \node at (0,-3) (f1) {$f_1^1$};
          \node at (0,-4) (delta1b) {$\delta$};
          \node at (0,-5) (delta2) {$\delta$};
          \node at (0,-6) (delta3) {$\delta$};
          \node at (0,-7) (f2) {$f_2^1$};
          \node at (0,-8) (delta4) {$\delta$};
          \node at (0,-10) (bc) {};
          \node at (-2,-9) (mu) {$\mu_\wAlg^\bullet$};
          \node at (-1,-5) (mub) {$\mu_\wAlg^\bullet$};
          \node at (-2,-10) (bl) {};
          \node at ($ (delta2)!.3!(delta3) $) (d23avg) {};
          \node at ($ (d23avg)!.5!(mu) $) (thead) {};
          \draw[dmoda] (tc) to (delta1);
          \draw[dmoda] (delta1) to (delta1a);
          \draw[dmoda] (delta1a) to (f1);
          \draw[dmoda] (f1) to (delta1b);
          \draw[dmoda] (delta1b) to (delta2);
          \draw[dmoda] (delta2) to (delta3);
          \draw[dmoda] (delta3) to (f2);
          \draw[dmoda] (f2) to (delta4);
          \draw[dmoda] (delta4) to (bc);
          \draw[taa, bend right=10] (delta1) to (mu);
          \draw[taa, bend right=7] (delta2) to (mu);
          \draw[taa, bend left=5] (delta3) to (mu);
          \draw[taa] (delta1a) to (mub);
          \draw[taa] (delta1b) to (mub);
          \draw[alga] (f1) to (mub);
          \draw[alga, bend left=5] (f2) to (mu);
          \draw[taa, bend left=5] (delta4) to (mu);
          \draw[alga] (mu) to (bl);
          \draw[alga, bend right=3] (mub) to (mu);
          \draw[->,thick,>=stealth] (mu) to (thead);
        \end{tikzpicture}
      }
      +
      \mathcenter{
        \tikzsetnextfilename{prop-huD-id-fig-7}
        \begin{tikzpicture}[smallpic]
          \node at (0,0) (tc) {};
          \node at (0,-1) (delta1) {$\delta$};
          \node at (0,-2) (f1) {$f_1^1$};
          \node at (0,-3) (delta2) {$\delta$};
          \node at (0,-4) (delta3) {$\delta$};
          \node at (0,-5) (delta3a) {$\delta$};
          \node at (0,-6) (f2) {$f_2^1$};
          \node at (0,-7) (delta3b) {$\delta$};
          \node at (0,-8) (delta4) {$\delta$};
          \node at (0,-10) (bc) {};
          \node at (-2,-9) (mu) {$\mu_\wAlg^\bullet$};
          \node at (-1,-8) (mub) {$\mu_\wAlg^\bullet$};
          \node at (-2,-10) (bl) {};
          \node at ($ (delta2)!.3!(delta3) $) (d23avg) {};
          \node at ($ (d23avg)!.5!(mu) $) (thead) {};
          \draw[dmoda] (tc) to (delta1);
          \draw[dmoda] (delta1) to (f1);
          \draw[dmoda] (f1) to (delta2);
          \draw[dmoda] (delta2) to (delta3);
          \draw[dmoda] (delta3) to (delta3a);
          \draw[dmoda] (delta3a) to (f2);
          \draw[dmoda] (f2) to (delta3b);
          \draw[dmoda] (delta3b) to (delta4);
          \draw[dmoda] (delta4) to (bc);
          \draw[taa, bend right=10] (delta1) to (mu);
          \draw[taa, bend right=8] (delta2) to (mu);
          \draw[taa, bend left=5] (delta3) to (mu);
          \draw[taa] (delta3a) to (mub);
          \draw[taa] (delta3b) to (mub);
          \draw[alga, bend right=8] (f1) to (mu);
          \draw[alga] (f2) to (mub);
          \draw[taa, bend left=5] (delta4) to (mu);
          \draw[alga] (mu) to (bl);
          \draw[alga] (mub) to (mu);
          \draw[->,thick,>=stealth] (mu) to (thead);
        \end{tikzpicture}
      }
    \]
    \caption[Generic case of the proof of Proposition~\ref{prop:hu-D-is-hu-Ainf-cat}]{\textbf{The generic case of the proof of Proposition~\ref{prop:hu-D-is-hu-Ainf-cat}.} The case of $\mu_D(\corolla{1\uparrow1})$ is shown. The terms in the first row are of types~\ref{item:huD-lhs-stump},\ref{item:huD-rhs-stump},~\ref{item:huD-lhs-stump},\ref{item:huD-rhs-split-2},~\ref{item:huD-lhs-split},\ref{item:huD-rhs-split-other}, and~\ref{item:huD-lhs-split},\ref{item:huD-rhs-split-other}. The term in the second row is of types~\ref{item:huD-lhs-split-after},\ref{item:huD-rhs-split-4}. The terms in the third row are of types~\ref{item:huD-lhs-split-after},\ref{item:huD-rhs-split-3}. The remaining terms obtained from the weighted $\Ainf$-relations cancel by the type $D$ structure relations.}
    \label{fig:prop-gen-case}
  \end{figure}
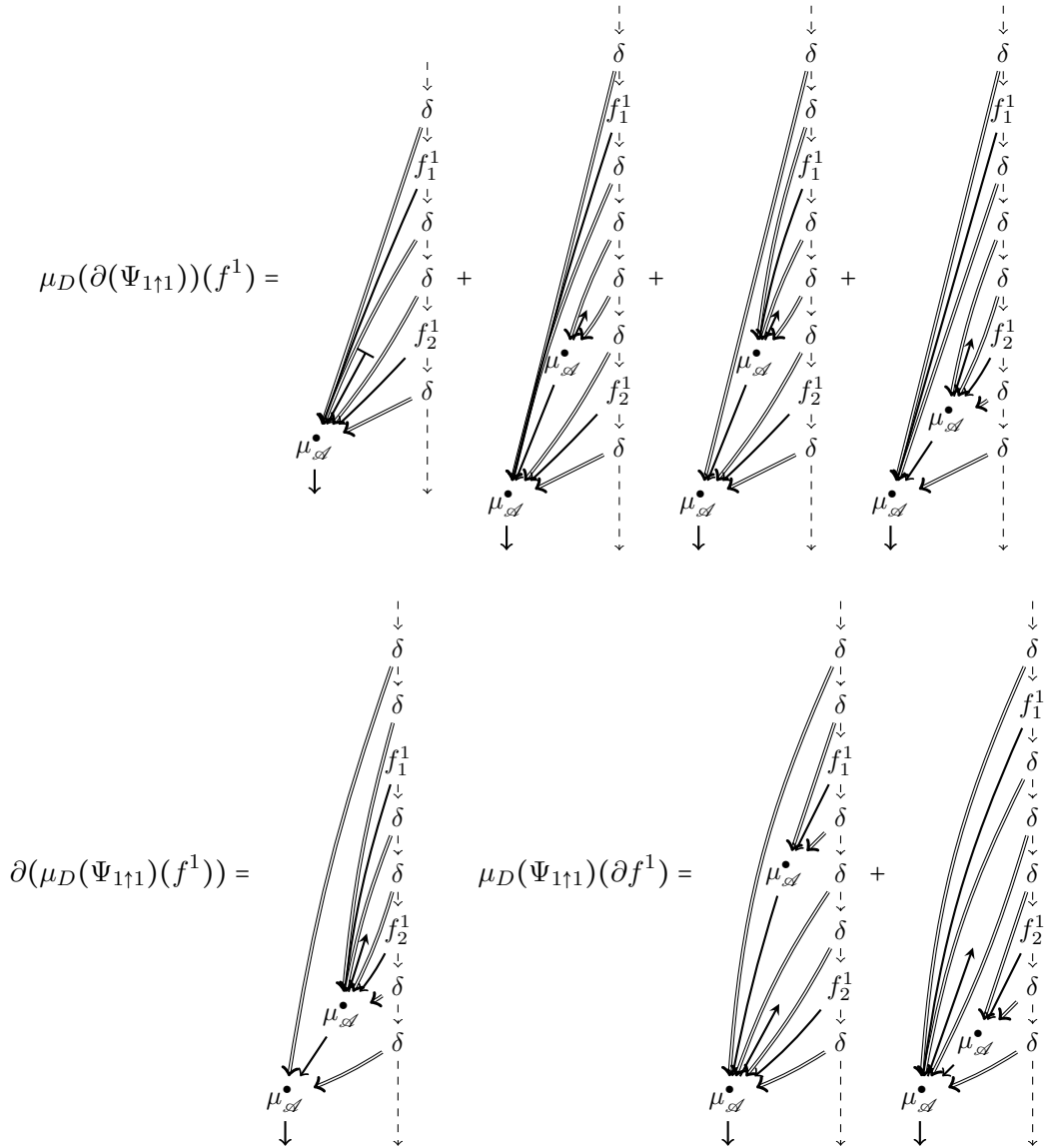

  In the $\mu_D \partial$ part of the left side of
  Equation~\eqref{eq:huD-big-check}, from the definition of the
  differential on
  $\uTreesCx{n}$, we have terms that
  \begin{enumerate}[label=(\arabic*)]
  \item\label{item:huD-lhs-stump} turn a thorn into a stump; or
  \item\label{item:huD-lhs-split} split the corolla into two corollas.
  \end{enumerate}
  We also have
  \begin{enumerate}[resume,label=(\arabic*)]
  \item\label{item:huD-lhs-split-after} terms in $\partial\mu_D$,
    i.e., coming from the differential on the space of type $D$
    morphisms.
  \end{enumerate}

  On the right side of Equation~\eqref{eq:huD-big-check}, the terms in $\partial\mu_\wAlg$ can be viewed as
  splittings of a vertex where one of the vertices involved is
  2-valent; we will usually not distinguish these terms from other
  types of splittings. Instead, we divide up the terms a different way.
  When we split a
  vertex, we label two resulting vertices as ``inner'' and ``outer'',
  with the output of the inner vertex feeding into the outer vertex.
  We consider what kinds of inputs feed into
  these two vertices: in addition to the strand from the inner vertex
  to the outer vertex, we have \emph{$\delta$ strands}, from one of the
  $\delta^{m_i}$ terms, \emph{morphism} strands from one of
  the~$\phi_i$, or thorns. Then on the right side we have terms that
  \begin{enumerate}[label=(\Alph*)]
  \item\label{item:huD-rhs-stump} turn a thorn into a stump;
  \end{enumerate}
  or splittings with
  \begin{enumerate}[resume,label=(\Alph*)]
  \item\label{item:huD-rhs-split-1} only $\delta$ strands into the inner vertex;
  \item\label{item:huD-rhs-split-2} $\delta$ strands and one thorn on the inner vertex;
  \item\label{item:huD-rhs-split-3} $\delta$ strands and one morphism strand on the inner vertex;
  \item\label{item:huD-rhs-split-4} only $\delta$ strands on the outer vertex (in addition to the
    strand from the inner vertex); or
  \item\label{item:huD-rhs-split-other} anything else.
  \end{enumerate}
  (Terms in $\partial\mu_\wAlg$ can only appear in
  types~\ref{item:huD-rhs-split-1} or~\ref{item:huD-rhs-split-4}.)

  These terms cancel as follows.
  \begin{itemize}
  \item Terms of type~\ref{item:huD-lhs-stump} cancel
    with terms of type~\ref{item:huD-rhs-stump} and~\ref{item:huD-rhs-split-2}.
  \item Terms of type~\ref{item:huD-rhs-split-1} form relations in the
    type~$D$ structures.
  \item Terms of type~\ref{item:huD-lhs-split-after} cancel with terms
    of type~\ref{item:huD-rhs-split-3}.
    and~\ref{item:huD-rhs-split-4}
  \item Terms of type~\ref{item:huD-lhs-split} cancel with terms of
    type~\ref{item:huD-rhs-split-other}; since a term of
    type~\ref{item:huD-rhs-split-other} does not fall into one of the
    other cases, it must correspond to a valid splitting of the vertex
    in the differential on $\uTreesCx{n}$.\qedhere
  \end{itemize}
\end{proof}

\begin{definition}
  Given a homotopy unital $\Ainf$-category $\Cat$, the \emph{homology
    category $H_*\Cat$} of $\Cat$ is the category with the same
  objects as $\Cat$ and $\Hom_{H_*\Cat}(X,Y)=H_*\Mor^{nu}_\Cat(X,Y)$,
  where $\Mor^{nu}_\Cat(X,Y)$ denotes the space of morphisms in the
  non-unital $\Ainf$-category associated to $\Cat$ (by ignoring
  morphisms with a positive number of thorns or stumps). The identity
  map in $\Hom_{H_*\Cat}(X,Y)$ is the image of the stump.
\end{definition}

\begin{lemma}\label{lem:homol-of-hu-cat}
  Given a homotopy unital $\Ainf$-category $\Cat$, the homology
  category $H_*\Cat$ is an ordinary category.
\end{lemma}
In other words, Lemma~\ref{lem:homol-of-hu-cat} asserts that a
homotopy unital $\Ainf$-category is homologically unital.
\begin{proof}
  We need to check that composition is associative (i.e.,
  $(f\circ_2 g)\circ_2 h = f\circ_2(g\circ_2 h)$) and that the image
  of the stump is an identity (i.e.,
  $f\circ_2\Id_X=\Id_Y\circ_2 f=f$). Both are immediate from the
  definitions.
\end{proof}

\begin{definition}
  Let $\Cat$ and $\Dat$ be a homotopy unital $\Ainf$-categories. A
  \emph{homotopy unital functor} $F\co\Cat\to\Dat$ consists
  of:
  \begin{itemize}
  \item For each object $X\in\ob(\Cat)$, an object $F(X)\in\ob(\Dat)$,
    and
  \item For each sequence of objects $X_0,\dots,X_n\in\ob(\Cat)$, a
    chain map
    \[
      F\co \uTransCx{n}\to
      \Mor(\Mor_\Cat(X_0,X_1)\otimes\cdots\otimes\Mor_\Cat(X_{n-1},X_n),\Mor_\Dat(F(X_0),F(X_n)))
    \]
  \end{itemize}
  so that for $S\in\uTreesCx{m}$ and $T\in \uTransCx{n}$
  \[
    F(T\circ_i S)=F(T)\circ_i\mu_\Cat(S)
  \]
  and for $S_i\in\uTransCx{m_i}$, $i=1,\dots,n$ and
  $T\in\uTreesCx{n}$,
  \[
    F(T\circ(S_1,\dots,S_n)) = \mu_\Dat(T)\circ(F(S_1),\dots,F(S_n)).
  \]

  Given a homotopy unital functor $F\co\Cat\to\Dat$, $F$ restricts to
  a functor of underlying non-unital $\Ainf$-categories.
\end{definition}

\begin{lemma}
  Given a homotopy unital functor $F\co\Cat\to\Dat$, there is an
  induced honest functor $F\co H_*\Cat\to H_*\Dat$ of ordinary
  categories.
\end{lemma}
\begin{proof}
  Again, this is immediate from the definitions.
\end{proof}

\begin{definition}
  A homotopy unital functor $F\co\Cat\to\Dat$ is a
  \emph{quasi-equivalence} if the induced functor of ordinary
  categories is an equivalence.
\end{definition}

Given strictly unital $\Ainf$-categories $\Cat$, $\Dat$, an $\Ainf$-functor
$F\co\Cat\to\Dat$ is \emph{strictly unital} if for each $X\in\ob(\Cat)$,
$F(\Id_X)=\Id_{F(X)}$ and for each $X_1,\dots,X_n\in\ob(\Cat)$, $n\geq
2$, and $f_i\in\Mor_\Cat(X_i,X_{i+1})$, if some $f_i=\Id_{X_i}$ then
$F(f_{n-1},\dots,f_1)=0$.
We note the following analogue of Proposition~\ref{prop:qi-unital}:
\begin{proposition}\label{prop:strictify-hu-quasi-equiv}
  If $\Cat$ and $\Dat$ are strictly unital $\Ainf$-categories which
  are homotopy unitally quasi-isomorphic then $\Cat$ and $\Dat$ are
  quasi-isomorphic via a strictly unital functor.
\end{proposition}
\begin{proof}
  The proof is similar to the proof of
  Proposition~\ref{prop:qi-unital} or Theorem~\ref{thm:hu-rectify}
  or~\cite[Lemma 2.1]{SeidelBook} (see also~\cite[Remark
  2.2]{SeidelBook}) or~\cite[Theorem 3.2.2.1]{LefevreAInfinity}, and
  is left to the reader.
\end{proof}

\begin{definition}\label{def:hu-induct-func}
  Fix homotopy unital weighted $\Ainf$-algebras $\wAlg$ and $\wBlg$
  over $\Ground$, an element $X\in\Ground$ acting centrally on $A$ and
  $B$,
  and a homotopy unital homomorphism $F\co\wAlg\to\wBlg$. Then there is a
  homotopy unital functor $F_*\co \lsup{\wAlg,X}\ModCat\to
  \lsup{\wBlg,X}\ModCat$ between categories of type $D$ structures with
  charge $X$ defined as follows. On objects, $F_*$ sends a type $D$
  structure $(P,\delta^1_P)$ to the same $\Ground$-module $P$ with
  operation
  \[
    \delta^1_{F_*(P)}=
    \sum_{n=0}^\infty
    \tikzsetnextfilename{hu-induct-func-1}
    \mathcenter{
      \begin{tikzpicture}[smallpic]
        \node at (0,0) (tc) {};
        \node at (0,-1) (delta) {$\delta^n_P$};
        \node at (0,-3) (bc) {};
        \node at (-1,-2) (muL) {$F^\bullet_n$};
        \node at (-1,-3) (bl) {};
        \draw[dmoda] (tc) to (delta);
        \draw[dmoda] (delta) to (bc);
        \draw[taa] (delta) to (muL);
        \draw[alga] (muL) to (bl);
      \end{tikzpicture}
    }
    =
    \sum_{(n,w)\neq(0,0)}X^w
    \tikzsetnextfilename{hu-induct-func-2}
    \mathcenter{
      \begin{tikzpicture}[smallpic]
        \node at (0,0) (tc) {};
        \node at (0,-1) (delta) {$\delta^n_P$};
        \node at (0,-3) (bc) {};
        \node at (-1,-2) (muL) {$F^w_n$};
        \node at (-1,-3) (bl) {};
        \draw[dmoda] (tc) to (delta);
        \draw[dmoda] (delta) to (bc);
        \draw[taa] (delta) to (muL);
        \draw[alga] (muL) to (bl);
      \end{tikzpicture}.
    }
  \]
  On sequences of morphisms, $F_*$ is defined by:
  \begin{multline}\label{eq:huD-F}
    F_*(\pcorolla{n_1\uparrow\cdots\uparrow n_k})(\phi_1,\dots,\phi_N)
    = \sum_w X^w \sum_{m_0,\dots,m_N=0}^\infty \,\,\sideset{}{'}\sum_{M_1,\dots,M_k=0}^\infty
    \Bigl(\bigl(F^w_{n_1+M_1\uparrow n_2+M_2\uparrow\cdots\uparrow n_k+M_k}\otimes\Id_{P_N}\bigr) \\
    \circ
    (\Id \otimes \delta^{m_N}_N) \circ (\Id \otimes \phi_N) \circ\dots\circ
    (\Id \otimes \delta^{m_1}_1) \circ (\Id \otimes \phi_1) \circ
    \delta^{m_0}_0\Bigr)
  \end{multline}
  where as in Equation~\eqref{eq:huD-mu} we constrain the sum over
  $M_i$ to make the arrows properly interspersed. That is, the map
  $F_*$ is obtained by
  feeding the outputs of the type $D$ structure morphisms into the corolla,
  interspersed with an arbitrary number of $\delta_j^1$ outputs, in all ways
  compatible with the thorns. (For a graphical representation of a few cases,
  see Figure~\ref{fig:hu-D-comps}, and replace the vertex $\mu_{\wAlg}^\bullet$
  with $F^\bullet$.)
\end{definition}

\begin{proposition}\label{prop:hu-D-quasi-equiv}
  Definition~\ref{def:hu-induct-func} defines a homotopy unital
  functor. Further, if $F$ is an isomorphism then $F_*$ is a
  quasi-equivalence.
\end{proposition}
\begin{proof}
  The proof of the first statement is similar to the proof of
  Proposition~\ref{prop:hu-D-is-hu-Ainf-cat}, and is left to the
  reader.

  For the second statement, it suffices to verify that for each pair
  of type $D$ structures $\lsup{\wAlg}P$, $\lsup{\wAlg}Q$, the map
  \[
    F_*(\pcorolla{1})\co \Mor(\lsup{\wAlg}P,\lsup{\wAlg}Q)\to\Mor(\lsup{\wBlg}F_*(P),\lsup{\wBlg}F_*(Q))
  \]
  induces an isomorphism on homology; in fact, we will show that if
  $G$ is the inverse of $F$ then
  $F_*(\pcorolla{1})$ is an isomorphism
  of chain complexes, with inverse given by $G_*(\pcorolla{1})$.
  (Essential surjectivity of
  the map of homology categories induced by $F_*$ is clear from
  this and invertibility of $F$.)

  By definition,
  \[
    F_*(\pcorolla{1})(\phi^1)=
    \tikzsetnextfilename{induct-D-1}
    \mathcenter{
      \begin{tikzpicture}[smallpic]
        \node at (0,0) (tc) {};
        \node at (0,-1) (delta1) {$\delta$};
        \node at (0,-2) (phi) {$\phi^1$};
        \node at (0,-3) (delta2) {$\delta$};
        \node at (-1,-4) (F) {$F^\bullet$};
        \node at (-1,-5) (bl) {};
        \node at (0,-5) (bc) {};
        \draw[dmoda] (tc) to (delta1);
        \draw[dmoda] (delta1) to (phi);
        \draw[dmoda] (phi) to (delta2);
        \draw[dmoda] (delta2) to (bc);
        \draw[taa] (delta1) to (F);
        \draw[taa] (delta2) to (F);
        \draw[alga] (phi) to (F);
        \draw[alga] (F) to (bl);
      \end{tikzpicture}
    }.
  \]
  Let
  \[
    \underline{F}^\bullet=\sum_{m=1}^\infty
    (F^\bullet\otimes\cdots\otimes F^\bullet)\co \bigoplus_n A^{\otimes n}\to
    \bigoplus_m B^{\otimes m},
  \]
  i.e.,
  \[
    \tikzsetnextfilename{induce-F-on-tens-1}
    \mathcenter{
      \begin{tikzpicture}
        \node at (0,0) (tc) {};
        \node at (0,-1) (F) {$\underline{F}^\bullet$};
        \node at (0,-2) (bc) {};
        \draw[taa] (tc) to (F);
        \draw[tbb] (F) to (bc);
      \end{tikzpicture}
    }
    =
    \sum
    \tikzsetnextfilename{induce-F-on-tens-2}
    \mathcenter{
      \begin{tikzpicture}
        \node at (0,0) (t1) {};
        \node at (2,0) (t2) {};
        \node at (.5,0) (t1b) {};
        \node at (1.5,0) (t2b) {};
        \node at (0,-1) (F1) {$F^\bullet$};
        \node at (.5,-1) (F1b) {$F^\bullet$};
        \node at (1,-1) (cdots) {$\cdots$};
        \node at (2,-1) (F2) {$F^\bullet$};
        \node at (1.5,-1) (F2b) {$F^\bullet$};
        \node at (1,-2) (bc) {};
        \draw[taa] (t1) to (F1);
        \draw[taa] (t2) to (F2);
        \draw[blga] (F1) to (bc);
        \draw[blga] (F2) to (bc);
        \draw[taa] (t1b) to (F1b);
        \draw[taa] (t2b) to (F2b);
        \draw[blga] (F1b) to (bc);
        \draw[blga] (F2b) to (bc);
      \end{tikzpicture}
    }.
  \]  
  The fact that $G$ is the inverse to $F$ implies that the composition
  \[
    \tikzsetnextfilename{comp-is-id}
    \begin{tikzpicture}
      \node at (0,0) (tc) {};
      \node at (0,-1) (F) {$\underline{F}^\bullet$};
      \node at (0,-2) (G) {$G^\bullet$};
      \node at (0,-3) (bc) {};
      \draw[taa] (tc) to (F);
      \draw[tbb] (F) to (G);
      \draw[alga] (G) to (bc);
    \end{tikzpicture}
  \]
  is the identity map $A\to A$ (and vanishes on $A^{\otimes m}$ for
  $m>1$).  So,
  \[
    G_*(\pcorolla{1})\bigl(F_*(\pcorolla{1})(\phi^1)\bigr)=
    \tikzsetnextfilename{induct-D-2}
    \mathcenter{
      \begin{tikzpicture}[smallpic]
        \node at (0,1) (tc) {};
        \node at (0,0) (delta0) {$\delta$};
        \node at (0,-1) (delta1) {$\delta$};
        \node at (0,-2) (phi) {$\phi^1$};
        \node at (0,-3) (delta2) {$\delta$};
        \node at (0,-4) (delta3) {$\delta$};
        \node at (-1,-1) (F1) {$\underline{F}^\bullet$};
        \node at (-1,-4) (F2) {$F^\bullet$};
        \node at (-1,-5) (F3) {$\underline{F}^\bullet$};
        \node at (0,-7) (bc) {};
        \node at (-2,-6) (G) {$G^\bullet$};
        \node at (-2,-7) (bl) {};
        \draw[dmoda] (tc) to (delta0);
        \draw[dmoda] (delta0) to (delta1);
        \draw[dmoda] (delta1) to (phi);
        \draw[dmoda] (phi) to (delta2);
        \draw[dmoda] (delta2) to (delta3);
        \draw[dmoda] (delta3) to (bc);
        \draw[taa] (delta0) to (F1);
        \draw[taa] (delta1) to (F2);
        \draw[taa] (delta2) to (F2);
        \draw[taa] (delta3) to (F3);
        \draw[tbb] (F1) to (G);
        \draw[tbb] (F3) to (G);
        \draw[alga] (phi) to (F2);
        \draw[blga] (F2) to (G);
        \draw[alga] (G) to (bl);
      \end{tikzpicture}
    }
    =
    \mathcenter{
      \tikzsetnextfilename{induct-D-3}
      \begin{tikzpicture}[smallpic]
        \node at (0,0) (tc) {};
        \node at (0,-1) (phi) {$\phi^1$};
        \node at (0,-2) (bc) {};
        \node at (-1,-2) (bl) {};
        \draw[dmoda] (tc) to (phi);
        \draw[dmoda] (phi) to (bc);
        \draw[alga] (phi) to (bl);
      \end{tikzpicture}
    }.
  \]
  So, $G_*(\pcorolla{1})\circ F_*(\pcorolla{1})$ is the
  identity map of $\Mor(\lsup{\wAlg}P,\lsup{\wAlg}Q)$. Similarly,
  $F_*(\pcorolla{1})\circ G_*(\pcorolla{1})$ is the identity map of
  $\Mor(\lsup{\wBlg}F_*(P),\lsup{\wBlg}F_*(Q))$. This proves the result.
\end{proof}

\begin{definition}
  Let $\wAlg$ and $\wBlg$ be homotopy unital $\Ainf$-algebras over
  $\Ground_1$ and~$\Ground_2$,
  $X_{12}\in\Ground$, and $\uwDiag{*}{*}$ a homotopy
  unital algebra diagonal. Then the \emph{homotopy unital category of
    type \DD\ structures over $\wAlg$ and $\wBlg$ with charge $X_{12}$}
  is the homotopy unital category of type $D$ structures over
  $\wAlg\ADtp[\uwDiagNS]\wBlg$ with charge $X_{12}$.
\end{definition}

\begin{corollary}\label{cor:hu-DD-indep-diag}
  Let $\wAlg$ and $\wBlg$ be homotopy unital $\Ainf$-algebras over
  $\Ground_1$ and~$\Ground_2$ and let
  $X_{12}\in\Ground$.  Let $\uwDiag{*}{*}_1$ and
  $\uwDiag{*}{*}_2$ be homotopy unital algebra diagonals with the same
  seed. Then the categories of type \DD\ structures over $\wAlg$ and
  $\wBlg$ with charge $X_{12}$ with respect to $\uwDiag{*}{*}_1$ and
  $\uwDiag{*}{*}_2$ are quasi-equivalent.
\end{corollary}
\begin{proof}
  By Theorem~\ref{thm:hu-tens-prod-props}, there is a homotopy unital
  isomorphism
  $\wAlg\wADtp[\uwDiagNS_1]\wBlg\to \wAlg\wADtp[\uwDiagNS_2]\wBlg$. So,
  this follows from Proposition~\ref{prop:hu-D-quasi-equiv}.
\end{proof}

We finally reach the goal of this discussion, a result not mentioning
homotopy unitality at all:
\begin{theorem}\label{thm:DD-indep-diag}
  Given split unital weighted $\Ainf$-algebras $\wAlg$ and $\wBlg$
  over $\Ground_1$ and~$\Ground_2$, respectively;
  $X_{12}\in\Ground$; weighted diagonals $\wDiag{*}{*}_1$ and
  $\wDiag{*}{*}_2$; and strictly unital weighted algebras $\wClg_i$
  isomorphic to $\wAlg\wADtp[\wDiagNS_i]\wBlg$, the (strictly
  unital) categories of type $D$ structures $\lsup{\wClg_1,X_{12}}\ModCat$
  and $\lsup{\wClg_2,X_{12}}\ModCat$ with charge $X_{12}$ are
  quasi-equivalent.
\end{theorem}
\begin{proof}
  Let $\uwDiag{*}{*}_i$ be a homotopy unital diagonal extending
  $\wDiag{*}{*}_i$, so $\wAlg\wADtp[\uwDiagNS_i]\wBlg$ is a homotopy
  unital algebra with underlying non-unital algebra
  $\wAlg\wADtp[\wDiagNS_i]\wBlg$.  By Theorem~\ref{thm:hu-rectify},
  there is a strictly unital algebra $\Clg'_i$ homotopy unitally
  isomorphic to $\wAlg\wADtp[\uwDiagNS_i]\wBlg$. By
  Theorem~\ref{thm:hu-tens-prod-props},
  $\wAlg\wADtp[\uwDiagNS_1]\wBlg$ and $\wAlg\wADtp[\uwDiagNS_2]\wBlg$
  are homotopy unitally isomorphic. So, by
  Proposition~\ref{prop:hu-D-quasi-equiv}, the categories of type $D$
  structures over $\Clg_1'$, $\wAlg\wADtp[\wDiagNS_1]\wBlg$,
  $\wAlg\wADtp[\wDiagNS_2]\wBlg$, and $\Clg_2'$ are all
  quasi-equivalent. On the other hand, by
  Proposition~\ref{prop:w-qi-unital}, $\Clg_i$ and $\Clg'_i$ are
  strictly unitally quasi-isomorphic, so the categories of type $D$
  structures over $\Clg_i$ and $\Clg'_i$ are quasi-equivalent.  Thus,
  the categories of type $D$ structures over $\Clg_1$ and $\Clg_2$ are
  homotopy unitally quasi-equivalent. By
  Proposition~\ref{prop:strictify-hu-quasi-equiv}, this, in turn,
  implies these categories of type $D$ structures are strictly
  unitally quasi-equivalent.
\end{proof}

%%% Local Variables: 
%%% mode: latex
%%% TeX-master: "AbstractDiagonal.tex"
%%% TeX-command-extra-options: "-shell-escape"
%%% End: 

\newcommand\InputSet{I}
\section{The associaplex}\label{sec:Associaplex}

\newcommand\cHHH{\overline{\HHH}}
\newcommand\oHHH{\HHH^o}
\newcommand{\Associaplex}[2]{{\mathbf X}^{#1,#2}}
\newcommand\Prod[1]{\widetilde{Q}^{#1}}
\newcommand\oCDisk{{\mathbb D}^{\circ}}
\newcommand\goesto{\rightarrow}
\newcommand\cCDisk{\mathbb D}
\newcommand{\SymQ}{Q}
\newcommand{\lQ}[1]{\widetilde{Q}^{#1}}
\newcommand\lAssoc[2]{\widetilde{\mathbf{X}}^{#1,#2}}

\begin{figure}
  \centering
  \tikzsetnextfilename{assoc-conv}
  \[
    \begin{array}{ccc}
      \mathcenter{\begin{tikzpicture}
    \coordinate (z1) at (0,-3);
    \node at ($ (z1) + (0,-.35) $) (z1lab) {$z_1$};
    \coordinate (z2) at (-1,-1);
    \node at ($ (z2) + (0,-.35) $) (z2lab) {$z_6$};
    \coordinate (z2b) at ($ (z2) +.5*(0,1) + .2*(-2,0) $);    
    \coordinate (z2c) at ($ (z2) +.75*(0,1) + .4*(-2,0) $);    
    \coordinate (z2d) at ($ (z2) +.875*(0,1) + .6*(-2,0) $);    
    \coordinate (z2e) at ($ (z2) +.9375*(0,1) + .8*(-2,0) $);    
    \coordinate (z3) at (-1,1);
    \node at ($ (z3) + (0,.35) $) (z3lab) {$z_5$};
    \coordinate (z3b) at ($ (z3) -.5*(0,1) + .2*(-2,0) $);    
    \coordinate (z3c) at ($ (z3) -.75*(0,1) + .4*(-2,0) $);    
    \coordinate (z3d) at ($ (z3) -.875*(0,1) + .6*(-2,0) $);    
    \coordinate (z3e) at ($ (z3) -.9375*(0,1) + .8*(-2,0) $);
    \coordinate (z4) at ({3*cos(45)},{3*sin(45)});
    \node at ($ (z4) + .35*({cos(45)},{cos(45)}) $) (z4lab) {$z_3$};
    \coordinate (z5) at ({3*cos(90)},{3*sin(90)});
    \node at ($ (z5) + .35*(0,1) $) (z5lab) {$z_4$};
    \coordinate (z5b) at ({3*cos((90+45)/2)},{3*sin((90+45)/2)});
    \coordinate (z5c) at ({3*cos((90+2*45)/3)},{3*sin((90+2*45)/3)});
    \coordinate (z5d) at ({3*cos((90+4*45)/5))},{3*sin((90+4*45)/5))});
    \coordinate (z5e) at ({3*cos((90+12*45)/13))},{3*sin((90+12*45)/13))});
    \coordinate (z6) at ({3*cos(-45)},{3*sin(-45)});
    \node at ($ (z6) + .35*({cos(-45)},{sin(-45)}) $) (z6lab) {$z_2$};
    \coordinate (z6b) at ({3*cos((3*-45+2*45)/5)},{3*sin((3*-45+2*45)/5)});
    \coordinate (z6c) at ({3*cos((-45+2*45)/3)},{3*sin((-45+2*45)/3)});
    \coordinate (z6d) at ({3*cos((-45+4*45)/5)},{3*sin((-45+4*45)/5)});
    \coordinate (z6e) at ({3*cos((-45+7*45)/8)},{3*sin((-45+7*45)/8)});
    \coordinate (z7) at (1,1.25);
    \node at ($ (z7) + .35*(0,1) $) (z7lab) {$z_7$};
    \coordinate (z7b) at ($ (z7)!.3!(z4) $);
    \coordinate (z7c) at ($ (z7)!{1-.7*.7}!(z4) $);
   \coordinate (z7d) at ($ (z7)!{1-.7*.7*.7}!(z4) $);
   \coordinate (z7e) at ($ (z7)!{1-.7*.7*.7*.7}!(z4) $);
   \coordinate (z8) at (1.5,-.5);
    \node at ($ (z8) - .35*(0,1) $) (z8lab) {$z_8$};
    \coordinate (z8b) at ($ (z8)!.3!(z4) $);
    \coordinate (z8c) at ($ (z8)!{1-.7*.7}!(z4) $);
   \coordinate (z8d) at ($ (z8)!{1-.7*.7*.7}!(z4) $);
   \coordinate (z8e) at ($ (z8)!{1-.7*.7*.7*.7}!(z4) $);
    \draw (0,0) circle (3cm); %Outer circle
    \filldraw (z1) circle (.1cm);  %z1
    \filldraw[fill=blue!100!white] (z2) circle (.1cm);
    \filldraw[fill=blue!80!white] (z2b) circle (.1cm);
    \filldraw[fill=blue!60!white] (z2c) circle (.1cm);
    \filldraw[fill=blue!40!white] (z2d) circle (.1cm);
    \filldraw[fill=blue!20!white] (z2e) circle (.1cm);
    \filldraw[fill=green!100!white] (z3) circle (.1cm);
    \filldraw[fill=green!80!white] (z3b) circle (.1cm);
    \filldraw[fill=green!60!white] (z3c) circle (.1cm);
    \filldraw[fill=green!40!white] (z3d) circle (.1cm);
    \filldraw[fill=green!20!white] (z3e) circle (.1cm);
    \filldraw[fill=white] (z4) circle (.1cm);
    \filldraw[fill=red!100!white] (z5) circle (.1cm);
    \filldraw[fill=red!80!white] (z5b) circle (.1cm);
    \filldraw[fill=red!60!white] (z5c) circle (.1cm);
    \filldraw[fill=red!40!white] (z5d) circle (.1cm);
    \filldraw[fill=red!20!white] (z5e) circle (.1cm);
    \filldraw[fill=cyan!100!white] (z6) circle (.1cm);
    \filldraw[fill=cyan!80!white] (z6b) circle (.1cm);
    \filldraw[fill=cyan!60!white] (z6c) circle (.1cm);
    \filldraw[fill=cyan!40!white] (z6d) circle (.1cm);
    \filldraw[fill=cyan!20!white] (z6e) circle (.1cm);
    \filldraw[fill=yellow!100!white] (z7) circle (.1cm);
    \filldraw[fill=yellow!80!white] (z7b) circle (.1cm);
   \filldraw[fill=yellow!60!white] (z7c) circle (.1cm);
   \filldraw[fill=yellow!40!white] (z7d) circle (.1cm);
   \filldraw[fill=yellow!20!white] (z7e) circle (.1cm);
    \filldraw[fill=brown!100!white] (z8) circle (.1cm);
    \filldraw[fill=brown!80!white] (z8b) circle (.1cm);
    \filldraw[fill=brown!60!white] (z8c) circle (.1cm);
    \filldraw[fill=brown!40!white] (z8d) circle (.1cm);
    \filldraw[fill=brown!20!white] (z8e) circle (.1cm);
  \end{tikzpicture}} & \longrightarrow & 
  \tikzsetnextfilename{assoc-conv-2}
  \mathcenter{\begin{tikzpicture}
    \draw (0,0) circle (1cm); %Root circle    
    \draw (-2,0) circle (1cm); %Left circle    
    \draw ({sqrt(2)},{sqrt(2)}) circle (1cm); %up-right circle    
    \draw ({2*sqrt(2)},{2*sqrt(2)}) circle (1cm);
    \coordinate (z1) at (0,-1);
    \node at ($ (z1) + (0,-.35) $) (z1lab) {$z_1$};
    \filldraw (z1) circle (.1cm);  %z1
    \coordinate (z23) at (-2,0);
    \node at ($ (z23) + (0,-.35) $) (z23lab) {$\{z_5,z_6\}$};
    \filldraw[fill=blue!50!green] (z23) circle (.1cm);
    \coordinate (z6) at ({sqrt(2)},{sqrt(2)-1});
    \filldraw[fill=cyan!100!white] (z6) circle (.1cm);
    \node at ($ (z6) + (0,-.35) $) (z6lab) {$z_2$};
    \coordinate (z7) at ({sqrt(2)-.25},{sqrt(2)+.15});
    \filldraw[fill=yellow!100!white] (z7) circle (.1cm);
    \node at ($ (z7) + (0,-.35) $) (z7lab) {$z_7$};
    \coordinate (z8) at ({sqrt(2)+.35},{sqrt(2)+.05});
    \filldraw[fill=brown!100!white] (z8) circle (.1cm);
    \node at ($ (z8) + (0,-.35) $) (z8lab) {$z_8$};
    \coordinate(z5) at ({2*sqrt(2)},{2*sqrt(2)+1});
    \filldraw[fill=red!100!white] (z5) circle (.1cm);
    \node at ($ (z5) + (0,-.35) $) (z5lab) {$z_4$};
    \coordinate(z4) at ({2*sqrt(2)+1},{2*sqrt(2)});
    \filldraw[fill=white] (z4) circle (.1cm);
    \node at ($ (z4) + (-.35,0) $) (z4lab) {$z_3$};
  \end{tikzpicture}}\\
      \tikzsetnextfilename{assoc-conv-tree}
      \mathcenter{\begin{tikzpicture}
        \draw[-] (0,0) to (0,-1);
        \draw[-] (-1,1) to (0,0);
        \draw[-] (0,1) to (0,0);
        \draw[-] (1,1) to (0,0);
        \filldraw[fill=white] (0,0) circle (.25cm);
        \node at (0,0) (weight) {$4$};
      \end{tikzpicture}}
    & &\tikzsetnextfilename{assoc-conv-tree-2}
      \mathcenter{\begin{tikzpicture}
        \draw[-] (0,0) to (0,-1);
        \draw[-] (-1,1) to (0,0);
        \draw[-] (1,1) to (0,0);
        \draw[-] (0,2) to (1,1);
        \draw[-] (2,2) to (1,1);
        \draw[-] (1,3) to (0,2);
        \draw[-] (-1,3) to (0,2);
        \filldraw[fill=white] (0,0) circle (.25cm);
        \node at (0,0) (weight1) {$0$};
        \filldraw[fill=white] (-1,1) circle (.25cm);
        \node at (-1,1) (weight2) {$2$};
        \filldraw[fill=white] (1,1) circle (.25cm);
        \node at (1,1) (weight3) {$2$};
        \filldraw[fill=white] (0,2) circle (.25cm);
        \node at (0,2) (weight4) {$0$};
      \end{tikzpicture}}
    \end{array}
  \]
  \caption[Convergence in the associaplex]{\textbf{Convergence in the
      associaplex.} Top-left: a sequence of points in
    $\Sym^{4,4}(\cCDisk)$. Subsequent terms in the sequence appear in
    lighter colors. Top-right: the limit in (one stratum of) the
    associaplex $\SymQ^T(\cCDisk)$. Note that the points $z_5,z_6$
    would converge to distinct points in a sphere bubble in the Deligne-Mumford
    compactification, and this sphere has been collapsed here.
    Bottom: the weighted trees corresponding to these strata.}
  \label{fig:associaplex-conv}
\end{figure}

The present section provides a geometric interpretation of the chain
complex of weighted trees.  Specifically, we define here a CW complex,
the \emph{weighted associaplex} $\Associaplex{n}{w}$, with one cell
$Q^T$ for each stably-weighted tree $T$, and whose differential agrees
with the differential on the weighted trees complex
$\wTreesCx{n}{w}$. In the
special case $w=0$, the associahedron serves as such a CW complex (see
for example~\cite{Devadoss98:tessellations}): indeed, the
associahedron is a polytope. For $w\geq 1$, however, the weighted
associaplex is not a polytope, though its underlying topological space
is homeomorphic to a disk.  Although the material here is logically
independent of the rest of the paper, the reader may find this
topological perspective useful as motivation.

We start with a brief description of the weighted associaplex, based
on the Deligne-Mumford compactification of the moduli space of points
on a disk, and then give a more detailed explanation using basic
complex analysis. 
Consider the configuration space of $w$ distinct,
unordered points in the interior of $D^2$ and $n+1$ points in
$\bdy D^2$, one of which is distinguished. The Deligne-Mumford compactification of this space has
strata corresponding to disks and spheres glued together at nodes
(which are additional marked points), so that:
\begin{itemize}
\item at each boundary node, two disks are glued together;
\item at each interior node either two spheres are glued together, or a sphere is glued to a disk;
\item smoothing all of the nodes gives a single topological disk (so,
  for instance, the disks are glued together according to some tree);
  and
\item each sphere has at least 3 special points (marked points or
  nodes) and each disk has at least three boundary special points or
  an interior special point and a boundary special point.
\end{itemize}
The stratum corresponding to such a configuration is the product over
the disks and spheres $\Sigma$ of the Deligne-Mumford moduli space of
$\Sigma$. Then the weighted associaplex with weight $w$ and $(n+1)$
boundary marked points is obtained by collapsing each factor
corresponding to a sphere $\Sigma$ to a single point. It is
straightforward to see that this collapsing respects the gluing
between strata. Roughly, this means that in the associaplex, we do not
keep track of the directions in which interior marked points collide.

Turning to the more detailed description, let $\HHH$ denote the lower
half plane $\{a+bi\mid a,b\in\RR,\ b\leq 0\}\subset \CC$, and let
$\oHHH$ denote its interior.
The conformal automorphism group $\Aut(\HHH)$ is the group
of affine transformations $z\to a z + b$, where $a,b\in\RR$ and $a>0$.
(Under the conformal identification of $\HHH$ with
the unit disk, this automorphism group corresponds to
the set
of M\"obius transformations that preserve $-i$.)
Let
$\Sym^{n,w}(\HHH)$ denote the product of the
space of order-preserving embeddings
$\{1,\dots,n\}\hookrightarrow \RR=\partial\HHH$
with the 
$w\th$ symmetric product
of $\oHHH$.  Our goal is to
define, for integers $n\geq 0$ and $w\geq 0$ with $n+2w\geq 2$, a
useful compactification of $\SymQ^{n,w}\coloneqq
\Sym^{n,w}(\HHH)/\Aut(\HHH)$. We call this compactification the {\em weighted
  associaplex} or just the \emph{associaplex}. 

Strata in the associaplex are indexed by stably-weighted trees.
Precisely, let $T$ be a stably-weighted tree with $n$ inputs and
weight $w$.
Let $\valence(v)$ denote the
valence of an internal vertex $v$ and let $w(v)$ denote the weight of
$v$. The open stratum corresponding to $T$ is
\[
  \SymQ^T\coloneqq \prod_{v\in\Vertices(T)}\SymQ^{\valence(v)-1,w(v)}.
\]

The space $Q^{n,w}$ is an open subspace of the
quotient of a space of ordered $(n+w)$-tuples of points,
\begin{equation}\label{eq:lQ}
  \lQ{n,w}\subset \left(\prod_{i=1}^n \partial \HHH\right)
  \times \left(\prod_{i=1}^w \oHHH\right) \subset \prod^{n+w}\HHH,
\end{equation}
by $\Aut(\HHH)$ and then by the action of the symmetric group
on $w$ letters. In defining $\lQ{n,w}$ we assume that the $n$
points in the boundary come first in the ordering, and that their
ordering induced by $\RR$ agrees with the chosen ordering;
correspondingly, in Equation~\eqref{eq:lQ}, a point
$(z^1,\dots,z^{n+w})$ satisfies $z^1<z^2<\cdots<z^n\in\RR$.

We will denote elements of $\prod_{i=1}^{n+w}\HHH$ with bold letters.
Let $z^i\in\HHH$ denote the $i^{th}$ component of
$\mathbf{z}\in\lQ{n,w}$.  We will also consider, for each weighted tree $T$,
the space
\[
  \lQ{T} =\prod_{v\in\Vertices(T)} \lQ{\valence(v)-1,w(v)}.
\]

Often, we will find it convenient to work with the one-point compactification
$\cHHH$ of $\HHH$.
Given a point $p\in \oHHH$, let $\phi^p\in\Aut(\HHH)$ denote the map
with $\phi^p(p)=-i$. Explicitly, if $p=a+bi$ then 
\[
  \phi^p(z)=\frac{a-z}{b}.
\]
Let
\[
  \Phi^{p}\co \prod^{n+w}\cHHH\to\prod^{n+w}\cHHH
\]
denote the induced map on the Cartesian product.

\begin{definition}
  \label{def:SliceCondition}
  A  tuple of points $(z_1,\dots,z_\ell)$ in $\HHH$,
  satisfies the {\em slice condition} if the following two properties hold:
  \begin{itemize}
  \item $\sum_{i=1}^{\ell}z_i$ have vanishing real part and
  \item  $\sum_{i=1}^{\ell}|z_i|=\ell$.
  \end{itemize}
\end{definition}
If $(z_1,\dots,z_\ell)$ contains at least two distinct entries,
or one entry which is not real, 
then there is a unique $p\in\oHHH$ so that $(\phi^p(z_1),\dots,\phi^p(z_\ell))$ satisfies the
slice condition.

\begin{definition}
  Given a sequence $\seq{z^1_j}$ in $\oHHH$, a {\em stabilizing
    sequence} $\seq{p_j}$ is a sequence in $\oHHH$ with the
  property that the sequence $\seq{\phi^{p_j}(z^1_j)}$ converges to a
  point in $\oHHH$.  Given a sequence $\seq{(z^1_j,z^2_j)}$ of pairs of points
  in $\HHH$  that are required to be distinct
  if both are on $\partial \HHH$, a
  stabilizing sequence $\seq{p_j}$ is a sequence in $\oHHH$ with the
  property that $\seq{(\phi^{p_j}(z^1_j),\phi^{p_j}(z^2_j))}$
  converges to pair of points $(z^1,z^2)$ in $\HHH$, which are
  distinct if both are on $\partial \HHH$.
\end{definition}

Any sequence $\seq{z^1_j}$ in $\oHHH$ has a stabilizing sequence:
let $p_j=z^1_j$. For pairs of points, we have the following:

\begin{lemma}
  \label{lem:ConstructStabilizingSequence}
  A sequence $\seq{(z^1_j,z^2_j)}$ of pairs of points in $\HHH$ that
  are distinct if they are both in $\partial\HHH$ has a subsequence
  that admits a stabilizing sequence.
\end{lemma}
\begin{proof}
  For each $j$, 
  choose $q_j$ so that for $(\phi^{q_j}(z^1_j),\phi^{q_j}(z^2_j))$
  satisfies the slice condition.
  By compactness of the unit ball in $\HHH$, we can find a subsequence
  of $\seq{(\phi^{q_j}(z^1_j),\phi^{q_j}(z^2_j))}$
  that converges to a pair of points $(z^1,z^2)$ in $\HHH$
  that also satisfy the slice condition. 
\end{proof}

\begin{lemma}
  \label{lemma:UniqueAutomorphism}
  If $\seq{p_j}$ and $\seq{q_j}$ are two stabilizing sequences
  for the same sequence $\seq{z_j}$ with $z_j\in\oHHH$
  or for the same sequence $\seq{(z^1_j,z^2_j)}$ of pairs
  of points in $\HHH$,
  then $\phi^{p_j}\circ (\phi^{q_j})^{-1}$ converges to some automorphism of
  $\HHH$.
\end{lemma}

\begin{proof}
  Suppose that $\seq{p_j}$ and $\seq{q_j}$ are stabilizing
  sequences for the same sequence of pairs $\seq{(z^1_j,z^2_j)}$.
  It suffices to consider the case that $\seq{q_j}$ is the stabilizing sequence from
  the proof of 
  Lemma~\ref{lem:ConstructStabilizingSequence}; i.e. 
  $(\phi^{q_j}(z^1_j),\phi^{q_j}(z^2_j))$ satisfies the slice condition
  for each $j$.
  Replacing $z^k_j$ with $\phi^{q_j}(z^k_j)$ for $k=1,2$, we
  can assume that the pair $(z^1_j,z^2_j)$ satisfies the slice
  condition and $\phi^{q_j}=\Id$ for all $j$.
  
  Suppose first that $\seq{p_j}$ converges to a point $p\in\cHHH$.
  If $p\in\oHHH$, then $\seq{\phi^{p_j}}\goesto\phi^p$. If $p\not\in\oHHH$, we claim that $\seq{p_j}$ is not
  a stabilizing sequence. Specifically, if $p=\infty$, then it is 
  easy to see that
  either
  \[
  \lim_{j\goesto\infty}\phi^{p_j}(z^1_j)=\lim_{j\goesto\infty}\phi^{p_j}(z^2_j)\in\RR \cup\infty. \]
  (The limit is real when  $\seq{p_j}$ converges.)
  Suppose next that $p\in\partial\HHH$.  If $p\neq \lim_{j\goesto\infty}\phi^{p_j}(z^k_j)$ for $k=1,2$, then it is
  easy to see that $\phi^{p_j}(z^k_j)\goesto \infty$.

  Finally, if the
  sequence $\seq{p_j}$ has more than one limit point then the argument
  above shows that each limit point must be in $\oHHH$ and, since the
  sequence $\seq{(\phi^{p_j}(z^1_j),\phi^{p_j}(z^2_j))}$ converges, in
  fact all the limit points must be the same.
  
  A similar argument applies in the case where $\seq{p_j}$ and $\seq{q_j}$
  are stabilizing sequences for the same sequence $\seq{z_j}$ with $z_j\in\oHHH$.
\end{proof}

A \emph{stable configuration} is a point in $\prod^{w+n}\HHH$ that
either contains at least one component in the interior of $\HHH$ or
at least two components that are distinct points on the boundary.
Fix a sequence $\seq{\mathbf{z}_j}$ in $\lQ{n,w}$.  A {\em stabilizing
  sequence} for $\seq{\mathbf{z}_j}$
is a sequence $\seq{p_j}$ of points in $\oHHH$, with the
following properties:
\begin{itemize}
\item The sequence of points $\seq{\Phi^{p_j}(\mathbf{z}_j)}$ converges to
  some point $\z\in\prod^{n+w}\cHHH$.
\item The limit point $\z$ is a stable configuration.
\end{itemize}
Elements $\z$ that arise as $\lim \seq{\Phi^{p_j}(\mathbf{z}_j)}$ for
some choice of stabilizing sequence $\seq{p_j}$, are called {\em stable limits}.

Call two stabilizing sequences for $\seq{\mathbf{z}_j}$ {\em equivalent} if
$\lim_{j\goesto  \infty}\Phi^{p_j}(\z_j)$ and $\lim_{j\goesto  \infty}\Phi^{q_j}(\z_j)$ differ by $\Aut(\HHH)$ or, 
more precisely, if there is a $\psi\in\Aut(\HHH)$ inducing
a map $\Psi$ on $\prod^{n+w}\cHHH$ such that
\[
\Psi\left(\lim_{j\goesto\infty} \Phi^{p_j}(\mathbf{z}_j)\right)
=\lim_{j\goesto\infty} \Phi^{q_j}(\mathbf{z}_j).\]
We denote the equivalence class of a stabilizing sequence $\seq{p_j}$ by
$\llbracket \seq{p_j}\rrbracket$.

\begin{definition}
  Fix a sequence of points $\seq{\mathbf{z}_j}$ in $\lQ{n,w}$.  Let
  $\seq{p_j}$ be any stabilizing sequence $\seq{p_j}$ for
  $\seq{z^i_j}$, for some $i\in\{n+1,\dots,n+w\}$, or for
  $\seq{(z^i_j,z^k_j)}$ for some pair of $i,k\in\{1,\dots,n+w\}$.  The
  sequence $\seq{\mathbf{z}_j}$ is called {\em pre-convergent} if for
  any $\seq{p_j}$ as above,
  $\seq{\Phi^{p_j}({\mathbf z}_j)}$ converges in $\prod^{n+m}\cHHH$.
\end{definition}

We will show that each pre-convergent sequence in $\lQ{n,w}$ converges
to a well-defined limit point in $\lQ{T}$ for some tree $T$.

\begin{lemma}\label{lem:equiv-stable-seqs}
  Let $\seq{\mathbf{z}_j}$ be a sequence in $\lQ{n,w}$, 
  suppose $\seq{p_j}$ and $\seq{q_j}$ are two stabilizing sequences
  for $\seq{\mathbf{z}_j}$, and suppose one
  of the following two conditions holds:
  \begin{itemize}
  \item For some $i\in\{1,\dots,n+w\}$,
    $\lim_{j\goesto\infty} \phi^{p_j}(z^i_j)$ and
    $\lim_{j\goesto\infty} \phi^{q_j}(z^i_j)$ both lie in $\oHHH$; or
  \item There are two distinct $i,k\in\{1,\dots,n+w\}$
    so that
    $\lim_{j\goesto\infty} \phi^{p_j}(z^i_j)$ and
    $\lim_{j\goesto\infty} \phi^{p_j}(z^k_j)$ are two distinct points
    in $\partial\HHH$, as are $\lim_{j\goesto\infty} \phi^{q_j}(z^i_j)$ and
    $\lim_{j\goesto\infty} \phi^{q_j}(z^k_j)$.
  \end{itemize}
  Then, $\seq{p_j}$ and $\seq{q_j}$ are equivalent.
\end{lemma}

\begin{proof}
  This follows from Lemma~\ref{lemma:UniqueAutomorphism}.
\end{proof}

\begin{definition}
  \label{def:InputSet}
  The {\em input set} of a stabilizing sequence $\seq{p_j}$,
  denoted $\InputSet(\seq{p_j})$, is the set
  of $i\in\{1,\dots,n+w\}$ with the property that
  the sequence $\seq{\phi^{p_j}(z^i_j)}$ converges to
  some point in $\HHH$ (i.e., the sequence does not diverge to infinity).
\end{definition}

\begin{lemma}
  \label{lem:InputSet}
  Let $\seq{p_j}$ and $\seq{q_j}$ be two stabilizing sequences
  for some sequence $\seq{{\mathbf z}_j}$, and suppose
  that there is $r\in\partial\HHH$ so that (at least) one of the following
  two conditions holds:
  \begin{itemize}
  \item for some $i\in\InputSet(\seq{p_j})$,
    \begin{equation}
      \label{eq:Hypotheses0}
      \lim_{j\goesto\infty}\phi^{p_j}(z^i_j)\in\oHHH
      \qquad\text{and}\qquad
      \lim_{j\goesto\infty}\phi^{q_j}(z^i_j)=r
    \end{equation}
    \item for some $i,k\in\InputSet(\seq{p_j})$,
      \begin{equation}
        \label{eq:Hypotheses1}
        \lim_{j\goesto\infty}\phi^{p_j}(z^i_j)\neq
        \lim_{j\goesto\infty}\phi^{p_j}(z^k_j); \qquad{\text{and}}\qquad
        \lim_{j\goesto\infty}\phi^{q_j}(z^i_j)=
        \lim_{j\goesto\infty}\phi^{q_j}(z^k_j)=r.
      \end{equation}
  \end{itemize}
  Then, for all $\ell\in I$,
  \begin{equation}
    \label{eq:AllConvergeToR}
    \lim_{j\goesto\infty} \phi^{q_j}(z^\ell_j)=r.
  \end{equation}
  Thus, any two stabilizing sequences with the same input set are equivalent.
\end{lemma}

\begin{proof}
  Write $\phi^{q_j}\circ (\phi^{p_j})^{-1}(z)=a_j z + b_j$.
  Either of the conditions from 
  Equations~\eqref{eq:Hypotheses0}
  and~\eqref{eq:Hypotheses1} implies that
  $\seq{a_j}\goesto 0$ and $\seq{b_j}\goesto r$,
  which in turn implies that if $\seq{\phi^{p_j}(z_j)}$ converges in $\HHH$,
  then $\seq{\phi^{q_j}\circ (\phi^{p_j})^{-1}(z_j)}\goesto r$,
  verifying Equation~\eqref{eq:AllConvergeToR}.

  Suppose the input set $I$ for $\seq{p_j}$ is contained in the input
  set $J$ of $\seq{q_j}$, but the sequences are inequivalent. Then, by
  Lemma~\ref{lem:equiv-stable-seqs}, at least one of the two
  conditions from Equations~\eqref{eq:Hypotheses0}
  and~\eqref{eq:Hypotheses1} holds.  Equation~\eqref{eq:AllConvergeToR}
  and the definition of a stable configuration ensures that $I$ is a
  proper subset of $J$.
\end{proof}

\begin{corollary}\label{cor:fin-many-stab-seq}
  There are only finitely many different equivalence classes of
  stabilizing sequences for a given sequence $\seq{\z_j}$.
\end{corollary}
\begin{proof}
  There are only finitely many possible input sets and, by
  Lemma~\ref{lem:InputSet}, the equivalence class of a stabilizing
  sequence is determined by its input set.
\end{proof}

\begin{lemma}
  \label{lem:PreConvergentSequence}
  Let $\seq{{\mathbf z}_j}$ be a pre-convergent sequence
  and $\seq{q_j}$ be a stabilizing sequence for a subsequence 
  $\seq{{\mathbf z}_{n_j}}$. Then, there is a stabilizing sequence
  $\seq{p_j}$ for $\seq{{\mathbf z}_j}$ with the same input set as
  $\seq{q_j}$.
\end{lemma}

\begin{proof}
  Suppose that $\seq{q_j}$ is a stabilizing sequence for
  $\seq{{\mathbf z}_{n_j}}$.  This means we can find either some
  $i\in\{n+1,\dots,n+w\}$ so that $\seq{\phi^{q_j}(z^i_{n_j})}$
  converges to a point in $\oHHH$, or we can find distinct
  $i,k\in\{1,\dots,n+w\}$ so that
  $\seq{(\phi^{q_j}(z^i_{n_j}),\phi^{q_j}(z^k_{n_j}))}$ converges to a
  pair of distinct points in $\partial \HHH$.  Let $\seq{p_j}$ be a
  stabilizing sequence for $\seq{z^i_j}$ (in the first case)
  or for 
  $\seq{(z^i_j,z^k_j)}$ (in the second case).
  Lemma~\ref{lemma:UniqueAutomorphism} shows that $\seq{q_{n_j}}$ and
  $\seq{p_{n_j}}$ are equivalent sequences.  In particular, they have
  the same input sets. For a pre-convergent sequence,
  the input set of a
  stabilizing sequence $\seq{p_j}$ is the same as the input set of any
  of its subsequences.
\end{proof}

\begin{lemma}
  \label{lem:IncomingEdge}
  Suppose that $\seq{q_j}$ is a stabilizing sequence for some sequence
  $\seq{{\mathbf z}_j}$ in $\lQ{n,w}$, and fix $r\in\partial \HHH$
  with the property that
  \[
    I=\{i\in\{1,\dots,n+w\}\mid\lim_{j\goesto\infty}\phi^{q_j}(z^i_j)=r\}
  \]
  contains at least two elements in $\{1,\dots,n\}$ (i.e., for which
  $z^i_j\in\partial \HHH$) or at least one element of
  $\{n+1,\dots,n+w\}$ (i.e., for which $z^i_j\in\oHHH$). Then there is
  a subsequence $\seq{{\mathbf z}_{n_j}}$ which has a stabilizing
  sequence $\seq{p_j}$ whose input set is $I$.  In particular, if the
  original sequence $\seq{{\mathbf z}_j}$ is pre-convergent, then
  there is another stabilizing sequence with input set $I$.
\end{lemma}

\begin{proof}
  For each $j$, let $I_j= \{z^i_j\mid i\in I\}$.  Choose $t_j$ so that
  $\phi^{t_j}(I_j)$ satisfies the slice condition.
  By compactness of the unit ball in $\HHH$, we can find
  a subsequence $\seq{t_{n_j}}=\seq{p_j}$ that is a stabilizing sequence
  for $\seq{{\mathbf z}_{n_j}}$.
  Let $J$ be the input set of
  $\seq{q_j}$.  By construction, $I\subseteq J$. By
  Lemma~\ref{lem:InputSet}, $J\subseteq I$; so $J=I$.
  The final statement in the lemma now follows from Lemma~\ref{lem:PreConvergentSequence}.
\end{proof}

\begin{lemma}
  \label{lem:UniqueInput}
  Let $\seq{{\mathbf z}_j}$ be a pre-convergent sequence in $\lQ{n,w}$.
  Given $i\in\{1,\dots,n\}$ (i.e., for which $z^i_j\in\partial \HHH$),
  there is a stabilizing sequence $\seq{p_j}$ with the
  property that for all $k\neq i$ with $k\in \{1,\dots,n+w\}$,
  \begin{equation}
    \label{eq:iIsUnique}
    \lim_{j\goesto\infty}\phi^{p_j}(z^i_j)\neq
    \lim_{j\goesto\infty}\phi^{p_j}(z^k_j).
  \end{equation}
  Moreover, $\seq{p_j}$ is
  uniquely characterized up to equivalence by that property.
 \end{lemma}
 
\begin{proof}
  To construct the desired stabilizing sequence, we argue as follows.
  Let ${\mathcal P}$ be the set of equivalence classes of stabilizing
  sequences whose input set contains $i$.

  First, we argue that ${\mathcal P}$
  is
  non-empty. Consider the sequence
  $\seq{t_j}$ so that
  $(\phi^{t_j}(z^1_j),\dots,\phi^{t_j}(z^{n+w}_j))$ satisfies the
  slice condition. As in the proof of Lemma~\ref{lem:IncomingEdge},
  this has a subsequence that is a stabilizing sequence for a subsequence
  of $\seq{{\mathbf z}_j}$ with input set
  $\{1,\dots,n+w\}$. Lemma~\ref{lem:PreConvergentSequence} now gives
  a stabilizing sequence for all of ${\mathbf z}_j$ with input set $\{1,\dots,n+w\}$.

  Next, given $v_1,v_2\in {\mathcal P}$, we declare
  $v_1\prec v_2$ if the input set of $v_1$ is contained in the input set of
  $v_2$.  This induces a partial ordering on the elements of
  ${\mathcal P}$. Lemma~\ref{lem:IncomingEdge} implies that
  any element $\llbracket \seq{p_j}\rrbracket$
  of ${\mathcal P}$ that is minimal with respect to this ordering
  satisfies Equation~\eqref{eq:iIsUnique} for all $k\neq i$.
  
  For uniqueness, we argue as follows.
  Let $\seq{p_j}$ and $\seq{q_j}$ be stabilizing sequences
  satisfying the hypothesis, with input sets $I$ and $J$ respectively.
  Then, we can find $k\in I$ with
  \[ \lim_{j\goesto\infty}\phi^{p_j}(z^i_j)\neq
  \lim_{j\goesto\infty}\phi^{p_j}(z^k_j),\]
  and $\ell\in J$ so that
  \[ \lim_{j\goesto\infty}\phi^{q_j}(z^i_j)\neq
  \lim_{j\goesto\infty}\phi^{q_j}(z^\ell_j).\]
   Choose a stabilizing sequence
   $\seq{s_j}$ so that
   $(\phi^{s_j}(z^i_j),\phi^{s_j}(z^k_j),\phi^{s_j}(z^\ell_j))$
   satisfies the slice condition.
   
   By the slice condition, at least two of 
   \[\lim_{j\goesto\infty}\phi^{s_j}(z^i_j),
   \qquad
   \lim_{j\goesto\infty}\phi^{s_j}(z^k_j),\qquad\text{and}
   \qquad
   \lim_{j\goesto\infty}\phi^{s_j}(z^\ell_j) \]
   are distinct.
   If 
   \begin{equation}
     \label{eq:TwoDistinct}
     \lim_{j\goesto\infty}\phi^{s_j}(z^i_j) \neq
     \lim_{j\goesto\infty}\phi^{s_j}(z^k_j),
   \end{equation}
   then by Lemma~\ref{lem:equiv-stable-seqs}, $\seq{p_j}$ and
   $\seq{s_j}$ are equivalent. In particular, $\ell\in I$; so
   Equation~\eqref{eq:iIsUnique} (for $k=\ell$) gives
   $\lim_{j\goesto\infty}\phi^{p_j}(z^i_j)\neq
   \lim_{j\goesto\infty}\phi^{p_j}(z^\ell_j)$. By
   Lemma~\ref{lem:equiv-stable-seqs}, we conclude that $\seq{p_j}$ is
   equivalent to $\seq{q_j}$, provided Equation~\eqref{eq:TwoDistinct}
   holds.  An analogous argument proves that $\seq{p_j}$ and
   $\seq{q_j}$ are equivalent when Equation~\eqref{eq:TwoDistinct} is
   replaced by
   \[
     \lim_{j\goesto\infty}\phi^{s_j}(z^i_j) \neq
     \lim_{j\goesto\infty}\phi^{s_j}(z^\ell_j).
   \]
   This completes the proof.
 \end{proof}

Fix a pre-convergent sequence $\seq{\mathbf{z}_j}$.
Construct a weighted planar graph $T$ as follows.  The vertices of $T$
correspond to equivalence classes of stabilizing sequences for $\seq{\mathbf{z}_j}$.
The weight of a vertex is the number of
$i\in\{n+1,\dots,n+w\}$ with the property that
$\seq{\phi^{p_j}(z^i_j)}$ converges to some point in $\oHHH$.  The
edges of $T$ coming into a vertex $\llbracket \seq{p_j}\rrbracket$
correspond to the distinct $q\in\partial \HHH$ that arise as limit
points of $\seq{\phi^{p_j}(z^i_j)}$.  The planar ordering of the edges
into each vertex is the one induced by the ordering of these points in
$\RR$.  The input edges into $T$
whose endpoint is at the vertex $\llbracket \seq{p_j}\rrbracket$
correspond to 
$i$ for which
$\lim_{j\goesto\infty}{\phi^{p_j}(z^i_j)}\in\partial\HHH$ is distinct from all
$\lim_{j\goesto\infty}{\phi^{p_j}(z^k_j)}$ for all $k\neq i$.
For each edge pointing into the vertex $\llbracket\seq{p_j}\rrbracket$
that is not an input edge,
Lemma~\ref{lem:IncomingEdge} constructs the other endpoint of
that edge.

Orient the edges of $T$ from
the vertex with smaller input set to the one with larger input set.

\begin{definition}
  An {\em increasing path} is a path between two vertices of $T$ obtained by
  juxtaposing edges oriented so that the input set is always increasing.
\end{definition}

\begin{lemma}
  \label{lem:UniquePath}
  Let $\seq{{\mathbf z}_j}$ be a pre-convergent sequence, and let $T$ be
  the associated weighted graph as above.
  Let $v_1$  and $v_2$ be two vertices of $T$
  whose stabilizing sequences have input sets $I$ and $J$ respectively,
  with $I\subset J$.
  Then, there is a unique increasing path from $v_1$ to $v_2$.
\end{lemma}

\begin{proof}
  The construction and proof of uniqueness is inductive on
  $|J\setminus I|$.  If $I=J$, the sequences are equivalent by
  Lemma~\ref{lem:equiv-stable-seqs} and the path is necessarily
  trivial.
  Otherwise,
  Lemma~\ref{lem:InputSet}, provides an $r\in\RR$ so that
  for all $i\in I$, $\phi^{q_j}(x^i_j)=r$. This clearly determines the
  final edge $e$ in any increasing path from $\seq{p_j}$ to
  $\seq{q_j}$, since all other edges into
  $\llbracket\seq{q_j}\rrbracket$ point from vertices that do not
  contain $I$ in their input set.
  Lemma~\ref{lem:IncomingEdge} supplies the initial endpoint
  $\llbracket\seq{s_j}\rrbracket$ of the edge $e$.
  The inductive hypothesis constructs uniquely the path from $\seq{p_j}$ to $\seq{s_j}$;
  the desired path is obtained by appending the edge $e$.
\end{proof}

\begin{proposition}
  For a pre-convergent sequence $\seq{{\mathbf z}_j}$, the associated
  weighted graph $T$ is  a finite, planar tree with $n$ inputs, weight
  $w$, and with no valence two, weight-0 internal vertices.
\end{proposition}

\begin{proof}
  Corollary~\ref{cor:fin-many-stab-seq} ensures that $T$ has
  finitely many vertices.  Orient each edge from
  $v_1=\llbracket\seq{q_j}\rrbracket$ to
  $v_2=\llbracket\seq{p_j}\rrbracket$ if the input set of $v_1$ is
  properly contained in the input set of $v_2$.
  If $I(v_1)\subset I(v_2)$ then Lemma~\ref{lem:UniquePath} constructs
  an oriented path from $v_1$ to $v_2$, and that path is unique.
  
  The trunk vertex $\trunkvertex$ (i.e., the predecessor to the root) is the vertex whose input set is all of
  $\{1,\dots,n+w\}$, which can be realized by the stabilizing sequence that
  achieves the slice condition for $(z^1_j,\dots,z^{n+w}_j)$.
  Lemma~\ref{lem:UniquePath} gives a path from any vertex to this trunk vertex;
  thus, the graph is connected.
  Uniqueness of the path, together with Lemma~\ref{lem:UniqueInput},
  shows that each input points into a unique internal vertex,
  and ensures that the graph is a tree.

  The stability condition eliminates the possibility of valence two,
  weight-0 internal vertices.

  To verify that the total weight is $w$, observe that for any $i$ so
  that $\seq{z^i_j}$ has $z^i_j\in \oHHH$ for all $j$, the sequence
  $\seq{z^i_j}$ is a stabilizing sequence. By
  Lemma~\ref{lem:equiv-stable-seqs}, any other stabilizing sequence
  such that $\seq{\phi^{p_j}(z^i_j)}$ converges to some point in
  $\oHHH$ is equivalent to this stabilizing sequence, so there is a
  unique vertex where $\seq{z^i_j}$ contributes to the weight.
  
  The fact that the tree has $n$ inputs follows from
  Lemma~\ref{lem:UniqueInput}. Any ordering of the edges into each
  vertex gives a tree a planar structure so, in particular, the
  ordering specified above does. It is immediate from the construction
  that the induced ordering of the inputs agrees with the ordering
  $1,\dots,n$.
\end{proof}

Having identified the stratum $T$ that a pre-convergent sequence
$\seq{{\mathbf z}_j}$ converges
to, we explain how to find the actual limit point.  Each vertex
$v\in \Vertices(T)$ corresponds to an equivalence class of stabilizing
sequence $\seq{p_j}$. Thought of as points in $\prod^{w+n}\cHHH$, the
sequence $\phi^{p_j}(\z_j)$ converges to some
$\z^v\in\Prod{w+n}(\cHHH)$, exactly $\valence(v)-1$
many components of which
are distinct points on $\partial\HHH$.
Thus, it represents a point in
$\lQ{\valence(v)-1,w(v)}$.  Clearly, if $\seq{q_j}$
is equivalent to $\seq{p_j}$,
then $\phi^{q_j}(\z_j)$ represents the same point
modulo the action of $\Aut(\HHH)$.
Extracting such limits for all vertices $v$, we
obtained the promised configuration ${\mathbf z}\in \lQ{T}$.

Let 
\begin{equation}
  \lAssoc{n}{w}=\bigcup_{T\in\wTrees{n}{w}}\lQ{T}.
\end{equation}
Given a vertex $v$ of $T$, let $\pi_v\co \lQ{T}\to \lQ{\valence(v)-1,w(v)}$
denote projection.
Define a sequence $\seq{\z_j}$ in $\lQ{T}$ to be convergent
if for each vertex $v$ of $T$, the sequence $\seq{\pi_v(\z_j)}$ is
pre-convergent. The limit of such a convergent sequence $\seq{\z_j}$ has
underlying tree $T'$ obtained by splicing in the tree for
$\lim \seq{\pi_v(\z_j)}$ at the vertex $v$ (for each $v$), and the limit of $\seq{\z_j}$ in
$\lQ{T'}$ is given by the limits of the points $\seq{\pi_v(\z_j)}$ in the
obvious way. Finally, call an arbitrary sequence $\seq{\z_j}$ in $\lAssoc{n}{w}$
\emph{convergent} if the intersection of $\seq{\z_j}$ with each stratum $\lQ{T}$
is either finite or a convergent sequence, and all of these convergent sequences
converge to the same limit.

A topological space $X$ is \emph{sequentially Hausdorff} if any convergent
sequence in $X$ has a unique limit.
\begin{proposition}\label{prop:sequential-top}
  The above notion of convergence gives the point set
  \[\lAssoc{n}{w}=\bigcup_{T\in\wTrees{n}{w}}\lQ{T}\]
  the structure of a sequentially compact, sequentially Hausdorff
  topological space. Further, the topology on each subspace
  $\lQ{T}$ is the standard topology on the open disk.
\end{proposition}

\begin{proof}
  We first show that $\lAssoc{n}{w}$ is sequentially Hausdorff, i.e.,
  that any convergent sequence as defined above has a unique limit. It
  follows that the convergent sequences induce a topology on
  $\lAssoc{n}{w}$. We then verify that this topology is sequentially
  compact.

  For a pre-convergent sequence in $\lQ{n,w}$, the above construction
  gives a tree $T$ and a unique limit point in $\lQ{T}$. In
  particular, convergent sequences in $\lQ{n,w}$ have unique limits. This implies
  immediately that any convergent sequence contained in a single open stratum has
  a unique limit. Since the limit of a sequence that intersects
  several open strata is the common limit of the subsequences in each
  stratum, such a sequence also has a unique limit.
  
  Now, define a topology on $\lAssoc{n}{w}$ by declaring that a subset
  $C\subset \lAssoc{n}{w}$ is closed if every sequence
  $\seq{\z_j}\subset C$ which converges to a point
  $\z\in \lAssoc{n}{w}$ has $\z\in C$. It is immediate from the
  construction that $\emptyset$ and $\lAssoc{n}{w}$ are closed and
  that arbitrary intersections of closed sets are closed. The fact
  that finite unions of closed sets are closed follows from sequential
  Hausdorffness and the fact that, for our collection of convergent
  sequences, a subsequence of a convergent sequence converges.
  It also follows from the fact that subsequences of convergent
  sequences are convergent that $U\subset\lAssoc{n}{w}$ is open (i.e.,
  $\lAssoc{n}{w}\setminus U$ is closed) if and only if every sequence
  which converges to a point in $U$ is eventually contained in $U$.
  
  To prove this topology is sequentially compact, we need to show that
  if $\seq{{\mathbf z}_j}$ is any sequence in $\lQ{n,w}$, then it has
  a pre-convergent subsequence.  For each
  $\{i,k\}\subset\{1,\dots,n+w\}$,
  Lemma~\ref{lem:ConstructStabilizingSequence} gives a stabilizing
  sequence $\seq{p_j}$ for $\seq{(z^i_j,z^k_j)}$. By compactness of
  $\cHHH$, there is a subsequence so that
  $\seq{\Phi^{p_{n_j}}({\mathbf z}_{n_j})}$ converges in
  $\prod^{n+w}(\cHHH)$. Similarly, for each $i\in\{n+1,\dots,n+w\}$,
  there is a stabilizing sequence $\seq{p_j}$
  for $\seq{z^i_j}$ (which we can take to be $p_j=z^i_j$).
  Again, we can find a subsequence so that
    $\seq{\Phi^{p_{n_j}}({\mathbf z}_{n_j})}$ converges in
  $\prod^{n+w}(\cHHH)$. Take a common subsequence of the $\seq{{\mathbf z}_j}$
  over the finitely many choices
  for $\{i,k\}\subset \{1,\dots,n+w\}$ and $i\in\{n+1,\dots,n+w\}$.
  This is pre-convergent because by Lemma~\ref{lemma:UniqueAutomorphism},
  any other stabilizing sequence for $\seq{{\mathbf z}_j}$ is equivalent to
  one of the finitely many stabilizing sequences considered above.

  Finally, we show that the subspace topology on $\lQ{T}$ agrees with
  the standard topology on the open disk.  If $U$   
  is open in $\lQ{T}$ as a subspace of $\lAssoc{n}{w}$, so $U=V\cap \lQ{T}$ where
  $V\subset \lAssoc{n}{w}$ is some open set, then any sequence in
  $\lQ{T}$ which converges to a point in $U$ is eventually contained
  in $U$, hence $U$ is open in the standard topology on $\lQ{T}\cong
  D^{\dim(T)}$. Conversely, if $U$ is open in the standard topology
  then let $V$ be the union of $U$ and all of the open strata whose
  closures strictly contain $U$ (i.e., all $\lQ{S}$ where $T$ is
  obtained by inserting edges in $S$). It is immediate from the
  definition of convergence that $V$ is open. Hence, $U$ is open in
  the subspace topology.
\end{proof}

Let $T\in\wTrees{n}{w}$. Let $E$ be a subset of the internal edges of
$T$ and let
$K_{E}(T)\in\wTrees{n}{w}$ be the tree obtained from $T$ by
contracting all the edges in $E$. Let $(0,\infty)^{E}$ denote the set
of maps from $E$ to $(0,\infty)$.

\begin{lemma}
  \label{lem:lAssociaplexCorners}
  Given any $T\in\wTrees{n}{w}$ and any subset $E$ of the internal edges of $T$,
  there is a map
  \[
    \varphi\co (0,\infty)^{E} \times \lQ{T}  \to \lQ{K_E(T)}
    \]
    which is a homeomorphism onto its image and 
  which extends continuously to a map
  \[
    {\overline\varphi}\co (0,\infty]^{E}  \times \lQ{T} \to \lAssoc{n}{w}
      \]
      which is also a homeomorphism onto its image,
      with the property that
      ${\overline\varphi}_{\{\infty\}^{E}\times \lQ{T}}$ is
      the inclusion of $\lQ{T}\subset \lAssoc{n}{w}$.
\end{lemma}

\begin{proof}
  For $E=\Edges(T)$, we write a formula for $\varphi$.
  The case that $E$ is a proper
  subset of $\Edges(T)$ is obtained by applying this formula to each
  sub-tree spanned by the edges in $E$.
  
  A point in $(0,\infty)^E\times \lQ{T}$ can be thought of as consisting of a
  pair:
  \begin{itemize}
  \item A function $r\co E \to (0,\infty)$;
  \item For each vertex $v$, a tuple
    $(z_i^{(v)}) \in \prod_{i\in \InputSet(v)} \HHH$ satisfying the
    slice condition.
  \end{itemize}
  Note that $z_i^{(v)}\in \partial \HHH$ unless $i\in\{n+1,\dots,n+w\}$
  and $v$ is the vertex which is furthest from the trunk vertex $\trunkvertex$,
  among all vertices containing $i$ in their input set.

  The image under $\varphi$ of such a pair $(r,(z_i^{(v)})_{v\in\Vertices(T),\
  i\in\InputSet(v)})$ is an $(n+w)$-tuple of points in $\HHH$
  whose $i^{th}$ component is defined as follows.  Let
  $v_1,\dots,v_\ell$ be the vertices with the property that $i\in
  \InputSet(v_j)$; order $v_1,\dots,v_\ell=\trunkvertex$ in decreasing distance to
  the output; i.e. $\InputSet(v_1)\subset\dots\subset
  \InputSet(v_\ell)$.  Then the $i^{th}$ component is given by
  \begin{equation}
    \label{eq:DefVarPhi}
    z_i=\sum_{j=1}^{\ell}\left( z^{(v_j)}_{i} /\prod_{k=j}^{\ell-1}
    r(e_k)\right),
    \end{equation}
    where $e_k$ is the edge from $v_k$ to $v_{k+1}$.
  
    We consider next the case where all the $\seq{r_j(e_k)}\to\infty$,
    to verify that the limit $\varphi(r_j,q)$ as $j\goesto \infty$ is
    the inclusion of $\lQ{T}\subset \lAssoc{n}{w}$.  For such
    sequences with $z_i$ is as in Equation~\eqref{eq:DefVarPhi}, we
    have that $\lim_{j\goesto\infty} z_i=z^{(\trunkvertex)}_i$. This shows that
    the $\seq{\z_j}$ is a sequence which has
    $(z_i^{(\trunkvertex)})_{i\in \InputSet(\trunkvertex)}$ as a stable limit.
    Next, we can pick any edge $e$ that points into
    the trunk vertex, and let $v_1$ denote its other vertex on
    $e$. Observe that $k=z^{(\trunkvertex)}_i$ is a real number that is independent
    of the choice of $i\in\InputSet(v_1)$.  Choose
    $\lambda_j\in\Aut(\HHH)$ defined by $\lambda_j(z)=(z-k)r_j(e)$.
    It is easy to see that
    \[
      \lim_{j\goesto\infty} \lambda_j(z_i)=
      \begin{cases}
        z_i^{(v_1)} &{\text{if $i\in\InputSet(v_1)$}} \\
        \infty &{\text{otherwise,}}
      \end{cases}
    \]
    exhibiting $(z_i^{(v_1)})_{i\in\InputSet(v_1)}$ as a stable limit.
    Proceeding in a similar manner for other vertices, we can find
    stable limits realizing $(z_i^{(v)})_{i\in\InputSet(v)}$
    for all vertices $v$, verifying that
    ${\overline\varphi}_{\{\infty\}^{E}\times
      \lQ{T}}$ is the inclusion of $\lQ{T}\subset \lAssoc{n}{w}$.
    
    To check that $\varphi$ is injective, 
    fix a  tree $T$, and let $E=\Edges(T)$.
    We construct a continuous map
    \[
      \psi_T\co \lQ{n,w}\to \lQ{T}
    \]
    as follows.  The input to $\psi_T$ is an element
    $\z=(z_1,\dots,z_{n+w})\in \lQ{n,w}$.  Fix a vertex $v$ of $T$, with
    $d$ incoming edges $e_1,\dots,e_d$ and weight $\omega$. For $j=1,\dots,d$,
    the edge $e_j$
    into $v$
    corresponds to a non-empty subset $S_j$ of the input set $\InputSet(v)$.
    Define $(Z_1,\dots,Z_{d+\omega})\in \prod_{i=1}^{d+\omega}\HHH$, as follows.
    Let $|S_j|$ denote the number of elements in $S_j$, and let
    \[
    Z_j=\sum_{i\in S_j}\frac{{\mathrm{Re}}(z_i)}{|S_j|}.
    \]
    If $k$ is any of the remaining elements $d+1,\dots,d+\omega$ of the input set for $v$
    let
    \[
    Z_k=z_k.
    \]
    The $v$ component of $\psi_T(\z)$ is
    the element $(z^{(v)}_1,\dots,z^{(v)}_{d+\omega})\in \lQ{d,\omega}$
    satisfying the slice condition that is in the same
    $\Aut(\HHH)$-orbit as $(Z_1,\dots,Z_{d+\omega})$.
    Letting $v$ vary, we obtain the desired element
    $\psi_T(z_1,\dots,z_{n+w})\in\lQ{T}$.

    We claim that $\psi_T(\varphi(r,q))=q$. To verify this identity,
    we argue as follows. Fix a vertex $v$ for the tree $T$.  Let
    $S_k\subset \InputSet(v)$ be the set coming from a parent vertex
    of $v$.  Equation~\eqref{eq:DefVarPhi} specifies the components of
    $\varphi(r,q)$. Each $z_i$ is a sum of terms, which are of the
    following three kinds: terms involving $z_i^{(v_j)}$, where
    $\InputSet(v_j)\subsetneq \InputSet(v)$; the
    term involving $z_i^{(v_j)}$ with $v_j=v$; and the terms involving
    $z_i^{(v_j)}$ with $\InputSet(v_j)\supsetneq\InputSet(v)$. We call the
    contributions of the first kind {\em ancestral
      contributions}; the contribution of $v$ the {\em
      contemporary contribution}; and terms of the third kind
    {\em descendant contributions}.  Observe that the descendant
    contributions are all the same, independent of $i\in\InputSet(v)$, and each
    contribution is real. We denote the sum of these contributions
    $K$.  Next, observe that in the definition of $Z_j$, the ancestral
    contributions cancel, because of the slice condition on the
    parents. Finally, for each of the $S_j\subset \InputSet(v)$, the
    contemporary contribution to $z_i$ is independent of the choice of
    $i\in S_j$.
    Let $d$ denote the number of inputs to $v$ and $\omega$ denote
    the weight of $v$.
    Comparing the labels $\{1,\dots, n+w\}$ of the factors
    of $\lQ{n,w}$ with the labels $\{1,\dots,d+\omega\}$
    of the corresponding factors of
    $\lQ{d,\omega}$ in $\lQ{T}$
    specifies a surjection
    \[
      \{1,\dots,n+w\}\supset \InputSet(v)\stackrel{\Pi}{\longrightarrow} \{1,\dots,d+\omega\}.
    \]
    Let $\sigma$ be a section of $\Pi$.
    For any $j\in\{1,\dots,d\}$, we have that
    \begin{equation}
      \label{eq:DetermineZ}
      Z_j=\frac{z_{\sigma(j)}^{(v)}}{R}  + K,
    \end{equation}
    where $R$ is the product of $r(e)$ over all edges from $v$ to the
    trunk vertex (in particular, it is independent of $j$).  For
    $j\in\{d+1,\dots,d+\omega\}$, there are no ancestral contributions, and
    $j$ is uniquely determined by $\sigma(j)$;
    Equation~\ref{eq:DetermineZ} holds for such $j$, as well.  It
    follows that $\Aut(\HHH)$ orbit of $(Z_1,\dots,Z_{d+\omega})$ contains
    $(z^{(v)}_{\sigma(1)},\dots,z^{(v)}_{\sigma(d+\omega)})$.  Since this
    holds for all $v\in T$, it follows that
    $q=\psi_T\bigl(\varphi(r,q)\bigr)$.

    To show $\varphi$ is
    injective, we extract the function $r$ from the element of
    $\varphi(r,q)$.  We determine $r$ on all edges inductively on the
    distance from the edge to the trunk vertex. Specifically, suppose
    that $e_1,\dots,e_\ell$ is an oriented path of edges to the trunk
    vertex, connecting vertices $v_1,\dots,v_{\ell+1}=\trunkvertex$.
    By stability, there is some $i$ so that $z_i^{(v_1)}\neq 0$.
    For that choice of $i$,
    Equation~\eqref{eq:DefVarPhi} determines $r(e_1)$ uniquely in
    terms of the real numbers $r(e_2),\dots,r(e_\ell)$, the complex
    numbers $z^{(v_1)}_i,\dots,z^{(v_\ell)}_i$ and $z_i$. Note that,
    for $i=1,\dots,\ell$, $z^{(v_i)}$ are specified by
    $q=\psi_T(\varphi(r,q))$, and $z_i$ a component of
    $\varphi(r,q)$.

    Since $\varphi$ is an injective map between Euclidean spaces of
    the same dimension, it follows from Invariance of
    Domain~\cite[Theorem~2B.3]{Hatcher02:book} that $\varphi$ is a
    homeomorphism onto its image.

    By construction, ${\overline\varphi}$ has the property that if $F\subset E$,
    then
    ${\overline\varphi}((0,\infty)^{F} \times \{\infty\}^{E\setminus F} \times
    \lQ{T}) \subset\lQ{K_F(T)}$.
    Indeed, the restriction of $\overline{\varphi}$ to
    $(0,\infty)^F\times\{\infty\}^{E\setminus F}\times\lQ{T}$
    is the Cartesian product of the maps
    obtained by applying $\varphi$ to the various components of
    $T\setminus (E\setminus F)$. 

    It follows that ${\overline\varphi}$ takes the boundary of
    $(0,\infty]^{E}\times \lQ{T}$ to the boundary of $\lAssoc{n}{w}$; and the
    description of ${\overline\varphi}$ in terms of $\varphi$ shows that this
    map, too, is a homeomorphism onto its image.  By construction,
    ${\overline\varphi}|_{(\infty)^{E}\times \lQ{T}}$ is the inclusion map of
    $\lQ{T}$ into $\lAssoc{n}{w}$.
\end{proof}

\begin{corollary}\label{cor:OrderedAssociaplexMfld}
  The space $\lAssoc{n}{w}$ is second-countable and Hausdorff,
  hence a topological manifold-with-boundary, the boundary of which is
  \begin{equation}\label{eq:ordered-assoc-top-bdy}
    \partial\lAssoc{n}{w}=\bigcup_{T\in\wTrees{n}{w}\setminus
      \{\wcorolla{n}{w}\}} Q^T.
  \end{equation}
  Further, $\lAssoc{n}{w}$ is compact.
\end{corollary}
\begin{proof}
  Lemma~\ref{lem:lAssociaplexCorners} gives a countable (in fact,
  finite) open cover of $\lAssoc{n}{w}$ by open subsets of
  $(0,\infty]^{n+2w-2}$. Hence, $\lAssoc{n}{w}$ is second-countable and locally
  Euclidean (in the sense of manifolds with boundary), and the boundary of $\lAssoc{n}{w}$ is given by
  Equation~\eqref{eq:ordered-assoc-top-bdy}.

  By Proposition~\ref{prop:sequential-top}, $\lAssoc{n}{w}$ is sequentially
  Hausdorff. For locally Euclidean spaces, sequential Hausdorffness implies
  Hausdorffness. Indeed, given points $\z$ and $\z'$, fix Euclidean charts
  around $\z$ and $\z'$. If the  balls $B_{1/n}(\z)$ and $B_{1/n}(\z')$
  intersect for all $n$, there is a sequence $\seq{\z_n}$ in their
  intersection converging to both $\z$ and $\z'$; so
  $\z=\z'$ by sequential
  Hausdorffness.

  For manifolds-with-boundary, compactness is equivalent to sequential
  compactness, which again was verified in Proposition~\ref{prop:sequential-top}. 
\end{proof}

Let
\[
  \Associaplex{n}{w}=\lAssoc{n}{w}/S_{w},
\]
the quotient by the symmetric group action permuting the labels of the
interior marked points.

\begin{proposition}\label{prop:AssociaplexMfld}
  The space $\Associaplex{n}{w}$ is a compact manifold-with-boundary, the boundary of which is
  \begin{equation}\label{eq:assoc-top-bdy}
    \partial\Associaplex{n}{w}=\bigcup_{T\in\wTrees{n}{w}\setminus
      \{\wcorolla{n}{w}\}} Q^T.
  \end{equation}
\end{proposition}
\begin{proof}
  First we verify that $\Associaplex{n}{w}$ is locally Euclidean, i.e., locally
  modeled on open subsets of $[0,1)^{n+2w-2}$.
    Recall that the symmetric product
  $\Sym^k(D^2)$ is homeomorphic to $D^{2k}$. It follows that the interior of
  each $Q^T$ is homeomorphic to an open ball. By
  Lemma~\ref{lem:lAssociaplexCorners}, a neighborhood of the interior of each
  $Q^T$ in $\lAssoc{n}{w}$ is given by $\lQ{T}\times(0,\infty]^k$. The action of
  the symmetric group on the $(0,\infty]^k$-factor is trivial: the group acts by
  permuting components of $\lQ{T}$ (each of which is a ball) and by acting on
  each component. Hence, the quotient space $\Associaplex{n}{w}$ is locally
  Euclidean.

  The space $\lAssoc{n}{w}$ is compact and the equivalence relation induced by
  the $S_w$-orbits is a closed subset of
  $\lAssoc{n}{w}\times\lAssoc{n}{w}$. Hence, the quotient space
  $\Associaplex{n}{w}$ is Hausdorff (see, e.g.,~\cite[Theorem
  4.57]{Lee11:top-mflds}).
  Of course, compactness is inherited by quotients, so by
  Proposition~\ref{prop:sequential-top}, $\Associaplex{n}{w}$ is
  compact. Similarly, second countability is inherited by quotients which are
  manifolds (see, e.g.,~\cite[Proposition 3.56]{Lee11:top-mflds}). Hence
  $\Associaplex{n}{w}$ is a manifold with boundary.  Since the action of $S_w$
  respects the decomposition in Equation~\eqref{eq:ordered-assoc-top-bdy} and,
  for each stratum of the boundary and each edge, the action of $S_w$
  respects the local collar neighborhoods from
  Lemma~\ref{lem:lAssociaplexCorners}, the boundary of $\Associaplex{n}{w}$ is given
  by Equation~\eqref{eq:assoc-top-bdy}.
\end{proof}

We will use the following folklore result:
\begin{lemma}\label{lem:ball}
  If $X$ is an $n$-dimensional compact topological
  manifold-with-boundary whose interior is homeomorphic to an open
  ball then $X$ is homeomorphic to a closed ball.
\end{lemma}
\begin{proof}
  It follows from Lefschetz duality that $H_i(X,\bdy X)=0$ for $i<n$
  and $H_n(X,\bdy X)\cong\ZZ$. Hence, the long exact sequence for the
  pair $(X,\bdy X)$ shows that $\bdy X$ is a homology sphere. 
  By the topological collar neighborhood theorem~\cite[Proposition
  3.42]{Hatcher02:book} (see also~\cite{Brown62:collar,Connelley71:collar}), the
  boundary of $X$ has a collar neighborhood.  It follows that $\partial X$ is
  simply-connected, since Euclidean space is simply-connected at infinity.
  Thus, $\bdy X$ is a homotopy sphere.  By the topological generalized
  Poincar{\'e}
  conjecture~\cite{Newman66:top-Poincare,Freedman82:Poincare,Perelman1,Perelman2,Perelman3}
  we obtain a homeomorphism from $\partial X$ to a sphere.  Further, the collar
  neighborhood gives an embedding of $X$ into the interior of $X$, i.e., into
  $\RR^n$. So, it follows from Brown's Schoenflies
  theorem~\cite{Brown60:Schoenflies} that $X$ is a closed ball.
\end{proof}

\begin{theorem}\label{thm:assoc-ball}
  The space $\Associaplex{n}{w}$ is homeomorphic to a closed ball,
  expressed as a $CW$ complex whose $k$-dimensional cells $Q^T$
  correspond to trees $T\in\wTreesCx[k]{n}{w}$ with $k=n+2w-2$, and
  whose cellular chain complex is identified with $\wTreesCx{n}{w}$ (over $\ZZ$).
\end{theorem}

\begin{proof}
  We prove by induction on the dimension $n+2w-2$ that
  $\Associaplex{n}{w}$ is a CW complex and that its cellular chain
  complex is $\wTreesCx{n}{w}$. For the base cases
  $(n,w)\in\{(2,0), (0,1)\}$, $\Associaplex{n}{w}$ is a single point.
  
  For the inductive step, by Proposition~\ref{prop:AssociaplexMfld} and the
  inductive hypothesis, Equation~\eqref{eq:assoc-top-bdy} expresses
  $\partial\Associaplex{n}{w}$ as a CW complex, with cells
  corresponding to the trees
  $T\in\wTrees{n}{w}\setminus\{\wcorolla{n}{w}\}$, and whose differential is
  identified with the corresponding subcomplex of $\wTreesCx{n}{w}$.
  Since the interior of $\Associaplex{n}{w}$ is an open ball,
  Lemma~\ref{lem:ball} implies that $\Associaplex{n}{w}$ is homeomorphic to
  a closed ball. In particular, the topological collar neighborhood theorem then
  implies that $\Associaplex{n}{w}$ is a CW complex.
  
  Finally, we verify that the cellular chain
  complex of $\Associaplex{n}{w}$ is isomorphic to
  $\wTreesCx{n}{w}$. Of course, both $\cellC{*}(\Associaplex{n}{w})$
  and $\wTreesCx{n}{w}$ have generating sets in bijection with the
  stably weighted trees with $n$ inputs and weight $w$.
  To verify the complexes are isomorphic we start by determining the signs in the
  identification between the generating sets. Fix arbitrarily an orientation of the (open)
  top-dimensional cell $D$ of $\Associaplex{n}{w}$. For any other
  $Q^T$, fix an ordering $\omega$ of the $e$ internal edges of $T$,
  i.e., an orientation of $T$ in the sense of
  Section~\ref{sec:w-signs}. By Lemma~\ref{lem:lAssociaplexCorners},
  this ordering induces an inclusion $(0,\infty)^e\times Q^T\into
  D$. The chosen orientation of $D$ and the standard orientation of
  $(0,\infty)^e$ then induce an orientation of $Q^T$. If we change
  $\omega$ by an even (respectively odd) permutation then this
  orientation is unchanged (respectively changed). Hence, we obtain a
  well-defined map
  $\cellC{*}(\Associaplex{n}{w})\stackrel{\cong}{\longrightarrow}
  \wTreesCx{n}{w}$, sending the generator corresponding to $T$ with
  orientation induced by $\omega$ to $(T,\omega)$.

  We must verify that this isomorphism is a chain map. From the
  description of the strata of $\Associaplex{n}{w}$, $Q^{T'}$ appears
  in the topological boundary of $Q^T$ if and only $T$ is obtained
  from $T'$ by collapsing an internal edge $e$. Fix such $T$, $T'$,
  and $e$, and let $\omega$ be an orientation of $T$. Then
  $e\wedge \omega$ is an orientation of $T'$. If we use $\omega$ to
  orient $Q^T$ and $e\wedge \omega$ to orient $Q^{T'}$ then, by
  Lemma~\ref{lem:lAssociaplexCorners} (taking $E=\{e\}$), the edge $e$ contributes $1$ to
  the coefficient of $Q^{T'}$ in the cellular boundary of $Q^T$. This
  agrees with the differential on $\wTreesCx{n}{w}$
  (Equation~\eqref{eq:signed-diff}), proving the result.
\end{proof}

%%% Local Variables: 
%%% mode: latex
%%% TeX-master: "AbstractDiagonal.tex"
%%% TeX-command-extra-options: "--shell-escape"
%%% End: 

\section{On boundedness}\label{sec:boundedness}
The conditions of being bonsai or bounded is used in the following places:
\begin{itemize}
\item For the type $D$ (respectively \DD) structure equation to make
  sense (Definitions~\ref{def:wD} and~\ref{def:wDD}), we need
  the algebra(s) to be bonsai or the type $D$ (or \DD) structure to be
  operationally bounded.
\item Similarly, to define the morphism complexes between type $D$
  (respectively \DD) structures and the composition maps on them
  (Section~\ref{sec:wD} and Definition~\ref{def:wDD}), we need the
  algebra(s) to be bonsai or the type $D$ (or \DD) structure to be
  operationally bounded.
\item For the box tensor products to be well-defined:
  \begin{itemize}
  \item For the box tensor product of a module and a type $D$
    structure to be defined, and an $\Ainf$-bifunctor, either the
    module must be bonsai or the type $D$ structure operationally
    bounded (Lemmas~\ref{lem:w-DT-sq-0} and~\ref{lem:w-DT-bifunc}).
  \item For the triple box tensor product
    $[\wMod_{\wAlg}\DT \lsup{\wAlg}P^{\wBlgop}\DT
    \lsub{\wBlgop}\wNod]_{\wModDiagNS}$
    (Definition~\ref{def:wtriple-prod} and
    Theorem~\ref{thm:wTripleTensorProduct}) to be well-defined, either
    the algebras and modules must be bonsai or the type \DD\ structure
    must be operationally bounded. (This is a special case of the
    previous point.)
  \item For the one-sided box tensor product
    $\wMod\DT^{\wTrPMDiagNS}P$ to be defined
    (Definition~\ref{def:w-one-sided-DT} and
    Proposition~\ref{prop:one-sided-DT-works}), 
    either the algebras and module must be bonsai or the type \DD\
    structure must be operationally bounded. Further, in the latter case, the
    resulting type $D$ structure is operationally bounded.
  \item For associativity of the box tensor product
    (Proposition~\ref{prop:one-sided-DT-works}), either all the
    algebras and modules involved must be bonsai or the type \DD\
    structure must be operationally bounded.
  \item For functoriality of the one-sided box tensor products
    (Lemma~\ref{lem:w-one-side-DT-id-chain-map} and
    Proposition~\ref{prop:w-one-side-DT-bifunc}), either the type \DD\
    structure must be operationally bounded or the algebras, modules,
    and module morphisms must be bonsai.
  \item For associativity of the box tensor product of morphisms
    (Lemma~\ref{lem:w-part-mod-map-diag-DT} and
    Corollary~\ref{cor:w-morph-DT-assoc}), again either the type \DD\
    structure must be operationally bounded or the algebras, modules,
    and module morphism must be bonsai.
  \end{itemize}
\end{itemize}

Further, these conditions need to satisfy:
\begin{itemize}
\item The (external) tensor product of two bonsai weighted algebras,
  modules, algebra homomorphisms, or module morphisms is bonsai
  (Lemmas~\ref{lem:wADtp-grading} and~\ref{lem:wMDtp-grading}).
\item A weakly unital bonsai algebra is isomorphic to a strictly
  unital bonsai algebra via a bonsai isomorphism
  (Theorem~\ref{thm:UnitalIsUnitalW} and
  Proposition~\ref{prop:w-qi-unital}).
\item A weakly unital bonsai module is isomorphic to a strictly unital
  bonsai module via a bonsai isomorphism
  (Theorem~\ref{thm:UnitalIsUnitalMw}).
\item The complex of bonsai morphisms between strictly unital bonsai
  modules is homotopy equivalent to the subcomplex of strictly unital
  bonsai morphisms (Proposition~\ref{prop:w-mod-mor-cx-unital}).
\end{itemize}

For our application to bordered Floer theory, we will need to weaken
these hypotheses slightly.

\begin{definition}\label{def:filtered-wAlg}
  Let $\wAlg=(A,\{\mu_n^w\})$ be a weighted $\Ainf$-algebra. A
  \emph{filtration} of $\wAlg$ is a sequence of subspaces
  $A=\Filt^0A\supset\Filt^1A\supset\Filt^2A\supset\cdots$ of $A$ so that
  if $a_i\in\Filt^{m_i}A$, $i=1,\dots,n$ then
  \[
    \mu_n^w(a_1,\dots,a_n)\in\Filt^{m_1+\cdots+m_n}A.
  \]
  The algebra $\wAlg$ is \emph{complete} with respect to the
  filtration $\Filt$ if the vector space $A$ is complete with respect
  to $\Filt$.

  A filtration of a weighted $\Ainf$-module $\wMod$ over a filtered,
  weighted $\Ainf$-algebra $\wAlg$ is a sequence of subspaces
  $M=\Filt^0M\supset\Filt^1M\supset\cdots$ so that if
  $x\in \Filt^{m_0}M$ and $a_i\in \Filt^{m_i}A$ then
  \[
    \mu_{1+n}^w(x,a_1,\dots,a_n)\in\Filt^{m_0+\cdots+m_n}M.
  \]
  The module $\wMod$ is complete with respect to the filtration if the
  underlying vector space $M$ is.

  If $(\wAlg,\Filt)$ and $(\wBlg,\Filt)$ are filtered weighted
  $\Ainf$-algebras then 
  a homomorphism $f\co \wAlg\to\wBlg$ is
  \emph{filtered} if for all $a_1,\dots,a_n\in A$ with $a_i\in \Filt^{m_i}A$,
  \[
    f_n^w(a_1,\dots,a_n)\in\Filt^{m_1+\cdots+m_n}B.
  \]

  If $(\wMod,\Filt)$ and $(\wNod,\Filt)$ are filtered weighted
  $\Ainf$-modules over $(\wAlg,\Filt)$ then a morphism 
  $f\co \wMod\to\wNod$ is \emph{filtered} if for all
  $x\in \Filt^{m_0}M$ and $a_i\in \Filt^{m_i}A$ we have
  \[
    f_{1+n}^w(x,a_1,\dots,a_n)\in\Filt^{m_0+\cdots+m_n}N.
  \]
\end{definition}

It is immediate from the definitions that the filtered morphisms form
a subcomplex $\Mor\bigl((\wMod,\Filt),(\wNod,\Filt)\bigr)$ of
$\Mor(\wMod,\wNod)$. Note also that if $T$ is a stably weighted tree
with $n$ inputs and $(\wAlg,\Filt)$ is a filtered $\Ainf$-algebra then
it follows from the definitions that
\[
  \mu(T)\co \Filt^{m_1}A\kotimes{\Ground}\cdots\kotimes{\Ground} \Filt^{m_n}A\to \Filt^{m_1+\cdots+m_n}A
\]
(and similar statements for modules, homomorphisms, and
morphisms).

\begin{example}
  If $\wAlg$ is a weighted $\Ainf$-algebra over the power series ring
  $\FF_2[[U]]$, and $A$ is a finitely-generated free module over
  $\FF_2[[U]]$, then there is a
  $U$-power filtration on $\wAlg$, $\Filt^m=U^mA$, and $\wAlg$ is
  complete with respect to this filtration.
\end{example}

Observe that if $(\wAlg,\Filt)$ is a filtered weighted $\Ainf$-algebra
then each quotient $A/\Filt^mA$ inherits the structure of a weighted
$\Ainf$-algebra $\wAlg/\Filt^m\wAlg$; and if $(\wMod,\Filt)$ is a
filtered weighted $\Ainf$-module over $(\wAlg,\Filt)$ then each
quotient $M/\Filt^mM$ inherits the structure of a weighted
$\Ainf$-module $\wMod/\Filt^m\wMod$ over
$\wAlg/\Filt^m\wAlg$. Homomorphisms (respectively morphisms) of
filtered weighted algebras (respectively modules) induce homomorphisms
(respectively morphisms) of these quotients

\begin{definition}\label{def:filtered-bonsai}
  \begin{itemize}
  \item A filtered weighted $\Ainf$-algebra $(\wAlg,\Filt)$ is
    \emph{filtered bonsai} if $\wAlg$ is complete with respect to
    $\Filt$ and each quotient $\wAlg/\Filt^m\wAlg$ is bonsai.
  \item A filtered weighted $\Ainf$-module $(\wMod,\Filt)$ over a
    filtered bonsai $(\wAlg,\Filt)$ is \emph{filtered bonsai} if
    $\wMod$ is complete with respect to $\Filt$ and each
    $\wMod/\Filt^m\wMod$ is bonsai.
  \item A filtered homomorphism $f\co (\wAlg,\Filt)\to(\wBlg,\Filt)$
    between filtered bonsai weighted $\Ainf$-algebras is
    \emph{filtered bonsai} if the induced map $\wAlg/\Filt^m\wAlg\to
    \wBlg/\Filt^m\wBlg$ is bonsai for each $m$.
  \item A filtered morphism $f\co (\wMod,\Filt)\to(\wNod,\Filt)$
    between filtered bonsai weighted $\Ainf$-modules is
    \emph{filtered bonsai} if the induced map $\wMod/\Filt^m\wMod\to
    \wNod/\Filt^m\wNod$ is bonsai for each $m$.
  \end{itemize}
\end{definition}

\begin{definition}\label{def:filtered-bounded}
  A weighted type $D$ or \DD\ structure with charge $X$ over a
  filtered weighted $\Ainf$-algebra $(\wAlg,\Filt)$ is \emph{filtered
    operationally bounded} if the operation $\delta^n$ vanishes for
  all sufficiently large $n$ and the action of $X$ on
  $\wAlg/\Filt^m\wAlg$ is nilpotent for each $m$ (i.e., for each $m$
  there exists an $N$ so that for $n>N$ and any $a\in A$,
  $X^n\cdot a\in\Filt^mA$).
\end{definition}

\begin{lemma}
  If $(\wAlg=(A,\{\mu_n^w\}),\Filt)$ and
  $(\wBlg=(B,\{\nu_n^w\}),\Filt)$ are filtered bonsai algebras then so
  is $\wAlg\wADtp\wBlg$. The analogous statement holds for filtered
  bonsai modules, filtered bonsai algebra homomorphisms, and filtered
  bonsai module morphisms.
\end{lemma}
\begin{proof}
  We explain the case of filtered bonsai algebras; the remaining cases
  are similar. Define the filtration on $A\rotimes{\Ring} B$ by declaring that
  \[
    \Filt^n(A\rotimes{\Ring} B)=\sum_{p+q=n}\Filt^p(A)\rotimes{\Ring}\Filt^q(B).
  \]
  To see this defines a filtration on $\wAlg\wADtp\wBlg$, note that if
  $T$ is a stably weighted trees with $n$ inputs then
  $\mu\otimes_{\wDiagNS}\nu(T)$ is a linear combination of pairs of
  trees with $n$ inputs, so defines a map
  \[
    \bigl(\Filt^{p_1}A\kotimes{\Ground_1}\cdots\kotimes{\Ground_1}\Filt^{p_n}A\bigr)\rotimes{\Ring}\bigl(\Filt^{q_1}B\otimes\cdots\otimes\Filt^{q_n}B\bigr)\to
    \bigl(\Filt^{p_1+\cdots+p_n}A\bigr)\rotimes{\Ring}\bigl(\Filt^{q_1+\cdots+q_n}B\bigr)\subset
    \Filt^{p_1+\cdots+p_n+q_1+\cdots+q_n}(A\rotimes{\Ring} B).
  \]
  For each $m$ and $p$, the composition
  \[
    (A\rotimes{\Ring} B)^{\kotimes{\Ground} n}\xrightarrow{\mu_n(T)}
    (A\rotimes{\Ring} B)\to (A\rotimes{\Ring} B)/\Filt^m(A\rotimes{\Ring} B)\to (A/\Filt^pA)\rotimes{\Ring}(B/\Filt^{m-p}B)
  \]
  vanishes for $\dim(T)$ sufficiently large. It follows that for each
  $m$, the map 
  \[
    (A\rotimes{\Ring} B)^{\kotimes{\Ground} n}\xrightarrow{\mu_n(T)}
    (A\rotimes{\Ring} B)\to (A\rotimes{\Ring} B)/\Filt^m(A\rotimes{\Ring} B)
  \]
  vanishes for $\dim(T)$ sufficiently large. Hence, $\wAlg\wADtp\wBlg$
  is filtered bonsai, as desired.
\end{proof}

The results about unitality also hold in the filtered bonsai case:
\begin{lemma}
  The bonsai statements in Theorem~\ref{thm:UnitalIsUnitalW},
  Proposition~\ref{prop:w-qi-unital},
  Theorem~\ref{thm:UnitalIsUnitalMw}, and
  Proposition~\ref{prop:w-mod-mor-cx-unital} hold for filtered bonsai
  weighted algebras, modules, and algebra homomorphisms, and module
  morphisms.
\end{lemma}
\begin{proof}
  We explain (briefly) the case of Theorem~\ref{thm:UnitalIsUnitalW};
  the other cases are similar. If suffices to check that each
  modification of the weighted operations respects the filtration, and
  the induced map of quotients $A/\Filt^mA$ is bonsai. The first
  statement follows from the fact that each $\phi_n^v$ and $\psi_n^v$
  used to perform the modification is filtered, and composition of
  filtered maps are filtered. The second follows from the fact that
  the new multiplication is a linear combination of operation trees
  for the old multiplication of the same dimension, with some inputs
  $\unit$.
\end{proof}

\begin{proposition}
  Throughout this paper, the hypothesis of being bonsai can be
  replaced by the hypothesis of being filtered bonsai, and the
  hypothesis of being operationally bounded can be replaced by the
  hypothesis of being filtered operationally bounded.
\end{proposition}
\begin{proof}
  This is straightforward, and is left to the reader.
\end{proof}

%%% Local Variables:
%%% mode: latex
%%% TeX-master: "AbstractDiagonal"
%%% End:

% \input{bialgebras}
\appendix
\newpage
\raggedbottom
\section{Terms in diagonals}\label{sec:pictures}
For the reader's convenience, we include the first few terms in the weighted diagonals of various types introduced in the text. To save space, we draw weight-0 vertices without the circle labelled 0.

\begin{center}
  %Font is 18 point.
  \includegraphics[scale=.55556]{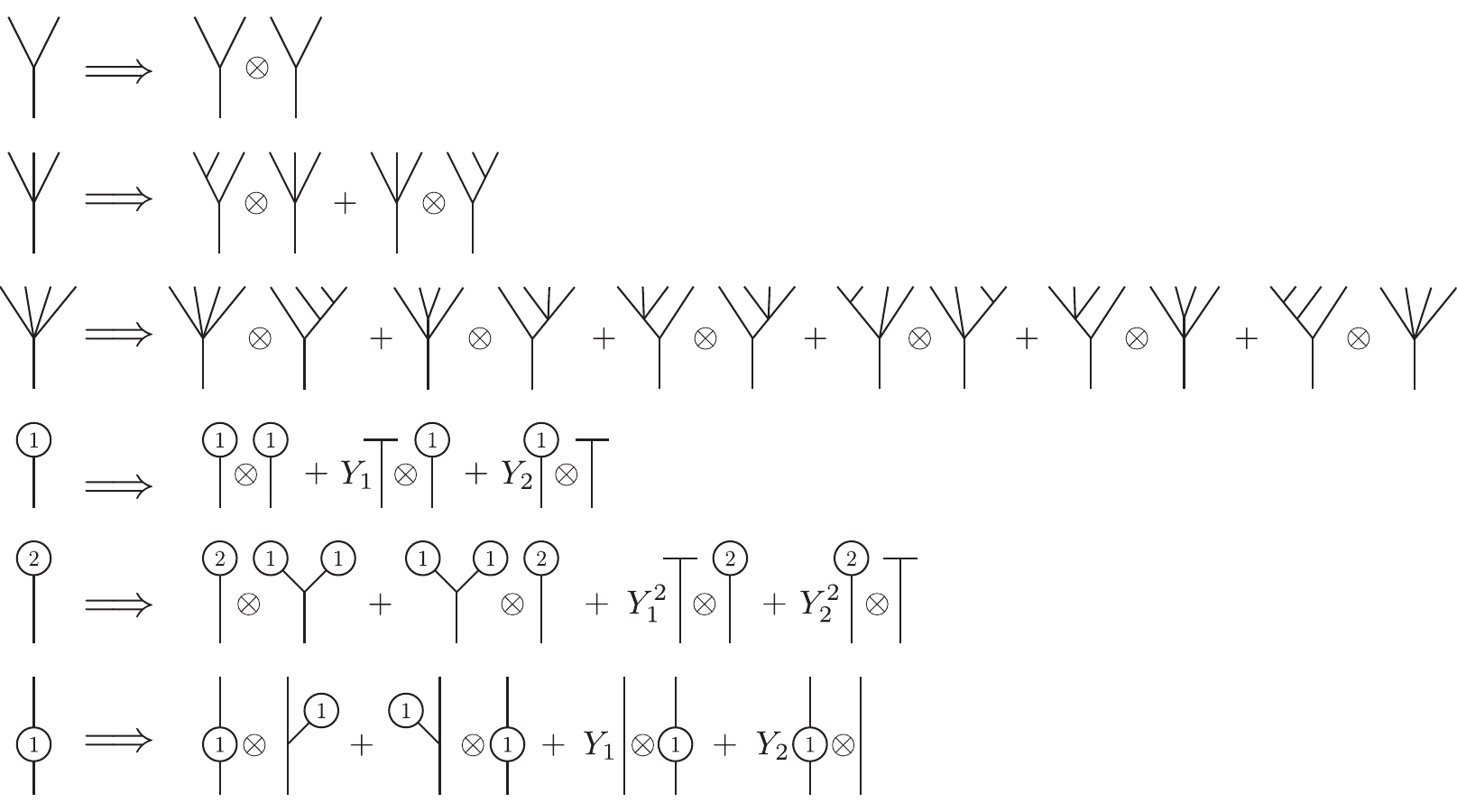}\\
  \textsc{A weighted algebra diagonal $\wDiag{*}{*}$.}
\end{center}
\bigskip\bigskip
\begin{center}
  %Font is 18 point.
  \includegraphics[scale=.55556]{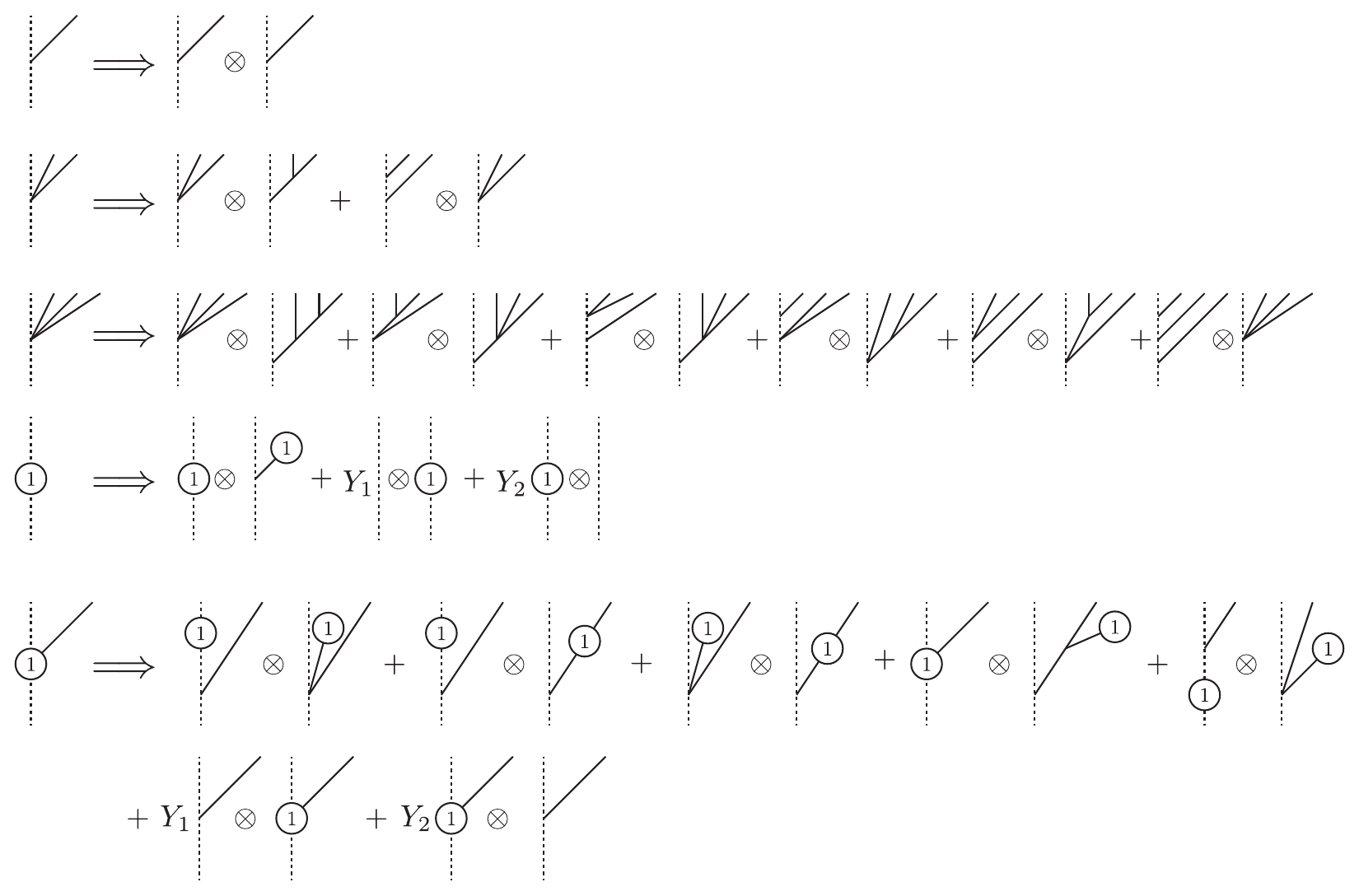}\\  
  \textsc{A weighted module diagonal $\wModDiag{*}{*}$ compatible with $\wDiag{*}{*}$.}
\end{center}
\bigskip\bigskip
\begin{center}
  %Font is 18 point.
  \includegraphics[scale=.55556]{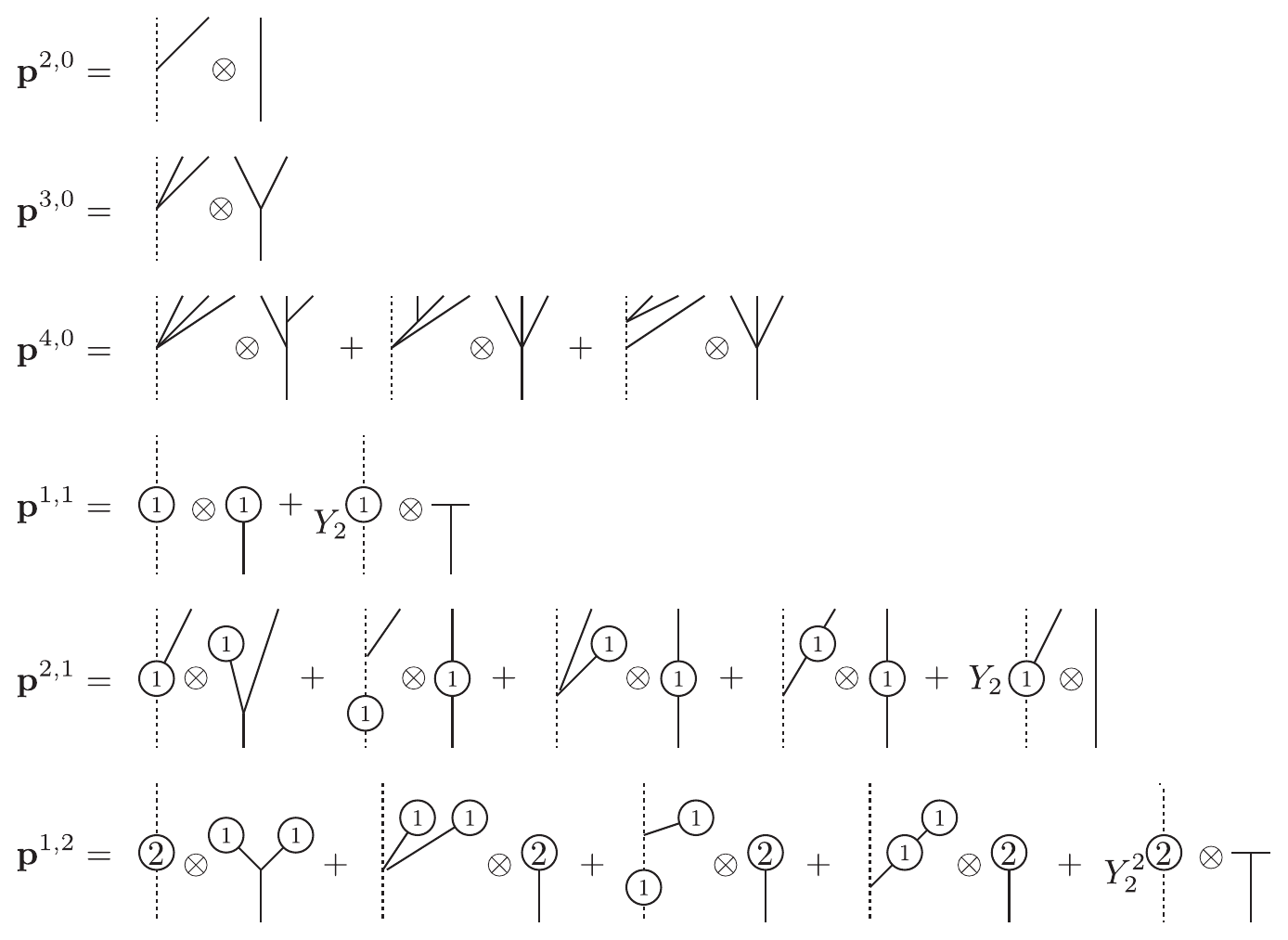}\\  
  \textsc{A weighted module diagonal primitive $\wTrPMDiag{*}{*}$ compatible with $\wDiag{*}{*}$ and inducing $\wModDiag{*}{*}$.}
\end{center}
\bigskip\bigskip
\begin{center}
  %Font is 12 point.
  \includegraphics[scale=.833333]{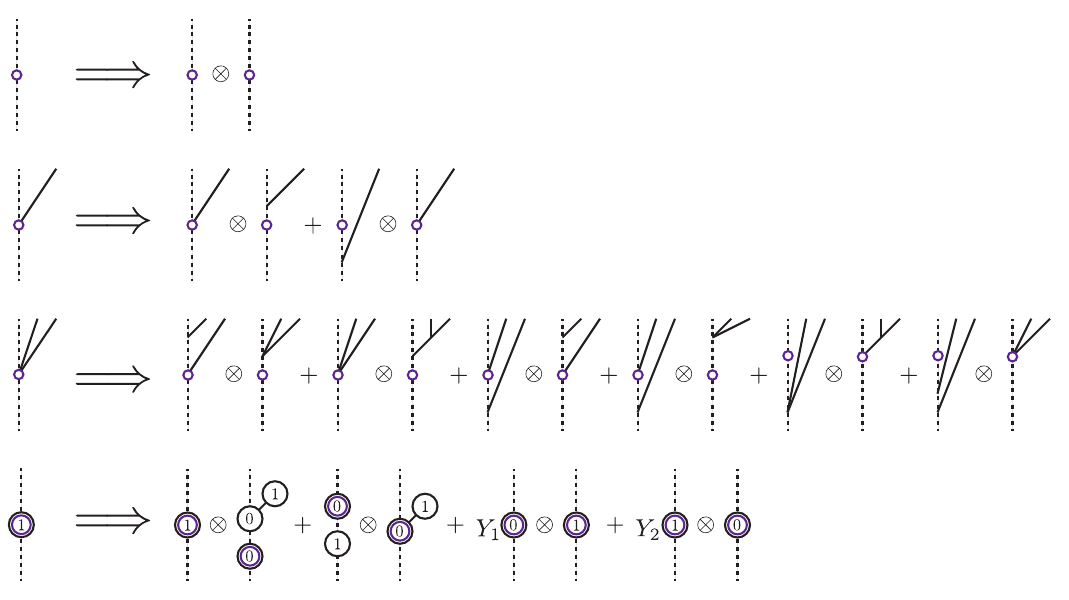}\\  
  \textsc{A weighted module-map diagonal $\wModMulDiag{*}{*}$ compatible with $\wModDiag{*}{*}$ and $\wModDiag{*}{*}$.}
\end{center}
\bigskip\bigskip
\begin{center}
  %Font is 12 point.
  \includegraphics[scale=.8333333]{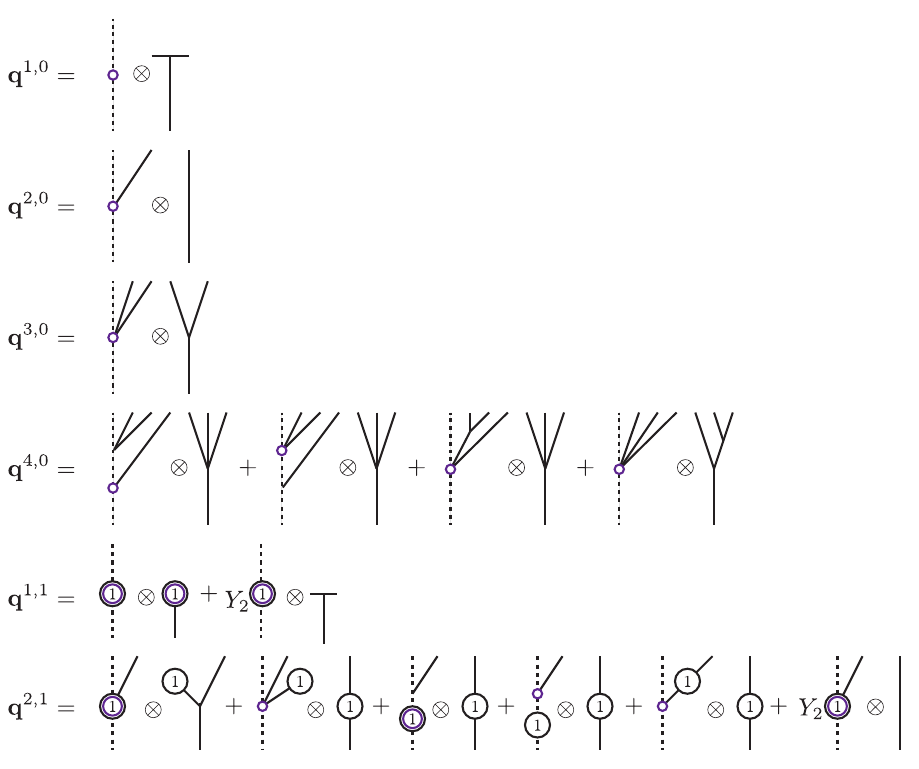}\\  
  \textsc{A weighted module-map primitive $\wTrPMorDiag{*}{*}$ compatible with $\wTrPMDiag{*}{*}$ and $\wTrPMDiag{*}{*}$.}
\end{center}
\bigskip\bigskip
\begin{center}
  %Font is 10 point.
  \includegraphics{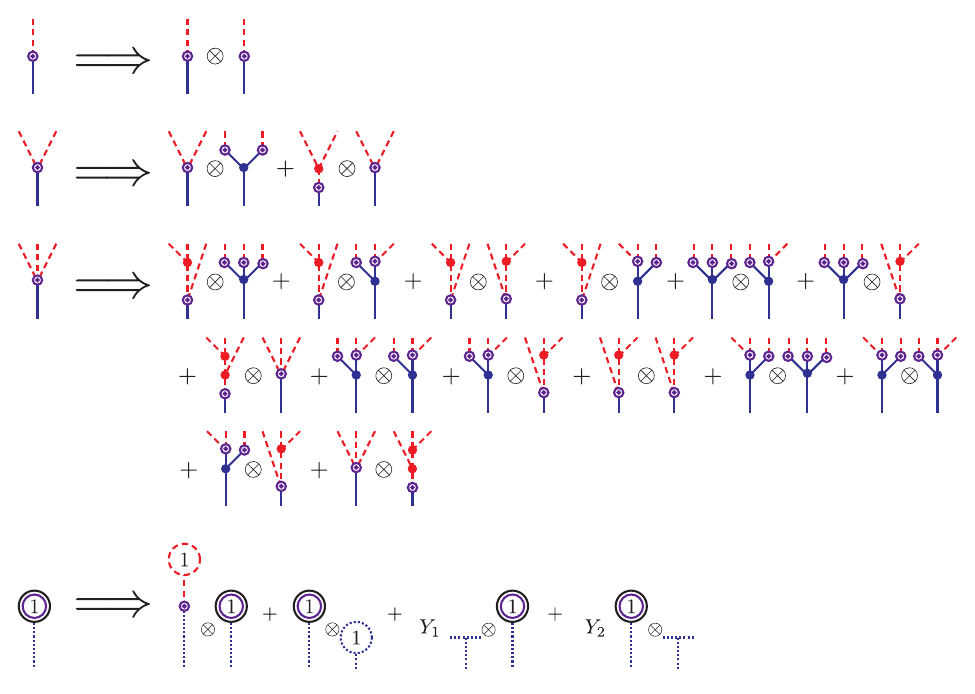}\\  
  \textsc{A weighted (algebra) map diagonal $\wMulDiag{*}{*}$ compatible with $\wDiag{*}{*}$ and $\wDiag{*}{*}$.}
\end{center}
\bigskip\bigskip
\begin{center}
  %Font is 12 point.
  \includegraphics[scale=.8333333]{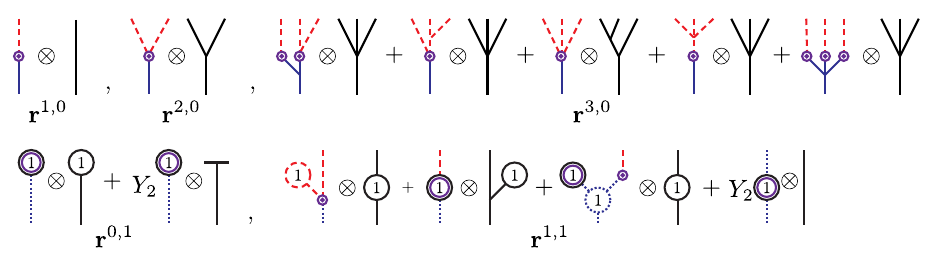}\\  
  \textsc{A weighted DADD diagonal $\wTrDADD{*}{*}$ compatible with $\wDiag{*}{*}$ and $\wDiag{*}{*}$.}
\end{center}

%%% Local Variables: 
%%% mode: latex
%%% TeX-master: "AbstractDiagonal.tex"
%%% End: 

\section{A convention of tables}\label{sec:tables}

%\begin{table}
{\small
\begin{center}
  \begin{tabular}{p{1.5in}lllp{2.5in}}
    \toprule
    Object & Notation & Definition & Terms & Use \\
    \midrule
    Associahedron diagonal & $\AsDiag^n$, $\TrDiag_n$ & \ref{def:AssociahedronDiagonal},~\ref{def:DiagonalCells} & Fig.~\ref{fig:diag-cells} & Tensor product of $\Ainf$-algebras, definition of \DD\ structure. \\
    Module diagonal & $\MDiag^n$, $\TrMDiag_n$ & \ref{def:ModuleDiagonal} & Fig.~\ref{fig:diag-cells} & Tensor product of $\Ainf$-modules, triple box tensor product.\\
    Module diagonal primitive & $\TrPMDiag_n$ & \ref{def:M-prim} & Fig.~\ref{fig:mprim-terms} & One-sided box tensor product.
    \\
    Multiplihedron diagonal & $\MulDiag^n$, $\TrMulDiag_n$ &
                                                             \ref{def:multiplihedron-diag} & Fig.~\ref{fig:mult-diag} & Tensor product of morphisms of $\Ainf$-algebras.\\
    Module-map diagonal & $\ModMulDiag_n$, $\TrModMulDiag_n$ & \ref{def:mod-map-diag} & Fig.~\ref{fig:m-tree-diag} & Tensor product of morphisms of $\Ainf$-modules.\\
    Module-map primitive & $\TrPMorDiag_n$ & \ref{def:mod-map-prim} & Fig.~\ref{fig:m-map-prim} & One-sided box product of identity map of a \DD\ structure with a morphism of $\Ainf$-modules.\\
    Partial module-map diagonal & $\PartTrModMulDiag_n$ & \ref{def:part-mod-map-diag} & & Tensor product of a morphism of $\Ainf$-modules with the identity map of another $\Ainf$-module.
    \\
    \midrule
    Weighted algebra diagonal & $\wDiag{n}{w}$, $\wDiagCell{n}{w}$& \ref{def:w-alg-diag},~\ref{def:wDiagCells} & Fig.~\ref{fig:wdiag-terms} & Tensor product of $w$-algebras, definition of weighted \DD\ structure.\\
    Weighted module diagonal & $\wModDiag{n}{w}$,
                               $\wModDiagCell{n}{w}$ &
                                                       \ref{def:w-mod-diag} & Fig.~\ref{fig:wmod-diag-terms} & Tensor product of $w$-modules, triple box tensor product of a weighted \DD\ structure and two $w$-modules.\\
    Weighted module diagonal primitive & $\wTrPMDiag{n}{w}$ &
                                                              \ref{def:wM-prim} & Fig.~\ref{fig:wmod-prim-diag-terms} & One-sided box tensor product of a weighted \DD\ structure and a $w$-module.\\
  Weighted map diagonal  & $\wMulDiag{n}{w}$, $\wTrMulDiag{n}{w}$ &
                                                                    \ref{def:w-map-diag},~\ref{def:w-map-diag-cell} & Fig.~\ref{fig:w-map-diag-base-case} & Tensor product of morphisms of $w$-algebras.\\
    Weighted module-map diagonal & $\wModMulDiag{n}{w}$,
                                   $\wTrModMulDiag{n}{w}$ &
                                                            \ref{def:w-mod-map-diag} & Fig.~\ref{fig:w-mod-map-diag-eg} & Tensor product of morphisms of $w$-modules.\\
    Weighted module-map primitive  & $\wTrPMorDiag{n}{w}$ & \ref{def:w-mod-map-prim} & Fig.~\ref{fig:w-map-prim}& One-sided box product of identity map of a \DD\ structure with a morphism of $\Ainf$-modules.\\
    Partial weighted module-map diagonal & $\wPartTrModMulDiag{n}{w}$ & \ref{def:w-partial-mod-map-diag} & & Tensor product of a morphism of $w$-modules with the identity map of another $w$-module.\\
    \bottomrule\\
  \end{tabular}
\end{center}}
\begin{center}
  \textsc{Table 1.} Things and their uses.
\end{center}
% \caption{Things and their uses.}
% \end{table}

% \begin{table}
\begin{center}
  \begin{tabular}{ll}
    \toprule
    Convention & Summary \\ \midrule
    \ref{conv:diag-tens-F2},~\ref{conv:wDiags-tens-Ring}
    & Tensor products in Sections~\ref{sec:diags} and~\ref{sec:wDiags} are over
                                         $\Ring$ by default.\\
    \ref{conv:Ring}, %~\ref{conv:Ground},
    \ref{conv:wDiagApps-ground}
    &The ground ring notation in
      Sections~\ref{sec:algebra},~\ref{sec:wAinfty},~\ref{sec:wDiagApps}.
    Tensor products over $\Ground$ by default.\\
    \ref{conv:trees} & Terminology for trees.\\
    % \ref{conv:leaves} & Leaves of weighted trees are viewed
    %                                as vertices.\\
    \ref{conv:cellC} & Cellular chain complex are with $\Ring$-coefficients.\\
    \ref{conv:grading-shift} & Grading shifts.\\
    \ref{conv:d-bounded} & Boundedness hypotheses for type $D$ structures.\\
    \ref{conv:uncurved} & Weighted $\Ainf$-algebras are uncurved.\\
    \ref{conv:su-mods} & Modules over strictly unital weighted algebras assumed strictly unital.\\
    \ref{conv:boundedness-suppressed} & Boundedness hypotheses suppressed in the rest of this section.\\
    \bottomrule
  \end{tabular}
\end{center}
\begin{center}
  \textsc{Table 2.} Table of conventions, to help the reader find which
  conventions are in force.
\end{center}
% \caption{Table of conventions, to help the reader find which
%   conventions are in force.}
% \end{table}

%%% Local Variables:
%%% mode: latex
%%% TeX-master: "AbstractDiagonal"
%%% End:

% \input{extraneous}

\bibliographystyle{hamsalpha}\bibliography{heegaardfloer}
\end{document}